%% file: main.tex
\documentclass[
12pt, % The default document font size, options: 10pt, 11pt, 12pt
english, % ngerman for German
singlespacing, % Single line spacing, alternatives: onehalfspacing or doublespacing
headsepline, % Uncomment to get a line under the header
]{MastersDoctoralThesis} % The class file specifying the document structure

\usepackage[utf8]{inputenc} % Required for inputting international characters
\usepackage[T1]{fontenc} % Output font encoding for international characters
\usepackage{subcaption}
\usepackage{mathpazo} % Use the Palatino font by default
\usepackage{enumerate}
\usepackage[export]{adjustbox}
\usepackage{mathtools}

\usepackage{epigraph}

\usepackage[backend=bibtex,natbib=true,style=alphabetic]{biblatex} % Use the bibtex backend with the authoryear citation style (which resembles APA)
\usepackage[autostyle=true]{csquotes} % Required to generate language-dependent quotes in the bibliography

\usepackage{amsmath,amssymb,amsthm,latexsym,bbm}
\usepackage{stmaryrd}

\theoremstyle{definition}
\newtheorem{lemma}{Lemma}
\newtheorem{conjecture}{Conjecture}
\newtheorem{definition}{Definition}
\newtheorem{theorem}{Theorem}
\newtheorem{remark}{Remark}

\newtheorem{proposition}{Proposition}

\newtheorem{question}{Question}

\newcommand{\normord}[1]{\vcentcolon\mathrel{#1}\vcentcolon}
\providecommand{\vcentcolon}{\mathrel{\mathop{:}}}

\DeclareMathOperator{\Tr}{Tr}

\DeclareMathOperator{\sgn}{sgn}
\DeclareMathOperator{\id}{id}
\DeclareMathOperator{\Cat}{Cat}

\newcommand{\cG}{\mathcal{G}}
\newcommand{\cT}{\mathcal{T}}

\newcommand{\cS}{\mathcal{S}}
\newcommand{\cI}{\mathcal{I}}
\newcommand{\cC}{\mathcal{C}}

\newcommand{\TJ}[6]{ \begin{pmatrix}
  #1 & #2 & #3 \\
  #4 & #5 & #6 
\end{pmatrix}}

\newcommand{\SJ}[6]{ \begin{Bmatrix}
  #1 & #2 & #3 \\
  #4 & #5 & #6 
 \end{Bmatrix}}

\thesistitle{D\'{e}compositions combinatoire pour des familles de cartes d\'{e}form\'{e}es ou d\'{e}cor\'{e}es} % Your thesis title, this is used in the title and abstract, print it elsewhere with \ttitle
\supervisor{Dr. James \textsc{Smith}} % Your supervisor's name, this is used in the title page, print it elsewhere with \supname
\examiner{} % Your examiner's name, this is not currently used anywhere in the template, print it elsewhere with \examname
\degree{Doctor of Philosophy} % Your degree name, this is used in the title page and abstract, print it elsewhere with \degreename
\author{John \textsc{Smith}} % Your name, this is used in the title page and abstract, print it elsewhere with \authorname
\addresses{} % Your address, this is not currently used anywhere in the template, print it elsewhere with \addressname

\subject{Biological Sciences} % Your subject area, this is not currently used anywhere in the template, print it elsewhere with \subjectname
 % Keywords for your thesis, this is not currently used anywhere in the template, print it elsewhere with \keywordnames
\university{\href{http://www.university.com}{University Name}} % Your university's name and URL, this is used in the title page and abstract, print it elsewhere with \univname
\department{\href{http://department.university.com}{Department or School Name}} % Your department's name and URL, this is used in the title page and abstract, print it elsewhere with \deptname
\group{\href{http://researchgroup.university.com}{Research Group Name}} % Your research group's name and URL, this is used in the title page, print it elsewhere with \groupname
\faculty{\href{http://faculty.university.com}{Faculty Name}} % Your faculty's name and URL, this is used in the title page and abstract, print it elsewhere with \facname

\AtBeginDocument{
\hypersetup{pdftitle=\ttitle} % Set the PDF's title to your title
\hypersetup{pdfauthor=\authorname} % Set the PDF's author to your name
\hypersetup{pdfkeywords=\keywordnames} % Set the PDF's keywords to your keywords
\hypersetup{linkcolor=Blue}
\hypersetup{citecolor=Cerulean}
}

\begin{document}

\frontmatter % Use roman page numbering style (i, ii, iii, iv...) for the pre-content pages

\pagestyle{plain} % Default to the plain heading style until the thesis style is called for the body content

%----------------------------------------------------------------------------------------
%	TITLE PAGE
%----------------------------------------------------------------------------------------

\include{Chapters/UBX_Garde}

%----------------------------------------------------------------------------------------
%	ABSTRACT PAGE
%----------------------------------------------------------------------------------------

\include{Chapters/UBX_Resume}

%----------------------------------------------------------------------------------------
%	ACKNOWLEDGEMENTS
%----------------------------------------------------------------------------------------
\include{Chapters/Remerciements}

%\begin{acknowledgements}
%\addchaptertocentry{\acknowledgementname} % Add the acknowledgements to the table of contents
%The acknowledgments and the people to thank go here, don't forget to include your project advisor\ldots
%\end{acknowledgements}

%----------------------------------------------------------------------------------------
%	LIST OF CONTENTS/FIGURES/TABLES PAGES
%----------------------------------------------------------------------------------------

\setcounter{tocdepth}{2}
\tableofcontents % Prints the main table of contents

%\listoffigures % Prints the list of figures

%\listoftables % Prints the list of tables

%----------------------------------------------------------------------------------------
%	SYMBOLS
%----------------------------------------------------------------------------------------

%\begin{symbols}{lll} % Include a list of Symbols (a three column table)

%$a$ & distance & \si{\meter} \\
%$P$ & power & \si{\watt} (\si{\joule\per\second}) \\
%Symbol & Name & Unit \\

%\addlinespace % Gap to separate the Roman symbols from the Greek

%$\omega$ & angular frequency & \si{\radian} \\

%\end{symbols}

%----------------------------------------------------------------------------------------
%	DEDICATION
%----------------------------------------------------------------------------------------

%\dedicatory{For/Dedicated to/To my\ldots} 

%----------------------------------------------------------------------------------------
%	THESIS CONTENT - CHAPTERS
%----------------------------------------------------------------------------------------

\mainmatter % Begin numeric (1,2,3...) page numbering

\pagestyle{thesis} % Return the page headers back to the "thesis" style

\include{Chapters/Intro_VF}
\include{Chapters/Intro}
\include{Chapters/Matrix_model}

\include{Chapters/Constellations} 
\include{Chapters/Tensor_model}

\include{Chapters/Double_scaling}

\include{Chapters/Melons_tens_vect}
\include{Chapters/Conclusion}

%----------------------------------------------------------------------------------------
%	THESIS CONTENT - APPENDICES
%----------------------------------------------------------------------------------------

\appendix % Cue to tell LaTeX that the following "chapters" are Appendices

% Include the appendices of the thesis as separate files from the Appendices folder
% Uncomment the lines as you write the Appendices

\include{Appendices/App_cubic_constraint}
\include{Appendices/App_UN2_OD}
\include{Appendices/App_prop_3j}

%----------------------------------------------------------------------------------------
%	BIBLIOGRAPHY
%----------------------------------------------------------------------------------------

\printbibliography[heading=bibintoc]

%----------------------------------------------------------------------------------------

\end{document}

%% file: Chapters/UBX_Garde.tex
\pagestyle{empty}
\begin{center}
    % Les logos en tête, ici une seule image incluse
%\includegraphics[scale=1, height=1.7cm]{images/logo0.png}
%\hfill
%\includegraphics[scale=1, height=1.7cm]{images/logo1.png}
\includegraphics[scale=1, height=1.7cm]{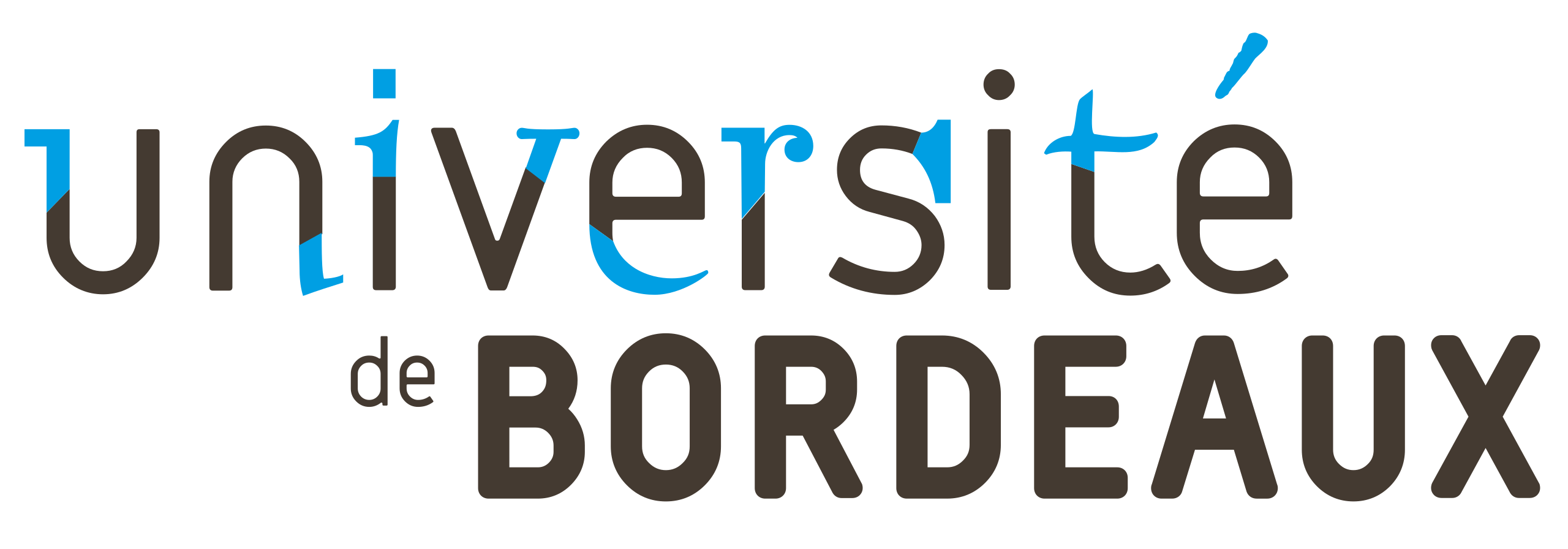}
%\includegraphics[scale=1, height=1.7cm]{Figures/brdx.png}
%\hfill
%\includegraphics[scale=1, height=1.7cm]{images/logo4.jpg}
\end{center}
\begin{center}
\doublespacing
\begin{large}

THÈSE PRÉSENTÉE\\ POUR OBTENIR LE GRADE DE \\
{\LARGE \textbf{DOCTEUR\\DE L'UNIVERSITÉ DE BORDEAUX} } \\
\vspace{0.55cm}
ECOLE DOCTORALE MATHÉMATIQUES ET INFORMATIQUE\\
{\normalsize SPÉCIALITÉ: INFORMATIQUE} \\
\vspace{0.55cm}
Par \textbf{Victor NADOR} \\
\vspace{0.55cm}
{\Large D\'{e}compositions combinatoire pour des familles de cartes d\'{e}form\'{e}es ou d\'{e}cor\'{e}es}
\end{large}
\vspace{0.55cm}
\begin{normalsize}
\begin{singlespace}
Sous la direction de : \textbf{Adrian TANASA}\\
Co-directeur : \textbf{Valentin BONZOM}
\end{singlespace}
\end{normalsize}
\end{center}
\vfill
{\large Soutenue le 21 Septembre 2023, devant le jury composé de: }
\vfill
\begin{table}[b]
\centering
\makebox[\textwidth]{%
\begin{tabular}{lllr}
M. Adrian TANASA & Professeur & Université de Bordeaux & Directeur  \\
M. Valentin BONZOM & Professeur & Université Gustave-Eiffel & Co-Directeur     \\
M. Guillaume CHAPUY & Directeur de Recherche & Université Paris Cité & Rapporteur  \\
M. Razvan GURAU  & Professeur & Heidelberg Universit\"at & Rapporteur  \\
Mme. Marie ALBENQUE & Directrice de Recherche & Université Paris Cité & Examinatrice\\
M. Ga\"{e}tan BOROT & Professeur & Humboldt U. zu Berlin & Examinateur \\
M. Wenjie FANG  & Ma\^{i}tre de conférences & Université Gustave Eiffel & Examinateur  \\
M. Thomas KRAJEWSKI & Ma\^{i}tre de conférences & Université Aix-Marseille & Examinateur \\
\end{tabular}
}
\end{table}

%% file: Chapters/UBX_Resume.tex
\thispagestyle{empty}
\vspace*{0pt}
\vfill
\begin{small}
\begin{center}
\textbf{D\'{e}compositions combinatoire pour des familles de cartes d\'{e}form\'{e}es ou d\'{e}cor\'{e}es}
\end{center}    
\textbf{R\'{e}sum\'{e} :} 
Le d\'eveloppement perturbatif en graphes de Feynman des th\'eories des champs tensorielles peut s'interpreter comme s\'eries g\'{e}n\'{e}ratrices pond\'{e}r\'{e}es de certaines vari\'{e}t\'{e}s lin\'{e}aires par morceaux. Cette remarque \'{e}tabli un lien entre deux domaines à priori tr\`{e}s diff\'{e}rents: la combinatoire des vari\'{e}t\'{e}s discr\`{e}tes d'une part et les th\'{e}ories de tenseurs al\'{e}atoires d'autre part. Dans cette th\`{e}se, nous nous int\'{e}ressons à diff\'{e}rentes propri\'{e}t\'{e}s gravitant autour de ce lien entre combinatoire et th\'{e}orie des champs.

\medskip

Dans un premier temps, nous \'{e}tudions certains mod\`{e}les de constellations. Ces objets g\'{e}n\'{e}ralisent les cartes, ce qui en fait des candidats naturels pour comprendre la $b$-d\'{e}formation, une d\'{e}formation de l'alg\`{e}bre des fonctions sym\'{e}triques dont l'interpr\'{e}tation combinatoire est à l'origine de plusieurs conjectures sur la combinatoire des cartes. Nous \'{e}tudierons les contraintes satisfaites par la s\'{e}rie g\'{e}n\'{e}ratrice des constellations $b$-d\'{e}form\'{e}es que nous qualifions de cubiques. \`{A} partir de l'\'{e}quation d'\'{e}volution satisfaite par cette serie g\'{e}n\'{e}ratrice, nous parvenons à extraire un ensemble de contraintes valides pour toute valeur du param\`{e}tre $b$.

\medskip

Dans un second temps, nous analysons la double limite d'\'{e}chelle de certains mod\`{e}les de tenseurs d'ordre $3$. Lorsque l'ordre est sup\'{e}rieur à deux, la nature du d\'{e}veloppement en $\frac{1}{N}$ - où $N$ est la taille du tenseur- est qualitativement tr\`{e}s diff\'{e}rente du cas matriciel d'ordre $2$. En particulier, seul l'ordre dominant du d\'{e}veloppement est connu explicitement. Malgr\'{e} celà, il est possible d'identifier les graphes sous-dominant contribuant à la double limite d'\'{e}chelle en impl\'{e}mentant la d\'{e}composition en sch\'{e}mas pour les graphes de Feynman. Une analyse des singularit\'{e}s de la s\'{e}rie g\'{e}n\'{e}ratrice de ces sch\'{e}mas nous permet ensuite de caract\'{e}riser pr\'{e}cis\'{e}ment les graphes de la double limite d'\'{e}chelle, et de donner une expression explicite de la fonction à deux points dans cette limite.

\medskip

Enfin, nous nous int\'{e}ressons à une connection entre les th\'{e}ories des champs tensorielles et vectorielles admettant une limite melonique. Nous montrons qu'il est possible d'obtenir la th\'{e}orie vectorielle de Amit-Roginski en consid\'{e}rant des perturbations particuli\`{e}res autour de solutions classiques aux \'{e}quations du mouvement du mod\`{e}le tensoriel de Boulatov. Nous donnons des conditions suffisantes sur ces solutions pour que l'action effective des perturbations prennent la forme de l'action du mod\`{e}le de Amit-Roginski.

\bigskip

\textbf{Mots-cl\'{e}s :} mod\`{e}les de matrices, b-d\'{e}formation, tenseurs aléatoires, théories meloniques, algèbre de contraintes, double limite d'\'{e}chelle\\
%\noindent\rule{\textwidth}{1pt}
\noindent\makebox[\linewidth]{\rule{\textwidth}{0.4pt}}

\newpage

\selectlanguage{english}
\begin{center}
\textbf{Combinatorial decompositions for deformed or decorated classes of maps}
\end{center}
\textbf{Abstract:} 
The perturbative expansion of tensorial field theories in Feynman graphs can be interpreted as weighted generating series of some piecewise linear varieties. This simple fact establishes a link between two a priori distinct fields: the combinatorics of discrete manifolds on one hand and tensorial field theories on the other hand. In this thesis, we study different aspects revolving around this connection between combinatorics and field theory.

\medskip

Firstly, we consider constellations model. These objects generalize maps and their algebraic properties. This makes them suited to probe the $b$-deformation, a deformation of the algebra of symmetric functions which has been conjectured to have a combinatorial interpretation. We will study the constraints satisfied by the generating series of cubical $b$-deformed constellations. Starting from a general evolution equation satisfied by this generating series, we extract the set of constraints satisfied by their generating series for all values of $b$.

\medskip

Secondly, we analyze the double scaling limit of particular tensor models of order $3$. For tensor of order greater than two, the nature of the $\frac{1}{N}$-expansion - where $N$ is the size of the tensor - is qualitatively different from the matrix case of order $2$. In particular, only the leading order graphs are fully characterized. Despite this fact, it is possible to identify graphs of subleading orders contributing to the double scaling limit by implementing the scheme decomposition for Feynman graphs of these theories. An analysis of the singularity of the schemes then allows us to give a complete characterization of the graphs contributing to the double scaling limit. This further enables an explicit computation of the two-point function in this limit.

\medskip

Finally, we investigate a particular link between a tensor and a vector field theory which both admit a melonic limit. Namely, we will show that we can obtain the vectorial Amit-Roginski model by considering perturbations around a classical solution of the Boulatov model, a tensorial theory. We give sufficient conditions on the classical solution so that the effective action for the perturbation around this solution takes the form of the Amit-Roginski action.

\bigskip

\noindent\makebox[\linewidth]{\rule{\textwidth}{0.4pt}}

\vfill
\selectlanguage{french}
\begin{center}
    \textbf{Unit\'{e} de recherche}\\
UMR 5800 Universit\'{e}, 33000 Bordeaux, France.
\end{center}
\end{small}
\vfill

%% file: Chapters/Remerciements.tex
\chapter*{Remerciements}

L'écriture de ce manuscrit me donne l'opportunité de remercier ceux grâce à qui j'ai pu trouver des repères et m'orienter dans le dédale des possibles pour finalement arriver ici aujourd'hui. La tâche n'est pas si facile car ils ont été nombreux ! Je vais faire de mon mieux et j'espère que ceux que je pourrais oublier voudront bien pardonner mon ingratitude.

Tout d'abord, je me dois de remercier Adrian et Valentin, mes deux directeurs de thèse. C'est avec eux que j'ai découvert le vaste et passionnant domaine des tenseurs et des cartes aléatoires.  Bénéficier de leur expérience et de leurs conseils tout au long de ces trois années a été une chance précieuse.

Je tiens à remercier Guillaume Chapuy et Razvan Gurau qui ont accepté la tâche méticuleuse de rapporter ce manuscrit, ainsi que tous les membres du jury: Marie Albenque, Gaëtan Borot, Wenjie Fang et Thomas Krajewski. Les modestes résultats de ce manuscrit s'inscrivent dans la prolongation de certains de leurs travaux de recherche que j'ai eu le plaisir de découvrir durant ma thèse. Je suis honoré qu'ils aient accepté de faire parti de mon jury.

Je remercie également Jean-François Marckert et David Dean pour leurs précieux conseils lors de nos réunions de suivi annuelles.

Je voudrais remercier tous ceux avec qui j'ai pu partagé un intérêt commun pour la recherche. Je pense tout d'abord à mes collaborateurs Daniele, Xiankai et Yili à Munich dont les visites à Bordeaux et les nombreuses réunions sur Zoom nous auront finalement permis de faire aboutir notre projet. Je pense également aux membres des différentes communautés que j'ai rencontré lors de conférences et de séminaires pour les nombreuses discussions passionnantes que j'ai pu avoir avec eux -sur des sujets scientifiques ou non- et plus spécifiquement à Stéphane, Luca, Sabine, Sylvain, Nicolas D., Romain, Léonard, Thomas M., Joseph, Fabien, Dario, Vincent R., Reiko, Roukaya, Dina, Davide, Hannes, Carlos, Ariane, Séverin, Baptiste, Luis, Harriet, Zéphyr, Houcine, Nicolas B., \'{E}ric, Vincent D., Andréa.  J'ai pris un grand plaisir à évoluer dans le même domaine de recherche qu'eux grâce à leur enthousiasme contagieux qui rend la recherche vivante.

Je souhaiterais aussi remercier tous ceux avec qui j'ai partagé le quotidien de la recherche en travaillant dans le même laboratoire ou le même bureau qu'eux. Au LABRI, Dimitri, Alexandre et Tobias qui m'ont accueilli quand je suis arrivé en stage de Master. Maxime, pour les verres lors de ses passages sporadiques à Bordeaux, et Imed dont la présence régulière même pendant le COVID était d'un grand soutien, avec qui j'ai partagé le bureau B367 bureau pendant presque 2 ans. Clément et Marie, pour les multiples séances de grimpe. Géraud et son opiniatreté pour que je rédige une leçon claire et accessible. 

Au LIPN, j'aimerais très chaleureusement tous les membres de la B103 que j'ai côtoyé. Ghazi, pour notre passion commune des memes et de m'avoir expliqué ce qu'étaient les NLP avant même ChatGPT. Francesco e Francesco, pour avoir ravivé mes notions d'italien par leurs discussions. J'espère avoir l'occasion de plus pratiquer pour un jour me joindre à vos échanges. Alex, pour les nombreuses discussions politiques jamais dénuées d'humour ce qui rendrait presque appréciable qu'on ne soit pas souvent d'accord. Dasha, pour ses impressionnantes figures et maîtrise de TeX. Mathieu, pour son énergie débordante, ses incitations répétées à faire des tâches administratives qu'on a veut toujours reporter et surtout de m'avoir initié à la photo. Théo, pour toutes ses anecdotes culturelles et les débriefs de séminaire ici, mais aussi à ALEA, aux JCB et jusqu'à Cracovie. Francis, pour sa bonhomie constante pendant quasiment 6 mois. Grâce à vous tous, l'ambiance au bureau a été au beau fixe pendant tout le temps que j'ai passé dans ce bureau et c'était un vrai plaisir de venir ici. Je remercie également Carole, pour les pauses au soleil salutaires en période de rédaction. Gaël, d'avoir toujours accepté de partager ses bons fruits. Nour, pour sa bonne humeur permanente et ses histoires invraisemblables. Will, pour toutes les discussions mathématico-philosophico-sociétales qu'on a pu avoir.

Durant ces trois années, j'ai aussi fait de belles rencontres en dehors des laboratoires. Je me dois de remercier mes amis de Bordeaux: Justine, Joanna, Federico, Léa, Camille, Emeric et Margaux. Arriver dans une nouvelle ville en pleine pandémie n'était pas forcément facile mais j'ai eu la chance de vous rencontrer ! Merci pour tout les bons moments qu'on a passé ensemble depuis !

C'est aussi l'occasion pour moi de remercier tous ceux qui me connaissent depuis un certain nombres d'années. Je pense à Peio, Calvin, Benjamin, Baptiste des années passées à l'ENS de Lyon. \`{A} Alexis, BEB, Corentin et Paul des années de prépa. \`{A} mes amis du lycée, Alexia, Clara, Jules, Lucas, Baptiste, Lisa, Baulte, Lilian, Guignard, Laurencin, Milet, Cécile, Pippo, Samuel, Max, Jérôme, Zola, Rozzette, Lorène, Juliane, Dylan. Et enfin à ceux qui me connaissent depuis mon plus jeune âge: Chaban, Loïc, Guillaume et Victor. Il ne fait aucun doute que le temps passé à vos côtés a eu une grande influence sur moi. Merci d'avoir toujours été là pour suivre mes péripéties et m'inclure dans les vôtres ! Et surtout, merci d'avoir supporté mes soliloques récurrents sur les maths, la physique, les sciences et leur importance, ou tout autre sujet qui peut me tenir à coeur. 

Je me dois aussi de remercier les membres de ma famille sans qui je n'aurais jamais pu arriver là. Mon cousin, Damien, qui m'a appris très tôt à bidouiller sur l'ordinateur de mes parents. Ma tante, Lucette, qui a éveillé ma curiosité scientifique en étant disponible pour répondre aux pléthores de questions de l'enfant que j'étais. Mes parents, Mireille et Philippe, qui m'auront transmis un peu malgré eux qu'il est possible de vivre de sa passion et qui m'ont toujours supporté dans mes choix.

Enfin, une pensée toute particulière pour Zelda. Par un curieux mélange de gentillesse, de second degré et de cynisme, tu parviens à rendre le quotidien à tes côtés (pas) si (dés)agréable.

%% file: Chapters/Intro_VF.tex
\selectlanguage{french}

\chapter*{Introduction {\small (en français)}}
\addcontentsline{toc}{chapter}{Introduction}  
\label{Chap:intro}

\epigraph{\itshape "Si tout ce que vous avez est un marteau, tout resssemble \`{a} un clou."}{\textup{A. Maslow}, Loi de l'instrument}

En publiant les \emph{Principia Mathematica} en $1687$, Newton a fondé les bases de la mécanique classique à travers ses \emph{lois du mouvement}, des équations mathématiques qui gouvernent le mouvement de tout objet, des planètes jusqu'aux pommes. Grâce à ces équations abstraites, il parvenait à expliquer la loi empirique de Kepler décrivant les trajectoires des planètes dans le Système solaire. En ancrant les observations physiques dans un cadre mathématique abstrait, Newton a initié la mathématisation de la physique. Depuis, ce processus n'a cessé de se poursuivre et est une des pierres angulaires de notre paradigme scientifique moderne. Mais l'\oe{}uvre de Newton a eu un impact tout aussi important en mathématiques. Pour formuler ses lois, Newton a dû développer de nouveaux outils mathématiques pour l'époque. Ces techniques ont été raffinées par Leibniz peu de temps après pour donner naissance au calcul différentiel avec lequel tout scientifique d'aujourd'hui est familier. Depuis, le développement de la physique et des mathématiques ont été intimement liés. Les physiciens emploient des notions de probabilités en physique statistique, utilisent les groupes de Lie pour décrire les interactions entre particules élémentaires à travers les théories de jauge, appliquent la géométrie riemannienne avec la relativité générale. En retour, certains domaines des mathématiques ont émergé grâce à des problèmes physiques. La transformée de Fourier a été inventée pour résoudre l'équation de la chaleur, la géométrie symplectique est apparue à travers l'étude de la formulation hamiltonienne de la mécanique classique, le développement de l'analyse fonctionnelle a été influencée par les besoins des théories quantiques. La physique emprunte des notions mathématiques et les applique pour former ses théories et en retour, certains champs des mathématiques apparaissent ou grandissent par la nécessité de disposer d'outils adaptés pour décrire notre monde à travers nos théories physiques.

\bigskip

Non sans ironie, la gravité est toujours au c\oe{}ur des interactions entre mathématiques et physique théorique plus de $300$ ans après Newton à traver le problème de la \emph{gravité quantique}. Ce problème peut essentiellement se résumer à la question suivante : comment la théorie de la relativité générale doit-elle être quantifiée ? Depuis plus de $70$ ans, de nombreux mathématiciens et physiciens se sont attaqués à cette question, conduisant a l'émergence de nouvelles idées dans les deux disciplines. Cependant, obtenir une réponse complète à cette question reste un des enjeux de la physique théorique moderne aujourd'hui. Même si la gravité quantique ne sera jamais mentionnée explicitement dans cette thèse, elle est présente en toile de fond comme l'une des principales motivations derrière les outils mathématiques qui seront discutés dans ce manuscrit. En tant que tel, il convient de dire quelques mots sur ce problème pour souligner certaines des difficultés qu'il présente.

\bigskip

Dans la première moitié du $XX$\textsuperscript{ème} siècle, notre compréhension des lois physiques a subi deux bouleversements par la découverte de la mécanique quantique d'une part, et celle de la relativitié générale d'autre part. La mécanique quantique nous a montrés que les lois physiques avec lesquelles nous sommes familiers à notre échelle cessent de fonctionner à l'échelle microscopique des atomes. Une des modifications majeures de la physique quantique~\footnote{D'autres aspects de la mécanique quantique sont tout aussi fondamentaux comme l'intrication ou le problème de la mesure mais il n'est pas nécessaire de les évoquer pour la discussion qui suit.} est que les objets ponctuels, typiquement les particules élémentaires comme l'électron, doivent être remplacé par une \emph{fonction d'onde} $\psi(x)$ dont le module $\vert \psi(x) \vert^2$ donne la probabilité que l'objet soit observé à un point $x$ de l'espace-temps. En conséquence, des quantités comme la position ou le moment cinétique de l'objet sont promus en opérateurs agissant sur la fonction d'onde et les observables macroscopiques correspondent aux valeurs propres de l'opérateur associé. Ce changement a conduit à la formulation de la \emph{quantification canonique}, une procédure permettant de quantifier une théorie, c'est à dire de construire l'analogue quantique d'une théorie hamiltonienne classique. Cette procédure a ensuite été raffinée en autorisant le nombre de particules d'un système à fluctuer grâce à la \emph{seconde quantification} qui conçoit les particules élémentaires comme excitations de champs et en rendant la théorie covariante dans l'espace-temps de Minkowski, c'est-à-dire compatible avec la théorie de la relativité restreinte. L'ensemble de ces notions a conduit à la formulation de la \emph{théorie quantique des champs} et sa formulation comme \emph{intégrale de chemins} par Feynman en $1948$. Dans ce formalisme, la probabilité de transition d'un champ entre deux configurations est donnée par l'ensemble des trajectoires possibles entre les deux configurations, incluant les trajectoires non-physiques, et en pondérant celles-ci par l'action classique $S[\phi]$ de la trajectoire. La \emph{fonction de partition} joue un rôle central car toutes les fonctions de corrélations des champs peuvent être calculées à partir de celle-ci. Pour un champ $\phi(x,t)$, sa fonction de partition $\mathcal{Z}$ est définie comme
\begin{equation}
    \label{eq:part_FPI_fr}
    \mathcal{Z} = \int \mathcal{D}\phi \exp\left(\frac{i}{\hbar}S\left[\phi\right]\right),
\end{equation}
où $\mathcal{D}\phi$ signifie que l'on intègre sur toutes les trajectoires possibles et $S\left[\phi\right]$ est l'action associée au champ $\phi(x,t)$. Ce formalisme rend apparent le lien entre la mécanique quantique et les processus stochastiques\footnote{Cette déclaration est un tant soit peu abusive car des nombres complexes apparaissent dans la définition de la fonction de partition. Les physiciens contournent ces difficultés en faisant une \emph{rotation de Wick} $t \rightarrow -it$ permettant d'obtenir un ensemble statistique bien défini.}. Dans la limite classique où $\hbar$ est petit devant l'action, l'intégrale de chemin est dominé par les trajectoires à phases stationnaires, c'est-à-dire les solutions des équations d'Euler-Lagrange pour l'action $S[\phi]$ et on retrouve alors les solutions de la théorie classique. Cela fournit un cadre efficace pour calculer les amplitudes de différents phénomènes physiques. Le modèle standard de la physique des particules, la théorie qui décrit toutes les particules élémentaires connues à ce jour et leurs interactions entre elles par les forces électromagnétique, forte et faible, est basée sur ce formalisme. Elle n'a cessé d'être confirmé expérimentalement par de nombreuses expériences dans les accélérateurs de particules au cours de la deuxième moitié du $XX$\textsuperscript{ème} siècle jusqu'à aujourd'hui.

\bigskip

La théorie de la relativité générale a mis en évidence que la gravité est de nature différente des trois autres forces connues. Cette théorie ne décrit pas comment des particules interagissent entre elles, mais comment les particules de matière interagissent avec l'espace-temps même. En relativité générale, l'espace-temps n'est plus un cadre fixe dans lequel les événements se produisent mais un objet dynamique dont la géométrie est caractérisée par son contenu en matière (plus précisément, en masse et en énergie) et sa répartition. La dynamique de l'espace-temps est gouvernée par l'action de Einstein-Hilbert\footnote{En général, un terme constant $\Lambda$ est ajouté dans cette action pour prendre en compte l'accélération de l'expansion de l'Univers, mais ce terme ne joue aucun rôle dans la discussion qui suit.}
\begin{equation}
    \label{eq:EH_action_fr}
    S_{EH}[g] = \int \frac{1}{2\kappa}R  \sqrt{-g} d^4x,
\end{equation}
où l'espace-temps a pour signature $(-,+,+,+)$, $g=\det (g_{\mu\nu})$ est le déterminant de la métrique $g_{\mu\nu}$, $R$ est le scalaire de Ricci de cette métrique et $\kappa$ est la constante gravitationnelle d'Einstein. Cette action peut-être couplée à la matière en incluant une action décrivant les interactions de la matière de la forme
\begin{equation}
    \label{eq:SM_act_fr}
     S_{m}[\phi,g] = \int \mathcal{L}_m(\phi) \sqrt{-g} d^4x,
\end{equation}
où $\mathcal{L}_m(\phi)$ est le lagrangien décrivant la dynamique des champs $\phi$. Cette action dépend elle-même de la métrique et induit donc à un couplage entre matière et espace-temps dans les équations d'Euler-Lagrange de l'action totale $S_{EH}+S_m$. Tout comme le modèle standard, la théorie de la relativité générale n'a cessé d'être confirmée expérimentalement, de l'observation de l'avancée du périhélie de Mercure et la prédiction des effets de lentilles gravitationnelles en $1915$ à la détection des ondes gravitationnelles à l'image du trou noir au centre de la galaxie M$87$ par l'Event Horizon Telescope en $2018$.

\bigskip

A partir des deux piliers que sont aujourd'hui la mécanique quantique et la relativité générale, on peut tenter de dresser un portrait heuristique des propriétés attendues d'une théorie de la gravité quantique. Cette théorie a pour but d'unifier ces deux mondes. Elle doit donc posséder les aspects géométriques de la relativité générale ainsi que les aspects probabilistes de la mécanique quantique. Elle doit donc donner du sens à une fonction de partition qui prendrait la forme
\begin{equation}
    \label{eq:QG_action_fr}
    \mathcal{Z} = \int \mathcal{D}\phi\mathcal{D}g \exp\left(\frac{i}{\hbar}S_m[g,\phi]+S_{EH}[g]\right),
\end{equation}
où $S_m$ décrit le contenu en matière de l'espace-temps (par exemple le modèle standard de la physique) et où on intègre maintenant sur toutes les métriques possibles en même temps que sur toutes les configurations possible pour les champs. Bien sûr, donner un sens précis et formel à une telle intégrale est une tâche qui couvre de nombreuses difficultés. Par exemple, la théorie de la relativité générale n'est pas quantifiable à travers la seconde quantification car cette approche donne des théories qui ne sont pas renormalisables. De plus, les champs dépendent explicitement des coordonnées de l'espace-temps et l'espace-temps de Minkowski ne peut pas être considéré comme un arrière-plan fixe autour duquel on pourrait former une théorie perturbative décrivant la courbure locale de l'espace-temps. Depuis environ $70$ ans, les physiciens s'attellent à formuler une théorie de la gravité quantique qui répondrait à ces difficultés. Les deux théories les plus connues sont la théorie des cordes~\cite{Pol_String,GSW_String} qui propose que les particules sont des cordes unidimensionnelles qui évoluent dans un espace-temps à $10$ (ou $26$ en l'absence de supersymétrie) dimensions à l'échelle de Planck, et la gravitation quantique à boucles~\cite{RV_LQG} qui propose un espace-temps discret à l'échelle de Planck, constitué d'"atomes d'espace-temps" possédant leurs propres degrés de liberté de jauge. D'autres approchent existent comme la géométrie non-commutative~\cite{Connes_NCG} où les triangulations causales dynamiques~\cite{AJL_CDT}. Chacune de ces théories fait face à son propre lot de difficultés. Ainsi, il n'existe pas encore de théorie de la gravité quantique qui répondent à toutes ces difficultés à ce jour.

\bigskip

Pour mieux sonder les propriétés de ce que serait une théorie de la gravité quantique, une méthode est d'étudier une théorie analogue dans un cadre plus simple que celui de notre espace-temps à $4$ dimensions incluant la matière. Faire ce pas en arrière en considérant des théories moins physiques (c'est-à-dire sans se soucier de leur conformité aux observations de notre univers) permet de concevoir puis d'affiner des outils mathématiques appropriés pour parvenir à surmonter cet immense défi que représente la définition d'une théorie de la gravité quantique. Étant donné le succès du modèle standard et les difficultés intrinsèques à la quantification de la relativité générale, un choix naturel est d'ignorer dans un premier temps la matière de la fonction de partition~\eqref{eq:QG_action_fr}, la transformant ainsi en une théorie de \emph{géométrie aléatoire}. Comme nous l'avons évoqué, une théorie de la gravité quantique doit parvenir à allier les aspects géométriques de la relativité générale avec les aspects probabilistes de la mécanique quantique. Il faut donc disposer de notions mathématiques nous permettant de donner un sens rigoureux à des questions que l'on s'attend à rencontrer en mélangeant ces deux domaines. Par exemple, quelle est la probabilité qu'une variété "tirée au hasard" ait une propriété donnée ? Ou encore, quelles sont les propriétés attendues d'une variété aléatoire "moyenne" ? Dans cette thèse, nous nous intéresserons à une des approches créant des outils qui permettent de donner du sens et de répondre à ce type de questions. Cette approche consiste à utiliser des variétés \emph{linéaires par morceaux} (PL-variétés, d'après l'anglais \emph{piecewise-linear manifold}) en lieu et place des variétés usuelles. L'idée principale est de remplacer les variétés par une notion analogue discrète. Une PL-variété de dimensions $d$ est constituée de blocs atomiques i.e. indivisibles qui sont assemblés ensemble de sorte à former un espace ressemblant fortement aux variétés continues. Pour des PL-variétés, leurs propriétés géométriques sont encodées dans la manière dont leurs constituants atomiques sont assemblés entre eux. Travailler avec des objets discrets conduit à substituer l'intégrale sur les métriques possibles de la fonction de partition~\eqref{eq:QG_action_fr} par une somme sur les différentes configurations possibles. Les distributions de probabilité sur l'ensemble des variétés discrètes sont donc définies en assignant un poids à chacune d'entre-elles, comme on le fait de manière usuelle en probabilité discrète. Ainsi, cette approche fournit un cadre simple pour associer des notions probabilistes et géométriques ensemble. L'étude des propriétés des PL-variétés et de leur distribution de probabilité se trouve à l'interface entre la physique (prenant inspiration dans les théories de gravité quantiques), des mathématiques (par ses liens avec la géométrie aléatoire) et la combinatoire (énumérant des objets générés à partir de constituants atomiques) et a été un sujet de recherche foisonnant dans les trois disciplines au cours de ces $40$ dernières années. 

\bigskip

Les PL-variétés en dimensions deux sont celles qui ont le plus été étudiées dans la littérature. Elles sont les \emph{cartes combinatoires} (simplement appelées \emph{cartes} ici après) introduite par Tutte dans ses travaux pionniers~\cite{Tutte1,Tutte2} au cours des années $1960$. Vingt ans plus tard, les formules obtenues par Tutte ont suscité l'intérêt de l'école française de combinatoire qui étudiait les décompositions des langages basés sur leur grammaire. Ces derniers ont obtenu des bijections~\cite{CV81} permettant de donner une interpretation combinatoire explicite aux formules obtenues par Tutte. Depuis, l'étude des cartes a été un sujet très actif en combinatoire, donnant lieu à une grande variété de résultats dans des directions différentes. Plusieurs autres décompositions ont été introduites, e.g.~\cite{ChaMa} permettant d'établir des bijections avec diverses classes d'objets~\cite{CV81,Schaeffer_Phd,BeFu,AlPo,BDFG}. Ces méthodes bijectives sont important tant d'un point de vue théorique que pratique. Ces bijections permettent de disposer d'un encodage numérique efficace qui est un ingrédient crucial pour la génération aléatoire et les autres implémentations algorithmique de ces objets. Certaines mettent en évidence certaines propriétés algébriques de leurs fonctions génératrices~\cite{BC_ratio,AL_ratio} et d'autres permettent de décrire les géodésiques entre deux points de la PL-variété associée~\cite{CV81,Schaeffer_Phd,AlPo,BDFG}. Cependant, il existe certaines formules obtenues par des techniques différentes, par exemple~\cite{GJ08,CaCha} reposant sur la hiérarchie de Kadomtsev-Petviashvili, qui n'admettent pas d'interprétation bijective à ce jour. En tant qu'espace métrique, l'étude de la limite continue des cartes (grossièrement, lorsque leur nombre d'"atomes" tend vers l'infini) a permis la définition d'un espace métrique aléatoire, la \emph{sphère brownienne}~\cite{MM_Brownian_map,LG_Brownian_map}. Dans la série d'articles~\cite{MS_LQG1,MS_LQG2,MS_LQG3}, il a été montré que la sphère brownienne était liée à la théorie de Liouville de la gravité. Ce résultat illustre l'importance des cartes pour une théorie de la gravité quantique en dimensions deux.

\bigskip

Cependant, dans la communauté de physique théorique, les cartes ne sont pas historiquement apparues comme un modèle-jouet pour la gravité quantique mais à travers les théories de jauges prenant valeur dans un espace de matrices. Le développement perturbatif en graphes de Feynman de ces théories est donné par la série génératrice des cartes organisée par genre. Ce résultat établit un pont entre la combinatoire des cartes et les matrices aléatoires permettant une confluence fructueuse entre les deux objets, principalement en exploitant les propriétés d'intégrabilité de ces objets. Cela a conduit au développement de la \emph{récursion topologique}~\cite{Ey_TR1,EyO_TR2,Eynard_book}, une méthode énumérative basée sur la structure récursive des \emph{équations de boucles} qui a depuis été implémenté pour plusieurs autres modèles~\cite{BE_loop,BCCGF22} et généralisé au delà de sa formulation originelle~\cite{BE_glob_loc,BBCCD_W}. Cela a également permis d'importer des méthodes issues de la théorie des champs pour l'énumération de cartes. Un tel exemple est donné par la formule de récursion de~\cite{GJ08} sur les triangulations, laquelle est la formule la plus efficace dont on dispose à ce jour pour énumérer ces objets numériquement. Cette formule a été obtenue en exploitant le fait que la série génératrice correspondante est une solution de la hiérarchie de Kadomtsev-Petviashvili et n'admet pas d'interprétation combinatoire à ce jour.

\bigskip

Le cas de la dimension deux est singulier pour deux raisons. Les notions de variétés et de PL-variétés coïncident et l'ensemble des informations sur la topologie de la variété est contenu dans un seul entier positif (le genre). En dimensions $d>2$, moins de choses sont connues sur les propriétés des PL-variétés de dimensions $d$, et leurs séries génératrices ne semblent pas (encore ?) disposer de la même richesse mathématique que le cas bidimensionnel. Toutefois, un mécanisme similaire à celui reliant matrices aléatoires et cartes existe toujours. La matrice aléatoire est alors remplacée par un tenseur aléatoire d'ordre $d$ (c'est-à-dire avec $d$ indices). La série perturbative associée au modèe de tenseurs correspond alors a la série génératrice de certaines classes de PL-variété de dimensions $d$~\cite{Amb_tens,Sasa_tens,Gross_tens}. Cette série admet un développement similaire au cas matriciel~\cite{Gu_exp1,Gu_exp2,FeVa}, indexé par un nombre rationnel positif mais qui n'est toutefois plus un paramètre topologique de la PL-variété. Ce fait rend l'identification des graphes de Feynman d'un ordre quelconque difficile. Il n'existe à ce jour aucune méthode permettant de caractériser un ordre arbitraire. Dans la plupart des modèles, seul l'ordre dominant est connu et est en bijection avec certaines familles d'arbres~\cite{GuRy} ce qui place ces modèles de tenseurs dans la classe d'universalité des polymères branchés contrairement au cas matriciel. Comme nous le verrons dans ce manuscrit, la double limite d'échelle -qui est reliée à la limite continue- ne permet pas de sortir de cette classe d'universalité~\cite{DaGuRi,TaGu,BeCa,BNT_DS_TM,BNT_DS_MM,KMT_DS}. En dehors des défis techniques que proposent les tenseurs aléatoires du point de vue de la géométrie aléatoire, ces modèles ont suscité un fort intérêt de la part de la communauté de physique des hautes énergies ces dernières années. Ces modèles possèdent des propriétés similaires~\cite{Wi_SYK}  au modèle de Sachdev-Ye-Kitaev qui a été fortement étudié dans la dernière décennie pour ses propriétés comme modèle-jouet pour la dynamique des trous noirs~\cite{BH_SYK} de par ses propriétés holographiques  (quasi-)AdS/(quasi-)CFT~\cite{Ros_SYK,GrRo_SYK,MaSta_SYK}. De plus, la limite large $D$ des équations d'Einstein est liée à certaines classes de modèles de tenseurs~\cite{Emparan1,Emparan2}. Cela a motivé l'introduction de diverses théories des champs tensorielles~\cite{GKT_Boson,GKPPT_Prism,BCGS_GN} et l'étude de leurs propriétés, notamment autour de la renormalisabilité de ces théories~\cite{BGHS_Unit,BGH_Line,BGHL_F,BN_FP}.

\paragraph{Organisation du manuscrit\\}

Dans le reste du manuscrit, nous étudions les propriétés de certaines familles de PL-variétés. En dimension $2$, nous nous intéressons aux propriétés de la série génératrice de certaines classes de cartes déformées, et en dimensions $3$ où nous implémentons la double limite d'échelle pour certains modèles de tenseurs. Le reste du manuscrit est organisé comme suit,

\begin{itemize}
    \item Le chapitre~\ref{Chap:mat_mod} présente en détail certaines des idées évoquées ici dans le cas bidimensionnel. Nous explicitons le lien entre les cartes et le développement perturbatif de la fonction de partition des matrices aléatoire~\cite{FGG86,Eynard_book}. Nous introduisons également la hiérarchie de Kadomtsev-Petviashvili~\cite{JM83,AZ13}. Cette hiérarchie donne une infinité dénombrable d'équations aux dérivées partielles toutes satisfaites par la fonction de partition des matrices aléatoires. Cette présentation sert deux objectifs. Le premier est d'illustrer les bénéfices qui peuvent être retirés en exploitant le lien entre combinatoire des cartes et systèmes intégrables. Le second est de motiver les résultats du chapitre suivant.
    \item Dans le chapitre~\ref{Chap:const}, nous introduisons deux généralisations aux cartes du chapitre précédent. La première est la $b$-déformation, une déformation à un paramètre de l'algèbre des fonctions symétriques qui permet d'interpoler entre les séries génératrices des cartes orientables et des cartes non-orientées. Cette déformation est la source de plusieurs conjectures autour de l'énumération de cartes~\cite{GJ_b_conj}. La seconde est l'introduction des constellations~\cite{LZ_book} qui généralisent les cartes biparties vues comme recouvrement ramifié de la sphère de Riemann. En combinant ces deux généralisations, nous obtenons les constellations $b$-déformées introduites dans~\cite{CD22}. Nous étudions ensuite les propriétés de certins modèles particuliers de constellations $b$-déformées que nous appelons \emph{cubiques} et nous montrons que leur série génératrice satisfait une infinité de contraintes qui forment une algèbre sous commutation. Ces résultats devraient être publiés bientôt dans~\cite{BN23}.
    \item Le chapitre~\ref{Chap:tens_mod} étend la construction des PL-variétés en dimensions supérieures à deux. Nous y introduisons le lien entre les PL-variétés de dimensions $d$ et les modèles de tenseurs d'ordre $d$ en présentant certains résultats maintenant classiques sur le sujet~\cite{Book_Gu,Book_Ta,GuRy2,Ta1}. Nous nous concentrons sur les modèles de tenseurs sans mélange d'indices. En particulier, nous présentons leur développement en $\frac{1}{N}$ et caractérisons l'ordre dominant de ce développement. Enfin, ces résultats sont étendus au cas multi-matriciels~\cite{FeVa}. Ce chapitre sert d'introduction aux modèles de tenseurs et permet d'introduire toutes les notions intervenant dans le chapitre suivant.
    \item Dans le chapitre~\ref{Chap:DScale}, nous implémentons la double limite d'échelle dans trois modèles de tenseurs différents. Bien que la méthode pour y parvenir soit commune aux trois modèles, son implémentation dépend de détails spécifiques de chacun des modèles considérés. Un outil combinatoire clé pour l'obtention de cette limite et la décomposition en schémas des graphes de Feynman de ces modèles de tenseurs. Cette décomposition permet d'analyser les singularités du développement perturbatif de la fonction de partition de ces modèles à tout ordre du développement en $\frac{1}{N}$ et d'identifier celles qui sont pertinentes pour la double limite d'échelle. Cela nous permet d'identifier les contributions de l'ordre dominant de la double limite d'échelle, ce qui nous permet enfin d'implémenter celle-ci en resommant les contributions dominantes. L'exemple prototypique est celui du modèle $O(N)^3$ étudié dans~\cite{TaCa} avec interactions quartiques. Nous montrons comment ces résultats peuvent être transférés au modèle multi-matrices $U(N)^2\times O(D)$. Enfin, nous appliquons une méthode similaire dans un modèle multi-matrices $U(N)\times O(D)$ complexe quartique. Ce modèle est plus compliqué du point de vue combinatoire car impliquant du mélange d'indice, mais la décomposition en schémas peut tout de même être implémentée, permettant par la suite d'obtenir l'ordre dominant de la double limite d'échelle de la même manière que dans les deux cas précédents. Ce chapitre est une version éditée des deux articles~\cite{BNT_DS_TM} and~\cite{BNT_DS_MM}
    \item Le chapitre~\ref{Chap:melons_th} établit un lien entre deux théories des champs de nature différente présentant un développement en $\frac{1}{N}$ dominé par des graphes meloniques. Plus précisément, nous montrons que l'action du modèle de Amit-Roginsky~\cite{AR_OG,AR_BeDe}, une théorie des champs vectorielle, peut être obtenue comme action effective pour des perturbations de solutions classiques du modèle de Boulatov~\cite{Boul_OG}, une théorie tensorielle avec dépendance sur un groupe de Lie. Ce résultat devrait être publié bientôt dans~\cite{BoMun}.
    \item Enfin, le chapitre~\ref{Chap:concl} donne des perspectives possibles pour prolonger les travaux présentés dans ce manuscrit. 
\end{itemize}

\selectlanguage{english}
%%%%%%%%%%%%%%%%%%%%%%%%%%%%%%%%%%%%%%%%%%%%%%%%%%%%%%%%%%%%%%%%%%%%%%%%%%%%%%%%%%%%%%%%%%%%%%%%%%%%%%%%%%%%%%%%%%%%%%%%%%%%%%%%%%%%%%%%%%%%%%

%% file: Chapters/Intro.tex
% Chapter Template

\chapter*{Introduction {\small (in english)}}

\epigraph{\itshape "If all you have is a hammer, everything looks like a nail."}{ \textup{A. Maslow}, Law of the instrument}

In the \emph{Principia Mathematica} in $1687$, Newton laid the foundation of classical mechanics through its \emph{law of motions}, mathematical equations governing the movement of any object from planets to apples. With these abstract equations, he could explain Kepler's empirical law for planetary motion. It grounded physical observations in an abstract mathematical framework, initiating the mathematization of physics that shaped our modern scientific paradigm to this day. But the \emph{Principia} also revolutionized mathematics as Newton developed new mathematical methods to formulate his laws. Shortly after, Leibniz refined these tools giving rise to differential calculus. Since then, the development of mathematics and physics have been intimately intertwined. Physicists employ probabilities in statistical physics, Lie groups to describe the interactions of particles through gauge theories, and use Riemannian geometry in general relativity. Mathematicians have discovered the Fourier transform to solve the heat equation. Symplectic geometry raised from the Hamiltonian formulation of classical mechanics and quantum mechanics influenced the development of functional analysis. Physics borrows and applies mathematical concepts to shape its theories and in return, some fields of mathematics grow from the need for new tools to describe our physical world through our theories.

\bigskip

With a touch of irony, more than $300$ years after Newton, gravity is once again at the core of the interplay between mathematics and theoretical physics through \emph{quantum gravity}. The problem of quantum gravity can be boiled down to the following question: how does general relativity have to be quantized? Over the last $70$ years, it has received the efforts of many mathematicians and physicists and has led to the emergence of new ideas in both fields. However, a complete answer to this question has yet to be given. Although the explicit mention of quantum gravity may be absent from this thesis, it lurks in the background as one of the main motivations behind the mathematics to be discussed in the rest of this manuscript. As such, it deserves a brief overview outlining the main issues surrounding this problem.

\bigskip

In the first half of the $XX$\textsuperscript{th} century, our understanding of physics underwent two major transformations through the emergence of quantum mechanics and general relativity. Quantum mechanics showed that the laws of physics we are familiar with at our human scale break down at the microscopic scale of atoms. One of the key change~\footnote{Other aspects of quantum mechanics like entanglement or the problem of measure are equally fascinating but need not be discussed here.} of quantum mechanics is that a point-like object (e.g. particles like the electron) has to be substituted by wavefunction $\psi(x)$ such that its modulus $\vert \psi(x) \vert^2$ encodes the probability of observing the object at a point $x$ in space (at a fixed time). As a consequence, the position and momentum of the object are promoted to operators acting on the wavefunction and macroscopic observables are given by eigenvalues of the corresponding operator. It led to the procedure of \emph{canonical quantization} which gives a method to quantize a theory, that is to construct the quantum mechanical analog from the Hamiltonian of a classical theory. These ideas were refined by allowing the number of particles in a system to fluctuate through \emph{second quantization} which sees particles as excitations of fields and requiring the theory to be covariant in Minkowski spacetime i.e. including special relativity. Altogether, it led to the formulation of \emph{quantum field theory} and its \emph{path integral} formulation by Feynman in $1948$. In the path integral formulation, the probability of transition of a field between two configurations is given by the sum over all possible trajectories - including a priori non-physically relevant ones- weighted by the classical action $S[\phi]$ associated with the trajectory. A key object is the \emph{partition function}. For a field $\phi(x,t)$, it is the generating function of all correlation functions between field configurations and is roughly defined as
\begin{equation}
    \label{eq:part_FPI}
    \mathcal{Z} = \int \mathcal{D}\phi \exp\left(\frac{i}{\hbar}S\left[\phi\right]\right),
\end{equation}
where $\mathcal{D}\phi$ means that we integrate over all possible trajectories. This formalism shows the connection between quantum physics and stochastic processes\footnote{This claim is slightly abusive since the partition function involves complex amplitudes. Physicists circumvent this in computation by performing \emph{Wick rotation} $t \rightarrow -it$ to get a well-defined statistical ensemble.}. In the classical limit where $\hbar$ can be neglected, the path integral is dominated by trajectories with stationary phases i.e. solutions of the Euler-Lagrange equations for the action $S[\phi]$ and thus the classical solutions are recovered. The path integral formulation provides an effective framework on which the amplitudes of events can be computed. Over the $XX$\textsuperscript{th} century, the standard model of particle physics -the theory describing all elementary particles and their interactions through the electromagnetic, weak and strong forces- has been described and confirmed by numerous experiments in particle colliders through the use of the path integral formulation.

\bigskip

General relativity established that gravity was of a different nature than the three other forces. It is not a theory describing how particles of matter interact together but a theory of the interactions between matter and spacetime. Spacetime is no longer a fixed background in which events occur but it is a dynamic object whose geometry is influenced by the local content in matter (more precisely in mass and energy). The dynamics of spacetime is described by the Einstein-Hilbert action\footnote{It is possible to include a constant $\Lambda$ in this action to account for the observed accelerated expansion of the Universe.} 
\begin{equation}
    \label{eq:EH_action}
    S_{EH}[g] = \int \frac{1}{2\kappa}R  \sqrt{-g} d^4x,
\end{equation}
where spacetime has signature $(-,+,+,+)$, $g=\det (g_{\mu\nu})$ is the determinant of the metric tensor $g_{\mu\nu}$, $R$ the Ricci scalar of the metric and $\kappa$ is Einstein gravitational constant. It can be coupled to matter by adding to it a matter action for the matter
\begin{equation}
    \label{eq:SM_act}
     S_{m}[\phi,g] = \int \mathcal{L}_m(\phi) \sqrt{-g} d^4x,
\end{equation}
where $\mathcal{L}_m(\phi)$ is the Lagrangian describing the dynamics of the matter fields $\phi$. This action depends on the spacetime metric, leading to a coupling between matter and spacetime in the Euler-Lagrange equations for the total action $S_{EH}+S_m$. This theory has also been confirmed experimentally throughout the last century by various observations from the perihelion precession of Mercury and gravitational lensing in $1915$ to the detection of gravitational waves and the image of the black hole at the center of the galaxy $M87$ by the Event Horizon Telescope in $2018$.

\bigskip

From the two pillars of modern physics which are quantum mechanics and general relativity, we can sketch a heuristic picture of what properties should be expected from a quantum theory of gravity. It should incorporate the geometric features of general relativity as well as the probabilistic aspects of quantum mechanics and make sense of a partition function of the form
\begin{equation}
    \label{eq:QG_action}
    \mathcal{Z} = \int \mathcal{D}\phi\mathcal{D}g \exp\left(\frac{i}{\hbar}S_m[g,\phi]+S_{EH}[g]\right),
\end{equation}
where $S_m$ describes the content in matter of spacetime (e.g. the action of the standard model) and we now integrate over all possible metric fields as well as all possible fields configurations. Of course, making sense of such a partition function is a difficult task that contains a lot of caveats. For example, general relativity is non-quantizable through the standard quantization procedure, as it yields a non-renormalizable theory. Moreover, the fields explicitly depend on the spacetime coordinate and the Minkowskian spacetime cannot be taken as a background around which to consider perturbation allowing for local curvature. In the last $70$ years, physicists have come up with various theoryies which aim at overcoming these issues. The two most notorious ones are string theory~\cite{Pol_String,GSW_String} which suggests that particles are one-dimensional strings that live in a $10$ (or $26$ in the absence of supersymmetry) dimensional space at the Planck scale, and loop quantum gravity~\cite{RV_LQG} which takes spacetime to be discrete at the Planck scale, thus forming "atoms of spacetime" with gauge degrees of freedom. Other approaches include for example non-commutative geometry~\cite{Connes_NCG} or causal dynamical triangulations~\cite{AJL_CDT}. All of these theories are facing their own issues and a fully fleshed physical theory of quantum gravity has still to be obtained to this day.

\bigskip

To gain a deeper insight into the features of quantum gravity, one strategy is to study it within simplified frameworks. Given the remarkable success of the standard model and the challenges related to quantization internal to general relativity, one natural simplification is to ignore the matter content in the partition function~\eqref{eq:QG_action}, thereby turning it into a theory of \emph{random geometry}. Taking a step back and considering theories that are physically less relevant is needed to shape the appropriate mathematical tools to address the immense challenge that finding a theory of quantum gravity represents. To obtain a full theory of quantum gravity, we need to have mathematical notions that allow us to make sense of questions we would expect to encounter in random geometry. For example, what is the probability that a space "drawn at random" has a given property or what are the properties expected of an "average" random space? In this thesis, we mostly study one attempt at shaping such notions which focuses on \emph{piecewise-linear manifolds} (PL-manifolds). The starting idea is to replace a $d$-dimensional manifold by a discretized counterpart, a ($d$-dimensional) PL-manifold, made of ($d$-dimensional) atomic building blocks suitably assembled. The geometric properties of the PL-manifold are encoded in the way its atomic constituents are glued together. Considering PL-manifolds instead of usual ones effectively turns the integral over all possible metric structures of Equation~\eqref{eq:QG_action} into a sum. Therefore, probability distributions over PL-manifolds can simply be defined by enumerating them and assigning a weight to each of them, similarly to discrete probability. Thus, this approach makes it easier to associate geometric and probabilistic notions altogether. The study of PL-manifolds lies right at the interface between physics (drawing inspiration from quantum gravity), mathematics (via its connection to random geometry) and combinatorics (enumeration of objects generated from atomic ones) and has been an active topic of research in all three fields over the last $40$ years. 

\bigskip

The two-dimensional PL-manifolds are the ones that have been studied the most in the literature. They are the \emph{combinatorial maps} introduced by Tutte through its pioneering work~\cite{Tutte1,Tutte2} in the $1960$'s. The formulas obtained by Tutte raised interest in the French school of computer scientists who were studying decomposition based on grammar in the $1980$'s. They found bijections that allowed them to give an explicit combinatorial interpretation to the formulas of Tutte~\cite{CV81}. Maps have been a very active topic of research since then, leading to a wide variety of results. Several other decompositions~\cite{ChaMa} and bijections with other combinatorial objects~\cite{CV81,Schaeffer_Phd,BeFu,AlPo,BDFG} have been obtained. These bijective methods are key both from the practical and theoretical standpoint. They allow for efficient numerical encodings of maps which are crucial for random sampling and other algorithmic implementations of these objects. They can be used to shed light on some of the algebraic properties of their generating functions~\cite{BC_ratio,AL_ratio}. Finally, they are also a key ingredient in the context of random geometry as some of them~\cite{CV81,Schaeffer_Phd,AlPo,BDFG} allow to keep track of the geodesic distance between points of the map. However, interestingly, some recursion formulas e.g.~\cite{GJ08,CaCha} can be obtained through different means and do not have a bijective interpretation to this day. As metric spaces, the study of their continuum limit led to the definition of a random metric space known as the \emph{Brownian map}~\cite{MM_Brownian_map,LG_Brownian_map}. The Brownian map has been shown to be related to the Liouville theory of quantum gravity in a recent series of articles~\cite{MS_LQG1,MS_LQG2,MS_LQG3}, further shedding light on the relevance of maps in the context of quantum gravity.

\bigskip 

Nevertheless, in the physics literature, maps did not make their first appearance as a toy model of quantum gravity but through the study of matrix field theory~\cite{HOFT74,BIPZ78}. The perturbative expansion in Feynman graphs of matrix field theory corresponds to the genus generating series of maps. This established a link between the combinatorics of maps and random matrix models which led to fruitful interactions between the two fields, mainly revolving on the integrability properties of these objects. It led to \emph{topological recursion}~\cite{Ey_TR1,EyO_TR2,Eynard_book}, an enumerative method based on the recursive structure of the \emph{loop equations} that has now been applied in a variety of other models~\cite{BE_loop,BCCGF22} and generalized beyond its original formulation~\cite{BE_glob_loc,BBCCD_W}. It also brought field theoretic methods to map enumeration which led to unique results. One such example is the recursion formula of~\cite{GJ08}, which is the most efficient formula to count map numerically. It has been obtained exploiting the fact that the generating series of triangulations is a solution of the KP hierarchy and cannot be obtained by standard combinatorial techniques to this day.

\bigskip

The two-dimensional case is quite special as the notions of PL-manifolds and manifolds coincide and the topological information of the surface is contained in a single non-negative integer (the genus). In dimension $d >2$, much less is known on the properties of PL-manifold in dimensions and it does not (yet ?) exhibit the same mathematical richness as in the two-dimensional case. Still, a similar mechanism as between matrix models and $2$ dimensional PL-manifolds can be obtained in $d$ dimensions by replacing the matrix with a tensor of order $d$ (that is, a tensor with $d$ indices). The perturbative expansion of the tensor field theory can then be interpreted as the generating function for certain classes of piecewise linear manifolds~\cite{Amb_tens,Sasa_tens,Gross_tens}. A $\frac{1}{N}$-expansion can then be obtained by a suitable choice of scaling with the size $N$ of the tensor for the terms of the action~\cite{Gu_exp1,Gu_exp2,FeVa}. However, contrary to the matrix case it is not a topological expansion. This renders the identification of Feynman graphs associated with subleading order difficult, and no general method to characterize them is known to this day. Moreover, its leading order graphs are in bijection with families of trees~\cite{GuRy} which puts tensor models in the universality class of branched polymers, contrary to the matrix case. As we shall see in this manuscript, the double scaling limit- related to a continuum limit- doesn't allow to escape this universality class~\cite{DaGuRi,TaGu,BeCa,BNT_DS_TM,BNT_DS_MM,KMT_DS}. Aside from the technical challenges they offer in the context of random geometry, tensor models have attracted a lot of interest from the high-energy physics community. They yield similar properties as the Sachdev-Ye-Kitaev model~\cite{Wi_SYK} which has been widely studied in the past decade as a toy model for black-hole dynamics~\cite{BH_SYK} due to its near AdS/CFT holographic properties~\cite{Ros_SYK,GrRo_SYK,MaSta_SYK} and some class of tensor models have been shown to be related to the large $D$ limit of Einstein equations~\cite{Emparan1,Emparan2}. This motivated the introduction of various tensor field theories~\cite{GKT_Boson,GKPPT_Prism,BCGS_GN} and the study of renormalization properties of such theories~\cite{BGHS_Unit,BGH_Line,BGHL_F,BN_FP}.

%Random tensor models were first introduced in the $90$'s~\cite{Amb_tens,Sasa_tens,Gross_tens} as a generalization of matrix models in dimensions $d\geq 3$. Similarly to matrix model, each interaction term in the action can be interpreted as a simplicial complex of dimension $d$. Therefore,  their perturbative expansion is given by a sum over possible gluings of these simplicial complexes and corresponds to the generating function for certain classes of piecewise linear manifolds. However, the perturbative expansion lacked a $\frac{1}{N}$-expansion for the next $20$ years until it was implemented for \emph{colored} tensor models~\cite{Gu_exp1,Gu_exp2,GuRi}.

\paragraph{Organisation of the manuscript\\}

In the rest of the manuscript, we will study some properties of PL-manifolds. In dimension $2$, we will study the properties of the generating function of some class of maps, and in dimension $3$ where we will study a particular limit of some tensor models. It is organized as follows,
\begin{itemize}
    \item Chapter~\ref{Chap:mat_mod} presents in detail the main ideas that were sketched to illustrate in the two-dimensional case the existing connection between maps and random matrix model, seen as a perturbative field theory~\cite{FGG86,Eynard_book}. We also introduce the KP hierarchy~\cite{JM83,AZ13}, an infinite set of partial differential equations that the generating series of the Hermitian matrix model satisfies. This presentation serves two purposes. First, it is an example of the fruitful interaction between combinatorics of maps and integrable systems. Second, it illustrates the integrable properties of both objects which was one of the reasons motivating the computation of constraints in the models of the next chapter.
    \item In Chapter~\ref{Chap:const}, we introduce two different generalizations to the objects of the previous chapter. The first one is the so-called $b$-deformation, a one-parameter deformation of the algebra of symmetric functions which allows to interpolates between the generating series of orientable and non-oriented maps and is the source of conjectures in map enumeration~\cite{GJ_b_conj}. The second is the introduction of constellations~\cite{LZ_book} which generalize bipartite maps seen as ramified coverings of the Riemann sphere. Merging the two generalizations gives $b$-deformed constellations introduced in~\cite{CD22}. We study some particular models of $b$-deformed constellations which we call \emph{cubic} and show that their generating series satisfies an infinite set of PDEs that form an algebra under commutation. These results are soon to be published in~\cite{BN23}.
    \item Chapter~\ref{Chap:tens_mod} generalizes the construction of the two-dimensional case to arbitrary dimension, reviewing classical results that can be found in the literature on tensor models~\cite{Book_Gu,Book_Ta,GuRy2,Ta1}. We establish the connection between tensor models of order $d$ and $d$-dimensional pseudo manifold, principally focusing on tensors with no mixed index symmetry. We also discuss the condition under which a $\frac{1}{N}$-expansion for these models exists, and how these results extend in the multi-matrix case~\cite{FeVa}. This chapter serves as a preliminary to introduce all the notions required in the next chapter.
    \item  In Chapter~\ref{Chap:DScale}, we study the double scaling limit~\cite{TaGu,BeCa} for three different tensor models. The strategy is the same in all three models, despite that its implementation depends on the details of the model considered. We first introduce a scheme decomposition~\cite{ChaMa,TaFu,GuSch} for the Feynman graphs of the model. Then, we study the generating function of the leading order graphs of the $\frac{1}{N}$-expansion and identify the relevant singularities. Finally, we identify the graphs contributing at the leading order of the double scaling limit and sum their contribution. The prototypical case is the quartic $O(N)^3$ tensor model. We then adapt this case to establish a similar result in the $U(N)^2\times O(D)$ multi-matrix model. We show that a similar strategy can be employed in some tensor models with mixed index symmetry through the example of the $U(N)\times O(D)$ multi-matrix model. This Chapter is an edited version of articles~\cite{BNT_DS_TM} and~\cite{BNT_DS_MM}.
    \item Chapter~\ref{Chap:melons_th} exhibits a connection between two field theories which have a $\frac{1}{N}$-limit dominated by melonic graphs. More precisely, we show that the Amit-Roginsky model~\cite{AR_OG,AR_BeDe} -a vector field theory- can be obtained as an effective action for perturbations around a classical solution to the equation of motions of the Boulatov model~\cite{Boul_OG} - a tensor field theory with additional group data. This result should be published soon in~\cite{BoMun}.
    \item Finally, Chapter~\ref{Chap:concl} gives some outlooks on possible directions for future work on the topics touched on in this manuscript.
\end{itemize}

%% file: Chapters/Matrix_model.tex
\chapter{Maps, random matrices and integrability} % Main chapter title
\label{Chap:mat_mod}

In two dimension, the notions of PL-manifold and surfaces coincides and the topological properties of surfaces are simple enough to be encoded by a single topological invariant, the genus. It provides a nice framework to illustrate the interplay between the combinatorics of PL-manifolds -combinatorial maps in this case- and tensor-matrix here- field theory. This chapter reviews classical results pertaining to maps and random matrices. As they are rich and active research topics on their own, each of these two objects is dedicated its own section. A remarkable feature of the generating functions of (some of) these objects is their integrability. Integrability manifests itself in various yet distinct forms, reflecting the impressive diversity of mathematical notions that revolve around map enumeration. It appears in conformal field theory, through topological recursion, or via integrable hierarchies, which will be its avatar studied here. Therefore the third section is dedicated to the Kadomtsev–Petviashvili hierarchy, a set of infinitely many partial differential equations satisfied by (some) generating functions of maps and random matrices.

\section{Maps}
\label{sec:maps}

\subsection{Different definitions for maps}
\label{ssec:def_maps}

The principal object we will work within this chapter is (orientable) \emph{combinatorial maps} (maps for short). Behind their apparent simplicity, maps are a very rich mathematical object. As such, they enjoy several equivalent definitions which reflects the diversity of the mathematical tools that can be employed to study them. We present here three equivalent definitions, each bringing a new perspective on maps.

\subsubsection{As embedded graphs}
\label{sssec:emb_graph}

A graph $G$ is a finite set of vertices $V$ and a finite multiset $E$ of couples in $E$.\footnote{Note that we allow multiple edges between vertices and loops, this definition corresponds to the usual notion of multi-graph in graph theory.} Usually, one represents vertices as points and edges as lines between those points, i.e. as copies of the segment $\left[0,1\right]$.

\begin{definition}[Maps as embedded graphs~\cite{MT01}]
\label{def:emb_graph}
An orientable map $\mathcal{M}$ is an embedding of a graph $G$ onto a surface $\mathbb{S}$ such that $\mathbb{S}\backslash \mathcal{M}$ is a union of disks. Each of the connected components of $\mathbb{S}\backslash \mathcal{M}$ is a \emph{face} of the map $\mathcal{M}$. The faces inherit the orientation of $\mathbb{S}$ if the latter is oriented. Two maps are isomorphic if there is an orientation preserving homeomorphism between them.
\end{definition}

Orientable surfaces without boundaries are characterized by a non-negative integer called the~\emph{genus} and denoted $g$. Intuitively, the genus counts the number of holes of the surface. If a map $\mathcal{M}$ is embedded on a surface of genus $g$, then we say the map $\mathcal{M}$ has genus $g$. For a connected map, the genus can be recovered from its number of faces, edges and vertices via the Euler characteristic
\begin{equation}
    \label{eq:Euler}
    \chi(\mathcal{M}) = 2-2g(\mathcal{M}) = V(\mathcal{M})-E(\mathcal{M})+F(\mathcal{M})
\end{equation}

\begin{remark}
    We will represent maps of genus $0$ in the plane rather than on the sphere by choosing a point inside any face as the point at infinity in the plane. This is harmless, as long as we keep in mind the existence of this outer face.
\end{remark}

\begin{figure}
\hfill
\begin{subfigure}{.4\textwidth}
  \centering
  \includegraphics[scale=0.17]{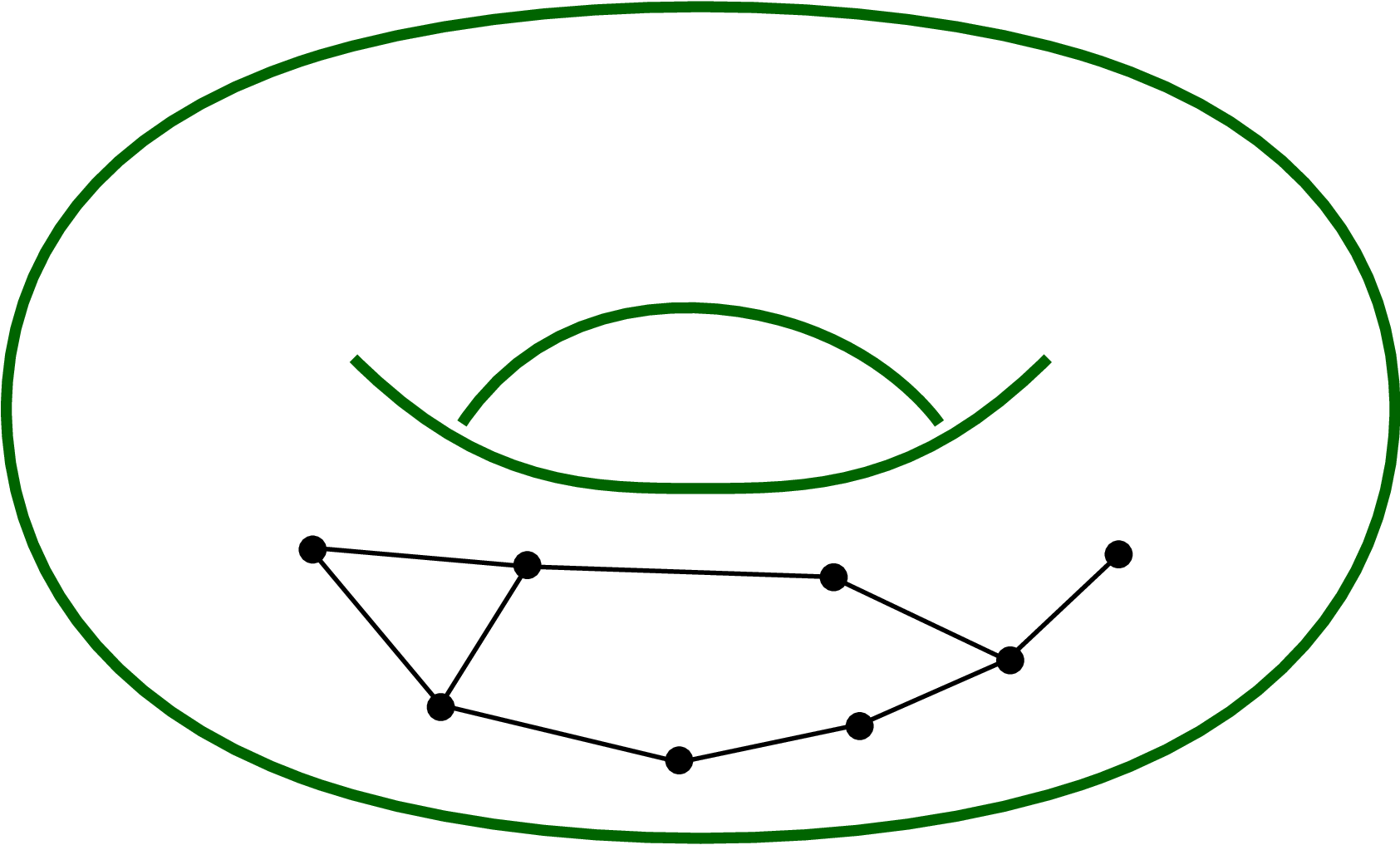}
  \caption{A graph embedding that is not a map.}
  \label{fig:not_map}
\end{subfigure}
\hfill
\begin{subfigure}{.4\textwidth}
  \centering
  \includegraphics[scale=0.32]{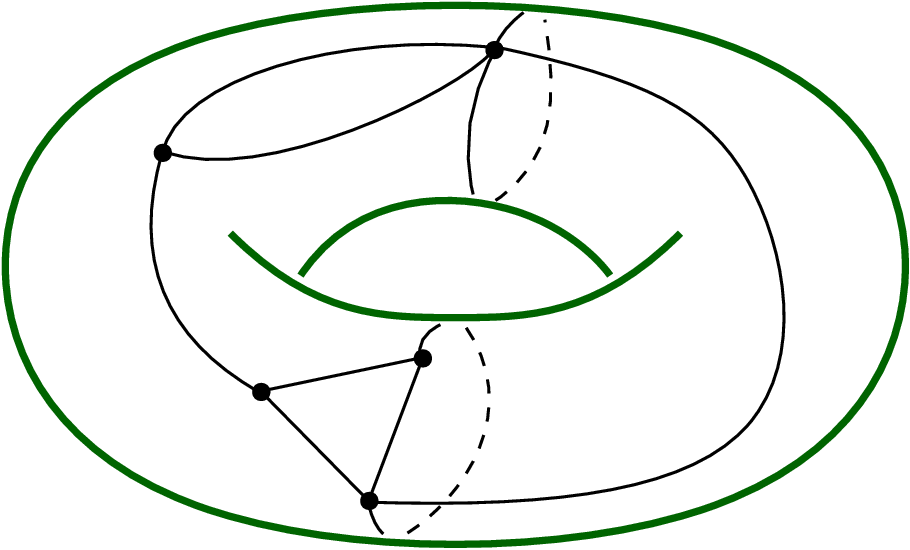}
  \caption{A map as an embedded graph.}
  \label{fig:map_emb}
\end{subfigure}
\hfill
\caption{Two graphs embeddings but only one map.}
\label{fig:example_maps}
\end{figure}
The graph embedded on the left of Figure~\ref{fig:example_maps} is not a map because one of the connected components of $\mathbb{S}\backslash \mathcal{M}$ is not homeomorphic to a disk but to a surface of genus $1$ with one boundary. On the contrary, the graph embedding on the right is a valid map. 

\begin{remark}
Note that for $g \geq 1$, there exist infinitely many isomorphic but non-isotopic maps. Such representations are generated by~\emph{Dehn twists}, which are obtained by cutting the map along any non-contractible of the surface (which exists in genus $g\geq 1$), rotating one of the borders by a multiple of $2\pi$ and gluing them back.
\begin{figure}[!ht]
    \centering
    \includegraphics[scale=0.25]{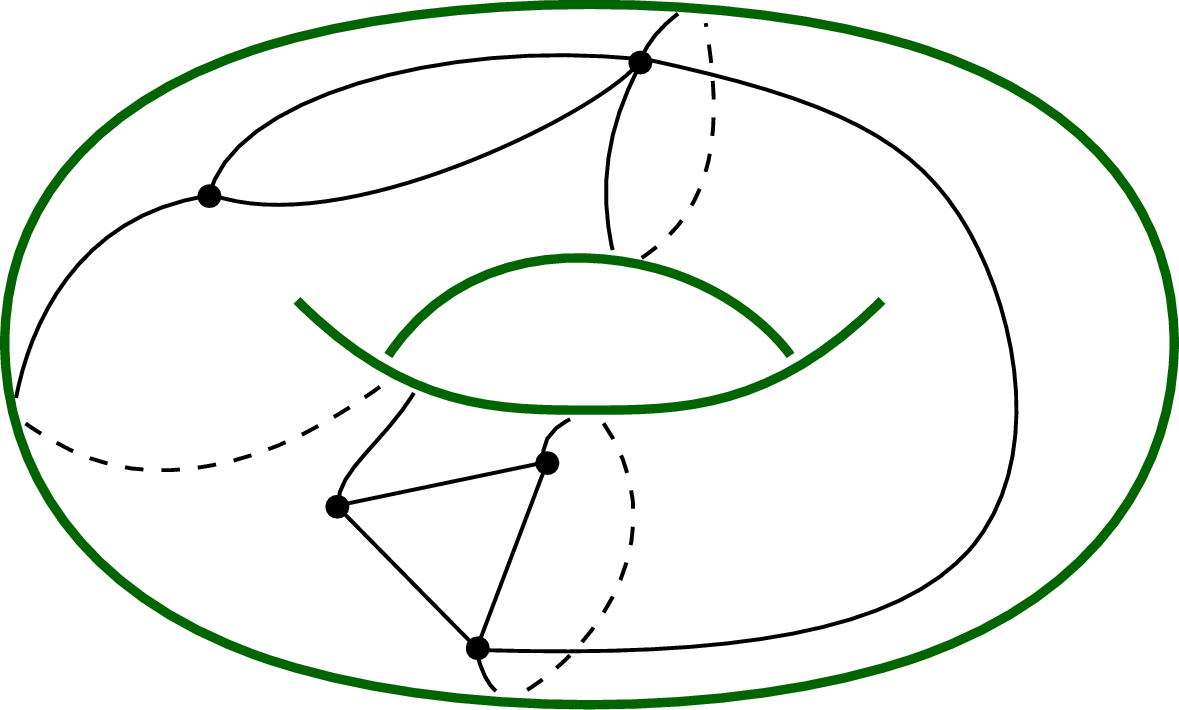}
    \caption{This map is isomorphic to the map of Figure~\ref{fig:map_emb} but cannot be obtained from continuous deformation of this map.}
    \label{fig:map_dehn}
\end{figure}
\end{remark}

The size of the map $\mathcal{M}$ is its number of edges $E(\mathcal{M})$. An edge is made of two half-edges obtained by removing the middle point of an edge (seen as the segment $\left[0,1\right]$). All half-edges are connected to a unique vertex. A vertex has degree $k$ if it has $k$ half-edges incident to it. The boundary of a face $F$ in a map $\mathcal{M}$ is a cycle in the original graph $G$. Starting from any vertex of this cycle, its edges can be enumerated in cyclic order following the orientation inherited from $\mathbb{S}$. The \emph{length} of a face is the length of this cycle. The total degree of the vertices is equal to the total length of the faces, both being given by twice the number of edges. The list of the length of the faces of a map is called its \emph{face profile} and similarly, the list of the degree of its vertices is called its \emph{vertex profile}. When ordered in decreasing order, the face and vertex profiles are given by two \emph{partitions} of $2 E(\mathcal{M})$ (see Def.~\ref{def:part} hereafter). 

\medskip

Closely related to maps are \emph{rooted maps}. A map $\mathcal{M}$ is rooted if it has a distinguished corner. Fixing an orientation of the surface gives a correspondence between corner and half-edges by associating a corner to the half-edge that borders it with respect to that orientation. The face (resp. the vertex) incident to the root corner is called the root face (resp. the root vertex). It will be convenient to think of the root face as a boundary of the surface. This extends to maps with $k$ distinguished half-edges carried on distinct faces which can be thought of as maps embedded on surfaces with $k$ boundaries.

\medskip

This definition shows the main properties of maps. They are graphs that have inherited an orientation from the surface they are embedded onto. This embedding fixes a cyclic order of edges around the vertices, which can be exploited to describe the faces of the map. 

\subsubsection{As a gluing of polygons}
\label{sssec:poly_glu}

All the information necessary to reconstruct a map is contained in the length of its faces and the incidence relations between these faces. This leads to the following definition of maps as a gluing of oriented polygons.

\begin{definition}[Maps as gluings of oriented polygons~\cite{MT01}]
\label{def:poly_glu}
    An orientable map $\mathcal{M}$ is a finite collection of $k$ oriented polygons $\{\mathcal{P}_i\}_{1 \geq i \geq k}$ such that $\mathcal{P}_i$ is a $r_i$-gon with $r_i$ sides in cyclic order $E^{(i)}=\left( e^{(i)}_1, \hspace{2pt}..\hspace{4pt}, e^{(i)}_{r_i}\right)$ and a matching $m$ of the set of all sides $E = \bigcup\limits_{i=1}^k E^{(i)} $.
\end{definition}

The total number of edges must be even as matchings only exist for set of even cardinality. The embedding of the graph is recovered by gluing every pair of sides $(e,e')$ such that $m(e) = e'$ while preserving the orientations of the polygon(s). One still has to check that the local gluing process does not produce topological or orientability defects when reconstructing the embedding globally. We refer to~\cite{MT01} for details of this construction. This definition stresses that maps are piecewise-linear manifolds, which are equivalent to manifolds in dimension $2$.

\begin{figure}[!ht]
\hfill
\begin{subfigure}{.5\textwidth}
  \centering
  \includegraphics[scale=0.65]{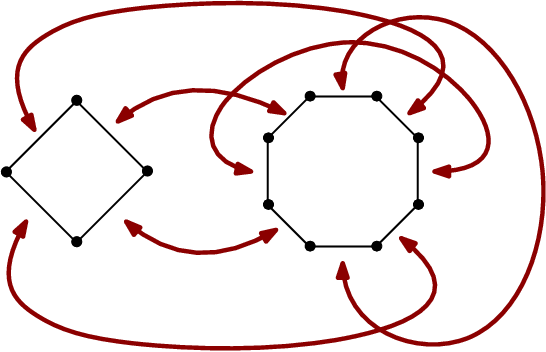}
  %\caption{The description of a map as a polygon gluing.}
  \label{fig:poly_gluing}
\end{subfigure}
\hfill
\begin{subfigure}{.4\textwidth}
  \centering
  \includegraphics[scale=0.42]{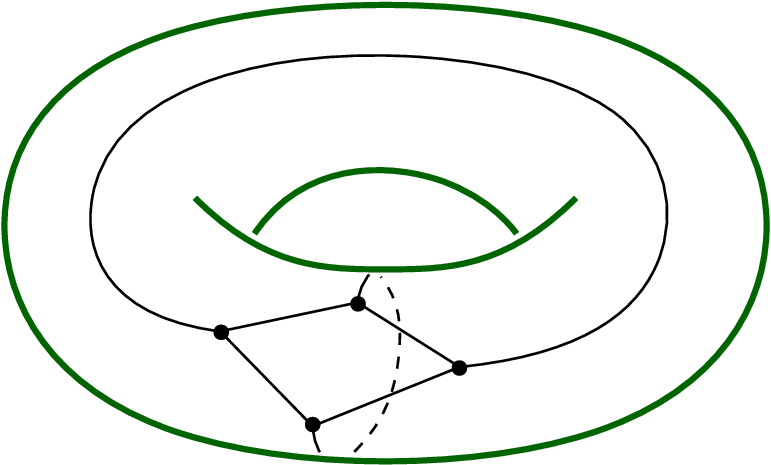}
  %\caption{The corresponding graph embedding.}
  \label{fig:poly_glu_map}
\end{subfigure}
\hfill
\caption{A description of a map as a polygon gluing and the corresponding graph embedding. The matching is represented by the red arrows.}
\label{fig:poly_gluing_fig}
\end{figure}

\subsubsection{As a rotation system}
\label{sssec:rot_sys}

\begin{definition}[Maps as factorization of permutation~\cite{JV90}]
\label{def:rot_sys}
    An orientable map $\mathcal{M}$ is a triple of permutation $(\sigma,\alpha,\phi)\in\mathfrak{S}_{2n}$ such that
    \begin{itemize}
        \item $\alpha$ is an involution without fixed point,
        \item $\sigma\alpha = \phi$.
    \end{itemize}
\end{definition}

The map $\mathcal{M}$ is reconstructed from the triple of permutation as follows. Label the half-edges of $\mathcal{M}$ from $1$ to $2n$. The cycles of $\alpha$ all have length $2$ and form the edges of $\mathcal{M}$. The cycles of $\sigma$ give the cyclic orientation of half-edges around a vertex. Cycles of $\sigma\alpha$ alternate between corners and edges and thus correspond to faces of the map, encoded by $\phi$. 

\begin{figure}[!ht]
    \centering
    \includegraphics[scale=0.6]{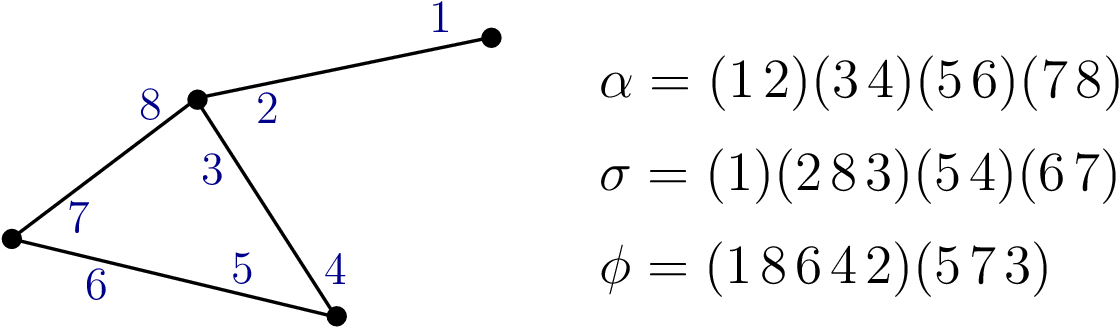}
    \caption{A description of a map as a rotation system.}
    \label{fig:rot_sys_map}
\end{figure}

\begin{remark}
    This definition allows for generalization to other classes of maps. Relaxing the hypothesis on $\alpha$ leads to an algebraic definition of bipartite maps. Increasing the number of permutations $m$ factorizing the identity leads to the notion of constellations that we will encounter in Chapter~\ref{Chap:const}.
\end{remark}

The map $\mathcal{M}$ is connected if the action of $\left<\sigma,\alpha\right>$ on $\left\{1,..,2n\right\}$ is transitive. Otherwise, its number of connected components is given by the number of orbits for this action. Its Euler characteristics is given by
\begin{align}
\label{eq:Euler'}
    \chi(\mathcal{M}) = 2c-2g=V-E+F
\end{align}
where $V=\ell(\sigma),F=\ell(\phi)$ and $E=n$ with $\ell(\rho)$ the number of cycles of the permutation $\rho\in\mathfrak{S}_{2n}$ and $c$ the number of connected components of the map.

\medskip

A triple of permutation $(\sigma,\alpha,\phi)\in\mathfrak{S}_{2n}$ actually encodes a \emph{labeled map}, that is a map endowed with a labeling of its $2n$ half-edges. Clearly, two labeled maps that differ only from the labeling of their half-edges are identical. Said differently, acting by conjugation of an element $\rho\in\mathfrak{S}_{2n}$ on the triple $(\sigma,\alpha)$ leaves an unlabeled map invariant. Therefore the automorphism group of a map $\mathcal{M}$ is
\begin{equation}
    \label{eq:aut_map}
    \text{Aut}(\mathcal{M}) = \{ \rho\in\mathfrak{S}_{2n} \vert \rho\sigma \rho^{-1} = \sigma, \rho\alpha \rho^{-1} = \alpha \}
\end{equation}
This hints at the connection between maps and conjugacy classes of $\mathfrak{S}_{2n}$, which allows to study maps using algebraic tools via the representation theory of the symmetric group $\mathfrak{S}_{2n}$, which we will study in subsection~\ref{ssec:map_sym_group}.

\medskip

To obtain a rooted map, only one half-edge has to be distinguished. Without loss of generality, we can think of the half-edge labeled $1$ as carrying the root~\footnote{It is sometimes convenient to take the opposite convention and think of the highest label half-edge as carrying the root. For example to define an iterable root-deletion procedure.}. This gives a mapping between labeled and rooted maps which consists of forgetting all labels but the root. Since any relabeling of the other $2n-1$ edges doesn't change the rooted map, this mapping is clearly $(2n-1)!$ to $1$. At the combinatorial level, enumerating rooted maps is much easier than unrooted ones. The reason for this is that rooting the map kills all of its symmetries, as stated by the following Lemma.

\begin{lemma}
\label{lem:root_aut}
    Connected rooted maps $\mathcal{M}$ have trivial automorphism group. 
\end{lemma}

\begin{proof}
Fix a labeling of half-edges of $\mathcal{M}$ such that the root carries label $1$ and consider $\rho\in \text{Aut}(\mathcal{M})$. Since $\left<\sigma,\alpha\right>$ acts transitively on half-edges, for all $i$ in $\{1,..,2n\}$ there exist $\rho_i \in \left<\sigma,\alpha\right>$ such that $\rho_i\cdot1=i$. Since $\rho$ fixes the root we have $\rho\cdot1=1$ and therefore
\begin{align}
    \rho\cdot i = \rho\rho_i \cdot 1 = \rho\rho_i\rho^{-1} \cdot 1 = \rho_i \cdot 1 = i 
\end{align}
Hence $\rho$ is the identity and $\text{Aut}(\mathcal{M})$ is trivial.
\end{proof}

On the contrary, to count unrooted maps one has to weight each contribution by the size of its automorphism group~\eqref{eq:aut_map} to make sure each map is counted with a unitary weight. Therefore, the relation between the enumeration of labeled (or rooted) maps and unlabeled ones is non-trivial.  

\subsection{Map enumeration via the Tutte equation}
\label{ssec:map_gen_Tutte}

\subsubsection{Generating series}
\label{sssec:gen_fct}

In combinatorics, one typically wants to count the number of objects with certain properties and give an asymptotic estimation when their size grows to infinity. To tackle these questions, the most ubiquitous tool is \emph{generating functions}. Combinatorial objects come with a notion of size and admit some relations allowing to generate larger objects recursively from smaller ones. The generating function is a formal power series whose coefficients count the number of objects for certain size parameters. Thus, the recursion relation satisfied by the objects studied can in turn be expressed as an equation satisfied by its generating function. In many situations, this equation is an algebraic equation, therefore the combinatorial problem of enumerating objects with given properties can be turned into an (almost entirely) algebraic question thanks to generating functions. While most enumerative applications use $\mathbb{K}=\mathbb{Z}$ or $\mathbb{K}=\mathbb{Q}$, we give the formal definition for an arbitrary ring $\mathbb{K}$.

\begin{definition}[Generating function]
\label{def:gen_fct}
    The generating function $f_{\mathcal{O}}\in \mathbb{K}\left[[t\right]]$ associated with the combinatorial objects $\mathcal{O}$ is the formal power series $f_{\mathcal{O}}(t)= \sum\limits_{k\geq 1} a_k t^k$ such that $a_k$ is the number of objects $O\in\mathcal{O}$ with size $k$. This definition extends naturally to coefficients keeping track of several parameters.
\end{definition}

This definition does not impose any restriction on the coefficients $a_k$. In particular, there is a priori no reason for the series to be summable. Actually, it is common in combinatorics to encounter coefficients with factorial (or exponential) growth with respect to their size parameter, leading to divergent series. One typical such case is the enumeration of labeled objects. Even in that case, we can use usual algebraic operations on the generating series, seen as encoding relations between its coefficients. This set of techniques is called the \emph{symbolic method}. The book~\cite{Flaj09} shows the depth of this simple idea in a wide variety of applications.

\medskip

Throughout this manuscript, we will work with maps or some of their various generalizations. We will mostly be interested in exact enumeration for some type of maps. Typically, we will enumerate maps according to their number of edges, face or vertex profile or genus (not necessarily simultaneously). Since faces of a map can be arbitrarily large, the associated generating function must have countably many indeterminate variables $\{t_i\}_{i\geq 1}$ to keep track of the face profile of a map. The associated ring of polynomials is the \emph{projective limit} of the rings $\left(\mathbb{K}\left[[t_1,\dotsc,t_n\right]]\right)_{n\geq 1}$ and consists of (infinite) sums of monomials in finitely many variables. This fits our purpose since every map has finitely many faces and thus contributes as a finite monomial. When working with infinite families of variables $\{t_i\}_{i\geq 1}$, we will use boldface characters $\textbf{t}$ to represent the countable family $\{t_i\}_{i\geq 1}$ and their weights. For example $\textbf{t}^n = \prod\limits_{i \geq 1} t_i^{n_i}$ (which is well defined since each contribution to the generating function is a monomial in finitely many variables).

\subsubsection{The Tutte equation}
\label{sssec:Tutte_deriv}

We now show how to derive algebraic equations satisfied by the generating function for planar maps counted by their number of edges via purely combinatorial means. This historical result is due to Tutte~\cite{Tutte68}.

\begin{theorem}[Tutte~\cite{Tutte68}]
\label{thm:Tutte}
    The number $m_n$ of rooted planar maps with $n$ edges is 
    \begin{equation}
        \label{eq:plan_map_Tutte}
        m_n = \frac{2\cdot 3^n}{n+2}\Cat_n,
    \end{equation}
\end{theorem}

In Equation~\eqref{eq:plan_map_Tutte}, the coefficients $\Cat_n$ are the famous Catalan numbers. They are the coefficients of the Taylor series expansion of the function $\frac{1-\sqrt{1-4x}}{2x}$ around $x=0$ and can be expressed as $Cat_n = \frac{1}{n+1}\binom{2n}{n}$. 

We denote $m_{n,k}$ the number of rooted planar maps of size $n$ with root face of length $k$ and the corresponding generating function $M_e$ and $\mathfrak{M}^{(e)}_{k,n}$ the set of rooted planar maps with $n$ edges and a root face of length $k$. 

\begin{equation}
    \label{eq:plan_root}
    M_e(t,x) = \sum\limits_{k,n\geq 1} \sum\limits_{\mathcal{M}\in\mathfrak{M}^{(e)}_{k,n}} t^n x^k = \sum\limits_{n,k} m_{n,k}t^n x^k.
\end{equation}

As we will see shortly after, keeping track of the size of the root face turns out to be necessary to write algebraic relations on the generating series. Ultimately, we can forget the introduction of this parameter by specializing its associated variable to one. Such variables are called \emph{catalytic variables}. Therefore, computing the coefficients of the generating series $M_e(t,x)$ allows to compute $m_n$ by extracting the coefficients of $t^n$ and taking the specialization $x=1$.

\medskip

We want to find some recursion allowing us to construct a rooted planar map of size $n$ from smaller maps of size $n'<n$. A natural operation to consider is the root edge deletion. When deleting the root-edge, we are necessarily in one of the following cases:
\begin{itemize}
    \item[(Join)] The same face borders both sides of the root edge. In this case, removing the root edge disconnects the map into two maps of smaller sizes.
    \item[(Cut)] Two distinct faces are bordered by the root edge. When the root edge is deleted, both faces are merged into a single larger face.
\end{itemize}

\begin{figure}[!ht]
\begin{subfigure}{.38\textwidth}
  \centering
  \includegraphics[scale=0.28]{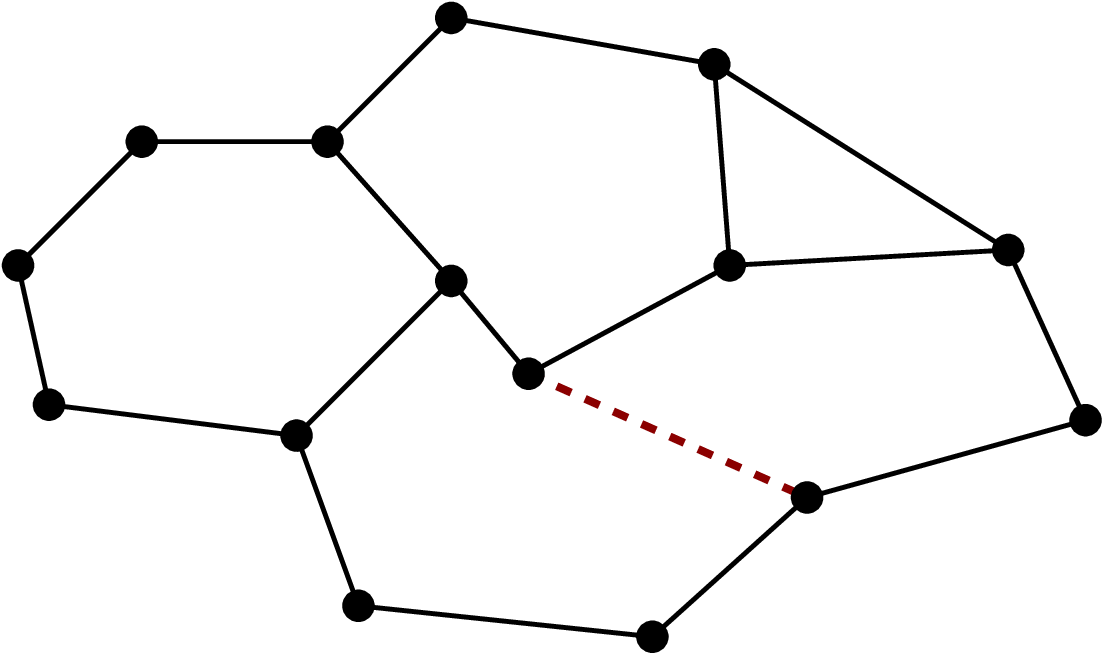}
  \caption{Join case.}
  \label{fig:join_case}
\end{subfigure}
\hspace{20pt}
\begin{subfigure}{.52\textwidth}
  \centering
  \includegraphics[scale=0.30]{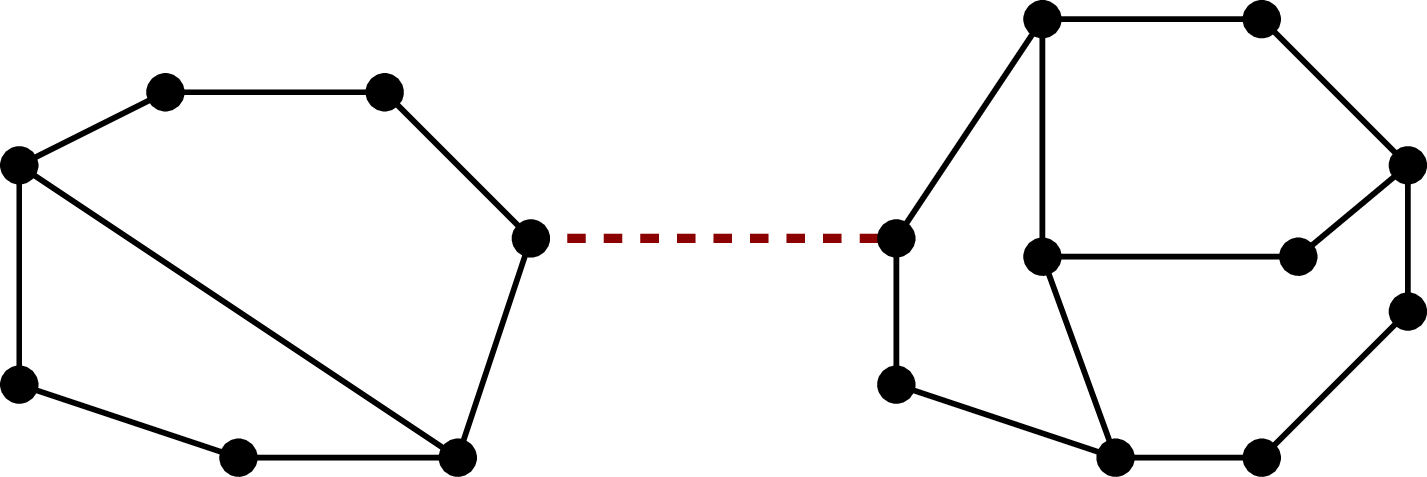}
  \caption{Cut case.}
  \label{fig:cut_case}
\end{subfigure}
\hfill
\caption{Illustrating the two situations that can occur when deleting an edge.}
\label{fig:cut_and_join}
\end{figure}

Now, we reverse this process and consider how many ways we can add an edge in smaller maps for each of these two cases.

\begin{itemize}
    \item For the join case, starting from two rooted planar maps, there is a unique way to add an edge between them while respecting the root face orientation. Hence, for each pair of maps contributing to monomials $t^{n_1}x^{k_1}$ and $t^{n_2}x^{k_2}$ there is a unique way to form a map contributing to the monomial $t^{n_1+n_2+1}x^{k_1+k_2+2}$. This contribution \textit{joins} two faces to form a new one.
    \item In the cut case, the number of ways we can split the root face into two smaller faces depends on the size of the root face. This is the reason why it is necessary to keep track of the root face to write an algebraic relation on $M_e(t,x)$. If the root face has length $k$, then for each $k'\in\llbracket 1,k+1\rrbracket$ there is a unique way to add an edge and obtain a root face of length $k'$. Therefore each map contributing to the monomial $t^nx^k$ generates a unique map contributing to the monomial $t^{n+1}x^{k'}$. This contribution \textit{cuts} one faces into two distinct ones.
\end{itemize}

Finally, one has to include the trivial map with weight $1$, which cannot be obtained from a recursive procedure. The above reasoning translates algebraically as follows

\medskip

Algebraically, the reasoning above leads to the following equation

\begin{align}
    \label{eq:Tutte_rel}
    M_e(t,x) &= 1 + \underbrace{tx^2\left(\sum\limits_{k_1,n_1\geq 1} \sum\limits_{\mathcal{M}_1\in\mathfrak{M}^(e)_{k_1,n_1}} t^{n_1}x^{k_1} \right)\left(\sum\limits_{k_2,n_2\geq 1} \sum\limits_{\mathcal{M}_2\in\mathfrak{M}^{(e)}_{k_2,n_2}} t^{n_2}x^{k_2} \right)}_{\text{Join}} \\ 
    &+ \underbrace{\sum\limits_{k,n\geq 1} \sum\limits_{\mathcal{M}\in\mathfrak{M}^{(e)}_{k,n}} t^{n+1)}\left(\sum\limits_{k'=1}^{k+1} u^{k'}\right)}_{\text{Cut}} \nonumber\\
    &= 1 + tx^2M_e(t,x)^2 + tx\frac{uM_e(t,x)-M_e(t,1)}{x-1}. \nonumber
\end{align}
 
This equation is polynomial in $M_e(t,x)$ and involves the catalytic variable $x$ and the specialization $M_e(t,1)$. At first sight, it is not at all obvious how to approach this equation and solve it. Historically, Tutte proceeded to guess and check the form of the solution. Since then, a general method has been developed by Bousquet-M\'{e}lou and Jehanne~\cite{BMJ06} to solve equations involving one catalytic variable. Solving this equation leads to the expression of $m_n$ as given in Theorem~\ref{thm:Tutte}.

\medskip

This root-edge deletion procedure can be applied to a wide variety of classes of maps, the main constraint being that the class of map should be stable under edge deletion. When applied, it leads to the above distinction of the cut and join cases depending on how the root edge is bordered. For this reason, the type of equations obtained are sometimes referred to as \emph{cut-and-join equations}. 

\subsection{Map enumeration via representation theory}
\label{ssec:map_sym_group}

In this subsection, we show how to apply the representation theory of the symmetric group for the enumeration of weighted maps.

\subsubsection{Representation of the symmetric group}
\label{sssec:rep_sym_grp}

Maps can be described as a triple of permutation factorizing the identity, up to a relabeling of its half-edges. Using this description of maps, we can make use of the representation theory of the symmetric group to deduce information on maps and their generating series. While most results presented here are general representation theory results, we will narrow down our presentation to the relevant case of the symmetric group.

\medskip

For $G$ a finite group, a representation of $G$ is a pair $(V,\rho)$ such that $V$ is a (complex) finite dimensional vector space and $\rho:G\rightarrow GL(V)$ is a group homomorphism. A representation is \emph{irreducible} if it cannot be written as a direct product of two representations. For an element $g \in G$, we define its conjugacy class $C_g$ as the set of elements that are conjugate to $g$
\begin{equation}
    \label{eq:conj_class}
    C_g = \{ hgh^{-1} \vert h \in G \},
\end{equation}

The character $\chi^\rho$ of an element $g$ is given by the trace of $\rho(g)\in GL(V)$ i.e. $\chi^\rho(g) = \Tr \rho(g)$. Since the trace is cyclically invariant, the characters are \emph{class functions}: they are  constant on conjugacy classes $\chi^\rho(hgh^{-1})=\chi^\rho(g)$. They completely characterize representations: two representations are isomorphic if and only if they have the same characters~\cite{Serre_grp}. Note that $\chi^\rho(1_G) = \dim V$, therefore we will sometimes abusively call $\rho$ a representation, omitting the vector space $V$ it acts on (which amounts to considering a representation up to isomorphism of that vector space). A character $\chi^\rho$ is irreducible if its representation $\rho$ is irreducible. The vector space of class functions is equipped with a scalar product given by 
\begin{equation}
\label{eq:scal_char}
    \langle \chi_1,\chi_2 \rangle = \frac{1}{\vert G\vert}\sum\limits_{g\in G} \chi_1(g)\chi_2(g^{-1}).
\end{equation}

In particular, a representation is irreducible if and only if $\langle \chi,\chi \rangle =1$~\cite{Serre_grp}.

\medskip

In the following, a key role will be played by the group algebra $\mathbb{C}G$. It is the vector space with basis labeled by elements of $G$ $V=\text{span}\{e_g \vert g\in G \}$. By a slight abuse of notation, we will often identify the vector basis $e_g$ with the group element $g$ itself. The law of $G$ induces a natural group action on the group algebra via $g\cdot e_h = e_{gh}$, which defines a product law for the group algebra. The character $\chi$ of this representation corresponds to the size of the stabilizer of an element $h\in G$. Since the action of $g$ on itself has no non-trivial fixed point we have $\chi(g) = \begin{cases} \vert G \vert \hspace{5pt}\text{if}\hspace{5pt}g=1_G \\ \vert 0 \vert \hspace{5pt}\text{otherwise}\hspace{5pt}\end{cases}$. It follows that if $\chi^n$ is the irreducible character associated with a representation of dimension $n$ then we have $\langle \chi, \chi^n \rangle = \chi^n(1_G) = n$. Therefore, the group algebra $\mathbb{C}G$ contains $n$ copies of any irreducible representation of dimension $n$. In particular, this implies that there are finitely many irreducible representations of $G$ since $\mathbb{C}G$ has finite dimension $\vert G\vert$. We will admit that irreducible characters form a basis of class functions and that the number of irreducible characters is equal to the number of conjugacy classes of $G$.

\medskip

We will now specialize to the symmetric group $\mathfrak{S}_n$. For an element $\sigma\in\mathfrak{S}_n$, we call its \emph{cycle type} the multi-set $\lambda = \{\lambda_i\}_i$ given by the length of its cycles. Since each element of $\llbracket 1,n \rrbracket$ appears exactly once in $\sigma$, we have $\sum\limits_i \lambda_i = n$. Ordering $\lambda$ such that the $\lambda_i$ are weakly decreasing, cycle types are in bijection with \emph{partitions} of $n$.

\begin{definition}[Partition]
\label{def:part}
    A partition $\lambda$ is a finite non-increasing sequence of positive integer $(\lambda_i)_{i\geq 1}$. Each $\lambda_i$ is called a part of length $\lambda_i$ of the partition. The number of parts is called the length of the partition and is denoted $\ell(\lambda)$. We say $\lambda$ is a partition of $n$, denoted $\lambda \vdash n$ if $\sum\limits_i \lambda_i = n$. If $\lambda \vdash n$ has parts $m_i$ parts of length $i$, then we write $\lambda = \left[ 1^{m_1} \dotsc n^{m_n} \right]$ omitting the length not appearing in $\lambda$.
\end{definition}

Partitions can be represented graphically using Ferrer diagrams. A partition is represented by rows of boxes such that the $k$-th row has $\lambda_k$ boxes and the rows are aligned to the left. In this manuscript, unless otherwise stated, we will use the French convention and will order rows from bottom to top, i.e. the largest part being at the bottom of the Ferrer diagram.

\begin{figure}[!ht]
\centering
\includegraphics[scale=0.4]{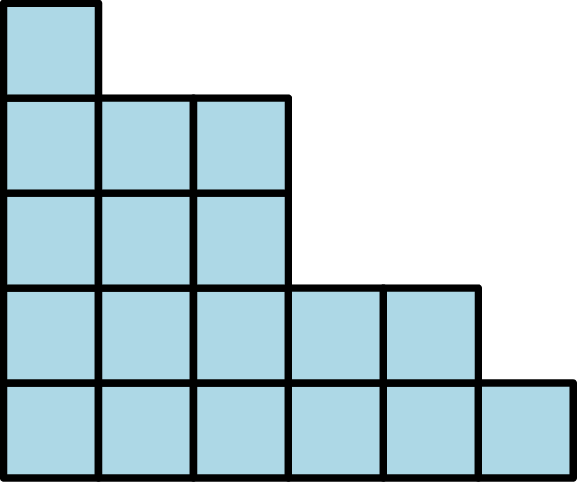}
\caption{The Ferrer diagram of the partition $\lambda =\left[6,5,3,3,1\right]$.}
\label{fig:ex_partition}
\end{figure}

Acting by conjugation on an element $\sigma\in \mathfrak{S}_n$ amounts to a relabeling of elements of $\llbracket 1,n \rrbracket$ hence it preserves the cycle type of a permutation. Therefore, conjugacy classes of $\mathfrak{S}_n$ can be labeled by partitions of $n$. We denote for $\lambda \vdash n$
\begin{equation}
\label{eq:conj_class_sym}
    C_\lambda = \{ \sigma\in \mathfrak{S}_n \vert \sigma \text{    has cycle type    } \lambda \}.
\end{equation}

For $\lambda = \left[ 1^{m_1} \dotsc n^{m_n} \right]$ and $\sigma \in C_\lambda$, the centralizer of $\sigma$ is the set of elements of $\mathfrak{S}_n$ commuting with $\sigma$, i.e. elements leaving $\sigma$ invariant when acting by conjugation. These elements must map cycles of length $\lambda_i$ between them. Therefore, such an element can act in two separate ways on each cycle of $\sigma$:
\begin{itemize}
    \item Rotate all the labels of a cycle, yielding a factor $i$ for each cycle of length $i$.
    \item Exchange all labels of cycles of the same length between them, yielding a factor $m_i!$ for all $i$.
\end{itemize}
It follows that the cardinality of the centralizer of $\sigma$ is $z_\lambda = \prod\limits_{i} i^{m_i}m_i!$. Therefore, the cardinality of the conjugacy class $C_\lambda$ is $\frac{n!}{z_\lambda}$.
\vspace{15pt}

With one more slight abuse of notation, a conjugacy class $C_\lambda$ can be identified with the sum of its elements as an element of the group algebra
\begin{equation}
    \label{eq:conj_grp_alg}
       C_\lambda \equiv \sum\limits_{\sigma\in C_\lambda} \sigma \in \mathbb{C}G.
\end{equation}

We can already get a glimpse of how conjugacy classes can be used to count maps. A labeled map with vertex profile $\lambda \vdash 2n$ and face profile $\nu \vdash 2n$ is characterized by $\sigma \in C_\lambda, \phi^{-1} \in C_\nu$, and $\alpha\in C_{\left[2^n\right]}$ a product of $n$ disjoint transposition under the condition $\sigma\alpha\phi^{-1} = \id_{\mathfrak{S}_n}$. Thus, if the product of two conjugacy classes $C_{\lambda_1}C_{\lambda_2}$ can be expanded over conjugacy classes again, labeled maps with face profile $\lambda$ and face profile $\phi$ are counted by the coefficients of the expansion of $C_{\id_{\mathfrak{S}_n}}$ in the product $C_\lambda C_{\left[2^n\right]} C_\nu$. To finally derive such a relation for the product of conjugacy classes, we need one last result from representation theory.

\begin{proposition}[\cite{Stanley}]
The conjugacy classes $C_\lambda$ generate the center of the group algebra $\mathcal{Z}(\mathbb{C}\mathfrak{S}_n)$
\end{proposition}

\begin{proof}
For any conjugacy class $C_\lambda$ and any element $g$ we have $gC_\lambda g^-1 = C_\lambda$, therefore $C_\lambda \subset \mathcal{Z}(\mathbb{C}\mathfrak{S}_n)$. It follows that any linear combination of $C_\lambda$ is also in the center of the group algebra. Conversely, any element $x\in \mathcal{Z}(\mathbb{C}G)$ can be expanded over the group algebra basis as $x = \sum\limits_{\sigma \in \mathfrak{S}_n} x_\sigma \sigma$. Since $x\in \mathcal{Z}(\mathbb{C}\mathfrak{S}_n)$, for any $\rho \in \mathfrak{S}_n$ we have $\rho x \rho^{-1} = x $ i.e. $x_\sigma=x_{\rho\sigma\rho^{-1}}$. Therefore, coefficients of $x$ in the same conjugacy classes are equal. It follows that $x$ is a linear combination of $C_\lambda$.
\end{proof}

A second basis of the center $\mathcal{Z}(\mathbb{C}\mathfrak{S}_n)$ is given by the following elements.
\begin{proposition}[\cite{Stanley}]
    The set $\{ P_\lambda \vert \lambda \vdash n \} $ is a basis of $\mathcal{Z}(\mathbb{C}\mathfrak{S}_n)$ where $P_\lambda = \frac{d_\lambda}{n!} \sum\limits_{\sigma \in \mathfrak{S}_n} \chi^\lambda(\sigma^{-1}) \sigma$ with $d_\lambda$ the dimension of the irreducible representation $\lambda$.
\end{proposition}
\begin{proof}

Since characters are invariant by conjugation, for any $\rho \in \mathfrak{S}_n$ we have the following relation $\rho P_\lambda \rho^{-1} = \frac{d_\lambda}{n!} \sum\limits_{\sigma \in \mathfrak{S}_n} \chi^\lambda(\rho\sigma^{-1}\rho^{-1}) \sigma\rho^{-1} = P_\lambda$ as any conjugacy class $C_\lambda$ is stable under conjugation. Moreover, we have $\sum\limits_{\lambda \vdash n} P_\lambda = \frac{1}{n!}\sum\limits_{\substack{\lambda \vdash n\\ \sigma \in \mathfrak{S}_n}} d_\lambda \chi^\lambda(\sigma)\sigma $. Since $\sum\limits_{\lambda \vdash n} d_\lambda \chi^\lambda(\sigma)$ is the character of the group algebra, it vanishes unless $\sigma= \id $ where it cancels the factor $\frac{1}{n!}$. Therefore $\sum\limits_{\lambda \vdash n} P_\lambda = \id_{\mathfrak{S}_n}$. It follows that any element of the center $x$ can be expanded as $x = \sum\limits_\lambda xP_\lambda$.
\end{proof}
This basis is particularly convenient for computations as it is idempotent. Indeed, for $\lambda,\nu \vdash n$ we have
\begin{align}
    \label{comp:idem_basis}
    P_\lambda P_\nu &= \frac{d_\lambda d_\nu}{n!^2}\sum\limits_{\substack{\sigma_1,\sigma_2}} \chi^\lambda(\sigma_1^{-1})\chi^\nu(\sigma_2^{-1}) \sigma_1\sigma_2 \\
    &= \frac{d_\lambda d_\nu}{n!^2}\sum\limits_{\sigma} \sum\limits_{\sigma_1} \chi^\lambda(\sigma_1^{-1})\chi^\nu(\sigma_1)\chi^\nu(\sigma^{-1}) \sigma \nonumber \\
    &= \delta_{\mu\nu} P_\lambda \nonumber \qquad \text{by orthonormality of irreducible characters}
\end{align}

The change of basis between $C_\lambda$ and $P_\nu$ is straightforward by factorizing conjugacy classes in the expression of $P_\nu$. Conversely, it can be checked that
\begin{equation}
\label{eq:conj_to_proj}
    C_\lambda = \frac{n!}{z_\lambda} \sum\limits_{\nu \vdash n} \frac{\chi^\nu_\lambda}{d_\nu} P_\nu
\end{equation}
where $\chi^\nu_\lambda$ is the character $\chi^\nu$ evaluated on any element of $C_\lambda$.
Now, we know that $C_\lambda$ form a basis of $\mathcal{Z}(\mathbb{C}\mathfrak{S}_n)$ and via the idempotent basis $\{P_\lambda\}$ we have for two conjugacy classes $C_\lambda C_\nu$
\begin{align}
\label{eq:prod_rule}
    C_\lambda C_\nu &= \frac{n!^2}{z_\lambda z_\nu} \sum\limits_{\mu \vdash n} \frac{\chi^\mu_\lambda \chi^\mu_\nu}{d_\mu^2 } P_\mu \\
    &= \frac{n!}{z_\lambda z_\nu} \sum\limits_{\mu \vdash n} \frac{\chi^\mu_\lambda \chi^\mu_\nu}{d_\mu}  \sum\limits_{\theta \vdash n} \chi^\mu_\theta C_\theta \nonumber
\end{align}
which generalized straightforwardly to any product of conjugacy classes. We can now apply this result to map enumeration. The number $M_{\lambda,\mu}$ of labeled maps with face profile $\lambda \vdash 2n$ and vertex profile $\mu\vdash 2n$ is given by the number of factorization of $\id_{\mathfrak{S}_n} = \sigma \alpha \phi^{-1}$ such that $\sigma \in C_\lambda, \phi^{-1} \in C_\nu$, and $\alpha\in C_{\left[2^n\right]}$. Therefore
\begin{align}
\label{eq:map_count}
 M_{\lambda,\nu} &= \left[C_{\left[1^{2n}\right]}\right] C_\lambda C_{\left[2^n\right]}C_\nu \\
 &= \frac{n!}{z_\lambda z_\mu 2^n} \sum\limits_{\mu \vdash n} \frac{\chi^\mu_\lambda \chi^\mu_{\left[2^n\right]} \chi^\mu_\nu}{d_\mu} \nonumber
\end{align}

Note that we counted \textit{all} maps, and not only connected ones. The number of connected maps can be computed from this expression since the generating function of connected and non-necessarily connected maps are related by a simple relation, the connected generating function being the logarithm of the non-connected one. 

\medskip

\begin{remark}
    The formula~\eqref{eq:map_count} is a particular case of a general result of the representation theory of finite group known as the \emph{Frobenius formula} (see e.g.~\cite[Sec. $A.1.3$]{LZ_book}). The Frobenius formula can be applied to other type of maps. For example, relaxing the condition on $\alpha$ gives the number of bipartite maps with specified vertices and face profiles. Another generalization is to increase the number of permutation factorizing to the identity, which counts \emph{constellations} we will encounter in Chapter~\ref{Chap:const}.
\end{remark}

\subsubsection{Application to weighted maps}
\label{sssec:weighted_map}

Using the representation theoretic tools of the last section, weights can be added to the generating function of maps to keep track of some parameters like the face or vertex profiles. As an example, we consider a generating function for (non-necessarily connected) rooted maps where we track
\begin{itemize}
    \item the number of faces of length $k$ via variables $t_k$
    \item the number of vertices of degree $i$ via variables $q_i$
    \item the number of edges via the variable $t$
\end{itemize}
Therefore, the generating function $M$ for these weighted maps writes
\begin{equation}
    \label{eq:gen_fct_maps}
    M(t,\textbf{p},\textbf{q}) = \sum\limits_{\mathcal{M}\in\mathfrak{M}_r} t^{E(\mathcal{M}) } \textbf{p}^{ F(\mathcal{M})}\textbf{q}^{V(\mathcal{M})},
\end{equation}
where $\mathcal{M}\in\mathfrak{M}_r$ is the set of rooted maps and $\textbf{p}^{F(\mathcal{M})} = \prod\limits_{i \geq 1} p_i^{F_i(\mathcal{M})}$.

\medskip

Therefore, if $\mathcal{M}$ has $n$ edges and $\sigma$ has cycle type $\lambda = \left[ 1^m_1 \dotsc 2n^{m_{2n}} \right]$ then the weight of $\mathcal{M}$ in variables $\textbf{p}$ is $\prod\limits_{i=1}^{2n} p_i^{m_i}$. To make the connection with symmetric functions, we will work with power-sum variables. We introduce variables $\textbf{x}$ and $\textbf{y}$ such that 
\begin{equation}
    \label{eq:change_power_sum}
    p_k = \sum\limits_{i \geq 1} x_i^k, \qquad q_k = \sum\limits_{j \geq 1} y_j^k.
\end{equation}
In the variables $\textbf{x}$ a map with a face profile $\lambda$ now has weight $p_\lambda(\textbf{x}) = \prod_{i=1}^{\ell(\lambda)} \sum\limits_{r \geq 1} x_r^{\lambda_i} $. While the reason for this choice might seem arbitrary at the combinatorial level, we will see that it is natural from the point of view of matrix models (see Sec.~\ref{sec:rand_mat}).

\medskip

A labeled map $\mathcal{M}$ is uniquely defined by a triple of permutations $(\sigma,\alpha,\phi)$ factorizing to the identity. The face and vertex profile are encoded in the cycle type of $\sigma$ and $\phi$ respectively. Since the number of maps with face profile $\lambda$ and vertex profile $\nu$ has been computed in Equation~\eqref{eq:map_count}, we get
\begin{equation}
    \label{eq:map_sym_fct_PS}
    M(t,\textbf{x},\textbf{y}) = \sum\limits_{n \geq 0} t^n \sum\limits_{\lambda,\nu \vdash n} M_{\lambda,\nu} p_\lambda(\textbf{x}) p_\nu(\textbf{y}).
\end{equation}
Making the coefficients $M_{\lambda,\nu}$ explicit amounts to performing the change of basis between the power-sum $p_\lambda$ to the Schur functions $s_\lambda$. Indeed, the Schur functions in variables $\textbf{x}$ are expressed in that basis as
\begin{equation}
    \label{eq:Schur_fct}
    s_\lambda(\textbf{x}) = \sum\limits_{\substack{ n \geq 0 \\ \mu \vdash n}} \frac{\chi^\lambda_\mu}{z_\mu} p_\mu(\textbf{x}).
\end{equation}
Therefore the generating series of face and vertex weighted maps can be expanded over the Schur functions as
\begin{equation}
    \label{eq:map_sym_fct}
    M(t,\textbf{x},\textbf{y}) = \sum\limits_{n \geq 0} t^n \sum\limits_{\mu \vdash n} \frac{ \chi^\mu_{\left[2^n\right]} }{n_\mu} s_\mu(\textbf{x})s_\mu(\textbf{y}).
\end{equation}

The Schur functions will play an important role in Section~\ref{sec:KP_hier} as they are related to the integrability properties of these generating series.

\subsubsection{The Murnaghan-Nakayama rule}
\label{sssec:MN_rule}

Equation~\eqref{eq:map_count} gives us a formula for the number of maps with a given vertex and face profile. To make use of this formula, we must compute characters of any irreducible representation of $\mathfrak{S}_n$. We present here a combinatorial method to compute the value of the character of an irreducible representation $\chi^\lambda$ on a conjugacy class $\mu$ known as the Murnaghan-Nakayama rule. We refer to~\cite[Section $7.17$]{Stanley} for a proof of this relation. The Murnaghan-Nakayama rule allows to compute characters evaluated on a given conjugacy class by successive insertions of \emph{border strips}.

\begin{definition}[Skew-partition]
    For $n\geq 1$ and two partitions $\lambda,\nu \vdash n$ such that for all $i>0$  $\lambda_i \leq \nu_i$, the skew-partition $\nu\backslash\lambda $ is the box-diagram of boxes of $\nu$ that are not in $\lambda$. 
\end{definition}

\begin{definition}[Border strip]
    A skew-partition $\nu\backslash\lambda $ is a border-strip if it is connected and doesn't contain any $2\times 2$ diagram. The height of a border strip is the number of rows minus one. The length of a border strip is its total number of boxes.
\end{definition}

\begin{figure}
\hfill
\begin{subfigure}{.4\textwidth}
  \centering
  \includegraphics[scale=0.75]{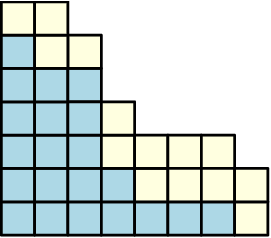}
  \caption{A skew-partition.}
  \label{fig:skew_part}
\end{subfigure}
\hfill
\begin{subfigure}{.4\textwidth}
  \centering
  \includegraphics[scale=0.75]{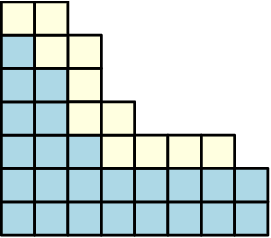}
  \caption{A border strip.}
  \label{fig:border_strip}
\end{subfigure}
\hfill
\caption{Example of Ferrer diagram of skew-partition and border strip, given by yellow boxes.}
\label{fig:example_skew_border}
\end{figure}

The characters $\chi^\lambda$ evaluated on a conjugacy class $\nu$ can then be computed recursively on the length of $\nu$ as follows.

\begin{theorem}[Murnaghan-Nakayama rule]
    For $n\geq 0$ two partitions $\lambda,\nu \vdash n$ we have
\begin{equation}
    \label{eq:MN_rule_rec}
    \chi^\lambda_\nu = \sum\limits_{\xi \in BS(\lambda,\nu_1)} (-1)^{\text{ht}(\xi)}\chi^{\lambda\backslash\xi}_{\nu\backslash\nu_1}
\end{equation}
where $BS(\lambda,\nu_1)$ is the set of border strip of $\lambda$ with length $\nu_1$ and $\chi^\emptyset_\emptyset = 1$.
\end{theorem}

This relation computes $\chi^\lambda_\nu$ by successive removal of the parts of $\nu$ as border strips in $\lambda$. It can be recast in a non-recursive way by applying the relation~\eqref{eq:MN_rule_rec} until the empty partition is reached after $\ell(\nu)$ iterations. To keep track of the successive ribbon insertions, the boxes of the $i$-th ribbon are labeled $i$. By construction, the labels of the boxes of $\lambda$ are weakly increasing in columns and rows. Here, we required $\nu$ to be a partition which led to adding ribbons of decreasing size, but it can be shown that the global sign factor is independent of the order of insertion of ribbons. Hence this decomposition is actually valid for any sequence of ribbon insertions.

\section{Random matrix models}
\label{sec:rand_mat}

Random matrix theory could at first sight seem unrelated to maps. And indeed, the first application of random matrices in theoretical physics dates to Wigner~\cite{Wigner55} in 1928 where he used random matrix to model the interaction of heavy atom nuclei. The connection with maps only became apparent almost $50$ years later through the \emph{Feynman path integral} of matrix field theory. The Feynman path-integral is an approach of quantum field theory that describes the correlation of a field at different points in space-time. In this formalism, the central object of a field that encodes all of its correlation functions is the partition function, similar to the one usually encountered in statistical physics, albeit more complicated due to the functional aspects related to the space-time coordinates dependency. The integral appearing in the perturbative expansion of the partition function (and thus of its correlators) can be described graphically by diagrams decorated to encode all necessary information. In their article, 't Hoft~\cite{HOFT74} and shortly after Brézin, Itzykson, Parisi and Zuber~\cite{BIPZ78} apply the path integral formalism to the $SU(N)$ gauge theory and notice that the diagrams contributing at leading order in $N$ are exactly planar maps. This leads to a renewed interest in random matrix models in the following $20$ years as a \emph{zero dimensional field theory}, a field theory on a space-time reduced to a single point i.e. a probability distribution on a matrix ensemble. Its perturbative expansion can be interpreted as a partition function for a quantum theory of spacetime in $2$ dimensions~\cite{DFGiZJ}. Conversely, this connection brought methods from field theory that can in turn be applied to the enumeration of maps. This section shows some of the connections between maps and random matrices and revisits some of the combinatorial notions of the previous section through a field theoretic lens via matrix models. It is based on the reviews on matrix models~\cite{Zvon97,EKR_Review,DFGiZJ}. We refer to~\cite{Ox_RMT} for applications of random matrix theory beyond random geometry.

\subsection{Gaussian Hermitian matrices}
\label{ssec:Gauss_Herm}

Free fields are a central object in quantum field theory. They are among the rare cases where the path integral of the system can be computed exactly~\footnote{The other one being integrable systems which can often be brought to a free field by a clever change of variable}. Interacting fields are not as tractable and we must often content ourselves with a perturbative expansion around free fields. While the path-integral is much tamer in zero dimensions since there are no functional involved, this still holds true in this framework as there is no obvious way to evaluate the corresponding integrals. In this subsection, we define the Gaussian measure on the space of Hermitian matrices and give the corresponding Wick's theorem.

\medskip

We denote $\mathcal{H}_N$ the space of $N\times N$ Hermitian matrix model. A matrix $\mathcal{H}\in H_N$ is described by $N^2$ real parameters: $N$ real diagonal entries $h_{i}$ and $\frac{N(N-1)}{2}$ complex entries $h_{ij}=x_{ij}+iy_{ij}$ off-diagonal. Therefore $\mathcal{H}_N$ is isomorphic to $\mathbb{R}^{N^2}$ and its Lebesgue measure is given by
\begin{equation}
    \label{eq:Leb_Herm}
    dH = \prod\limits_{i=1}^N dh_{ii}\prod\limits_{1\leq i<j\leq N} dx_{ij}dy_{ij}.
\end{equation}

The quadratic form on $\mathcal{H}_N$ is 
\begin{align}
    \Tr(H^2) &= \sum\limits_{1\leq i,j\leq N} h_{ij}h_{ji} = \sum\limits_{1\leq i \leq N} h_{i}^2 + 2\sum\limits_{1\leq i,j\leq N} x_{ij}^2 + y_{ij}^2,
\end{align}

and we define the corresponding Gaussian integral as
\begin{align}
    \label{eq:Gauss_Herm}
    Z^{(0)}_{\mathcal{H}_N} &= \int_{\mathcal{H}_N} dH\exp\left(-\frac{N}{2}\Tr(H^2)\right) \\
                            &= (N\pi)^{-\frac{N^2}{2}}2^{N}.\nonumber
\end{align}

The measure 
\begin{equation}
    \label{eq:norm_meas_Herm}
    d\mu(H) = \frac{1}{Z^{(0)}_{\mathcal{H}_N}} dH,
\end{equation}
is normalized and positive for any $H\in \mathcal{H_N}$, therefore, it is a probability distribution on the space of Hermitian matrices. For a function $f:\mathcal{H}_N \rightarrow \mathbb{C}$, we denote $<f>$ its expectation value for the probability measure $\mu$. A quick computation shows that for quadratic polynomials in matrix entries, we have

\begin{equation}
    \label{eq:2pt_Gauss}
    <h_{ij}h_{kl}> = \frac{1}{N}\delta_{il}\delta_{kj}.
\end{equation}

Monomials of degree greater than $2$ can be computed inductively via integration by parts. This leads to Wick's Theorem.

\begin{theorem}[Wick's Theorem]
\label{thm:Wick}
The expectation value of the monomial $h_{i_1j_1}..h_{i_nj_n}$ is
\begin{equation}
\label{eq:Wick}
    <h_{i_1j_1}..h_{i_nj_n}> = \begin{cases} 0 &\text{  if   } $n$ \text{  odd} \\ \frac{1}{N}\delta_{il}\delta_{kj} &\text{  if    } $n=2$ \\ \sum\limits_{M\in \mathcal{M}_k} \prod\limits_{\substack{(i,j)(k,l)\\ \text{s.t.} M\cdot(i,j)=(k,l)}} <h_{ij}h_{kl}>  &\text{  if   } $n=2k$ \end{cases},
\end{equation}
where $\mathcal{M}_k$ is the set of matchings on the set with $2k$ elements.
\end{theorem}

By linearity, this Theorem allows us to compute the expectation value of any polynomial function for a Gaussian measure.

\subsection{The Hermitian matrix model}
\label{ssec:Herm_model}

Anticipating on later computations, we rescale the kinetic term coupling constant as $N \rightarrow \frac{N}{t}$ with $t$ a positive real number. This rescaling has no incidence beyond a straightforward change of the normalization constant $Z^{(0)}_{\mathcal{H}_N}$. We now consider the following non-Gaussian measure
\begin{equation}
    d\nu(H) = \frac{1}{Z^{(0)}_{\mathcal{H}_N}}\exp\left(-\frac{N}{2t}\Tr H^2+ N V(H)\right),
\end{equation}

with $V(H)$ the \emph{potential} defined as 
\begin{equation}
    V(H) = \sum\limits_{i=3}^d \frac{t_k}{k}  \Tr H^k,
\end{equation}
for some integer $d\geq 1$, where the $t_k$ are the coupling constants of the model. This measure leads to the partition function

\begin{equation}
\label{eq:part_herm}
    Z\left[t,\bf{t_k}\right] = \frac{1}{Z^{(0)}_{\mathcal{H}_N}}\int_{H_N} dH\exp\left(-\frac{N}{2t}\Tr H^2 + N\sum\limits_{k=3}^d \frac{t_k}{k} \Tr H^k\right).
\end{equation}

The nature of this partition function is very different than the Gaussian one. Not only there is no clear way on how this integral could be computed, but it is non-analytic around zero. This argument was pointed out by Dyson in the context of quantum electrodynamics~\cite{Dyson52} and is inherent to interacting field theories. For example, if we take only $t_4$ to be non-zero, then the integral is convergent for $t_4<0$ but diverges for $t_4<0$ and therefore has zero radius of convergence around zero, even in the scalar case $N=1$ ! 

\medskip

Physicists circumvent this issue somewhat boldly by considering the \emph{perturbative expansion} around the free field, that is, by expanding in a power series all the terms of the potential, leaving only the quadratic term in the measure, commuting the sum and the integral and applying Wick's Theorem to compute the expectation values for each term of the power series. Essentially, this amounts to performing a Taylor expansion around the point $t_i=0$, where we in all rigor couldn't because the partition function isn't analytic around this point. And indeed, as we shall see through the example of the Hermitian matrix model hereafter, the perturbative series is divergent and thus cannot be naively resummed to recover the partition function we started from. How to make sense of the perturbative expansion and how to reconstruct (some information on) the original partition has been among the most crucial questions surrounding quantum field theory as it is one of the only generic methods available to perform computations in a quantum field theory. This turns out to be possible in some cases via \emph{resurgence theory} which gives tools to resum divergent series. We refer the interested reader to the introductions~\cite{Marino12,Dorigoni14} for more details on resurgence theory. In the rest of this manuscript, we will leave these issues aside and the matrix integrals should be understood as perturbative expansions around the free field.

\medskip

We are now ready to expand the Hermitian matrix model perturbatively. The potential develops as
\begin{equation}
\label{eq:V_exp}
    \exp(NV(H)) = \sum\limits_{n\geq 0} \frac{N^n}{n!}\sum\limits_{k_1,..k_n \in\{1,..,d\}} \frac{t_{k_1}...t_{k_n}}{k_1...k_n} \Tr H^{k_1}...\Tr H^{k_n}.
\end{equation}

Therefore the perturbative expansion $Z^{(\text{pert})}_{\mathcal{H}_N}$ of the Hermitian matrix model is
\begin{equation}
\label{eq:pert_herm}
    Z^{(\text{pert})}_{\mathcal{H}_N}\left[t_k\right] = \sum\limits_{n\geq 0} \frac{N^n}{n!}\sum\limits_{k_1,..k_n \in\{1,..,d\}} \frac{t_{k_1}...t_{k_n}}{k_1...k_n} \left\langle\Tr H^{k_1}...\Tr H^{k_n}\right\rangle.
\end{equation}

\begin{figure}
\hfill
\begin{subfigure}{.4\textwidth}
  \centering
  \includegraphics[scale=0.4]{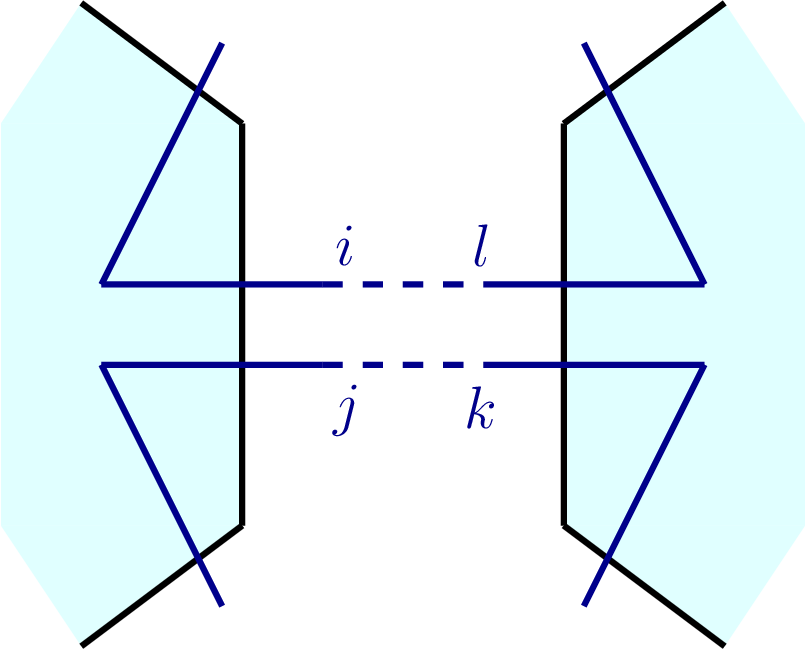}
  \caption{Wick Theorem identifies pairs of indices together.}
  \label{fig:mat_propa}
\end{subfigure}
\hfill
\begin{subfigure}{.4\textwidth}
  \centering
  \includegraphics[scale=0.5]{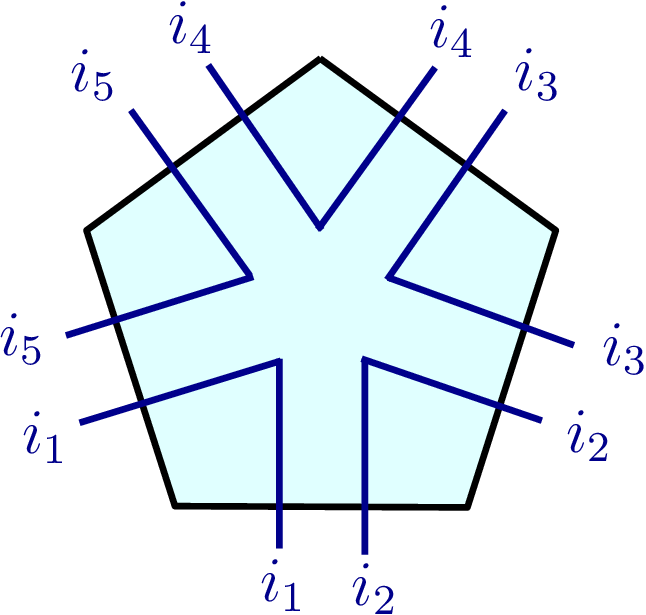}
  \caption{A vertex of the matrix model.}
  \label{fig:mat_vert}
\end{subfigure}
\hfill
\caption{The propagator and the vertex $\Tr H^5$ of the Hermitian matrix model.}
\label{fig:matrix_map}
\end{figure}

Each of the terms in this expansion can be represented diagrammatically via \emph{Feynman diagrams}. Each trace $\Tr H^k$ is represented as a vertex with $k$ pairs of half-edges and carries a weight $Nt_k$. The pair of half-edges correspond to copies of the matrix entries $H_{ij}$ and are labeled with the corresponding matrix indices $(i,j)$. Observe that the $k$ pairs are cyclically ordered due to the contraction pattern of the indices in the trace, and two consecutive indices in different half-edges are identical. This vertex can be represented as a $k$-gon such that each side carries a pair of half-edges. This is illustrated in Figure~\ref{fig:mat_vert}. We now apply Wick's Theorem~\ref{thm:Wick} to compute the expectation value $\left\langle\Tr H^{k_1}...\Tr H^{k_n}\right\rangle$. This expectation value is given by the sum of the possible pairings between the half-edges $(i,j)$ with a weight $\frac{t}{N}$ for each pair. Moreover, when two half-edges are paired their indices are identified according to~\eqref{eq:Wick} as represented in Figure~\ref{fig:mat_propa}. This identification of the indices amounts to gluing the corresponding sides of the polygons together, as in Definition~\ref{def:poly_glu}. Therefore, a term of the Wick expansion is nothing but a particular choice of gluing for the polygons of the map and the sum over pairings of the Wick Theorem in the monomial $\Tr H^{k_1}...\Tr H^{k_n}$ corresponds to a sum over gluings leading to a map with face profile $(k_1,...,k_n)$.

\medskip

We have exactly one free index $i\in\{1,...,n\}$ per face of the map, which yields a factor $N$ after the sum on that index is performed. Adding the contribution from the pairings and the vertices, we get that the scaling with $N$ of a map $\mathcal{M}$ is given by its Euler characteristics $\chi(\mathcal{M})$. Finally, the cardinality of the automorphism group of a map with $m_i$ faces of length $i$ such that $\sum_{i} m_i = n$ each carrying a label $\ell \in \llbracket 1,N \rrbracket$ is given by the automorphism factor leaving invariant the labels of each map, i.e. 

\begin{align}
    \label{eq:aut_grp}
    \vert Aut(\mathcal{M}) \vert &= \binom{n}{m_1\hspace{2pt}\dotsc\hspace{2pt}m_d}\prod\limits_{i} i^{m_i}i! \\
    &= n!\prod\limits_{i} i^{m_i} \nonumber.
\end{align}

It follows that the numerical prefactor that appears for a map $\mathcal{M}$ in the expansion in perturbative series is exactly $\vert Aut(\mathcal{M}) \vert^{-1} $. Therefore we have
\begin{align}
    \label{eq:matrix_maps}
    Z^{(\text{pert})}_{\mathcal{H}_N}\left[t,t_k\right] &= \sum\limits_{\mathcal{M}\in\mathfrak{M}} N^{\chi(\mathcal{M})}t^{E(\mathcal{M})}\prod\limits_{i=1}^d t_k^{V_i(\mathcal{M})},
\end{align}
where $\mathfrak{M}$ is the set of maps and $E(\mathcal{M})$ is the number of edges of $\mathcal{M}$ and $V_i(\mathcal{M})$ its number of vertices of degree $i$.

\medskip 

We can see that the perturbative expansion of the Hermitian matrix model corresponds to the genus series of (unrooted) maps with weights $t_k$ per vertex of degree $k$ and $t$ per edge. By duality, it is also the genus series of unrooted maps where the weights $t_k$ per face of degree $k$ and $t$ per edge. In the dual representation, each vertex represents a $k$-gon and the Wick Theorem corresponds to the sum over possible gluings as in Definition~\ref{def:poly_glu}.

\medskip 

This sum is over all maps, including non-connected ones. To restrict the sum to connected maps, we consider the logarithm of the partition function, called the \emph{free energy} $F$ in field theoretic terms. To see combinatorially why taking the log restricts the sum to connected maps, it is more convenient to inverse this relation: a map with $k$ connected components consists of $k$ connected maps up to a permutation of these $k$ components. This is exactly the $k$-th order of the power series of $\exp(F\left[t,t_k\right])$.

\begin{align}
    \label{eq:log_part}
    F\left[t,t_k\right] &= \log Z^{(\text{pert})}_{\mathcal{H}_N}\left[t,t_k\right] \\
    &= \sum\limits_{\substack{\mathcal{M}\in\mathfrak{M}\\\mathcal{M} \text{   connected}}} N^{\chi(\mathcal{M})}t^{E(\mathcal{M})}\prod\limits_{i=1}^d t_k^{V_i(\mathcal{M})}.\nonumber
\end{align}

For a connected map $\mathcal{M}$, its Euler characteristics is $\chi(\mathcal{M})=2-2g(\mathcal{M})$. Therefore the free energy admits a \emph{topological expansion} in $\frac{1}{N^2}$ where the $g$-th order of the expansion is given by the sum over maps of genus $g$.

\begin{align}
\label{eq:genus_exp}
       F\left[t,t_k\right] &= \sum\limits_{g \geq 0} N^{2-2g} \sum\limits_{\substack{\mathcal{M}\in\mathfrak{M}_g\\\mathcal{M} \text{   connected}}} t^{E(\mathcal{M})}\prod\limits_{i=1}^d t_k^{V_i(\mathcal{M})} \\
       &= \sum\limits_{g \geq 0} N^{2-2g} F_g\left[t,t_k\right], \nonumber
\end{align}
where $\mathfrak{M}_g$ denotes the set of maps of genus $g$. In particular, in the large $N$ limit, only the leading order in $N$ survives: the contribution of maps with non-zero genus is suppressed and we are left with a sum over planar maps. 

\medskip

For an operator $O$ of the matrix entries $M_{ij}$, its expectation value is defined asas 
\begin{align}
    \label{eq:exp_val_mat}
    \langle O \rangle = \frac{1}{Z^{(\text{pert})}_{\mathcal{H}_N}}\int_{H_N} O\exp\left(-\frac{N}{2t}\Tr H^2 + N\sum\limits_{k \geq 3} \frac{t_k}{k} \Tr H^k\right)dH,
\end{align}
taken as a perturbative expansion around the Gaussian free field. In the following, an important role is played by (products of) polynomial invariant of the matrix $H$ given by $\Tr H^k$ for $k>0$. We introduce the \emph{correlation functions} $\hat{W}_n(x_1,\dotsc,x_n)$ which are generating series for products of $n$ such polynomials as 
\begin{align}
    \hat{W}_n(x_1,\dotsc,x_n) &= \sum\limits_{k_1,\dotsc ,k_n \in \mathbb{N}} \langle \prod\limits_{i=1}^n x_i^{-1-k_i} \Tr H^{k_i} \rangle \\
    &= \langle \prod\limits_{i=1}^n \Tr \frac{1}{x_i-H} \rangle. \nonumber
\end{align}

For example, the $1$-point correlation function $\hat{W}_1(x_1)$ contains all the information on the expectation value of the eigenvalues of the matrix. It is related to the eigenvalue density $\rho(x) = \langle \frac{1}{N}\sum\limits_{i=1}^N \delta(x-\lambda_i) \rangle$ via a Stieljes transform
\begin{equation}
    \label{eq:W1_ev}
    \hat{W}_1(x) = \int_{\text{supp  } \rho} \frac{\rho(x')dx'}{x-x'}.
\end{equation}

Terms $\Tr H^k$ for some $k\geq 3$ are linearly coupled to some (appropriately named) coupling constant $t_k$. Thus any expectation value of the form $\langle \prod\limits_{i=1}^n \Tr H^{k_i} \rangle $ can be represented as the action of a differential operator on the partition function $Z^{(\text{pert})}_{\mathcal{H}_N}\left[t,t_k\right]$.
\begin{equation}
\label{eq:diff_op}
    \langle \prod\limits_{i=1}^n  \Tr H^{k_i} \rangle \equiv \frac{1}{Z^{(\text{pert})}_{\mathcal{H}_N}}\frac{k_1\dotsm k_n}{N^n}\partial_{t_{k_1}} \dotsm \partial_{t_{k_n}} Z^{(\text{pert})}_{\mathcal{H}_N}\left[t,t_k\right].
\end{equation}

Since $Z^{(\text{pert})}_{\mathcal{H}_N}\left[t,t_k\right]$ is the generating function for maps, acting with $k_i\partial_{t_{k_i}}$ on the generating series $Z^{(\text{pert})}_{\mathcal{H}_N}\left[t,t_k\right]$ can be interpreted as marking one of the $k_i$ edges in a face of length $k_i$. Therefore, these products admit an expansion over maps with $n$ distinct boundaries of respective length $k_1,..,k_n$ and each of the boundaries have exactly one marked edge. When treating these faces as boundaries of the surface, the scaling with $N$ of a map is still given by the Euler characteristic of the corresponding surface. The power of the variable $x_i$ simply keeps track of the length of the $i$-th boundary (up to a shift of $1$). Hence they also admit a topological expansion over maps with $n$ boundaries. However, we are so far working with potentially non-connected maps, thus the boundaries may be carried on different connected components of the map. This can be circumvented by the introduction of the connected correlation function $W_n(x_1,..,x_n)$. For $n=1$, the two coincide. For $n=2$ it is defined as
\begin{equation}
    \label{eq:cumul_n2}
    W_2(x_1,x_2) = \hat{W}_2(x_1,x_2)-W_1(x_1)W_1(x_2).
\end{equation}

And it is defined similarly from an inclusion-exclusion formula for $n \geq 3$.

\subsection{The loop equations}
\label{ssec:Vir_eq}

Provided that the integrand vanishes on the boundary of the integration domain, which we shall admit here, the matrix integral can be used to generate relations on its correlation functions via the Stokes Theorem. For example, we have for $k>0$ and $i,j \in \llbracket 1,N\rrbracket$:
\begin{equation}
    \label{eq:ex_loop}
    0 = \int dH\frac{\partial}{\partial H_{ij}} \left( (H^k)_{ij} \exp\left(-\frac{N}{2t}\Tr H^2 + N\sum\limits_{k \geq 3} \frac{t_k}{k} \Tr H^k\right) \right).
\end{equation}

Via matrix calculus, we have $\frac{\partial}{\partial H_{ij}} H^k _{ij} = \sum\limits_{\ell = 0}^{k-1} (H^\ell)_{ii}(H^{k-1-\ell})_{jj}$ and $\frac{\partial}{\partial H_{ij}} \Tr H^k = k(H^{k-1})_{ji}$. Setting $i=j$ and summing over $i$, Equation~\eqref{eq:ex_loop} leads to the following relation between expectation values
\begin{equation}
    \label{eq:ex_loop2}
    \langle \Tr H^{k+1} \rangle + N \sum\limits_{i \geq 3} t_i \langle \Tr H^{i+k-1} \rangle = \sum\limits_{\ell = 0}^{k-1} \langle \Tr H^\ell \Tr H^{k-1-\ell} \rangle
\end{equation}

This gives a first example of a \emph{loop equation}. Seeing the second term of~\eqref{eq:ex_loop2} as associated with the coupling constant $t_2 = \frac{-1}{t}$, this equation can be recast as a differential operator as 
\begin{equation}
    \label{eq:Vir_eq}
    \underbrace{\left( \frac{1}{N^2}  \sum\limits_{\ell = 0}^{k-1} \ell(k-1-\ell) \partial_\ell \partial_{k-1-\ell} + \sum\limits_{i \geq 2} t_i(i+k-1)\partial_{k+i-1} \right)}_{V_k} Z_{\mathcal{H}_N}\left[t,t_k\right].
\end{equation}

The differential operator $V_k$ is the $(k-2)$-th Virasoro constraint. A brief computation shows that these operators form a Lie algebra known as (half of) the \emph{Virasoro algebra}. Its commutator is given by
\begin{equation}
    \label{eq:bracket_Vir_alg}
    \left[V_i,V_j\right] = (i-j)V_{i+j},\qquad i,j \geq 1.
\end{equation}
These operators play a special role in the Hermitian matrix model as the set of Virasoro constraints $\{V_k\}_{k \geq 1} $ determines the partition function $Z_{\mathcal{H}_N}\left[t,t_k\right]$ uniquely. 

\medskip

Similar relations can be obtained by adding products of traces to the initial function $(H^k)_{ij}$. The relations thus obtained are called \emph{loop equations}.  For $p\geq 1$ and $\mu_i \in \mathbb{N}$, the function $(H^{\mu_1})_{ij} \prod\limits_{i=1}^p \Tr H^{\mu_i}$ leads to the following loop equation~\cite{EKR_Review}
 \begin{align}
     \label{eq:loop_eq}
     0 = &\sum\limits_{\ell = 0}^{\mu_1-1} \langle \Tr H^\ell \Tr H^{\mu_1-1-\ell}\prod\limits_{i=2}^p \Tr H^\mu_i\rangle + \sum\limits_{j=2}^p \langle \Tr H^{\mu_1+\mu_j -1 } \prod\limits_{\substack{2 \leq i \leq p \\ i \neq j}} \Tr H^\mu_i\rangle \\
     &- N \sum\limits_{k \geq 2} t_k \langle \Tr M^{\mu_1+k-1} \prod\limits_{i=2}^p \Tr H^\mu_i\rangle.  \nonumber
 \end{align}

These equations have a similar combinatorial interpretation to the Tutte equations of Section~\ref{ssec:map_gen_Tutte}. They encode what can happen when the marked edge of the first boundary is deleted. Either the edge was bordered by the same face on both sides and we end up with two new boundary faces of length $\ell$ and $\mu_1-1-\ell$, or it was bordered by another boundary face, resulting merging of the two boundaries in a single one of length $\mu_1+\mu_j-1$, or it was bordered by a face of length $k$ with weight $t_k$. Algebraically, these three cases correspond respectively to the first, second and third terms of Equation~\eqref{eq:loop_eq}.

\subsection{Eigenvalue decomposition}
\label{ssec:ev_decomp}

Any Hermitian matrix $H$ can be diagonalized by an unitary matrix $U\in U(N)$ and has real eigenvalues $\Lambda=\left(\lambda_1,..,\lambda_N\right)$

\begin{equation}
\label{eq:diag_herm}
    H=U\Lambda U^{\dagger}.
\end{equation}

Note that this decomposition is not unique. Multiplying $U$ to the right by a diagonal matrix of the form $\left(e^{i\theta_1},...,e^{i\theta_N}\right)$ or by a permutation matrix leaves the decomposition unchanged. Therefore the space of Hermitian matrices $\mathcal{H}_N$ decomposes as

\begin{equation}
    \mathcal{H}_N \simeq \left(U_N/\left(U_1\right)^N \times \mathbb{R}^N\right)/\mathfrak{S}_N. \nonumber
\end{equation}

Therefore the measure writes
\begin{align}
    \label{eq:meas_ev}
    dH &= dU\Lambda U^\dagger + Ud\Lambda U^\dagger + U\Lambda dU^\dagger \\
       &= \left[dUU^\dagger,\Lambda U^\dagger\right] + Ud\Lambda U^\dagger
\end{align}
where we used that $dU^\dagger = -U^\dagger dU U^\dagger$ since $U$ is unitary. The Jacobian of this change of variable is $\prod\limits_{i>j} (\lambda_i-\lambda_j)^2$. We recognize the square of the Vandermonde determinant associated with the matrix $\left(\lambda_i^{j-1})\right)_{1\leq i,j \leq N}$. As we will see later in section~\ref{sec:KP_hier} this is related to the integrability of random Hermitian matrices and its connection with the KP hierarchy. More generally, being able to reduce the matrix integral to its eigenvalues allows for powerful methods to solve the model e.g. orthogonal polynomials~\cite{EKR_Review}. The partition function~\eqref{eq:part_herm} can be written as an integral of the eigenvalues as
\begin{equation}
\label{eq:ev_pert_fct}
    Z^{(\text{pert})}_{\mathcal{H}_N}\left[t,t_k\right] = \int \prod\limits_{i=1}^n d\lambda_i \prod\limits_{i>j} (\lambda_i-\lambda_j)^2 \exp\left(\sum\limits_{i=1}^N -\frac{N}{2t}\lambda_i^2 + N\sum\limits_{k\geq 3} \frac{t_k}{k} \lambda_i^k \right)
\end{equation}
where we used the fact that for the Haar measure $\int_{U_N} dU = 1$. 

\subsection{The double scaling limit}
\label{ssec:double_scaling_limit}

In the previous section, we have seen that the partition function of the matrix model can be seen as the partition function of a toy model for $2$-dimensional quantum theory of gravity. The $\frac{1}{N}$-expansion organizes the generating function of connected maps by their genus as given by~\eqref{eq:genus_exp}. Each of the coefficients $F_g\left[\bf{t_k}\right]$ can then be expanded as a function of the coupling constant. For clarity of the presentation, we shall restrict our maps to quadrangulations, but the argument holds in full generality~\cite{Eynard_book}. At the level of the partition function, this is achieved by setting all coupling constant but $t_4$ to zero. Each coefficient of the free energy can be expanded as a power series of the coupling constant $t_4$.
\begin{equation}
\label{eq:genus_coupling_exp}
    F_g\left[t_4\right] = \sum\limits_{n \geq 0} F_{g,n}t_4^n
\end{equation}
The coefficients $F_{g,n}$ is the number of maps of genus $g$ made of exactly $n$ quadrangles. This series is convergent for each genus $g$ and its radius of convergence is independent of the genus~\cite{BC_ratio, AL_ratio} and denoted $t^{(c)}_4$. The average number of quadrangles in a map of genus $g$ can be computed as
\begin{equation}
    \label{eq:avg_face}
    \left\langle n \right\rangle_g = t_4\frac{\partial \ln F_g\left[t_4\right]}{\partial t_4}
\end{equation}
This number becomes infinite when $t_4 \rightarrow t^{(c)}_4$. Therefore, sending $t_4$ to its critical value can be interpreted as taking a continuum limit with infinitely many quadrangles. Around $ t^{(c)}_4$, the leading singular part of $F_g\left[t_4\right]$ behaves like $\left(t_4-t^{(c)}_4\right)^{2-\gamma}$ for some $\gamma \in \mathbb{R}$. Therefore when approaching $t^{(c)}_4$, the number of quadrangles behave like $\left(t_4-t^{(c)}_4\right)^{-\frac{1}{2}(2-\gamma)\chi}$ with $\chi = 2-2g $. If we impose each quadrangle to have area $\epsilon$, then the area is fixed when sending $t_4 \rightarrow t^{(c)}_4$, $\epsilon\rightarrow 0$ while holding their ratio fixed. Hence sending $t_4$ to its critical value can be thought of as taking a continuum limit, refining the mesh size of the discretized surface. The generating functions $F_g$ behave similarly for all genuses, they all have the same singular value $t^{(c)}_4$ and critical exponent $\gamma$. Therefore sending $N \rightarrow \infty $ and  $t_4 \rightarrow t^{(c)}_4$, the generating series expands at leading order in $(t_4-t^{(c)}_4)$ as
\begin{equation}
     F\left[t_4\right] = \sum\limits_{g \geq 0} N^{2-2g}(t-t_c)^{\frac{1}{2}(2-\gamma)\chi}\tilde{F_g}(t_4)
\end{equation}
for some coefficients $\tilde{F_g}(t_4)$. These coefficients are rational functions of the coupling constants~\cite{BC_ratio, AL_ratio} but they are not singular at $t_4^{(c)}$. If the limits are taken such that $\kappa^{-1} = N\left(t_4-t_4^{(c)}\right)^{\frac{1}{2}(2-\gamma)}$ is held fixed then the generating admits an expansion as
\begin{equation}
    \label{eq:ds_map_exp}
    \tilde{F}\left[t_4\right] = \sum\limits_{g \geq 0} \kappa^{2-2g}\tilde{F_g}(t_4).
\end{equation}
The generating series $\tilde{F_g}\left[t_4\right]$ \emph{double scaling limit} of $F$. Its correlation functions $\tilde{W}_n$ correspond to the double scaling limit of the correlation function $W_n$ of $F\left[t_4\right]$. 

\medskip

A similar mechanism holds when the generating functions have several coupling constants. In this case, there are several different loci where the generating function is singular and they may intersect, leading to a more involved phase diagram with several different critical manifolds and multi-critical points at their intersection~\cite{Eynard_book}. 

\medskip

The double scaling limit of the Hermitian matrix model is equivalent to topological gravity described as intersection theory on the moduli space of connected Riemann surfaces~\cite{DFGiZJ}. This equivalence was conjectured by Witten~\cite{W90_conj} from the fact that both series satisfy a reduction of the KP hierarchy called the Korteweg–de Vries (KdV) hierarchy~\footnote{To quote Kontsevich in~\cite{Kont92}: "The reason for this conjecture is an irrational idea (for mathematicians), that gravity is unique".}. The conjecture was subsequently proven by Kontsevich~\cite{Kont92} via a cellular decomposition of the moduli space of Riemann surfaces via \emph{Strebel graphs} and establishing that these graphs are the Feynman graphs of a matrix model now known as the Kontsevich model given by the partition function
\begin{equation}
    \label{eq:Konts_mod}
    Z_N\left[\Lambda\right] = \int_{H_N} \exp(i\Tr H^3 - \Tr X\Lambda)dH,
\end{equation}
where $\Lambda$ is a diagonal matrix with real eigenvalues $\left(\lambda_i\right)_{1 \leq i \leq N}$. Note that this model is of a different nature than the Hermitian matrix model since it involves an external matrix $\Lambda$. 

\medskip

This equivalence connects the double-scaling limit of the matrix model with finite representations of the conformal group. These representations have been classified by Kac~\cite{Kac90} and are characterized by two coprime integers $(p,q)$. The Kontsevich model corresponds to the case $(3,2)$. Other values of $(p,q)$ are related to different statistical physics models on a random lattice e.g. the Ising model for $(p,q)=(4,3)$ or 3-Potts for $(p,q)=(6,5)$. However, it is not known if any representation $(p,q)$ of the conformal group can be represented as the double scaling series of a matrix model.

\medskip

Parallelly, planar maps have been shown to converge to the Brownian map in the continuum limit~\cite{MM_Brownian_map,LG_Brownian_map}. The Brownian map is a continuous surface with spectral dimension $2$ but Hausdorff dimension $4$. From a field theoretic perspective, it corresponds to a mesoscopic limit and describes the properties of an "average" large random map. In a series of papers by Sheffield and Miller~\cite{MS_LQG1,MS_LQG2,MS_LQG3}, it has been shown to be related to a particular value of Liouville conformal field theory coupled to gravity. While of a different nature, these results strengthen the connection between the continuum limit of map generating series and the representation of the conformal group.

\section{The Kadomtsev–Petviashvili hierarchy}
\label{sec:KP_hier}

The Kadomtsev–Petviashvili (KP) equation is a non-linear partial differential equation originally introduced to describe the propagation of waves of small amplitudes in shallow mediums~\cite{KP_OG}. Sato~\cite{Sato82} showed that the original KP equation is actually the first equation of a countable set of non-linear partial differential equations. Each of these equations can be expressed as a Lax equation, and their set forms the KP hierarchy. In the same article, Sato established the connection between the KP equations and the infinite dimensional Grassmanian, which gives a geometric origin to the KP hierarchy and gives an effective framework to show that certain functions satisfy the KP hierarchy. This section gives a self-contained derivation of the KP hierarchy and shows how it applies to map enumeration and random matrices. It is based on the articles~\cite{Sato82,JM83,Kac90,AZ13}.

\subsection{The Grassmanian of finite dimensional vector space}
\label{ssec:finite_Grass}

Let $V$ be a vector space of dimension $n$ over $\mathbb{C}$ and fix a basis $(e_i)_{1 \geq i \geq n}$ of $V$. The Grassmanian $\text{Gr}(k,V)$ is the set of all $k$-dimensional vector subspaces of $V$. The group $GL(n,\mathbb{C})$ has a natural action on $V$ that describes changes of basis in $V$. In particular, it acts transitively on any $k$-dimensional subspace of $V$. Thus, the Grassmanian $\text{Gr}(k,V)$ can be described as $GL(n,\mathbb{C})$ up to transformations leaving any subspace of dimension $k$ invariant. This gives the Grassmanian the structure of a smooth manifold. To compute its dimension, consider $V_k$ a $k$-dimensional subspace of $V$. It has a basis $(v_i)_{1 \leq i \leq k}$ with $v_i = \sum\limits_{i=1}^n a_{ij}e_j$ where $(a_{ji})_{1 \leq i \leq n,1 \leq j \leq k}$ is a $n\times k$ matrix. Clearly, any element of $GL(k,\mathbb{C})$ acting on $(v_i)_{1 \leq i \leq k}$ leaves $V$ invariant. Therefore, by performing Gaussian elimination we can find an element of $GL(k,\mathbb{C})$ such that the matrix $(a_{ji})_{\substack{1 \leq i \leq n\\1 \leq j \leq k}}$ is
\begin{equation}
\label{eq:loc_coord_grass}
    \left[\begin{array}{cccc|ccc}
1&0&\cdots &0&a_{1,k+1} & \cdots & a_{1,n}\\
0&1&\cdots &0&a_{2,k+1} & \cdots & a_{2,n}\\
\vdots & &\ddots &\vdots\\
0&0&\cdots &1&a_{k,k+1} & \cdots & a_{k,n}\\
\end{array}\right].
\end{equation}
We are left with $n(k-n)$ coefficients, which shows that the Grassmanian $\text{Gr}(k,V)$ has dimension $k(n-k)$.

\subsubsection{The exterior algebra and the Plücker embedding}
\label{ssec:ext_alg}

The quotient algebra $\bigwedge(V)$ is the quotient of the tensor algebra $T(V)=\bigoplus\limits_{k \geq 0} V^{\otimes k}$ by the two-sided ideal $I=\{ x \otimes x \vert x \in V \}$. We denote the product on $\bigwedge(V)$, called the wedge product, by $\wedge$. By construction, it is alternating, (and therefore antisymmetric). The natural grading of $T(V)$ extends to $\bigwedge(V)$ as $\bigwedge(V)=\bigoplus\limits_{k \geq 0} \bigwedge^{k}(V)$ where 
\begin{equation}
    {\bigwedge}^{k}(V) = \text{span}\{v_1 \wedge .. \wedge v_k \vert (v_i)_{1 \leq i \leq k} \in V\}. \nonumber
\end{equation} 
The set ${\bigwedge}^{k}(V)$ is called the $k$-th exterior product of $V$ and its elements of the form $v_1 \wedge \dotsc \wedge v_k$ are called $k$-vectors. Using the antisymmetry of the wedge product, we can see that for any permutation $\sigma \in \mathfrak{S}_k$ and elements $v_1,\dotsc,v_k \in V$ we have $v_{\sigma(1)} \wedge \dotsc \wedge v_{\sigma(k)} = (-1)^{\vert\sigma\vert} v_1 \wedge \dotsc \wedge v_k$. Hence, for a basis $\left(e_i\right)_{1 \leq i \leq n}$ of $V$, the set
\begin{equation}
    \mathcal{B}^{(k)}= \{e_{i_1} \wedge .. \wedge e_{i_k} \vert 1 \leq i_1 < ... < i_n \leq n \}, \nonumber
\end{equation}
forms a basis of $\bigwedge^{k}(V)$. The coordinates of an element of $\bigwedge^{k}(V)$ in the basis $\mathcal{B}^{(k)}$ are the \emph{Plücker coordinates} (with respect to the basis $(e_i)_{1\leq i \leq n}$ of $V$).  Given $k$ vectors $v_i = \sum\limits_{j=1}^n a_{ij} e_j$ with $i \in \llbracket 1, k\rrbracket$ then the Plücker coordinates of $v_1 \wedge \dotsc \wedge v_k$ are given by the $k\times k$ minors of the matrix $A=(a_{ij})_{1 \leq j \leq n,1 \leq i \leq k}$. We denote $\Delta_{J}$ the $k\times k$ minor of the matrix $A$ with columns $J=\{j_1,\dotsc,j_k\}$ where the indices  $j_i \in \llbracket 1,n \rrbracket$ are two by two distinct. Consequently, $\bigwedge^{k}(V)$ has dimension $\binom{k}{n}$ and $\bigwedge^{k}(V) = \{0\}$ for $k>n$. 

\medskip

The following two Lemmas make the connection between the $k$-th exterior power $\bigwedge^{k}(V)$ and the Grassmanian $\text{Gr}(k,V)$.

\begin{lemma}
\label{lem:lin_dep}
    $(v_1,\dotsc,v_k) \in V^k$ are linearly dependent \hspace{3pt} $\Longleftrightarrow$ \hspace{3pt} $ v_1 \wedge \dotsc \wedge v_k = 0$.
\end{lemma}

\begin{proof}
    \begin{itemize}
        \item[$\Longrightarrow$] {If $(v_1,..v_k) \in V^k$ are linearly dependent then we have $v_k = \sum\limits_{i=1}^{k-1} c_i v_i$ with $c_i \in \mathbb{C}$. By linearity of the wedge product
        \begin{equation}
            v_1 \wedge \dotsc \wedge v_k = \sum\limits_{i=1}^{k-1} c_i v_1 \wedge \dotsc \wedge  v_{k-1} \wedge v_i
        \end{equation}
        and each term of the sum vanishes since the wedge product is alternating.}
        \item[$\Longleftarrow$] {If $v_1 \wedge \dotsc \wedge v_k = 0$, then all the $k\times k$ minors of the matrix $A$ are zero. Hence $A$ doesn't have maximal rank $k$, there exists a linear relation between rows of $A$, i.e. the elements $v_i$ are linearly dependent.}
    \end{itemize}
\end{proof}

\begin{lemma}
\label{lem:lin_span}
     If $(v_1,\dotsc,v_k)$ and $(v'_1,\dotsc,v'_k)$ span the same $k$-dimensional vector space then there exists $c$ such that $v_1 \wedge \dotsc \wedge v_k = c v'_1 \wedge \dotsc \wedge v'_k $.
\end{lemma}

\begin{proof}
    If $(v_1,\dotsc,v_k)$ and $(v'_1,\dotsc ,v'_k)$ span the same vector space, then there is a matrix $M = (m_{ij})_{1 \leq i,j \leq k} \in GL(k)$ such that $v'_i = \sum\limits_{j=1}^k m_{ij}v_j$. Denoting respectively $A$ and $B$ the matrix coordinates of $(v_i)_{1 \leq i\leq k}$ and $(v'_i)_{1 \leq i\leq k}$ in the basis $(e_i)_{1 \leq i\leq n}$ of $V$ we have $B = MA$. Therefore, the coordinate of $v_1 \wedge \dotsc \wedge v_k$ and $v'_1 \wedge \dotsc \wedge v'_k$ in the basis $\mathcal{B}^{(k)}$ are given by the $k \times k$ minors of $A$ and $MA$ respectively. Computing the minors of $MA$, it follows that their coordinates only differ by a global factor $\det M$.
\end{proof}

From these two lemmas, it follows that any $k$-vector $v_1 \wedge \dotsc \wedge v_k$ in $\bigwedge^{k}(V)$ which is not zero spans a unique $k$-dimensional vector space of $V$ with basis $(v_i)_{1 \leq i \leq k}$, and two different bases of a $k$-dimensional vector space map to the same element of $\bigwedge^{k}(V)$ up to a scalar factor. This gives a natural embedding  of the Grassmanian $\text{Gr}(k,V)$ into the projective exterior power $\mathbb{P}\bigwedge^k(V)$ called the \emph{Plücker embedding} given by
\begin{align}
    \label{eq:def_plucker}
    \iota \colon \text{Gr}(k,V) &\to \mathbb{P}{\bigwedge}^k(V) \\
    \text{span}\{v_1,\dotsc,v_k\} &\to \left[v_1 \wedge \dotsc \wedge v_k\right]. \nonumber
\end{align}
This mapping is injective. For a point $\left[v_1 \wedge \dotsc \wedge v_k\right] \in \text{Im}(\iota)$, we can recover the corresponding $k$-dimensional subspace $V_k$ via Lemma~\ref{lem:lin_dep} as $\{ v\in V \hspace{2pt}\vert \hspace{2pt} v \wedge v_1 \wedge \dotsc \wedge v_k = 0\}$ which determines $V_k$. The coordinates of $\text{span}\{v_1,\dotsc,v_k\}$ in $\mathbb{P}\bigwedge^k(V)$ are the Plücker coordinates of $\left[v_1 \wedge \dotsc \wedge v_k\right]$. 

\medskip

In general the Grassmanian $\text{Gr}(k,V)$ is not a projective space itself. The image of the Plücker embedding is exactly the $k$-vectors of $\mathbb{P}\bigwedge^k(V)$ while $\mathbb{P}\bigwedge^k(V)$ contains all linear combinations of $k$-vectors. The points of $\mathbb{P}\bigwedge^k(V)$ that correspond to points of the Grassmanian are characterized by the \emph{Plücker relations}.

\subsubsection{The Plücker relations}
\label{ssec:Plucker_relations}

The Plücker relations characterize the points of the Grassmanian $\text{Gr}(k,V)$ by a set of quadratic relations, thereby showing that the Grassmanian $\text{Gr}(k,V)$ is an algebraic variety. Any $k$-dimensional subspace of $V$ admits a basis $(v_i)_{1 \geq i \geq k}$ such that the matrix coordinates of $(v_i)_{1 \geq i \geq k}$ with respect to the basis $(e_j)_{1\leq j\leq n}$ takes the form~\eqref{eq:loc_coord_grass}. The parameters $(a_{ij})_{\substack{1 \leq i \leq k\\k+1 \leq j \leq n}}$ can be thought of as local coordinates on $\text{Gr}(k,V)$. This choice of basis for amounts the choice of a representative of $\left[v_1 \wedge \dotsc \wedge v_k\right]$ such that the Plücker coordinate $\Delta_{1,\dotsc,k}$ is $1$. The local coordinates $(a_{ij})$ are related to the Plücker coordinates of $\left[v_1 \wedge \dotsc \wedge v_k\right]$ via
\begin{equation}
    \label{eq:plucker_local}
    \Delta_{1,\dotsc \hat{\imath}\dotsc ,k,j} =
    \begin{vmatrix} 1 & \cdots & 0 & a_{1j} \\
    0 & \ddots & 0 & \vdots \\
    0 & \cdots & 1  & a_{kj}
    \end{vmatrix} = (-1)^{k+i} a_{ij}
\end{equation}
where $1 \leq i \leq k, k+1 \leq j \leq n$ and $\hat{\imath}$ means that we omitted the $i$-th column. This shows that we can identify the $k(n-k)$ Plücker coordinates corresponding to the above case with our choice of local coordinate on the Grassmanian (up to a sign factor). But since $\mathbb{P}\bigwedge^k(V)$ has dimension $\binom{n}{k}-1$, we have many more Plücker coordinates that can be computed as function of the $a_{ij}$, thus giving us relations between the Plücker coordinates themselves. Since they are minors of the matrix $A$, they can be computed recursively on the number of indices $r$ of omitted indices in $\llbracket 1,k\rrbracket$ by developing the minors with respect to lines corresponding to one of these omitted indices. For $r=2$, assuming $i_1<i_2$ and $j_1<j_2$ we get
\begin{equation}
\label{eq:plucker_local2}
\Delta_{1,\dotsc,\hat{\imath}_1,\hat{\imath}_2,\dotsc,k,j_1,j_2} = \begin{vmatrix}
    1 & 0 &\cdots & 0 & a_{1j_1} & a_{1j_2} \\
    0 & \ddots & \cdots &0 & \vdots & \vdots\\
    0 & \cdots & \cdots &0 & a_{i_1j_1} & a_{i_1j_2}\\
    0 & \cdots & \ddots &0 & \vdots & \vdots\\
    0 & \cdots & \cdots &0 & a_{i_2j_1} & a_{i_2j_2}\\
    0 & \cdots & 0 &1  & a_{kj_1} & a_{kj_2}
    \end{vmatrix} = (-1)^k \begin{vmatrix} a_{i_1j_1} & a_{i_1j_2}\\ a_{i_2j_1} & a_{i_2j_2}\end{vmatrix}.
\end{equation}
 This relation can be homogenized by restoring the dependency on $\Delta_{1,\dotsc ,k}$ on the left hand side. This gives us a first set of quadratic relations between the Plücker coordinates
 \small{
 \begin{equation}
     \label{eq:plucker_rel_delta}
     \Delta_{1,\dotsc ,k} \Delta_{1,\dotsc,\hat{\imath}_1,\hat{\imath}_2,\dotsc,k,j_1,j_2} = (-1)^k \left(\Delta_{1,\dotsc,\hat{\imath}_1,\dotsc,k,j_1}\Delta_{1,\dotsc,\hat{\imath}_2,\dotsc,k,j_2} - \Delta_{1,\dotsc,\hat{\imath}_1,\dotsc,k,j_2} \Delta_{1,\dotsc,\hat{\imath}_2,\dotsc,k,j_1}\right).
 \end{equation}}
 A similar computation still holds for $r\geq 2$. We can see that the $(r-1)\times (r-1)$ minors that appear when developing the minor is a Plücker coordinate itself (up to a sign factor). While we worked here on the open subset of $\mathbb{P}\bigwedge^k(V)$ where $\Delta_{1,\dotsc,k} \neq 0 $ and thus obtained relations involving $\Delta_{1,\dotsc,k}$, similar relations can be derived on any other open subset of $\mathbb{P}\bigwedge^k(V)$. The set of these quadratic relations are the \emph{Plücker relations}. In general, for two sequences of two-by-two distinct integers $(i_1,\dotsc,i_{k-1})$ and $(j_0,\dotsc,j_k)$ they write
 \begin{equation}
     \label{eq:plucker_relations}
     \sum\limits_{p=0}^{k} (-1)^p \Delta_{i_1,\dotsc,i_{k-1},j_p}\Delta_{j_0,\dotsc,\hat{\jmath}_p,\dotsc,j_k} = 0
 \end{equation}
 where coefficients vanish when the same integer appears twice.

\begin{remark}
     Note that when specifying an open subset of $\mathbb{P}\bigwedge^k(V)$ i.e. when working with local coordinates $(a_{ij})$ of the Grassmanian, the Plücker coordinates corresponding to minors with $r$ columns omitted in $\llbracket 1,k\rrbracket $ are not quadratic but of order $r$ as a function of the local coordinates $(a_{ij})$.
\end{remark}

\begin{remark}
    The Plücker relations can be given in a coordinate-free approach by working with the natural isomorphism $\bigwedge^k V \cong \bigwedge^{n-k} V^{*}$ and the bilinear map $\bigwedge^k V \times \bigwedge^{n-k} V \rightarrow \bigwedge^n V$ given by $(v,w) \mapsto v \wedge w$.
\end{remark}

The Plücker relations characterize points of the projective space $\mathbb{P}\bigwedge^k(V)$ which lie in the image of the Plücker embedding, that is they characterize the $k$-vectors of $\mathbb{P}\bigwedge^k(V)$. 

\subsection{The Sato Grassmanian and the semi-infinite wedge}
\label{ssec:Sat_and_wedge}

\subsubsection{The Sato Grassmanian}
\label{sssec:infinite_Grass}

For two vector spaces $V$ and $V'$, respectively of dimension $k+n$ and $k'+n'$ with basis $(e_i)_{1 \leq i \leq k+n}$ and $(e'_{i'})_{1 \leq j \leq k'+n'}$ such that $k\leq k'$,$\hspace{3pt}n\leq n'$, the Grassmanians $\text{Gr}(k,V)$ can be embedded in the Grassmanian $\text{Gr}(k',V)$. Before defining this embedding, we start by embedding the group $GL(k+n,V)$ in $GL(k'+n',V)$ via the following mapping
\begin{align}
    \label{eq:GL_map}
    GL(k+n,V) &\to \qquad GL(k'+n',V) \nonumber\\
    G \hspace{10pt}&\mapsto \left(
\begin{array}{c|c|c}
\text{id}_{k'-k}& 0 & 0 \\ \hline
0 & G & 0 \\ \hline
0 & 0 & \text{id}_{n'-n}
\end{array}\right)
\end{align}
Note that we could have chosen a different embedding in $GL(k'+n',V)$. In particular, it could seem more natural to work with $M$ only acting on the first $k+n$ vectors of the basis of $V'$. The reason behind this particular choice is that we are anticipating the following paragraph where both $k$ and $n$ will be sent to infinity to construct the infinite dimensional Grassmanian. 

\medskip

Now, we can describe an embedding of the Grassmanian that is compatible with the embedding of $GL(k+n,V)$. Using matrix coordinates $A = (a_{ij})_{\substack{1 \leq i \leq k \\ 1 \leq j \leq n+k}}$ to define the element $V_k$ of $\text{Gr}(k,V)$ with basis $\left(\sum\limits_{j=1}^{k+n} a_{ij} e_j\right)_{1 \leq i \leq k}$, we define it as
\begin{equation}
    \label{eq:map_grass}
    A 
    \mapsto
    \left(
\begin{array}{c|c|c}
\text{id}_{k'-k}& 0 & 0 \\ \hline
0 & A & 0 
\end{array}\right).
\end{equation}
This embedding is well-defined on the Grassmanian $\text{Gr}(k,V)$ since the action of $GL(n+k,V)$ embedded as an element of $GL(n',V)$ leaves the first $(k'-k)$ vector invariant. Using this embedding, the Grassmanian $\text{Gr}(k,V)$ can be embedded into any Grassmanian $\text{Gr}(k',V')$ as long as $k\leq k'$ and $n \leq n'$.

\medskip

The \emph{Sato Grassmanian} is obtained as the closure of the inductive limit for this embedding as $k$ and $n$ go to infinity. While we will not give the details of this construction here, let us sketch an intuitive picture of this object. The Sato Grassmanian is obtained by starting from the smallest possible Grassmanian $\text{Gr}(1,1) = \mathbb{P}^0 \cong \{1\}$, embedding it into all possible "slightly larger" Grassmanians and iterating the process countably many times. Suppose we start with Grassmanian $\text{Gr}(k,k+n)$ for some $k,n \in \mathbb{N}$. There are two "slightly larger" Grassmanians corresponding to a change in the parameter $n \rightarrow n+1$ or $k \rightarrow k+1$ with corresponding embeddings $\text{Gr}(k,k+n)\hookrightarrow\text{Gr}(k,k+n+1)$ and $\text{Gr}(k,k+n)\hookrightarrow \text{Gr}(k+1,k+1+n)$. The first case simply embeds the Grassmanian in a larger vector space. For an element of $\text{Gr}(k,k+n)$ associated with a matrix $\left[A\right]$, it adds a column of zero to the right of $\left[A\right]$. The second case corresponds to an embedding into a $k+1$-vector space while also increasing the dimension of the vector space, it extends the size of the identity matrix in $\left[A\right]$ by one. It is convenient to relabel indices of the matrix so that the first $k$ lines are labeled with negative indices. By doing so, we can think of $k$ and $n$ as independent parameters of the Grassmanian where increasing $k$ creates a new line and column with negative indices while increasing $n$ adds a new column with positive indices. Additionally, at each step of the process, new contributions from $\text{Gr}(k,k+n+1)$ or $\text{Gr}(k+1,k+1+n)$ are collected, that will in turn be embedded in larger Grassmanians. The Sato Grassmanian is obtained as the limit of this process as $k$ and $n$ go to infinity. 

\begin{definition}[Sato Grassmanian~\cite{Sato82}]
\label{def:Sato}
    The Sato Grassmanian $\mathcal{S}\text{Gr}$ is the set of matrices $A=(a_{ij})_{i \in \mathbb{Z}^-, j \in \mathbb{Z}}$ such that
    \begin{itemize}
        \item[$(i)$] $\exists m \in \mathbb{Z}^- \hspace{2pt}$ such that $a_{ij} = \delta_{ij}\hspace{2pt}$ for $i<m$.
        \item[$(ii)$] The $m$ line vector of $A$ for $i\in\llbracket -m,-1\rrbracket$ are linearly independent
    \end{itemize}
    quotiented by the group action induced by left-multiplication of matrices $(g_{ij})_{i,j \in \mathbb{Z}^-}$ satisfying $(i)$ and $(ii)$.
\end{definition}

\subsubsection{The semi-infinite wedge and Fermionic Fock space}
\label{sssec:sem_inf_wedge}

Now that we have constructed the Sato Grassmanian, we can embed it in an infinite dimensional projective space similarly as in the finite dimensional case. In order to do that, let us fix a basis $(e_i)_{i\in\mathbb{Z}}$ of a vector space of infinite dimension. An element of $\mathcal{S}\text{Gr}$ represented by matrix $A$ represents the infinite dimensional subspace spanned by vectors $\left(v_{i} = \sum\limits_{j\in\mathbb{Z}} a_{ij}e_j\right)_{i \in \mathbb{Z}^-}$. Following the final dimensional case, $\text{span}\{v_{i}\}_{i \in \mathbb{Z}^-}$ should be mapped to $\bigwedge\limits_{i=-\infty}^{-1} v_{-i}$. From condition $(i)$, there exists $m \in \mathbb{N}$ such that only the first $m$ elements $v_{-i}$ are actually different from $e_{-i}$. Therefore this wedge product can be written as $\bigwedge\limits_{i=-\infty}^{-m-1} e_{i} \wedge \bigwedge\limits_{i=-m}^{-1} v_{-i}$. Since the wedge product is alternating, when trying to expand the remaining $m$ vectors $v_{-i}$ over the basis $(e_{i})_{i\in\mathbb{Z}}$, only contribution from $e_k$ with $k\geq -m$ will be non-zero. This motivates the following definition of an infinite dimensional analog of the exterior power called the \emph{semi-infinite wedge}. 

\begin{definition}[Semi-infinite wedge]
\label{def:semi_inf_wedge}
    The semi-infinite wedge $\bigwedge^{\frac{\infty}{2}}_0$ is defined as
    \begin{equation}
        {\bigwedge}^{\frac{\infty}{2}}_0 \hspace{4pt}= \text{span}\{\bigwedge\limits_{i=-\infty}^{-m-1} e_{i} \wedge \bigwedge\limits_{k=1}^{m} e_{i_k} \vert m\in\mathbb{N}, i_k \in\llbracket-m+1,+\infty\llbracket \hspace{3pt}\text{s.t.}\hspace{3pt} i_1<\dotsc<i_m \}.
    \end{equation}
    Moreover, this set forms a basis of $\bigwedge^{\frac{\infty}{2}}_0$.
\end{definition}

The semi-infinite wedge is said to be \emph{frozen at $-\infty$} as its elements omit finitely many $v_k$ with negative index. Moreover, a basis element $b$ of $\bigwedge^{\frac{\infty}{2}}_0$ contains as many omitted elements $e_i$ with negative indices as elements $e_i$ carrying a positive label. This allows us to describe the basis of $\bigwedge^{\frac{\infty}{2}}_0$ given in Definition~\ref{def:semi_inf_wedge} compactly by labeling them with integer partitions. To do that, we introduce the \emph{Frobenius coordinates} of an integer partition. For a partition $\lambda$, we define its transposition $\lambda'$ as the partition with parts $\lambda'_i = \vert \{ j \vert \lambda_j \geq i \} \vert$. Graphically, $\lambda'$ is obtained by taking the mirror image of the Ferrer diagram of $\lambda$ along the diagonal axis $x=y$.

\medskip

\begin{figure}
    \centering
    \includegraphics[scale=0.5]{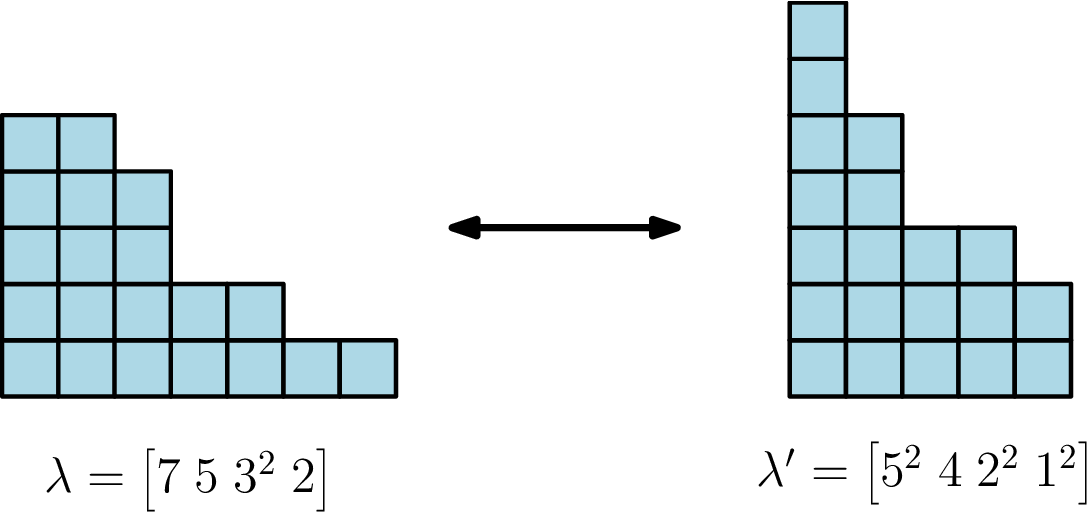}
    \caption{A partition $\lambda$ and its transpose $\lambda'$.}
    \label{fig:part_transpose}
\end{figure}

For an integer partition $\lambda =(\lambda_i)$, we define $d(\lambda)=\max\{i \vert \lambda_i-i \geq 0\}$. Observe that $d(\lambda)=d(\lambda')$.

\begin{definition}[Frobenius coordinates]
     For $i\in\llbracket1,l\rrbracket$, we define $\alpha_i = \lambda_i -i$ and $\beta_i = \lambda'_i-i$. They are a decreasing sequence of integers of length $d(\lambda)$. The two vectors $\vec{\alpha}=(\alpha_i)$ and $\vec{\beta}=(\beta_i)$ are the Frobenius coordinates of the partition.
\end{definition} 

Since $d(\lambda) = d(\lambda')$, the sequence $\vec{\alpha}$ and $\vec{\beta}$ have the same length. These two vectors completely characterize $\lambda$. Graphically they encode how many blocks $\lambda$ contains to the right and above the diagonal axis respectively, while their length gives the number of blocks on that diagonal. We write $\lambda = (\vec{\beta}\vert \vec{\alpha})$.

\medskip

We associate the partition $\lambda = (\vec{\beta}\vert \vec{\alpha})$ to the basis element $\vert \lambda \rangle$ with omitted negative indices $e_{-\beta_i-1}$ and including positives indices $e_{\alpha_i}$

\begin{equation}
\label{eq:basis_wedge_part}
    \vert \lambda \rangle = \dotsc \wedge \hat{e}_{-\beta_1-1} \wedge \dotsc \wedge \hat{e}_{-\beta_{d(\lambda)-1}} \wedge \dotsc \wedge e_{\alpha_{d(\lambda)}} \wedge \dotsc \wedge e_{\alpha_1}.
\end{equation}
In particular, we have
\begin{equation}
\label{eq:def_vac}
    \vert \emptyset \rangle = \bigwedge\limits_{i=-\infty}^{-1} e_{i}.
\end{equation}

Using this labeling for this basis of $\bigwedge^{\frac{\infty}{2}}_0$, the basis element $\vert \emptyset \rangle$ can be thought of as a reference point in the infinite dimensional Grassmanian. The other basis elements are characterized by how they can be reached from this element. To make this idea more precise, we need to define operators that allow us to add or remove wedge with basis vector $e_i$ in an infinite product. However, the semi-infinite wedge is not stable under these operations. To circumvent this issue, we first introduce a larger space, the \emph{Fermionic Fock space}, by relaxing the balance constraints between positive and negative labels. 

\begin{definition}[Fermionic Fock space]
    The Fermionic Fock space $\bigwedge^{\frac{\infty}{2}}$ is defined as the vector space
    \begin{equation}
        \text{span}\{\bigwedge\limits_{k=1}^{+\infty} e_{i_k}  \vert (i_k)_{k\in\mathbb{N}} \text{\hspace{5pt}strictly decreasing \hspace{5pt}s.t.\hspace{5pt}}\exists m\in\mathbb{N} \hspace{4pt}\forall m'>m \hspace{4pt}i_{m'+1} = i_{m'}-1 \}.
    \end{equation}
    Moreover, this set forms a basis of $\bigwedge^{\frac{\infty}{2}}$.
\end{definition}

The Fermionic Fock space contains all infinite wedge products of $e_i$ frozen at $-\infty$. Therefore we can now define operators adding and removing wedge products on elements of that space. For any $i \in \mathbb{Z}$, we define $\psi_i$ as
\begin{align}
\label{eq:def_crea}
    {\bigwedge}^{\frac{\infty}{2}} &\to {\bigwedge}^{\frac{\infty}{2}} \nonumber\\
    \psi_i\colon \hspace{12pt}v \hspace{8pt}&\mapsto e_i \wedge v.
\end{align}
The vector space $\text{span}\{e_i\}_{i \in \mathbb{Z}}$ admits a dual whose basis are linear forms $(e^{*}_i)_{i \in \mathbb{Z}}$ such that $e^{*}_i(e_j) = \delta_{ij}$. This is equivalent equips the vector space with a scalar product, and we will denote $e^{*}_i(e_j)= \langle e_i \vert e_j \rangle $. This definition extends to the Fermionic Fock space by linearity of the wedge product. For two basis elements $b_i$ and $b_j$ of ${\bigwedge}^{\frac{\infty}{2}}$ corresponding to sequences $i=(i_k)$ and $j=(j_k)$ we have $\langle b_i \vert b_j \rangle = b^{*}_i(b_j) = \prod\limits_{k \geq 1} \delta_{i_kj_k}$. For $i \in \mathbb{Z}$, we define $\psi^{*}_i$ the adjoint of the operator $\psi_i$ by the relation
\begin{equation}
    \label{eq:def_ann}
    \langle v' \vert \psi^{*}_i(v) \rangle = \langle \psi_i(v') \vert v \rangle \qquad \forall v,v' \in {\bigwedge}^{\frac{\infty}{2}}.
\end{equation}
Explicitly, for $v \in {\bigwedge}^{\frac{\infty}{2}}$, $\psi^*_i(v)$ corresponds to the action of $e_i^{*}$ on $v$ induced by linearity over the wedge product. These operators satisfy fermionic anticommutation relations
\begin{equation}
    \label{eq:anticomm_rel}
    \{ \psi_i, \psi_j \} = 0,\qquad \{ \psi^*_i, \psi^*_j \} = 0, \qquad \{ \psi_i, \psi^*_j \} = \delta_{ij}
\end{equation}
Hence the algebra generated by operators $\{\psi_i,\psi_i^*\vert i \in \mathbb{Z}\}$ is a Clifford algebra.

\medskip

\begin{figure}[!ht]
    \centering
    \includegraphics[scale=0.85]{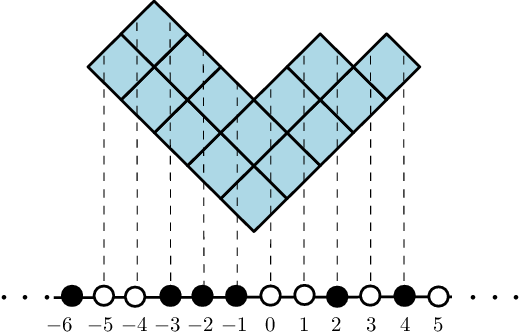}
    \caption{The Maya diagram associated with the basis element with partition $\lambda =  \left[5,4,2,2,2\right] $}
    \label{fig:wedge_state}
\end{figure}

Now, the basis element $\vert \lambda \rangle $ with $\lambda = (\vec{\alpha}\vert \vec{\beta})$ can be described via its relation to $\vert \emptyset \rangle$ as 
\begin{equation}
    \vert \lambda \rangle = \prod_{i=1}^{d(\lambda)} \psi^*_{-\alpha_i+1}\psi_{\beta_i} \vert \emptyset \rangle
\end{equation}
These operators can be used to describe basis elements of the Fermionic Fock space as states of a system of infinitely many fermionic particles, or equivalently in terms of \emph{Maya diagram}. The operators $\psi^*_i$ are interpreted as the deletion of a fermionic particle in position $i$, while $\psi_i$ corresponds to its creation. The element $\vert \emptyset \rangle$ defined in~\eqref{eq:def_vac} can thus be seen as a vacuum state - a Dirac sea - where all the negative modes are filled and the states $ \vert \lambda \rangle$ correspond to excited states. The correspondence between a state with partition $\lambda$ and Maya diagrams is illustrated in Figure~\ref{fig:wedge_state}.

The annihilation operators (with respect to the vacuum $\vert \emptyset \rangle$) are defined as the annihilator of $\vert \emptyset \rangle$. They are given by $\psi_i$ with $i<0$ and $\psi^*_i$ with $i \geq0$. The creation operators are its complement i.e. they correspond to $\psi_i$ with $i\geq 0$ and $\psi^*_i$ with $i <0$. 

\begin{definition}[Normal ordering]
\label{def:norm_ord}
    The normal order $\normord{O}$ (with respect to the vacuum $\vert \emptyset \rangle$) of a linear combination of finite product of operators $O$ is obtained by pushing all annihilation operators to the right while changing sign whenever two operators are exchanged in each monomial of $O$.
\end{definition} 

Annihilation/creation operators and normal ordering with respect to any other basis state $\vert n,\lambda \rangle$ can be defined similarly with respect to their own creation and annihilation operators, but this will not be needed in this manuscript.

\medskip

Observe that the set of creation operators of the vacuum $\vert \emptyset \rangle$ corresponds to annihilation operators of the dual state $\langle \emptyset \vert$ and vice-versa, therefore the expectation value of a normally ordered operator on the vacuum is always well-defined. This allows us to consider infinite sums of monomials in the fermions that would otherwise be ill-defined. This can be illustrated via the charge operator
\begin{equation}
\label{eq:def_charge}
    Q = \sum\limits_{k\in \mathbb{Z}} \normord{\psi_k\psi_k^*}.
\end{equation}
It counts particles in negative position with weight $-1$ and absence of particles in positive position with weight $1$. Then clearly we have $\langle Q \rangle = 0$ while any monomial $\psi_k\psi_k^*$ with $k<0$ would contribute as $1$ without the normal ordering, leading to a divergent expression. 

\medskip

For a basis element $e_{(i_k)}=\bigwedge\limits_{k=1}^{+\infty} e_{i_k}$\hspace{2pt}, there is a unique integer $n$ such that $e_{(i_k)}$ contains as many basis element with label $i_k \geq n$ as it miss elements $e_i$ with label $i<n$. This integer $n$ is called the charge and it is the eigenvalue of $e_{(i_k)}$ for the operator $Q$. In particular, the semi-infinite wedge corresponds to the subspace spanned by elements with zero charge. The shifted vacua of charge $n$, denoted $\vert n \rangle$ are defined as
\begin{equation}
    \vert n \rangle = \vert n, \emptyset \rangle = \begin{cases} \psi_{n-1}\dotsc\psi_0 \vert \emptyset \rangle \hspace{5pt}\text{if}\hspace{5pt} n>0 \\ \psi^*_{n}\dotsc\psi^*_1\hspace{5pt}\text{if}\hspace{5pt} n<0\end{cases}
\end{equation}

The vacuums are related by a simple relation, hence we can see that all subspaces of charge $n$ are isomorphic. As a consequence, basis elements with charge $n$ can be labeled via a partition describing their relation to their shifted vacua, similarly as in the zero charge case. Therefore the set
\begin{equation}
    \label{eq:basis_FFs}
   \{ \vert n,\lambda \rangle = \prod_{i=1}^{d(\lambda)} \psi^*_{-\beta_i-1+n}\psi_{\alpha_i+n} \vert n \rangle \vert n\in \mathbb{Z},\hspace{5pt} \lambda \hspace{5pt} \text{integer partition} \}
\end{equation}
forms a basis of the Fermionic Fock space. The dual basis $\langle n, \lambda \vert$ is defined similarly by exchanging the labels of $\psi^*$ and $\psi$ and acting to the right of $\langle n\vert$. This equips the Fermionic Fock space with a scalar product such that $\langle n, \lambda \vert$ form an orthonormal basis
\begin{equation}
    \label{eq:scal_FFS}
    \langle m, \nu \vert n, \lambda \rangle = \delta_{nm} \delta_{\lambda\nu}.
\end{equation}

For $O$ a linear combination of finite product of operators $\psi$ and $\psi^*$, its expectation value (with respect to the vacuum $\vert \emptyset \rangle$) as 
\begin{equation}
    \label{eq:exp_val}
    \langle O \rangle = \langle \emptyset \vert O \vert \emptyset \rangle
\end{equation}

\subsection{The Kadomtsev–Petviashvili equations}
\label{ssec:KP_equ}

In the previous paragraph, we have constructed the semi-infinite wedge (and the Fermionic Fock space) in which the Sato Grassmanian can be embedded. Instead of directly manipulating the Sato Grassmanian, it is more convenient to give an equivalent description of its embedding in the semi-infinite wedge (or the Fermionic Fock space) as the fermionic operators $\psi$ and $\psi^*$ give a very efficient framework on which to perform computations. As in the finite dimensional case, the Fermionic Fock space is a much larger space than the Sato Grassmanian and the embedding of the Sato Grassmanian in the semi-infinite wedge is characterized by an infinite set of Plücker relations which are related to the Kadomtsev–Petviashvili (KP) equations. We will show two different characterizations of the solution to the KP equations. One obtained directly in the semi-infinite wedge by the infinite dimensional generalization of the Plücker relations. The second will use the invariance of the Sato Grassmanian under the action of the infinite dimensional Lie algebra $GL(\infty)$.

\subsubsection{Plücker embedding of the Sato Grassmanian}
\label{sssec:Pluck_inf_case}

In this paragraph, we show briefly how the Plücker relations generalize in the semi-infinite wedge. Since the space is infinite dimensional, one has to be careful when manipulating expressions such as determinants or matrices to make sure they are well-defined. Here, we will only sketch how to write the Plücker relations in the semi-infinite wedge and refer to~\cite{Segal:1985aga} for technical details. 

\medskip

Recall that the Sato Grassmanian was obtained as the inductive limit of the embedding of finite dimensional Grassmanian $Gr(k,k+n)$. From the embedding~\eqref{eq:map_grass}, we can see that the embedding of $Gr(k,k+n)$ in the sato Grassmanian corresponds to infinite matrices $A=(a_{ij})_{i \in \mathbb{Z}^-, j \in \mathbb{Z}}$ with $a_{ij} = \delta_{ij} $ for $i\leq -k-1$ and $a_{ij} = 0$ for $j>n$, up to the action of $GL(k)$ on the last $k$ lines of $A$. Consider such an element $\left[A\right] = \text{span}\{v_{-i} = \sum\limits_{j=-k}^{n-1} a_{-ij}e_j \vert i \in \llbracket1,k\rrbracket\}$. Then all minors of $\left[A\right]$ but finitely many are trivial. Indeed, there is only one non-zero coefficient in columns with $j\leq -k+1$, and all columns with $j>n$ are identically zero. Therefore for an infinite minor $\Delta_{(i_k)}$ indexed by a decreasing sequence $(i_k)_{k\in \mathbb{N}}$ such that $p+1$ is the first index with $i_{p+1}\leq -k+1$ we have $\Delta_{(i_k)} = \Delta_{i_1,\dotsc,i_p}$. Moreover $\Delta_{i_1,\dotsc,i_p} = 0$ if $i_1 > n$. Therefore, the Plücker coordinates of $\left[A\right]$ in the semi-infinite wedge are totally determined from the Plücker coordinates $\Delta_{\lambda}$ associated with partitions $\lambda = (\vec{\alpha}\vert \vec{\beta})$ with $\beta_1+1\leq k$ and $\alpha_1 \leq n$. Observe that going back to the finite dimensional case, these partitions can be used to label the Plücker coordinates of $Gr(k,k+n)$ in the exterior algebra ${\bigwedge}^{k}$ (see Section~\ref{ssec:Plucker_relations}). Therefore, the set of Plücker relations corresponding to these partitions is nothing but the Plücker relations satisfied by the Grassmanian $Gr(k,k+n)$. When labeled with partition $\lambda = (\vec{\alpha}\vert \vec{\beta})$, they write~\cite{AZ13}
\begin{equation}
    \label{eq:plucker_part}
    \Delta_{(\vec{\alpha}\vert \vec{\beta})} \Delta_{(\vec{\alpha}_{\hat{\imath}\hat{\jmath}}\vert \vec{\beta}_{\hat{\imath}\hat{\jmath}})} =  \Delta_{(\vec{\alpha}_{\hat{\imath}}\vert \vec{\beta}_{\hat{\imath}})}\Delta_{(\vec{\alpha}_{\hat{\jmath}}\vert \vec{\beta}_{\hat{\jmath}})} -\Delta_{(\vec{\alpha}_{\hat{\imath}}\vert \vec{\beta}_{\hat{\jmath}})}\Delta_{(\vec{\alpha}_{\hat{\jmath}}\vert \vec{\beta}_{\hat{\imath}})}
\end{equation}
for all $i,j \in \llbracket 1,d(\lambda)\rrbracket$, where $(\vec{\alpha}_{\hat{\imath}}\vert \vec{\beta}_{\hat{\jmath}})$ corresponds to the partition $\lambda$ where the $i$-th element of $\vec{\alpha}$ and $j$-th element of $\vec{\beta}$ have been removed (and similarly removing two hooks when two indices are omitted). 

\medskip

By considering increasingly large Grassmanians, we get more Plücker relations corresponding to partitions with larger values of $\alpha_1$ and $\beta_1$. An element of the semi-infinite wedge $\sum\limits_{\lambda} \Delta_\lambda \vert \lambda \rangle$ is in the image of the Sato Grassmanian if and only if $(\Delta_\lambda)_{\lambda}$ satisfies Plücker relations~\eqref{eq:plucker_part} for all partitions $\lambda$.

\subsubsection{The \texorpdfstring{$\mathfrak{gl}(\infty)$}{gl(inf)} algebra}
\label{sssec:gl_inf}

In the finite dimensional case, the set of points of a Grassmanian $\text{Gr}(k,k+n)$ can be generated as the orbit of the transitive action of $GL(k+n)$ on any point of $\text{Gr}(k,k+n)$. We give here a similar construction for the Fermionic Fock space. The bilinear combinations $\psi_i \psi_j^*$ with $i,j \in \mathbb{Z}$ preserve the charge and satisfy commutation relations
\begin{equation}
    \label{eq:comm_bilin}
    \left[\psi_i \psi_j^*,\psi_k \psi_l^*\right] = \delta_{jk}\psi_i \psi_l^* - \delta_{il}\psi_j \psi_k^*,
\end{equation}
which is the same commutation relation as matrix elements $E_{ij}$. Provided their action on the semi-infinite wedge is well-defined, they can be seen as an endomorphism of the semi-infinite wedge. The simplest way to ensure that the action of an element of $\text{span}\{\psi_i \psi_j^*\}$ is to restrict to finite linear combination in the bilinear. This algebra $\mathfrak{gl}(\infty)$ is defined as

\begin{equation}
    \label{eq:def_gl_inf}
    \mathfrak{gl}(\infty) = \{ \sum\limits_{i,j\in \mathbb{Z}} a_{ij} \psi_i \psi_j^* \vert \hspace{2pt}\text{finitely many } a_{ij} \text{  are non-zero  } \},
\end{equation}

and the corresponding group is obtained by exponentiating these elements i.e. as
\begin{equation}
    \label{def:GL_grp}
    GL(\infty) = \{ \exp(g) \vert \hspace{2pt} g \in \mathfrak{gl}(\infty)\}.
\end{equation}
Elements $G\in GL(\infty)$ are exactly the elements such that their adjoint action on fermionic operators $\psi_i$ and $\psi^*_i$ leave the linear space spanned by $\psi_i$ and $\psi^*_i$ unchanged. For these elements, there is a matrix $R$ such that for all $i\in\mathbb{Z}$
\begin{equation}
\label{eq:bilinear_cond_proto}
G\psi_i^*G^{-1} = \sum\limits_{j \in \mathbb{Z}}  \psi^*_j R_{ji}, \qquad G\psi_iG^{-1} = \sum\limits_{j\in \mathbb{Z}} R'_{ij} \psi_j.
\end{equation}
Since they are invertible, we have $R^{-1} = R'$. These elements correspond to points that lie in the image of the Sato Grassmanian. They exactly characterize polynomial solutions to the KP hierarchy~\cite{Kac90}. However, in most applications in theoretical physics or combinatorics, the relevant objects are not polynomials but formal power series, thus their generating series do not fall in this class of solutions to the KP hierarchy. Therefore, we would like to have a similar construction allowing for countable sums in the bilinears $\psi_i \psi_j^*$. This can be achieved thanks to the the normal ordering (see Def. ~\ref{def:norm_ord}) which ensures that the action of infinite sum in the bilinear on Fermionic Fock space is well-defined. The algebra thus obtained is~\cite{JM83}

\begin{equation}
   \label{eq:def_gl_hat}
    \hat{\mathfrak{gl}}(\infty) = \{ \sum\limits_{i,j\in \mathbb{Z}} a_{ij} \normord{\psi_i \psi_j^*} \vert \hspace{2pt}\exists k \in \mathbb{N} \text{  s.t.  } a_{ij}=0 \text{  for  } \vert i - j \vert > k \}\oplus \mathbb{C}\cdot 1,
\end{equation}
and the corresponding group $\hat{GL}(\infty)$ is obtained by exponentiation of elements of the algebra. The definition of $\hat{\mathfrak{gl}}(\infty)$ allows for infinitely many non-zero coefficients $a_{ij}$ via the introduction of the normal ordering. This leads to the existence of elements $G'$ that are non-invertible but still satisfy similar commutation relations as~\eqref{eq:bilinear_cond_proto}. For those elements, only one of the two matrices between $R$ and $R'$ exists, and these relations become (in the case where $R$ is well-defined)
\begin{equation}
\label{eq:bilinear_cond_noninv}
    G'\psi^*_i = \sum\limits_{j\in \mathbb{Z}} R_{ji}\psi^*_iG', \qquad \psi_iG' = \sum\limits_{j\in \mathbb{Z}} R_{ij}G'\psi_j.
\end{equation}
These elements satisfy the following commutation relation.
\begin{equation}
    \label{eq:BBC}
    \left[G\otimes G, \sum\limits_{i\in\mathbb{Z}} \psi_i \otimes \psi^*_i\right] = 0.
\end{equation}
All these elements correspond to group-like~\cite{AZ13} elements associated with an element of $\hat{\mathfrak{gl}}(\infty)$ and are such that their action stabilizes the image of the Sato Grassmanian. Actually, Equation~\eqref{eq:BBC} is sufficient to ensure that an element $G$ can be used to construct a solution to the KP hierarchy, which is what we will exploit in the rest of this section. This allows us to go one step beyond the presentation sketched here to consider elements that satisfy the commutation relation~\eqref{eq:BBC} but that are \textit{not} in $\hat{GL}(\infty)$. In particular, for these elements, neither $R$ nor $R'$ matrix exists. An example of such an element is given by the projector operators
\begin{equation}
\label{eq:def_proj}
    P^+ = \prod_{i < 0} \psi_i \psi^*_i, \qquad P^- = \prod_{i \geq 0} \psi^*_i \psi_i,
\end{equation}
which are projectors on the space of positive or negative charge respectively. These operators will be relevant for the application of the KP hierarchy to matrix models (see section~\ref{sssec:matrix_KP}). This class of solution contains all the group-like elements of $\hat{GL}(\infty)$ and is still stable under multiplication~\cite{AZ13}. Finally, we shall admit that all elements $G$ satisfying the condition~\eqref{eq:BBC} have a definite charge $n$, i.e. they satisfy $\left[Q,G\right]=nG$.

\subsubsection{The boson-fermion correspondence}
\label{sssec:current_alg}

To construct solutions to the KP hierarchy, we finally have to make the connection between the previous constructions and generating series. This is achieved through the \emph{boson-fermion} correspondence. 

\medskip

Among the elements of $\hat{\mathfrak{gl}}(\infty)$, a key role is played by the current operators $J_k$. They are defined for $k\neq 0$ as 
\begin{equation}
    \label{eq:def_current}
    J_k = \sum\limits_{i\in \mathbb{Z}} \normord{\psi_i \psi^*_{i+k}}.
\end{equation}

As the current are bilinear in the fermions, they are bosonic operators. Their commutator is
\begin{equation}
     \label{eq:comm_current}
     \left[ J_k,J_l \right] = k\delta_{k,-l},
\end{equation}
therefore they form an infinite dimensional Heisenberg algebra.

\medskip

\begin{figure}[!ht]
    \centering
    \includegraphics[scale=0.75]{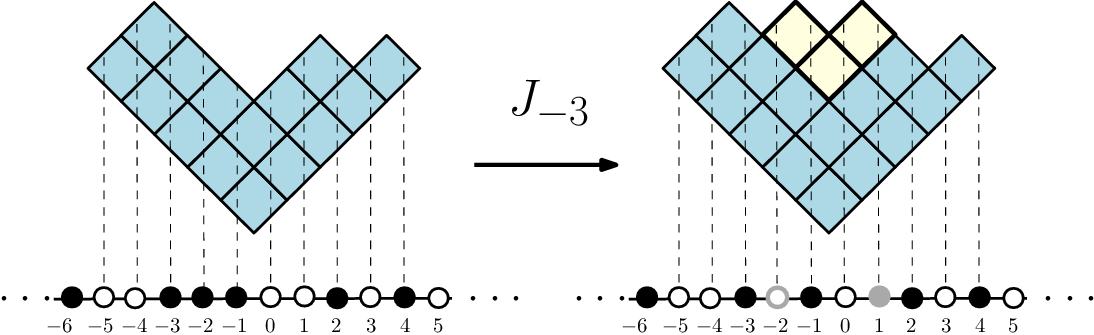}
    \caption{One of the states obtained by acting with the current operator $J_{-3}$ on the state of partition $\lambda = \left[5,4,2,2,2\right]$.}
    \label{fig:wedge_current}
\end{figure}

The action of a current $J_k$ on a state $\vert \lambda \rangle$ has a simple interpretation. On the fermions, everywhere possible it moves one particle of the state $\vert \lambda \rangle$ by $-k$ slots to the right while taking a sign factor for each fermion between the two positions. Combinatorially, if $k<0$ then it adds a strip of length $k$ wherever possible in the partition $\lambda$ and conversely if $k>0$ it deletes a strip of length $k$. This action is illustrated in Figure~\ref{fig:wedge_current}.

\medskip

In both cases, the sign factor is given by the height of the transposed strip. As we have seen in Section~\ref{sssec:MN_rule}, the effect of adding a strip in a partition is described by the Murnaghan-Nakayama rule. Tracking the sign carefully with respect to our definition of $\vert \lambda \rangle$ leads to a sign factor $(-1)^{\epsilon(\lambda)}$ with $ \epsilon(\lambda) = \sum\limits_{i=1}^{d(\lambda)} \beta_i +1$. Therefore we have
    \begin{equation}
        \label{eq:MN_exp}
        J_{-\mu_1}\dotsc J_{-\mu_k} \vert \emptyset \rangle = \sum\limits_{\lambda} (-1)^{\epsilon(\lambda)}\chi^\lambda(\mu) \vert \lambda \rangle,
    \end{equation}
and for any $k>0$  we have $J_{k} \vert \emptyset \rangle = 0$.

\medskip

We introduce times $(t_k)_{k \in \mathbb{Z}}$ which can be seen as weights associated with operators $J_k$ and the generating functions
\begin{equation}
    \label{eq:gen_current}
    J_+(\textbf{t}_+) = \sum\limits_{k\geq 1 } t_k J_k, \qquad J_-(\textbf{t}_-) = \sum\limits_{k\leq 1 } t_k J_k.
\end{equation}

From relation~\eqref{eq:MN_exp} and the definition of Schur functions~\eqref{eq:Schur_fct}, it follows that
\begin{equation}
    \label{eq:Schur_exp}
    \exp(J_-(\textbf{t}_-))\vert \emptyset \rangle = \sum\limits_\lambda (-1)^{\epsilon(\lambda)} s_\lambda(\textbf{t})\vert \lambda \rangle,
\end{equation}
and similarly for $\exp(J_+(\textbf{t}_+))$ acting to the right of the dual state $\langle \emptyset \vert$.

\begin{remark}
    The action of $\exp(J_-(\textbf{t}_-))$ on a state $\vert \lambda \rangle$ admits a similar expansion over the \emph{skew-Schur} functions $s_{\lambda\backslash \mu}$~\cite{AZ13}.
\end{remark}

Note that since $\exp(J_-(\textbf{t}_-))$ satisfies the condition~\eqref{eq:BBC} as an element of $\hat{GL}(\infty)$, it follows that $\exp(J_-(\textbf{t}_-))\vert \emptyset \rangle$ is in the image of the Sato Grassmanian, therefore the Schur functions satisfy the Plücker relations~\eqref{eq:plucker_part}.

\medskip

The operator $\exp(J_-(\textbf{t}_-))$ can be used to obtain an isomorphism between formal (normally ordered) series in the fermionic operators and formal power series in $\mathbb{C}[[q,\textbf{t}]]$. This isomorphism is known as the boson-fermion correspondance~\cite{JM83}. It is given by the following mapping
\begin{align}
    \label{eq:boson-fermion}
    \Psi \mapsto \bigoplus_{k=-\infty}^\infty q^k\langle k \vert \exp(J_+(\textbf{t}_+)) \Psi \vert 0 \rangle 
\end{align}

The action of any operator on the Clifford algebra can be equivalently described by the action of a differential operator acting on the space of formal power series $\mathbb{C}[[q,\textbf{t}]]$. In particular, observe that if $\Psi$ has charge $n$, then all terms are zero but the monomial in $q^n$. We can now ready to define the $\tau$-functions of the KP hierarchy.

\begin{definition}[$\tau$-function]
    If an element $\Psi$ with zero charge of the Clifford algebra satisfies Equation~\eqref{eq:BBC}, then the formal power series
    \begin{equation}
    \label{eq:tau_fct}
        \tau(\textbf{t}) = \langle 0 \vert \exp(J_+(\textbf{t}_+)) \Psi \vert \emptyset \rangle 
    \end{equation}
    is a KP $\tau$-function.
\end{definition}

A KP $\tau$-function can be expanded over the Schur polynomials $s_\lambda(\textbf{t})$. Inserting the identity $\sum\limits_\lambda \vert \lambda \rangle \langle \lambda \vert $ we get due to the orthogonality of the states $\vert \lambda \rangle$
\begin{equation}
    \label{eq:tau_schur}
    \tau(\textbf{t}) = \sum\limits_\lambda (-1)^{\epsilon(\lambda)}s_\lambda(\textbf{t}) \Delta_\lambda(\Psi)
\end{equation}
where $\Delta_\lambda(\Psi)$ are the Plücker coordinate of $\Psi$. As $\Psi$ satisfies~\eqref{eq:BBC}, $\Psi\vert 0\rangle$ satisfies the Plücker relations~\cite[Sec. $3.3.3$]{AZ13}.

\medskip

Other types of $\tau$-functions can be defined similarly. The $\tau$-functions of the $n$-th modified KP hierarchy are defined similarly with respect to the shifted vacuum $\vert n \rangle$. Also, since if $\Psi$ satisfies~\eqref{eq:BBC}, then so does $\exp(J_+(\textbf{t}_+)) \Psi$ and $\Psi \exp(J_-(\textbf{t}_-)) $, the $\tau$-function of the 2D-Toda lattice (2DTL) are obtained as the expectation value of $\exp(J_+(\textbf{t}_+)) \Psi \exp(J_-(\textbf{t}_-)) $. Fixing $\textbf{t}_+$, it is a KP $\tau$-function in variables $\textbf{t}_-$ and vice-versa.

\subsubsection{The KP equations}
\label{sssec:KP_equa}

In this subsection, we derive the infinite set of partial differential equations satisfied by the KP $\tau$-functions by further exploiting the boson-fermion correspondence. We first need to introduce the fermions generating functions
\begin{equation}
    \label{eq:ferm_gen_fct}
    \psi (z) = \sum\limits_{k\in \mathbb{Z}} \psi_k z^k, \qquad \psi^* (z) = \sum\limits_{k\in \mathbb{Z}} \psi^*_k z^{-k}.
\end{equation}
The function $\psi(z)$ is an eigenvector for the adjoint action of $\exp(J_+(\textbf{t}))$:
\begin{align}
    \label{eq:adj_J_ferm}
    \exp(J_+(\textbf{t})) \psi(z) \exp(-J_+(\textbf{t})) &= \exp(\xi(\textbf{t},z)) \psi(z), \\
    \exp(J_-(\textbf{t})) \psi(z) \exp(-J_-(\textbf{t})) &= \exp(\xi(\textbf{t},z^{-1})) \psi(z),
\end{align}
where $\xi(\textbf{t},z) = \sum\limits_{k \geq 1} t_k z^k$. Similarly $\psi^*(z)$ is an eigenvector for this action with eigenvalue $\exp(-\xi(\textbf{t},z^{\pm}))$.

\medskip

We have the following bosonization formulas
\begin{align}
\label{eq:ferm_to_boson}
    \psi(z) \vert n \rangle &= z^n\exp(J_-([z]))\vert n+1 \rangle, \\
    \psi^*(z) \vert n \rangle &= z^{-n+1}\exp(-J_-([z]))\vert n-1 \rangle,
\end{align}
where $[z] = \left(z,\frac{z^2}{2},\frac{z^3}{3},\dotsc\right)$. Similar relations holds on the dual states $\langle n \vert $.

\medskip

For any $k \in \mathbb{Z}$, the state $\vert 0 \rangle $ is killed by either $\psi_k$ or $\psi^*_k$. Therefore, for an element $\Psi$ satisfying condition~\eqref{eq:BBC}, we get by projecting on states $\langle 1 \vert \exp(J_+(\textbf{t})) \vert \otimes \langle -1 \vert \exp(J_+(\textbf{t}')) \vert$ on the left and $\vert 0 \rangle \otimes \vert 0 \rangle $ on the right
\begin{equation}
\label{eq:bilin_fermion}
    \sum\limits_{k \in \mathbb{Z}} \langle 1 \vert \exp(J_+(\textbf{t})\psi_k\Psi \vert 0 \rangle \langle -1 \vert \exp(J_+(\textbf{t}')\psi^*_k\Psi \vert 0 \rangle = 0.
\end{equation}

The left hand side of this equation corresponds to the residue in $z^{-1}\psi(z)\psi^*(z)$. Using the relations~\eqref{eq:adj_J_ferm} and~\eqref{eq:ferm_to_boson}, it can be recast as a relation on $\tau$-functions as 
\begin{equation}
    \label{eq:bilin_tau_fct}
    \oint_{\mathcal{C}_\infty} \exp(\xi(\textbf{t}-\textbf{t}',z))\tau(\textbf{t}-[z^{-1}])\tau(\textbf{t}'+[z^{-1}])dz = 0.
\end{equation}

Similar equations can be derived for the $\tau$-functions of the MKP and 2DTL hierarchies. This equation contains infinitely PDE equations which are the KP equations. They are extracted by setting $\textbf{t} \rightarrow \textbf{t}+\textbf{a}$, $\textbf{t}' \rightarrow \textbf{t}'-\textbf{a} $ and Taylor expanding around $\textbf{a}=\textbf{0}$. This leads to the following expression.
\begin{equation}
    \label{eq:KP_eq_all}
    \sum\limits_{j \geq 0}h_j(-2\textbf{a})h_{j+1}(\tilde{\partial}_a) \tau(\textbf{t}+\textbf{a})\tau(\textbf{t}-\textbf{a}) = 0.
\end{equation}

Each monomial in $\textbf{a}$ is a PDE  bilinear in $\tau(\textbf{t})$ that must vanish. For example, the first non-trivial relation is obtained from the coefficients in $a_1^4$ and the corresponding equation is
\begin{equation}
    \label{eq:first_KP}
    \tau_{}\tau_{1111} -4\tau_{1}\tau_{111} +3 \tau_{11}\tau_{11} + 3 \tau_{}\tau_{22} - 3 \tau_{2}\tau_{2} -4 \tau_{}\tau_{13} + 4 \tau_{1}\tau_{3} = 0.
\end{equation}
where $\tau_i = \partial_{t_i} \tau$. It is often recast for $F = \log \tau$ where it takes the form
\begin{equation}
    \label{eq:first_KP_F}
    -F_{3,1}+F_{2,2}+\frac{1}{2}F_{1,1}^2+\frac{1}{12}F_{1,1,1,1} = 0,
\end{equation}

which is the KP equation usually found in the literature.

\subsection{Maps and Hermitian matrices as KP \texorpdfstring{$\tau$}{tau}-function}
\label{ssec:maps_matrix_KP}

Now that the stage is set, we show that the partition function of Hermitian matrices is a KP $\tau$-function, and similarly for face and vertex weighted maps generating series.

\subsubsection{Random Hermitian matrix as KP \texorpdfstring{$\tau$}{tau}-function}
\label{sssec:matrix_KP}

Start from the eigenvalue decomposition of the partition function~\eqref{eq:ev_pert_fct}. We will include coupling constants $t_1$ and $t_2$ and omit the scaling with $N$ for now as it can be restored by a rescaling of the coupling constants $t_k \rightarrow Nt_k$. Apply Andreev's formula to express the partition function as the determinant
\begin{equation}
    Z^{(\text{pert})}_{\mathcal{H}_N}\left[t,t_k\right] = N! \det\left( \int d\lambda  \lambda^{i+j-2}\exp\left(-\frac{1}{2t}\lambda^2 + \sum\limits_{k\geq 1} \frac{t_k}{k} \lambda_i^k \right)\right)_{1\leq i,j \leq N}.
\end{equation} The columns can be rearranged by $i \rightarrow N-i$ without changing the sign of the determinant, and we shift $j \rightarrow j-1$ so that both index range in $\llbracket0,N-1\rrbracket$. The potential expands as a series over completely symmetric polynomials $h_k(\textbf{t})$ as
\begin{equation}
    \exp(\sum\limits_{k\geq 1} \frac{t_k}{k} \lambda^k) = \sum\limits_{k \geq 1} h_k(\textbf{t})\lambda^k.
\end{equation}
We can now apply the Cauchy-Binet formula to separate the $x$ and $\textbf{t}$ variables. Via the Jacobi-Trudi formula, we have $\det(h_{\mu_i-i+j}(\textbf{t})) = s_\mu(\textbf{t})$, which leads to
\begin{equation}
\label{eq:Schur_Herm_mat}
     Z^{(\text{pert})}_{\mathcal{H}_N}\left[t,t_k\right] = \sum\limits_{\mu} \det \left(\int d\lambda \lambda^{\mu_i+N-i+j}\exp\left(-\frac{1}{2t}\lambda^2\right) \right)_{0\leq i,j \leq N-1} s_\mu(\textbf{t}).
\end{equation}

Now, we simply have to check that $\det \left(\int d\lambda \lambda^{\mu_i+N-i+j}\exp\left(-\frac{1}{2t}\lambda^2\right) \right)_{0\leq i,j \leq N-1} $ satisfy the Plücker relations, which is straightforward.

\medskip

Alternatively, the matrix model can be shown to satisfy the KP equations by exhibiting a bilinear element in the fermionic operator satisfying~\eqref{eq:BBC} that yields the partition function of the Hermitian matrix model. This element is given in ~\cite[Sec. $4.3.4$]{AZ13} as 

\begin{equation}
    \label{eq:G_matrix_KP}
    \Psi = \normord{\hspace{2pt}\exp\left(\int \psi_+(\lambda)\psi_+^*(\lambda)d\lambda - \sum\limits_{i \geq 0} \psi_i\psi_i^* \right)}.
\end{equation}

\subsubsection{Maps generating function and KP \texorpdfstring{$\tau$}{tau}-function}
\label{sssec:maps_KP}

The expression of the generating function~\eqref{eq:map_sym_fct} involving Schur functions could make it seem like it satisfies the KP hierarchy in both sets of variables $\textbf{x}$ and $\textbf{y}$. However, this generating function is \emph{not} a $\tau$-function of the KP hierarchy~\cite{GJ08}. The solution to the KP hierarchies (or more precisely the $2$DTL) only allow us to control two infinite sets of parameters and, in maps, we are implicitly controlling a third one as the permutation associated with edges must have only transpositions. A natural extension is to relax the condition on $\alpha$ and to introduce a third set of indeterminates $\textbf{t}$ to weight its cycle's length, thus treating the three permutations on equal foot. This gives the generating function for bipartite maps with weights on both sets of vertices and faces
\begin{equation}
    \label{eq:bip_gen_all_weights}
    B(\textbf{x},\textbf{y},\textbf{t}) = \sum\limits_{\substack{n\geq 0 \\ \mu \vdash n}}  s_\mu(\textbf{x})s_\mu(\textbf{y})s_\mu(\textbf{t}).
\end{equation}
Specializing $t_i = t\delta_{i,2}$ forces all cycles to have length $2$ and therefore gives back the generating function of maps.

\medskip

Still, it is possible to keep track of more statistics, but more roughly, via the \emph{content product} of partitions $\lambda$. For a partition $\lambda$ represented via its Ferrers diagram, the content of a box $\Box\in \lambda $ on the $i$-th row and $j$-th column is $c(\Box) = j-i$. For a third set of indeterminates $\textbf{t}$ we have 
\begin{theorem}[\cite{GJ08}]
\label{thm:content_prod}
    $\{\prod\limits_{\Box \in \lambda }t_{c(\Box)}s_\lambda(\textbf{x})) \vert \hspace{2pt}n\geq 0, \hspace{2pt} \lambda \vdash n\}$ \hspace{3pt} satisfies Plücker relations.
\end{theorem}
The $\tau$-functions of the KP hierarchy that can be written under this form are called \emph{hypergeometric}~\cite{OS01}. This allows to construct solutions to the KP hierarchy by specializing the Schur functions of one of the three sets of indeterminates to obtain a content product via the \emph{principal specialization} of the Schur functions~\cite{Stanley}
\begin{equation}
    \label{eq:princ_spec}
    s_\lambda(\textbf{z})_{\vert_{p_i(\textbf{z} = u}} = \frac{n!}{d_\lambda}\prod_{\Box \in \lambda } \left(u+ c(\Box)\right).
\end{equation}
 Setting $u_j = u$ for an indeterminate $t$, we obtain an hypergeometric $\tau$-function corresponding to the principal specialization $s_\lambda(\textbf{t})_{\vert_{p_i = t^{-1}}}$ allowing to track the number of cycles of the permutation $\sigma$ with weight $t^{n-\ell(\sigma)}$:
 \begin{align}
 \label{eq:bip_map_KP}
     \tilde{B}(\textbf{x},\textbf{y},t) &= \sum\limits_{\substack{n\geq 0 \\ \mu \vdash n}}  s_\mu(\textbf{x})s_\mu(\textbf{y}) \prod\limits_{\Box \in \lambda} (1+tc(\Box)) \\
     &= \sum\limits_{\substack{n\geq 0 \\ \mu \vdash n}}  \sum\limits_{\nu \vdash n}  \frac{n!\chi^\lambda_{\nu}}{z_\nu n_\mu} t^{n-\ell(\nu)}s_\mu(\textbf{x})s_\mu(\textbf{y})
 \end{align}

Now, to get back to the generating function of maps, we simply need to require that $\alpha$ is an involution without fixed points. This cannot be achieved via the content product since it also contains partitions which have length $n$ while some parts have size one. Hence we need to specialize either $\textbf{x}$ (or $\textbf{y}$) such that only $p_2(\textbf{x})$ does not vanish. It follows that the only map generating series that is a $\tau$-function of the KP hierarchy is the one where either the face or vertex profile is tracked, and one other parameter can be given a coarser weight via the content product.

%% file: Chapters/Constellations.tex
\chapter{The \texorpdfstring{$b$}{b}-deformation and constellations}
\label{Chap:const}

The interplay between combinatorics of maps, random matrices and the KP hierarchy seen in the previous chapter extends to at least one other case. Similar relations exist between non-oriented maps, symmetric matrices and the BKP hierarchy. On the side of maps, the two generating series can be interpolated via the Jack polynomials~\cite{Ja70} which are a deformation of the Schur functions by a parameter $b$. This interpolation also appears in random matrix models where it is related to $\beta$-ensembles where $\beta$ is related to the parameter $b$ through a simple relation. This suggests that the interplay between the combinatorics of maps, random matrices and symmetric functions could be more general than these two particular cases and hold for arbitrary values of the parameter $b$. However, establishing these connections relies on the representation theory of the symmetric group on the combinatorial side or on the existence of well-known matrix ensembles on the random matrix side, and there is no known analogous structure for a generic value of $b$. On the combinatorial side, the existence of a combinatorial interpretation of the $b$-deformation as a statistic on non-oriented maps has been conjectured by Goulden and Jackson~\cite{GJ_b_conj} via the $b$-conjecture~\cite{LaC_PhD}. Several partial results have been obtained in the past two decades~\cite{LaC_PhD,DF16,Dol17,BD21,BD22} but the conjecture is still open to this day.

\medskip

\emph{Constellations} are a different type of generalization of orientable maps that raised a lot of interest in the past decade. They are combinatorial object generalizing bipartite maps that count ramified covering of the sphere up to homeomorphism and are related to various types of Hurwitz numbers. As with maps, their generating function can be expressed over the Schur functions and taking the principal specialization on all but two sets of indeterminates gives a solution to the KP hierarchy~\cite{Hurwitz_KP_Oko,OS00}. Also, constellations have recently been proven to satisfy topological recursion~\cite{BCCGF22,BDKS22}. Therefore, they are a natural object to study the properties of the $b$-deformation beyond the case of maps.

\medskip

This chapter introduces both the $b$-deformation and constellations to finally arrive at $b$-deformed constellations~\cite{CD22}. We will then study the constraints satisfied by the partition function of the two types of cubic (in a sense we will precise later) $b$-constellations and show that they form an algebra for any value of $b$. These constraints also show that the operator $i\frac{\partial}{\partial p_i}$ can be interpreted as rooting an edge in a face of degree $i$ even for a generic value of $b$ where it is not immediate since the weight in $b$ of the map depends on the choice of root corner. 

\section{The \texorpdfstring{$b$}{b}-deformation}
\label{sec:b_deform}

The algebra of symmetric functions admits several different deformations which retain some of the properties of the Schur functions, and allow us to interpolate between different families of symmetric polynomials. One of these deformations is the so-called $b$-deformation, associated with \emph{Jack polynomials}~\cite{Ja70}. For $b=0$, Jack polynomials are simply the Schur polynomials (up to a rescaling), appearing in orientable map enumeration, random Hermitian matrices and related to the KP hierarchy. At $b=1$ they are the zonal polynomials, which is a family of polynomials related to the enumeration of maps on non-oriented surfaces~\cite{GJ_b_conj,BCD22_bmaps} and orthogonal matrices~\cite{HSS_matrix,BCD22_ON_HCIZ}. While often not directly expanded over the zonal polynomials, these problems have been shown to be related to an integrable hierarchy called the BKP hierarchy~\cite{KVdL98,VdL01}. Yet, aside from these two specific values (and $b=-\frac{1}{2}$ which is related to the case $b=1$ by a duality $ b \longleftrightarrow \frac{-b}{1+b}$) far less is known on the structure of the $b$-deformation for a generic value of $b$. After performing numerical checks, Goulden and Jackson~\cite{GJ_b_conj} noticed that when replacing the Schur functions with the Jack polynomials, the coefficients are polynomials in $b$ with non-negative coefficients. This led them to propose two conjectures about the existence of a combinatorial interpretation of the $b$-deformation. These conjectures are known as the \emph{$b$-conjecture} and \emph{matching Jack conjecture}. Despite several partial results~\cite{LaC_PhD,DF16,Dol17,BD21,BD22}, the two conjectures are still open to this day. In matrix models, the $\beta$-ensembles~\cite{EKR_Review,Desro09,BE09} give an interpolation between well-known matrix ensembles associated with different values of $\beta$ e.g. Hermitian matrices for $\beta =2$ and symmetric matrices at $\beta = 1$. The $\beta$-ensembles admit an expansion over the Jack polynomials~\cite{Desro09} and the parameter $\beta$ is related to $b$ through a simple relation.

\medskip

In this section, we first generalize the algebraic techniques of the previous chapter to the enumeration of non-oriented maps to motivate the introduction of the $b$-deformation. We then present the main properties of the $b$-deformation and some of the difficulties that arise when considering a generic value of $b$.

\subsection{Enumeration of non-oriented maps}
\label{ssec:non_oriented_maps}

\begin{figure}[!ht]
\hfill
\begin{subfigure}{.35\textwidth}
  \centering
  \includegraphics[scale=0.38]{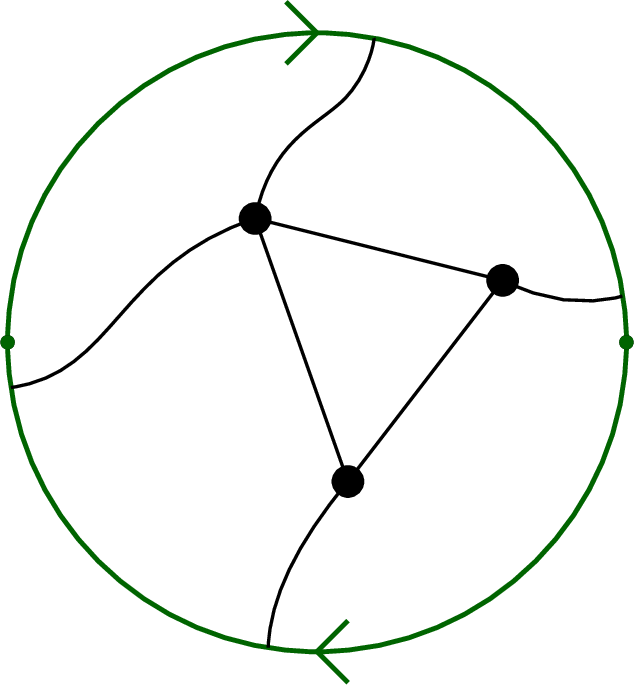}
  %\caption{The projective plane is obtained by identifying the two green edge together.}
  \label{fig:proj_plane_map}
\end{subfigure}
\hfill
\begin{subfigure}{.6\textwidth}
  \centering
  \includegraphics[scale=0.52]{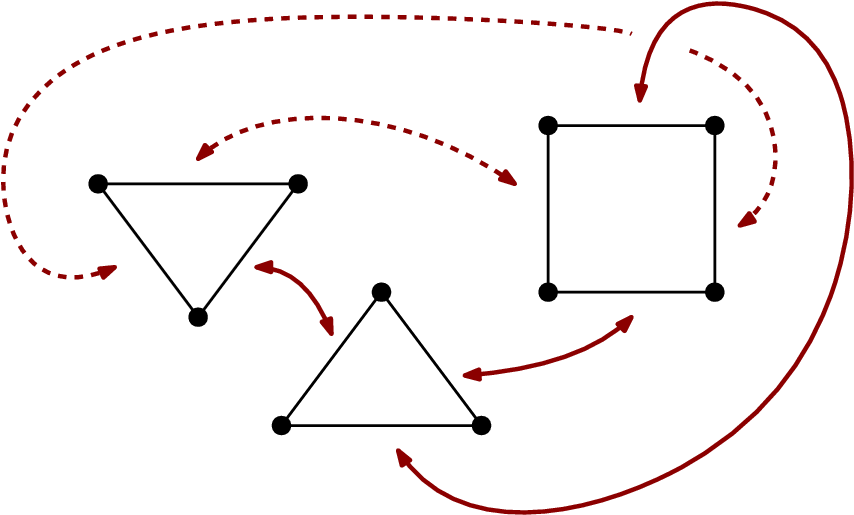}
  %\caption{Dashed arrows represent sides that have to be glued with opposite orientation.}
  \label{fig:proj_plane_map_glu}
\end{subfigure}
\caption{An example of map on the projective plane and one of its polygon gluing. Dashed arrows represent edges glued with opposite orientation.}
\label{fig:proj_plane_ex}
\end{figure}

In the previous chapter, we have seen that orientable maps can be defined as particular graph embedding on an orientable surface. This definition is straightforward to generalize to non-orientable surfaces: a non-orientable map is a graph embedded in a non-orientable surface, such that its complement is a union of disks. For example, the left side of Figure~\ref{fig:proj_plane_ex} shows a map embedded on the projective plane $\mathbb{P}^2$. Cutting the map along its edges, we still obtain a decomposition of the map, but some sides are now glued with opposite orientations with respect to the orientation of their polygons. One possible polygonal gluing is represented on the right of Figure~\ref{fig:proj_plane_ex}, where the dashed arrow represents sides that are glued with opposite orientation. Observe that even for a fixed matching of the sides of the polygons, different gluings can give the same map. For example, turning all dashed arrows into non-dashed ones and conversely, we get a different gluing which gives the same non-orientable map.

\medskip

The definition of maps as a rotation system implicitly relies on the orientability of the underlying surface and thus also needs to be modified. At first sight, it is not obvious how one should generalize this definition to account for non-orientable maps. It turns out to be possible to construct a sub-algebra of the center of the group algebra of $\mathfrak{S}_n$ whose connection coefficients count all maps, both oriented and non-oriented~\cite{HSS_matrix,GJ96_zonal} which allows us to extend to all of the group theoretic tools available in the orientable case to include maps on non-orientable surfaces as well. We give this construction hereafter and show how it relates to zonal polynomials. 

\medskip

When considering non-oriented maps, matchings (i.e. involution without fixed points) play a role similar to permutations in the orientable case. We denote $\mathfrak{M}_n$ the set of matchings on $\llbracket 1,2n \rrbracket$. Given two elements $m_1,m_2 \in \mathfrak{M}_n$, we denote $G(m_1,m_2)$ the graphs with edges the matchings of $m_1$ and $m_2$. By construction, the graph $G(m_1,m_2)$ has cycles of even length and therefore is associated with a partition $2\lambda$ of $2n$. We define $\Lambda(m_1,m_2) = \lambda$ the partition of half-length of cycles of $G(m_1,m_2)$. A (non-necessarily oriented) map is encoded via a \emph{corner system}.

\begin{definition}[Corner system~\cite{GJ96_zonal}]
    A corner system is a triple of matchings \\ $(m_e,m_s,m_c) \in \mathfrak{M}_{4n}$ such that $\Lambda(m_e,m_s) = \left[ 2^n \right]$.
\end{definition}

With this definition, each element of the set $\llbracket 1 , 4n \rrbracket$ should be thought of as labeling one quarter of an edge. Then, the matching $m_e$ assembles quarters that are on the same end of an edge, while $m_s$ assembles those sharing the same side. Therefore, the condition $\Lambda(m_e,m_s) = \left[ 2^m \right]$ is nothing but the condition that these two matchings give the usual notion of edges. Finally, the third matching $m_c$ encodes the structure at the corner of the vertices. Following this correspondence, $\Lambda(m_e,m_c)$ is the face profile of the map while $\Lambda(m_e,m_c)$ gives its vertex profile.

\begin{remark}
    Since a pair of matchings only have cycles of even length, $G(m_e,m_c)$ (or $G(m_s,m_c)$) can always be given a bipartite coloring of its vertices. This coloring encodes the local orientability of the faces (or vertices) of the map.
\end{remark}

To make contact with the representation theory of the symmetric group, we now have to connect matchings with elements of the symmetric group i.e. permutations. This is achieved via the \emph{chain decomposition} of a permutation. We denote $\epsilon_n$ the matching $\{(1 \hspace{2pt} 2) \dotsc (2n-1\hspace{2pt} 2n)$ on the ordered set $\llbracket 1, 2n \rrbracket$.

\begin{definition}[Chain decomposition~\cite{GJ95_comb}]
    For a permutation $\sigma \in \mathfrak{S}_{2n}$, a chain of length $m$ is a list of ordered pairs $((i_1,j_1),\dotsc,(i_{m},j_{m}))$ such that $\sigma(i_{k}) = j_{k}$ and $(i_{2k-1},i_{2k}),(j_{2k},j_{2k+1}) \in \epsilon_n $ for $k\in \llbracket 1,m\rrbracket$ (with $i_{2m+1} = i_1$). The set of chains $c(\sigma)$ of a permutation $\sigma$ is the chain decomposition of $\sigma$.
\end{definition}

By construction, each of the $n$ transpositions of $\epsilon_n$ appear once as a right-link and once as a left-link, both occurring in the same chain. Therefore all chains have even lengths. We define $\kappa(\sigma)$ as the partition of $n$ given by the half-length of its chains ordered in decreasing order. To make the connection with the corner decomposition, notice that to a permutation $\sigma$ we can associate the matching $m(\sigma)$ such that $(i \hspace{2pt}j)\in m(\sigma) \hspace{2pt} \Longleftrightarrow (\sigma(i) \hspace{2pt} \sigma(j)) \in \epsilon_n$.

\begin{figure}
    \centering
    \includegraphics[scale=0.4]{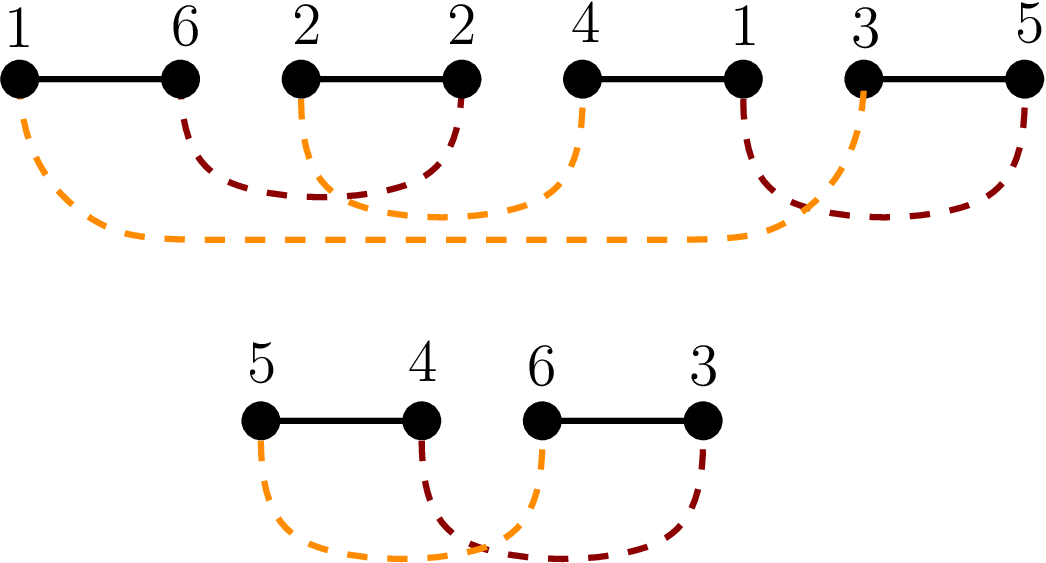}
    \caption{The chain decomposition of $\sigma = (2)(1\hspace{2pt}6\hspace{2pt}3\hspace{2pt}5\hspace{2pt}4)$ with left-links in orange and right-links in red.}
    \label{fig:my_label}
\end{figure}

\medskip

The hyperoctahedral group $B_n$ is the set of elements of $\mathfrak{S}_{2n}$ stabilizing an arbitrary matching. In particular, it can be embedded in $\mathfrak{S}_{2n}$ as the centralizer of $\epsilon_n$. Multiplying to the right of a permutation $\sigma$ with $\delta \in B_n$ relabels the indices $i_k$, while multiplying to the left relabels the indices $j_k$. Therefore, we have the following proposition.
\begin{proposition}[\cite{GJ95_comb,HSS_matrix}]
\label{prop:chain_equiv}
    For any $\sigma,\sigma' \in \mathfrak{S}_{2n}$
    \begin{equation}
        \kappa(\sigma) = \kappa(\sigma') \Longleftrightarrow \exists\hspace{2pt} \delta_1,\delta_2 \in B_n \hspace{4pt}\text{such that} \hspace{4pt} \sigma = \delta_1 \sigma' \delta_2
    \end{equation}
\end{proposition}

This proposition tells us that the role of the conjugacy classes $C_\lambda$ is played by the double cosets of $B_n$ in $\mathfrak{S}_{2n}$ in the non-oriented case who are labeled by partitions of $n$. They are defined as
\begin{align}
    \label{eq:double_coset}
    K_\lambda &= \{ \sigma \in \mathfrak{S}_{2n} \hspace{2pt}\vert\hspace{2pt} \kappa(\sigma) = \lambda \} \\
    &= \{ \sigma \in \mathfrak{S}_{2n} \hspace{2pt}\vert\hspace{2pt} G(\epsilon_n,m(\sigma)) \hspace{4pt} \text{has cycle type} \hspace{4pt} 2\lambda \} \nonumber
\end{align}

We identify $K_\lambda$ with their sum in the group algebra of $\mathfrak{S}_{2n}$. The elements $K_\lambda$ form a commutative subalgebra called the Hecke algebra associated with the Gelfand pair $(B_n,\mathfrak{S}_{2n})$. Its product reads
\begin{equation}
    \label{eq:double_coset_prod}
    K_\mu K_\nu = \sum\limits_{\lambda \vdash n} b^\lambda_{\mu\nu} K_\lambda
\end{equation}
Choosing a representative $[\sigma_\lambda] \in K_\lambda$, the coefficient $b^\lambda_{\mu \nu}$ is the number of pairs $(\rho,\phi )\in K_\mu \times K_\nu$ such that $\rho \phi = \sigma_\lambda$~\cite{HSS_matrix}. They can also be interpreted using the matchings
\begin{lemma}[\cite{HSS_matrix}]
\label{lemma:coeff_b_match}
For any representative $\sigma_\lambda \in K_\lambda$ with matching $m_\lambda = m(\sigma_\lambda)$ we have
\begin{equation}
    \label{eq:coeff_b_match}
    b^\lambda_{\mu\nu} = 2^nn! \vert\{ \sigma \vert G(\epsilon_n,m(\sigma))=\mu, \hspace{3pt} G(m(\sigma),m_\lambda)=\nu \} \vert
\end{equation}
\end{lemma}

A similar interpretation holds when replacing $\epsilon_n$ by any matching $m\in \mathfrak{M}_n$. This allows us to make the connection with the definition of non-oriented maps as a corner system. The number of non-oriented labeled maps with vertex profile $\mu \vdash 2n$ and face distribution $\nu \vdash 2n $ is given by the coefficients $b^{[2^n]}_{\mu\nu}$. The Hecke algebra admits an idempotent basis which allows us to get an expression of these coefficients as characters of $\mathfrak{S}_{2n}$~\cite[Lemma $3.3$]{HSS_matrix}:
\begin{equation}
    b^{\lambda}_{\mu\nu} = \frac{1}{\vert K_\lambda \vert} \sum\limits_{\beta \vdash n } \frac{1}{H_{2\beta}}\phi^{\beta}(\lambda)\phi^{\beta}(\mu)\phi^{\beta}(\nu),
\end{equation}
where $\phi^{\beta}(\lambda) = \sum\limits_{\sigma \in K_\lambda} \chi^{2\beta}(\sigma)$ and $H_{2\beta}$ is an explicit numerical factor.

\medskip

Adding the weights $p_\mu(\textbf{x})$ and $p_\nu(\textbf{y})$ of faces and vertices respectively, we can see that we get a factor of the form
\begin{equation}
    \label{eq:zonal_poly}
    Z_\lambda(\textbf{x}) =  \sum\limits_{\mu \vdash \lambda} \phi^{\lambda}(\mu) p_\mu(\textbf{x}).
\end{equation}
appearing in the expression of the generating series. Up to a factor $\vert B_n \vert^{-1}$, it characterizes the power-sum expansion of zonal polynomial $Z_\lambda$~\cite{James_zonal}, a well-known family of symmetric functions. Therefore, the generating series for non-oriented maps with face distribution $\mu \vdash n$ weighted with $\textbf{x}$, vertex profile $\nu \vdash n$ with weight $\textbf{y}$ and weight $t$ becomes
\begin{equation}
    \label{eq:map_no_gen}
    \tilde{M}(t,\textbf{x},\textbf{y}) = \sum\limits_{n \geq 0} t^n\sum\limits_{\beta \vdash 2n} \frac{\phi^\beta([2^n])}{H_{2\beta}} Z_\beta(\textbf{x})Z_\beta(\textbf{y}).
\end{equation}

This generating series is strikingly similar to the one for orientable maps obtained in the previous chapter~\eqref{eq:map_sym_fct}. They only differ from the type of symmetric functions involved and the numerical prefactor. It generalizes to bipartite maps on non-oriented surfaces by relaxing the condition on $\lambda$, thus treating all three sets of variables symmetrically, leading to the generating series

\begin{equation}
    \label{eq:bip_map_no_gen}
    \tilde{M}(t,\textbf{x},\textbf{y},\textbf{z}) = \sum\limits_{n \geq 0} t^n \sum\limits_{\beta \vdash 2n} Z_\beta(\textbf{x})Z_\beta(\textbf{y})Z_\beta(\textbf{z}).
\end{equation}

For bipartite maps, the difference in the generating series between the orientable and non-oriented case only lies in the numerical prefactor and the family of symmetric function appearing in the generating series.

\subsection{Extension to arbitrary \texorpdfstring{$b$}{b}}
\label{ssec:b_def}

The similitudes in the form of the generating series of orientable~\eqref{eq:bip_gen_all_weights} and non-oriented~\eqref{eq:bip_map_no_gen} bipartite maps hint at the existence of a more general structure encompassing both cases. This generalization is obtained through the introduction of the \emph{Jack polynomials} at the level of the algebra of symmetric functions. This allows for an interpolation between the generating series of orientable and non-orientable maps (up to a global rescaling). Yet at the combinatorial level, there is no family of maps associated with the Jack polynomials for arbitrary values of $b$. This has been the topic of conjectures by Goulden and Jackson~\cite{GJ_b_conj}. Similarly, the Jack polynomials have appeared in matrix models through the $\beta$-ensembles. Less is known on the side of integrable hierarchy, where the generating function of non-oriented maps has been shown to satisfy the BKP hierarchy~\cite{VdL01}. This section introduces Jack polynomials and details what is known about the $b$-deformation from all three points of view.

\subsubsection{The Jack polynomials}
\label{sssec:jack_pol}

The algebra of symmetric functions in variables $\textbf{x}$ can be equipped with a scalar product called the \emph{Hall scalar product} defined in the power-sum basis as 
\begin{equation}
    \label{eq:Hall_scal}
    \langle p_\lambda \vert p_\mu \rangle = \delta_{\lambda\mu} z_\lambda.
\end{equation}

Thus, the Schur functions $s_\lambda$ form an orthogonal basis for the Hall scalar product. Additionally, the expression of the Schur functions in the monomial symmetric function is upper triangular
\begin{equation}
    \label{eq:monom_to_Schur}
    s_\lambda = m_\lambda + \sum\limits_{\mu < \lambda} K_{\lambda\mu} s_\mu,
\end{equation}
for some coefficients $K_{\lambda\mu}$ called the Kostka numbers~\cite{Stanley}. These two properties characterize the Schur functions. For a different choice of scalar product than the Hall scalar product, these two properties define a different basis of symmetric functions. For example, setting
\begin{equation}
    \label{eq:zonal_scalar_prod}
    \langle p_\lambda \vert p_\mu \rangle = \delta_{\lambda\mu} 2^{\ell(\lambda)} z_\lambda 
\end{equation}
defines uniquely the zonal polynomials $Z_\lambda$. The Jack polynomials $P_\lambda^{(b)}$ generalize these two cases, they are defined for the scalar product
\begin{equation}
    \label{eq:scalar_prod_b}
    \langle p_\lambda \vert p_\mu \rangle_b = \delta_{\lambda\mu}(1+b)^{\ell(\lambda)} z_\lambda
\end{equation}
for a parameter $b \in \left]-1;\infty\right[$.

\begin{definition}[Jack polynomials]
\label{def:Jack_poly}
    The Jack polynomials $P_\lambda^{(b)}(\textbf{x})$ are the unique polynomials with coefficients in $\mathbb{Q}(b)$ such that
    \begin{itemize}
        \item {Their expansion over symmetric monomial $(m_\nu)_\nu$ is triangular superior 
        \begin{equation}
            P_\lambda^{(b)}(\textbf{x}) = m_\lambda(\textbf{x}) + \sum\limits_{\mu < \lambda} c_{\lambda\mu} m_\mu(\textbf{x})
        \end{equation}}
        \item {They are orthogonal w.r.t the $b$-deformed scalar product $\langle \hspace{2pt}\vert \hspace{2pt} \rangle_b$.}
    \end{itemize}
\end{definition}

Alternatively, they can be defined by replacing the orthogonality condition with the requirement that they are the eigenvectors of the Laplace-Beltrami operator
\begin{align}
    \label{eq:Laplace_Beltrami}
    D^{(b)} = \sum\limits_{i \geq 1} (1+b) x_i^2\partial_i^2 + 2 \sum\limits_{j \neq i} \frac{x_ix_j}{x_i-x_j}\partial_i 
\end{align}

With this choice of normalization, the Jack polynomials at $b=0$ are rescaled Schur functions $H_\lambda s_\lambda(\textbf{x})$ while they are the zonal polynomial at $b=1$. Introducing the $b$-deformation leads to modification of the usual statistics on partitions. For example, for a box $\Box$ of a partition $\lambda$ with coordinate $(i,j)$, its content becomes $c_b(\Box)= (1+b)(i-1)+(j-1)$, which follows from the principal specialization of Jack polynomials $P_\lambda^{(b)} (\textbf{1}) = \prod_{\Box\in\lambda} (n + c_b(\Box))$. We define two deformed hook length $h_b(\Box)$ and $h'_b(\Box)$ are defined as
\begin{align}
    h_b(\Box) &= (1+b)a(\Box) + \ell(\Box) + 1 ,\\
    h'_b(\Box) &= (1+b)a(\Box) + \ell(\Box) + 1+b.
\end{align}
The corresponding hook products are defined as the product over boxes of $\lambda$ of the hook length i.e. $H_b(\lambda) = \prod\limits_{\Box \in \lambda} h_b(\Box) $ and $H'_b(\lambda) = \prod\limits_{\Box \in \lambda} h'_b(\Box) $. 
The hook products compute the norm $j_\lambda^{(b)}$ of the Jack polynomials for the scalar product~\eqref{eq:scalar_prod_b}
\begin{equation}
    \label{eq:norm_jack_b}
    j_\lambda^{(b)} = \langle P^{(b)}_\lambda \vert P^{(b)}_\lambda \rangle_b = H_b(\lambda)H'_b(\lambda).
\end{equation}
In particular, the scaling factor of the Schur functions is nothing but the usual hook product at $b=0$.

\subsubsection{The \texorpdfstring{$b$}{b}-conjecture(s)}
\label{sssec:b-conj}

Generalizing the generating series of orientable and non-oriented bipartite maps, we define the following generating series
\begin{equation}
    \label{eq:Jack_gen_series}
    \phi(\textbf{x},\textbf{y},\textbf{z}) = \sum\limits_{n \geq 0} t^n \sum\limits_{\lambda \vdash n} \frac{1}{\langle P_\lambda^{(b)},P_\lambda^{(b)} \rangle_b}  P_\lambda^{(b)}(\textbf{x})P_\lambda^{(b)}(\textbf{y})P_\lambda^{(b)}(\textbf{z})
\end{equation}

This generating series interpolates between the (rescaled) generating series for orientable bipartite maps at $b=0$ and the generating series of non-oriented bipartite maps at $b=1$. Changing basis from the Jack polynomials to the power sum, this series can be expanded as
\begin{equation}
    \label{eq:b_series_PS}
    \phi(\textbf{x},\textbf{y},\textbf{z}) = \sum\limits_{n \geq 0} t^n \sum\limits_{\lambda \vdash n} \frac{c^\lambda_{\mu\nu}(b)}{(1+b)^{\ell(\lambda)}z_\lambda}  p_\lambda(\textbf{x})p_\mu(\textbf{y})p_\nu(\textbf{z})
\end{equation}

The coefficients $c^\lambda_{\mu\nu}(b)$ should a priori be rational functions of $b$. However, in the explicit cases where these coefficients have been computed in~\cite{GJ_b_conj} (and later in subsequent work e.g.~\cite{LaC_PhD}), they turn out to be integer valued polynomials in $b$. This led Goulden and Jackson to conjecture that the $b$-deformation has a combinatorial interpretation in terms of matchings

\begin{conjecture}[Matching Jack conjecture~\cite{GJ_b_conj}]
There exists non-negative weight function $\text{wt}_\lambda$ such that, for all $\lambda,\mu,\nu \vdash n$
\begin{equation}
    \label{eq:match_jack_coef}
    c^\lambda_{\mu\nu}(b) = \sum\limits_{m} b^{\text{wt}_\lambda(m)}
\end{equation}
where the sum is over matchings such that $G(\epsilon_n,m) = \mu$ and $G(\delta,m_\lambda) = \nu $ for any matching $m_\lambda$ associated with a representative $[\sigma_\lambda] \in K_\lambda$
\end{conjecture}

A similar conjecture -called the $b$-conjecture~\cite{GJ_b_conj,LaC_PhD}- is made at the level of the connected generating series. It is not a direct implication of the matching Jack conjecture because of the difference in the factors in $b$ appearing in both series, but both conjectures are expected to be combinatorially related. One of the difficulties that arise in the study of this conjecture is that to this day, there is no known way to use the representation theory of the symmetric group for a generic value of $b$, which would help make the connection with the power-sum expansion of Jack polynomials~\cite[Cor. $4.3.1$]{HW_Jack_exp}.

%(Brief survey of the results on the $b$-conjecture ?) (Details on the computation of $b$-weight ?)

\subsubsection{The matrix \texorpdfstring{$\beta$}{B}-ensembles}
\label{ssec:beta_ens}

In random matrix models, the Gaussian $\beta$-ensembles~\cite{Meh04,EKR_Review} interpolates between different diagonalizable matrix models with real eigenvalues. Its partition function is defined as
\begin{equation}
\label{eq:Gauss_beta}
    Z[\beta] = \int_{\mathbb{R}^N} \prod\limits_{i=1}^N d\lambda_i \prod\limits_{i>j} \vert \lambda_i - \lambda_j \vert^{\beta} \exp(- \sum\limits_{i=1}^N x_i^2+ V(\lambda_i))
\end{equation}
where $\beta$ is a positive real number and $V(\lambda_i)$ is a potential (a polynomial in most cases). This partition function can be interpreted as the partition function of a $1D$ classical Coulomb gas of $N$ particles in a potential $V$. For $\beta=2$, we recognize the eigenvalue decomposition of random Hermitian matrices given by Equation~\eqref{eq:ev_pert_fct}. For $\beta = 1$ it corresponds to the partition function of symmetric matrices, which can be diagonalized by an orthogonal matrix $O\in O(N)$. We have seen that the Hermitian matrix model admits an expansion over Schur functions in~\eqref{eq:Schur_Herm_mat}. Similarly, integrals over symmetric matrices can be expanded over zonal polynomials~\cite{James_zonal,BCD22_ON_HCIZ}. The $\beta$-ensembles give an interpolation between the Hermitian and symmetric random matrices which is akin to the relation between the generating series of face-weighted orientable and non-oriented maps and the parameter $\beta$ of matrix model is related to the parameter $b$ of symmetric functions via the relation $\beta = \frac{2}{(1+b)}$~\cite{Desro09,BCD22_ON_HCIZ}. 

\begin{remark}
Note that on the combinatorial side, the maps generating series only carry weights on their faces. To this day, there is no known analog of the beta ensembles which would correspond to the generating series of maps with weights on both faces and vertices. 
\end{remark}

The $\beta$-ensembles are defined for any positive value of $\beta$ and they admit an expansion over Jack polynomials~\cite{Desro09}. Outside the particular values $\beta=1,2,4$ corresponding respectively to the eigenvalue decomposition of symmetric, (complex) Hermitian and quaternionic Hermitian matrices, there is no known explicit random matrix ensemble that admits an eigenvalue decomposition of the form~\eqref{eq:Gauss_beta}. These three particular values of $\beta$ correspond to matrices with real eigenvalues which are diagonalized by compact Lie groups (respectively orthogonal, unitary and symplectic) and it is unclear what properties should still hold for a generic value of $\beta$ where it is a priori no longer the case. Yet, it is possible to generalize some results for arbitrary $\beta$. For example the Harish-Chandra-Itzykson-Zuber integral
\begin{equation}
\label{eq:IZ_int}
    I(X,Y) = \int_{U_N} dU \exp\left(-\Tr(XUYU^{-1})\right)
\end{equation}
can be extended to an arbitrary value of $\beta$~\cite{BE09}. It is then defined as an integral over abstract Lie groups $U_N^\beta$ on which integration can be performed using an eigenvalue equation known as the Calogero-Moser equation.

\subsubsection{\texorpdfstring{$b$}{b}-deformation and integrable hierarchies}
\label{sssec:b_hierarch}

In the previous paragraph, we have shown that the $b$-deformation has a natural interpretation as interpolating between different classes of maps or random matrices and this correspondence also holds for integrable hierarchies. The generating function of symmetric matrices, which is related to non-oriented maps, also satisfies an infinite set of PDE known as the BKP hierarchy~\cite{VdL01} which generalizes the KP hierarchy. However, contrary to matrices and maps, there is no clear candidate structure to generalize the integrable hierarchies to arbitrary values of $b$. It is unclear whether the integrability properties of the generating series of maps and random matrices hold for a generic value of $b$. Fortunately, integrability has several different facets which can be used to fill this gap. In the orientable case, all the examples above satisfy topological recursion, and their generating functions satisfy constraints that form a particular type of $\mathcal{W}$-algebra known as the Virasoro algebra. These other aspects can be used to try to probe to robustness of the integrability properties for a generic value of $b$. For example, it is known that the $b$-deformed generating functions of maps satisfy Virasoro constraints for all values of $b$~\cite{AvM01}.

\section{Constellations and their \texorpdfstring{$b$}{b}-deformation}
\label{sec:const}

Algebraically, constellations count ramified covering of the sphere up to homeomorphism~\cite{LZ_book}. Combinatorially, they generalize the algebraic properties of maps. As such, their generating function admits an expansion over Schur functions and they yield similar integrability properties. Their generating function has specializations which are $\tau$-function of the KP hierarchy and they have been recently been shown to satisfy topological recursion~\cite{BCCGF22,BDKS22}. Therefore, they are a natural object to probe deeper into the properties of the $b$-deformation. The corresponding combinatorial object -which we will still call constellations unless we explicitly want to make the distinction with the orientable case explicit- has been defined in~\cite{CD22}.

\subsection{Orientable constellations}
\label{ssec:orient_const}

Orientable constellations (which we will simply call constellations) generalize the definition of maps as a rotation system to an arbitrary number $k+2$ of permutations.

\begin{definition}[Constellations~\cite{LZ_book}]
\label{def:const}
A (labeled and orientable) $k$-constellation of size $n$ is a $(k+2)$-uple of permutations $\left(\sigma_0,\dotsc,\sigma_{k},\phi\right)\in \mathfrak{S}_n$ such that the product $\phi \sigma_0\dotsc\sigma_{k}$ is a factorization of  $\id_{\mathfrak{S}_n}$. A constellation is connected if the action of $\langle\sigma_0,\dotsc,\sigma_{k},\phi\rangle$ on $\{1\dotsc,n\}$ is transitive.
\end{definition}

In particular, a $1$-constellation of size $n$ corresponds to bipartite maps with $n$ edges. In the general case, the corresponding combinatorial object is obtained as follows. As for bipartite maps, the cycle type $\lambda^{(i)}$ of each of the $m$ permutations $\sigma_i$ encodes the profile of vertices carrying color $i$ and the cycle type $\lambda^{-1}$ of $\phi$ encodes the face profile of the constellation. All faces have length multiples of $k$ and the successor in clockwise order of a vertex of color $c$ has color $c+1$ modulo $k+1$. The main difference comes from the objects that carry the labels $\{1,\dotsc,n\}$ corresponding to edges in the case of bipartite maps. For constellations, they become \emph{hyperedges}. A hyperedge is a face of length exactly $m$ whose vertices have color $0$ to $k$ in clockwise order. When $m=2$ a hyperedge is simply a digon with vertices colored $0$ and $1$ and collapsing the two edges of the digon gives back the notion of edge we are familiar with.

\begin{figure}[!ht]
    \centering
    \includegraphics[scale=0.5]{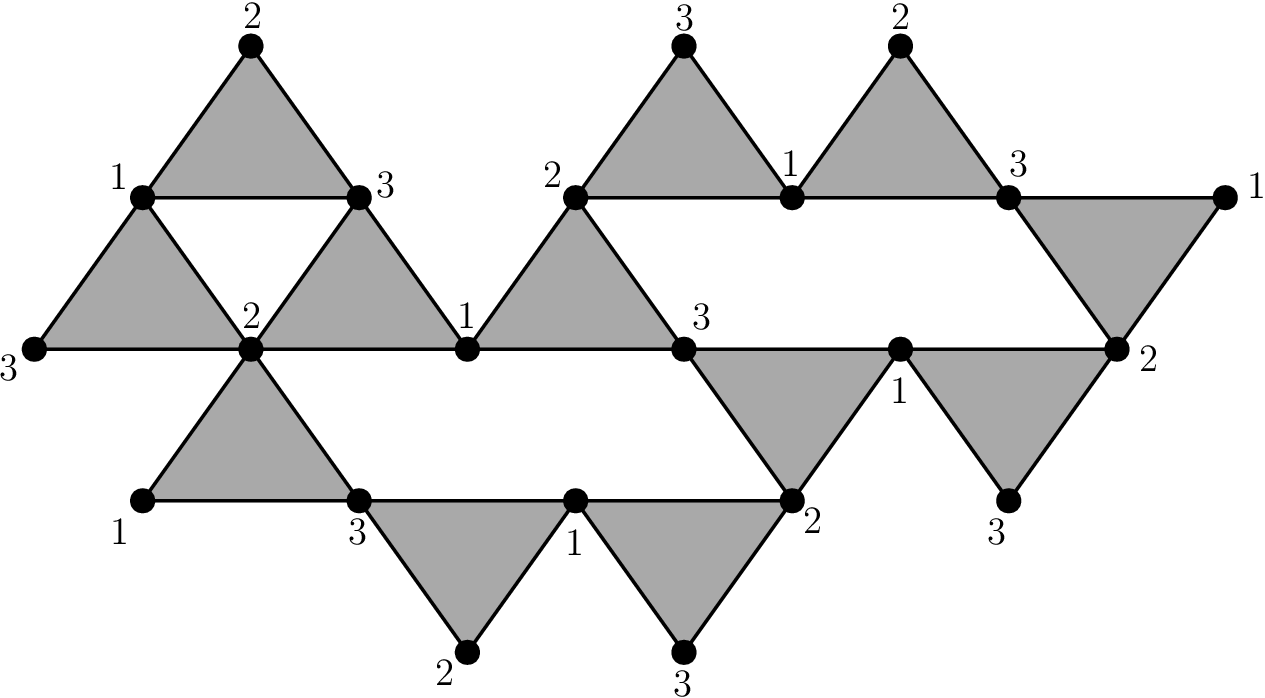}
    \caption{An example of orientable $3$-constellation. Hyperedges are represented in grey.}
    \label{fig:3const}
\end{figure}

\medskip

\begin{figure}[!ht]
    \centering
    \includegraphics[scale=0.7]{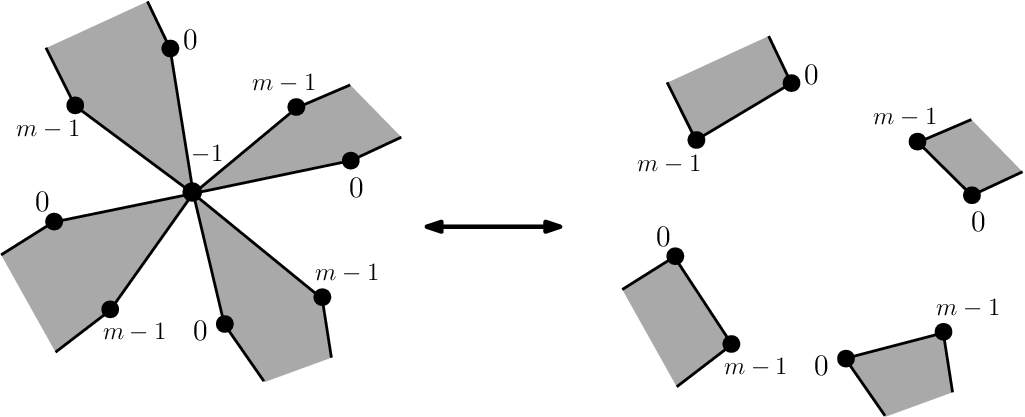}
    \caption{The bijective move between pure $k+1$-constellations and $k$-constellations.}
    \label{fig:pure_to_const}
\end{figure}

Constellations have a natural interpretation in enumerative geometry as they count ramified coverings of the Riemann sphere $\mathbb{CP}^1$. Let $X$ be a Riemann surface and $f:X\to \mathbb{CP}^1$ a ramified covering of degree $n$ of $\mathbb{CP}^1$ with $(k+1)$ ramification points $y_{-1},\dotsc,y_{k} \in \mathbb{CP}^1$. Up to homeomorphism, we can assume the points $\{y_i\}$ lie on the equator $\Gamma$ of the complex sphere and are encountered with increasing labels (modulo $k+2$) when following the equator in clockwise order. Both hemispheres contain no ramifications points, therefore $f^{-1}(\Gamma)$ splits $X$ into $2n$ distinct faces, $n$ in which vertices are visited with increasing colors and $n$ with decreasing colors in clockwise order. Each of these $n$ faces corresponds to one of the $n$ preimage of the Northern (or Southern) hemisphere of $\mathbb{CP}^1$. The preimages of the Southern hemisphere correspond to hyperedges of Definition~\ref{def:const} which are labeled from $1$ to $n$ and the preimages $f^{-1}(y_i)$ encode the ramification profile of $y_i$. The graph thus obtained is a $(k+1)$-constellation where all faces are required to have length $k+2$, which are called \emph{pure} constellations. To obtain a constellation as given in Definition~\ref{def:const}, we have to transform vertices of color $-1$ and degree $p$ into faces of length $p(k+1)$. To do so, we cut each of the $n$ clockwise faces along their vertices of color $k$ and $0$. Thus, each vertex of color $-1$ becomes the center of a face and if the vertex initially had degree $p$, then $p$ faces are merged together and the resulting face has length $p(k+1)$. This operation is represented in Figure~\ref{fig:pure_to_const}. 

This process is bijective. The vertices of color $-1$ can be recovered by adding a vertex of color $-1$ at the center of each face and adding edges between $-1$ and vertices of color $1$ and $k$. Thus we have a bijection between $k$-constellations and $(k+1)$-pure constellations.

\medskip

When the constellation is connected, the genus of the initial Riemann surface can be recovered from the Riemann-Hurwitz formula. The pure constellations has exactly $2n$ faces, $(k+2)n$ edges and $\ell(\lambda^{(i)})$ vertices of color $i$ ($\ell(\phi)$ for vertices of color $-1$). Therefore
\begin{proposition}[Riemann-Hurwitz formula]
\begin{equation}
\label{eq:RH_formula}
    \chi(X) = 2 -2g(X) = 2n - (n-\ell(\phi)) - \sum\limits_{i=0}^{k} \left(n-\ell(\lambda^{(i)})\right).
\end{equation}
\end{proposition}

As ramified coverings, we require $y_i$ to be ramification points. This implies that vertices (resp. faces for $y_{-1}$) cannot have only degree $1$ (resp. all length $1$). This condition is relaxed in our combinatorial definition which amounts to allowing for trivial "ramification points".  By relaxing this condition, $k$-constellations can be embedded as $(k+1)$-constellations by adding a trivial "ramification point" which enables to define the projective limit of constellations with arbitrary (but finitely many) ramifications points. Via the Riemann-Hurwitz formula, it is possible to enumerate other families of constellations by placing constraints on the values of $\ell$ rather than fixing the number of ramification points $k+1$, as we shall see examples of in Section~\ref{ssec:Hurwitz_nb}.

\medskip 

The algebraic method presented in the previous chapter for maps (see Sec.~\ref{ssec:map_sym_group}) can be extended straightforwardly to the case of constellations. We just need to add more conjugacy classes to the computation of~\eqref{eq:map_count}. In particular, the number of $k$-constellations of size $n$ with vertex profile $(\lambda_0,\dotsc,\lambda_{k+1}) \vdash n$ and face profile $\nu \vdash n$ is given by
\begin{align}
    \label{eq:const_count}
    K_{\lambda_0\dotsc\lambda_{k}\nu} &= \left[C_{[1^n]}\right] C_{\lambda_0} \dotsc C_{\lambda_{k}}C_{\nu} \\
    &= \frac{(n!)^{k+1}}{z_{\lambda_0}\dotsc z_{\lambda_{k}}z_{\nu}} \sum\limits_{\mu \vdash n} \frac{\chi^\mu_{\lambda_0}\dotsc\chi^\mu_{\lambda_{k}}\chi^\mu_{\nu}}{d_\mu^{k}}\nonumber
\end{align}

Associating to the vertex of color $i$ indeterminates $\textbf{t}^{(i)}$ and faces to indeterminates $\textbf{p}$ and an indeterminate $t$ to the size of the constellation, the generating function of weighted $k$-constellations admits an expansion over Schur functions as
\begin{equation}
    \label{eq:const_weight}
    \tilde{Z}_k(t,\textbf{q}^{(i)},\textbf{p}) = \sum\limits_{n \geq 0} \frac{t^n }{n!} \sum\limits_{\mu \vdash n} \frac{(n!)^{k+1}}{d_\mu^{k}} s_\mu(\textbf{q}^{(0)})\dotsc s_\mu(\textbf{q}^{(k)})s_\mu(\textbf{p})
\end{equation}

As with bipartite maps, this function is not a KP $\tau$-function, but $\tau$-function can be obtained from the principal specialization of $k-1$ sets of indeterminates e.g. $p_k(\textbf{t}^{(i)}) = u_i$ for vertices of color $i \geq 1$. This generating series is a KP $\tau$-function since it can be expressed as a content product via Theorem~\ref{thm:content_prod}.

\begin{equation}
    \label{eq:const_KP}
    Z_k(t,\textbf{q},\textbf{p},u_1,\dotsc,u_{k}) = \sum\limits_{n \geq 0} t^n \sum\limits_{\mu \vdash n} \frac{(n!)^{k}}{n_\mu^{k}} s_\mu(\textbf{q})s_\mu(u_1)\dotsc s_\mu(u_{k})s_\mu(\textbf{p})
\end{equation}

\subsection{Relation with Hurwitz numbers}
\label{ssec:Hurwitz_nb}

One of the main interest surrounding constellations lies in their link with \emph{Hurwitz numbers}. The Hurwitz numbers count factorization of the identity with any number of permutations while controlling the cycle type of some of them. Therefore, they correspond combinatorially to constellations where the number of color $m$ is not held fixed. Different conditions on the factorization lead to different families of Hurwitz numbers, but a wide variety of Hurwitz numbers share the same properties. The prototypical example of Hurwitz numbers is the \emph{simple} Hurwitz numbers (see their definition below). It is the case understood best and showcases all of the expected properties of Hurwitz numbers: they admit a combinatorial expansion through their link with constellations but they can also be computed as an integral over the moduli space of complex curves via the ELSV formula~\cite{ELSV}. Finally, they are integrable as they satisfy their generating function and satisfy the KP hierarchy~\cite{Hurwitz_KP_Oko}. Thus, the Hurwitz number establishes a link between constellations and intersection theoretic quantities. In the last two decades, these results have been extended to other families of Hurwitz numbers. We illustrate their relevance for integrable hierarchies and matrix models through the following three examples. 

\paragraph{Simple Hurwitz number}

The simple Hurwitz numbers $h_{\mu,k}$ count the number of factorizations of a permutation $\sigma$ with cycle type $\mu$ in a product of $k$ transpositions. Weighting cycle length of $\mu$ with variables $\textbf{p}$ and $u$ for each transposition, the corresponding generating function is
\begin{equation}
\label{eq:Hurwitz_gen}
    H^d(\textbf{p},u) = \sum\limits_{\substack{n \geq 0 \\ \mu \vdash n}} \sum\limits_{m \geq 0} h_{\mu,m}\textbf{p}_\mu \frac{u^m}{m!}
\end{equation}
Via the Riemann Hurwitz formula~\eqref{eq:RH_formula}, we have
\begin{equation}
\label{eq:RH_Hurwitz}
    \chi_\mu = n + \ell(\mu) - m.
\end{equation}
Hence for a fixed partition $\mu$, increasing the number of transposition factorizing $\mu$ increases the genus of the constellation. This generating series is similar to the constellation generating series in that the variables $\textbf{p}$ keep track of the ramification profile of the branched covering, except that we now allow $m$ to fluctuate and impose that all $m$ permutations are simple instead of controlling directly the genus. Despite these differences in nature between the two generating series, the generating series of Hurwitz numbers also satisfies the KP hierarchy and its description as a matrix model was established in~\cite{BEMS_Hurwitz_MM}. 

\paragraph{Double Hurwitz number}

As we have seen in the previous Chapter, the KP hierarchy is contained in the $2D$-Toda lattice, which allows us to keep track of two infinite sets of indeterminates $\textbf{p}$ and $\textbf{q}$ and satisfying KP equations in both sets of variables. Thus the previous generating function~\eqref{eq:Hurwitz_gen} can be generalized to track of the cycle type of a second permutation while satisfying the 2D Toda equations~\cite{Hurwitz_double_OKO}. The corresponding Hurwitz number are called the double Hurwitz numbers and denoted $h_{\mu,\nu,m}$. They are defined similarly to the simple ones, but now factorizing two permutations $\mu$ and $\nu$ of size $n$ simultaneously as a product of transposition. Their generating function is
\begin{equation}
    \label{eq:double_Hurw}
    H^s(\textbf{p},\textbf{q},u) = \sum\limits_{\substack{n \geq 0 \\ \mu,\nu \vdash n}}  \sum\limits_{m \geq 0} h_{\mu,\nu,m}\textbf{p}_\mu \textbf{q}_\nu \frac{u^m}{m!}
\end{equation}

While the double Hurwitz numbers admit no representation as a matrix model to this day, they yield all other features present for simple Hurwitz numbers~\cite{Hurwitz_double_BDKLM}.

\paragraph{Monotone Hurwitz number}

Consider the HCIZ integral~\eqref{eq:IZ_int}. Its moments are
\begin{equation}
\label{eq:HCIZ_moment}
I(X,Y)^{(d)} = N^d \int_{U_N} dU \left(\Tr(XUYU^{-1})\right)^d.  
\end{equation}
Using invariance by right and left multiplication of the Haar measure of $U(N)$, they expand over the eigenvalues $(x_i)_{1 \leq 1 \leq N}$,$(y_i)_{1 \leq 1 \leq N}$ of $X$ and $Y$ as~\cite{GGPN_Mono_HCIZ} 
\begin{equation}
    I(X,Y)^{(d)} = N^d \sum\limits_{\sigma,\rho} x_{\sigma(1)}\dotsc x_{\sigma(d)}y_{\rho(1)}\dotsc y_{\rho(d)} \int_{U_N} \vert U_{\sigma(1)\rho(1)}\dotsc U_{\sigma(d)\rho(d)} \vert^2 dU
\end{equation}
where the sum is over functions $\rho,\sigma: \llbracket,1,d\rrbracket \rightarrow \llbracket,1,N\rrbracket$
This integral over the unitary group is called the \emph{Weingarten function} and can be computed via \emph{Weingarten calculus}~\cite{CS06}. Weingarten calculus recasts the computation of this integral as a combinatorial problem of counting factorization (with suitable weight in $N$) of the permutation $\rho^{-1}\sigma \in \mathfrak{S}_d $ in transpositions with the additional constraints that the sequence of transposition must be monotonous, that is for a finite product of transposition $(a_1 \hspace{2pt} b_1)\dotsc (a_k \hspace{2pt} b_k)$ with $a_i < b_i$ for all $i$ it must satisfy $b_i \leq b_{i+1}$ for $1\leq i \leq k-1$. Hence the situation is similar to simple Hurwitz numbers with this additional monotony constraint, which defines the monotone Hurwitz numbers. These numbers also satisfy the KP hierarchy~\cite{ALS_Mono_Hurw}.

\subsection{The \texorpdfstring{$b$}{b}-deformation for constellation}
\label{ssec:b_def_const}

Since the generating function of constellations~\eqref{eq:const_KP} admits a Schur expansion, it can naturally be $b$-deformed by replacing the Schur functions with the Jack polynomials,
\begin{equation}
\label{eq:b_const_KP}
    Z_k^{(b)}(t,\textbf{q},\textbf{p},u_1,\dotsc,u_{k}) = \sum\limits_{n \geq 0} t^n \sum\limits_{\mu \vdash n} \frac{1}{j^{(b)}_\mu} P^{(b)}_\mu(\textbf{q})P^{(b)}_\mu(u_1)\dotsc P^{(b)}_\mu(u_{k})P^{(b)}_\mu(\textbf{p}).
\end{equation}

However, contrary to constellations which have a clear combinatorial interpretation, it is (at this stage) not clear what combinatorial object this generating function counts. The combinatorial structure of constellations makes it amenable to combinatorial techniques like Tutte-like equations or slice decomposition~\cite{BCCGF22} to study its generating function. The lack of a combinatorial interpretation for this generating function essentially reduces the set of tools available to study its properties to algebraic techniques. Fortunately, this gap was partially filled in~\cite{CD22} where the notion of \emph{generalized branched coverings} was introduced. The $b$-deformed constellations are the combinatorial object associated with the generalized branched coverings and their (specialized) generating function is given by~\eqref{eq:b_const_KP}. Yet some difficulty inherent to the $b$-deformation (e.g. dependency on the rooting~\cite{LaC_PhD}) still makes it difficult to apply usual combinatorial methods to these objects. We follow here the presentation of $b$-constellations of~\cite{CD22}.

\subsubsection{Generalized branched coverings and constellations}
\label{sssec:b_def_const}

As we have seen through the example of maps at $b=1$, the $b$-deformation is related to the non-orientability of the surfaces. For a Riemann surface $X$ and a point $x \in X$, a local orientation at $x$ is a choice of a direct basis $\{e^{(x)}_i\}$ for the tangent space $T_x X$. Any other basis $\{{e'_i}^{(x)}\}$ of $T_x X$ has the same orientation as $\{e^{(x)}_i\}$ if the transition matrix between the two bases has positive determinant. Hence for each point $x\in X$ there exist two different orientations. A choice of orientation at $x$ extends to some neighbourhood $U_x$ of $X$ by requiring that for all $y \in U_x$, the determinant of the transition function between tangent space $T_x X$ and $T_y X$ in their respective basis $\{e^{(x)}_i\}$ and $\{e^{(y)}_i\}$ has positive determinant. When it extends to the whole of the Riemann surface $X$, then $X$ is orientable. We denote $\tilde{X}$ its orientation double cover. It is the set of points of $X$ equipped with a local orientation. When $X$ is orientable, its double cover $\tilde{X}$ is exactly two disconnected copies of $X$ where each copy corresponds to a choice of local orientation at any point $x$. Conversely, when $X$ is non-orientable then $\tilde{X}$ is connected. The orientation double cover comes with a $2$-to-$1$ projection to $X$ noted $\pi$ which simply forgets the choice of orientation. While constellations count branched coverings of the sphere, their $b$-deformed counterpart counts \emph{generalized branched coverings} where the initial Riemann surface $X$ is allowed to be non-orientable. We denote $\mathbb{CP}_+^1$ and $\mathbb{CP}_-^1$ the Northern and Southern hemispheres of the complex Riemann sphere $\mathbb{CP}^1$. We denote $p:\mathbb{CP}^1 \rightarrow \mathbb{CP}_+^1$ the projection identifying both hemispheres. 

\begin{definition}[Generalized branched coverings~\cite{CD22}]
\label{def:gen_branch}
    A continuous map \\ $f:X \rightarrow \mathbb{CP}_+^1$ is a generalized branched covering of the sphere if there exists a branched covering $\tilde{f}:\tilde{X}\rightarrow \mathbb{CP}^1$ such that
    \begin{equation}
        \label{eq:def_gen_branch}
        f \circ \pi = p \circ\tilde{f}
    \end{equation}
    Two generalized branched coverings $f$ and $f'$ are equivalent if $\tilde{f}$ and $\tilde{f}'$ are equivalent as (usual) branched coverings.
\end{definition}

Since $\tilde{X}$ is a double cover of $X$, it follows that if $f$ has degree $n$ then the ramification profile of $\tilde{f}$ over a ramification point $\tilde{y}\in \mathbb{CP}^1$ will be given by a partition $\tilde{\lambda} = (\lambda_1, \lambda_1, \dotsc, \lambda_{\ell(\lambda)},\lambda_{\ell(\lambda)}) \vdash 2n $. The ramification profile of $f$ over $y \in \mathbb{CP}_+^1$ is then given by $\lambda$.

\medskip

If $f$ has $(k+2)$ ramification points $y_{-1},y_0,\dotsc,y_{k}$ along the equator $\partial\mathbb{CP}_+^1$ encountered in that order, the corresponding constellations are given by $f^{-1}(\partial\mathbb{CP}_+^1)$ which gives a pure $(k+1)$-constellation. As in the orientable case, there is a bijection between pure $(k+1)$-constellations and $k$-constellations which is obtained by considering the preimage of the segment $\left[y_0,y_{k-1}\right] \subset \partial\mathbb{CP}_+^1$ instead. The generalized constellations are represented using ribbon graphs via the following definition. 

\begin{definition}[Generalized constellations~\cite{CD22}]
    \label{def:const_ribbon_graph}
    A (generalized) $k$-constellation $\mathcal{C}$ is a map equipped with a coloring of its vertices in $\llbracket 0,k\rrbracket$ such that
    \begin{itemize}
        \item Vertices of color $0$ (resp. $k$) have neighbours of color $1$ (resp $k-1$),
        \item Each corner of a vertex of color $i \in \llbracket 1,k-2\rrbracket$ separates between a vertex of colour $i-1$ and a vertex of color $i+1$.
    \end{itemize}
\end{definition}

The degree of a face is given by its number of vertices of color $0$ (or equivalently $k-1$ or half the number of vertices in any other color). The size of a constellation is the sum of the degree of its faces. The profile of the constellations is the $(k+2)$-uple of partitions $(\lambda^{(-1)},\lambda^{(0)},\dotsc,\lambda^{(k)})$ such that $\lambda^{(-1)}$ encodes the face distribution and $\lambda^{(i)}$ the profile of vertices of colour $i$. A labeled constellation of size $n$ has its corner of color $0$ labeled with unique integers in $\llbracket 1,n \rrbracket$ and an orientation. Finally, note that the Riemann-Hurwitz formula~\eqref{eq:RH_formula} is still valid for generalized branched coverings. The complete proof of the one-to-one correspondence between equivalence classes of generalized branched coverings and constellations can be found in~\cite[Prop. $2.3$]{CD22}. 

\medskip

Maps can be decomposed by successive deletion of edges via the Tutte decomposition, and constellations admit a similar decomposition. The edges are the objects that carry the labels in a labeled map and this role is played by the corners of color $0$ in constellations. The difference with maps is that, to ensure that the resulting object is still a constellation (so that the process can be iterated), several edges of the constellation must be deleted. 

\begin{definition}[Right-path~\cite{CD22}]
    The right-path of an oriented corner $c$ of a vertex of color $0$ is the sequence of $k+1$ successive edges $e_0,\dotsc,e_{k}$ that follow the orientation of $c$.
\end{definition}

As a consequence of the coloring constraints of vertices of constellations, $e_c$ connects vertices of color $c$ and $c+1$. Therefore deleting a right-path in a constellation of size $n$ results in a constellation of size $n-1$. The right-paths play a role akin to edges in maps and it is possible to give a combinatorial decomposition of a constellation by successive deletion of right-paths. In the case of orientable constellations, such a procedure has been described in~\cite{Fang_PHD} while tracking the degree of the faces of the constellations and the combinatorial decomposition is quite similar to maps. Here, to identify the generating series~\eqref{eq:b_const_KP} with constellations, we need to assign a suitable weight to constellations. This is harmless for the weights associated with face and/or vertices degrees but the weight in $b$ (which we will call $b$-weight in the following) has to be collected while performing the combinatorial decomposition which induces additional difficulties due to the dependency of the $b$-weight on the rooting of the maps as we shall explain in the next paragraph.

\subsubsection{The \texorpdfstring{$b$}{b}-weight and measure of non-orientability}
\label{sssec:b_weight_MON}

The $b$-weight of a constellation can be tracked by successive root-edge deletion following axioms which define a \emph{measure of non-orientability}. 

\begin{definition}[Measure of non-orientability~\cite{CD22}]
\label{def:MON}
    A measure of non-orientability is a function $\rho$ that associates to a map $\mathcal{M}$ and an edge $e \in \mathcal{M}$ a value in $\mathbb{Q}\left[b\right]$ such that $\rho(\mathcal{M},e)$ satisfies the following properties
    \begin{itemize}
        \item {For $\mathcal{N}=\mathcal{M}\setminus \{e\}$ and $c_1,c_2$ the corners of $\mathcal{N}$ such that $\mathcal{M}$ is obtained from $\mathcal{N}$ by inserting an edge between $c_1$ and $c_2$,
        \begin{itemize}
            \item[$(a)$] If $c_1,c_2$ are in different connected component of $\mathcal{N}$ then $\rho(\mathcal{M},e) = 1$.
            \item[$(b)$] If $c_1,c_2$ are in same connected component but different faces, let $\tilde{e}$ the other edge that could be added between $c_1$ and $c_2$ and $\tilde{\mathcal{M}}$ the corresponding map. Then, $\rho(\mathcal{M},e)+\rho(\tilde{\mathcal{M}},\tilde{e})=1+b$.
            \item[$(c)$] If $c_1,c_2$ are in the same face, then $\rho(\mathcal{M},e)=1$ if $\mathcal{N}$ has one face less than $\mathcal{M}$, and $\rho(\mathcal{M},e)=b$ if $\mathcal{N}$ has the same number of faces than $\mathcal{M}$.
        \end{itemize}}
        \item The value of $\rho(\mathcal{M},e)$ depends only on the connected component of $\mathcal{M}$ containing $e$.
    \end{itemize}
    For an ordered list of $p$ edges $(e_1,\dotsc,e_p)$ we denote
    \begin{equation}
        \rho(\mathcal{M},e_1,\dotsc,e_p) = \rho(\mathcal{M},e_1)\dotsc \rho(\mathcal{M}\setminus \{e_1,\dotsc,e_{p-1}\},e_p)
    \end{equation}
\end{definition}

Note that these axioms do not fully specify the measure of non-orientability because of condition $(b)$. Different explicit measures of non-orientability can be obtained by requiring that they satisfy specific properties. However it turns out that the decomposition equation given in the sequel of the manuscript does not depend on the choice of measure of non-orientability, therefore it is sufficient for our purpose to work with the axioms of Definition~\ref{def:MON}. We refer to~\cite[Sec $3.1$]{CD22} for more details on these aspects and examples of different measures of non-orientability.

\begin{definition}[$b$-weight~\cite{CD22}]
For a constellation $\mathcal{C}$ of size $n$ and $n$ disjoint right-paths $R^{(i)}= \{e^{(i)}_1, \dotsc, e^{(i)}_{k}\}$, the $b$-weight of a constellation is defined as
\begin{equation}
    \label{eq:b_weight}
   \rho(\mathcal{C},e^{(1)}_1, \dotsc, e^{(1)}_{k-1},\dotsc , e^{(n)}_1, \dotsc, e^{(n)}_{k})
\end{equation}
\end{definition}

\begin{figure}[!ht]
    \centering
    \includegraphics[scale=0.3]{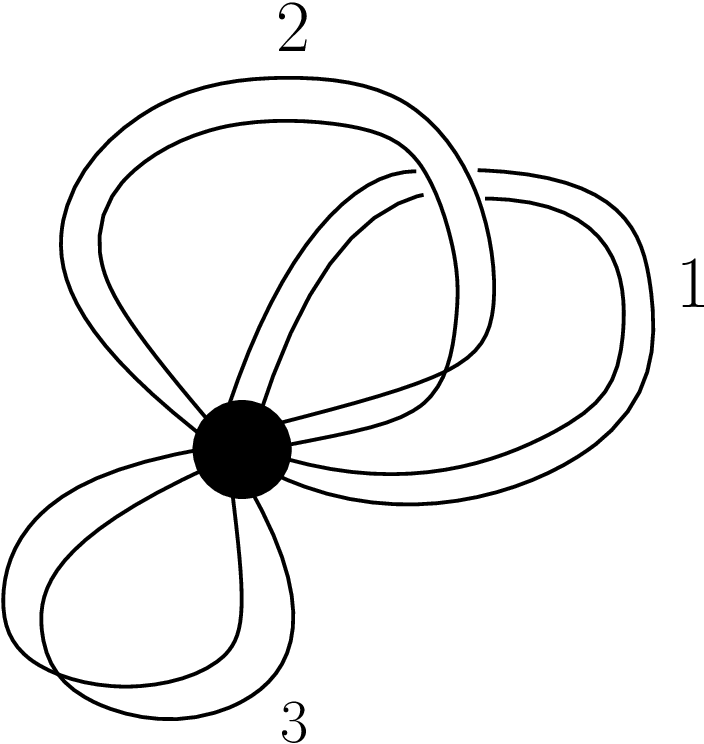}
    \caption{A labeled non-orientable map of genus $\frac{3}{2}$.}
    \label{fig:example_bweight}
\end{figure}

The $b$-weight of a constellation depends on the particular order of edge deletion (and therefore also on the order of right-paths deletion). To illustrate in the simple framework of maps, consider the non-orientable labeled map represented in Figure~\ref{fig:example_bweight}. Deleting edges in order $(3,1,2)$ gives $b$-weight $b^3$ while deleting edges in order $(2,3,1)$ gives a $b$-weight $b^2$ (Deleting edges in any other order also give a $b$-weight of $b^2$ or $b^3$). 

This inductive nature of the $b$-weight makes manipulation of generating series of $b$-deformed constellations more difficult to work with than in the orientable case and renders combinatorial decomposition inaccessible. Instead, most results are obtained by generating larger constellations from smaller ones and exploiting the algebraic properties of the operators involved in the process to identify the corresponding generating series. This is how contact is made between the series~\eqref{eq:b_const_KP} and the generating series of $b$-deformed constellations, which is the main result of~\cite{CD22}.

\begin{theorem}{\cite[Theorem $1.1$]{CD22}}
    The series $(1+b) \frac{t\partial }{\partial t}\ln Z_m^{(b)}(t,\textbf{p},\textbf{q},u_1,\dotsc,u_{k})$ is the generating series of rooted $b$-deformed constellations where $t$ is associated with the size of the constellation, $\textbf{p}$ its face profile, $\textbf{q}$ the profile of vertices of color $0$ and $u_i$ the number of vertex of color $i$ for $i\in \llbracket 1, k-1 \rrbracket $.
\end{theorem}

\subsubsection{Evolution equation for constellations}
\label{sssec:evo_eq_const}

The evolution equation reverses the decomposition process and proceeds by inserting right-path instead. This allows us to control precisely the $b$-weight associated with each edge insertion, thus circumventing the difficulty related to the deletion order. In particular, it avoids the ambiguity related to case $(b)$ of Definition~\ref{def:MON} since both edges will be inserted, leading to a weight $1+b$ independently of the choice of measure of non-orientability. It also makes it possible to get a clear combinatorial interpretation of the terms appearing in this equation. The evolution equation essentially amounts to successive insertion of vertices of color $0$ with arbitrary degree. To add a vertex of color $0$ in the constellation, we create as many right-paths for this new vertex as needed. In the process, we need to keep track of the degree of the face of the right-path we are adding and thus need to introduce catalytic variables $\textbf{y}=(y_i)_{i \geq 0}$ which are replaced by the usual $\textbf{p}$ at the end of the process. We define the following operators:
\begin{itemize}
    \item For $k \geq 2$, we define $Y_k = \sum\limits_{i \geq 0} y_{i+k}\frac{\partial}{\partial y_i}$ which increases the size of the root-face by $k$. The operator $Y_1$ incrementing the size of the root face by one is denoted $Y_+$ instead.
    \item $\Theta_Y = \sum\limits_{i \geq 1} p_i\frac{\partial}{\partial y_i}$ which replaces the catalytic variable with the corresponding $\textbf{p}$ variable.
    \item {$\Lambda_Y$ is the cut-and-join operator. It is given by
    \begin{equation}
        \label{eq:cut_join_b}
        \Lambda_Y = (1+b)\sum\limits_{i,j \geq 1} y_{i+j-1}i\frac{\partial^2}{\partial p_i \partial y_{j-1}} + 
   \sum\limits_{i,j \geq 1} p_iy_{j-1} \frac{\partial}{\partial y_{i+j-1}} + b\sum\limits_{i \geq 1} y_i i \frac{\partial}{\partial y_i}
   \end{equation}
Recall that for colors $c \in \llbracket 1,k-2\rrbracket$ a face of length $2p$ visits $2p$ vertices of color $c$. When fixing an orientation for the face, half of them are visited in increasing order (i.e. the colors of the vertices when following the orientation of the face is $c-1,c,c+1$) and the other half in decreasing order. Three cases can occur when adding an edge in the constellations, which corresponds to the three terms appearing in $\Lambda_Y$.
\begin{itemize}
\item {The first term correspond to adding an edge that connects the root face with a non-root face. If the vertex has color $c$ and is in an increasing (resp. decreasing) order, then it has weight $1$ if it connects to a vertex of color $c+1$ in decreasing (resp. increasing) order and weight $b$ otherwise.}
\item {The second term correspond to adding an edge that decreases the size of the root face. It corresponds to the usual cut case. It can only occur if the two vertices are visited with the same order, and hence has weight $1$.}
\item {The third term corresponds to adding the edge cutting the root face without increasing the root size. It occurs when connecting two vertices with opposite orders. and has weight $b$.}
\end{itemize}}
  
\end{itemize}

We can now give the evolution equation satisfied by the generating function of $b$-deformed rooted connected $k$-constellations. 

\begin{theorem}{~\cite[Theorem $5.10$]{CD22}}
The generating series of $b$-deformed rooted connected $k$-constellations $(1+b) \frac{t\partial}{\partial t} \ln Z^{(b)}_k$ satisfies the evolution equation
    \begin{equation}
    \label{eq:decomp_eq_const}
    (1+b) \frac{t\partial}{\partial t} \ln Z^{(b)}_k = \Theta_Y \sum\limits_{m \geq 1} t^mq_m\left(Y_+\prod\limits_{i=1}^{k} (1+u_i\Lambda_Y)\right)^m \frac{y_0}{1+b} Z^{(b)}_k.
\end{equation}
\end{theorem}

The l.h.s of this equation simply corresponds to the sum of rooted constellations. The r.h.s. encodes how they can be generated from smaller ones by inserting vertices of color $0$ successively. Initially, we start from a root vertex of color $0$ with degree $0$. The operator $\prod\limits_{c=1}^k (u_c+\Lambda_Y)$ creates a right-path for this root vertex by creating new vertices of color $c \geq 1$ and connecting to other elements of $Z^{(b)}_k$. After a right-path is created, we added a corner of color $0$ and a right-path and therefore the total size of the constellation increases by one. This is accounted for by the operator $Y_+$ when re-rooting the constellations at the next corner of the root-vertex. The process is iterated at this new rooted corner as many times as needed. Finally, the catalytic variable associated to the size of the root face is replaced by its corresponding weight.

\section{Constraints for \texorpdfstring{$b$}{b}-deformed constellations}
\label{sec:constraints}

In most map models, the constraints are obtained via a combinatorial decomposition of the maps -typically via cut and join equations \textit{a la Tutte}- and the evolution equation is obtained by summing over the constraints with appropriate weights. The evolution equation then ensures the existence and uniqueness of the partition function since it can be understood as a set of linear relations between the coefficients of the partition function in $t$. When working with the $b$-deformation, despite the similar looking aspect of the equations, its combinatorial interpretation is lost and equations \textit{a la Tutte} cannot be written directly. However, in some situations, it is possible to reverse this construction and to extract constraints from the evolution equation. This can be achieved using the following Lemma, which is a generalization of~\cite[Lemma $2.6$]{BCD22_ON_HCIZ}.

\begin{lemma}
\label{lemma:constr_from_ev}
    Let $Z\in \mathbb{Q}[\textbf{p}][[t]]$ and let $r\geq 1$ be an integer. Assume that 
  \begin{enumerate}[(i).]
  \item\label{enum:Homogeneity} $[t^n]Z$ is a homogeneous polynomial of degree $n$ in the variables $p_i$s.% (where the degree of $p_i$ is $i$).
   \item\label{enum:Evolution} There exist operators $M^{(1)}, \ldots, M^{(r)} \in \mathbb{Q}[[p_1, p_2, \dotsc, \frac{\partial}{\partial p_1}, \frac{\partial}{\partial p_2}, \dotsc]]$ which are independent of $t$ such that
   \begin{equation} \label{EvolutionEquation}
   \frac{\partial Z}{\partial t} = \sum_{m=1}^r t^{m-1} M^{(m)} Z
   \end{equation}
   %\item\label{enum:Initial} for all $i>0$, $[t^0]\frac{\partial}{\partial p_i} Z=0$.
   \item\label{enum:Algebra} For $m \in \llbracket 1,r \rrbracket$, $M^{(m)}$ can be rewritten as $M^{(m)} = \sum\limits_{i \geq 1} p_i M^{(m)}_i$, such that the operators
   \begin{equation}
   L_i \coloneqq i\frac{\partial}{\partial p_i} - \sum_{m=1}^r t^m M^{(m)}_i \ \in \mathbb{Q}[[p_1, p_2, \dotsc, \frac{\partial}{\partial p_1}, \frac{\partial}{\partial p_2}, \dotsc]][t]
   \end{equation}
    for $i\geq 1$ satisfy the following two conditions
    	\begin{itemize}
    	\item $[t^n]L_i F$ is a homogeneous polynomial of degree $n-i$ in the variables $p_i$s,
    	\item the operators $L_i$ are closed under commutation in the sense that there exist $D_{ij}^k\in \mathbb{Q}[[p_1, p_2, \dotsc, \frac{\partial}{\partial p_1}, \frac{\partial}{\partial p_2}, \dotsc]][t]$, for $i, j, k\geq 1$, such that
   		\begin{equation} \label{CommutationRelation}
   		\left[L_i,L_j\right] = t \sum_{k \geq 1} D_{ij}^k L_k
   		\end{equation}
   		\end{itemize}
   %\item $[t^n] L_i Z$ is a homogeneous polynomial in the $p_j$s of degree $n-i$,
 \end{enumerate}
    then $Z$ satisfies the constraints
    \begin{equation}
        \label{eq:constraints}
        L_i Z = 0 \qquad \text{for all $i\geq1$.}
    \end{equation}
\end{lemma}

Notice that the factor $t$ in the RHS of \eqref{CommutationRelation} is crucial. Without it, the commutation relation $\left[L_i,L_j\right] = \sum_{k \geq 1} D_{ij}^k L_k$ has the simple solution $D_{ij}^k = \delta_{j,k} L_i - \delta_{i,k} L_j$, which is obviously not of interest for us. As we will see, the difficulty in applying the above lemma is identifying the $M_i^{(m)}$s and actually proving the existence of the commutation relation \eqref{CommutationRelation}.

As another technical remark, notice that in general there exist several ways\footnote{For instance, for some fixed $k, l$, one can always add a term $p_k$ to $M^{(m)}_l$ and $- p_l$ to $M^{(m)}_k$ without changing $M^{(m)}$.} of writing $M^{(m)}=\sum_j p_j M^{(m)}_j$ but one looks for a set of $M^{(m)}_j$s such that condition \ref{enum:Algebra} (homogeneity of $[t^n]L_iZ$ and commutation relation \eqref{CommutationRelation}) holds which is a very much non-trivial requirement. 

\begin{proof}
The proof is based on the proof of Lemma $2.6$ in~\cite{BCD22_ON_HCIZ} which we adapt to our conditions. We proceed by induction on the coefficients of the the formal power series $Z$ in $s$. Since $[t^n]Z$ is homogeneous of degree $n$, $[t^0]Z$ is independent of the $p_i$s, hence $[t^0]L_i Z = 0$ (in fact the assumption that $[t^n]L_i F$ is homogeneous of degree $n-i$ gives that it vanishes for $n<i$). Let $n>0$ and assume that for all $n'<n$ we have $\left[t^{n'}\right]L_i Z = 0$. Since $Z$ satisfies the evolution equation, we have
\begin{equation}
    \left[ t\frac{\partial}{\partial t} - \sum_{m=1}^r t^m M^{(m)}, L_i\right] Z = \left( t\frac{\partial}{\partial t} - \sum_{m=1}^r t^m M^{(m)}\right) L_i Z.
\end{equation}
We extract the coefficient of $t^n$ and use $[t^n] t\frac{\partial F}{\partial t} = n[t^n] F$ for any formal power series $F$, as well as the induction hypothesis $[t^{n-m}]L_i Z =0$ for $m=1, \dotsc, \min(n,r)$, so that
\begin{equation} \label{eq:sfwd_exp}
[t^n]\left[ t\frac{\partial}{\partial t} - \sum_{m=1}^r t^m M^{(m)}, L_i\right] Z = n[t^n] L_iZ.
\end{equation}

Let us calculate the commutator on the LHS independently. We start with $[t\frac{\partial}{\partial t}, L_i] Z$ which is evaluated by noting that
\begin{itemize}
\item since $[t^n]Z$ is homogeneous of degree $n$, $t\frac{\partial Z}{\partial t} = \sum_{j\geq 1} jp_j \frac{\partial Z}{\partial p_j}$.
\item since $[t^n]L_i Z$ is homogeneous of degree $n-i$, $t\frac{\partial L_i Z}{\partial t} = \sum_{j\geq 1} jp_j \frac{\partial L_i Z}{\partial p_j} + iL_i Z$.
\end{itemize}
Hence,
\begin{equation}
\left[t\frac{\partial}{\partial t}, L_i\right] Z = iL_i Z + \sum_{j\geq 1} j\left[p_j\frac{\partial}{\partial p_j}, L_i\right]Z
\end{equation}
Writing $M^{(m)} = \sum_{j\geq 1} p_j M^{(m)}_j$,
\begin{equation}
\begin{aligned}
\left[ t\frac{\partial}{\partial t} - \sum_{m=1}^r t^m M^{(m)}, L_i\right] Z &= iL_i Z + \sum_{j\geq 1} \left[p_j L_j, L_i\right]Z \\
& = iL_i Z + \sum_{j\geq 1} p_j \left[ L_j, L_i\right]Z 
+ \sum_{j\geq 1} \left[ p_j, L_i\right] L_jZ
\end{aligned}
\end{equation}
The commutator $\left[ L_j, L_i\right]$ is given by the lemma assumption. As for the other commutator,
\begin{equation}
\left[ p_j, L_i\right] = -i\delta_{ij} - \sum_{m=1}^r t^m B^{(m)}_{ji}
\end{equation}
where $B^{(m)}_{ji} \coloneqq [p_j, M^{(m)}_i]$ is independent of $t$. Overall one finds
\begin{equation}
\left[ t\frac{\partial}{\partial t} - \sum_{m=1}^r t^m M^{(m)}, L_i\right] Z = t\sum_{j,k\geq 1} p_j D_{ji}^k L_k Z - \sum_{m=1}^r t^m \sum_{j\geq 1} B^{(m)}_{ji} L_jZ
\end{equation}
We now extract the coefficient of $t^n$ to get
\begin{equation}
[t^n]\left[ t\frac{\partial}{\partial t} - \sum_{m=1}^r t^m M^{(m)}, L_i\right] Z = \sum_{j,k\geq 1} p_j [t^{n-1}] D_{ji}^k L_k Z - \sum_{j\geq 1} \sum_{m=1}^r B^{(m)}_{ji} [t^{n-m}] L_jZ.
\end{equation}
Our induction hypothesis implies that $[t^{n-m}]L_k Z = 0$ for all $m= 1, \ldots, n$, and therefore the RHS is found to be zero. Equating with the RHS of \ref{eq:sfwd_exp}, we find $[t^n] L_iZ=0$.
\end{proof}

\medskip

Before diving into particular models, let us introduce some operators to describe how we can make use of this lemma. We define the evolution operator $M^{(k)}$ as
\begin{align}
\label{eq:def_all_modes}
    M^{(k)} &= \Theta_Y \sum\limits_{m \geq 1} t^m q_m\left(Y_+\prod\limits_{c=1}^{k} (u_c+\Lambda_Y)\right)^m \frac{y_0}{1+b} \\
    \label{eq:def_modes}
     M^{(k)}_i &= \left[y_i\right] \sum\limits_{m \geq 1} t^m q_m\left(Y_+\prod\limits_{c=1}^{k} (u_c+\Lambda_Y)\right)^m \frac{y_0}{1+b}
\end{align}
such that $M^{(k)}= \sum\limits_{i \geq 1} p_i  M^{(k)}_i $. Hence the operators $M^{(k)}_i$ are the candidate operators to extract constraints in constellations. However, Lemma~\ref{lemma:constr_from_ev} cannot be applied immediately. Indeed~\eqref{eq:decomp_eq_const} doesn't satisfy condition~\ref{EvolutionEquation}. Thus, to make use of the Lemma, we shall specialize $\textbf{q}$ such that $q_{m} = 0$ when $m>p$ for some $p\geq 1$. Combinatorially, this specialization comes to restricting the degree of the vertices of color $0$ of the constellation. We denote $Z^{(b)}_{k,p}$ the corresponding generating function. Under these additional conditions, the evolution operators write
\begin{align}
\label{eq:def_mode_spec}
    M^{(k,m)}_i &= \left[y_i\right]\left(Y_+\prod\limits_{c=1}^{k} (u_c+\Lambda_Y)\right)^m \frac{y_0}{1+b}
\end{align}

We introduce the $b$-deformed currents $J^{(b)}_i$. They are defined for $i \neq 0$ as
\begin{equation}
    \label{eq:currents}
    J^{(b)}_i = \begin{cases} p_i \qquad &\text{for} \hspace{2pt}i<0 \\ (1+b)i\frac{\partial}{\partial p_i} \qquad &\text{for}\hspace{2pt}i>0 \end{cases},
\end{equation}
and satisfy the commutation relation
\begin{equation}
    \left[J^{(b)}_i,J^{(b)}_j\right] = (1+b)i\delta_{i,-j}.
\end{equation}

With these notations, we have
\begin{align}
\label{eq:def_lambda_J}
\left(u_j + \Lambda_Y \right)&= \sum\limits_{i \in \mathbb{Z}\setminus\{0\}} Y_i J^{(b)}_i + b \sum\limits_{i \geq 0} iy_i \frac{\partial}{\partial y_i}
\end{align}

Now, to apply Lemma~\ref{lemma:constr_from_ev}, it remains to check that the modes are closed under commutation. Since we have an explicit expression for the modes, this can be done through a direct computation of their commutators. The complexity of this computation grows with the number of times the operator $\left(u^{-1}+\Lambda_Y\right)$ is applied. We will say that a constellation model is of order $n$ if its modes involve $n$-th power of the operator $\left(u^{-1}+\Lambda_Y\right)$. Hence constellations with $k$-colors where vertices of color $0$ have degree $m$ (whose modes are given by $M^{(k,m)}_i$) are of order $km$. There are two quadratic constellation models which correspond to maps ($k=1,m=2$) and bipartite maps ($k=2,m=1$). The constraints for these models are known to be the Virasoro constraints~\cite{AvM01}. Our goal is to derive constraints for the case of cubic constellations models, corresponding to the two cases $k=3,m=1$ and $k=1,m=3$. Therefore, in the following, we will specialize to evolution operators with either $k=1$ or $m=1$. While the computations are similar in both families, the expansion of its evolution operators is different and thus they should be treated separately.
\begin{itemize}
\item For evolution operators with $m=1$, we will expand the evolution operators over elementary symmetric polynomials.
\item For evolution operators with $k=1$, we have a single variable $u_1=u$. Thus, we can introduce another current operator $J_0 = u$ which leaves the commutation relation unchanged and allow to recast $\left(u_j + \Lambda_Y \right)$ as a sum of currents. 
\end{itemize}

\subsection{Expansion of evolution operators with $m=1$}

The operators 
\begin{equation}
M^{(k,m)}_i = [y_i] \left(Y_+ \prod_{c=1}^k (\Lambda_Y+u_c)\right)^m \frac{y_0}{1+b}
\end{equation}
are symmetric polynomials in $u_1, \ldots, u_k$. They can be expanded over the elementary symmetric polynomials $e_{p}(u_1,\ldots,u_{k})$ for $0 \leq p \leq k-1$ to organize our calculation. For $i\geq m$ and integers $s_1, \dotsc, s_m\geq 0$, let
\begin{equation}
A^m_i(s_1, \ldots, s_m) \coloneqq [y_i] \prod_{l=1}^m Y_+\Lambda_Y^{s_l}\ \frac{y_0}{1+b},
\end{equation}
so that
\begin{equation}
M^{(k,m)}_i = \sum_{s_1, \ldots, s_m =0}^k A^m_i(s_1, \ldots, s_m) \prod_{l=1}^m e_{k-s_l}(u_1, \ldots, u_k).
\end{equation}

When $m=1$, there is a single contribution to each elementary symmetric function. As we will see shortly after, this allows to give simple conditions for the commutators of evolution operators $M^{(k,m)}_i$ -and in turn for the candidate constraints- to close. We now focus on the following case, relevant for $k$-constellations with $q_m = \delta_{m,1}$. 
\begin{equation}
A_i(s)\equiv A_i^1(s) = [y_i]Y_+\Lambda_Y^s \frac{y_0}{1+b},
\end{equation}
for $i\geq 1, s\geq 0$.
We form the candidate constraints
\begin{equation}
L^{(k)}_i \equiv L^{(k,1)}_{i|q_m = \delta_{m,1}} = -p_i^* + t M^{(k,1)}_i = -p_i^* + t \sum_{p=0}^k e_{k-p}(u_1, \dotsc, u_k) A_i(p),
\end{equation}
for which the following proposition provides a recursion.
\begin{proposition}
\label{prop:mod_rec_no_ui}
One has $A_i(0) = \frac{\delta_{i,1}}{1+b}$ for all $i\geq 1$. Then for all $i\geq1, p\geq 0$
\begin{equation} \label{eq:mod_rec_no_ui}
    A_i(s+1) = \sum\limits_{\substack{n \geq 1 }} J^{(b)}_{i-n} A_n(s) + b(i-1) A_i(s) 
\end{equation}
\end{proposition}

\begin{proof}
By definition, $A_i(0) = [y_i] Y_+\frac{y_0}{1+b} = \frac{\delta_{i,1}}{1+b}$ for all $i\geq 0$. %(correct in particular for $i=1$ since $J_0 = \alpha=0$, but not correct for $i=0$ since then $A_0(1) = 0$). 
Then
\begin{equation}
    A_i(s+1) = \left[y_i\right] Y_+ \Lambda_Y^{s+1} \frac{y_0}{1+b} = \left[y_{i-1}\right] \Lambda_Y \sum\limits_{j \geq 1} y_{j-1} \underbrace{\left[y_{j-1}\right]  \Lambda_Y^{s} \frac{y_0}{1+b}}_{A_j(s)}
\end{equation}
and using~\eqref{eq:def_lambda_J} gives~\eqref{eq:mod_rec_no_ui}.
\end{proof}

Now, the following Lemma allows to reconstruct commutator of the constraints from those o operators $A_i(s)$.
\begin{lemma} \label{lemma:SimplifiedCommutators}
Assume that there exist differential operators $(D_{ij, l}(s))_{i,j,l\geq 1, s\geq 0}$ in the variables $p_i$s such that
\begin{equation} \label{Commutatorp*A}
[p_i^*, A_j(s)] - [p_j^*, A_i(s)] = \sum_{l\geq 1} D_{ij, l}(s) p_l^*
\end{equation}
and
\begin{equation} \label{eq:SimplifiedCommutators}
\begin{aligned}
    \left[A_i(s), A_j(s)\right] &= \sum_{l\geq 1} D_{ij,\ell}(s) A_l(s), \\
    \left[A_i(s), A_j(s')\right] -  \left[A_j(s), A_i(s')\right] &= \sum_{l\geq 1} D_{ij,l}(s') A_l(s) + D_{ij,l}(s) A_l(s').
\end{aligned}
\end{equation}
Then the operators $L^{(k)}_i = -p_i^* + t \sum_{s=0}^k e_{k-s}(u_1, \dotsc, u_k) A_i(s)$ satisfy
\begin{equation}
\label{eq:CommutConstr}
\left[L^{(k)}_i, L^{(k)}_j\right] = t \sum_{l\geq 1} D_{ij, l}^{(k)} L_l^{(k)},
\end{equation}
with
\begin{equation}
D_{ij, l}^{(k)} = \sum_{s=0}^k e_{k-s}(u_1,\ldots,u_k) D_{ij, l}(s).
\end{equation}
\end{lemma}

\begin{proof}
Recall that $M^{(k,1)}_i = \sum_{p=0}^k e_{k-p}(u_1,\ldots,u_k) A_i(p)$, which combined with~\eqref{eq:SimplifiedCommutators}, gives directly
\begin{equation}
\left[M^{(k,1)}_i, M^{(k,1)}_j\right] = \sum_{l\geq 1} D_{ij, l}^{(k)} M_l^{(k)}.
\end{equation}
In addition, the constraints read $L^{(k)}_i = -p_i^* + tM^{(k,1)}_i$, hence
\begin{equation}
[L^{(k)}_i, L^{(k)}_j] = -t\left( \left[p_i^*,M^{(k,1)}_j \right] + \left[M^{(k,1)}_i, p_j^* \right] \right) 
+ t^2\left[M^{(k,1)}_i, M^{(k,1)}_j\right].
\end{equation}
We use \eqref{Commutatorp*A} and get directly~\eqref{eq:CommutConstr}.
\end{proof}

As it turns out, for any set of operators $A_i(s)$, the equations \eqref{eq:SimplifiedCommutators} always admit the solution $D_{ij, l}(s) = \delta_{j,l} A_i(s) - \delta_{i,l} A_j(s)$, and therefore any operators of the form $M^{(k,1)}_i = \sum_{s=0}^k e_{k-s}(u_1,\ldots,u_k) A_i(s)$ satisfy $[M^{(k,1)}_i, M^{(k,1)}_j] = \sum_{l\geq 1} D_{ij, l}^{(k)} M_l^{(k)}$. The non-triviality of Lemma \ref{lemma:SimplifiedCommutators} comes the fact that the constraints read $L^{(k)}_i = -p_i^* + tM^{(k,1)}_i$, and that we assume that the coefficients $D_{ij, l}(s)$ satisfy Equation \eqref{Commutatorp*A}. For the same reason, the assumption that $[L_i, L_j] = t \sum_{k\geq 1} D_{ij,k} L_k$ in the Lemma \ref{lemma:constr_from_ev} would be trivial without the factor $t$ on the RHS, since $D_{ij,k}$ could be a combination of the $L_i$s themselves.

Therefore for those models, our strategy is to prove the relations \eqref{eq:SimplifiedCommutators} from which the commutators~\eqref{eq:CommutConstr} and condition~\ref{enum:Algebra} follow.

\subsection{Expansion of evolution operator with $k=1$}

For these evolution operators model, a simplification appears compared to the previous section: since there is a single variable $u$, there is in fact no need to expand $L^{\text{bip}\leq 3}_i$ on elementary symmetric functions of $u_1, \dotsc, u_k$. Instead, it is enough to introduce another current $J^{(b)}_0  = u$ in \eqref{eq:currents}. This leaves the commutation relations unchanged since  $J^{(b)}_0$ commutes with all other currents. This trick allows to bypass the computation of commutators $A^m_i(s_1,\dotsc,s_m)$ to work directly with the evolution operators $M^{(1,m)}_i$. 

\medskip
For those models, we can give a recursion relation on operators  $M^{(1,m)}_i$ as follows.

\begin{proposition}
One has $M^{(1,1)}_i = \frac{J^{(b)}_{i-1}}{1+b}$ for all $i\geq 1$. Then for all $i\geq1, m\geq 1$
\begin{equation} \label{eq:mod_rec_no_ui'}
    M^{(1,m+1)}_i = \sum\limits_{n \geq 1} J^{(b)}_{i-n-1} M^{(1,m)}_n + b(i-1) M^{(1,m)}_{i-1}.
\end{equation}
(using $M^{(1,m)}_0 = 0$ for all $m\geq 1$ when necessary).
\end{proposition}

\begin{proof}
Let $Y_i$, $i\in\mathbb{Z}$, be the operator defined as $Y_i y_j = y_{i+j} \delta_{i+j\geq 0}$. Then,
\begin{equation} \label{Lambda+u}
\Lambda_Y+u = \sum_{i\in\mathbb{Z}} Y_i J^{(b)}_i + b \sum\limits_{i \geq 0} iy_i [y_i].
\end{equation}
By definition, $M^{(1,m)}_i = [y_i] (Y_+(\Lambda_y+u))^m \frac{y_0}{1+b}$ for all $i, m\geq 0$, which gives directly $M^{(1,1)}_i = \frac{J^{(b)}_{i-1}}{1+b}$ for all $i\geq 0$, and $M^{(1,m)}_0 = \frac{\delta_{m,0}}{1+b}$ for all $m\geq 0$. Then for $m\geq 1$
\begin{equation}
M^{(1,m+1)}_i = \sum_{n\geq 0} [y_i] Y_+(\Lambda_Y+u) y_n [y_n] (Y_+(\Lambda_y+u))^m \frac{y_0}{1+b} = \sum_{n\geq m} [y_i] Y_+(\Lambda_Y+u) y_n M^{(1,m)}_n
\end{equation}
and we conclude with \eqref{Lambda+u}.
\end{proof}

The structure of the calculations is completely similar to what was done in the previous section. The analog of Lemma \ref{lemma:SimplifiedCommutators} is as follows.
\begin{lemma} \label{lemma:SimplifiedCommutators2}
Assume that there exist differential operators $(D_{ij, l}^{(m)})_{i,j,l\geq 1, p\geq 0}$ in the variables $p_i$s such that
\begin{equation} \label{Commutatorp*M}
[p_i^*, M^{(1,m)}_j] - [p_j^*, M^{(1,m)}_j] = \sum_{l\geq 1} \tilde{D}_{ij, l}(m) p_l^*
\end{equation}
and
\begin{equation} \label{eq:SimplifiedCommutators2}
\begin{aligned}
    \left[M^{(1,m)}_i, M^{(1,m)}_j\right] &= \sum_{l\geq 1} \tilde{D}_{ij, l}(m) M^{(1,m)}_l, \\
    \left[M^{(1,m)}_i, M^{(1,m')}_j\right] -  \left[M^{(1,m)}_j, M^{(1,m')}_i\right] &= \sum_{l\geq 1} \tilde{D}_{ij, l}(m') M^{(1,m)}_l + \tilde{D}_{ij, l}(m) M^{(1,m')}_l.
\end{aligned}
\end{equation}
Then the operators $L^{\text{bip}\leq p}_i = -p_i^* + \sum_{m=1}^p q_m t^m M^{(1,m)}_i$ satisfy
\begin{equation}
\label{eq:CommutConstr2}
\left[L^{\text{bip}\leq p}_i, L^{\text{bip}\leq p}_j\right] = t \sum_{l\geq 1} \tilde{D}_{ij, l}^{(p)} L_l^{\text{bip}\leq p},
\end{equation}
with
\begin{equation}
\tilde{D}_{ij, l}^{(p)} = \sum_{m=1}^r q_m t^{m-1} \tilde{D}_{ij, l}(m).
\end{equation}
\end{lemma}

\begin{remark}
The operators $M^{(k,m)}$ could also be expanded over elementary symmetric polynomials as in the previous case, yielding similar computations than in the previous case. When $m > 1$, there are several contributions to each elementary symmetric polynomial. This leads to weaker conditions for the commutators of operators $M^{(k,m)}$ and increases the number of commutators of operators $A^m(s_1,\dotsc,s_m)$ that must be computed. The introduction of a current $J_0=u$ allows to bypass these difficulty to work directly at the level of evolution operators $M^{(k,m)}$.
\end{remark}

\subsection{Quadratic constellations}
\label{ssec:quad_const}

There are two different quadratic constellation models. The first one is $m=2,k=2$ and corresponds to face-weighted $b$-deformed maps. The second one is given by $m=3,k=1$, the vertices of color $0$ are leaves and can be bijectively deleted leading to $b$-deformed bipartite maps. The constraints for $b$-deformed maps were known to form a Virasoro algebra via~\cite{AvM01} through the connection between the $b$-deformation and the $\beta$-ensembles of random matrix theory. The derivation we present here leads to the same result but it is more direct since it relies only on the generating series and its evolution equation. The computation of the commutator is very similar in both models and leads to the same algebra. We will only present here the computation of the commutator for the bipartite maps. 
\medskip

To compute the commutator $M^{(2,1)}$, we have to compute all commutators of operators $A_i(s)$ with $s\leq 2$. We compute them by increasing $s$ and using Proposition~\ref{prop:mod_rec_no_ui}.

\subsubsection{Computation of $\left[ A_i(s), A_j(s') \right]$ for $s, s' =0, 1$}
We have $A_i(0) = \frac{\delta_{i,1}}{1+b}$ and $A_i(1) = \frac{J^{(b)}_{i-1}}{1+b}$ for $i\geq 1$. Therefore
\begin{align}
[A_i(0), A_j(0)] =[p_i^*, A_j(0)] &=  0, \\
 [A_i(1), A_j(1)] =  [A_i(0), A_j(1)]  &= [p_i^*, A_j(1)] = 0.
\end{align}

For later purposes, we also give for $r\in\mathbb{Z}$
\begin{equation}
    \left[ J_j^{(b)}, A_i(1) \right] = j \delta_{i+j,1}.
\end{equation}

\subsubsection{Computation of $\left[A_i(s), A_j(s') \right]$ for $s, s' \leq 2$} \label{sec:A2A2}
We use the relation~\eqref{eq:mod_rec_no_ui} to compute the various commutators involving $A_i(2)$, starting with
\begin{equation}
\begin{aligned}
    \left[J_j^{(b)} , A_i(2) \right] &= \sum_{n\geq 1} J_{i-n}^{(b)} [J^{(b)}_j, A_n(1)] + [J^{(b)}_j, J^{(b)}_{i-n}]A_n(1) + b(i-1)[J^{(b)}_j, A_i(1)]\\
    &= jJ^{(b)}_{i+j-1}\left( \delta_{j \leq 0} + \delta_{i+j \geq 1} \right) - b(i-1)^2 \delta_{i+j,1}.
\end{aligned}
\end{equation}
A direct application gives the commutator $\left[A_i(1) , A_j(2) \right]=(i-1) A_{i+j-1}(1)$, and $\left[p_i^*, A_j(2) \right]=i p_{i+j-1}^*$. Therefore
\begin{equation}
\begin{aligned}
\left[A_i(1) , A_j(2) \right] - \left[A_j(1) , A_i(2) \right] &= (i-j) A_{i+j-1}(1),\\
\left[p^*_i , A_j(2) \right] - \left[p^*_j , A_i(2) \right] &= (i-j) p^*_{i+j-1}.
\end{aligned}
\end{equation}

Moreover, %after a lengthier but classical calculation, one finds the half-Virasoro algebra
\begin{equation}
\begin{aligned}
    \left[A_i(2) , A_j(2) \right] &= \left[\sum_{n\geq 1} J^{(b)}_{i-n} A_n(1) + b(i-1) A_i(1), A_j(2)\right]\\
    &\begin{multlined} = \biggl(\sum_{n\geq j} (n-j) + \sum_{n\geq i} (i-n) + \sum_{n=1}^{i+j-1} (i-n)\biggr) J^{(b)}_{i+j-n-1} A_n(1) \\ + b(i-j)(i+j-2) A_{i+j-1}(1)\end{multlined}
\end{aligned}
\end{equation}
On one hand,
\begin{multline}
\biggl(\sum_{n\geq j} (n-j) + \sum_{n\geq i} (i-n)\biggr) J^{(b)}_{i+j-n-1} A_n(1) \\ 
= \biggl(\sum_{n\geq M} (i-j) + \sgn(i-j) \sum_{n=\mu}^{M-1} (n-\mu)\biggr) J^{(b)}_{i+j-n-1} A_n(1).
\end{multline}
In the second sum, we notice that the sum can be started at $n=\mu+1$, we write $A_n(1) = J^{(b)}_{i-1}/(1+b)$, perform the change of summation index $n\to i+j-n$ and commute the two $J^{(b)}$s, so that
\begin{equation}
\begin{aligned}
\sum_{n=\mu+1}^{M-1} (n-\mu)J^{(b)}_{i+j-n-1} A_n(1) &= \sum_{n=\mu+1}^{M-1} (M-n) J^{(b)}_{i+j-n-1} A_n(1)\\
&= \frac{M-n}{2} \sum_{n=\mu+1}^{M-1} J^{(b)}_{i+j-n-1} A_n(1)
\end{aligned}
\end{equation}
the last being obtained by taking the half-sum of both expressions. This gives
\begin{multline}
\biggl(\sum_{n\geq j} (n-j) + \sum_{n\geq i} (i-n)\biggr) J^{(b)}_{i+j-n-1} A_n(1) = (i-j) \biggl(\sum_{n\geq M}  + \frac{1}{2} \sum_{n=\mu+1}^{M-1} \biggr) J^{(b)}_{i+j-n-1} A_n(1).
\end{multline}
On the other hand, using the same tricks as above, one finds
\begin{equation}
\sum_{n=1}^{i+j-1} (i-n) J^{(b)}_{i+j-n-1} A_n(1) = \frac{i-j}{2} \sum_{n=1}^{i+j-1} J^{(b)}_{i+j-n-1} A_n(1)
\end{equation}
which is decomposed as
\begin{equation}
\frac{i-j}{2} \sum_{n=1}^{i+j-1} J^{(b)}_{i+j-n-1} A_n(1) = \frac{i-j}{2} \biggl(\sum_{n=1}^{\mu} + \sum_{n=\mu+1}^{M-1} + \sum_{n=M}^{i+j-1}\biggr) J^{(b)}_{i+j-n-1} A_n(1).
\end{equation}
The first and third sum can be shown to be equal, by using $A_n(1) = J^{(b)}_{n-1}/(1+b)$, performing the change $n\to i+j-n$ and commuting the two $J^{(b)}$s. All in all, it comes that
\begin{equation}
\left[A_i(2) , A_j(2) \right] = (i-j) A_{i+j-1}(2).
\end{equation}

\medskip

The constraints of bipartite maps are given by
\begin{equation}
\begin{aligned}
L^{\text{bip}}_i &= -p_i^* + t(A_i(2) + e_1(u_1, u_2) A_i(1) + e_2(u_1, u_2))
\end{aligned}
\end{equation}

Now that we have computed the commutators for all operators $A_i(s)$, we can obtain the commutator of constraints immediately via Lemma~\ref{lemma:SimplifiedCommutators}.
\begin{equation}
[L^{\text{bip}}_i, L^{\text{bip}}_j] = t\sum_{l\geq 1} D_{ij, \ell}^{(2)} L^{\text{bip}}_\ell
\end{equation}
with $D_{ij, \ell}^{(2)} = (i-j)\delta_{\ell=i+j-1}$.
Therefore the candidate constraints  $L^{\text{bip}}_i$ are closed form a (shifted, half) Virasoro algebra. All the conditions to apply Lemma~\ref{lemma:constr_from_ev} are now satisfied.
 This leads to the following theorem.
\begin{theorem}
For all $b>-1$, the generating function $Z^{(b)}_{2,1}\left[t,q,u_1,u_2,\textbf{p}\right]$ satisfies constraints $L^{\text{bip}}_i \cdot Z = 0$ for all $i\geq 0$ where
\begin{equation}
\begin{aligned}
L^{\text{bip}}_i &= -p_i^* + t(A_i(2) + e_1(u_1, u_2) A_i(1) + e_2(u_1, u_2))
\end{aligned}
\end{equation}
\end{theorem}

In particular at $b=0$ (resp. $b=1$), they are nothing but the usual Virasoro constraints for orientable (resp. non-oriented) bipartite maps. However, this combinatorial interpretation of the constraints as a cut-and-join operation is valid only for these values of $b$ due to the dependency of the $b$-weight on the order of deletion of the edges.

\medskip

For the sake of completeness, we state the analogous result for the generating function of $b$-deformed maps $Z^{(b)}_{1,2}$.  Its evolution operator $M^{(1,2)}_i $ is given by
\begin{equation}
    M^{(1,2)}_i(u) = \sum\limits_{k \geq 1} J^{(b)}_{i-1-k}J^{(b)}_{k} + b(i-1)J^{(b)}_{i-2} + u J^{(b)}_{i-1}.
\end{equation}

The computation of the commutation relations for their modes can be found in Appendix~\ref{sec:comm_triv_bip}. They also form a shifted Virasoro algebra 
\begin{equation}
    \label{eq:Vir_22}
    \left[ L^{\text{map}}_i, L^{\text{map}}_j \right] = (i-j)L^{\text{map}}_{i+j-2}
\end{equation}
which leads to the following constraints.
\begin{theorem}[\cite{AvM01}]
\label{thm:constr_bmaps}
For all $b>-1$, the generating function $Z^{(b)}_{1,2}\left[t,q,u,\textbf{p}\right]$ satisfies constraints $L^{\text{map}}_i \cdot Z = 0$, for all $i\geq 0$, where
\begin{equation}
   L^{\text{map}}_i =  \frac{1}{t^2}\frac{i\partial}{\partial p_i}-M^{(1,2)}_i .
\end{equation}
\end{theorem}

\subsection{Cubic constellations}
\label{ssec:cubic_const}

We now turn to the case of cubic constellations. Again, there are two different cubic constellation model. The first one is given by $k=1,m=3$ and corresponds combinatorially to bipartite maps where one set of vertices have degree $3$ (which we will call trivalent bipartite maps hereafter). The second one is given by $k=3, m=1$ and corresponds to $3$-constellations. 

\begin{figure}[!ht]
\hfill
\begin{subfigure}{.4\textwidth}
  \centering
  \includegraphics[scale=0.45]{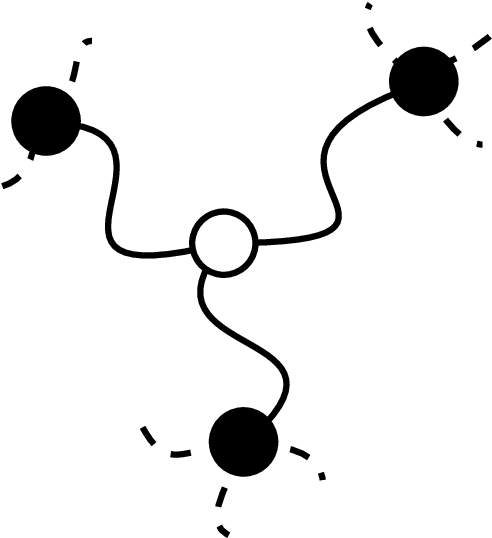}
  \label{fig:trivalent_bip_map}
\end{subfigure}
\hfill
\begin{subfigure}{.4\textwidth}
  \centering
  \includegraphics[scale=0.7]{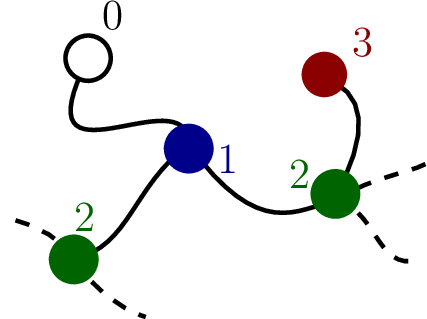}
  \label{fig:vert_3_const}
\end{subfigure}
\caption{Subgraphs of two different models of cubic constellation.}
\label{fig:cubic_const}
\end{figure} 

To derive the set of constraints for these models, we need to apply Lemma~\ref{lemma:constr_from_ev} as in the quadratic case. While the connection between maps and the Virasoro algebra was known to be valid for all values of $b$~\cite{AvM01}, there are no analogous results on constellations. Therefore, it is not obvious that the modes of these models should form an algebra. Yet the constraints satisfied by orientable constellations were known from~\cite{Fang_PHD} where they were obtained from a purely combinatorial decomposition. These constraints correspond to the one we would obtain by applying the Lemma~\ref{lemma:constr_from_ev} to our generating functions at $b=0$ which motivated this computation.

\medskip

The computation of the commutator of the modes for both models is lengthy and not particularly enlightening hence will only say a few words on the method behind the computation here. The detailed computation for each model can be found in Appendix~\ref{app:cubic_const}. The strategy is similar to the quadratic case. For each model, we expand the commutator of cubic constellations via~ Proposition\ref{lemma:SimplifiedCommutators} when $p=1$ or ~Proposition~\ref{lemma:SimplifiedCommutators2} when $k=1$. Each contribution can then be computed from the commutator of quadratic constellations. While direct, the number of terms appearing makes identification of cancellations across various terms difficult. It requires rewriting of each of them by further splitting terms, relabeling sums and using commutation relations. At the end of the day, we show that the candidate constraints are closed under commutation in both cases, as expected. In the following, we denote $\mu = \min (i,j)$ and $M = \max (i,j)$.

\subsubsection{\texorpdfstring{$3$}{3}-constellations}

For $3$-constellations, we obtain the following proposition.

\begin{proposition} \label{thm:SimplifiedCommutators}
Let $\mu = \min(i,j)$ and $M = \max(i,j)$. For $0\leq s, s'\leq 3$,
\begin{subequations} \label{SimplifiedCommutators}
\begin{alignat}{1}
    \left[A_i(s), A_j(s)\right] &= \sum_{l\geq 1} D_{ij, l}(s) A_l(s),\label{CommutatorA} \\
    \left[A_i(s), A_j(s')\right] -  \left[A_j(s), A_i(s')\right] &= \sum_{l\geq 1} D_{ij,l}(s') A_l(s) + D_{ij,l}(s) A_l(s'), \label{CommutatorAA}\\
    \left[p_i^*, A_j(s)\right] -  \left[p_j^*, A_i(s)\right] &= \sum_{l\geq 1} D_{ij,l}(s) p_l^* \label{CommutatorpA}
\end{alignat}
\end{subequations}
with $D_{ij, l}(0) = D_{ij, l}(1) = 0$ and $D_{ij, l}(2) = (i-j) \delta_{l,i+j-1}$ and
\begin{multline}
D_{ij, l}(3) = (i-j) (2\delta_{l\geq M} + \delta_{M\leq l\leq i+j-1}) J^{(b)}_{i+j-1-l} \\
+ \operatorname{sgn}(i-j)(2l-3\mu+1) \delta_{\mu\leq l\leq M-1} J^{(b)}_{i+j-1-l}
+ b(i-j)(i+j-2)\delta_{l, i+j-1}.
\end{multline}
\end{proposition}

This proposition shows that the assumptions of the Lemma \ref{lemma:SimplifiedCommutators} are satisfied by $3$-constellations. From Lemma~\ref{lemma:SimplifiedCommutators} we then get the commutator for candidate constraints
\begin{equation}
L^{\text{3-const.}}_i = -p_i^* + t(A_i(3) + e_1(u_1, u_2, u_3) A_i(2) + e_2(u_1, u_2, u_3) A_i(1) + e_3(u_1, u_2, u_3))
\end{equation}
which reads
\begin{equation}
\label{eq:comm_3const}
\left[ L^{\text{3-const.}}_i, L^{\text{3-const.}}_j \right] = t\sum_{\ell \geq 1} D_{ij, \ell}^{(3)} L^{\text{3-const.}}_\ell
\end{equation}

with\begin{equation*}
D_{ij, \ell}^{(3)} = D_{ij,\ell}(3) + e_1(u_1, u_2, u_3) D_{ij, \ell}(2) + e_2(u_1, u_2, u_3) D_{ij, \ell}(1) + e_3(u_1, u_2, u_3) D_{ij, \ell}(0).
\end{equation*}

It follows that Lemma~\ref{lemma:constr_from_ev} applies to $3$-constellations, allowing us to extract constraints.

\begin{theorem}
    For all $b>-1$, the generating function $Z^{(b)}_{3,1}\left[t,q,u_1,u_2,u_3,\textbf{p}\right]$ satisfies constraints $L^{3-\text{const}}_i \cdot Z = 0$ for all $i\geq 0$.
\end{theorem}

\subsubsection{Trivalent bipartite maps}

Details of the computations for this model can be found in Appendix~\ref{sec:comm_triv_bip}. The structure of the computation is similar to the previous model. We find.

\begin{proposition} \label{thm:SimplifiedCommutators2}
Let $\mu = \min(i,j)$ and $M = \max(i,j)$. For $1\leq m, m'\leq 3$,

\begin{align}
    \left[M^{(1,m)}_i, M^{(1,m)}_j\right] &= \sum_{l\geq 1} \tilde{D}_{ij, l}(m) M^{(1,m)}_l,\label{CommutatorM} \\
    \left[M^{(1,m)}_i, M^{(1,m')}_j\right] -  \left[M^{(1,m')}_i, M^{(1,m)}_j\right] &= \sum_{l\geq 1} \tilde{D}_{ij,l}(m') M^{(1,m)}_l + \tilde{D}_{ij,l}(m) M^{(1,m')}_l, \label{CommutatorMM}\\
    \left[p_i^*, M^{(1,m)}_j\right] -  \left[p_j^*, M^{(1,m)}_i\right] &= \sum_{l\geq 1} \tilde{D}_{ij,l}(m) p_l^* \label{CommutatorpM}
\end{align}
with $\tilde{D}_{ij, l}(1) = 0$ and $\tilde{D}_{ij, l}(2) = (i-j) \delta_{l,i+j-2}$ and
\begin{multline}
\tilde{D}_{ij, l}(3) = (i-j) (2\delta_{l\geq M-1} + \delta_{M-1\leq l\leq i+j-3}) J^{(b)}_{i+j-3-l} \\
+ \operatorname{sgn}(i-j)(2l-3\mu+3) \delta_{\mu-1\leq l\leq M-2} J^{(b)}_{i+j-3-l}
+ b(i-j)(i+j-3)\delta_{l, i+j-3}.
\end{multline}
\end{proposition}

This shows that the candidate constraints 
\begin{equation}
L^{\text{bip}\leq 3}_i \coloneqq L^{(1,3)} = -p_i^* + \sum_{m=1}^3 q_m t^m M^{(1,m)}_i
\end{equation}
satisfy commutation relations

\begin{multline}
\label{eq:comm_modes_tri_bip}
[L_i^{\text{bip}\leq 3}, L_j^{\text{bip}\leq 3}] = 2t^3q_3(i-j) \sum_{n\geq 1} p_n L^{\text{bip}\leq 3}_{i+j+n-3} + t^3 q_3b(i-j)(i+j-3) L^{\text{bip}\leq 3}_{i+j-3}\\
+ (1+b)t^3 q_3\Bigl(3(i-j) \sum_{n=1}^{\mu-2} + \operatorname{sgn}(i-j) \sum_{n=\mu-1}^{M-2} (2M-2n-\mu-3)\Bigr) p_n^* L_{i+j-3-n}^{\text{bip}\leq 3} \\
+ 3t^3q_3u(i-j) L^{\text{bip}\leq 3}_{i+j-3} + t^2 q_2 (i-j)L^{\text{bip}\leq 3}_{i+j-2}.
%[L_i^{\text{bip}\leq 3}, L_j^{\text{bip}\leq 3}] = 2t^3(i-j) \sum\limits_{n \geq M-1} J^{(b)}_{i+j-3-n} L_n^{\text{bip}\leq 3} +t^3 b(i-j)(i+j-3)L_{i+j-3}^{\text{bip}\leq 3} \\
%+ t^3 \sgn(i-j)\sum\limits_{n=\mu-1}^{M-2} (2n-3\mu+3)J^{(b)}_{i+j-3-n}L_n^{\text{bip}\leq 3} + t^3(i-j) \sum\limits_{n=M-1}^{i+j-3} J^{(b)}_{i+j-3-n}L_n^{\text{bip}\leq 3}
\end{multline}

Therefore the Lemma~\ref{lemma:constr_from_ev} applies we have the following theorem.
\begin{theorem}
    For all $b>-1$, the generating function $Z^{(b)}_{1,3}\left[t,q,u,\textbf{p}\right]$ satisfies constraints $L^{\text{bip}\leq 3}_i \cdot Z^{(b)}_{1,3}\left[t,q,u,\textbf{p}\right] = 0$ for all $i\geq 0$.
\end{theorem}

%As in the previous case, it can be checked that the constraints satisfy the same commutation relation as its modes~\eqref{eq:comm_modes_tri_bip}. The computation follows the same steps as in the case of $3$-constellation, starting from the commutation relation
%\begin{equation}
%\left[J^{(b)}_i, M^{(2,3)}_j \right] = i \sum\limits_{k \geq i-1} J^{(b)}_{i+j-3-k} J^{(b)}_k + i M^{(2,3)}_{i+j-1} + b(j-1)iJ^{(b)}_{i+j-3}.
%\end{equation}

\section{Conclusion}

\subsubsection{Analysis of the results }

Note that contrary to the quadratic case of maps, algebra formed by the modes $M^{(4,1)}$ and $M^{(3,2)}$ are not Lie algebra as the structure "constants" are, in fact, differential operators. in the variables $\textbf{p}$.

\medskip

Furthermore, comparing the commutators of the modes of the two models~\eqref{eq:comm_3const} and~\eqref{eq:comm_modes_tri_bip}, the two could seem to be related by a simple shift of variables, and a relabelling of some variables $u_i$. However this is a not the case. The difference in applications of the operator $Y_+$ induces non-trivial effects on the expression of the modes for models with similar order but different values of $k$ and $m$ as the operator $Y_i$ truncates modes that would yield a catalytic variable $y$ carrying negative index. For example, this induces different truncations in the range of the sum defining the recurrence relations~\eqref{eq:mod_rec_no_ui} and~\eqref{eq:mod_rec_no_ui'}. These differences prevent a direct identification of the modes of the two models.

\medskip

Finally, let us mention that in both cases studied here, the constraints of each model had already been derived at $b=0$ albeit through different methods. For $3$-constellation, they are a particular case of the results of~\cite{Fang_PHD} (albeit their commutator had not been computed) and for trivalent bipartite maps, it had been obtained in~\cite{MMM_Wtilde} using the Ward identities of the corresponding $2$-matrix model, though using a derivation that is not fully rigorous. The results obtained here are coherent with those obtained in both articles.

\subsubsection{Comments on the method}

The method employed here could in principle be generalized to other constellation models. However, the computation would get more complex than it is in the present case, both because of the inner complexity of the computation due to the increasing number of terms, and the sophistication of the recurrence relation if $m \neq 2$ or $k \neq 1$. Numerical computations don't seem to be easily accessible either for two reasons. First, the recurrence relations involve sums going to infinity which have to be cancelled altogether. Second, we had to use commutation relations and non-trivial rewriting of some terms to cancel some contributions. It appears that there is no direct way to give a prescription to ensure cancellations between the different terms.

\medskip

Similarly, since the constraints are expected to form an algebra as they allow us to characterize the generating function order by order in $t$, we could think of inferring the form of commutator of the constraints from the simpler (from a computational perspective) commutator $\left[J^{(b)}_i, M^{(m,k)}\right]$. Observe that in the cases above this computation does not match the form of the commutator of the modes $M^{(m,k)}$ immediately after it is antisymmetrized. Instead, we relied on our knowledge of the commutator of the evolution operator to match this expressions when computing the commutator of the constraints.

\subsubsection{Conjecture for the general case}

We conclude this Chapter with a conjecture on the expected behaviour of the constraints for evolution operators $M^{(k,m)}$. The operator $M^{(k,m)}(u_1,\ldots,u_k)$ can be expanded over elementary symmetric polynomials as
\begin{equation}
 M^{(k,m)}(u_1,\ldots,u_k) = \sum\limits_{\lambda \in \Lambda_{k,m}} e_{\lambda^{\perp}}(u_1,\ldots,u_k)\mathcal{A}^{(\lambda)}_i
\end{equation}
where $\Lambda_{k,m}$ is the set of partitions with at most $m$ parts and $\lambda_1 \leq k$ i.e. it corresponds to the set of partitions that fit in a $k \times m$ rectangle using Ferrers diagrams, $\lambda^{\perp}$ is the complementary of $\lambda$ in the $k \times m$ Ferrers diagram and
\begin{equation}
\mathcal{A}^{(\lambda)}_i = N_{\lambda}^{-1} \sum\limits_{\sigma \in \mathfrak{S}_m} A\left(\sigma\cdot \lambda\right)  
\end{equation}
where $\sigma$ acts on $\lambda$ by permuting its parts (including empty ones) and 
\begin{equation}
    N_{\lambda} = \prod\limits_{n \geq 0} m_n(\lambda)!
\end{equation}
with $m_n(\lambda)$ the number of parts of length $n$ in $\lambda$.

We make the following conjecture on the commutator of operators $\mathcal{A}^{(\lambda)}$.

\begin{conjecture}
For all partitions $\lambda,\mu$, there exists differential operators in variables $p_i$ such that
\begin{equation}
\left[\mathcal{A}^{(\lambda)},\mathcal{A}^{(\mu)} \right]_{[ij]} = \sum\limits_{ \ell \geq 1 } \frac{2-\delta_{\lambda,\mu}}{2} \left( D^{(\lambda)}_{ij,\ell} ,\mathcal{A}^{(\mu)}_{\ell} + D^{(\mu)}_{ij,\ell} ,\mathcal{A}^{(\lambda)}_{\ell} \right).
\end{equation}
\end{conjecture}

This conjecture implies that the commutator of evolution operators $M^{(k,m)}$ closes for all $k,m \geq 0$ with brackets given by
\begin{equation}
    D^{(k,m)}_{ij,\ell} = \sum\limits_{\lambda \in \Lambda_{k,m}} e_{\lambda^{\perp}}D^{(\lambda)}_{ij,\ell}. 
\end{equation}

Additionally, we conjecture that the commutator of the currents $J^{(b)}_i$ with $i>0$ with operators $\mathcal{A}^{(\lambda)}_j$ behaves such that the commutator of the candidate constraints $L_i$ is always closed
\begin{conjecture}
    For all partition $\lambda$, we have
    \begin{equation}
        \left[ J^{(b)},\mathcal{A}^{(\lambda)}\right]_{[i+1,j]} = \sum \limits_{\ell \geq 1} D^{(\lambda)}_{ij,\ell} J^{(b)}_{\ell}.
    \end{equation}
\end{conjecture}

Indeed, this conjecture implies that the commutator of operators $L^{(k,m)}$ closes with bracket
\begin{equation}
    \left[ L^{(k,m)}_i, L^{(k,m)}_j \right] = -\sum \limits_{\ell \geq 1} D^{(k,m)}_{ij,\ell} L^{(k,m)}_\ell.
\end{equation}

%% file: Chapters/Tensor_model.tex
\chapter{Random tensor models} % Main chapter title
\label{Chap:tens_mod} 

Random tensor models were first introduced in the $90$'s~\cite{Amb_tens,Sasa_tens,Gross_tens} as a generalization of matrix models in dimensions $d\geq 3$. Similarly to matrix models, each interaction term in the action can be interpreted as a particular simplicial complex of dimension $d$ and their perturbative expansion is given by a sum over possible gluings of these simplicial complexes and corresponds to the generating function for certain classes of piecewise linear manifolds. However, the perturbative expansion lacked a $\frac{1}{N}$-expansion for the next $20$ years until it was implemented for \emph{colored} tensor models~\cite{Gu_exp1,Gu_exp2,GuRi}. This expansion is organized via a parameter called the Gurau degree, but it is of a different nature than in the matrix case. In particular, it is not topological and the subleading orders are unknown aside from the first few orders in specific models. However, tensor models exhibit universality in the large $N$ limit as they are dominated by \emph{melonic graphs} which places them in the universality class of branched polymers~\cite{GuRy}. This mechanism opened the way to the study of many classes of tensor models beyond the colored case. It also holds for uncolored tensor models~\cite{BoGuRi}, multi-matrix models~\cite{Ferrari,FeVa} and irreducible tensor models~\cite{Ca18,CaHa21}. All these models are dominated by melonic graphs and thus fall in the universality class of branched polymers. This chapter draws from literature on the topic~\cite{Book_Gu,Book_Ta} to review these results. We first establish the connection between piecewise linear manifolds and random tensor models. We then show how we can construct an action for the random tensor such that it admits a $\frac{1}{N}$-expansion by an appropriate choice of scaling for the terms of the action. This chapter serves as a prelude before going into the study of the double-scaling limit of some particular random tensor models in the following chapter. 

\section{Simplicial complexes and tensor invariants}
\label{sec:invs_tensor_simpl_complex}

\subsection{Colored simplicial complexes}
\label{ssec:col_simp_com}

To extend the interpretation of the perturbative expansion of matrix models as a sum over classes of $2$D-piecewise linear manifolds, we first need to give an analogous description of PL-manifolds in dimension $d \geq 3$. To construct elementary building blocks of these manifolds, we first introduce \emph{colored simplicial complexes}.

\begin{definition}[Colored simplex]
\label{def:colo_simplex}
    A $d$-colored simplex is the simplex of dimension $d$ such that all of its $d+1$ facets are labeled with a unique color $c\in \llbracket 0,d \rrbracket$.
\end{definition}

\begin{figure}[!ht]
    \centering
    \includegraphics[scale=0.7]{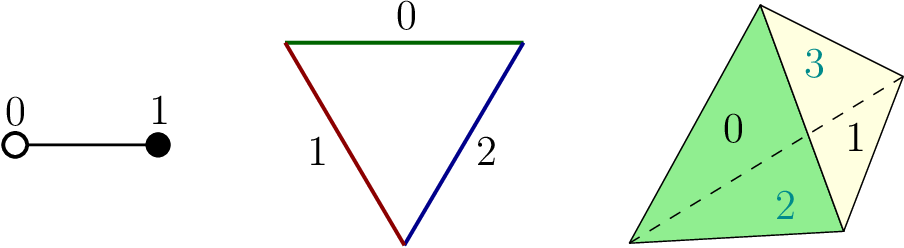}
    \caption{Colored simplices in dimension $1$, $2$ and $3$.}
    \label{fig:colored_simp}
\end{figure}

The coloring of the facets induces a bijection between $(d-k)$-subsimplices and subsets of $\llbracket 0,d \rrbracket$ of size $k$. Each $(d-k)$-subsimplex can be defined as the intersection of the facets that define it. For example, in dimension $3$, edges are given by pairs of distinct colors $(i,j)$ and vertices by triples $(i,j,k)$. The coloring allows to define a unique gluing between two $d$-colored simplices $s_1$ and $s_2$ by specifying the color $c\in\llbracket 0,d \rrbracket$ corresponding to the $(d-1)$-simplices along which the two simplices are glued. The gluing is uniquely defined by requiring that all lower dimensional simplices involving color $c$ of $s_1$ and $s_2$ are identified together. These gluings generate \emph{colored triangulations} by requiring that all $(d-1)$-subsimplices are shared by exactly two $d$-colored simplices.

\begin{figure}[!ht]
    \centering
    \includegraphics[scale=0.45]{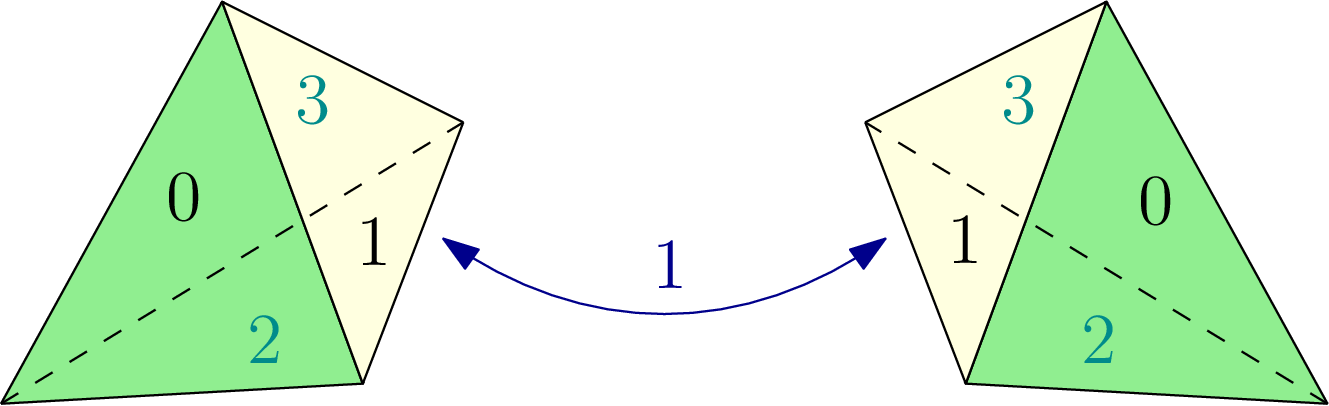}
    \caption{The unique gluing of two colored simplices of dimension $3$ along their facets of color $1$.}
    \label{fig:glu_simp}
\end{figure}

This condition ensures that the gluing between two simplices is locally homeomorphic to a ball $\mathbb{B}^d$ at the interior of the $(d-1)$-subsimplices. However, there is no reason that the neighborhood of a $(d-k)$-subsimplex should satisfy this condition. For this reason, the colored triangulations do not generate $d$-dimension manifolds but only \emph{pseudomanifolds} which allows for singularities of codimension at least $2$. This type of defect appears generically in PL-manifolds in dimension $d \geq 3$. The following theorem shows that colored triangulations can be used as building blocks of PL-manifolds.
\begin{theorem}[\cite{LiMa85}]
\label{thm:PL_triangu}
Every PL-manifolds admits a decomposition as a colored triangulation. Moreover, for every $d$ there are finitely many moves such that two homeomorphic decompositions are related by a finite sequence of those moves.
\end{theorem}

Since, in a colored triangulation, specifying the color of the facets along which two simplices are glued completely specifies the gluing between these two simplices, the colored triangulations are in bijection with \emph{properly $(d+1)$-colored graphs}.

\begin{definition}[Properly colored graph]
\label{def:colo_graph}
A graph $G=(V,E)$ is \emph{properly $d$-colored} if each vertex is incident to exactly one edge of color $c$ for all $c \in \llbracket 1,d\rrbracket$.
\end{definition}

Each $d$-simplex of a colored triangulation $\mathcal{T}^d$ corresponds to a vertex in the associated $d+1$-colored graph $G(\mathcal{T}^d)$ and the $d+1$ edges incident to a vertex encode the gluing with other $d$-simplices. Said differently, the $(d+1)$-colored graph is obtained as the $1$-skeleton of the dual triangulation. We stress that this compact encoding of the triangulation as $(d+1)$-colored graph is only possible because of the introduction of the colors on the simplices. This construction automatically extends to subsimplices of the colored triangulations: the $(d-k)$-subsimplices labeled with colors $(c_1,\dotsc,c_k)$ are given by connected components of the graphs $G(\mathcal{T}^d)_{\rvert c_1,\dotsc,c_k}$ obtained from $G(\mathcal{T}^d)$ by keeping only edges with color $(c_1,\dotsc,c_k)$ have been deleted. 

\medskip

However, note that in the two-dimensional case of matrix models, the elementary building blocks of the PL-manifolds are not colored triangles but uncolored polygons. To help clarify the role played by colored simplices in dimensions $d \geq 3$, let us first have a look at their role in the familiar two dimensional case, albeit in a slightly simpler setup. For any $n\geq 1$, a $2n$-gon can be divided into $2n$ triangles by adding a new vertex in its center and connecting each vertex of the $2n$-gon to this new vertex in a line. Assigning color $0$ to the edges of the original polygon, color $1$ to any newly created edge and coloring the other new edges such that each triangle has one edge of each color $\{0,1,2\}$, we obtain a decomposition of the $2n$-gon into $2n$ triangles, which corresponds to one possible triangulation of the $2n$-gon of Theorem~\ref{thm:PL_triangu}.

\begin{figure}
    \centering
    \includegraphics[scale=0.5]{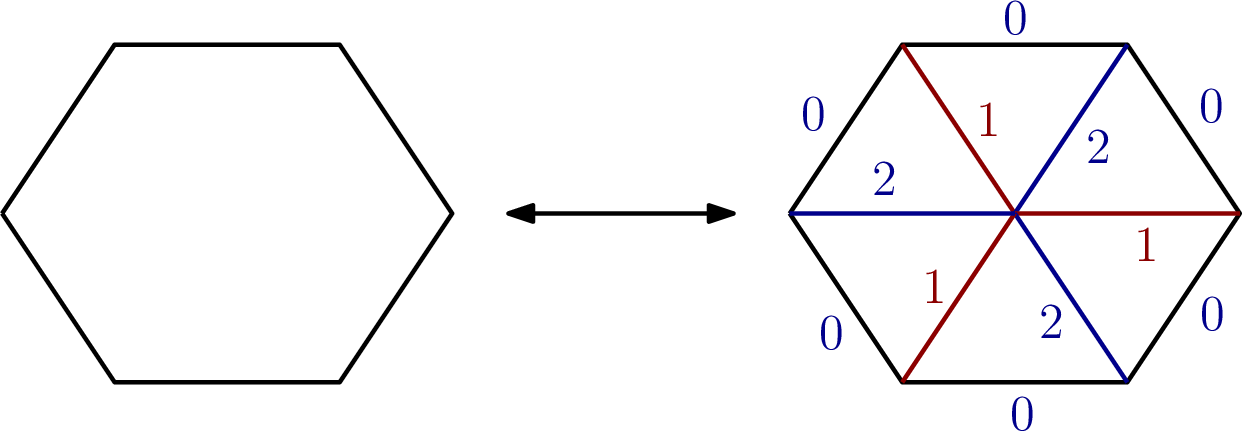}
    \caption{The hexagon, corresponding to a term $\Tr M^6$ in matrix model, can be decomposed in $6$ colored triangles.}
    \label{fig:poly_to_simpl}
\end{figure}

\medskip

On the random matrix side, the $2n$-gons are generated by including terms $\Tr M^{2n}$ in the action, which writes $\Tr M^{2n} = M_{i_1j_1}\dotsc M_{i_{2n}j_{2n}}$. Each copy of the matrix $M$ can be identified with a vertex of the colored graph so that two vertices are connected by an edge of color $c$ if their $c$-th index are contracted together. Finally, the edges of color $0$ of the $2n$-gon correspond to the identification of indices as obtained in the perturbative expansion via Wick's theorem.

\begin{remark}
\label{rem:bub_odd}
Note that the triangles alternate between triangles with edges colored $(0,1,2)$ in clockwise order and triangles with edges colored $(0,2,1)$ in the same order. In particular, this condition doesn't allow for non-even interactions. However, we have seen in Chapter~\ref{Chap:mat_mod} that such interactions do exist in the random Hermitian matrix model. This is because the Hermitian matrix model has \emph{mixed index symmetry} which allows for a larger class of interactions compared to the colored-graph framework, as we shall see below.
\end{remark}

\subsection{Bubbles and tensor invariants}
\label{ssec:bub_tens_inv}

The above construction generalizes to random tensor models. The building blocks of the $d$-dimensional PL-manifolds are given by connected $d$-colored graphs with color in $\llbracket 1,d \rrbracket$ called \emph{bubbles}. After the gluings of color $i \neq 0$ have been performed, we obtain a simplicial complex whose $(d-1)$-subsimplexes all carry color $0$ which can be seen as an uncolored simplicial complex. The corresponding polynomial in the tensor entries is obtained as follows. Each vertex of the graph corresponds to a copy of a tensor $T_{i_1\dotsc i_d}$ whose $c$-th index $i_c$ ranges in $\llbracket1,N_c\rrbracket$ for some integer $N_c$. If two vertices are connected by an edge of color $c$, then the $c$-th index of the corresponding copies of $T$ are contracted together. Hence, to each bubble with $n$ vertices corresponds a unique monomial of degree $n$ in the tensor entries $T_{i_1\dotsc i_d}$ with no free index. 

\begin{figure}
    \centering
    \includegraphics[scale=0.5]{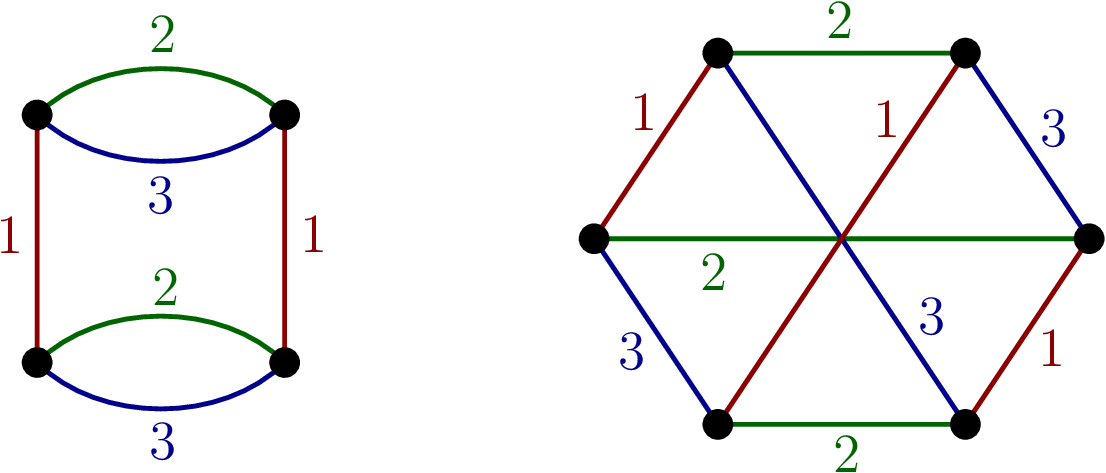}
    \caption{Two examples of bubbles of tensors of rank $3$.}
    \label{fig:ex_bubble}
\end{figure}

\medskip

In random tensor models, the set of possible interactions is restricted by requiring that the tensor is invariant under the action of some group $G$. Monomials in $T$ invariant under this action of $G$ are called \emph{tensor invariants} and correspond to the terms appearing in the action of the field theory. The tensor models whose tensor invariants are given by bubbles correspond to a particular type of model which have \emph{no index symmetry}.

\begin{definition}[No mixed index symmetry]
    A tensor model $T$ with the action of a group $G$ is said to have no index symmetry if each of the $d$ indices $i_k$ of the tensor $T_{i_1\dotsc i_k}$ transforms under the fundamental representation of $O(N_k)$ or $U(N_k)$. 
\end{definition}

These group actions are the only cases in which the action doesn't mix the indices of the tensor together. In that case, the $c$-th index of the tensor can be assigned color $c$ and we obtain a correspondence between bubbles and invariant for the corresponding group action. If the $c$-th index transforms under the unitary group, then the $c$-th index of $T$ must be contracted with the $c$-th index of its complex conjugate $\bar{T}$, which is encoded in the graph by assigning a color black or white to the vertex and requiring that edge of color $c$ connect black and white vertices. Up to this additional coloring of the vertices, there is a one-to-one correspondence between bubbles and tensor invariants. This bipartiteness condition translates the orientability of the associated PL-manifold as shown by the following proposition.

\begin{proposition}[\cite{FGG86}]
A PL-manifold is orientable if and only if its dual colored graph is bipartite.
\end{proposition}

Choosing all indices to transform under $U(N_c)$ ensures the bipartiteness of the colored graph and, in turn, the orientability of the associated simplicial complex, while having copies of $O(N_c)$ allows for non-orientable PL-manifolds. For these models, there is a unique quadratic invariant that plays the role of the kinetic term in the action. This corresponds to the Gaussian term for the tensor field $T_{i_1\dotsc i_d}$ where all $d$ indices of two copies of $T$ (or $T$ and $\bar{T}$) are contracted together. Via the Wick theorem, the perturbative expansion around the Gaussian can be interpreted as the gluing of the simplicial complexes along their $(d-1)$-simplices of color $0$, which form their boundaries. Note that while the Feynman graphs are regular $(d+1)$-colored graphs, the color $0$ plays a different role than the other ones as it identifies all $d$ indices of the copies of $T$ it connects while the remaining colors identify only one index. In the field theory terminology, the bubbles correspond to the vertices encoding the interactions while the edges of color $0$ correspond to free propagators. %When drawing tensor graphs we will represent edges of color $0$ using dashed edges to make this distinction explicit.

\medskip

Thus, the perturbative expansion of random tensors with no mixed index symmetry whose invariants are given by a set of $d$-colored graphs $\mathcal{B}$ can be interpreted as the generating function of $d$-dimensional PL-manifolds that can be obtained as gluings of the simplicial complexes associated with elements of $\mathcal{B}$. The corresponding Feynman graphs generated are $(d+1)$-colored graphs such that the connected components of the graph $G(\mathcal{T}^d)_{\rvert \hat{0}}$ obtained by deleting edges of color $0$ is a union of elements of $\mathcal{B}$.

\medskip

One major difference between the matrix case at $d=2$ and higher orders $d \geq 3$ lies in the growth of the number of tensor invariants. At $d=2$, there is a unique connected $2$-colored graph with $2n$ vertices obtained by forming one cycle of length $2n$ with alternating color $1$ and $2$, which corresponds to a term $\Tr (MM^{\dagger})^{n}$ in the matrix model. With tensors of rank $d\geq 3$, the number of invariants is given by the number of connected $d$-colored graphs. Since the graphs encode a contraction pattern of the indices, they can be described using matchings and thus can be counted using the representation theory of the symmetric group~\cite{ABGD_inv,BG_inv}. For example, when all indices transform under the unitary group the number $K_d(n)$ of (disconnected) invariants of order $n$ for a rank $d$ tensor is~\cite{BG_inv}
\begin{equation}
    \label{eq:nb_tens_inv}
    K_d(n) = \sum\limits_{\lambda \vdash n} z_\lambda^{d-2},
\end{equation}
where we recall that $z_\lambda = \prod\limits_{i} i^{m_i}m_i!$ with $m_i$ the number of parts of length $i$ in the partition $\lambda$.

\subsection{Some comments on tensor models with mixed index symmetry}

Tensor models with mixed index symmetry are obtained by considering a group action on $T$ of rank $d$ different than $O(N)^d$. They include for example irreducible tensor model~\cite{Ca18,CaHa21} and multi-matrix models~\cite{TaFe,Aze}. In full generality, a tensor model with mixed index symmetry does not admit an interpretation as the generating function of PL-manifold as their Feynman graphs cannot be described using colored graphs. Yet it may be the case, as we have seen in Chapter~\ref{Chap:mat_mod} with the Hermitian matrix model whose interactions terms $\Tr H^n$ are invariants of $H$ under action by conjugation of $U(N)$. Instead, the Feynman vertices for these models are represented using \emph{stranded graphs}. Each copy of a rank $d$ tensor is represented by a set of $d$ half-edges, and two half-edges are connected when the corresponding tensor index are contracted together. All bubbles can be represented as stranded graphs but the converse is not true since stranded graphs allow to contract together indices at different positions in their respective tensor copies. 

\begin{figure}
    \centering
    \includegraphics[scale=0.62]{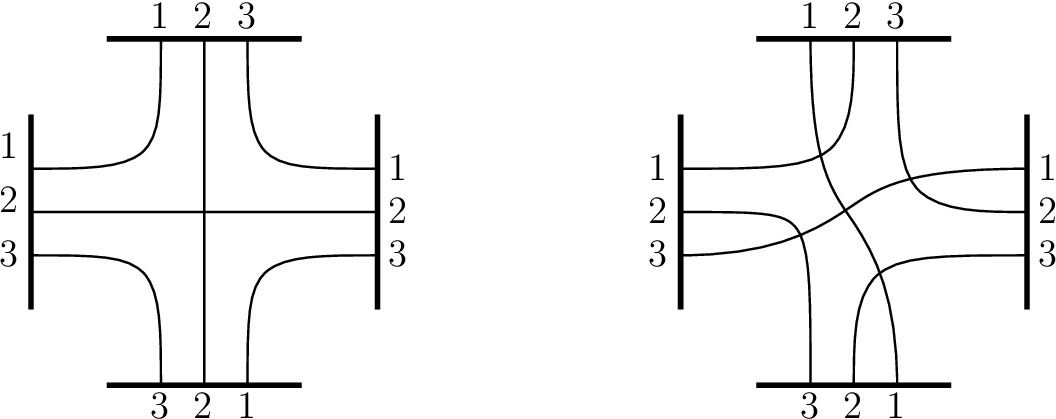}
    \caption{Two examples of stranded graphs. The right one cannot be represented using colored graphs.}
    \label{fig:stranded_graph}
\end{figure}

\medskip

The combinatorial properties of stranded graphs is more intricate than colored graphs. For example, it depends on the type of index mixing that is allowed, and the possible index contraction increases with the rank $d$ of the tensor, while colored graphs have a similar combinatorial structure for any $d$. Thus, far less is known on generic properties for this type of model, and their graphs are often studied via a different representation of its Feynman graphs (e.g. as loop decorated maps in the $U(N)\times O(D)$ tensor model of the Section~\ref{sec:DS_U(N)xO(D)}). However, in the cases studied in the literature~\cite{Ca18,CaHa21,TaFe}, they share similar properties than tensor models with no mixed index symmetry (albeit with additional necessary conditions on the structure of the tensor e.g. tracelessness~\cite{Ca18,BCGK_symm} or bipartiteness~\cite{TaFe}). The rest of this chapter focuses on the case of tensor models with no mixed index symmetry, but this shall be illustrated in the next chapter where we will implement the double scaling limit for the quartic bipartite $U(N)\times O(D)$ multi-matrix model (see Section~\ref{sec:DS_U(N)xO(D)}).

\section{Large \texorpdfstring{$N$}{N} expansion for tensor model}
\label{sec:tensor_large_N}

In the previous section, we have seen how to construct tensor invariants that can be interpreted as building blocks of PL-manifolds in dimension $d$, but we have yet to give an explicit action for a tensor model that admits a $\frac{1}{N}$-expansion. This mechanism was first established for colored tensor models in a series of paper~\cite{Gu_exp1,Gu_exp2,GuRi} and generalized shortly after to uncolored tensor models in~\cite{BoGuRi,DRT_exp_MO,TaCa} which will be framework described hereafter. While the previous subsection relied only on the encoding of simplicial complexes using colored graphs and did not involve the dimensions of the tensor $N_c$, a diagrammatic expansion organized by inverses of the dimensions is possible in a wide variety of models but the nature of this expansion changes when these dimensions differ. In this section, we will focus on the case where all $d$ dimensions $N_c$ of the tensor are equal. This mechanism will then be extended to the case of matrix-tensor models with indices in two different ranges in the next section. To illustrate the large $N$ expansion of tensor models, we consider a rank $d$ tensor $T=\left(T_{i_1\dotsc i_d}\right)_{1 \leq i_1,\dotsc ,i_d \leq N}$ invariant under the action of the group $O(N)^{\otimes d}$
\begin{equation}
    \label{eq:grp_action_tensor}
    T_{i_1\dotsc i_d} \rightarrow O^{(1)}_{i_1j_1}\dotsc O^{(d)}_{i_dj_d} T_{j_1\dotsc j_d}\hspace{5pt}, \qquad O^{(i)} \in O(N).
\end{equation}

The set of allowed bubbles $\mathcal{B}$ is exactly the set of $d$-colored graphs, which is the most general class that can be considered as the orthogonal group imposes no bipartiteness conditions on the edges of the graph. The cases of models where the $c$-th index of the tensor transforms under $U(N_c)$ can then be recovered by requiring that edges of color $c$ of the bubbles must respect the bipartite structure of the vertex of the graph. We further restrict ourselves here to connected $d$-colored graphs, which are analogous to single trace invariants in matrix models, but the arguments can be adapted for disconnected graphs with no major difficulty. For a bubble $B \in \mathcal{B}$, we denote the corresponding tensor invariant $I_B(T)$. We recall that there is a unique quadratic bubble $B_2$ which contracts all indices of the two copies together. We denote the corresponding invariant as $I_2(T) = \sum_{i_1,\dotsc,i_d} T_{i_1\dotsc, i_d} T_{i_1\dotsc, i_d}$. To all bubbles, we associate a coupling constant $t_B$ and a scaling factor $N^\rho(B)$. The scaling factor and coupling constant of the quadratic invariant can be absorbed in the definition of the tensor field $T$, therefore we set without loss of generality $t_2 = 1$ and $\rho(B_2) = 1$ and absorb a global factor $N^{d-1}$ in the definition of the action, anticipating on upcoming computation. The action for the random tensor model then writes
\begin{equation}
\label{eq:action_O(N)_model}
    S(T) = N^{d-1}\left(I_2(T) + \sum\limits_{b\in \mathcal{B}\setminus{I_2(T)}} t_b N^{\left(\rho(b)\right)}I_b(T)\right),
\end{equation}
and the corresponding partition function is
\begin{equation}
    \label{eq:part_tensor}
    Z[t_b] = \int dT \exp\left(N^{d-1}\left(I_2(T) + \sum\limits_{b\in \mathcal{B}\setminus{I_2(T)}} t_b N^{\rho(b)}I_b(T)\right)\right).
\end{equation}

The existence of a $\frac{1}{N}$-expansion for the partition function depends on the existence of an upper bound to the scaling in $N$ of the Feynman graphs and therefore on the choice of scaling function $\rho$. Moreover, when such an expansion exists, the scaling function $\rho$ also determines the structure of the expansion. To determine the possible choices of scaling function, let us compute the scaling exponent with $N$ of a Feynman graph of this model. Since the model has no mixed index symmetry, the contribution of each of the $d$ indices can be treated separately. The $c$-th index of the diagram propagates in the diagram by alternating between edges of color $0$ and color $c$, each cycle of alternating color $0$ and $c$ gives a free index $i_c$ is summed on and contributes as $N$. Even though they actually correspond to $(d-2)$-dimensional simplices in the PL-manifold, these cycles are called faces of color $l$. This is justified by the fact that they are faces of ribbon graphs associated with the graph, see Definition~\ref{def:jackets} below. The number of edges of the graph can be recovered from the valency $N_b$ of the bubbles $b$ as $\sum\limits_{b\in \mathcal{B}} \frac{n_bN_b}{2} $. The scaling exponent $s(G)$ with $N$ of a graph $G$ with $n_b(G)$ bubbles $b$ and $F_c(G)$ faces of color $c$ of the perturbative expansion is given by
\begin{equation}
\label{eq:scaling}
 s(G) = \sum\limits_{b\in \mathcal{B}} n_b \rho_b + \sum\limits_{c=1}^d F_c(G).
\end{equation}
Therefore the model admits a $1/N$ expansion only if the function $s$ is bounded from above. Under this form of the scaling factor, it is not clear whether the scaling factor is bounded for any graphs as this expression depends on the number of bubbles of the graph both explicitly and through the number of faces $F_c$. Therefore, we first rewrite this expression using \emph{jackets}~\cite{Gu_exp1} in order to count the number of faces of the graph.
\begin{definition}[Jackets~\cite{Gu_exp1}]
\label{def:jackets}
    For $\tau$ a cycle of length $d$ in $\mathfrak{S}_d$, the jacket $J_{\tau}(G)$ of a $(d+1)$-colored graph graph $G$ is the ribbon graph with all vertices and edges of $G$ faces of colors $(\tau^i(1),\tau^{i+1}(1))$ for $i\in\llbracket 0,d\rrbracket$ modulo the equivalence relation $\tau \sim \tau^{-1}$.
\end{definition}

\begin{figure}[!ht]
    \centering
    \includegraphics[scale=0.70]{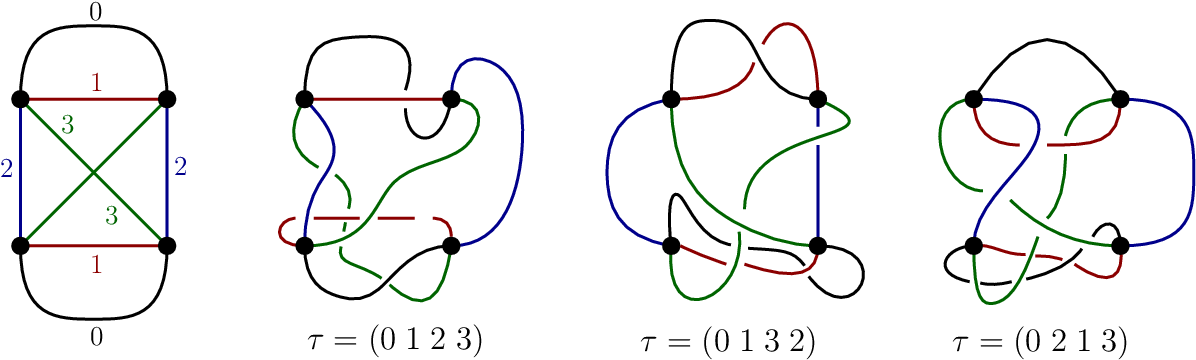}
    \caption{A $4$-colored graph and its three jackets.}
    \label{fig:jackets_OG}
\end{figure}

 A $d$-colored graph has exactly $\frac{(d-1)!}{2}$ different jackets. The choice of a cycle $\tau$ amounts to fixing one cyclic ordering of the colors for every vertex of the graph. The introduction of this canonical orientation for the vertices of the graph allows us to study $d$-colored graphs with the tools available for maps. Note that the jacket has the same connectivity as the original graph as it preserves all edges and vertices of the graph. The Euler characteristic of a jacket $J_\tau$ of a $d$-colored graph with $V=2v$ vertex is
\begin{equation}
\label{eq:euler_jacket}
2-2k_\tau = 2v - Dv + F_\tau,
\end{equation}
where $k_\tau$ is the non-orientable genus of the jacket and $F_\tau$ its number of faces. A jacket $J_\tau$ has faces of color $(i,j)$ if and only if $\tau(i)=j$ or $\tau(j)=i$. Hence there are exactly $(d-2)!$ jackets with faces of color $(i,j)$. It follows that the number of bicolored cycles of $G$ is
\begin{equation}
\label{eq:sum_euler_jacket}
    \sum\limits_{i,j} F_{ij} = (d-1) + \frac{(d-1)(d-2)}{2}v - \frac{2}{(d-2)!}\sum\limits_{\tau \in \mathfrak{C}_d} k_\tau,
\end{equation}
where $\mathfrak{C}_d$ denotes the set of cycles of length $d$ in $\mathfrak{S}_d$. We define the Gurau degree (hereafter called degree) of a $d$-colored graph $G$ as 
\begin{equation}
    \label{eq:G_degree}
    \omega(G) = \sum\limits_{\tau \in \mathfrak{C}_d} k_\tau(G).
\end{equation}
In two dimensions, it coincides with the genus of the unique jacket associated with the $3$-colored graphs and we recover the $\frac{1}{N}$-expansion of matrix models. However, for $d\geq 3$ it is no longer a topological invariant, and should be thought of as a linear bound between the number of bubbles of the graph and its number of faces.

\medskip

The number of faces of color $l$ can now be computed by noticing that faces of color $(i,j)$ with $i,j \neq 0$ occur only inside of bubbles. Therefore, in a graph $G$ we can apply formula~\eqref{eq:sum_euler_jacket} to $G$ as a $(d+1)$-colored graph and subtract contributions of each of its bubbles. Rather than performing this computation, we shall compute directly the amplitude of a Feynman graph of the model. We shift the scaling function $\rho$ writing it as $\rho(b) = \tilde{\rho}(b) -\frac{2}{(d-2)!}\omega(b)$ for some function $\tilde{\rho}$. We can see that our choice of normalization of the action and the redefinition of the scaling function kills the last two contributions of formula~\eqref{eq:sum_euler_jacket} when applied to the bubbles. Therefore, the scaling of a graph $G$ is given by 

\begin{equation}
    s(G) = D - \frac{2}{(d-1)!}\omega(G) + \sum\limits_{b \in G} \tilde{\rho}(b)
\end{equation}
where the sum is over all bubbles $b$ of $G$ counted with multiplicities. We can see that taking $\tilde{\rho}(b)=0$ ensures the existence of a $\frac{1}{N}$-expansion since the degree is a non-negative integer. This choice corresponds to the original scaling factor of~\cite{Gu_exp1,Gu_exp2,BoGuRi}. 

\subsection{Leading order of the \texorpdfstring{$\frac{1}{N}$}{1/N}-expansion}
\label{ssec:melons}

In all this subsection, we set $\tilde{\rho}(b)=0$. Now that the existence of the $\frac{1}{N}$-expansion for this choice of scaling factor has been established, we would like to characterize graphs at each order of this expansion. However, the combinatorics of the $(d+1)$-colored graphs for $d \geq 3$ are more intricate than in the matrix case $d=2$. In most tensor models, only the leading order of the $\frac{1}{N}$ are fully characterized. The leading order graphs are graphs with vanishing Gurau degree $\omega(G) = 0$. Using Euler relation~\eqref{eq:euler_jacket}, these graphs are the graphs that have a maximal number of bicolored cycles for a fixed number of vertices. 

\medskip

All faces appearing in the jackets of $G$ also appear in the jackets of some graphs $G^{(c)}$ obtained by deleting edges of color $c$ in $G$. Using equation~\eqref{eq:G_degree} for $G^{(c)}$ and~\eqref{eq:sum_euler_jacket} in both $G$ and $G^{(c)}$, we obtain the following expression for the degree of $G$ as a function of the degree of the subgraphs $G^{(c)}$
\begin{equation}
\label{eq:deg_rec_jacket}
\omega(G) = \frac{(d-1)!}{2}\left(d+v-\sum\limits_c C^{(c)} \right) + \sum\limits_{c} \omega(G^{(c)}),
\end{equation}
where $C^{(c)}$ is the number of connected components of $G^{(c)}$ and $\omega(G^{(c)})$ is the sum of the degree of the connected components of $G^{(c)}$. The first term of the r.h.s of this equation can be shown to be non-negative by induction on the number of vertices, and the second one is non-negative by definition of the degree. Therefore the graphs of vanishing degree are exactly the graphs with $v$ vertices such that
\begin{equation}
    \label{eq:char_LO}
    \forall c\in\llbracket 0,d \rrbracket \hspace{5pt} \omega(G^{(c)}) = 0 \hspace{5pt}\text{and} \hspace{5pt} d+v-\sum\limits_c C^{(c)} = 0.
\end{equation}

\begin{figure}
\hfill
\begin{subfigure}{.2\textwidth}
  \centering
  \includegraphics[scale=0.22]{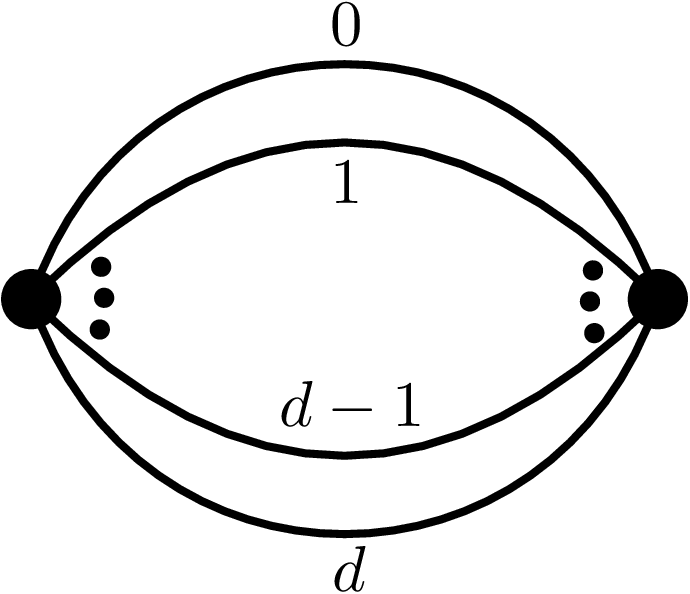}
  \label{fig:elem_melon}
  \caption{The elementary melon.}
\end{subfigure}
\hfill
\begin{subfigure}{.6\textwidth}
  \centering
  \includegraphics[scale=0.47]{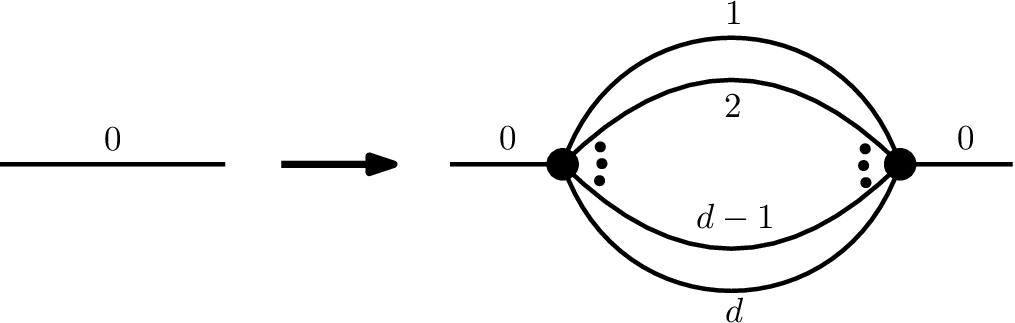}
  \caption{The melonic move (on an edge of color $0$).}
  \label{fig:mel_move}
\end{subfigure}
\caption{Melonic move inserts the elementary melon edge of any color (here the color $0$). In the graph, the root edge is marked with a square.}
\label{fig:melon_insert}
\end{figure}

When $d=2$, this condition characterizes connected planar maps and the degree coincides with the genus of the corresponding map. For $d \geq 3$, the graphs satisfying these conditions are known as \emph{melonic graphs} or melons.
\begin{definition}[Melonic graphs~\cite{Gu_exp2}]
    The elementary melon is the only $(d+1)$-colored graph with two vertices. The melonic graphs is the set of graph that can be obtained by successive iterations of the move described by figure~\ref{fig:mel_move} on the elementary melon.
\end{definition}
The melonic graphs can be thought of as "hyperplanar" since all their jackets are planar. Seeing the Euler relation as a linear relation between the number of vertices and faces, planar graphs can be characterized as graphs that have maximal number of faces when fixing the valency of its vertices. Similarly, the melonic are graphs with maximal number of faces i.e. bicolored alternating cycles in the colored graph at a fixed number of vertices. 

\medskip 

The elementary melon satisfies both conditions of~\eqref{eq:char_LO}, and the combinatorial move of Figure~\ref{fig:mel_move} preserves both properties. Therefore melons are leading order graphs. Conversely, we can prove that they are the only graphs satisfying both properties~\cite[Lemma 4]{Gu_exp2}. The set of melonic graphs is in bijection with planar $(d+1)$-ary trees. To construct this bijection and enumerate melonic graph, it is more convenient to work with rooted objects. We consider melonic graphs with a marked edge which we can without loss of generality assume to have color $0$. The elementary melon corresponds to the rooted tree with a single vertex and $d+1$ leaves. Then, performing the move of Figure~\ref{fig:mel_move} on an edge of color $c$ corresponds to the addition of a vertex and $d+1$ edges on the corresponding edge of color $c$. 

\begin{figure}
    \centering
    \includegraphics[scale=0.75]{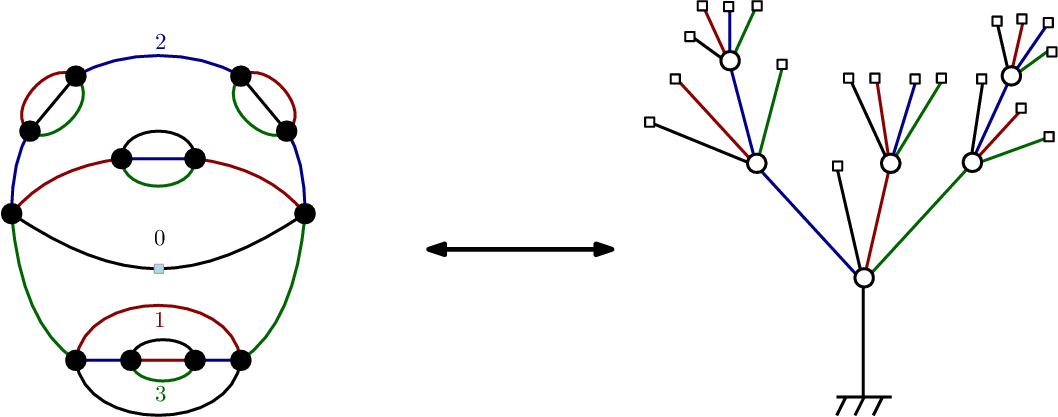}
    \caption{Example of the bijection between melonic graphs and trees for $d=3$.}
    \label{fig:melon_to_tree}
\end{figure}

\medskip 

Via this bijection, it is clear that the melonic graphs fall into the universality class of branched polymers. In terms of PL-manifold, they correspond to a particular class of triangulations of the $d$-sphere called \emph{stacked triangulations}~\cite{GuRy}. This is in stark contrast with the $d=2$ case which retained all PL-manifolds with the topology of the $2$-sphere at leading order of the $\frac{1}{N}$-expansion. It is possible to show that as long as a graph has a planar jacket, it corresponds to a triangulation of the $d$-sphere~\cite[Proposition $4.3$]{GuRy2}. As a consequence, there can exist triangulations of the $d$-sphere of arbitrarily large degree and the first graphs which do not have the topology of the sphere can only appear at order at least $\frac{d!}{2}$.

\subsection{Optimal scaling at \texorpdfstring{$d=3$}{d=3}}
\label{ssec:dim_3}

Due to conditions~\eqref{eq:char_LO}, we have in particular $\omega(G^{(0)}) = 0$ for leading order graphs. It implies that the bubbles of leading order graphs must be melonic graphs themselves. Conversely, any bubble that is not melonic cannot contribute at leading order of the $\frac{1}{N}$-expansion. To get a larger family of graphs -and thus a richer set of PL-manifold surviving in the large $N$ limit- we can tune $\tilde{\rho}(b)$ appropriately to enhance the contribution of some non-melonic bubbles. However, this might spoil the existence of the $\frac{1}{N}$-expansion depending on how the insertion of the bubbles $b$ in the graph changes its degree. For example, if $\tilde{\rho}(b) > 0 $ for a melonic bubbles then the existence $\frac{1}{N}$-expansion is immediately lost. Indeed, since these bubbles have only planar jackets, a graph made of these bubbles has a vanishing degree if and only if the graph obtained by contracting all edges of color $c \neq 0$ is planar. Clearly, these graphs exist and can be arbitrarily large. For these graphs, choosing $\tilde{\rho}(b)>0$ for each bubble would give a net positive contribution, leading to the scaling $s$ being unbounded. This leads to the following definition.
\begin{definition}[Optimal scaling~\cite{FeVa}]]
    A scaling factor $\tilde{\rho}$ is said to be \emph{optimal} for a set of bubbles $\mathcal{B}$ if it cannot be enhanced without spoiling the existence of the $\frac{1}{N}$-expansion.
\end{definition}
For example, from the argument above, the scaling factor $\tilde{\rho}(b)=0$ is optimal for melonic bubbles. Determining optimal scaling for generic values of $d$ and interactions $\mathcal{B}$ is a difficult problem as it requires a fine understanding of the combinatorics of the graph and its subgraphs. Thus the optimal scaling is only known for particular types of models or interactions. Essentially, it is known for \emph{maximally single trace} bubbles - which are bubbles with exactly one cycle of each pair of color $(i,j)$- for any $d$~\cite{FeVa}, and in the case $d=3$~\cite{TaCa}. 

\medskip 

To illustrate the previous discussion as well as introducing the $O(N)^3$ model, we present the case $d=3$ following~\cite{TaCa}. When $d=3$, the bubbles are $3$-colored graphs which simplifies the analysis. Since identifying how to choose the scaling factor $\tilde{\rho}$ such that it is optimal is difficult, we will work the other way around and present a different $\frac{1}{N}$-expansion which we will then identify as the optimal scaling for this model, follwing~\cite{TaCa}. The action of the model writes
\begin{equation}
\label{eq:action_O(N)3}
    S(T) = N^{\frac{3}{2}}\left(I_2(T) + \sum\limits_{b\in \mathcal{B}\setminus{I_2(T)}} t_b N^{\rho(b)}I_b(T)\right).
\end{equation}    
Therefore the scaling of a graph $G$ with $E(G)$ edges, $n_b$ bubbles $b$ and $F=\sum\limits_{c=1}^d F_d$ faces is
\begin{equation}
    \label{eq:scal_O(N)3}
    s(G) = -\frac{3}{2}E(G) +\sum\limits_{b \in \mathcal{B}} \left(\rho(b)-\frac{3}{2}\right)n_b - F(G).
\end{equation}

To study this expansion, we introduce a different notion of jackets than in the previous case.
\begin{definition}[Jackets~\cite{FeVa}]
\label{def:jackets_bis}
The jackets $\tilde{J}^{(c)}(G)$ of a $4$-colored graph $G$ are the graph obtained when deleting edges of color $c \in \{1,2,3\}$ in $G$ and contracting the cycles of the remaining non-zero color $(c',c'')$. Performing this contraction preserves the genus of the ribbon graph and maps cycles with color $(c',c'')$ to vertices and faces of color $(0,c')$ and $(0,c'')$ to faces.
\end{definition}

\begin{figure}[!ht]
    \centering
    \includegraphics[scale=0.5]{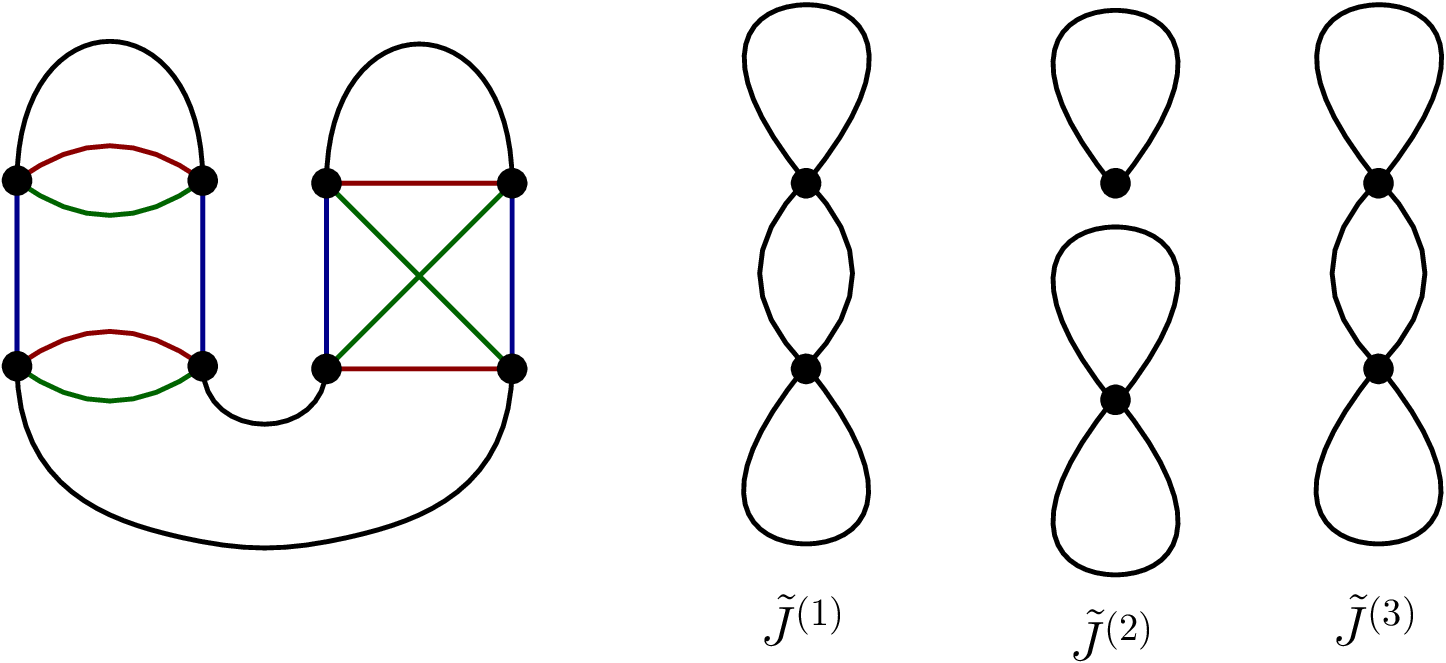}
    \caption{The second type of jackets for a $4$-colored graph.}
    \label{fig:jackets_NG}
\end{figure}

\medskip 

Contrary to the previous case, these jackets have fewer edges than the initial graph therefore they can have more connected components than $G$. We denote $C^{(c)}=\delta^{(c)}+1$ the number of connected components of the jacket $\tilde{J}^{(c)}(G)$ and similarly for a bubble $b$ we denote $C^{(l)}_b=\delta^{(c)}_b+1$ the number of connected components of the bubble after deletion of edges of color $l$. The Euler characteristic for the jacket $\tilde{J}^{(c)}(G)$ for a graph $G$ reads
\begin{equation}
    2C^{(c)} - 2k^{(c)} = \sum\limits_{b\in \mathcal{B}} n_bC^{(c)}_b - \frac{N_bn_b}{2} - F_{c'} - F_{c''},
\end{equation}
where $k^{(c)}(G)$ is the non-orientable genus of the jacket $\tilde{J}^{(c)}(G)$ and $N_b$ the valency of the bubble $b$.

\medskip 

The faces of color $c$ appear in $2$ of the jackets. Summing the Euler relations~\eqref{eq:euler_jacket} for all jackets we get
\begin{equation}
    2\sum\limits_{c=1}^d F_d(G) = \sum\limits_{c\in\{1,2,3\}} 2C^{(c)}(G)-k^{(l)}(G) + \sum\limits_{b\in \mathcal{B}} n_b\frac{d(d-1)N_b}{4}  - \sum\limits_{c} \sum\limits_{b\in \mathcal{B}} n_bC^{(c)}_b.
\end{equation}

With our choice of scaling for the quadratic bubble, the contribution of the propagator cancels the terms in $N_b$ and the scaling function $s(G)$ reads
\begin{align}
     s(G) &= 3 + \sum\limits_{c,c'} \frac{1}{2}\left( 2\delta^{(c)}(G)-k^{(c)}(G)\right) + \sum\limits_{b\in \mathcal{B}} n_b\left(\rho(b)-\frac{1}{2}\sum\limits_{c} \delta^{(c)}_b\right).
\end{align}

We now choose the following scaling factor $\rho_{opt}$
\begin{equation}
    \label{eq:scal_opt}
    \rho_{opt}(b) = -\frac{1}{2} \sum\limits_{c} \delta^{(c)}_b.
\end{equation}
Alternatively, it can be expressed as
\begin{equation}
\label{eq:scal_opt2}
    \rho_{opt}(b) = \frac{3}{2}-\frac{1}{2}\sum\limits_{c,c'} F_{cc'}(b),
\end{equation}
since the number of connected components of the bubble $b$ in the jacket $\tilde{J}^{(c)}(G)$ is given by its number of cycles of alternating color $c'$ and $c''$ in $b$, which allows for straightforward computation. This scaling factor leads to degree
\begin{equation}
\label{eq:opt_deg}
    \tilde{\omega}(G) = \sum\limits_{c} \frac{1}{2}k^{(c)} + \left(\sum\limits_{b\in \mathcal{B}} n_b\sum\limits_{c} \delta^{(c)}_b-\delta^{(c)} \right).
\end{equation}
The fact that the degree is bounded from below and thus leads to a well-defined $\frac{1}{N}$-expansion follows from the following Lemma.

\begin{lemma}[\cite{TaCa}]
For any graph $G$ and any jacket $\tilde{J}^{(c)}(G)$ we have 
\end{lemma}
\begin{equation}
    \label{eq:ineq_jack}
    \sum\limits_{b\in\mathcal{B}} n_b \delta_b^{(c)} \geq \delta^{(c)}.
\end{equation}
This inequality simply states that the jackets of a graph have strictly less connected components than the union of its bubbles which is clear since the graph can be recovered from its set of bubbles by adding edges of color $0$ appropriately. Moreover, there exist subsets of bubbles such that this bound can be reached for all three jackets simultaneously e.g. melonic bubbles. This leads to the following proposition

\begin{proposition}[\cite{TaCa}]
 The scaling factor $\rho_{opt}$ is optimal for the set of $3$-colored graphs.
\end{proposition}

\begin{figure}
\hfill
\begin{subfigure}{.45\textwidth}
  \centering
  \includegraphics[scale=0.6]{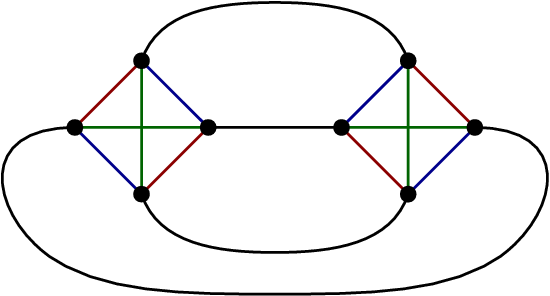}
  \label{fig:melon_enh1}
  %\caption{The elementary melon.}
\end{subfigure}
\hfill
\begin{subfigure}{.45\textwidth}
  \centering
  \includegraphics[scale=0.40]{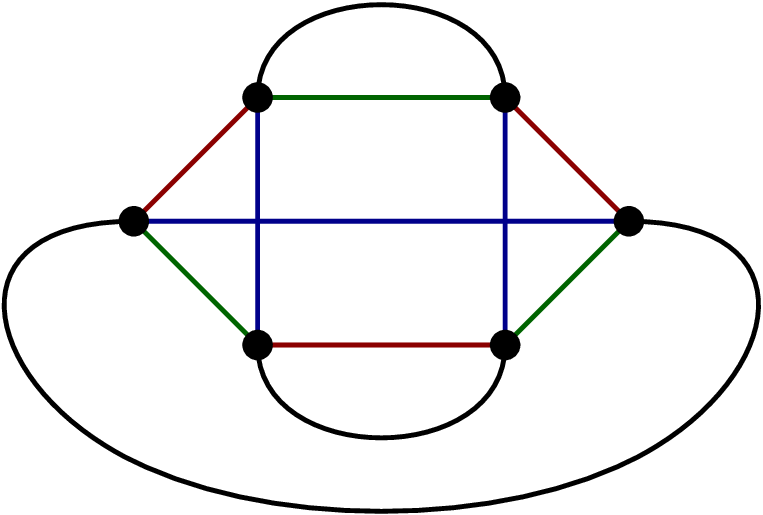}
  %\caption{The melonic move (on an edge of color $0$).}
  \label{fig:melon_enh2}
\end{subfigure}
\caption{Two examples of leading order graphs using non-melonic bubbles.}
\label{fig:melon_enh}
\end{figure}

This scaling factor enhances the scaling with $N$ of the non-melonic bubbles so that they contribute at leading order. For example, the following graph with bubbles given by the $4$-clique $K_4$ contributes at leading order while it was subleading before. An example of a leading order graph with this bubble is given on the left of Figure~\ref{fig:melon_enh}. Hence the expansion is a different $\frac{1}{N}$-expansion than in Section~\ref{sec:tensor_large_N} and includes the largest possible class of graph at leading order when including all possible bubbles.

\medskip 

\begin{figure}[!ht]
    \centering
    \includegraphics[scale=0.5]{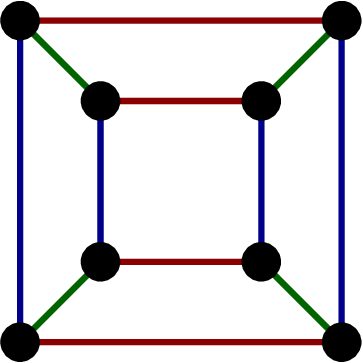}
    \caption{This bubble cannot appear in any leading order graph.}
    \label{fig:octa_bubble}
\end{figure}

Note that there still exists some bubbles that do not contribute in any leading order graphs. For example, consider the bubble associated with the octahedron, represented in Figure~\ref{fig:octa_bubble}. For a graph to contribute at leading order, inequality~\eqref{eq:ineq_jack} must be saturated for each of its three jackets. The inequality is saturated when the jackets of the graph $G$ have as many connected components as its bubbles. Therefore, for each of the three jackets, the edges of color $0$ must connect two vertices in the same connected component of a bubble. This is not possible for the octahedral bubble of Figure~\ref{fig:octa_bubble} since for that bubble, there is no pair of vertices that belong to the same connected component for all three jackets. It follows that this bubble cannot appear at leading order of the expansion. This type of interaction has been studied in~\cite{BoLi} where they obtained a formula for an enhanced degree so that graphs with this octahedric bubble contribute at leading order which is not compatible with the choice of scaling we work with here.

\subsection{Beyond leading order}
\label{ssec:beyond_LO}

Let us conclude this section by saying a few words on the subleading orders of the $\frac{1}{N}$-expansion of tensor models. In matrix models, the topological nature of the $\frac{1}{N}$-expansion gives an explicit characterization of the maps appearing at $g$-th order of the expansion: they are genus $g$ maps. For $d \geq 3$ the $\frac{1}{N}$-expansion is no longer topological, thus there is no such straightforward characterization of the subleading order graphs. They are solely defined by their degree i.e. by their combinatorial characteristics. For example, the conditions~\eqref{eq:char_LO} can be generalized to arbitrary order $\omega$ for the Gurau degree~\eqref{eq:deg_rec_jacket} by requiring that the sum of the two contributions of~\eqref{eq:char_LO} sum to $\omega$ i.e.
\begin{equation}
    \label{eq:char_NLO}
\hspace{5pt} d+v-\sum\limits_c C^{(c)} = K \hspace{5pt}\text{such that }\hspace{3pt} K + \sum\limits_{c=0}^{d} \omega(G)^{(c)} = \omega.
\end{equation}
Roughly speaking, this allows us to derive conditions on the graphs of a given order from the fact that they should not differ too much from graphs that are one order lower, as only one of the values $\omega(G^{(c)})$ or $K$ increases between two consecutive orders. However, the possible combination giving a graph of degree $\omega$ grows as the number of partitions of $\omega$ in $d+2$ parts which makes this approach difficult to carry through arbitrary order. In the best known cases, for the quartic $O(N)^3$ tensor model and the multi-orientable tensor model $U(N)^2\times O(N)$, their graphs have been characterized up to degree $\omega = \frac{3}{2}$ in~\cite{BNT_subleading}. However, it seems out-of-reach to characterize sets of graphs of arbitrary degree $\omega$ with this approach since the complexity of the combinatorial problem increases with the degree.

\section{Large \texorpdfstring{$N,D$}{N,D} expansion for multi-matrix model}
\label{ssec:large_D_mat}

Matrix-tensor models generalize the tensor models of the previous sections. They are tensor models of rank $d+2$ where two indices play a singular role with respect to the previous case. The last $d$ indices range in $\llbracket 1,D \rrbracket$ for some integer $D$ while the first two now range in $\llbracket 1,N \rrbracket$. This tensor can be interpreted as a tensor of dimension $d$ whose entries are matrices, i.e. as a set of $D^d$ matrices $X_{i_1,\dotsc,i_d}$. The faces of color $(1,2)$ correspond to traces of products of the matrices $X_{i_1,\dotsc,i_d}$, hence the bubbles can be interpreted as sums of products of those traces. The set of allowed bubbles doesn't depend on the dimensions of the tensor, therefore we will work with the largest class of bubbles and require the tensor to be invariant under the action of the group $O(D)^{\otimes d} \bigotimes O(N)^{\otimes 2}$. Again, we restrict here to connected bubbles for simplicity of the exposition but the arguments can be extended to non-connected bubbles.

\begin{equation}
    \label{eq:grp_action_tensor_mat}
    T_{i_1\dotsc i_d} \rightarrow O^{(1)}_{i_1j_1}\dotsc O^{(d)}_{i_dj_d} \tilde{O}^{(1)}_{aa'}\tilde{O}^{(2)}_{bb'}T_{a'b'j_1\dotsc j_d}, \qquad O^{(i)} \in O(N),\hspace{2pt}\tilde{O}^{(1,2)} \in O(D).
\end{equation}

While changing the range of the indices doesn't change the combinatorial structure of the graphs, it can induce changes to the structure of the perturbative expansion of the associated partition function. Still, when setting $D=N$ we must recover the $\frac{1}{N}$-expansion of random tensor models. We consider the action

\begin{equation}
    \label{eq:action_tens_mat}
    S(X) = ND^d \left(\Tr X_{i_1,\dotsc,i_d} \Tr X_{i_1,\dotsc,i_d} + \sum\limits_{b \in \mathcal{B}} N^{1-F_{12}(b)}D^{\rho(b)}I_b(X_{i_1,\dotsc,i_d})\right).
\end{equation}

This choice of scaling of the action coincides with~\eqref{eq:action_O(N)_model} when $D=N$. Since it doesn't depend on the choice of scaling factor, we first check that the $\frac{1}{N}$-expansion is well-defined. The scaling with respect to $N$ of a graph $G$ with $E(G)$ edges and $V(G)$ bubbles is given by 

\begin{equation}
    \label{eq:scaleN_TM}
    s_N(G) = V(G) - E(G) +\sum\limits_{b \in \mathcal{B}} F_{12}(b) + F_1(G) + F_2(G).
\end{equation}

All quantities appearing in this equation are encoded in the ribbon graph $G_{(012)}$ obtained by keeping only edges of color $c\in\{0,1,2\}$ in $G$ and then contracting all cycles of alternating color $(1,2)$. We denote with index $(012)$ all quantities associated with this graph. Its Euler characteristic is
\begin{equation}
\label{eq:euler_012}
    2c_{(012)}-2g_{(012)} = v_{(012)}-e_{(012)} + F_1(G)+F_2(G).
\end{equation}
Its number of vertices corresponds to the number of cycle of color $(1,2)$ of $G$ and it has $E(G)$ edges. Plugging this expression for $F_1(G)+F_2(G)$ into~\eqref{eq:scaleN_TM} gives
\begin{equation}
    \label{eq:scaleN_TM'}
    s_N(G) = -2g_{(012)} + 2c_{(012)} -2\sum\limits_{b \in \mathcal{B}} \left(F_{12}(b)-1\right)
\end{equation}
This scaling can be written as $s_N(G)=2-h(G)$ with
\begin{equation}
\label{eq:genus_TM}
    h(G) = 2g_{(012)} +2\left(1-c_{(012)} +\sum\limits_{b \in \mathcal{B}} \left(F_{12}(b)-1\right)\right).
\end{equation}
The genus $g_{(012)}$ is clearly non-negative. The second term is also non-negative as for a connected graph $G$, $G_{(012)}$ can only have as many connected components as there are multi-trace terms in the bubbles of $G$. Therefore $h(G)$ is non-negative and the $\frac{1}{N}$-expansion is well-defined. This is nothing but the usual $\frac{1}{N}$ expansion of matrix models when including multi-trace terms. Since it coincides with $g_{(012)}$ when all bubbles are single-trace terms, the parameter $h(G)$ is called the genus of the graph even though it doesn't correspond to the genus of the graph $G_{(012)}$ when there are multi-trace terms. 

\medskip 

We now turn to the structure of the $\frac{1}{D}$-expansion which depends on the choice of scaling factor. Consider first the scaling factor of tensor models of section~\ref{sec:tensor_large_N} associated with the Gurau degree. Since the bubbles are $d+2$-colored graph, this scaling factor reads for a bubble $b\in \mathcal{B}$ is
\begin{equation}
    \label{eq:scal_Gur_TM}
        \rho(b) = -\frac{2}{d!}\omega(b) + F_{12}(b).
\end{equation}
where we have substracted the contribution of faces of color $(1,2)$ which correspond to matrix traces and scale with $N$. Therefore the scaling with $D$ of a connected graph $G$ reads
\begin{equation}
    s_D(G) = dE(G)+2F_{12}(G)-\sum\limits_{b\in \mathcal{B}}\frac{2}{(d+1)!}\omega(b) +\sum\limits_{c =3}^{d+2} F_c(G).
\end{equation}
Using equation~\eqref{eq:genus_TM} and~\eqref{eq:sum_euler_jacket}, this can be written as
\begin{equation}
\label{eq:scaleD_TM}
s_D(G) = d +h(G) - \frac{2}{(d+1)!}\omega(G).
\end{equation}
For a graph $G$, we can obtain a lower bound for its degree from the degree of its subgraphs. In particular we have
\begin{equation}
    \label{eq:deg_bound}
    \omega(G) \geq \frac{1}{2}(d+2)! \omega(G_{(012)}).
\end{equation}
Since $h(G)$ and $g_{012}$ are related by~\eqref{eq:genus_TM} and $\omega(G)$ contains a term $g_{012}$, this inequality can be used to show that the contribution to the $\frac{1}{D}$-expansion of a graph $G$ is bounded by $d$ independently of $h(G)$ i.e. of its order in the $\frac{1}{N}$-expansion. Therefore, the $\frac{1}{D}$-expansion is always well-defined and the two expansions in $\frac{1}{N}$ and $\frac{1}{D}$ commute and can be taken in any order.

\medskip
Now consider the enhanced scaling factor
\begin{equation}
    \label{eq:scal_Gur_enh_TM}
    \rho'(b) = \frac{2}{d+1!}\omega(b) + F_{12}(b).
\end{equation}

The scaling with respect to $D$ of a graph $G$ reads
\begin{equation}
    s_D(G) = dE(G)+2F_{12}(G)+\sum\limits_{b\in \mathcal{B}}\frac{2}{d+1!}\omega(b) +\sum\limits_{c =3}^{d+2} F_c(G),
\end{equation}
and it can be expanded as~\cite{FeVa},
\begin{equation}
     s_D(G) = d+h(G)-\frac{2}{d+1}\sum\limits_{3\leq c,c'\leq d+2} h_{cc'}(G)
\end{equation}
where $h_{cc'}(G)$ is defined identically as $h(G)$ by~\eqref{eq:genus_TM} for the graph $G_{(0cc')}$.

\medskip

We define the grade $\ell(G)$ of a graph $G$ as 
\begin{equation}
    \label{eq:def_grade}
    \frac{\ell}{2} = \frac{2}{d+1}\sum\limits_{3\leq c,c'\leq d+2} h_{cc'}(G)
\end{equation}
As it is a sum of genuses of graphs, the grade is clearly a non-negative integer.

\medskip 

For this choice of scaling factor, the scaling with $D$ of a graph $G$ can be arbitrarily large depending on the genus of the graphs $G_{(012)}$. This implies that the $\frac{1}{D}$-expansion is ill-defined at fixed $N$. Contrary to the previous case, it is crucial to perform the $\frac{1}{N}$-expansion first to ensure the existence of an upper bound of the scaling with $D$ at each order of the $\frac{1}{N}$-expansion. This scaling factor generalizes the scaling factor~\eqref{eq:scal_O(N)3}. It is optimal for all \emph{maximally single trace interaction}~\cite[Proposition $3.2$]{FeVa} which are bubbles $b$ such that for every pair of colors $(i,j)$, the graph $b_{(ij)}$ where edges of $b$ of color $i$ and $j$ have been deleted is connected.

%% file: Chapters/Double_scaling.tex
\chapter{Double scaling limit in tensor models} % Main chapter title
\label{Chap:DScale}

As in the two dimensional case of matrix models, the double scaling limit of tensor models is related to a continuum limit and is attained by sending the coupling constants of the field theory to a critical value for its partition function while the size of the tensor goes to infinity. But only the leading order of the $\frac{1}{N}$-expansion is fully characterized in tensor models, while the double-scaling limit enhances the contributions of some graphs of subleading order, picking up contributions from all orders of the $\frac{1}{N}$-expansion. Yet it is possible to access these contributions using the \emph{scheme decomposition}. The scheme decomposition is a combinatorial tool initially introduced for maps~\cite{ChaMa}. The schemes of a family of maps (e.g. of maps genus $g$) is a "blueprint" that repackage together infinitely many maps sharing similar combinatorial structures. The scheme decomposition was first adapted to tensor models in~\cite{GuSch} where it was applied to quartic colored tensor model~\footnote{In the quartic case, the double scaling limit was also studied by an independent method in~\cite{DaGuRi} where it was derived using the Loop-Vertex-Expansion.}. The scheme decomposition gives information on the singularities of the generating function of graphs of subleading order which is sufficient to extract the graphs contributing at leading order in the double-scaling limit. This method has been used to derive the double scaling limit in several other tensor models e.g. for the multi-orientable tensor model~\cite{TaFu,TaGu} and for the $U(N)^2 \times O(D)$ model with tetrahedral interaction in~\cite{BeCa}.

\medskip 

In this chapter, we implement the double scaling in two different types of models. First, we implement this method to the $O(N)^3$ tensor model~\cite{TaCa} with all quartic interactions. Additional coupling constant complexifies the structure of the generating function of melonic graphs and could lead to different phases corresponding to different critical points. However, we find that this is not the case. There is a single phase associated with a particular combinatorial structure called \emph{broken chains}. The combinatorial tools introduced for this model can be extended to multi-matrix models as we show by extending the results of~\cite{BeCa} of the $U(N)^2 \times O(D)$ to include the pillow interactions. Finally, we show that these tools can also be used in models with mixed index symmetry and implement the double scaling mechanism in the $U(N)\times O(D)$ tensor model with tetrahedral interaction, a model whose Feynman graphs are not given by colored graphs but maps decorated with loops.

\medskip 

This chapter is an edited version of the articles~\cite{BNT_DS_TM} and~\cite{BNT_DS_MM}.

\section{Overview of the strategy}
\label{sec:plan_DS_limit}

While the computation depends on the details of the tensor model e.g. range of the indices and symmetry group, the method employed to access the double scaling limit consists of the same steps independently of the model. It is achieved through a three steps method: 
\begin{enumerate}
\item \label{enum:Classify} Classify the graphs according to their degree via the scheme decomposition.
\item \label{enum:WriteGF} Write the generating series of graphs at a fixed degree in terms of known series and identify its singularities.
\item \label{enum:DS} Describe the most singular contributions and resum them using a double scaling limit.
\end{enumerate}
We give below an overview of the tools used to implement this method for a given tensor model. 
\medskip

\paragraph*{Scheme decomposition.} There are infinitely many graphs contributing at each degree of the $\frac{1}{N}$-expansion. What we are looking for is packing the infinities into well controlled graphical objects, which will be the \emph{melons} and \emph{chains} (also known, in the theoretical physics literature, as ladders) whose generating function can be explicited. Melons can be eliminated or added on every edge without changing the degree, while chains can be extended or reduced without changing the degree either. The \emph{schemes} will then be defined as graphs without melons and with minimal chains. This gives the following result.
\begin{theorem}
\label{thm:graph-scheme}
Any 2-point graph can be reconstructed from a unique scheme by extending some chains and adding some melons on the edges.
\end{theorem}
In other words, each individual scheme represents an infinite family of graphs obtained by extending some chains and adding some melons. It is therefore possible to repackage the sum of all graphs of any given degree $\omega$ as a sum over schemes of the same degree. The 2-point function then reads
\begin{equation*}
G_\omega = \sum_{\substack{\text{Schemes $\mathcal{S}$}\\\text{ of degree $\omega$}}} P_{\mathcal{S}}(C(M), M),
\end{equation*}
where $P_{\mathcal{S}}$ is a polynomial, $C$ the generating series of chains and $M$ the generating series of melons. The quantity $P_{\mathcal{S}}(C(M), M)$ is the amplitude resulting from the sum over all graphs in the family of the scheme $\mathcal{S}$. The singularities of $G_\omega$ may then come from the series $C$, $M$ and from the sum over schemes of degree $\omega$ if there is an infinite number of them. However, this turns out not to be the case. We will prove the following result:
\begin{theorem}
\label{th:sch}
The set of schemes of a given degree is finite.
\end{theorem}

Enumerating all schemes of a given degree is still a hard combinatorial problem as it amounts to characterizing the graph of a given order in the $\frac{1}{N}$ expansion of tensor models. However, in the double scaling limit, we only need to identify a subset of schemes of any given degree which we will be able to obtain explicitly.

\medskip

\paragraph*{Dominant schemes and double scaling limit.} The double scaling limit consists of taking the large $N$ limit while sending the coupling constants $\lambda$ to a critical value of the partition function while maintaining a certain parameter $\kappa(N, \lambda)$ fixed. Note that this parameter $\kappa$ is defined such that the contribution of non-melonic graphs are enhanced in this limit so that graphs of arbitrarily large order can contribute. 

\medskip 

Since there is a finite number of schemes of any fixed degree $\omega$, the sum over schemes in $G_\omega$ is a finite sum. All singularities then come the series $C$ and $M$, and the double scaling limit can be obtained by finding the schemes for which $P_{\mathcal{S}}(C(M),M)$ is most singular. Those schemes are said to be \emph{dominant} because they are the most divergent when the coupling constant gets close to its critical value. As we will see in the sequel, the dominant schemes are the ones that have a maximal number of a specific type of chains, called \textit{broken} chains. This leads to the following Theorem

\begin{theorem}
\label{thm:dom_sch_tree}
The dominant schemes are in bijection with decorated rooted binary trees.
\end{theorem}

This Theorem is obtained through an explicit combinatorial decomposition of the dominant schemes. The decorations of the tree depend on the details of the model but can easily be explicited via this decomposition. It allows us to get an explicit expression for the generating function of the dominant schemes for any model where the theorem can be applied, ultimately allowing us to give the expression of the $2$-point function in the double scaling limit.

\section{Double scaling limit for the quartic \texorpdfstring{$O(N)^3$}{O(N)3} tensor model}
\label{sec:DS_O(N)^3}

We first apply this method to the case of the quartic $O(N)^3$ tensor model of~\cite{TaCa} and studied in the previous chapter of this manuscript (see Sec.~\ref{ssec:dim_3}). 

\subsection{Interactions and its \texorpdfstring{$\frac{1}{N}$}{1/N}-expansion}
\label{ssec:def}

As in~\cite{TaCa}, we limit ourselves to quartic interactions. This leaves one quadratic invariant bubble polynomial, which gives rise to the propagator of the model, and four invariant bubble polynomials which are the interaction terms of the model. These terms are respectively called the tetrahedral, and pillow of color $i=1,2,3$ (or quartic melonic bubble) interactions. The tetrahedral interaction corresponds to the only quartic maximally single trace interaction. The pillow of color $i$ is by convention the one that is disconnected when the edges of color $i$ are removed. The tetrahedral bubble is obviously invariant under color permutations, but the pillows are not and are instead swapped under color permutations. They correspond to the following invariants:
\begin{align}
I_k(\phi) &= \sum_{a, b, c} \phi_{abc}\phi_{abc} = 
\begin{array}{c}\includegraphics[scale=.35]{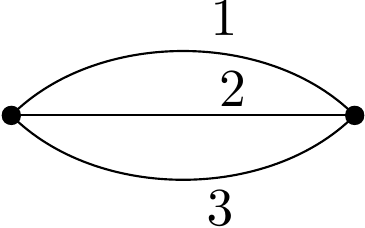}\end{array}\\
I_t(\phi) &= \sum_{a, a', b, b', c, c'} \phi_{abc}\phi_{ab'c'}\phi_{a'bc'}\phi_{a'b'c} = \begin{array}{c}\includegraphics[scale=.3]{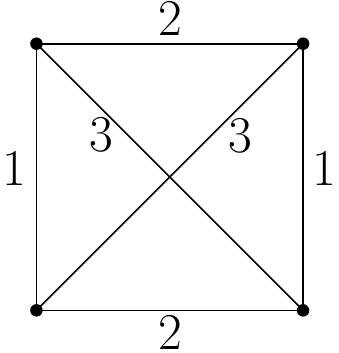}\end{array}\\
I_{p,1}(\phi) &= \sum_{a, a', b, b', c, c'} \phi_{abc}\phi_{a'bc}\ \phi_{ab'c'}\phi_{a'b'c'} = \begin{array}{c}\includegraphics[scale=.3]{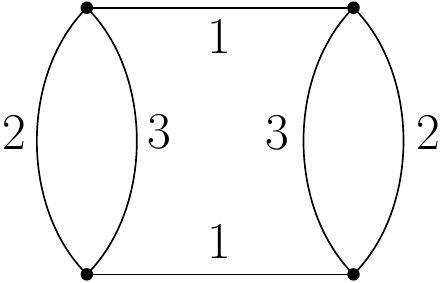}\end{array}\\
I_{p,2}(\phi) &= \sum_{a, a', b, b', c, c'} \phi_{abc}\phi_{ab'c}\ \phi_{a'bc'}\phi_{a'b'c'} = \begin{array}{c}\includegraphics[scale=.3]{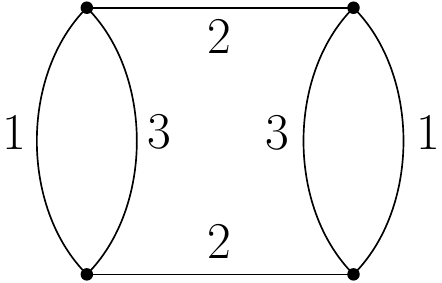}\end{array}\\
I_{p,3}(\phi) &= \sum_{a, a', b, b', c, c'} \phi_{abc}\phi_{abc'}\ \phi_{a'b'c}\phi_{a'b'c'} = \begin{array}{c}\includegraphics[scale=.3]{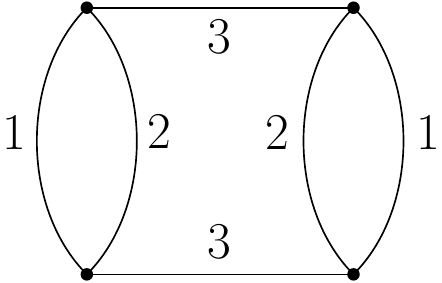}\end{array}
\end{align}

We will work with scaling factor~\eqref{eq:scal_opt} leading to the largest possible family of graphs contributing at leading order for this model where all interactions contribute at leading order. The action of the quartic $O(N)^3$-invariant tensor model thus writes

\begin{equation} \label{O(N)3Action}
S_N (\phi) = -\frac{N^{\frac{3}{2}}}{2}I_k(\phi) + N^{\frac{3}{2}} \frac{\lambda_1}{4}I_t(\phi) + N \frac{\lambda_2}{4}\Bigl(I_{p,1}(\phi) + I_{p, 2}(\phi) + I_{p, 3}(\phi)\Bigr)
\end{equation}
and its partition function is
\begin{equation}
Z_{N}(\lambda_1, \lambda_2) = \int \prod_{a,b,c=1}^N d\phi_{abc}\ e^{S_N(\phi)}.
\label{eq:part_fct_O(N)3}
\end{equation}

Thus, the $\frac{1}{N}$ expansion of the free energy reads \cite{TaCa}
\begin{equation}
F_{N}(\lambda_1, \lambda_2) = \ln Z_{N}(\lambda_1, \lambda_2) = \sum_{{ \bar{\cG}}\in{\bar{\mathbb{G}}_{}}} N^{3-\omega(\bar{\cG})} \mathcal{A}({\bar{\cG}}). %= \sum_{\omega \in \frac{\mathbb{N}}{2}} N^{3-\omega}\ \sum_{\substack{\cG\in\mathbb{G}_{}\\  \omega(\cG) = \omega}} \tilde{\mathcal{A}}(\cG)
\label{eq:largeN}
\end{equation}
The set $\bar{\mathbb{G}}$ is the set of connected $4-$regular properly-edge-colored graphs such that the subgraph obtained by removing all edges of color 0 is a disjoint union of tetrahedral bubbles and pillows. For quartic interactions, the degree~\eqref{eq:opt_deg} can be expressed as
\begin{equation}
\omega (\bar{\cG}) = 3 + \frac{3}{2}n_t(\bar{\cG}) + 2n_p(\bar{\cG}) - F(\bar{\cG}),
\label{eq:deg_O(N)3}
\end{equation}
where
\begin{itemize}
\item $n_t(\bar{\cG})$ and $n_p(\bar{\cG})$ are respectively the number of tetrahedral and pillow bubbles in the graph,
\item $F(\bar{\cG})$ is the number of faces of the graph. Here a face is defined as a cycle of alternating colors $\{0,i\}$ for $i \in \{1,2,3\}$. The \emph{degree} of a face is the number of edges of color 0 incident to it.
\end{itemize}

\medskip

\paragraph*{2-point graphs.} The 2-point function is 
\begin{equation}
\langle \phi_{abc} \phi_{a'b'c'}\rangle = \frac{1}{N^3} G_N(\lambda_1, \lambda_2)\ \delta_{aa'} \delta_{bb'} \delta_{cc'}
\end{equation}
with $G_N(\lambda_1, \lambda_2) = \left\langle \sum_{i,j,k} \phi_{ijk} \phi_{ijk}\right\rangle$ as a consequence of the $O(N)^3$-invariance. It has an expansion on 2-point graphs, similar to that of the free energy. We denote $\mathbb{G}$ the set of 2-point graphs. From a graph $\mathcal{G}\in\mathbb{G}_{}$, we can obtain a vacuum graph $\bar{\mathcal{G}}$ by connecting the two half-edges together. $\bar{\mathcal{G}}$ is moreover equipped with a marked edge -the one formed by connecting the two half-edges- called a \emph{root}. Rooted graphs are the Feynman graphs of the expansion of $G_N(\lambda_1, \lambda_2)$.

\medskip

Calculating the free energy from its Feynman expansion onto vacuum graphs requires taking into account the graph automorphisms which are difficult to track combinatorially. In contrast, rooted graphs have no symmetry. This makes the computation of the 2-point function more straightforward than that of the free energy (in fact, rooting objects is often the first step in this type of combinatorial problems).

\medskip

To extract the free energy out of the 2-point function, one can introduce a coupling constant for the quadratic part of the action, and then integrate the 2-point function with respect to it. This can be achieved by first rescaling the variables as follows. One rescales the coupling constants $\lambda_1, \lambda_2$ by $1/t^2$ for some parameter $t$, then rescales $\phi$ by $\sqrt{t}$ so that $t$ now only appears in front of the quadratic terms of the action. Denote $\tilde{\lambda}_{1,2} = t^2 \lambda_{1,2}$ and $\tilde{\phi} = \phi/t$, then
\begin{equation}
S_N(\phi) = -\frac{N^{\frac{3}{2}}t}{2} I_k(\tilde{\phi}) + N^{\frac{3}{2}} \frac{\tilde{\lambda}_1}{4}I_t(\tilde{\phi}) + N \frac{\tilde{\lambda}_2}{4}\Bigl(I_{p,1}(\tilde{\phi}) + I_{p, 2}(\tilde{\phi}) + I_{p, 3}(\tilde{\phi})\Bigr).
\end{equation}
In this normalization, the free energy is obtained by integrating the 2-point function with respect to $t$, at fixed $\tilde{\lambda}_1, \tilde{\lambda}_2$. The latter is
\begin{equation}
\langle \tilde{\phi}_{abc} \tilde{\phi}_{a'b'c'}\rangle = \frac{1}{t^2 N^3} G_N(\tilde{\lambda}_1/t^2, \tilde{\lambda}_2/t^2)\ \delta_{aa'} \delta_{bb'} \delta_{cc'}.
\end{equation}

\vspace{10pt}
The difference of exponents of $N$ between a 2-point graph $\cG$ and the associated vacuum graph $\bar{\cG}$ is just a $N^3$ due to opening/closing 3 faces in a systematic way. Therefore, we use as a convention for the degree of $\cG\in\mathbb{G}$ the degree $\omega(\bar{\cG})$ of the graph obtained by connecting the two half-edges.

\medskip

In the following, we will represent bubbles with partial coloring whenever possible i.e. when the omitted colors can be placed anyhow on the remaining pair of edges.

\subsubsection{Melons, dipoles, chains and schemes}
\label{sssec:MDCS_O(N)3}

We introduce here the different families of subgraphs which will play a key role in the upcoming analysis.

\paragraph{Melons\\}

Melons of the quartic $O(N)^3$-invariant tensor model were introduced and studied in~\cite{TaCa}. Their characterization is analogous to the one performed in Subsection~\ref{ssec:melons} but for the scaling factor~\ref{eq:scal_O(N)3}. We give here a short summary of the relevant results. 
\begin{definition}[Melons]
The connected melonic graphs -or melons- are the leading order graphs of the $\frac{1}{N}$-expansion of the free energy~\ref{eq:largeN}.
\label{def:mel}
\end{definition}
The structure of melonic graphs relies on the following elementary building blocks:
\begin{definition}[Elementary melons]
A melonic graph is elementary if it is a 2-point graph with no non-trivial melonic $2$-point subgraph. We also say that it is an elementary melon.
\label{def:elem_mel}
\end{definition}
Melonic graphs can then be proved to be the set of graphs obtained by recursively inserting elementary melons on arbitrary edges, starting from the elementary melon itself. 

\medskip

There is an elementary melon associated with each of the quartic bubbles. An elementary melon is said to be of type I if it comes from a tetrahedral bubble and of type II if it comes from a pillow bubble.
\begin{equation}
\text{Type I} = \includegraphics[scale=.5, valign=c]{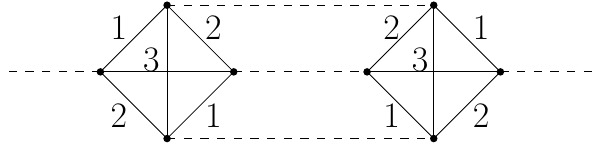} \qquad \text{Type II} = \includegraphics[scale=.5, valign=c]{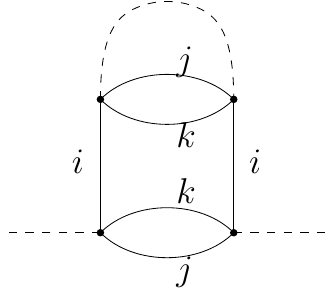}
\end{equation}
where $(i,j,k)$ is a cyclic permutation of $(1,2,3)$. Denote $M(\lambda_1, \lambda_2)$ the generating series of melonic 2-point graphs. From the recursive decomposition of melonic graphs, one has 
\begin{equation}
\includegraphics[scale=.5, valign=c]{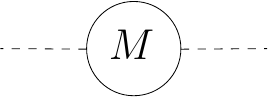} = \includegraphics[scale=.5,valign=c]{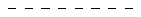} + \includegraphics[scale=.45, valign=c]{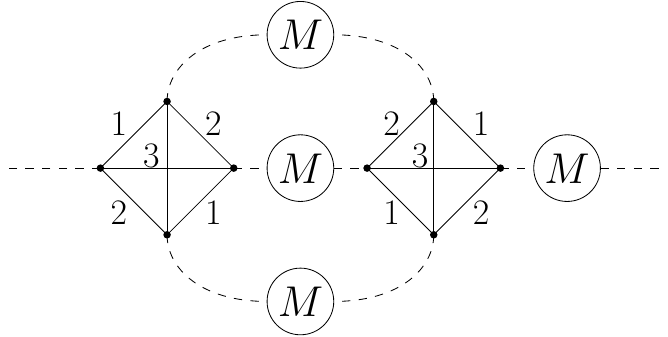} + \sum_{(i,j,k)} \includegraphics[scale=.45, valign=c]{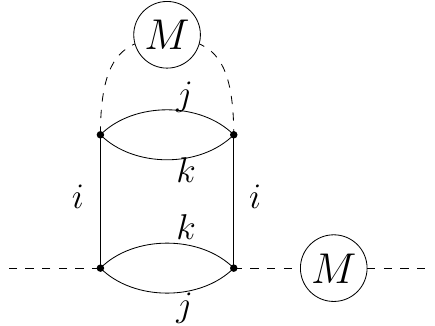}
\end{equation}
which leads to the following equation
\begin{equation}
M(\lambda_1,\lambda_2) = 1 + \lambda_1^2M(\lambda_1,\lambda_2)^4 + 3\lambda_2M(\lambda_1,\lambda_2)^2
\end{equation}
We change to the variables $(t,\mu) = (\lambda_1^2, \frac{3\lambda_2}{\lambda_1^2})$ (while retaining the notation $M$ for the melonic 2-point function), so that
\begin{equation}
M(t,\mu) = 1 + tM(t,\mu)^4 + t\mu M(t,\mu)^2
\label{eq:mel}
\end{equation}

which is the equation satisfied by the generating function of melons as studied in~\cite{TaCa}. With these variables, $t$ counts the number of melonic insertions performed and $\mu$ is the ratio between the number of melons of type I and II.

\medskip

In particular, its critical points are known. For a fixed value of $\mu \geq 0$, there is a single critical value $t_c(\mu)$ such that $(t_c(\mu),\mu)$ is a critical point of $M(t,\mu)$ and its behaviour near this critical point is
\begin{equation}
M(t,\mu) \underset{t \rightarrow t_c(\mu)}{\sim} M_c(\mu) + K(\mu)\sqrt{1 - \frac{t}{t_c(\mu)}}
\label{eq:crit_behav}
\end{equation}
where $M_c(\mu)$ is the unique positive real root of the polynomial equation 
\begin{equation} \label{eq:Mc_poly}
-3x^3+4x^2-\mu x +2\mu = 0,
\end{equation}
and $K(\mu) = \sqrt{\frac{M_c(\mu)^2\left(M_c(\mu)^2 + \mu \right)}{6M_c(\mu)^2+\mu}}$.

\paragraph{Dipoles\\}

The dipoles are obtained from the melons as follows.

\begin{definition}[Dipoles]
A dipole is a $4$-point graph obtained by cutting an edge in an elementary melon.
\label{def:dip}
\end{definition}

Starting from an elementary melon of type I, cutting two edges of color 0 leaves a single face of degree $2$ untouched. If this face has color $i$, we get a dipole of type I and color $i$,
\begin{equation}
\includegraphics[scale=.6,valign=c]{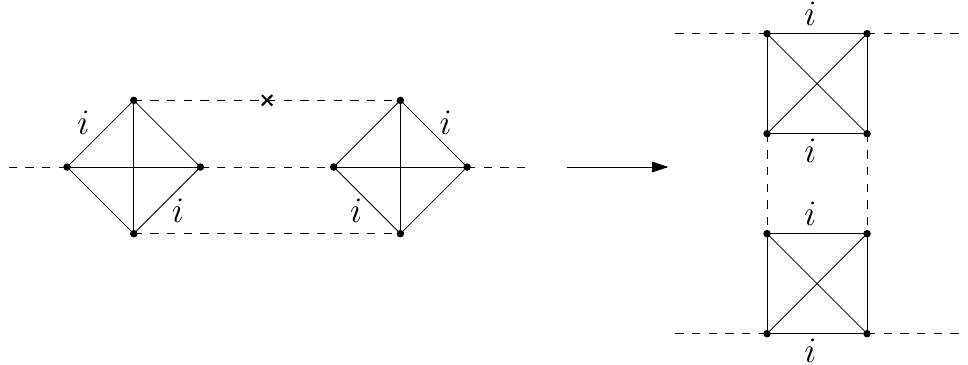}
\end{equation}
A dipole of type II is obtained by cutting the edge of color 0 of an elementary melon of type II. Dipoles of type II are simply given the color of their corresponding pillow interaction, 
\begin{equation}
\includegraphics[scale=.6,valign=c]{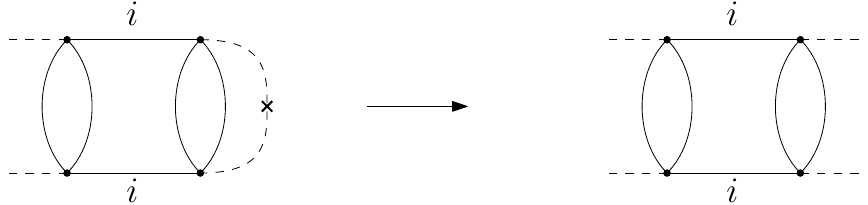}
\end{equation}
Note that dipoles of type II are always bubble-disjoint, meaning that two of them can never share a bubble. However, dipoles of type I may not be bubble-disjoint, in which case we say that they are \emph{non-isolated}. It turns out that non-isolated dipoles can only occur in the following subgraph,
\begin{equation}
    \includegraphics[scale=.5]{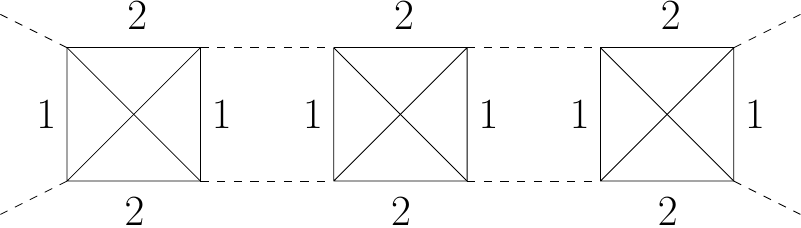}
\end{equation}
up to color permutations. Dipoles that are bubble-disjoint from others are said to be \emph{isolated}.

\medskip

It is easy to check that changing a dipole of type I and color $i$ in a graph $\cG$ for a dipole of type II of the same color does not change the degree of $\cG$. Therefore, we can replace dipoles of either type which are bubble-disjoint by a \emph{dipole-vertex} which represents any of them. In terms of generating series, it means that the dipole-vertex of color $i$ is the sum of the two dipoles of color $i$,
\begin{equation}
\includegraphics[scale=0.45,valign=c]{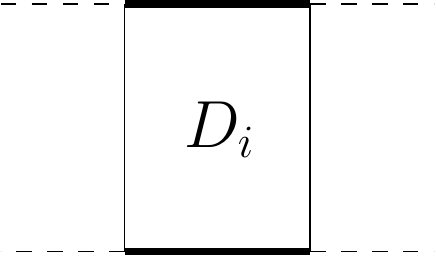} = \includegraphics[scale=0.55,valign=c]{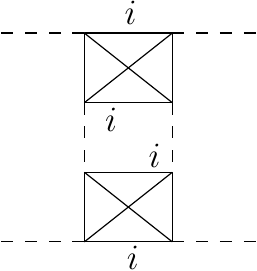} + \includegraphics[scale=0.60,valign=c]{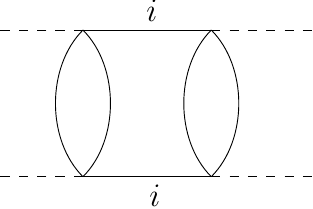}
\label{eq:dip}
\end{equation}
We have added fat edges on two opposite sides of the dipole-vertex in order to have the same symmetry as in the dipoles of type I and II. We can thus remember which external legs come from cutting the same edge in the original elementary melon: they sit on the same side of the fat edges.

\medskip

Since the $O(N)^3$ tensor model is color-symmetric, the generating function of dipoles is independent of its color. Inserting generating functions of melons on one side of the dipoles, the generating function for dipoles reads
\begin{align}
U(t,\mu) &= \bigl(\lambda_1^2 M(t,\mu)^2 + \lambda_2\bigr) M(t,\mu)^2 \\
 &=tM(t,\mu)^4 + \frac{1}{3}t\mu M(t,\mu)^2 \underset{\eqref{eq:mel}}{=} M(t,\mu)- \frac{2}{3}t \mu M(t,\mu)^2 - 1
\label{eq:dip_rec}
\end{align}
In particular, the critical points of $U(t,\mu)$ are the same as the critical points of $M(t,\mu)$.

\vspace{5pt}
\paragraph{Chains\\}
The last family of graphs that will appear in our decomposition are \emph{chains}.
\begin{definition}
A chain is either an isolated dipole or a $4$-point function obtained by connecting an arbitrary number of dipoles by matching one side of a dipole to another side of a distinct dipole.
\label{def:chains}
\end{definition}

The \emph{length} of a chain is defined as the number of dipoles that compose the chain. Note that a chain can have just a single dipole provided it is bubble-disjoint from other dipoles. Crucially, we can check that changing the length of a chain leaves the degree~\ref{eq:deg_O(N)3} unchanged. Notice that a chain of length $\ell$ contains subchains of all lengths $1\leq \ell'\leq\ell$. A chain is said to be \emph{maximal} in a graph $\bar{G}$ if it cannot be included in a longer chain in $\bar{G}$. Two different maximal chains are necessarily bubble-disjoint.
A chain is said to be of color $i$ if it only involves dipoles of color $i$,
\begin{equation}
\label{eq:chain}
\includegraphics[scale=0.45,valign=c]{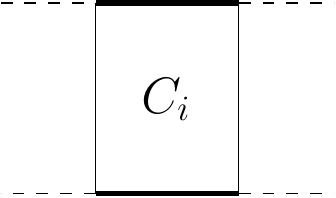} = \includegraphics[scale=0.45,valign=c]{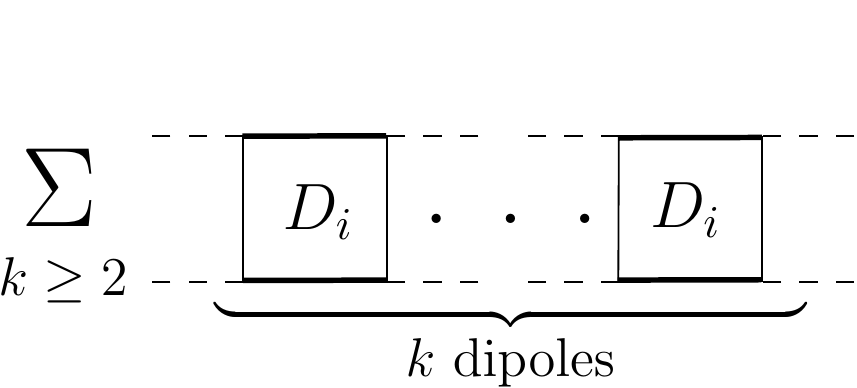}
\end{equation}
Otherwise, it is said to be \textit{broken}. Here we represent a chain of arbitrary length by a \emph{chain-vertex} (on the left hand side). In the following, chains of color $i$ will be denoted by a chain-vertex labeled $C_i$ and broken chains by a chain-vertex labeled $B$. A chain of color $i$ is a sequence of dipoles of length at least $1$, with melons inserted between the dipoles. Those melons have already been inserted on exactly one side of the generating series of the dipole-vertex $U(t,\mu)$. Therefore we have:
\begin{equation}
C_i(t,\mu) = U(t,\mu) \sum_{k\geq 0} U(t,\mu)^k = \frac{U(t,\mu)}{1-U(t,\mu)}
\end{equation}
Broken chains are obtained whenever a sequence of dipoles of any color does not lead to a colored chain. Therefore their generating function is that of all chains minus those of chains of colors $i\in \{1,2,3\}$:
\begin{align}
B(t,\mu) &= 3U(t,\mu) \sum_{k\geq 0} \bigl(3U(t,\mu)\bigr)^k -\sum_{i=1}^{3} C_i(t, \mu) \\
&= \frac{3U(t,\mu)}{1-3U(t,\mu)} - 3 \frac{U(t,\mu)}{1-U(t,\mu)} = \frac{6U(t,\mu)^2}{(1-3U(t,\mu))(1-U(t,\mu))} \nonumber
\end{align}
From this expression, it follows that critical points for chains are either critical points for $M(t,\mu)$, or points where $U(t,\mu) = 1$ or $U(t,\mu) = 1/3$.

\subsection{Schemes of the \texorpdfstring{$O(N)^3$}{O(N)3} model}
\label{ssec:scheme_O(N)^3}

The three families introduced in the last paragraph allow us to introduce the \emph{schemes} of the quartic $O(N)^3$ tensor model.
\begin{definition}
The scheme $\mathcal{S}$ of a 2-point graph $\mathcal{G}$ is obtained by first removing all melonic 2-point subgraphs, then replacing all maximal (broken and not broken) chains with chain-vertices of the same type.
\label{def:scheme}
\end{definition}

The scheme of a graph does not depend on the order of melon removals and is thus uniquely defined. Conversely, a graph $\mathcal{G}$ is uniquely obtained by taking its scheme, first re-extending the chain-vertices as chains and then adding melons. This is precisely the content of Theorem \ref{thm:graph-scheme}. Observe that a scheme has no pillow: every pillow in a graph $\cG$ is an isolated dipole and is, therefore, part of a maximal chain, the latter becoming a chain-vertex in a scheme.

\medskip

In this section, we will prove Theorem~\ref{th:sch} for the quartic $O(N)^3$ tensor model which takes the following form.
\begin{theorem}
    \label{thm:sch_ON3}
    The set of scheme $\mathcal{S}_\omega$ of degree $\omega \in \frac{\mathbb{N}}{2}$ is finite.
\end{theorem}

The result holds thanks to the following two lemmas, which ensure that there are finitely many schemes of any degree $\omega$.
 
\begin{lemma}
\label{lemma:fst}
A scheme of fixed degree $\omega$ has finitely many chain-vertices.
\end{lemma}

\begin{lemma}
\label{lemma:sec}
A vacuum graph $\bar{\cG}\in\bar{\mathbb{G}_{}}$ of degree $\omega$ with $k$ isolated dipoles has a bounded number of bubbles. There exists a constant  $n_{\omega,k}$ such that $n(\bar{\cG})\leq n_{\omega,k}$.
\end{lemma}
%These two lemmas ensures that the set $C_{\omega,min}$ is finite. Since the mapping is injective, it will follow that the set of scheme $S_{\omega}$ is finite as well, thus proving theorem~\ref{th:sch} in the case of the $O(N)^3$ model.

In the proof of Lemma \ref{lemma:fst} and in the analysis of the structure of dominant schemes in the double-scaling limit, a key role is played by its \emph{skeleton graph}.

\begin{figure}[!ht]
\centering
\includegraphics[scale=0.58]{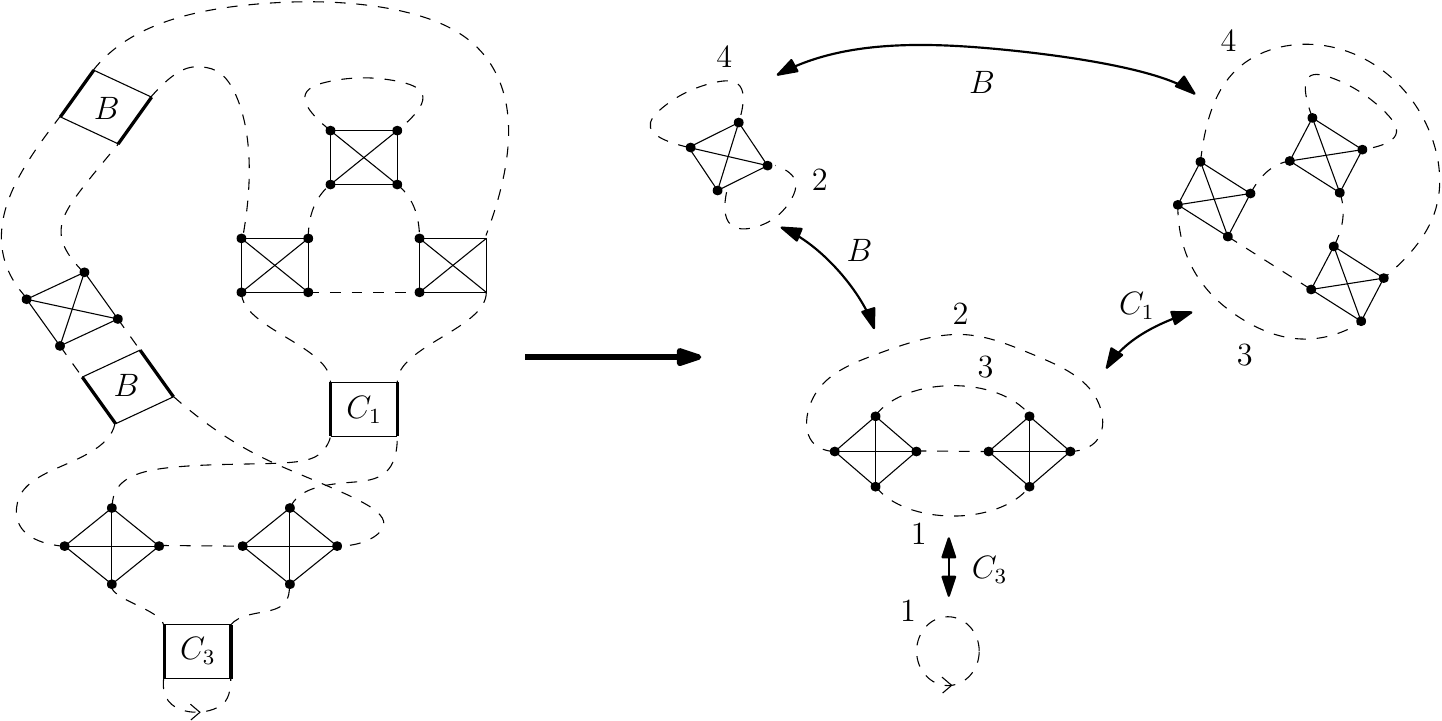}
\caption{A scheme $\cS$ and its skeleton graph $\cI(\cS)$.}
\label{fig:skeleton_graph}
\end{figure}

\begin{definition}
\label{def:skel_graph}
The \emph{skeleton graph} $\cI(\cS)$ of a scheme $\cS$ is the graph such that
\begin{itemize}
\item The vertex set of $\cI(\cS)$ is the set of connected components obtained by removing all chain-vertices of $\cS$.
\item There is an edge between two vertices in $\cI(\cS)$ if the two corresponding connected components are connected by a chain-vertex. This edge is labeled by the type of chain-vertex.
\end{itemize}
\end{definition}
In other words, $\cI(\cS)$ is the incidence graph between chain-vertices, and the disjoint connected components obtained after removing those. An example of a scheme and its skeleton graph is given in Figure~\ref{fig:skeleton_graph}. 

\subsection{Dipole and chain removals}
\label{ssec:dip_chain_rem_ON3}

Since our analysis will rely on skeleton graphs, we have to study how graphs behave under dipole and chain removal. In particular, it is important that those connected components obtained by removing the chain-vertices of a scheme can be seen as graphs from $\mathbb{G}$ (or $\bar{\mathbb{G}}$) and associated with a degree to keep track of the degree through dipole or chain removal. Formally, this operation is defined as follows.

\begin{definition}
A \emph{dipole removal} in $\cG$ consists of removing the dipole and reconnecting the half-edges which were connected on the same side of the dipole together, and similarly for maximal chain removals,
\begin{equation}
\includegraphics[scale=0.65, valign=c]{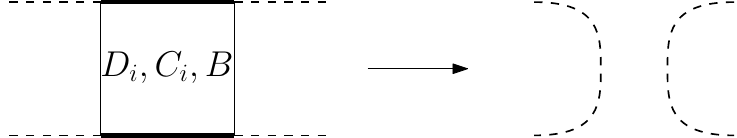}
\label{fig:dip_rem}
\end{equation}
A dipole or a chain is said to be \emph{separating} if its removal disconnects $\cG$, and \emph{non-separating} if it does not. Note that a separating dipole/chain removal in a graph $\mathcal{G}\in \mathbb{G}$ creates two connected components, one being a 2-point graph, the other being a vacuum graph.
\end{definition}

Note that dipole removals can be applied to non-isolated dipoles and not only dipoles that are part of chains.

\medskip

In a scheme, the removal of a chain-vertex does not in general lead to another scheme, because it can create new maximal chains. To circumvent this, a \emph{chain-vertex removal} is defined by replacing the chain-vertex with any chain of the same type, performing the removal, then re-identifying the maximal chains in the newly obtained graph and thus finding its corresponding scheme.

\begin{lemma} \label{thm:ChainRemoval}
Let $S$ be a scheme and consider a chain-vertex removal. If it is a non-separating chain, denote $S'$ the resulting scheme, then
\begin{equation} \label{NonSeparatingRemoval}
\omega(S)-3\leq \omega(S')\leq \omega(S)-1.
\end{equation}
If it is a separating chain or dipole, denote $S_1$ and $\bar{S}_2$ the two resulting schemes, then
\begin{equation} \label{SeparatingRemoval}
\omega(S) = \omega(S_1) + \omega(\bar{S}_2).
\end{equation}
\end{lemma}

\begin{proof}
It is enough to consider the case of graphs, and the case of schemes follows directly. Among all possible types of dipoles and chains, it is enough to consider only dipoles of type I as the two types of dipoles can be interchanged. Similarly, it is sufficient to consider maximal chains of length 2 (to account for broken chains) made of two dipoles of type I since the length of a chain leaves its degree invariant. Therefore we replace the dipole-vertex and chain-vertex as follows
\begin{equation}
\begin{aligned}
&\includegraphics[scale=.5, valign=c]{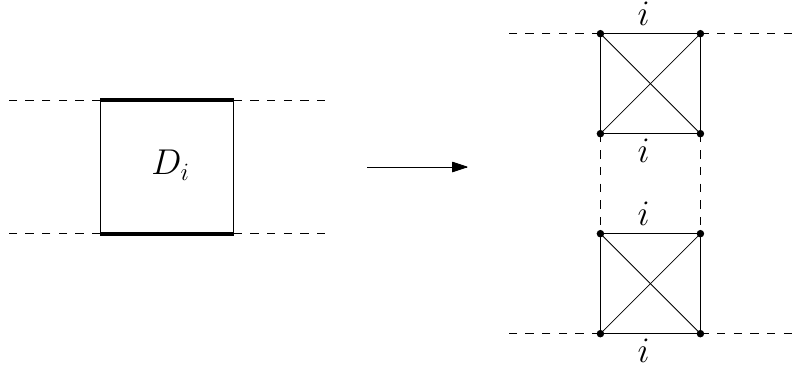}\\
&\includegraphics[scale=.5, valign=c]{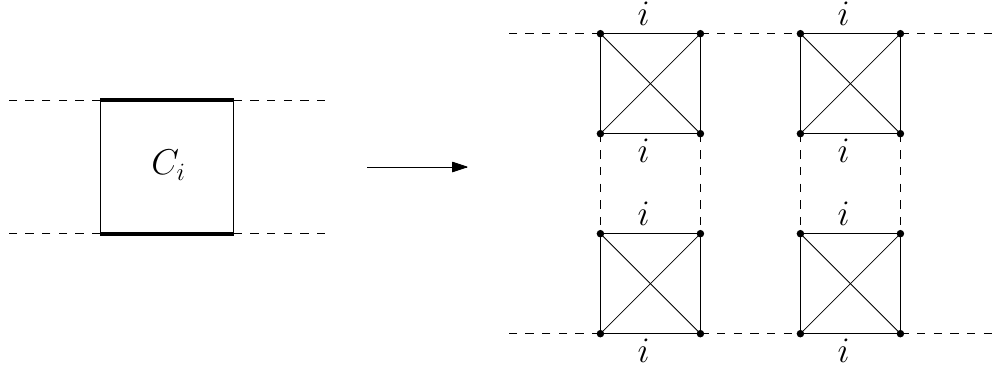}\\
&\includegraphics[scale=.5, valign=c]{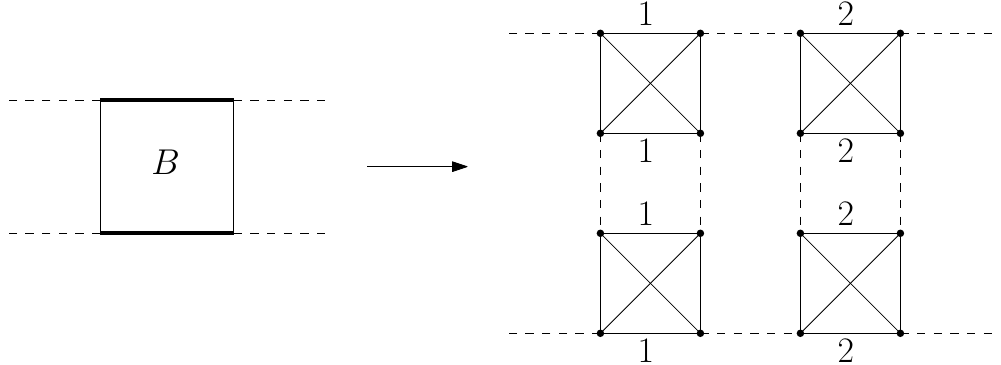}
\end{aligned}
\end{equation}

We now list the possible cases that can occur when removing a dipole or a chain.
\begin{description}
\item[Removal of a non-separating dipole of color $i$]
There are either one or two faces of color $i$ incident to the dipole in $\cG$, whose structures can be as follows (we represent the paths of the faces with dotted edges)
\begin{equation}
\includegraphics[scale=.5, valign=c]{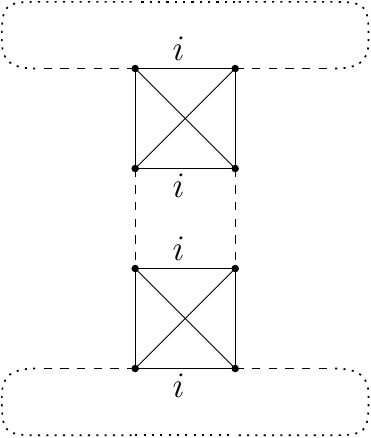} \qquad \includegraphics[scale=.5, valign=c]{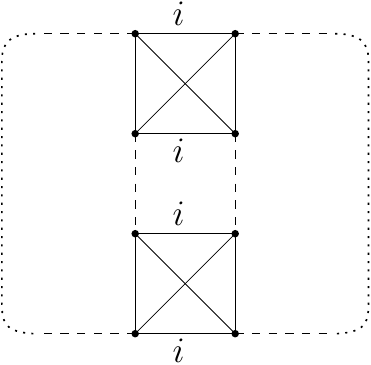}\qquad \includegraphics[scale=.5, valign=c]{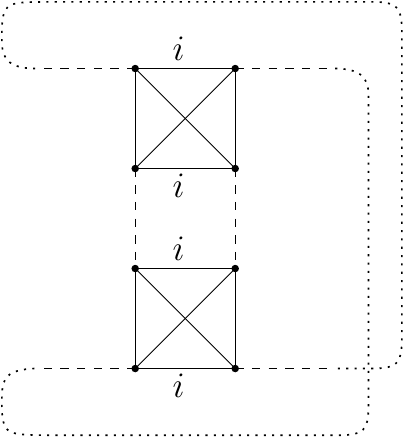}
\end{equation}
and the dipole has an internal face of color $i$. When removing the dipole, the internal face is deleted and one or two faces of color $i$ are formed. The faces of the other colors $j, k\neq i$ are unaffected by the move. The number of bubbles decreases by $2$. Thus deleting a non-separating dipole gives $ -1 \geq \Delta \omega \geq -3$, where $\Delta \omega$ is the variation of the degree.

\item[Removal of a non-separating chain]
We have to distinguish two cases depending on whether the chain is broken or not.
\begin{itemize}
\item If the chain is broken then the structure of the faces incident to the chain is unchanged by the removal. Recall that since $\cG$ is a minimal realization of a scheme, a broken chain has exactly two dipoles of type I, in which case we can check that the chain has exactly three internal faces. Therefore deleting the chain gives $\Delta \omega = -3$.
\item If the chain is not broken, the discussion is similar to the case of the dipole above, and we have $ -1 \geq \Delta \omega \geq -3$.
\end{itemize}

\item[Removal of separating chains and dipoles]
In the following, we tackle the case of separating chains, but the discussion is identical for separating dipoles. If a chain is separating, then the edges of color 0 on either side form a 2-edge-cut,
\begin{equation}
\includegraphics[scale=.5, valign=c]{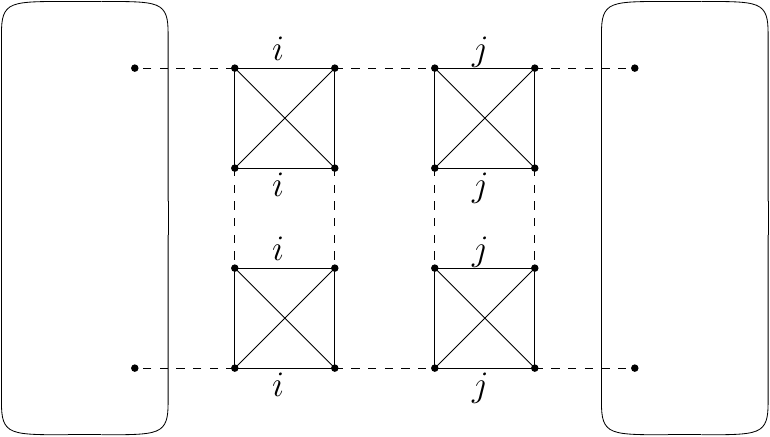}
\end{equation}
with $i,j\in\{1,2,3\}$. The chain removal produces two connected components $\cG_1$ and $\bar{\cG}_2$. It is straightforward to check that in either case $i=j$ and $i\neq j$, $\omega(\cG) = \omega(\cG_1) + \omega(\bar{\cG}_2)$.
\end{description}
\end{proof}

\subsection{Combinatorics of the skeleton graph \texorpdfstring{$\cI(\cS)$}{I(S)}}

To identify the vertex set of the skeleton graph $\cI(\cS)$, one first performs all dipole- and chain-vertex removals in $\cS$. This gives one 2-point graph $\cG^{(0)}$ and a collection of vacuum graphs $\bar{\cG}^{(1)}, \dotsc, \bar{\cG}^{(p)}$, which in turn gives a root vertex and $p$ additional vertices in $\cI(\cS)$.

\begin{lemma} \label{thm:SkeletonGraph}
For any scheme $\cS$, define by convention $\omega(\cI(\cS)) = \omega(\cS)$. The following properties hold.
\begin{enumerate}
\item\label{enum:Valency3} If $\bar{\cG}^{(r)}$ for $r\in\{1, \dotsc, p\}$ has vanishing degree, then the corresponding vertex in $\cI(\cS)$ has valency at least equal to 3.
\item\label{enum:RemoveNonSeparating} Let $\cT\subset \cI(\cS)$ be a spanning tree. Let $q$ be the number of edges of $\cI(\cS)$ which are not in $\cT$. Then, $\omega(\cT) \leq \omega(\cS)-q$.
\item\label{enum:SpanningTree} $\omega(\cT) = \omega(\cG^{(0)}) + \sum_{r=1}^p \omega(\bar{\cG}^{(r)})$, in other words, if $\cI(\cS)$ is a tree, then the degree of $\cS$ is the sum of the degrees of its components obtained by removing all chain-vertices.
\end{enumerate}
\end{lemma}

\begin{proof}
\ref{enum:Valency3}. If $\cG^{(r)}$ has degree zero and the corresponding vertex in $\cI(\cS)$ has valency 1, then $\cG^{(r)}$ is a melonic 2-point function, which cannot happen in schemes. If $\cG^{(r)}$ has degree zero and the corresponding vertex in $\cI(\cS)$ has valency 2, then $\cG^{(r)}$ is a chain, which is also impossible in a scheme as it would mean that a chain-vertex was used in place of a non-maximal chain.
\medskip
\ref{enum:RemoveNonSeparating}. The edges of $\cI(\cS)$ which are not in $\cT$ correspond to non-separating chain-vertices in $\cS$. One then concludes using Equation~\eqref{NonSeparatingRemoval} from Lemma~\ref{thm:ChainRemoval}.
\medskip
\ref{enum:SpanningTree}. Denote $\cS_{\cT}$ the scheme obtained after removing the above $q$ non-separating chain-vertices. Its skeleton graph is $\cT$, i.e. $\cI(\cS_{\cT}) = \cT$, meaning that all its chain-vertices are separating. One concludes with Equation \eqref{SeparatingRemoval} from Lemma \ref{thm:ChainRemoval}.
\end{proof}

%%%%%%%%%%%%%%%%%%%%%
\subsection{Proof of Lemma~\ref{lemma:fst}}

Let $\cS$ be a scheme of degree $\omega>0$ and consider the notations of Lemma~\ref{thm:SkeletonGraph}. Notice that a chain-vertex removal can decrease the number of chain-vertices by at most 3. Indeed, only two edges of color 0 are affected by the removal, which means that at most two pairs of chain-vertices can be joined after the removal to form new chain-vertices. This happens in the following situation,
\begin{equation}
\includegraphics[scale=.6, valign=c]{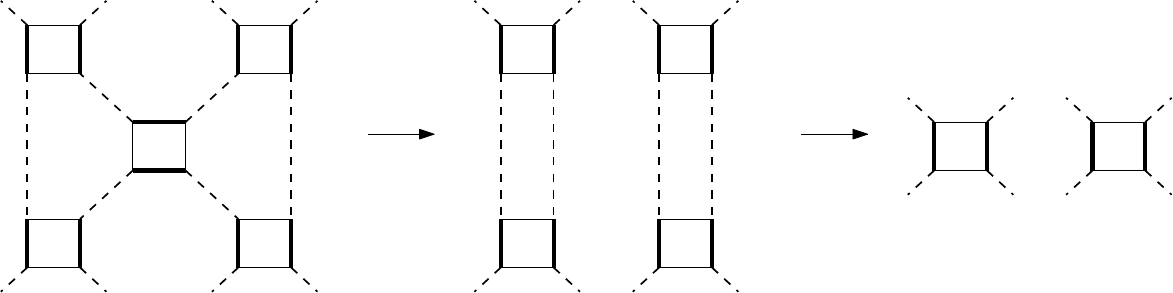}
\end{equation}
After the removal (first arrow), one obtains an object which is not a scheme, because two consecutive chain-vertices would account for two consecutive, hence non-maximal, chains. The correct scheme is obtained by replacing the maximal chains with chain-vertices, represented by the second arrow.

\medskip

Let $N(\cS)$ be the number of chain-vertices of $\cS$ and $N(\cT)$ the number of chain-vertices of the scheme $\cS_{\cT}$ obtained after removing the above $q$ non-separating chain-vertices from $\cS$. One finds $N(\cS) \leq N(\cT) + 3q$. A consequence of point~\ref{enum:RemoveNonSeparating} in Lemma~\ref{thm:SkeletonGraph} is that $q\leq \omega(\cS) - \omega(\cT) \leq \omega(\cS)$ and therefore 
\begin{equation}
\label{BoundSchemeTree}
N(\cS) \leq N(\cT) + 3\omega(\cS).
\end{equation}

For $r\in\{1, \dotsc, p\}$, we say that $\bar{\cG}^{(r)}$ is \emph{tracked} if it is incident to one of the non-separating chain-vertices (corresponding to the edges of $\cI(\cS)$ which are not in $\cT$). Denote $N_0^{t}$ the number of tracked components of degree 0, $N_0^{nt}$ the number of non-tracked components of degree 0, and $N_+$ the number of components of positive degree, all excluding $\cG^{(0)}$. By construction, there are fewer than $2q$ tracked components,
\begin{equation}
\label{eq:bound_C0}
N_0^t\leq 2q.
\end{equation}
Since components of positive degree have degree at least $1/2$, we find from \ref{enum:SpanningTree} in Lemma \ref{thm:SkeletonGraph} that
\begin{equation} \label{eq:bound_pos}
\omega(\cT) \geq \frac{1}{2}N_+
\end{equation}
Together with \ref{enum:RemoveNonSeparating} from Lemma \ref{thm:SkeletonGraph}, and \eqref{eq:bound_C0}, this leads to 
\begin{equation} \label{PartialDegreeBound}
N^t_0 + N_+ \leq 2\omega(\cT)
\end{equation}
Notice that $N(\cT) = p+q$ is also the number of edges of $\cT$. Since it is a tree, the number of edges is the number of vertices minus one, hence the following relation
\begin{equation}
N(\cT) = N^{nt}_0 + N^t_0 + N_+ %\leq N^{nt}_0 + 2\omega(\cT)
\label{eq:tree_euler}
\end{equation}
Finally, counting leaves and nodes of $\mathcal{T}$ weighted by their valency amounts to counting twice the number of edges of $\mathcal{T}$. Due to point \ref{enum:Valency3} in Lemma \ref{thm:SkeletonGraph}, we have
\begin{equation}
2N(\cT) \geq 3N^{nt}_0 + N^t_0 + N_+ + 1
\label{eq:tree_val}
\end{equation}
the additional one being due to the component $\cG^{(0)}$ (with valency at least 1). The previous two equations lead to $N(\cT) \leq 2(N_0^t+N_+)$. From~\eqref{BoundSchemeTree} together with~\eqref{PartialDegreeBound}, $N(\cT) \leq 4\omega(\cT)$, one gets
\begin{equation}
N(\cS) \leq 4\omega(\cT) + 3\omega(\cS) - 1 \leq 7\omega(\cS) -1,
\end{equation}
which proves Lemma~\ref{lemma:fst}.

\subsection{Proof of Lemma~\ref{lemma:sec}}
\label{ssec:proof_Lem_sec}

Here we prove Lemma~\ref{lemma:sec}, i.e. graphs of finite degree and with a finite number of isolated dipoles have a bounded number of bubbles. Our strategy is to show that there exist bounds on the number of faces of every degree, depending on the degree $\omega$ and the number of isolated dipoles $k$, i.e.
\begin{equation} \label{BoundFaceDegreeP}
F^{(p)}(\mathcal{G}) \leq \phi^{(p)}(\omega, k)
\end{equation}
where $F^{(p)}(\mathcal{G})$ is the number of faces of degree $p$.

\medskip

There is exactly one face of color $i=1,2,3$ along each edge of color $0$. Therefore, denoting $F_i^{(p)}$ the number of faces with color $i$ and degree $p$, we have for each color $i$
\begin{equation}
\sum\limits_{p \geq 1} pF_i^{(p)}(\mathcal{G}) = E_0(\cG) = 2n(\cG)
\end{equation}
Using the degree formula~\eqref{eq:deg_O(N)3}, $n(\cG)$ can be eliminated to give
\begin{equation} 
\sum_{p \geq 5} (p-4)F^{(p)}(\mathcal{G}) = 4(\omega-3) + 3 F^{(1)}(\mathcal{G}) + 2 F^{(2)}(\mathcal{G}) + F^{(3)}(\mathcal{G})
\label{eq:face_bound}
\end{equation}
We have written this equation so that both sides come with positive coefficients (except the irrelevant constant $-12$). Therefore if we can bound $F^{(1)}(\mathcal{G})$, $F^{(2)}(\mathcal{G})$, $F^{(3)}(\mathcal{G})$, i.e. if we can prove \eqref{BoundFaceDegreeP} for $p=1,2,3$, then it automatically gives the bound \eqref{BoundFaceDegreeP} for $p\geq 5$. It then remains to prove the bound independently for $p=4$, since $F^{(4)}$ does not appear in \eqref{eq:face_bound}. There is indeed a risk that graphs with an arbitrarily large number of faces of degree 4 exist, without having these faces affecting the degree. However, this is not the case, as we will prove. In our analysis, we will distinguish self-intersecting and non-self-intersecting faces. A face is said to be \emph{self-intersecting} if it visits the same bubble more than once (that is twice since our interactions are quadratic and thus only have $2$ edges in each color). Otherwise it is called \emph{non-self-intersecting} (n.s.i.).

\paragraph{Faces of degree $1$\\}
If a graph has a face of degree $1$, then it has the following structure:
\begin{equation}
\begin{array}{c} \includegraphics[scale=.3]{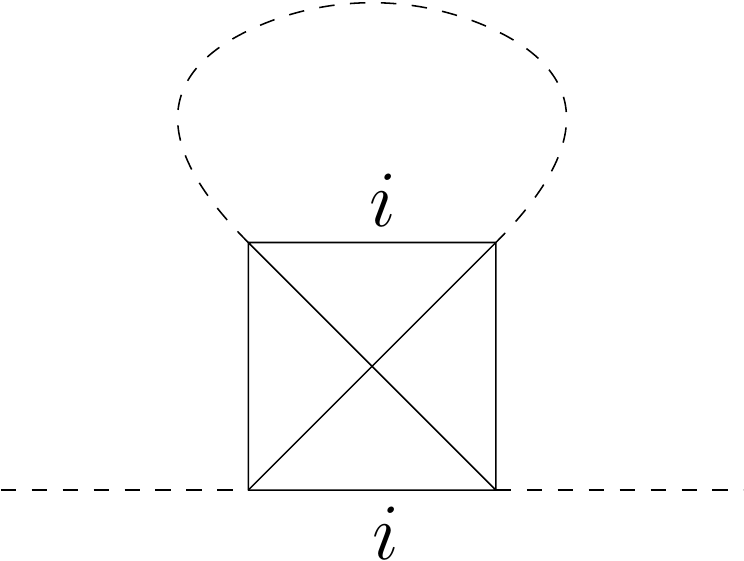} \end{array}
\end{equation}
When that bubble and the tadpole are removed and replaced with an edge of color 0, the degree decreases by $\frac{1}{2}$. Therefore, there are at most $2\omega$ tadpoles in a graph $\cC$ of degree $\omega$ i.e. $F^{(1)}(\cG)\leq 2\omega$.

\begin{figure}
    \centering
    \includegraphics[scale=0.4]{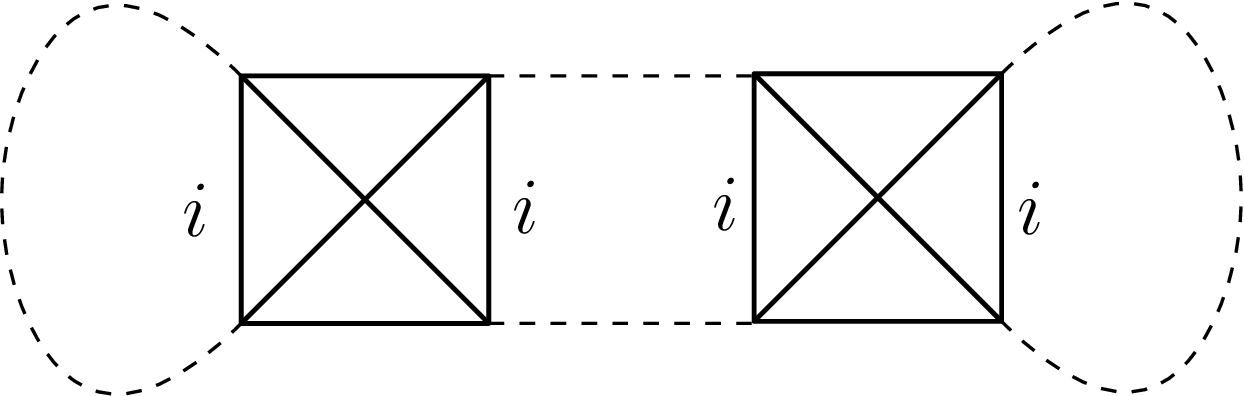}
    \caption{The unique graph with a self-intersecting face of length $2$.}
    \label{fig:double_tadpole}
\end{figure}

\paragraph{Faces of degree $2$\\}
First notice that there is a single graph that has a face of degree 2 which is self-intersecting, represented in Figure~\ref{fig:double_tadpole}.

If $\cG$ has a face of degree $2$ which is n.s.i., then it is a dipole. Either it is a dipole that is isolated and there are $k$ of them, or it is a non-isolated dipole. In the latter case, we perform the dipole removal,
\begin{equation} \label{NonIsolatedDipoleRemoval}
    \begin{array}{c} \includegraphics[scale=.45]{NonIsolatedDipoles.pdf}\end{array} \quad \to \quad \begin{array}{c} \includegraphics[scale=.45]{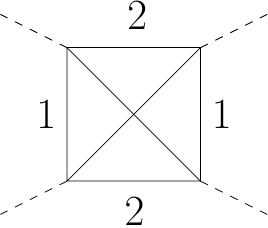}\end{array}
\end{equation}
A non-isolated dipole is always non-separating, therefore one obtains a connected graph $\cG'$ whose degree satisfies $\omega(\cG')<\omega(\cG)$. We thus proceed by induction on the degree.

\medskip

Denote $F^{(2)}_{\text{NI}}(\cG)$ the number of non-isolated dipoles of $\cG$. If $\cG$ is melonic, then $F^{(2)}_{\text{NI}}(\cG)=0$. Let $\omega>0$ and assume that for all $\cG'$ of degree $\omega'<\omega$ with $k'$ isolated dipoles there is a bound $F^{(2)}_{\text{NI}}(\cG) \leq \phi^{(2)}(\omega',k')$. Let $\cG$ of degree $\omega$ with $k$ isolated dipoles. We perform the move~\eqref{NonIsolatedDipoleRemoval} and obtain $\cG'$ of degree $\omega'<\omega$. We track the changes in the number of dipoles:
\begin{itemize}
    \item The number of non-isolated dipoles cannot increase (but it may remain unchanged).
    \item Moreover the bubble on the RHS of~\eqref{NonIsolatedDipoleRemoval} can belong to at most one isolated dipole, therefore the number of isolated dipoles of $\cG$ satisfies $k'\leq k+1$. 
\end{itemize}
From the induction hypothesis, we thus find
\begin{equation}
    F^{(2)}_{\text{NI}}(\cG) \leq F^{(2)}_{\text{NI}}(\cG') \leq \max_{\substack{\omega'<\omega\\ k'\leq k+1}} \phi^{(2)}(\omega',k').
\end{equation}
The RHS defines the bound $\phi^{(2)}(\omega,k)$.

\paragraph{Faces of degree $3$.\\}  We distinguish cases depending whether $\cG$ has a face of degree $3$ which is self-intersecting or not.

\subparagraph{Self-intersecting faces of degree 3\\}
If a face of degree $3$ is self-intersecting, then it necessarily has the following structure,
\begin{equation}
\begin{array}{c} \includegraphics[scale=0.45]{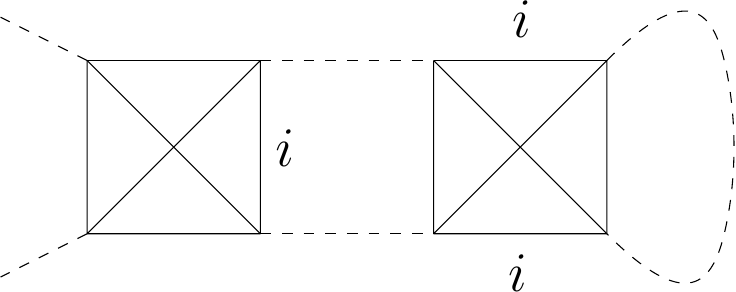} \end{array}
\end{equation}
It is a 2-point function whose removal decreases the degree by 1. Thus there are at most $\omega$ such faces in the graph, $F^{(3)}_{\text{s.i.}}(\mathcal{G})\leq \omega$.

\subparagraph{Non-self-intersecting faces of length $3$\\}
A n.s.i. face of degree 3 and color $1$ necessarily has the following structure,
\begin{equation} \label{F3nsi}
\begin{array}{c} \includegraphics[scale=0.4]{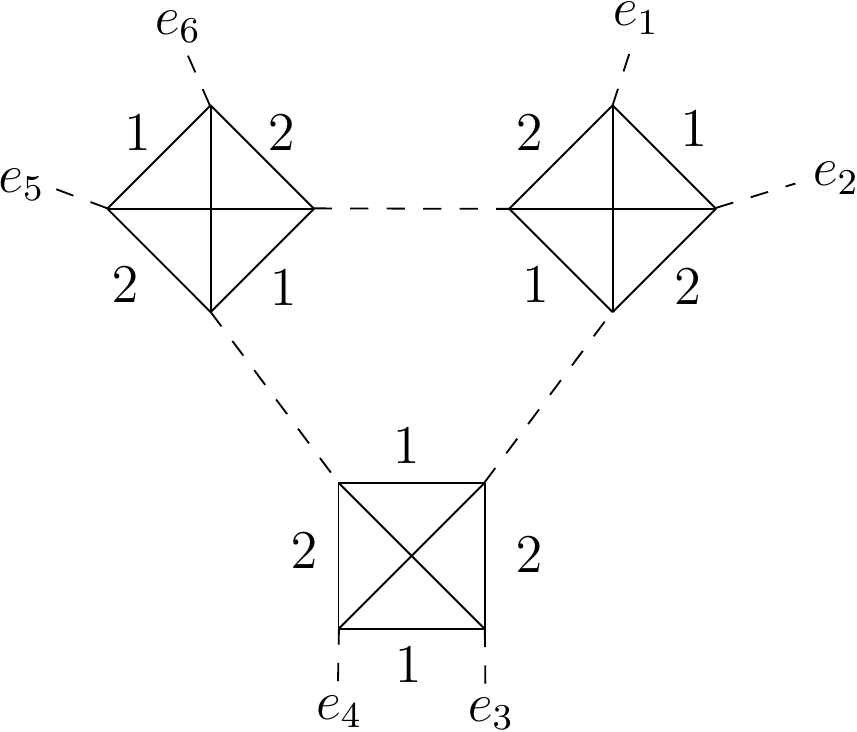} \end{array}
\end{equation}
There are at most $3$ different faces of each color passing along the edges $e_1, \dotsc, e_6$. We denote $f_1, f_2, f_3 \leq 3$ their numbers. Note that $2$ and $3$ play symmetric roles. We want to show that it is always possible to remove this subgraph from $\cG$ and reglue the $6$ half-edges while decreasing the degree of the graph. We will use the following lemmas.

\begin{lemma}
\label{thm:Variations}
The following two assertions hold.
\begin{enumerate}
\item Consider a n.s.i. face of degree 3 in $\cG\in\mathbb{G}_{O(N)^3}$, as in \eqref{F3nsi}, and a move on $\cG$ which removes the 3 bubbles of that face and reconnects the half-edges $e_1, \dotsc, e_6$ pairwise (the resulting graph may not be connected). Denote $\Delta F^{(d)}$ the variation of the number of faces of degree $d$ and $\Delta k$ the variation of the number of dipoles. Then
\begin{equation}
\Delta F^{(d)} \geq -9 \qquad \text{and} \qquad \Delta k \leq 9.
\end{equation}
\item Let $\cG\in\mathbb{G}_{O(N)^3}$ with $k$ dipoles. Let $e, e'$ be two edges of color 0 in $\cG$ forming a 2-edge-cut, and perform a flip
\begin{equation}
\begin{array}{c} \includegraphics[scale=.3]{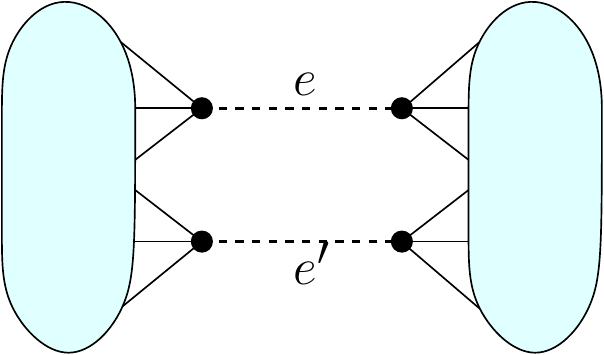} \end{array} \qquad \to \qquad \mathcal{G}_L = \begin{array}{c} \includegraphics[scale=.3]{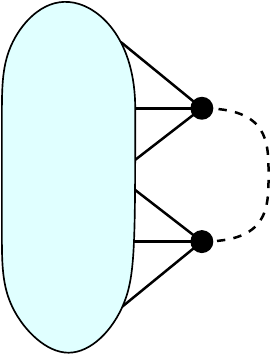} \end{array} \quad \cup \quad \mathcal{G}_R = \begin{array}{c} \includegraphics[scale=.3]{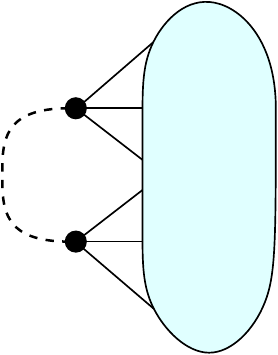} \end{array}
\end{equation}
which gives rise to two connected components $\cG_L$, $\cG_R$. Then the number of dipoles of each is bounded as $k(\cG_\alpha) \leq k+3$, $\alpha=L, R$, and the number of faces of degree $d$ as $F^{(d)}(\cG_L) + F^{(d)}(\cG_R) \geq F^{(d)}(\cG) - 3$.
\end{enumerate} 
\end{lemma}

This lemma contains only particular cases of a more general principle. If $\mathcal{H}$ is a subgraph, then the number of incident faces of each degree can be bounded as a function of $\mathcal{H}$ and not the graphs it is contained in. Then when one removes $\mathcal{H}\subset \cG$ in some way, the variations of the number of faces of each degree is bounded independently of $\cG$.

\begin{proof}
Along each edges there are exactly 3 faces, one of each color.
\begin{enumerate}
\item There are at most 9 different faces going along the edges $e_1, \dotsc, e_6$ before the move.. If they all have degree $d$ and none of this degree are created by the move we have $\Delta F^{(d)} =-9$, which is the extremal case.

After the move, the pairing of the half-edges gives rise to 3 edges of color 0. If they have no faces in common and all of them are dipoles, this gives 9 dipoles. If no dipoles are destroyed in the move, this gives $\Delta k = 9$.

\item For each color 1, 2, 3, it is the same face going along $e$ and $e'$ and the move splits each of them into two. If those three faces in $\cG$ were of degree $d$, and they are split in $\cG_L, \cG_R$ into faces of different degrees, then $F^{(d)}(\cG_L) + F^{(d)}(\cG_R) = F^{(d)}(\cG) - 3$ which is the extremal case.

\end{enumerate}
\end{proof}

\begin{lemma} \label{thm:TypeIMove}
Assume $\cG$ has a n.s.i. of degree 3 such that $f_2\leq 3$ and $f_3\leq 2$. Further assume that the following move gives a graph $\cG'$ connected
\begin{equation}
\begin{array}{c} \includegraphics[scale=.75]{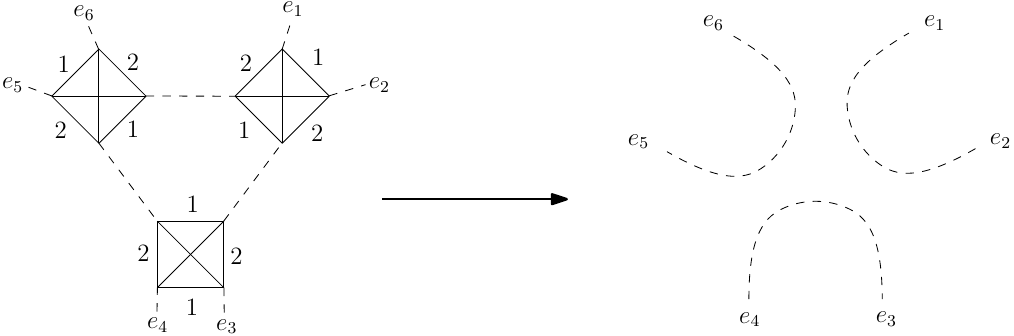} \end{array}
\end{equation}
then the degree changes as $\Delta \omega = \omega(\cG') - \omega(\cG) \leq -1/2$.
\end{lemma}

\begin{proof} 
The variation in the degree after the move is given by 
\begin{equation}
\Delta \omega = -\frac{7}{2} - \Delta f_1 - \Delta f_2 - \Delta f_3
\end{equation}
where $\Delta f_{i}$ is the variation of the number of faces of color $i$ which go along $e_1, \dotsc, e_6$ in $\cG$. In the lemma, the move does not change $f_1$. Moreover, $\Delta f_2\geq -2$ and $\Delta f_3 \geq -1$, since there is at least one face of color $2$ and one of color $3$ going along the edges on the configuration on the right. Thus $\Delta \omega \leq -\frac{1}{2}$.
\end{proof}

We can now prove \eqref{BoundFaceDegreeP} for n.s.i. faces of degree 3. We proceed by induction on the degree. At degree 0, all faces have even degree so $F^{(3)}(\cG)=0$.

\medskip

Let $\omega>0$ and assume that there exists a bound $F^{(3)}_{\text{n.s.i.}}(\cG')\leq \phi^{(3)}_{\text{n.s.i.}}(\omega', k')$ for all graphs $\cG'$ of degree $\omega'<\omega$. Let $\cG\in\mathbb{G}_{O(N)^3}$ have degree $\omega$.

\medskip

We first consider the cases where a pair of edges $e_i, e_j$ forms a 2-edge-cut. Without loss of generality, we will consider the pairs to be $\{e_1, e_j\}$ for $j=2, \dotsc, 6$. Clearly, the case of $\{e_1, e_5\}$ is a 2-cut is equivalent to $\{e_1, e_3\}$ by exchanging the colors 2 and 3. Same for $\{e_1, e_6\}$ which is equivalent to $\{e_1, e_4\}$. It is therefore enough to consider the cases where $\{e_1, e_j\}$ is a 2-cut for $j=2,3,4$.

\medskip

The first step is to perform the cut and obtain two connected components $\cG_1, \cG_2$ and $\omega(\cG) = \omega(\cG_1) + \omega(\cG_2)$. We denote $\cG_1$ the component which inherits the n.s.i. face of degree 3. From Lemma \ref{thm:Variations}
\begin{equation}
F^{(3)}_{\text{n.s.i.}}(\cG) \leq  F^{(3)}_{\text{n.s.i.}}(\cG_1) + F^{(3)}_{\text{n.s.i.}}(\cG_2) + 3
\end{equation}
$\cG_1$ is not melonic (it has a face of degree 3), so $\omega(\cG_1)>0$. This gives $\omega(\cG_2)<\omega(\cG)$ and thus by the induction hypothesis 
\begin{equation}
F^{(3)}_{\text{n.s.i.}}(\cG_2) \leq \phi^{(3)}_{\text{n.s.i.}}(\omega(\cG_2), k(\cG_2)) \leq \max_{k_2\leq k+3} \phi^{(3)}_{\text{n.s.i.}}(\omega(\cG_2), k_2)
\end{equation}
To bound $F^{(3)}_{\text{n.s.i.}}(\cG_1)$, we will notice from Lemma~\ref{thm:Variations} that $\cG_1$ has at most $k+3$ dipoles (hence isolated dipoles) and distinguish the following cases.

\begin{description}
\item[$\{e_1, e_2\}$ is a 2-cut (or form a single edge)] 
$\cG_1$ has a tadpole, 
\begin{equation}
\begin{array}{c} \includegraphics[scale=.6]{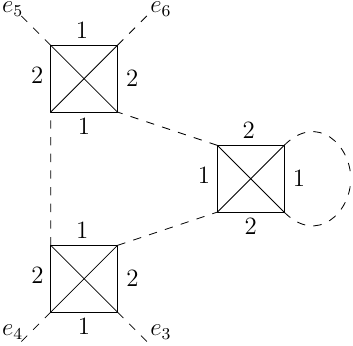} \end{array}
\end{equation}
which we remove to get $\cG_1'$ whose degree is $\omega(\cG_1') = \omega(\cG_1)-1/2$. From Lemma \ref{thm:Variations}, $\cG_1'$ has at most $k+6$ isolated dipoles ($k+3$ in fact holds because in this particular case, $\cG_1$ has at most all the dipoles of $\cG$). Thus,
\begin{equation}
F^{(3)}_{\text{n.s.i.}}(\cG_1) \leq F^{(3)}_{\text{n.s.i.}}(\cG_1') + 3 \leq \max_{k_1\leq k+6} \phi^{(3)}_{n.s.i.}(\omega-\frac{1}{2}, k_1) + 3
\end{equation}

\item[$\{e_1, e_3\}$ is a 2-cut (or form a single edge)] 
$\cG_1$ is as follows
\begin{equation}
\begin{array}{c} \includegraphics[scale=.6]{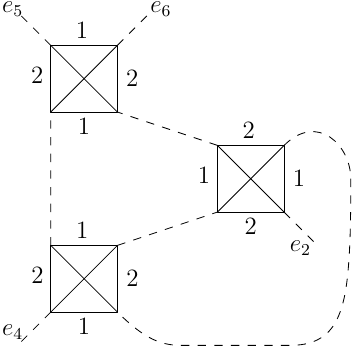} \end{array}
\end{equation}
so that $f_2, f_3\leq 2$. We can assume that $\{e_5, e_6\}$ is not a 2-edge-cut, since it would be equivalent to the case where $\{e_1, e_2\}$ is a 2-cut if it were. Therefore the move from Lemma \ref{thm:TypeIMove} does not disconnect $\cG_1$ and one obtains $\cG_1'$ of degree $\omega(\cG_1')\leq \omega(\cG_1)-1/2$. From Lemma \ref{thm:Variations}, it has at most 9 more isolated dipoles than $\cG_1$ and 9 n.s.i. faces of degree 3 less than $\cG_1$. This gives
\begin{equation}
F^{(3)}_{\text{n.s.i.}}(\cG_1) \leq F^{(3)}_{\text{n.s.i.}}(\cG_1') + 9 \leq \max_{\substack{\omega_1\leq \omega-1/2\\ k_1\leq k+12}} \phi^{(3)}_{\text{n.s.i.}}(\omega_1, k_1) + 9
\end{equation}

\item[$\{e_1, e_4\}$ is a 2-cut (or form a single edge)] 
$\cG_1$ is as follows
\begin{equation}
\begin{array}{c} \includegraphics[scale=.6]{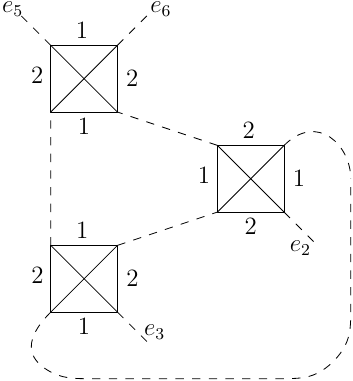} \end{array}
\end{equation}
so that $f_2\leq 2, f_3\leq 3$. The rest is therefore identical to the case where $\{e_1, e_3\}$ is a 2-cut above.
\end{description}

We now assume that no pair $\{e_i, e_j\}$ forms a 2-cut and distinguish two cases.
\begin{description}
\item[$f_2 \leq 3$ and $f_3\leq 2$ (or $f_3 \leq 3$ and $f_2\leq 2$)] 
We can apply the move of Lemma \ref{thm:TypeIMove}. The resulting graph $\cG'$ is connected and of degree $\omega(\cG')\leq \omega-1/2$. From Lemma \ref{thm:Variations} we further get
\begin{equation}
F^{(3)}_{\text{n.s.i.}}(\cG) \leq F^{(3)}_{\text{n.s.i.}}(\cG') + 9 \leq \max_{\substack{\omega'\leq \omega-1/2\\ k'\leq k+9}} \phi^{(3)}_{\text{n.s.i.}}(\omega', k') + 9
\end{equation}

\item[$f_2 = f_3 = 3$]
In this case, the move from Lemma \ref{thm:TypeIMove} would only guarantee $\Delta \omega \leq 1/2$ which is not enough for our purpose. Instead, we exchange some of the half-edges as follows
\begin{equation}
\begin{array}{c} \includegraphics[scale=0.65]{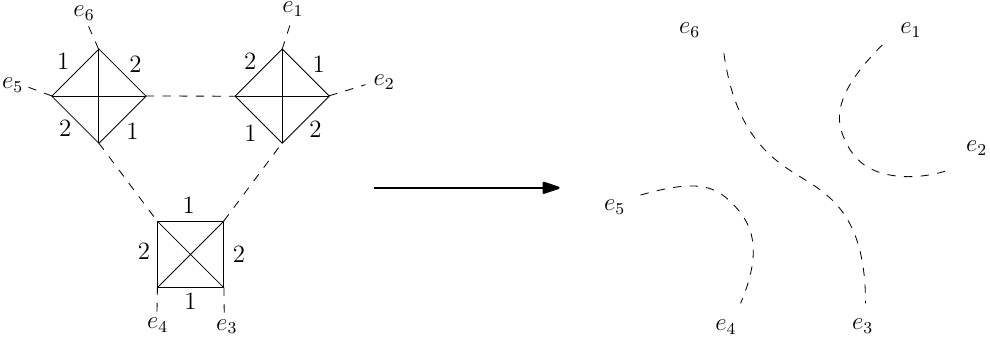} \end{array}
\end{equation}
which gives exactly $\Delta f_2 = \Delta f_3 = -1$. Moreover, if $f_1=3$, then $\Delta f_1=-1$ too, and if $f_1\geq 2$ then $\Delta f_1\geq -1$ because there is still at least one face on the configuration of the right hand side. Overall, $\Delta f_1 \geq -1$ and thus $\Delta \omega \leq -1/2$. We then conclude as previously.
\end{description}

\paragraph{Faces of degree $4$\\}

We have so far proved that at fixed degree $\omega$ and fixed number of isolated dipoles $k$, there is a bound on the number of faces of every degree $d$, except for $d=4$ which does not enter \eqref{eq:face_bound}. So {\it a priori}, there could be an unbounded number of faces of degree 4. However, for this to be possible they would have to be at an arbitrarily large distance of any faces of degree $d\neq 4$ (since a bubble has a finite number of other bubbles at a finite distance). Here the distance between two bubbles is the minimal number of edges of color 0 to go from any vertex of the first bubble to any vertex of the second bubble.

\medskip

Similarly as for faces of degrees $2$ and $3$, we distinguish two cases according to whether the face is self-intersecting or not.

\subparagraph{Self-intersecting faces of degree $4$\\}

Since a tetrahedral bubble has two edges of each color, a face can pass through a bubble at most twice. Therefore a self-intersecting face of degree $4$  has either 2 or 3 bubbles. If it is 2, then $\cG\in\mathbb{G}_{O(N)^3}$ has in fact only those two bubbles and there is a finite number of such graphs.

\medskip

A self-intersecting face of degree 4 which goes along three bubbles can have a subgraph which contain a subgraph with a face of length $1$ as follows (up to color permutations),
\begin{equation}
\begin{array}{c} \includegraphics[scale=.6]{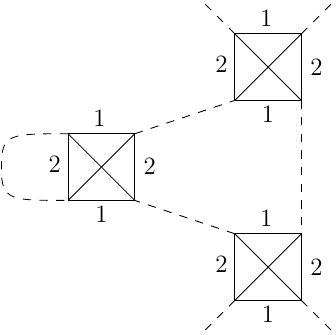} \end{array}
\end{equation}
the number of which we know to be bounded at fixed degree, or be as follows (up to color permutations)
\begin{equation}
\begin{array}{c} \includegraphics[scale=.45]{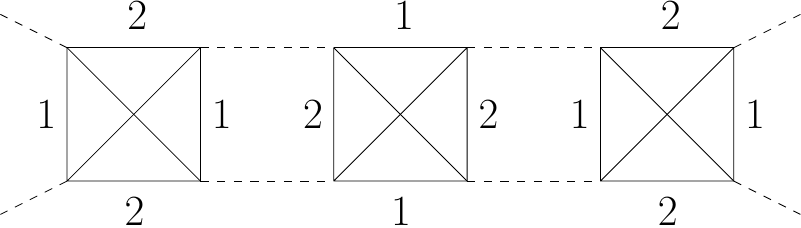} \end{array}
\end{equation}

\begin{itemize}
\item If the graph has a 2-cut as follows
\begin{equation}
\begin{array}{c} \includegraphics[scale=.45]{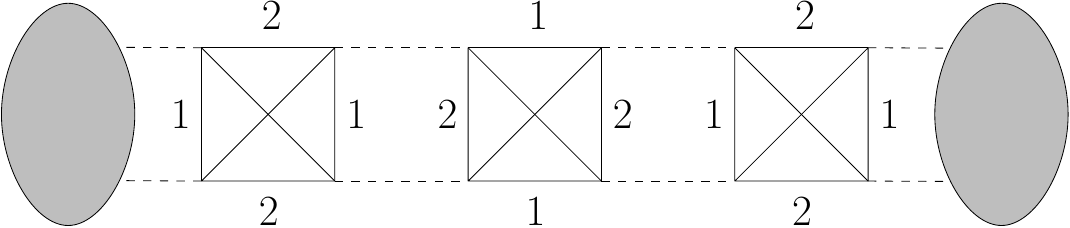} \end{array}
\end{equation}
we transform it to $\begin{array}{c} \includegraphics[scale=.25]{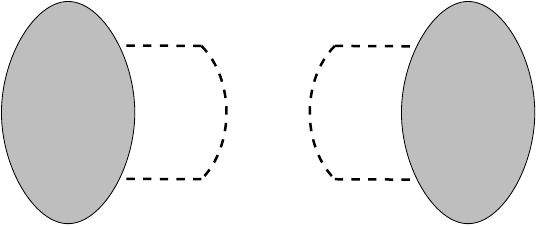} \end{array}$ and denote $\cG_L, \cG_R$ the two connected components. We find that
\begin{equation}
\omega(\cG) = \omega(\cG_L) + \omega(\cG_R) + \frac{3}{2}
\end{equation}
which means that $\omega(\cG_{L,R}) <\omega(\cG)$ so that the induction hypothesis can be applied to $\cG_{L, R}$. By applying Lemma \ref{thm:Variations} on both sides of the subgraph we find
\begin{equation}
F^{(4)}(\cG) \leq F^{(4)}(\cG_L) + F^{(4)}(\cG_R) + 7 \leq 2 \max_{\substack{\omega'\leq \omega-3/2\\ k'\leq k+3}} \phi^{(4)}(\omega',k') + 7.
\end{equation}

\item Else, we perform the same move but the new graph $\cG'$ is connected and we have $\omega(\cG') \leq \omega(\cG)-5/2$. By adapting the arguments of Lemma \ref{thm:Variations}, we also find $F^{(4)}(\cG') \geq F^{(4)}(\cG) - 10$ and $k(\cG') \leq k(\cG)+10$. This leads to
\begin{equation}
F^{(4)}(\cG) \leq \max_{\substack{\omega'\leq \omega-5/2\\k'\leq k+10}} \phi^{(4)}(\omega', k') + 10.
\end{equation}
\end{itemize} 

\subparagraph{Non-self-intersecting faces of degree $4$\\}
\label{par:nsi_f4}

Let $\cG\in\mathbb{G}_{O(N)^3}$ have a n.s.i. face of degree 4 and $B_{\text{exc}}\subset \cG$ the set of bubbles incident to at least one face of degree $d\neq 4$ or one self-intersecting face of degree 4. From the previous cases, we have shown that there is a bound $\beta(\omega, k)$ such that
\begin{equation}
|B_{\text{exc}}| \leq \beta(\omega, k).
\end{equation}
Now, we want to show that $\cG$ cannot have bubbles at an arbitrarily large distance from $B_{\text{exc}}$. Since the number of bubbles at a finite distance of $B_{\text{exc}}$ is itself finite (because bubbles have degree 4), we obtain that $\cG$ have a bounded number of bubbles which proves Lemma~\ref{lemma:sec}. To do so, we will show that there are only a finite number of graphs where a bubble is at a distance 3 or more from $B_{\text{exc}}$.

\medskip 

We denote $A, B, C, D$ the four bubbles of a n.s.i. face of degree 4
\begin{equation}
\begin{array}{c} \includegraphics[scale=.74]{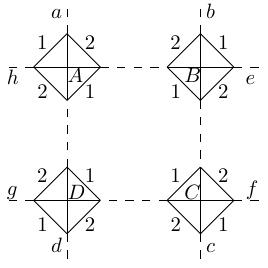} \end{array}
\end{equation}
and consider that \emph{A is at distance 3 of $B_{\text{exc}}$}. In particular, it is incident to only n.s.i. faces of degree 4. The half-edge $a$ cannot be connected to $b$ as it would form a dipole, nor to $d$ for the same reason. If $a$ is connected to $c$, then, in order for the faces of colors 2 and 3 to have degree 4, it is necessary to connect $b$ to $d$. This leaves a subgraph which we replace with a single bubble as follows
\begin{equation}
\begin{array}{c} \includegraphics[scale=.75]{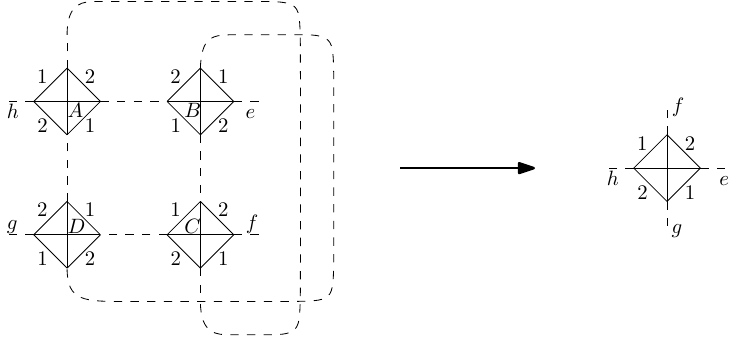} \end{array}
\end{equation}
This does not change any of the faces that go along $e, f, g, h$. The new graph $\cG'$ is connected, has only three faces less than $\cG$ and has degree $\omega(\cG') = \omega(\cG) -3/2$. Again adapting the arguments of Lemma \ref{thm:Variations}, we see that the number of faces of degree 4 cannot decrease by more than 6 (this is the largest number of faces which can go through $e,f,g,h$), the number of dipoles cannot increase by more than 6 (for the same reason). This gives
\begin{equation}
F^{(4)}(\cG) \leq \max_{k'\leq k+6} \phi^{(4)}(\omega-\frac{3}{2}, k') + 6
\end{equation}

We can now consider the case where $a$ is connected to another bubble,
\begin{equation}
\begin{array}{c} \includegraphics[scale=.7]{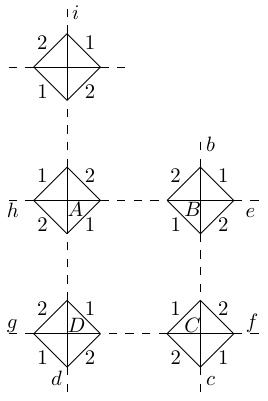} \end{array}
\end{equation}
The face of color 3 which goes along $d$ and $i$ must be of degree 4. The half-edges $d, i$ must therefore be connected to two vertices which are themselves connected by an edge of color 3. There is no such edge available in the subgraph, so a new bubble must be added,
\begin{equation}
\begin{array}{c} \includegraphics[scale=.7]{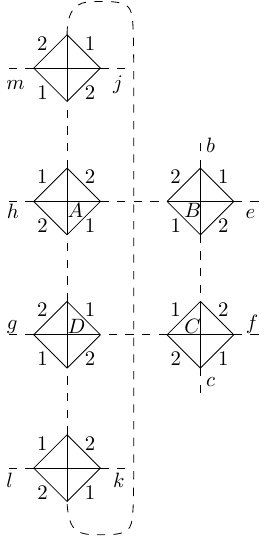} \end{array}
\end{equation}
The face of color 2 which goes along $b$ and $j$ must be of degree 4. The half-edges $b, j$ must therefore be connected to two vertices with an edge of color 2 between them. There are no such vertices available in the subgraph, so a new bubble must be added,
\begin{equation}
\begin{array}{c} \includegraphics[scale=.7]{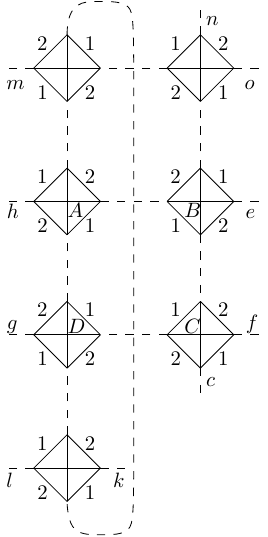} \end{array}
\end{equation}
The face of color 3 which goes along $c$ and $n$ must be of degree 4. The half-edges $c, n$ must therefore be connected to two vertices which are themselves connected by an edge of color 3. There is in our subgraph the edge of color 3 adjacent to the half-edges $k, l$ available to do so.
\begin{itemize}
\item If $c$ is connected to $k$ and $n$ to $l$, it creates a face of degree 3 (that of color 2 along $c$ and $k$), which is forbidden.
\item If $c$ is connected to $l$ and $n$ to $k$, it creates a face of degree greater than 4 (that of color 2 along $c$), which is forbidden.
\end{itemize}
Therefore another bubble must be added
\begin{equation}
\begin{array}{c} \includegraphics[scale=.7]{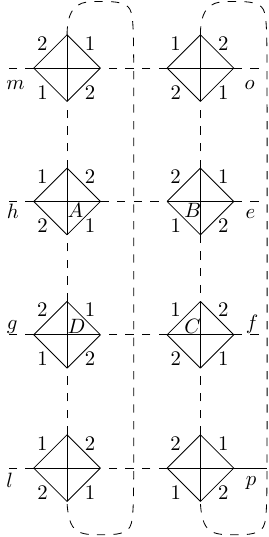} \end{array}
\end{equation}
We now consider the face of color 3 along $h$.
\begin{itemize}
\item If $h$ is connected to $g$, $m$ or $e$, this creates a dipole, which is forbidden.
\item $h$ can be connected to $f$, $l$, $m$, or $p$.
\item $h$ can be to a new bubble.
\end{itemize}

The situations of the second type are all treated similarly and each lead to a single possible graph.
\begin{description}
\item[$h$ to $l$] It then forces $e$ to $p$ so that the face of color 3 along $h$ has degree 4. It then forces $g$ to $m$ so that the face of color 1 along $h$ has degree 4. This in turn forces $f$ to $o$ so that the face of color 3 along $g$ has degree 4 (indeed the bubble labeled $D$ is at distance at least 2 of $B_{\text{exc}}$ so all its incident faces must have degree 4). This fully determines $\cG\in\mathbb{G}_{O(N)^3}$,
\begin{equation}
\begin{array}{c} \includegraphics[scale=.7]{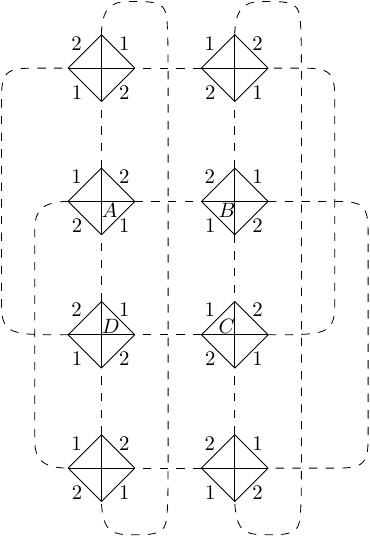} \end{array}
\end{equation}
\item[$h$ to $f$] then forces $g$ to $e$, then $p$ to $m$, then $o$ to $l$.
\item[$h$ to $o$] similar to the previous case by symmetry.
\item[$h$ to $p$] then forces $l$ to $e$, then $f$ to $m$, then $o$ to $g$.
\end{description}

We now consider the case where $h$ is connected to a new bubble. To close the face of color 3 along $h$, it is necessary to have two vertices connected by an edge of color 3, and they cannot belong to the newly added bubble or else the face of degree 4 would be self-intersecting. It is therefore necessary to add yet another bubble. We arrive at
\begin{equation}
\begin{array}{c} \includegraphics[scale=.7]{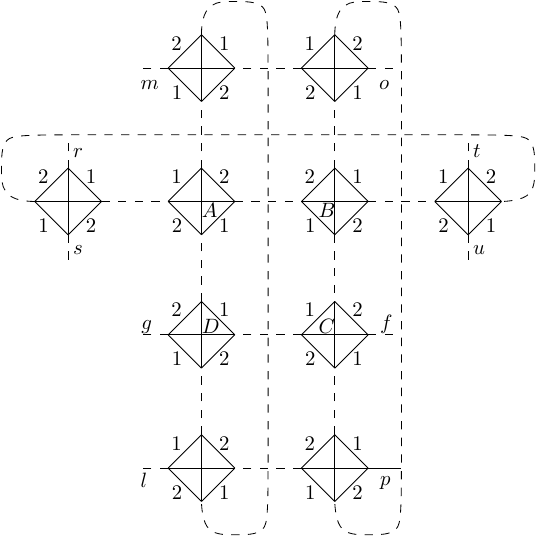} \end{array}
\end{equation}
The face of color 1 along $m$ and $r$ must have degree 4. The half-edges $m, r$ must therefore be connected to the vertices of an edge of 1. Since there is none available, a new bubble must be added. The same holds true for the face of color 2 along $g$ and $s$ (even with the previously added new bubble), so that we get
\begin{equation}
\begin{array}{c} \includegraphics[scale=.7]{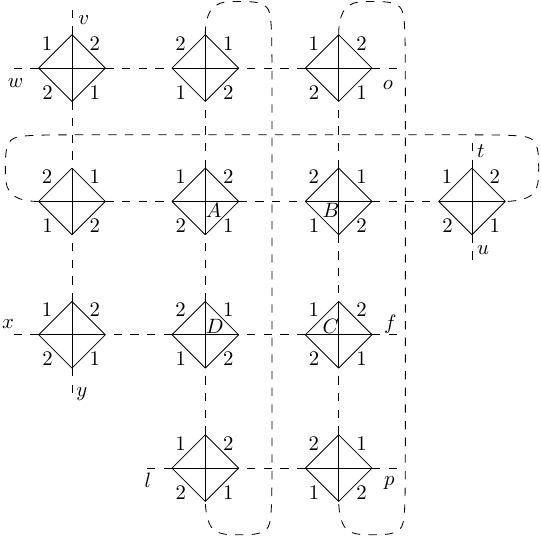} \end{array}
\end{equation}
To close the face of color 1 which goes along $o$ and $t$, one could try to use the edge of color 1 between $v$ and $w$. However, $o$ to $w$ would create a face of color 3 and degree 3, while $o$ to $v$ would create a self-intersecting face of color 3. We therefore need a new bubble. The same argument applies to the face of color 2 along $f$ and $u$ (even with the previously added new bubble since the latter has no edge of color 2 with both ends available), and we get
\begin{equation}
\begin{array}{c} \includegraphics[scale=.7]{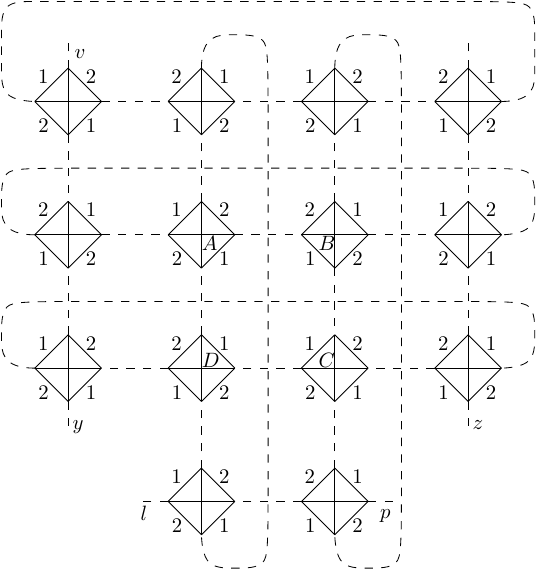} \end{array}
\end{equation}
The face of color 1 going along $l$ and $y$ must have degree 4 (since it is incident to the bubble labeled $D$ which is at a distance of at least 2 from $B_{\text{exc}}$). The half-edges $l, y$ must therefore be connected to both ends of an edge of color 1, but there is no such edge available in the subgraph. A new bubble must therefore be added. This also allows for closing the face of color 3 along $v$ and $y$ (which must also have degree 4 because the bubble to the left of $A$ is at a distance at least 2 of $B_{\text{exc}}$).

\medskip

The whole argument is then repeated one last time: the face of color 1 along $z$ and $p$ must have degree 4 (since it is incident to the bubble labeled $C$ which is at a distance at least 1 from $B_{\text{exc}}$). The half-edges $z, p$ must therefore be connected to both ends of an edge of color 1, but there is no such edge available in the subgraph. A new bubble must therefore be added. It allows for closing the face of color 3 which goes along $z$ (which must have degree 4 because the bubble to the right of $B$ is at a distance at least 1 from $B_{\text{exc}}$), and for closing the face of color 3 along $l$ and $p$ (which must have degree 4 because the bubble below $D$ is at a distance at least 1 from $B_{\text{exc}}$). This fully determines $\cG$ as
\begin{equation}
\begin{array}{c} \includegraphics[scale=.7]{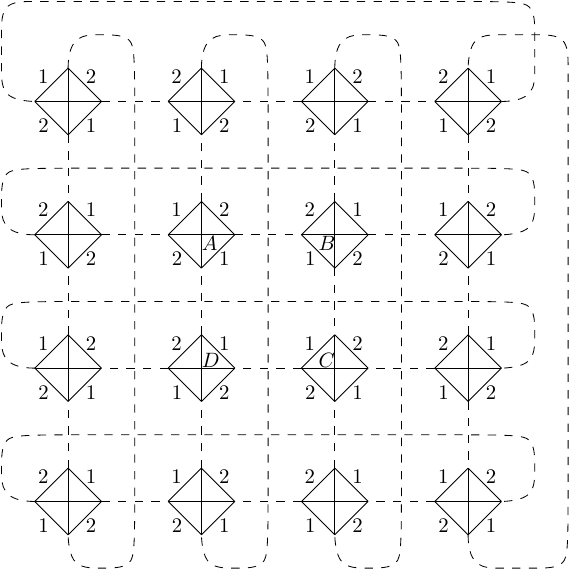} \end{array}
\end{equation}
This exhausts all the possibilities of having a bubble at a distance 3 from $B_{\text{exc}}$. They consist of a finite number of graphs. This concludes the proof of Lemma \ref{lemma:sec}.

%%%%%%%%%%%%%%%%%%%%%%%%%%%%%%%%%%%%%%%%%%%%%%%%%%%%%%%%%%%%%%%%%%%%%%%%%%%%%%%%%%%%%%%%%%%%%%%%%%%%
\subsection{Identification of the dominant schemes and double scaling}

We first analyze the singularity that can be attained and the corresponding combinatorial structure, which turns out to be broken chains. Thus we obtain additional information on the structure of the dominant schemes. This allows us to characterize exactly the dominant schemes at all degree $\omega$ even though schemes of degree $\omega$ are not.

\subsubsection{Identification of relevant singularities}
\label{sssec:RelevantSingularitiesO(N)3}

From Theorem \ref{th:sch}, we know that there is a finite number of schemes at fixed degree, so the singularities can only come from the generating functions of chains and dipoles and that of melonic 2-point graphs. Therefore, we have a priori three different types of singular points:
\begin{itemize}
\item Singular points of $M(t,\mu)$, which are also singular for $U(t,\mu)$ and $B(t,\mu)$.
\item Points such that $U(t,\mu) = 1$, which are singular for any type of chains.
\item Points such that $U(t,\mu) = \frac{1}{3}$ for broken chains only.
\end{itemize}

$M(t,\mu)$ is a generating series whose coefficients $[t^p\mu^q]M$ are the numbers of melonic 2-point graphs with $p$ melons in total and $q$ of type II, therefore its coefficients are all positive. This implies that the function $U_{\mu}: t \mapsto U(t,\mu)$ is an increasing function of $t$. Hence at fixed $\mu$, the point where $U(t,\mu) = \frac{1}{3}$ is always reached for a smaller value of $t$ than $U(t,\mu) = 1$.

\medskip

\begin{figure}
\begin{center}
\includegraphics[scale=0.6]{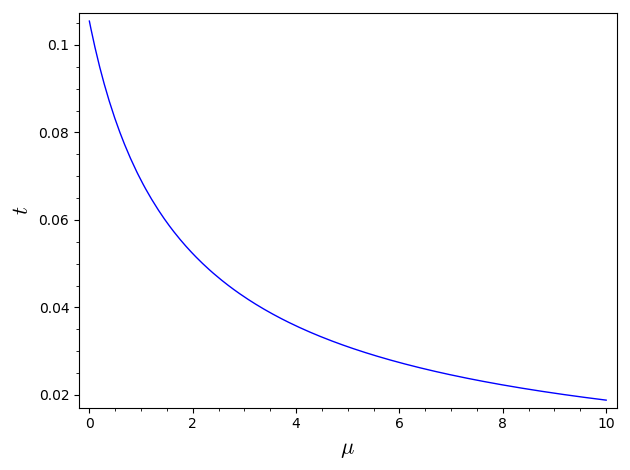}
\caption{Critical points for the generating function $M(t,\mu)$. They also correspond to points where $U(t,\mu)= \frac{1}{3}$ and correspond to critical points for $B(t,\mu)$ as well.}
\label{fig:plot_tc}
\end{center}
\end{figure}

Thus, we only have to know whether we first reach a value of $t$ where $U(t,\mu) = \frac{1}{3}$ or $t_c(\mu)$ such that $(t_c(\mu),\mu)$ is a singular point of $M(t,\mu)$. Using equation~\eqref{eq:mel}, we can express the variable $t$ as $ t = \frac{M(t,\mu)-1}{M(t,\mu)^4 + \mu M(t,\mu)^2}$. Therefore the equation $U(t,\mu) = 1/3$ can be written as
\begin{equation}
M(t,\mu) - \frac{4}{3} - \frac{2}{3} \frac{(M(t,\mu)-1)M(t,\mu)^2 \mu}{M(t,\mu)^4 + \mu M(t,\mu)^2} = 0.
\end{equation}
Clearing the denominator gives
\begin{equation} \label{eq:crit_eq}
-3M(t,\mu)^3 + 4t M(t,\mu)^2 -\mu M(t,\mu) + 2 \mu = 0,
\end{equation}
which actually coincides with the equation~\eqref{eq:Mc_poly} determining the critical values of $M(t,\mu)$. Thus, points where $U(t,\mu) = \frac{1}{3}$ are again exactly the points that are critical for $M(t,\mu)$ and they are called the \emph{dominant singularities} or \emph{critical curve}. It is plotted in Figure~\ref{fig:plot_tc}.

\medskip

Close to the critical curve, the behavior of  the generating series of melons $M(t,\mu)$ is given by equation~\eqref{eq:crit_behav}. Therefore the generating function for chains $U(t,\mu)$ can be expressed as
\begin{align}
U(t,\mu) \underset{t \rightarrow t_c(\mu)}{\sim} &\frac{1}{3} + \left(1 - \frac{4}{3}t_c(\mu)\mu M_c(\mu)\right) K(\mu)\sqrt{1-\frac{t}{t_c(\mu)}} \\ &- \frac{2}{3} t_c(\mu) \mu K(\mu)^2 \left(1-\frac{t}{t_c(\mu)}\right), \nonumber
\end{align}

from which we deduce the behavior of $B$ near the critical curve
\begin{align}
B(t,\mu) &\underset{t\rightarrow t_c(\mu)}{\sim} \frac{1}{\left(1 - \frac{4}{3}t_c(\mu)\mu M_c(\mu)\right) K(\mu)\sqrt{1-\frac{t}{t_c(\mu)}}}
\label{eq:B_crit}
\end{align}

\subsubsection{Identification of dominant schemes}

To perform the double-scaling limit, we identify the schemes at fixed $\omega$ which are the most singular at criticality, and call them \emph{dominant schemes}. From the above analysis, they are the schemes that maximize the number of broken chains. It results in Theorem \ref{thm:DominantSchemesO(N)3} below.

\medskip

Let us recall that a tree is a graph with no cycles. Vertices of valency 1 are called leaves and the others are called internal nodes. A rooted tree is a tree with a marked leaf. A binary tree is a tree whose internal nodes all have valency 3. A tree is said to be plane if it is embedded in the plane. 

\begin{theorem} \label{thm:DominantSchemesO(N)3}
The dominant schemes of degree $\omega>0$ are given bijectively by rooted plane binary trees with $4\omega-1$ edges, with the following correspondence
\begin{itemize}
\item The root of the tree corresponds to the two external legs of the 2-point function.
\item Edges of the tree correspond to broken chains.
\item The leaves are faces of length $2$: $\begin{array}{c} \includegraphics[scale=.55]{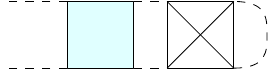}\end{array}$
\item There are two types of internal nodes,
\begin{equation} \label{InnerNodes}
\begin{array}{c} \includegraphics[scale=.6]{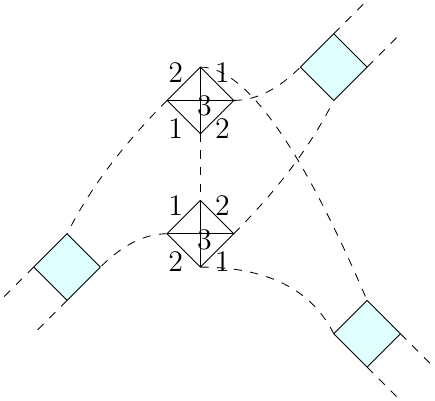}\end{array}, \qquad \begin{array}{c} \includegraphics[scale=.6]{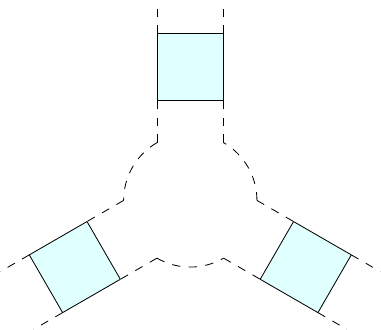}\end{array}
\end{equation}
\end{itemize}
\end{theorem}

The degree of a dominant scheme is thus entirely ``contained'' in the leaves of the tree at the end of the broken chains.

\begin{proof}
We use the notations of Lemma~\ref{thm:SkeletonGraph}. From this lemma, it is clear that if $\cS$ is dominant, then its skeleton graph is a tree, i.e. $\cI(\cS) = \cT$ (or $q=0$), with $N(\cT)=4\omega(\cT)-1$ edges, all corresponding to broken chains in $\cS$. Moreover
\begin{itemize}
\item all its internal nodes have valency 3 and correspond to components $\bar{\cG}^{(r)}$ of degree 0,
\item all leaves correspond to components $\bar{\cG}^{(r)}$ of degree 1/2,
\item the component $\cG^{(0)}$, which has the two external legs, gives rise to a root for $\cT$, and has degree 0.
\end{itemize}

We now have to identify which graphs can appear as nodes and leaves of the tree.
\begin{itemize}
\item Leaves correspond to graphs of degree $1/2$. These graphs have been identified in~\cite{TaCa}. They are tadpoles, graphs with one tetrahedral interaction and two edges of color 0. There are three different tadpoles, depending on the color of the two faces of length $1$. One edge of color 0 is cut and connected to one side of a chain in $\cS$.
\item An internal node of $\cT$ corresponds to a graph of degree 0, with three edges of color 0 cut in order to connect it to three chains. The resulting 6-point function must have no melons and no dipoles. One can check that it must have either $0$ bubbles, i.e. be a single propagator closed on itself, giving rise to the second type of vertices in \eqref{InnerNodes}, or 2 bubbles which is the elementary type I melon, giving rise to the first type of vertices in \eqref{InnerNodes}.
\end{itemize}
Clearly, the dominant schemes are fully encoded by their skeleton graphs which are rooted binary trees, except for the order of the chains meeting an internal node. We thus obtain a bijection between the dominant schemes and rooted binary trees with an order of the edges incident at every vertex, i.e. plane trees. This concludes the proof of the theorem.
\end{proof}

\subsubsection{Generating function for dominant schemes}

A rooted binary tree with $N$ edges has $\frac{N-1}{2}$ inner nodes and $\frac{N+1}{2}$ leaves (not counting the root). Here a leaf carries a weight $3t^{1/2}$ accounting for all three possible colors of the face of length $2$ of the leaf. An inner node receives a weight $1+6t$, the 1 being due to the second type of nodes in~\eqref{InnerNodes} and the $6t$ to the first type.

\medskip

A dominant scheme corresponds to a rooted, plane binary tree $\mathcal{T}$ with $4\omega-1$ edges and thus its generating function is
\begin{equation}
G_{\mathcal{T}}^{\omega}(t,\mu)  = (3t^{\frac{1}{2}})^{2\omega} (1+6t)^{2\omega-1} \frac{6^{4\omega-1}U^{8\omega-2}}{\left((1-U)(1-3U)\right)^{4\omega-1}}
\end{equation}
This function only depends on $\omega$ and not on the shape of $\mathcal{T}$. We can therefore easily sum over all trees and also add all melonic insertions at the root,
\begin{align}
G_{\text{dom}}^{\omega}(t,\mu) &= M(t, \mu) \sum_{\substack{\mathcal{T}\\ \text{$2\omega$ leaves}}} G_{\mathcal{T}}^{\omega}(t,\mu) \\ &
= \Cat_{2\omega-1} G_{\mathcal{T}}^{\omega}(t,\mu) \nonumber
\label{eq:fct_dom_scheme}
\end{align}
where $\Cat_{2\omega-1} = \frac{1}{2\omega}\binom{4\omega-2}{2\omega-1}$ is the number of rooted, plane, binary trees with $2\omega$ leaves.

%%%%%%%%%%%%%%%%%%%%%%
\subsection{Double scaling limit of the quartic \texorpdfstring{$O(N)^3$}{O(N)3} model}

Using Equation~\eqref{eq:crit_behav}, the generating function of dominant schemes behaves near singular points as
\begin{align}
G_{dom}^\omega(t,\mu) \underset{t\rightarrow t_c}{\sim} &M_c(\mu)\Cat_{2\omega-1}9^{\omega}t_c^{\omega}\left(1+6t_c\right)^{2\omega-1}\\
&\times\left(\frac{1}{\left(1 - \frac{4}{3}t_c(\mu)\mu M_c(\mu)\right) K(\mu)\sqrt{1-\frac{t}{t_c(\mu)}}}\right)^{4\omega-1} \nonumber
\end{align}
Since in the large $N$ expansion a graph $\cG$ of degree $\omega$ scales as $N^{3-\omega}$ we define the following double scaling parameter
{\small \begin{equation}
\kappa(\mu)^{-1} = N^{\frac{1}{2}}\frac{1}{3}\frac{1}{t_c(\mu)^\frac{1}{2}\left(1+6t_c(\mu)\right)}\left( \left(1 - \frac{4}{3}t_c(\mu)\mu M_c(\mu)\right) K(\mu) \right)^2\left(1-\frac{t}{t_c(\mu)}\right)
\label{eq:kappa}
\end{equation}}%
Using Equation~\eqref{eq:kappa} we get:
{\small \begin{equation}
\left[ \frac{ (1+6t_c(\mu))}{\left(1 - \frac{4}{3}t_c(\mu)\mu M_c(\mu)\right) K(\mu)\sqrt{1-\frac{t}{t_c(\mu)}}}\right]^{-1} = \kappa(\mu)^{-\frac{1}{2}}N^{-\frac{1}{4}}\frac{\sqrt{3}t_c(\mu)^\frac{1}{4}}{\left(1+6t_c(\mu)\right)^\frac{1}{2}}
\end{equation}}%

Therefore in the double scaling limit, the dominant graphs of degree $\omega > 0$ contribute as:
\begin{equation}
G_{dom}^\omega(\mu) = M_c(\mu)\frac{N^{\frac{11}{12}}}{\kappa(\mu)^\frac{1}{2}}\sqrt{3}\frac{t_c(\mu)^\frac{1}{4}}{\left(1+6t_c(\mu)\right)^\frac{1}{2}}\Cat_{2\omega-1}\kappa(\mu)^{2\omega}
\end{equation}
where one has to add the contribution of the graphs of degree $0$ i.e. the melons, which contribute simply as $M_c(t,\mu)$. Hence summing over contributions of all degree, the total contribution to $G_2^{DS}$ is:
\begin{align}
\label{eq:GDS_ON3}
G_{2}^{DS}(\mu) &= \sum\limits_{\omega\in\mathbb{N}/2} G_{dom}^\omega(\mu) \nonumber \\
&= M_c(\mu) + M_c(\mu)\frac{N^{\frac{11}{12}}}{\kappa(\mu)^\frac{1}{2}}\sqrt{3}\frac{t_c(\mu)^\frac{1}{4}}{\left(1+6t_c(\mu)\right)^\frac{1}{2}} \sum\limits_{n \in \frac{\mathbb{N}}{2} > 0 } \Cat_{2\omega-1}\kappa(\mu)^{2\omega} \nonumber \\
&= M_c(\mu) + M_c(\mu)\kappa(\mu)N^{\frac{11}{12}}\sqrt{3}\frac{t_c(\mu)^\frac{1}{4}}{\left(1+6t_c(\mu)\right)^\frac{1}{2}} \sum\limits_{n \in \mathbb{N}} \Cat_n \kappa(\mu)^{n} \nonumber \\
&= M_c(\mu) \left(1 + N^{\frac{11}{12}}\sqrt{3}\frac{t_c(\mu)^\frac{1}{4}}{\left(1+6t_c(\mu)\right)^\frac{1}{2}} \frac{1-\sqrt{1-4\kappa(\mu)}}{2\kappa(\mu)^\frac{1}{2}}\right).
\end{align}

\section{Double scaling limit for the quartic \texorpdfstring{$U(N)^2\times O(D)$}{U(N)2XO(D)} multi-matrix model}
\label{sec:UN2OD} 

The method employed to derive the double scaling limit of tensor models can also be applied in multi-matrix models. We illustrate this on the quartic $U(N)^2\times O(D)$ multi-matrix model. Despite the differences in the nature of the expansion of the two models, their Feynman diagrams are similar which allows us to exploit the combinatorial results obtained in the previous section for the quartic $O(N)^3$ for this model. Therefore this computation is similar to the case of the $O(N)^3$ tensor model and we will only explain here how the results of the quartic $O(N)^3$ tensor model can be applied and point out the differences that exist between the two cases. The details of the computation can be found in Appendix~\ref{app:DS_UN2OD}.

\subsection{Description of the model}
\label{ssec:Un2xOd_model}

The $U(N)^2 \times O(D)$ multi-matrix model is a model involving a vector of $D$ complex matrices of size $N \times N$, denoted $(X_\mu)_{\mu=1, \dotsc, D} = (X_1, \dotsc, X_D)$. The model is required to be invariant under unitary actions to the left and on the right of each matrix $X_\mu$,
\begin{equation}
X_\mu \rightarrow X'_\mu = U_1 X_{\mu} U_2^\dagger
\end{equation}
with $U_1,U_2 \in U(N)$ and under orthogonal transformations of the matrices $X_{\mu}$
\begin{equation}
X_\mu \rightarrow X'_\mu = \sum_{\mu'=1}^D O_{\mu\mu'}\ X_{\mu'} 
\end{equation}
for any $O \in O(D)$. 

\medskip

To describe the bubbles of the model, we associate the color $1$ and $2$ respectively to the matrix indices $i$ and $j$ of the matrices $X_\mu = \left((X_\mu)_{ij}\right)_{1 \leq i,j, N}$. The set of possible bubbles of this model corresponds to $3$-colored graph whose vertices are assigned a color, white or black corresponding respectively to copies of $X_\mu$ or $X_\mu^\dagger$, and such that edges of color $1$ and $2$ connect black to white vertices. Therefore, the set of bubbles for this model can be obtained from bubbles of the $O(N)^3$ tensor model by specifying a suitable coloring of its vertices. From the quartic bubbles of the $O(N)^3$ tensor model, we get the following quartic interactions for the $U(N)^2 \times O(D)$ multi-matrix model.

\medskip

The connected, quartic interactions are
\begin{align}
I_t(X, X^\dagger) = \sum_{\mu, \nu} \Tr \bigl(X_\mu X^\dagger_\nu X_\mu X^\dagger_\nu\bigr) &= \begin{array}{c}\includegraphics[scale=.28]{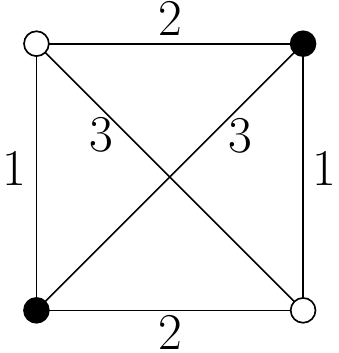}\end{array},\\
I_{1}(X, X^\dagger) = \sum_{\mu, \nu} \Tr \bigl(X_\mu X^\dagger_\mu X_\nu X^\dagger_\nu\bigr) &= \begin{array}{c}\includegraphics[scale=.28]{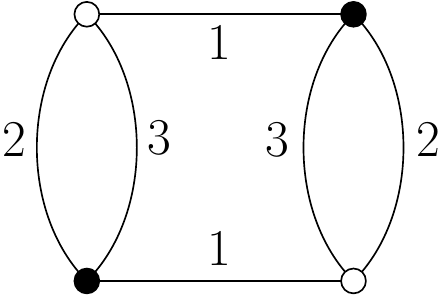}\end{array},\\
I_{2}(X, X^\dagger) = \sum_{\mu, \nu} \Tr \bigl(X_\mu X^\dagger_\nu X_\nu X^\dagger_\mu\bigr) &= \begin{array}{c}\includegraphics[scale=.28]{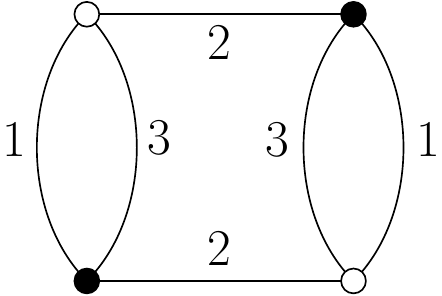}\end{array},\\
I_{3b}(X, X^\dagger) = \sum_{\mu, \nu} \Tr \bigl(X_\mu X^\dagger_\nu\bigr)\ \Tr \bigl(X_\nu X^\dagger_\mu\bigr) &= \begin{array}{c}\includegraphics[scale=.28]{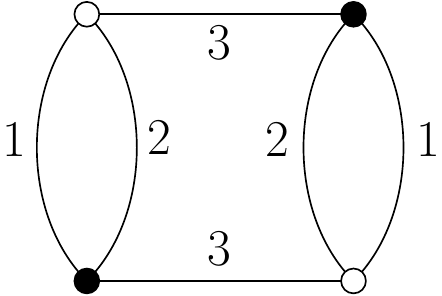}\end{array},\\
I_{3nb}(X, X^\dagger) = \sum_{\mu, \nu} \Tr \bigl(X_\mu X^\dagger_\nu\bigr)\ \Tr \bigl(X_\mu X^\dagger_\nu\bigr) &= \begin{array}{c}\includegraphics[scale=.28]{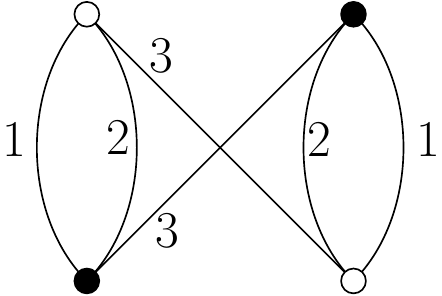}.\end{array}
\end{align}
For the tetrahedral interaction and the pillows of color $1$ and $2$, there is a unique coloring of its vertices which satisfies this bipartiteness condition. However, there are two possible coloring that can be associated with the pillow of color $3$, leading to an additional interaction compared to the $O(N)^3$ tensor model.

\medskip

Finally, the propagator of the model is given by the quadratic bubble
\begin{equation}
    I_k(X, X^\dagger) = \sum_{\mu=1}^D \Tr X_\mu X^\dagger_\mu,
\end{equation}
and associated with the color $0$. The set of connected Feynman graphs, denoted $\bar{\mathbb{G}}$ is the set of connected, 4-regular, properly-edge-colored graphs such that the subgraph obtained by removing all edges of color 0 is a disjoint union of quartic bubbles $I_{p;1}, I_{p;2}, I_{p;3b}, I_{p;3nb}, I_t$.

\medskip

The action of the model with quartic interaction and enhanced scaling factor~\eqref{eq:scal_Gur_enh_TM} writes 
\begin{align} \label{UN2ODAction}
S_{U(N)^2\times O(D)}(X_{\mu},X_{\mu}^{\dagger })= -&ND \sum_{\mu=1}^D \Tr(X_{\mu} X^\dagger_{\mu}) + N D^{3/2}\frac{\lambda_1}{2} I_t(X_{\mu},X_{\mu}^{\dagger }) \\
&+ ND\frac{\lambda_2}{2} \bigl(I_{p;1}(X_{\mu},X_{\mu}^{\dagger })+I_{p;2}(X_{\mu},X_{\mu}^{\dagger })\bigr) \nonumber\\
&+ D^2\frac{\lambda_2}{2} \bigl(I_{p;3b}(X_{\mu},X_{\mu}^{\dagger }) + I_{p;3nb}(X_{\mu},X_{\mu}^{\dagger })\bigr)\nonumber.
\end{align}
Observe that since the indices of the tensor now play different roles, the scaling of the three pillows now differ depending on their color. The partition function is
\begin{equation}
Z_{U(N)^2\times O(D)}(\lambda_1, \lambda_2) = \int \prod_{\mu=1}^D \prod_{a,b=1}^N d(X_\mu)_{ab} d(\overline{X_\mu})_{ab}\ e^{S_{N,D}(X_\mu, X^\dagger_\mu)}.
\end{equation}

Let $\bar{\cG}\in\bar{\mathbb{G}}$, then denote
\begin{itemize}
\item $n_t, n_{1}, n_{2}, n_{3b}, n_{3nb}$ the number of interactions given by $I_t, I_{1}, I_{2}, I_{3b}, I_{3nb}$ respectively,
\item $E$ the number of edges of color $0$ (which can be expressed as $E = 2(n_t+n_{1}+n_{2}+n_{3b}+n_{3nb})$),
\item $F_{a}$ for $a=1,2,3$, the number of bicolored cycles which alternate the colors $0, a$. The faces of color $3$ are associated with scaling in $D$ of the graph while the faces of color $1,2$ are associated with its scaling with $N$.
\end{itemize}

The free energy $1/N,1/D$-expansion is well-defined for this choice of scaling factor\cite{FeVa}. The $1/N$ expansion is governed by $h(\bar{\cG})$ given by~\eqref{eq:genus_TM} and the $1/D$ expansion is governed by the grade denoted $l(\bar{\cG})$ defined by~\eqref{eq:def_grade}. Therefore the free energy expands as
\begin{equation}\label{F_exp}
F_{U(N)^2\times O(D)}(\lambda_1, \lambda_2) = \sum_{\bar{\cG}\in\bar{\mathbb{G}}} \biggl(\frac{N}{\sqrt{D}}\biggr)^{2-2h(\bar{\cG})} D^{2-l(\bar{\cG})/2}\ \mathcal{A}_{\bar{\cG}}(\lambda_1, \lambda_2).
\end{equation}

\subsection{Scheme decomposition}
\label{ssec:scheme_UN2}

The scheme decomposition of the $U(N)^2\times O(D)$ is completely similar to the $O(N)^3$ tensor model. The melons, dipoles and chains allow us to repackage together infinitely many graphs of given genus $h$ and grade $l$ in a single scheme. However, the inclusion of the second pillow of color $3$ induces some slight changes to graphs occurring as dipoles and chains. Indeed we must now distinguish two dipole vertices of color $3$ according to the coloring of their vertices.
\begin{itemize}
\item One involving the tetrahedral interaction, and the non-bipartite pillow of color $3$,
\begin{equation}
D_3 = \left\{\includegraphics[scale=0.25,valign=c]{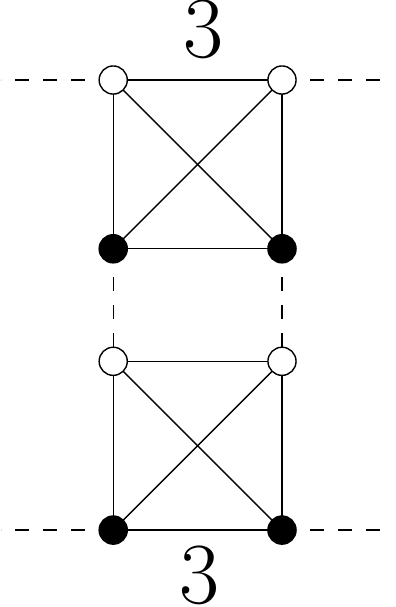}, \hspace{.5cm} \includegraphics[scale=0.26,valign=c]{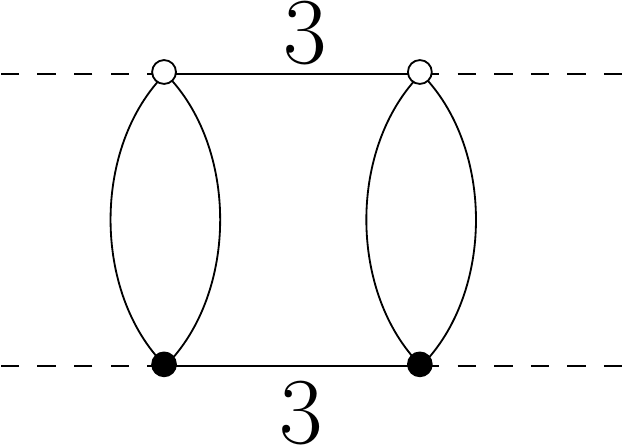}\right\}
\label{eq:d3}
\end{equation}
\item One for the bipartite pillow of color $3$, %denoted $D'_3$
\begin{equation}
D_3' = \left\{\includegraphics[scale=0.3,valign=c]{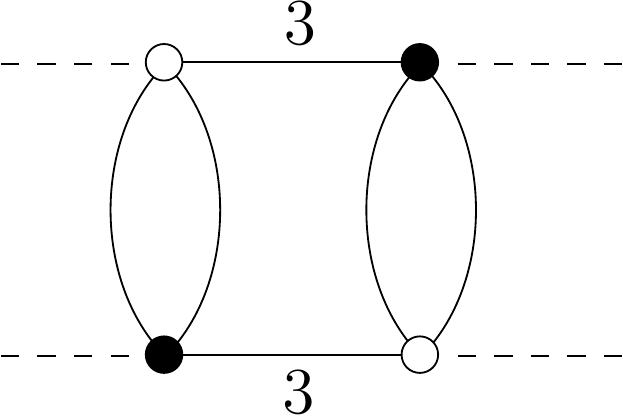}\right\}
\label{eq:d3'}
\end{equation}
\end{itemize}
While the situation is completely identical for dipoles of color $1$ and $2$. Denoting $U$ the generating function for the dipoles $D_{1,2,3}$ decorated with melons on one side and on the internal edges, and $V$ the generating function for $D'_3$ decorated with melons on one side. Then we have
\begin{align}
U(t,\mu) &= tM(t,\mu)^4 + \frac{1}{4}t\mu M(t,\mu)^2 \underset{\eqref{eq:mel_un2}}{=} M(t,\mu)-1-\frac{3}{4}t\mu M(t,\mu)^2, \\
V(t,\mu) &= \frac{1}{4}t\mu M(t,\mu)^2.
\end{align}

This change propagates to chains of color $3$ as the type of dipole induces changes in the coloring of the chain at its boundary. Therefore we now also distinguish two types of chains of color $3$ depending on the coloring at its boundary. Their respective generating functions are
\begin{align}
C_{3,\hspace{1pt}\includegraphics[scale=0.2,valign=c]{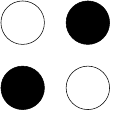}} &= \frac{U^2+V-V^2}{(1-V-U)(1-V+U)},
\end{align}
and
\begin{align}
C_{3,\hspace{1pt}\includegraphics[scale=0.2,valign=c]{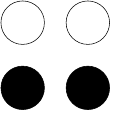}} &=  \frac{U}{(1-V-U)(1-V+U)}.
\end{align}

However, this difference is only apparent at the level of the generating series. When manipulating graphs, the two models are much more similar as they only differ from the coloring of the vertices in the $U(N)^2 \times O(D)$ model. These similitudes can be used to derive the following Theorem, which is analogous to Theorem~\ref{thm:sch_ON3} for this model.

\begin{theorem} \label{thm:FiniteSchemesI}
$\mathbb{S}_{h,l}$ is a finite set.
\end{theorem}

\begin{proof}
Since bubbles of the $U(N)^2 \times O(D)$ model can be obtained from those of the ones of the $O(N)^3$ model, there is a map between the graph of the two models
\begin{equation}\label{eq:theta}
\theta: \mathbb{G} \to \mathbb{G}_{O(N)^3}
\end{equation}
which simply consists of forgetting the vertex coloring of the graph in the $U(N)^2 \times O(D)$ model.

\medskip

If $\mathcal{G}\in\mathbb{G}_{U(N)^2\times O(D)}$, we define the degree of $\cG$ as $\omega(\mathcal{G}) := \omega(\theta(\mathcal{G}))$ the degree of the graph it is mapped to in the $O(N)^3$ model. Then
\begin{equation}
\omega(\mathcal{G}) = h(\mathcal{G}) + \frac{l(\mathcal{G})}{2}
\end{equation}
which can be established by setting $D=N$ in~\eqref{UN2ODAction}.

\medskip

The map $\theta$ maps melons to melons, dipoles to dipoles and chains to chains. We denote $\mathbb{S}_{O(N)^3}(\omega)$ be the set of schemes of degree $\omega$ in the $O(N)^3$ model. We can therefore descend the map $\theta$ to schemes by
\begin{equation}
\tilde{\theta}_{h,l}: \mathbb{S}_{h,l} \to \mathbb{S}_{O(N)^3}(h+l/2)
\end{equation}
from the schemes of genus $h$ and grade $l$ to those of the $O(N^3)$ model with degree $h+l/2$. It simply consists of forgetting the coloring of the vertices on the vertices, dipoles and chains.

\medskip

Theorem~\ref{th:sch} has been established for the quartic $O(N)^3$ tensor model (see Section~\ref{ssec:scheme_O(N)^3}). Therefore, it is sufficient for us to show that each scheme $\mathcal{S}\in \mathbb{S}_{O(N)^3}(h+l/2)$ has a finite fiber via $\tilde{\theta}_{h,l}$. This is clear since the fiber of $\mathcal{S}\in \mathbb{S}_{O(N)^3}(\omega)$ is found by considering all possible colorings of the vertices and chain-vertices. This proves Theorem \ref{thm:FiniteSchemesI}.
\end{proof}

\begin{figure}[!ht]
\begin{center}
\includegraphics[scale=0.55]{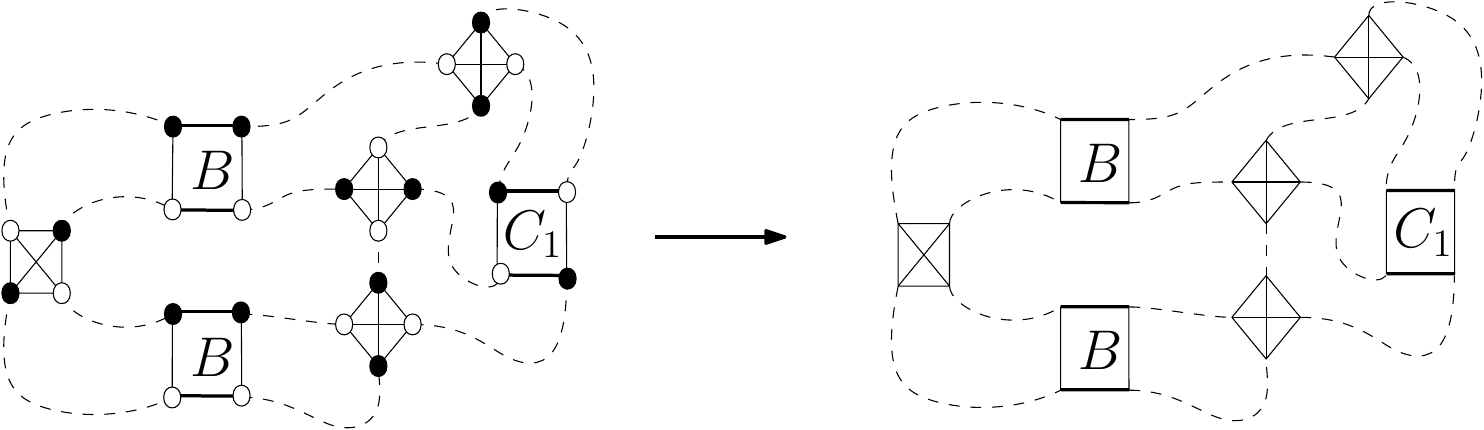}
\caption{A scheme of the $U(N)^2\times O(D)$-invariant model and its corresponding scheme on the side of the $O(N)^3$-invariant model.}
\end{center}
\end{figure}

The map $\theta$ can be used to identify to show that leading order graphs of the $U(N)^2 \times O(D)$ model are melonic graphs. It maps graphs of genus $h$ and grade $l$ to graphs of degree $\omega = g+ \frac{l}{2}$, therefore the leading order graphs of this model are a subset of the leading order graph of the $O(N)^3$ tensor model and conversely, any melons of the $O(N)^3$ tensor model can be assigned at least one coloring of its vertices satisfying the bipartiteness condition of this model.

\medskip

However, due to the difference in the nature of the expansion, it cannot be used to identify directly the dominant scheme for the $U(N)^2\times O(D)$ model from the dominant schemes of the $O(N)^3$ model. Indeed, the $\frac{1}{N},\frac{1}{D}$-expansion selects schemes with vanishing grade. Since the genus is an integer, we can already see that all dominant schemes of the $O(N)^3$ tensor model with half-integer degree cannot have vanishing grade and therefore cannot be dominant schemes for the $U(N)^2\times O(D)$ model. Still, it is possible to perform a similar combinatorial decomposition of the schemes via their skeleton graph which presents no particular difficulty. This leads to the following structure for the dominant scheme for this model.

\begin{proposition} 
\label{prop:dom_scheme_un2}
A dominant scheme of genus $h>0$ has $2h-1$ broken chain-vertices, all separating. Such a scheme has the structure of a rooted binary plane tree where
\begin{itemize}
\item Edges correspond to broken chain-vertices.
\item The root of the tree corresponds to the two external legs of the 2-point function.
\item The $h$ leaves are one of the two following graphs
\begin{equation}
\includegraphics[scale=0.5]{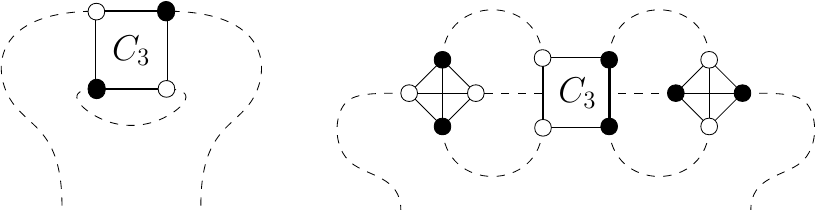}
\end{equation}
Observe that the chains of color $3$ have different boundary vertices in the two graphs.
\item Each internal vertex corresponds to one of four $6$-point subgraphs:
\begin{equation}
\includegraphics[scale=0.5]{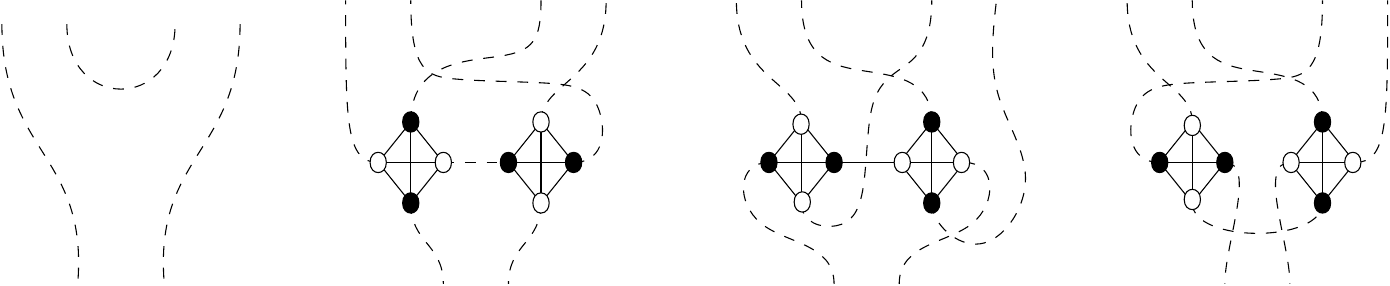}
\end{equation}
Here we have used the embedding convention that the edge coming out of the bottom left is connected to a white vertex.
\end{itemize}
\end{proposition}

As we can see, the structure of the dominant schemes for the $U(N)^2\times O(D)$ model is different from the one of Theorem~\ref{thm:DominantSchemesO(N)3} of the quartic $O(N)^3$ model because of the difference in the leaves of the tree, which here correspond to graphs of genus $1$ with vanishing grade. Also, note that there are now two types of leaves with different generating functions. Once the structure of the dominant scheme has been obtained, a simple summation of their contribution gives the expression of the $2$-point function in the double scaling limit. This limit is achieved by sending first the size parameters of the tensor $N$ and $D$ to infinity while holding the ratio $L = \frac{N}{\sqrt{D}}$ fixed, and then sending $L$ to infinity while holding the double scaling parameter $\kappa(\mu)$ fixed. This parameter is defined as
\begin{equation}
\kappa(\mu)^{-1} = L^2 \frac{1}{(1+6t)(C_{3,\hspace{1pt}\includegraphics[scale=0.2,valign=c]{wbbw.pdf}}(t,\mu)+tC_{3,\hspace{1pt}\includegraphics[scale=0.2,valign=c]{wwbb.pdf}}(t,\mu))}\left(\frac{1}{B(t,\mu)}\right)^2.
\label{eq:kappa_un2}
\end{equation}
In this limit, when summing the contributions of the dominant schemes of all degrees we obtain the $2$-point function of the $U(N)^2 \times O(D)$ model in the double scaling limit
{\small \begin{align}
G_{2}^{DS}(\mu) &= M_c(t_c(\mu),\mu)\left(1+ L \left(\frac{(C_{3,\hspace{1pt}\includegraphics[scale=0.2,valign=c]{wbbw.pdf}}+tC_{3,\hspace{1pt}\includegraphics[scale=0.2,valign=c]{wwbb.pdf}})}{1+6t}\bigg\rvert_{(t_c(\mu),\mu)} \right)^{\frac{1}{2}} \frac{1-\sqrt{1-4\kappa(\mu)}}{2\kappa(\mu)^\frac{1}{2}}\right). \label{eq:GDS_U(N)2}
\end{align}}%

\section{Double scaling limit for the \texorpdfstring{$U(N)\times O(D)$}{U(N)xO(D)} bipartite tensor model}
\label{sec:DS_U(N)xO(D)}

In tensor models with mixed indices symmetry, the mixing of indices prevents from encoding the possible interaction of the model in terms of colored graphs, and their propagator usually takes the form of a sum of monomials which can lead to non-trivial cancellations between Feynman graphs. Additional constraints (to ensure the irreducibility of the representation to which the group action corresponds) often have to be imposed to ensure the existence of the $\frac{1}{N}$-expansion for this model. This makes this type of model more difficult to study via purely combinatorial techniques, even to establish fundamental results such as melonic dominance in the $\frac{1}{N}$-expansion. However, in the few cases where it was established -for example in~\cite{Ca18} for a rank $3$ tensor model transforming in an irreducible representation of $O(N)$ or in~\cite{CaHa21} at rank $5$- these model seem to share similar features as models with no-mixed indices symmetry. In this section, we show that the tools used to derive the double scaling limit for the $O(N)^3$ model can be adapted and applied to derive the double scaling limit to the quartic $U(N)\times O(D)$ model with tetrahedral interaction.

\subsection{The model and its large \texorpdfstring{$N$}{N}, large \texorpdfstring{$D$}{D} expansion}

\subsubsection{\label{UNOD:graph_exp} Feynman graphs, genus and grade}

\paragraph{The model.\\}

The $U(N)\times O(D)$ multi-matrix model is a close cousin of the $U(N)^2\times O(D)$ tensor model. We consider a single quartic interaction and its complex conjugate, leading to the following action 
\begin{align} \label{UNODAction}
S_{U(N)\times O(D)}(X_\mu, X_\mu^\dagger) = &\frac{\lambda}{4} ND^{\frac{3}{2}} \sum_{\mu, \nu=1}^D \left( \Tr (X_\mu X_\nu X_\mu X_\nu) + \Tr (X^\dagger_\mu X^\dagger_\nu X^\dagger_\mu X^\dagger_\nu) \right) \\
& -ND\sum\limits_{\mu = 1}^{D} \Tr(X_{\mu}^{\dagger} X_{\mu}) \nonumber
\end{align}%
It is invariant under 
\begin{equation}
X_{\mu} \mapsto \sum_{\mu'=1}^D O_{\mu\mu'}\ UX_{\mu'}U^\dagger,
\end{equation}
where $O\in O(D)$ and $U\in U(N)$. There is now a single copy of the unitary group acting on the matrices $X_{\mu}$ by conjugation, which induces mixing of the tensor indices associated with indices of the matrices. Due to this change, the Feynman graphs for this model are not colored graphs, but can be represented via decorated maps instead.
\begin{itemize}
\item We represent an interaction $\sum_{\mu, \nu} \Tr (X_\mu X_\nu X_\mu X_\nu)$ as a white vertex of degree 4, and each matrix as an incident half-edge. These half-edges are cyclically ordered, say counter-clockwise, around each vertex, reflecting the structure of the trace.
\item The interaction $\sum_{\mu, \nu} \Tr (X_\mu^\dagger X_\nu^\dagger X_\mu^\dagger X_\nu^\dagger)$ is represented as a black vertex, also with 4 incident, cyclically ordered, half-edges.
\item Propagators glue half-edges incident to black vertices to half-edges incident to white vertices.
\end{itemize}
Feynman graphs are thus ribbon graphs which are bipartite and tetravalent. The faces of the map are the cycles obtained by following edges from black to white vertices and corners at vertices counter-clockwise.

\medskip

The propagator is $\frac{1}{ND} \delta_{ad} \delta_{bc} \delta_{\mu\nu}$. It therefore identifies the matrix indices of half-edges and as usual in matrix models, each face receives a weight $N$. Moreover, propagators also identify the vector indices of half-edges. There is therefore another notion of cycles that receive the weight $D$. They are cycles that follow edges and cross the vertices (i.e. leaving one half-edges on each side). They are similar to the bicolored cycles of color $\{0,3\}$ of the previous models, and we will call them \emph{straight cycles}. This is an additional structure compared to the usual case of matrix models which can be interpreted as loops carried on the edges of the map.

\medskip

Let $\bar{\mathbb{M}}$ be the set of vacuum, connected Feynman graphs. For $\bar{\mathcal{M}}\in\bar{\mathbb{M}}$ we denote $V$, $E$, $F$, $\phi$ the numbers of vertices, edges, faces and straight cycles.  

\medskip

The $1/N,1/D$ expansion of this model has been established in~\cite{TaFe}.

\begin{theorem}[\cite{TaFe}]\label{thm:free_energy_U(N)}
The free energy admits the following expansion on maps
\begin{equation}
F_{U(N)\times O(D)}(\lambda) = \sum_{\bar{\mathcal{M}}\in\bar{\mathcal{M}}} \biggl(\frac{N}{\sqrt{D}}\biggr)^{2-2g(\bar{\mathcal{M}})} D^{2-l(\bar{\mathcal{M}})/2}\ \mathcal{A}_{\bar{\mathcal{M}}}(\lambda)
\end{equation}
where $g(\bar{\mathcal{M}})$ is the genus of $\bar{\mathcal{M}}$ and $l(\bar{\mathcal{M}})$ is a non-negative integer called the \emph{grade},
\begin{equation} \label{GenusAndGradeU(N)}
\begin{aligned}
2-2g(\bar{\mathcal{M}}) &= F-E+V,\\
\frac{l(\bar{\mathcal{M}})}{2} &= 1 + g(\bar{\mathcal{M}}) + E - \frac{3V}{2} - \phi.
\end{aligned}
\end{equation}
\end{theorem}
Like in the previous model, the large $N$, large $D$ expansion is in fact an expansion in $D$ and $L = \frac{N}{\sqrt{D}}$.

\paragraph{2-point maps.\\}

Due to the symmetries of the model, the 2-point function has the form
\begin{equation}
    \langle (X_\mu)_{ab} (X^\dagger_\nu)_{cd}\rangle = \frac{1}{N^2D} G_{N, D}(\lambda) \delta_{\mu\nu} \delta_{ad} \delta_{bc},
\end{equation}
where $G_{N, D}(\lambda) = \langle \sum_{\mu=1}^D \Tr X_\mu X^\dagger_\mu\rangle$. It has an expansion on 2-point maps, whose set is denoted $\mathbb{M}^\circ$. A 2-point map is like a vacuum map with exactly one half-edge left incident on a black vertex and one half-edge incident to a white vertex. These two half-edges are necessarily part of the same straight path, hence the factor $\delta_{\mu\nu}$ above.

\medskip

From a 2-point map $\mathcal{M}^\circ\in\mathbb{M}^\circ$ one can obtain a unique vacuum graph $\bar{\mathcal{M}}\in \bar{\mathbb{M}}$ called its \emph{closure}, by connecting the two half-edges. This vacuum graph is furthermore equipped with a marked edge (that is obtained by connecting the two half-edges) called a \emph{root}. The set of rooted maps, denoted $\mathbb{M}$, is in fact the set of Feynman graphs for the expansion of $G_{N,D}(\lambda)$.

\medskip

Note that in contrast with the model of the previous section where we had defined equivalent objects, it is here not true that the amplitudes of $\mathcal{M}^\circ\in\mathbb{M}^\circ$ and of its closure $\bar{\mathcal{M}}$ are simply related by a factor $N^2 D$. This is because in the previous models, the faces on both sides of an edge could never be the same (they carry different colors) while here they may. As a consequence, cutting an edge of $\mathcal{M}$ might break either one or two faces. 

\medskip

However, it is clear that the amplitude of a rooted map is the same as that of the vacuum map obtained by removing the marking on the root. Although there is a clear bijection between 2-point maps and rooted maps, the expansions are thus a bit different. In the following, we will focus on the rooted function $G_{N,D}(\lambda)$.

\medskip

If $\mathcal{M}$ is a rooted map and $\bar{\mathcal{M}}\in\bar{\mathbb{M}}$ the one obtained by removing the marking, we define $h (\mathcal{M}):= h(\bar{\mathcal{M}})$ and $l(\mathcal{M}) := l(\bar{\mathcal{M}})$.

\subsubsection{Melons, dipoles and chains}

Similarly to other models, we define here the melons, dipoles and chains which are the building blocks of the scheme decomposition for this model.

\paragraph{Melons\\}

The elementary melon is the following 2-point map,
\begin{equation}
\includegraphics[scale=.45, valign=c]{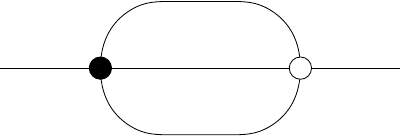}
\end{equation}
Melonic maps are defined starting from the elementary melon itself, and replacing any edge (or one of the half-edges) with another elementary melon, and so on. After inserting some elementary melons, a rooted map is obtained by gluing the two half-edges and marking this edge. 

\begin{proposition}
Melonic maps are the only maps with vanishing genus and grade.
\label{claim_mel}
\end{proposition}
This proposition will be proved below after the introduction of dipoles and their removal.

\medskip

The generating function $M(t)$ of melonic maps satisfies
\begin{equation}
M(t) = 1 +tM(t)^4
\label{eq:mel_unxod}
\end{equation}
where $t$ counts the number of melonic insertions in the map (equivalently half the number of vertices). The singular points and values of this function have been studied in~\cite{GuSch}. In particular its leading singularity is at $t_c=\frac{3^3}{4^4}$ where it has value $M(t_c)=\frac{4}{3}$.

\paragraph{Dipoles\\}

Dipoles are $4$-point functions obtained from the elementary melon by cutting an edge. The $4$ half-edges thus obtained form two pairs of half-edges: one pair formed by the original half-edges of the elementary melon, and one pair from the cutting of the edge to get a dipole. There are two different dipoles depending on the edge which is cut.
\begin{itemize}
\item Cutting the ``top'' or ``bottom'' edge of the elementary melon gives a dipole with a face of degree $2$, called a \emph{U-dipole}
\begin{equation}
\includegraphics[scale=.5, valign=c]{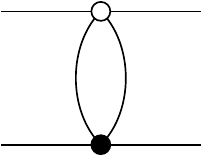}
\end{equation}
\item Cutting the ``middle'' edge of the elementary melon gives a dipole with a cycle of length $2$, called an \emph{O-dipole}
\begin{equation}
\includegraphics[scale=.5, valign=c]{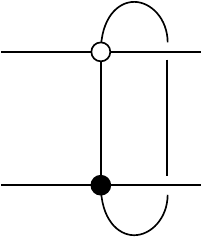}
\end{equation}
\end{itemize}
Notice that a U-dipole (respectively an O-dipole) is equivalent to a face of degree 2 (respectively a straight cycle of degree 2). We usually represent dipoles by drawing the half-edges of the same pair on the same ``side''.

\medskip

\begin{figure}
    \centering
    \includegraphics[scale=.5, valign=c]{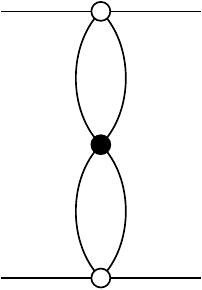}\hspace{2cm}\includegraphics[scale=.5,valign=c]{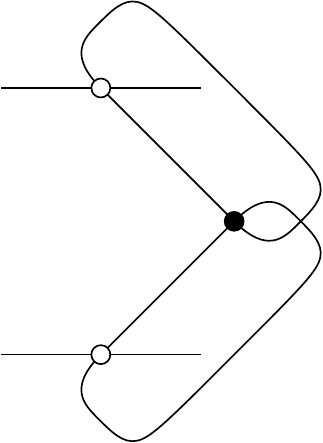}
    \caption{There are two possible minimal subgraphs containing non-isolated dipoles, up to vertex coloring.}
    \label{fig:VertexJointDipoles}
\end{figure}

Dipoles that are vertex-disjoint from other dipoles are said to be \emph{isolated}, else they are \emph{non-isolated}. A non-isolated dipole can in fact be part of only two possible minimal subgraphs, given in Figure~\ref{fig:VertexJointDipoles}.

Let us denote $D_O$ and $D_U$ the generating functions of O- and U-dipoles decorated with melons on their two (internal) edges and on two half-edges of a side (it does not matter which side; there is actually no way to distinguish them). We have
\begin{align}
D_O(t) &= U(t) = tM(t)^4 \underset{\eqref{eq:mel_unxod}}{=} M(t)-1 \\
D_U(t) &= 2U(t),
\end{align}
where the 2 in the last line comes from the two possible edges of the elementary melon which can be cut to obtain the U-dipole.

\medskip

\begin{figure}[!ht]
\begin{center}
\includegraphics[scale=0.6]{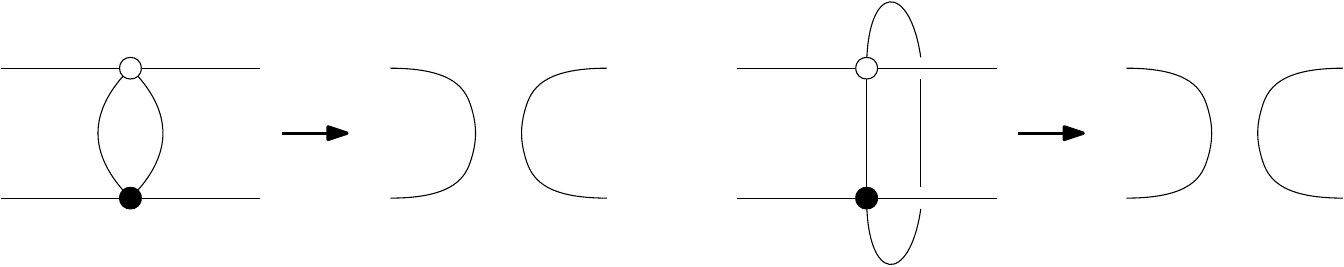}
\caption{On the left is the dipole removal of a U-dipole and on the right, of an O-dipole.}
\label{fig:dip_rem_unxod}
\end{center}
\end{figure}

The \emph{dipole removal} deletes a dipole and connects the half-edges of the same side. This move is represented in Figure~\ref{fig:dip_rem_unxod}. It is straightforward to see that dipole removals cannot increase the genus nor the grade.

\paragraph{Proof of Proposition~\ref{claim_mel}.\\}

First, it is easy to prove, by induction on the number of vertices, that melonic maps all have vanishing genus and grade by checking that the melonic insertion leaves the genus and grade unchanged and that the elementary melon has vanishing genus and grade.
\begin{equation}
\includegraphics[scale=0.45, valign=c]{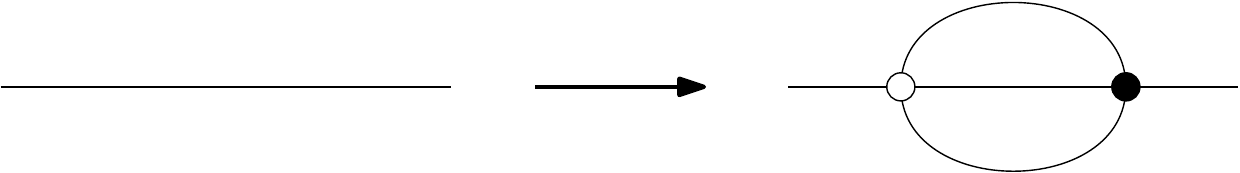}
\label{fig:melon_unxod}
\end{equation}

To prove that those are the only such maps, we also use an induction on the number of vertices $v$ of the map. It is obviously true at $v=2$. Let $v'> 2$ and assume that for all $v < v'$, there are no maps with vanishing genus and grade which are not melonic.

\medskip

Let $\phi_{2n}$ be the number of straight cycles of length $2n$ in a map $\mathcal{M}$ of genus $g$ and grade $l$. We first show that $\phi_{2}>0$ for maps of vanishing genus and grade. In an arbitrary map, the total number of straight cycles is $\phi = \sum_{n\geq 1} \phi_{2n}$ and since each edge belongs to exactly one cycle one has $E = \sum_{n\geq 1} 2n\phi_{2n}$. Moreover, $V=E/2$ since vertices have degree 4. Plugging those relations into the second Equation of \eqref{GenusAndGradeU(N)} gives
\begin{equation}
\sum\limits_{n\geq 1} (n-2)\phi_{2n} = l - 2 -2g.
\label{eq:cycle_length}
\end{equation}
This shows that such a map $\mathcal{M}$ has $\phi_2>0$, i.e. at least one O-dipole which we denote $c$. If $\mathcal{M}$ with $v'$ vertices has vanishing genus, it is planar and thus $c$ splits the map in two regions,
\begin{equation}
\mathcal{M} = \includegraphics[scale=.6, valign=c]{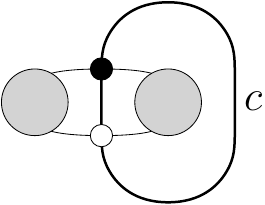}
\end{equation}
where each grey blob is {\it a priori} an arbitrary 2-point map. Removing the O-dipole disconnects the map into two connected components with fewer vertices, which are melonic by hypothesis. Reconstructing $\mathcal{M}$ as above, we find that it is melonic.

\paragraph{Chains.\\}

A chain is either an isolated dipole or a sequence of dipoles glued side by side, as in Equation~\eqref{eq:chain} where each dipole-vertex can now be a U-dipole or an O-dipole. The length of a chain is its number of dipoles. A chain of length $\ell$ contains chains of all lengths $1\leq \ell'\leq \ell$. A chain is said to be \emph{maximal} in $\mathcal{M}$ if it is not contained in a longer chain. Note that maximal chains are \emph{vertex-disjoint}. As a remark, this would not be true if non-isolated dipoles were allowed as chains as in Figure~\ref{fig:VertexJointDipoles}.

\medskip

We distinguish between chains made of O-dipoles only (resp. U-dipoles) named \emph{O-chains} (resp. \emph{U-chains}), which can have length greater than or equal to 1, and other chains which are called \emph{broken chains} and are made of at least two different dipoles. The following proposition is easily proved.

\begin{proposition}
Changing the length of a chain does not change the genus nor the grade (both in vacuum and 2-point graphs).
\end{proposition}

Let $C_O$, $C_U$ and $B$ be the generating functions of O-chains, U-chains and broken chains respectively, decorated with melons. They are given by
\begin{align}
C_O(t) &= \frac{U}{1-U}, \\
C_U(t) &= \frac{2U}{1-2U}, \\
B(t) &= \frac{3U}{1-3U} - \frac{U}{1-U} -\frac{2U}{1-2U}  = \frac{(4-6U)U^2}{(1-U)(1-2U)(1-3U)}.
\end{align}
In the following, we will further need to distinguish between O-chains of even and odd lengths. Their respective generating functions are
\begin{equation}
C_{O,e}(t) = \frac{U^2}{1-U^2}, \qquad
C_{O,o}(t) = \frac{U}{1-U^2}
\end{equation}

\subsubsection{Singularity analysis}

Since we only have one coupling constant, the singularity analysis is much simpler than in the previous model. The leading singularity of $M(t)$ and of $U(t) = M(t)-1$ both occur at $t_c = \frac{3^3}{4^4}$. It is such that $U(t_c) = \frac{1}{3}$. The denominators of the series $C_O(t)$ and $C_U(t)$ remain finite at this point. However, that of $B(t)$ blows up. Near $t_c$ it behaves like
\begin{equation}
B(t) \underset{t \rightarrow t_c}{\sim} \frac{1}{\sqrt{\frac{8}{27}}\sqrt{1-\frac{t}{t_c}}}.
\end{equation}

Therefore the relevant singularity is associated with broken chains as they are the only ones to diverge at the critical point.

%%%%%%%%%%%%%%%%%%%%%%%%%%
\subsection{Schemes of the \texorpdfstring{$U(N)\times O(D)$}{U(N)xO(D)} multi-matrix model}

\emph{Schemes} are maps obtained by
\begin{itemize}
\item Replacing every melonic submap with an edge,
\item Replacing every maximal chain with a \emph{chain-vertex}, i.e. a 4-point vertex where we forget the length of the chain. A chain-vertex keeps track of the pairing of the half-edges by drawing them on the same side of a box, and separating both sides with fat edges along the box,
\begin{equation}
\includegraphics[scale=.5]{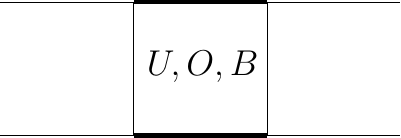}
\end{equation}
The labels $U,O,B$ correspond to chain-vertices replacing a U-chain, an O-chain or a broken chain.
\end{itemize}

\medskip

From this definition, it is clear that the reduction of a map to its scheme does not change its genus and grade. It thus makes sense to define $\mathbb{MS}_{h,l}$ the set of schemes of genus $h$ and grade $l$. The following Theorem is the specialization of Theorem~\ref{th:sch} to this model.

\begin{theorem} \label{thm:SchemesUNxOD}
The set $\mathbb{MS}_{h,l}$ is finite.
\end{theorem}

The proof of Theorem~\ref{thm:SchemesUNxOD} uses the similar $2$-step strategy as for the quartic $O(N)^3$ tensor model. We show that
\begin{enumerate}
    \item Schemes of fixed genus $g$ and grade $l$ have finitely many chain-vertices via Lemma~\ref{lem:max_ns}
    \item There are finitely many maps of genus $g$ and grade $l$ with $k$ isolated dipoles via Proposition~\ref{prop:sch_ns}.
\end{enumerate}
This proves Theorem \ref{thm:SchemesUNxOD} by observing that there is a bijection between schemes of genus $g$ and grade $l$ and maps of the same genus and grade whose chains all have length 1. One can thus apply the result of step $1$, then replace chain-vertices with chains of length $1$ and apply the result of step $2$.

\medskip

We start by showing that the number of chain-vertices is bounded by the genus and the grade. We have the following Lemma which can be found in~\cite[Lemma $4$]{BeCa} which can be adapted to the graph of this model.

\begin{lemma}
If a scheme $\cS$ has a non-separating chain-vertex, then it is possible to construct a scheme $\cS'$ of the same genus and grade with more chain-vertices.
\label{lem:max_ns}
\end{lemma}

\begin{proof}
Let $\mathcal{T}\subset \mathcal{I}(\cS)$ be a spanning tree in the skeleton graph of $\cS$. The edges in the complement $\mathcal{I}(\cS)\setminus\mathcal{T}$ correspond in $\cS$ to a set of non-separating chain-vertices. Removing each of them decreases the genus by 1 or the grade by 4. After all those removals, one obtains a scheme $\cS_r$ of genus $g-q_1$ and grade $l-4q_2$, such that $q_1+q_2$ is the number of non-separating chain-vertices which have been removed. It is then possible to attach to any (separating) remaining chain-vertex of $\cS_r$ a plane binary tree $\tilde{\mathcal{T}}$ with $q_1+q_2$ leaves and whose internal vertices correspond to components of vanishing genus and grade. Corresponding to the leaves of $\tilde{\mathcal{T}}$, one can choose, thanks to Lemma~\ref{thm:SkeletonGraph}, components that make up for the loss of genus and grade induced by the removals of the non-separating chain-vertices. More formally, if $\ell$ is a leaf of $\tilde{\mathcal{T}}$, we can choose a component $\mathcal{M}_\ell$ with genus $g_\ell$ and grade $\ell$, and using Lemma~\ref{thm:SkeletonGraph}, such that $\sum_\ell g_\ell = q_1$ and $\sum_\ell l_\ell = q_2$. This gives a new skeleton graph which is a tree and which we denote $\mathcal{T}'$ and a (non-unique) scheme $\cS'$, such that $\mathcal{I}(\cS') = \mathcal{T}'$, whose degree and grade are the same as those of $\cS$ but has more chain-vertices.

\medskip

As an illustration, we describe the attachment of $\tilde{\mathcal{T}}$ to $\mathcal{T}$ in the case $q_1+q_2=1$ (i.e. a single non-separating chain-vertex in $\cS$) on Figure~\ref{fig:non_sep_non_dom}.
\begin{figure}[!h]
    \centering
    \includegraphics[scale=0.6]{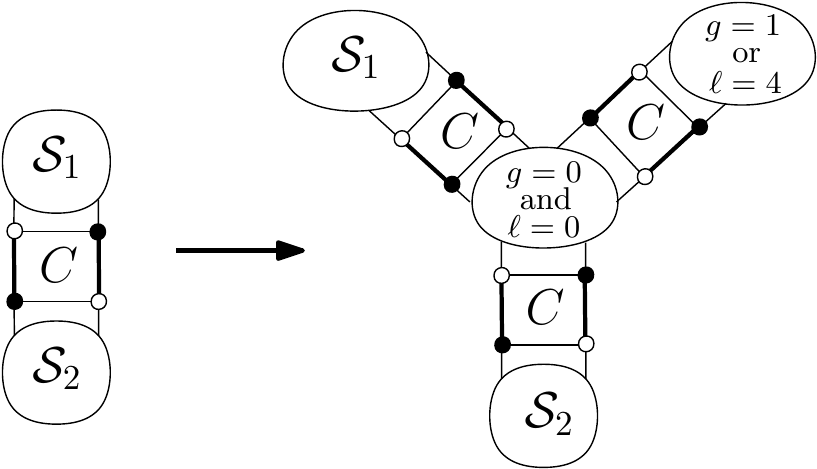}
    \caption{\label{fig:non_sep_non_dom} On the left is $\cS_r$ with a separating chain-vertex and a move attaching a binary vertex to it on the right hand side. It restores the genus and grade of $\cS$ but has more chain-vertices.}
\end{figure}
\end{proof}

From this lemma, we get the following bound on the number of chain-vertices of a scheme.
\begin{proposition}
A scheme of genus $g$ and grade $l$ has at most $2(g+l)-1$ chain-vertices. This bounded is saturated when all chain-vertices are separating.
\label{prop:sch_chain}
\end{proposition}

\begin{proof}
By contraposition of the Lemma~\ref{lem:max_ns}, schemes that have a maximal number of chains cannot have non-separating chains. Therefore their skeleton graph must be a tree $\mathcal{T}$. Using Lemma~\ref{lem:UNxOD}, we know that leaves of $\mathcal{T}$ cannot have vanishing genus and grade, therefore all leaves of $\mathcal{T}$ must have genus at least $1$ or grade at least $2$. Therefore, maximizing the number of chains in the scheme is a simple linear problem where one has to maximize the number of edges of a tree given a bound on the number of leaves given by $g+l$. The solution is to have as many leaves as possible, here $g+l$, and inner vertices of degree exactly 3 with vanishing genus and grade. A binary tree with $k$ leaves has $2k-1$ edges, thus the result follows.
\end{proof}

The second part of the proof of Theorem \ref{thm:SchemesUNxOD} comes from the following proposition.

\begin{proposition}
The set of maps of genus $g$ and grade $l$ with $k$ isolated O-dipoles is finite.
\label{prop:sch_ns}
\end{proposition}

\begin{proof}
The strategy is similar to the strategy adopted in Subsection~\ref{ssec:proof_Lem_sec} for the $O(N)^3$ tensor model. We want to show that there exist functions $b_{2n}(g,l,k)$ such that
\begin{equation}
    \phi_{2n} \leq b_{2n}(g,l,k),
\end{equation}
for all $n\geq 1$, where we recall that $\phi_{2n}$ is the number of straight cycles of length $2n$.

\medskip

Observe that equation \eqref{eq:cycle_length}, i.e. $\sum_{n\geq 1} (n-2)\phi_{2n} = l - 2 -2g$, bounds the number of straight cycles of length greater than or equal to $6$ for fixed values of $g$ and $l$ and $\phi_2$. Therefore, two bounds remain to be obtained.
\begin{itemize}
    \item Bounding $\phi_2\leq b_2(h,l,k)$ in terms of the genus, grade and number of isolated O-dipoles. This implies bounds $\phi_{2n}\leq b_{2n}(h,l,k)$ for all $n\geq3$ thanks to \eqref{eq:cycle_length}.
    \item Bounding $\phi_4\leq b_4(h,l,k)$ similarly. This must be done independently because $\phi_4$ does not appear in \eqref{eq:cycle_length} and there could thus, {\it a priori}, exist maps with an arbitrarily large number of straight cycles of length 4.
\end{itemize}

\paragraph{Straight cycles of length 2\\}
Those cycles are O-dipoles by definition.
\begin{itemize}
    \item Either it is an isolated O-dipole and there are $k$ of them at most.
    \item Either it is a non-isolated O-dipole.
\end{itemize}
%In the latter case, notice as shown in Figure \ref{fig:VertexJointDipoles} that two non-isolated O-dipoles form a topological minor of genus 1 and therefore there can be at most $g$ such disjoint configurations.
In the latter case, notice that a non-isolated dipole is non-separating, so its removal decreases the genus by $1$ and we can proceed by induction on the genus. If $\mathcal{M}$ has genus zero (i.e. it is a planar map) then it cannot have non-isolated O-dipoles (they would form a topological minor of genus 1, see the right of Figure \ref{fig:VertexJointDipoles}). Hence we can define $b_2(0,l,k) = k$, which satisfies $\phi_2 \leq b_2(0,l,k)$ when the map $\mathcal{M}$ is planar. Now if $\mathcal{M}$ has genus $g>0$ and a non-isolated O-dipole, the latter can be removed and this decreases the genus by one and gives a connected map $\mathcal{M}'$. This operation cannot increase the grade (in fact $l(\mathcal{M}') = l$ or $l(\mathcal{M}') = l-4$). It can form at most $2$ non-isolated O-dipoles (if the two new edges in $\mathcal{M}'$ belong to non-isolated dipoles) or $1$ isolated O-dipole. Since $\mathcal{M}'$ has at most $k+1$ isolated O-dipoles, we can define
\begin{equation}
b_2(g,l,k) = 2+\max_{\substack{k' \leq k+1\\ l'\leq l}} b_2(g-1,l',k').
\end{equation}
By construction, this function satisfies $\phi_2 \leq b_2(g,l,k)$ for a map of genus $g$ and grade $l$. This shows that the number of O-dipoles is bounded in terms of the genus, grade and number of isolated O-dipoles.

\medskip

This in turn proves that there is a finite number of straight cycles of length larger than or equal to 6.

\paragraph{Straight cycles of length $4$\\}

This case is proved via a similar method to the previous one but it is more involved. Instead of proving directly that the number of cycles of length $4$ is bounded, we show that any cycle of length $4$ is within bounded distance of a submap which carries a non-zero genus. We start by recalling some definitions and the Lemma~\ref{thm:TopMinors} which carries the main idea of the proof.

\begin{definition}[Graph distance]
In a map $\mathcal{M}$, two vertices $v$ and $v'$ are at distance $k$ if there is a sequence of edges $(e_1,..,e_k)$ such that $e_1$ starts at $v$, $e_k$ ends at $v'$ and for $i \in \{1,..k-1\}$ the end point of $e_i$ is the starting point of $e_{i+1}$.
\end{definition}

\begin{definition}[Topological minor]
A topological minor of genus $h$ of a map $\mathcal{M}$ is a submap $\mathcal{M}' \subset \mathcal{M}$ such that $\mathcal{M}'$ has genus $h$.
\end{definition}

We have the following useful lemma (which could have in fact been used above to deal with cycles of length 2 too).

\begin{lemma} \label{thm:TopMinors}
Let $\mathfrak{m}$ be a connected topological minor of genus $h>0$. Then the number of copies of $\mathfrak{m}$ which can occur in $\mathcal{M}$ is bounded linearly by the genus $g$ of $\mathcal{M}$.
\end{lemma}

\begin{proof}
We say that a set of copies of $\mathfrak{m}$ are \emph{separated} in $\mathcal{M}$ if they are vertex-disjoint in $\mathcal{M}$, and it is maximal if it cannot be properly contained in another independent set.

\medskip

Let $\operatorname{Sep}_{\mathfrak{m}}(\mathcal{M})$ be a maximal separated set for $\mathfrak{m}$ in $\mathcal{M}$. Since $h>0$ there can only be a finite number of disjoint copies of $\mathfrak{m}$, at most $\lfloor g/h\rfloor$, so $\lvert\operatorname{Sep}_{\mathfrak{m}}(\mathcal{M})\rvert\leq g/h$. Any other occurrence of $\mathfrak{m}$ shares a vertex with one of $\operatorname{Sep}_{\mathfrak{m}}(\mathcal{M})$. Since $\mathfrak{m}$ has a finite number of vertices $v(\mathfrak{m})$ and they have finite valency $4$, it comes that there is at most $4^{v(\mathfrak{m})-1} g/h$ copies of $\mathfrak{m}$ in $\mathcal{M}$.
\end{proof}

To treat the case of straight cycles of length 4, we are going to show that in all configurations where such a cycle can appear, we can find a submap of genus (at least) $1$ close to it. After counting how many non-separated sets of such a submap can appear, we can then give a very rough bound on the total number of such structures in the graphs using~Lemma~\ref{thm:TopMinors}, which is sufficient for our goal since we only have to show that it is finite. Let $\cC$ be such a cycle and let $V_{\text{exc}}$ the set of vertices that belong to at least one straight cycle of length 2 or 6 or larger, i.e. a cycle of a type which can only appear finitely many times in the graph.

\subparagraph{$\cC$ has a vertex in $V_{\text{exc}}$\\}
If $\cC$ has a vertex that belongs to $V_{\text{exc}}$, then it has all its vertices at a distance at most $2$ of $V_{\text{exc}}$. Since vertices all have valency $4$, for any vertex $v$ there are at most $17$ distinct vertices at a distance at most $2$ of $v$. Thus there are at most $17|V_{\text{exc}}|$ straight cycles of length $4$ with a vertex in $V_{\text{exc}}$.

\subparagraph{$\cC$ has no vertices in $V_{\text{exc}}$\\}
This means that $\cC$ only intersects straight cycles of length exactly $4$. In the following, we will assume $\cC$ is non-self-intersecting. The case where $\cC$ is self-intersecting can be tackled similarly. We denote  $v_1$ to $v_4$ the vertices of $\cC$ in any cyclic order. The cycle $\cC'$ incident to $v_1$ that is not $\cC$ also has length $4$. We distinguish cases depending on how $\cC$ and $\cC'$ intersect.

\subparagraph{$\cC$ and $\cC'$ intersect only at $v_1$.\\} In this case, $\cC\cup\cC'$ forms a topological minor of genus $1$. From Lemma \ref{thm:TopMinors}, there is a bound in the genus on the number of such pairs, hence on the number of cycles of length $4$ which intersect another cycle of length $4$ exactly once. 

\subparagraph{$\cC$ and $\cC'$ intersect at least at $v_1$ and $v_2$.\\} By symmetry, it means they intersect at two vertices of distinct colors.
\begin{itemize}
    \item If they intersect at $v_1$ and $v_2$ only, then they either form an isolated U-dipole and there are at most $k$ of them, or a topological minor of genus 1 and we conclude with Lemma \ref{thm:TopMinors}. This is illustrated in Figure \ref{fig:v1v2}.
    \begin{figure}
        \centering
        \includegraphics[scale=.65]{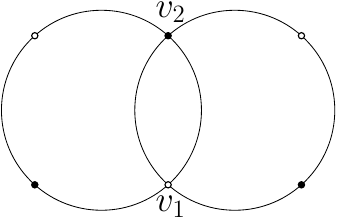}
        \hspace{1cm}
        \includegraphics[scale=.65]{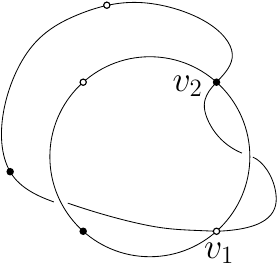}
        \caption{There are two possibilities for $\cC$ and $\cC'$ to intersect at exactly two vertices of different colors.}
        \label{fig:v1v2}
    \end{figure}
    \item If they intersect at $v_1$, $v_2$, $v_3$, or by symmetry at any triple of vertices, we have one of the situations depicted in Figure \ref{fig:v1v2v3}. In all of them, $\cC\cup \cC'$ forms a topological minor of genus 1. From Lemma \ref{thm:TopMinors}, we conclude that there is a bounded number of those minors.
    \begin{figure}
        \centering
        \includegraphics[scale=.65]{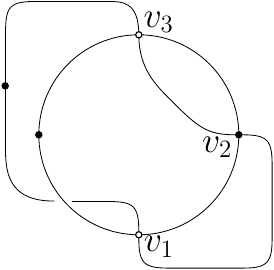}
        \hspace{1cm}
        \includegraphics[scale=.65]{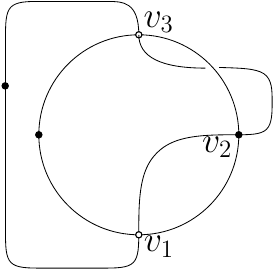}
        \hspace{1cm}
        \includegraphics[scale=.65]{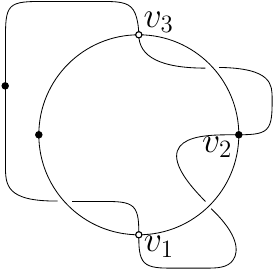}
        \caption{There are three possibilities for $\cC$ and $\cC'$ to intersect at exactly three vertices, up to symmetry.}
        \label{fig:v1v2v3}
    \end{figure}
    \item If they intersect on all four vertices, then only a finite number of maps can be drawn.
\end{itemize}

\subparagraph{$\cC$ and $\cC'$ intersect at $v_1$ and $v_3$ only.} By symmetry, it corresponds to the case where they intersect at exactly two vertices of the same color. We distinguish two cases depending on whether $\cC'$ is self-intersecting or not.
\begin{itemize}
\item If $\cC'$ is self-intersecting, we can check, see Figure \ref{fig:4nsi_with_si}, that $\cC\cup\cC'$ either forms a topological minor of genus 1 or of genus $2$. We conclude again from Lemma \ref{thm:TopMinors} that there is a finite number of such configurations at fixed genus. 

\begin{figure}
        \centering
        \includegraphics[scale=.85]{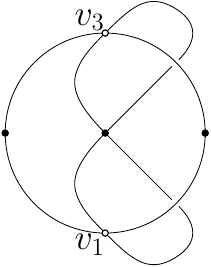}
        \hspace{1cm}
        \includegraphics[scale=0.45]{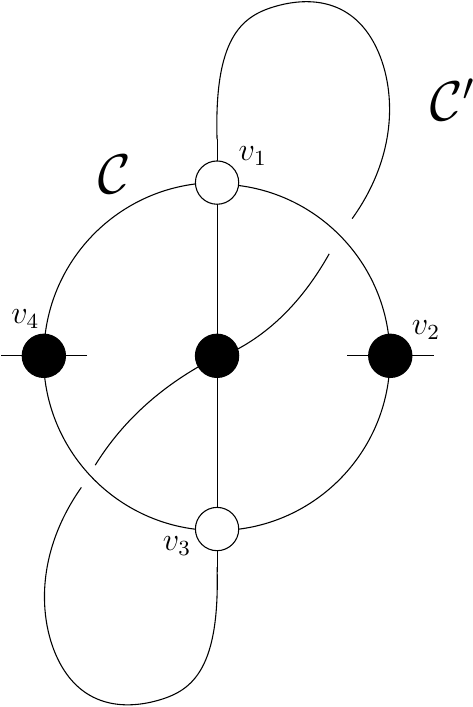}
        \caption{There are two possibilities for $\cC$ non-self-intersecting and $\cC'$ self-intersecting to intersect at two vertices of the same color, up to symmetry.}
        \label{fig:4nsi_with_si}
\end{figure}

\item The last situation we have yet to treat is when $\cC'$ is non-self-intersecting. Then denote $v_{13}$ and $v_{31}$ the two other vertices of $\cC'$,
    \begin{equation} \label{LastCase}
    \includegraphics[scale=.65,valign=c]{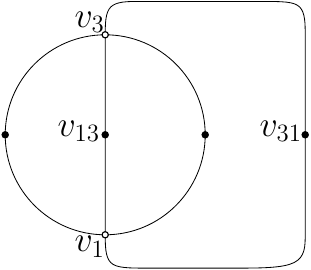} \hspace{1cm} \includegraphics[scale=.65,valign=c]{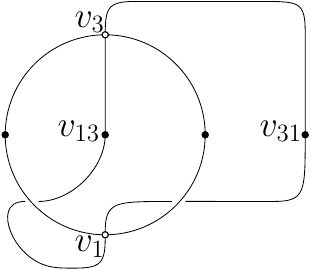}
    \end{equation}
    In the second case, $\cC\cup\cC'$ forms a topological minor of genus 1 and again from Lemma \ref{thm:TopMinors} there is a bounded number of such configurations.

    \medskip 
    
    We thus focus on the first case of \eqref{LastCase} and consider the straight cycle $\cC''$ through $v_{13}$. If $\cC'\cup \cC''$ forms a configuration which we have already bounded, so is $\cC\cup \cC'\cup \cC''$. Therefore the only situation to be treated is when $\cC'\cup\cC''$ has the same structure as in the first case of \eqref{LastCase}, i.e. $\cC''$ goes through $v_{13}$ and $v_{31}$ and $\cC'\cup\cC''$ is planar.
    This situation is represented in Figure~\ref{fig:4nsi_triple}.
    \begin{figure}[!h]
        \centering
        \includegraphics[scale=0.5]{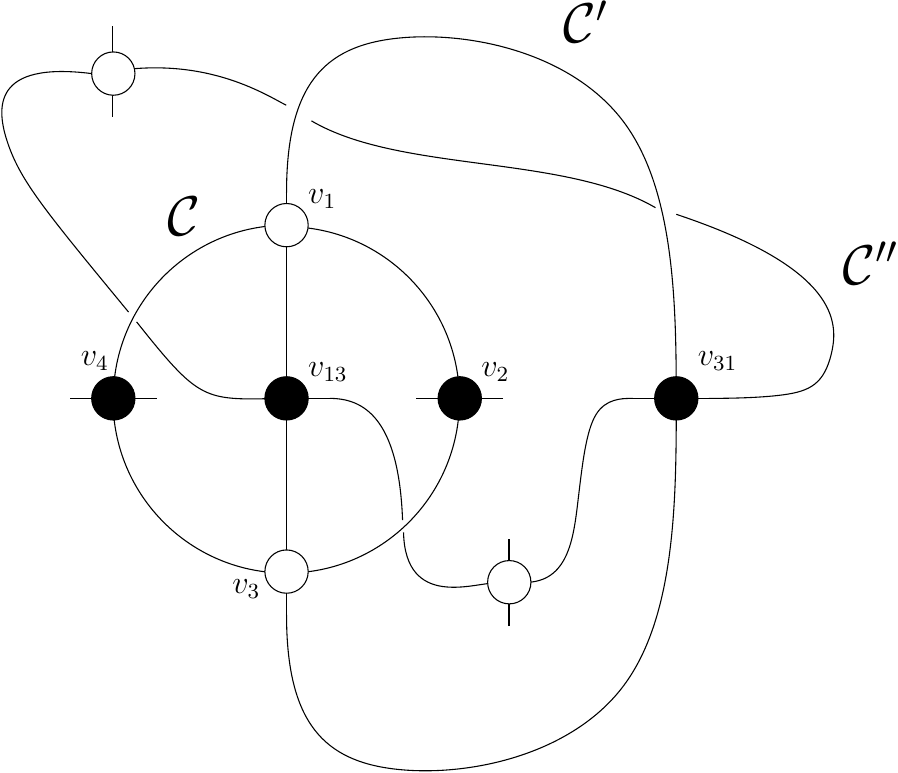}
        \caption{$\cC\cup\cC'\cup\cC''$ forms a topological minor of genus $2$.}
        \label{fig:4nsi_triple}
    \end{figure}
    In this configuration,  $\cC\cup\cC'\cup\cC''$ forms a topological minor of genus $2$, and it is concluded from Lemma \ref{thm:TopMinors} that there is a bounded number of such configurations.
\end{itemize}
We have exhausted all the possible cases for a straight cycle of length 4 and find that there is indeed a bound $\phi_{2n}\leq b_{2n}(g,l,k)$ for $n=2$. This was the last missing bound, so we have them for all $n\geq 1$, which proves Proposition \ref{prop:sch_ns}. Proposition~\ref{prop:sch_ns} together with Proposition~\ref{prop:sch_chain} finally proves Theorem~\ref{thm:SchemesUNxOD}. 
\end{proof}

\subsection{Skeleton graph of the schemes}

This subsection introduces the skeleton decomposition for the $U(N) \times O(D)$ model and establishes its main properties.

\subsubsection{Removals of chains and dipoles}

We say that a dipole in $\mathcal{M}$ is \emph{separating} if its removal disconnects $\mathcal{M}$ into two connected maps $\mathcal{M}_1, \mathcal{M}_2$ and \emph{non-separating} otherwise, and similarly for chain-vertices.

\vspace{10pt}
\paragraph{Dipole removals.\\}

We consider the dipole removals represented in Figure~\ref{fig:dip_rem_unxod}. The following analysis parallels that made in the previous model, see Section~\ref{ssec:dip_chain_rem_ON3}.

\subparagraph{Separating dipole removal.} Denote $\Delta Q = Q(\mathcal{M}_1) + Q(\mathcal{M}_2) - Q(\mathcal{M})$ the variation of any quantity $Q$ through the removal. Recall that the formula \eqref{GenusAndGradeU(N)} holds for each connected map. During the removal one finds for both O- and U-dipoles that $\Delta F=\Delta \phi=0$, $\Delta E=-4$, $\Delta V=-2$ and the variation of the number of connected components is $\Delta C=1$. This gives
\begin{equation}
    \Delta h = \Delta l = 0.
\end{equation}

\subparagraph{Non-separating O-dipole removal.} Let $\Delta Q = Q(\mathcal{M}') - Q(\mathcal{M})$ be the variation of $Q$ through the move. The number of faces is unchanged, i.e. $\Delta F=0$, while $\Delta \phi = 0,-2$ depending on the structure of the straight cycles which are incident to the dipole. This gives
\begin{equation} \label{ODipoleRemoval}
    \Delta h=-1, \qquad \Delta l = 0 \text{ or } -4.
\end{equation}

\subparagraph{Non-separating U-dipole removal.} In contrast with the previous case of O-dipoles, it preserves the straight cycles incident to it, i.e. $\Delta \phi=0$, while $\Delta F=0,-2$ depending on the structure of the faces which are incident to the dipole. This gives
\begin{equation} \label{UDipoleRemoval}
    \Delta h = -1\text{ or }0, \qquad \Delta l = -2\text{ or }-4,
\end{equation}

\paragraph{Chain removals.\\}

A chain removal corresponds to the following move
\begin{equation}
\includegraphics[scale=.56, valign=c]{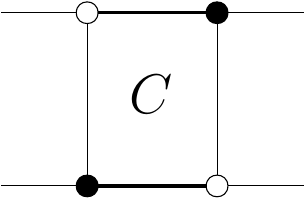} \qquad \to \qquad \includegraphics[scale=.55, valign=c]{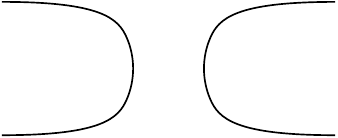}
\end{equation}
It can be studied by removing one dipole from the chain and then removing the resulting melonic 2-point functions. The removal of a separating chain is obviously the same as that of a separating dipole.

\subparagraph{Non-separating chains.}
O-chains and U-chains are sequences of isolated O- and U-dipoles. Thus their removal is exactly the same as the removal of a non-separating O- or U-dipole.

\subparagraph{Non-separating broken chains.} By definition a broken chain has at least one O-dipole and one U-dipole. Therefore the removal of a broken chain has to be a special case of both \eqref{ODipoleRemoval} and \eqref{UDipoleRemoval}. This directly gives
\begin{equation} \label{BrokenChainRemovalU(N)}
    \Delta h=-1,\qquad \Delta l=-4.
\end{equation}
This is the same result as for broken chain removals in the previous model \eqref{BrokenChainRemoval}.

\subsubsection{Skeleton graph of a scheme}

\begin{figure}
\centering
\includegraphics[scale=0.55]{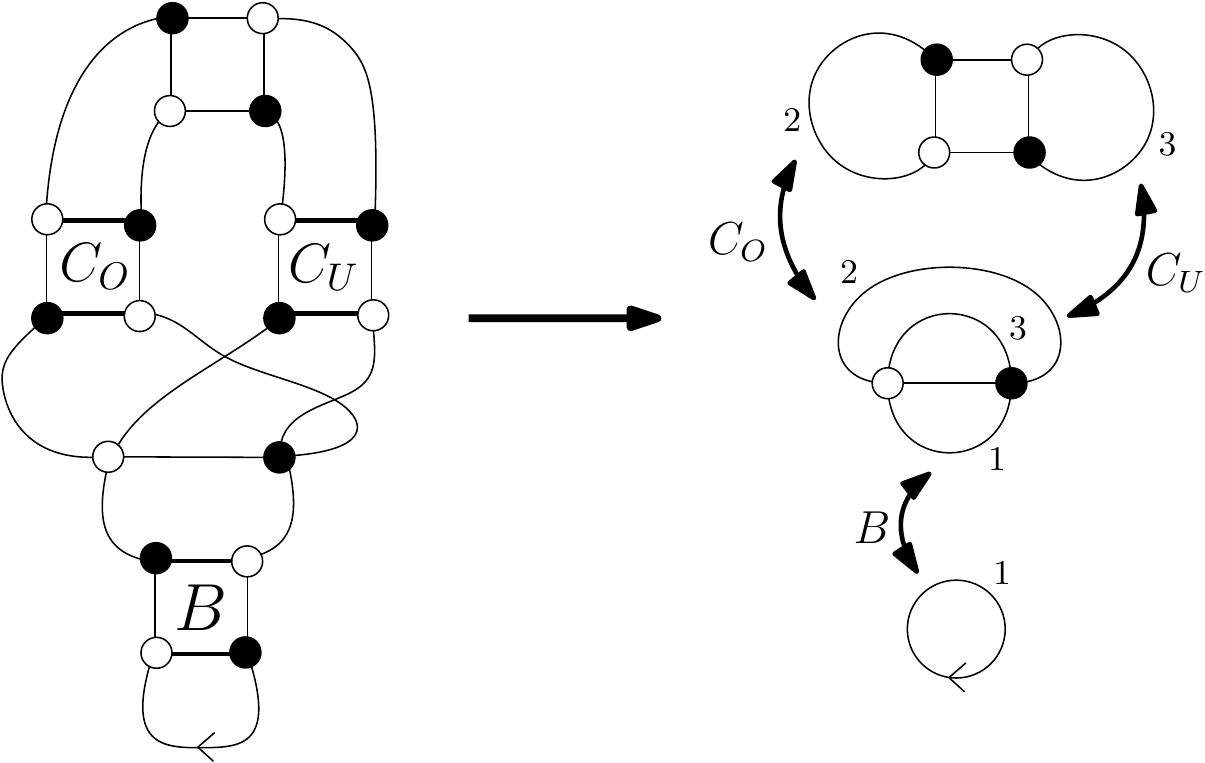}
\caption{A scheme of the $U(N) \times O(D)$ model and its skeleton graph.}
\label{fig:skel_unxod}
\end{figure}

Let $\cS$ be a scheme with genus $g$ and grade $l$. Similarly to what has been done for the $O(N)^3$ model, we introduce the \emph{skeleton graph} $\mathcal{I}(\cS)$ of a scheme $\cS$.
\begin{itemize}
\item We call the components of $\cS$ the connected components obtained after removing all chain-vertices of $\cS$. In each component, we mark the edges created by the removals.
\item The vertices of $\mathcal{I}(\cS)$ are the components of $\cS$.
\item Two vertices of $\mathcal{I}(\cS)$ are connected by an edge if the two components are connected by a chain-vertex in $\cS$. Each edge is labeled by the type of the chain connecting them.
\end{itemize}
An example of a skeleton graph is shown in Figure~\ref{fig:skel_unxod}.

\begin{lemma}\label{lem:UNxOD}
The skeleton graph $\mathcal{I}(\cS)$ of a scheme $S$ satisfies the following properties:
\begin{enumerate}
\item Any vertex of $\mathcal{I}(\cS)$ corresponding to a component of vanishing genus and grade, and not carrying the root edge of $\cS$, has degree at least $3$.
\item If $\mathcal{I}(\cS)$ is a tree, then the genus and grade of $\cS$ are split among the components.
\end{enumerate}
\end{lemma}
This lemma is completely analogous to Lemma \ref{thm:SkeletonGraph} for the $O(N)^3$ model.

\begin{proof}
\begin{enumerate}
\item If a component of vanishing genus and grade not carrying the root edge of $\cS$ has valency $1$ in $\mathcal{I}(\cS)$, then it is a melonic component, which is not possible for a scheme. Similarly, if it has valency $2$ in $\mathcal{I}(\cS)$ then it is a chain. This implies that the (two) chains incident to this component are not maximal, which is not possible in a scheme.
\item As we have shown, removing a separating chain-vertex in $\cS$ splits the degree among the components of $\cS$. Saying that $\mathcal{I}(\cS)$ is a tree is equivalent to saying that all chain-vertices of $\cS$ are separating, therefore the genus and grade of $\cS$ split among the components.
\end{enumerate}
\end{proof}

\subsection{\label{UNOD:dom_scheme} Identification of dominant schemes and double scaling limit}

Via Theorem~\ref{thm:SchemesUNxOD}, we know that for a given genus and grade $(g,l)$, there are finitely many schemes. Thus all singularities of the generating function for $\mathbb{MS}_{g,l}$ come from the generating series of melons and chains. As we have shown, the only divergent objects at the critical point are the broken chains. Moreover, we restrict attention to dominant schemes of vanishing grade.

\medskip

Proposition~\ref{prop:sch_chain} gives the structure of the dominant schemes of genus $g$. All their chain-vertices are separating, i.e. they are schemes whose skeleton graphs are plane trees. In order to maximize the singularity, all chain-vertices are broken chain-vertices. Given Lemma \ref{thm:SkeletonGraph}, maximizing the number of chain-vertices is a linear problem on the degrees of the internal vertices of $\mathcal{I}(\cS)$ and their genus and grades, and on the genus and grades of the leaves. The solution is for $\mathcal{I}(\cS)$ to be a plane binary tree, whose internal vertices correspond to components of $\cS$ of vanishing genus and grade, and whose leaves correspond to components of genus 1 (and vanishing grade).

\medskip

By adapting the analysis of~\cite[Proposition $1$]{BeCa} to identify the components of genus 1 and the components of vanishing grade and genus, the following proposition is obtained.

\begin{proposition} \label{prop:dom_scheme_un}
A rooted dominant scheme of genus $g\geq1$ has $2g-1$ broken chain-vertices which are all separating. Such a scheme has the structure of a rooted\footnote{The root is a marked leaf.} binary plane tree where
\begin{itemize}
\item Each edge corresponds to a broken chain-vertex.
\item The root of the tree corresponds to the root of the map.
\item Each of the $g$ leaves is one of the following two graph
\begin{equation}
\includegraphics[scale=0.4]{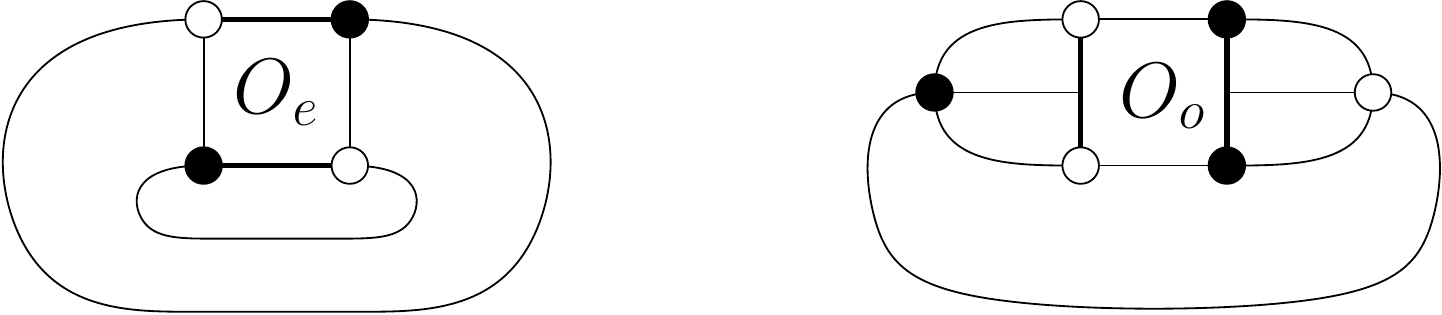} \nonumber
\end{equation}
\item {Each internal vertex corresponds to one of the four $6$-point submaps
\begin{equation}
\includegraphics[scale=0.80]{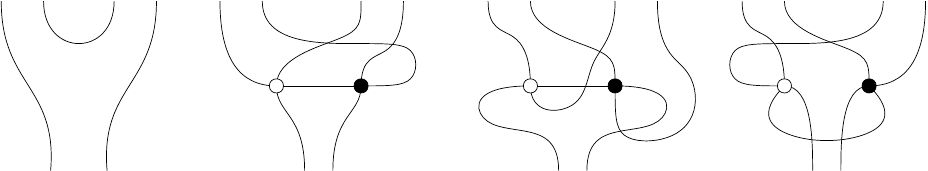} \nonumber
\end{equation}
with the convention that each bottom left vertex is related to a black vertex.}
\end{itemize}
\end{proposition}

From this proposition, the generating function of dominant schemes of genus $g$ is given by
\begin{align}
G_{dom}^g(t) = L^{2-2g}M(t) \Cat_{g-1} \left(C_{O,e}+t C_{O,o}\right)^g\left(1+6t\right)^{g-1} B(t)^{2g-1}
\end{align}
where we recall that $L=N/\sqrt{D}$ and that maps of genus $g$ scale like $L^{2-2g}$ in the large $N$, large $D$ expansion (see Theorem~\ref{thm:free_energy_U(N)}).

\subsubsection{Double scaling limit}

We recall that the critical point for $B(t)$ is $t_c = \frac{3^3}{4^4}$. Near this point we have:
\begin{align}
C_{O,e} = \frac{1}{8} + O(\frac{t}{t_c}), \qquad C_{O,o} = \frac{1}{24} + O(\frac{t}{t_c}) \\
U(t) \underset{t \rightarrow t_c}{\sim} \frac{1}{3} - \sqrt{\frac{8}{27}}\left(1-\frac{t}{t_c}\right)^{-\frac{1}{2}} 
\end{align}
Therefore it follows that
\begin{align}
B(t) \underset{t \rightarrow t_c}{\sim} \sqrt{\frac{27}{8}}\left(1-\frac{t}{t_c}\right)^{-\frac{1}{2}}
\end{align}
Thus near the critical point $t_c$ we find
{\small \begin{equation}
G_{\text{dom}}^g(t) \underset{t \rightarrow t_c}{\sim} L^{2-2g}\frac{4}{3}\Cat_{g-1} \left(\frac{1}{8}+t_c \frac{1}{24}\right)^g \left(1+6t_c\right)^{g-1} \left(\frac{27}{8}\frac{1}{\left(1-\frac{t}{t_c}\right)}\right)^{g-\frac{1}{2}}
\end{equation}}%
After introducing the double scaling parameter $\kappa$ defined as
\begin{equation}
\kappa^{-1} = L^2 \frac{1}{\left(\frac{1}{8}+t_c \frac{1}{24}\right)\left(1+6t_c\right)}\frac{8}{27}\left(1-\frac{t}{t_c}\right)
\label{eq:kappa_un}
\end{equation}
the following equality holds
\begin{align}
\left(1+3t\right)^{-1}\left(\frac{8}{27}\left(1-\frac{t}{t_c}\right)\right)^{\frac{1}{2}} = \frac{\kappa^{-\frac{1}{2}}}{L} \left(\frac{\frac{1}{8}+t_c \frac{1}{24}}{1+6t_c}\right)^{\frac{1}{2}}
\end{align}
so that in the double scaling limit, $G_{dom}^g(t)$, for $g>0$, contributes to the 2-point function as
\begin{equation}
G_{dom}^g(t) = \frac{4}{3}L\Cat_{g-1} \left(\frac{\frac{1}{8}+t_c \frac{1}{24}}{1+6t_c}\right)^{\frac{1}{2}}\kappa^{g-\frac{1}{2}}.
\end{equation}
For $g=0$, we simply have to account for melons thus contributing as $M(t_c) = \frac{4}{3}$.

\medskip 

Summing over the genus $g$, it is found that the series is convergent for $\kappa \leq \frac{1}{4}$. This gives the expression of the 2-point function in the double scaling limit,
\begin{align}
G_2^{DS}(t) &= L^{-2} \sum\limits_{g \geq 0}^{} G_{dom}^g(t) \\
&= 4/3 + \frac{4}{3}\frac{1}{L\kappa^{\frac{1}{2}}} \left(\frac{\frac{1}{8}+t_c \frac{1}{24}}{1+6t_c}\right)^{\frac{1}{2}} \sum\limits_g \Cat_{g-1} \kappa^g \nonumber \\
&= 4/3 + \frac{2}{3}\frac{1}{L}\left(\frac{\frac{1}{8}+t_c \frac{1}{24}}{1+6t_c}\right)^{\frac{1}{2}} \frac{1-\sqrt{1-4\kappa}}{\kappa^\frac{1}{2}}
\label{eq:GDS_U(N)xO(D)}
\end{align}

\section{Conclusion}
\label{sec:concl_DS}

The principal tool that allowed us to identify the dominant schemes is the implementation of a scheme decomposition in these tensor models. While they are dependent on precise details of the structure of the model, that is the symmetry group of the model and specific choice of bubbles, the scheme decomposition seems to hold for a wide variety of tensor models e.g. for some multi-matrix models and some tensor models with mixed index symmetry. Moreover, some features of the scheme decomposition are similar in all three models studied here.

\medskip

The double-scaling limit of tensor models is dominated by graphs that have a maximal number of broken chains. The dominant schemes are then given by trees whose edges correspond to broken chains. The contribution of the dominant schemes corresponds to a sum over decorated trees and can then be performed explicitly. Thus, the expression of the $2$-point function in the double scaling limit is convergent for $\kappa \leq \frac{1}{4}$ (where the definition double-scaling parameter $\kappa$ depends on the model), contrary to the two-dimensional matrix case where the double scaling series is not summable~\cite{Eynard_book}.

\medskip

The combinatorial methods employed to identify the double scaling limit of the $2$-point function can be generalized in a straightforward way to compute the double scaling limit of the $2r$-point function for $r \geq 1$, similarly as in~\cite{TaGu} where it was implemented for the multi-orientable tensor model.

%% file: Chapters/Melons_tens_vect.tex
\chapter{Amit-Roginsky from the Boulatov model} % Main chapter title
\label{Chap:melons_th}

In the last decade, alongside the development of $0$-dimensional tensor field theories, the Sachdev-Ye-Kitaev (SYK) model~\cite{Ki_SYK} attracted a lot of interest in the high energy physics community for its connection with holography~\cite{Ros_SYK,GrRo_SYK,MaSta_SYK} and black holes dynamics~\cite{BH_SYK}. This model is fermionic field theory in $0+1$-dimensions of $N$ Majorana fermions coupled by a random Gaussian coupling constant $J_{i_1\dotsc i_q}$ for some integer $q>2$. Albeit the tensor appearing in this model is of a different nature than in tensorial field theory, the SYK model admits a large $N$ expansion that is also dominated by melonic graphs~\cite{BNT_SYK}. The similarities in the nature of their $\frac{1}{N}$-expansion motivated the search for common properties between the SYK model and tensor field theories~\cite{Wi_SYK,KleTa_SYK}, establishing links between tensorial field theory and randomly coupled vector ones.

\medskip

In this chapter, we exhibit a similar connection between tensor and vector field theory with melonic leading order in their large $N$ expansion. More precisely, we show that the \emph{Amit-Roginski model} can be obtained as a perturbation over particular classical solutions of the \emph{Boulatov model} and give sufficient conditions on the form of the solution and perturbation to yield this property.

\section{The Boulatov model}
\label{sec:Boulatov_model}

The Group Field Theories (GFT) are field theories whose dynamical field depends on $n$ points $g_i$ of a Lie group $G$. The group elements $g_i$ can be interpreted as holonomies, that is as the parallel transport of a $G$-vector bundle. The Boulatov model~\cite{Boul_OG} is a $3$D GFT model whose field is a function $T:SU(2)^{3}\to \mathbb{C}$. The field is required to be invariant under global right multiplication by an element of $SU(2)$

\begin{equation}
    T(g_1h,g_2h,g_3h)=T(g_1,g_2,g_3)\hspace{10pt}\forall h\in SU(2),
    \label{right_inv}
\end{equation}
and to satisfy the reality condition
\begin{equation}
\label{reality}
    T(g_1,g_2,g_3)=\bar{T}(g_3,g_2,g_1),
\end{equation}
where $\bar{T}$ denotes the complex conjugate of $T$.
%The dynamics of the field $T$ is governed by a non-local quartic action~\cite{Boulatov:1992vp}
The action of the Boulatov model is non-local. It writes
\begin{align}
\label{eq:actionTgroup}
          S[T]&=\frac{\mu^2}{2}\int d g_1 d g_2 d g_3T(g_1,g_2,g_3)T(g_3,g_2,g_1)  \\
      &-\frac{\lambda}{4!}\int \prod_{i=1}^6 d g_i T(g_1,g_2,g_3) T(g_3,g_5,g_4)T(g_4,g_2,g_6)T(g_6,g_5,g_1), \nonumber
\end{align}
where $\mu$ is the mass of the field and $\lambda$ the coupling constant.

\medskip

Applying the Euler-Lagrange equation to this action, we obtain the classical equation of motion for the field $T(g_1,g_2,g_3)$. They are
\begin{equation}
      \mu^2T(g_3,g_2,g_1)=\frac{\lambda}{3!}\int d g_4d g_5d g_6T(g_3,g_5,g_4)T(g_4,g_2,g_6)T(g_6,g_5,g_1). \label{eq:boulatoveqTgroup}
\end{equation}

This equation is a non-linear Fredholm integral of the second kind. The study of this type of equations and of the properties of their solutions is an active research topic of its own (see e.g.~\cite{Brunner_book}). In particular, very few explicit solutions to this equation are known. To bring this equation in a more tractable form, we rely on the \emph{spin representation} of the Boulatov model. The spin representation relies on the Peter-Weyl theorem.

\begin{theorem}[Peter-Weyl Theorem~\cite{Wey_thm}]
\label{thm:PW}
    For a compact Lie group $G$, the matrix coefficients of the irreducible unitary representations form a basis of square-integrable function of $G$.
\end{theorem}

In the case of the group $SU(2)$, its irreducible representations are labeled by a non-negative half-integer $j$ called the \emph{spin}. The representation of spin $j$ has dimension $d_j = 2j+1$. The corresponding matrix elements are given by the \emph{Wigner matrices}~\cite{Wig31} $\left(D^j_{mn}(g)\right)_{-j \leq m,n \leq j}$. The indices of the Wigner matrices are called \emph{magnetic indices}. The Appendix~\ref{app:su2recoupling} contains some useful formulae on Wigner matrices and the $3j$ and $6j$ symbols of the group $SU(2)$.

\medskip

As a function of $SU(2)^{\otimes3}$, the field $T$ can be expanded in terms of Wigner matrices $D^{j_i}_{m_in_i}(g_i)$ via the Peter-Weyl theorem. Using the invariance~\eqref{right_inv} and the properties of the Wigner matrices (see Appendix~\ref{app:su2recoupling}), this decomposition further simplifies and takes the form 
\begin{equation}
\label{eq:Tdecompos}
      T(g_1,g_2,g_3)=\sum_{\{j,m,n\}}T^{m_1m_2m_3}_{j_1j_2j_3}\prod_{i=1}^{3}\sqrt{d_{j_i}}D^{j_i}_{m_in_i}(g_i)\TJ{j_1}{j_2}{j_3}{n_1}{n_2}{n_3}, 
\end{equation}
with $\displaystyle \TJ{j_1}{j_2}{j_3}{n_1}{n_2}{n_3}$ the $3j$ symbol of $SU(2)$ and where ${\{j\}}$ denotes the summation over the triples $j_1,~j_2$ and $j_3$ (and similarly for for $\{m\}$ and $\{n\}$). The coefficients $T^{m_1m_2m_3}_{j_1j_2j_3}$ of this expansion can be extracted using the orthogonality of Wigner matrices for the Haar measure as
\begin{equation}
      T^{m_1m_2m_3}_{j_1j_2j_3}=\int d g_1 d g_2 d g_3 \sum_{\{n\}}T(g_1,g_2,g_3)\prod_{i=1}^3\sqrt{d_{j_i}}\bar{D}^{j_i}_{m_in_i}(g_i)\TJ{j_1}{j_2}{j_3}{n_1}{n_2}{n_3}, \label{eq:Tdecomposcoeff}
\end{equation}
In the spin representation, the dependency on group elements is entirely carried by the Wigner matrices. Therefore the integral over the group elements can be performed explicitly and we are left with a series of polynomials in the coefficients $T^{m_1m_2m_3}_{j_1j_2j_3}$. The Boulatov action~\eqref{eq:actionTgroup} then rewrites as
\begin{equation}
\label{eq:actionTspin}
      S_B[T] = \sum_{j_1,j_2,j_3} \frac{\mu^2}{2}|T^{m_1,m_2,m_3}_{j_1,j_2,j_3}|^2 - \frac{\lambda}{4!} \sum_{j_1,..,j_6} \SJ{j_1}{j_2}{j_3}{j_4}{j_5}{j_6}  T^{4_{6j}}, 
\end{equation}
where the kinetic term is
\begin{equation}
|T^{m_1,m_2,m_3}_{j_1,j_2,j_3}|^2=\sum_{\substack{j_1,j_2,j_3\\m_1,m_2,m_3}}(-1)^{\sum_{i=1}^3(j_i-m_i)}T^{m_1,m_2,m_3}_{j_1,j_2,j_3}T^{-m_1,-m_2,-m_3}_{j_1,j_2,j_3},
\end{equation}
and the term $T^{4_{6j}}$ encodes the contraction of the magnetic indices $m_i$ of the four copies of the field as in the $6j$ symbol, that is
\begin{equation}
\label{eq:T6jdef}
      T^{4_{6j}}=\sum_{\{j,m\}}(-1)^{\sum_{i=1}^{6}(j_i-m_i)}T^{-m_1,-m_2,-m_3}_{j_1j_2j_3}T^{m_3,m_5,-m_4}_{j_3j_5j_4}T^{m_4,m_2,-m_6}_{j_4j_2j_6}T^{m_1,-m_5,m_1}_{j_6j_5j_1}. 
\end{equation}

In this form, the action of the Boulatov model looks like the action of a tensor field theory. The two main differences are the summation over spin labels $j$, thus the size of the tensor is not fixed, and that the coupling constant now explicitly depends on the spin of the tensor entries it couples through the $6j$-symbol of $SU(2)$. The Boulatov model admits a similar perturbative expansion as a sum over PL-manifold of dimension $3$ where the simplicial complexes carry additional group data~\cite{Boul_OG}. In the spin representation, the equation of motion~\eqref{eq:boulatoveqTgroup} now becomes
\begin{equation}
\label{eq:boulatoveqTspin}
      \mu^2T^{m_1,m_2,m_3}_{j_1,j_2,j_3} = \frac{\lambda}{3!} \sum_{j_4,j_5,j_6} \SJ{j_1}{j_2}{j_3}{j_4}{j_5}{j_6} T^{4_{6j}}_{\backslash \{m_1,m_2,m_3\}}, 
\end{equation}
where
\begin{equation}
      T^{4_{6j}}_{\backslash \{m_1,m_2,m_3\}}=\sum_{m_4,m_5,m_6}(-1)^{\sum_{i=4}^6(j_i-m_i)}T^{m_3,m_5,-m_4}_{j_3j_5j_4}T^{m_4,m_2,-m_6}_{j_4j_2j_6}T^{m_6,-m_5,m_1}_{j_6j_5j_1}.
\end{equation}

%The perturbative expansion of the Boulatov model
%has a clear interpretation as a discrete $3D$ pseudo-manifold
%\footnote{For the definition of pseudo-manifold and its differences with manifold, see for example \cite{David:1992jw,Krajewski:2012aw}} from this point of view. The interaction term encodes a tetrahedron whose vertices are copies of $T$ and whose edges are labeled by representations of $SU(2)$. The kinetic term then encodes how these tetrahedra are glued together. Thus, the Boulatov model can be viewed as a toy model of three dimensional lattice gravity.

\subsection{Classical homogeneous solutions to the Boulatov model}
\label{ssec:class_sol_Boul}

A one-parameter family of classical solutions to the Boulatov model  was proposed in~\cite{FaLi07}. These solutions are parametrized by normalised class functions $f:SU(2)\to \mathbb{C}$. The associated field $T_f$ is given by

\begin{equation}
     T_f(g_1,g_2,g_3)=\mu\sqrt{\frac{3!}{\lambda}}\int d h\delta(g_1h)f(g_2h)\delta(g_3h), \label{eq:Tfsolgroup}
\end{equation}
where $\delta(g)$ is the Dirac delta function over the group $SU(2)$ satisfying
\begin{equation}
      \int d h \delta(h)=1, \qquad \int d h\delta(h)f(h)=f(I),
\end{equation}
with $I$ the identity of $SU(2)$ group. 

\medskip

Applying the Peter-Weyl Theorem to the solution~\eqref{eq:Tfsolgroup} we obtain
\begin{equation}
  (T_f)_{j_1,j_2,j_3}^{m_1,m_2,m_3} = \mu\sqrt{\frac{3!}{\lambda}} \sqrt{d_{j_1}d_{j_3}} \sum_{l_2} f^{j_2}_{m_2,l_2} \TJ{j_1}{j_2}{j_3}{m_1}{l_2}{m_3}, \label{eq:Tfsolspin}
\end{equation}
where $f^{j}_{mn}$ is the coefficients of the Peter-Weyl decomposition of $f(g)$
\begin{equation}
     f^{j}_{mn}=\sqrt{d_j}\int d g f(g)\bar{D}^{j}_{mn}(g), \label{eq:fdecomposcoeff}
\end{equation}
and the normalization condition of $f$ reads $\sum_{j,m,n} (-1)^{m-n}f^j_{mn} f^j_{-m,-n}=1$. 

\medskip

There are two things to notice about this type of solution. Firstly, the solution~\eqref{eq:Tfsolgroup} is asymmetrical in the group elements $g_i$ since $g_2$ plays a preferential role through $f$. Secondly, the presence of Dirac delta function in the expression~\eqref{eq:Tfsolgroup} leads to several divergences. For example, the action~\eqref{eq:actionTspin} is divergent when evaluated on this class of solution. This can also be seen from the Peter-Weyl expansion of the Dirac delta
\begin{equation}
     \delta(g)= \sum_{j,m} d_j D^j_{mm}(g),
\end{equation}
as the coefficients of the Wigner matrices $D^j_{mm}(g)$ diverge as $j$ tends to infinity. Thus we need to regularize our solution. This can be achieved via different methods. For example, one possible solution could be to introduce a cut-off parameter $J$ in the Peter-Weyl expansion of $T(g_1,g_2,g_3)$~\eqref{eq:Tdecompos}. This would automatically make the action finite, and we can then send $J$ to infinity to recover our original model. For our purpose here, we will use a \emph{heat kernel regularization} to make all quantities well-defined, at the cost of only having an approximate solution to the equations of motion. To implement this regularization, we introduce a new real parameter $\varepsilon$. For any function $f$ of $SU(2)$ with coefficients $f^j_{mn}$ in its Peter-Weyl expansion, we define its heat kernel regularization as
\begin{equation}
    f_\epsilon(g) = \sum\limits_{j,m,n} \sqrt{d_j} f^j_{mn} D^j_{mn}(g) e^{-\varepsilon C_j},
\end{equation}
with $C_j= j(j+1)$ is the Casimir of the spin $j$ representation of $SU(2)$. This function is well-defined for any $\epsilon>0$ and its leading order when $\varepsilon \rightarrow 0$ is the original function $f$. Applying the heat kernel regularization to the Dirac delta function gives
\begin{equation}
    \delta_\varepsilon(g)=\sum_{j,m}d_j D^j_{mm}(g)e^{-\varepsilon C_j}.
\end{equation}

However, this function is no longer normalised. If we denote its norm as $N_\epsilon^{-2}$, the normalised function associated with $\delta_\varepsilon$ is 
\begin{equation}
    \Delta_\varepsilon(g) = \sum\limits_{j,m,n} (\Delta_\varepsilon)^j_{mn} D^j_{mn}(g) e^{-\varepsilon C_j},
\end{equation}
where the Peter-Weyl coefficients $(\Delta_\varepsilon)^j_{mn}$ has the form
\begin{equation}
     (\Delta_\varepsilon)^j_{mn}=N_\varepsilon d_j \delta_{mn} e^{-\varepsilon C_j}. \label{eq:Deltavarepsiloncoef}
\end{equation}

We can build now use the function $\Delta_\varepsilon(g)$ to construct a regularized and symmetric field $T_\varepsilon$ defined as 
\begin{align}
\label{eq:Tepsilonsolgroupreg}
     T_\varepsilon(g_1,g_2,g_3)&=\mu\sqrt{\frac{3!}{\lambda}}\int d h\delta_\varepsilon(g_1h)\Delta_\varepsilon(g_2h)\delta_\varepsilon(g_3h) \\
     &=\mu N_\varepsilon\sqrt{\frac{3!}{\lambda}}\int d h\delta_\varepsilon(g_1h)\delta_\varepsilon(g_2h)\delta_\varepsilon(g_3h). \nonumber
\end{align}

Due to the introduction of the heat-kernel regularization, $T_\varepsilon(g_1,g_2,g_3)$ is only an approximate solution of the equation of motion~\eqref{eq:boulatoveqTspin}. That is, it is a solution at first order in $\varepsilon$. Its Peter-Weyl coefficients are
\begin{equation}
    \left(T_\varepsilon\right)^{m_1m_2m_3}_{j_1j_2j_3}= \mu N_{\varepsilon}\sqrt{\frac{3!}{\lambda}}\prod_{i=1}^3\sqrt{d_{j_i}}e^{-\varepsilon C_{j_i}} \TJ{j_1}{j_2}{j_3}{m_1}{m_2}{m_3}. \label{eq:PW_eps}
\end{equation}
In particular, when $\varepsilon \rightarrow 0$ the coefficients of $T$ are given by the $3j$ symbol. This shows that the $3j$ symbol of $SU(2)$ is a (regularized) classical solution to the Boulatov model.

\medskip

Since in most situations the introduction of the heat kernel regularization doesn't lead to meaningful changes in computations but makes all expressions heavier, we will work with the non-regularized solution~\eqref{eq:Tfsolgroup} whenever possible.

\subsection{Matter degrees of freedom and spacetime coordinates}

An essential property of general relativity is the independence of the theory on coordinates used to describe the structure of spacetime. This suggests that in a non-perturbative theory of quantum gravity, the coordinates are non-physical degrees of freedom. Instead, local matter fields play the roles of rods and clocks and are used to form local reference frames to describe the dynamical system \emph{relationally}. A solution of the equation of motion of a theory of quantum geometry describes a classical spacetime geometry. These matter fields can then be used to localize events in that classical background, similarly to the spacetime coordinates we are familiar with. While different choice of matter fields can be used to fill that role, the simplest framework uses free massless scalar fields~\cite{Or_rel_coord1,Or_rel_coord2} for each spacetime dimension. 

\medskip

In our three dimensional case of the Boulatov model, this role is played by a vector field $\vec{\chi}=(\chi_1,\chi_2,\chi_3)$. Requiring the theory to be invariant under translations $\chi_i\to\chi_i+a_i$ allows for a term in the action~\eqref{eq:actionTgroup} defined as $\displaystyle \nabla =\left(\partial_{\chi_1},\partial_{\chi_2},\partial_{\chi_3}\right)$ and thus extends $T(g_1,g_2,g_3)$ to $T(g_1,g_2,g_3;\vec{\chi}):SU(2)^3\times\mathbb{R}^3\to\mathbb{C}$. This action new action in the spin representation writes
\begin{align}
\label{eq:actionTchispin}
        S_B[T(\vec{\chi})] &= \sum_{j_1,j_2,j_3}\int d^3\vec{\chi}\left[\frac{1}{2}\left|\nabla T^{m_1,m_2,m_3}_{j_1,j_2,j_3}(\vec{\chi})\right|^2+ \frac{\mu^2}{2}\left|T^{m_1,m_2,m_3}_{j_1,j_2,j_3}(\vec{\chi})\right|^2\right. \nonumber \\
    & \left.- \frac{\lambda}{4!} \sum_{j_1,..,j_6} \SJ{j_1}{j_2}{j_3}{j_4}{j_5}{j_6}  \int d^3\vec{\chi}T(\vec{\chi})^{4_{6j}}\right]. 
\end{align}
and the equation of motion becomes
\begin{align}
\label{eq:boulatoveqTchigroup} 
      &\nabla^2 T(g_3,g_2,g_1;\vec{\chi})+\mu^2T(g_3,g_2,g_1;\vec{\chi})\nonumber\\
  &=\frac{\lambda}{3!}\int d g_4d g_5d g_6T(g_3,g_5,g_4;\vec{\chi})T(g_4,g_2,g_6;\vec{\chi})T(g_6,g_5,g_1;\vec{\chi}).
\end{align}

At the level of the Boulatov model, the introduction of these scalar fields is harmless. For our purpose here, they are required to obtain a kinetic term for the perturbation around the classical solution in the effective action we will consider.

%The classical solution of the Boulatov model should be thought of as generating a $3$-dimensional classical background with coordinate $\chi$  in which the perturbation lives.

\section{Amit-Roginsky-like model as perturbations}
\label{sec:emergenceAR}

We now turn to the study of specific perturbations around the classical solution of the equation of motion constructed in the previous section. More precisely, we show that by further restraining the type of perturbations considered, it is possible to get the following effective action for the perturbation: 
\begin{align}
\label{eq:actionAR}
  S_{AR}[\phi]=&\int d^dx\left\{\frac{1}{2}\sum_m(-1)^{j-m}\left[(\nabla \phi^j_m)(\nabla \phi^j_{-m})+\mu\phi^j_m\phi^j_{-m}\right]\right. \nonumber \\
  &+\left.\sum_{m_1,m_2,m_3}\frac{\lambda}{3!}\sqrt{d_j}\TJ{j}{j}{j}{m_1}{m_2}{m_3}\phi^j_{-m_1}\phi^j_{-m_2}\phi^j_{-m_3}\right\}. 
\end{align}
This action is the action of the Amit-Roginsky (AR) model. The AR model is a cubic field theory of a vector field $\phi$ self-coupled through the $3j$ symbol for a fixed value of the spin $j$. It was originally introduced in~\cite{AR_OG} as an example of solvable field theory in the large $N=2j+1$ limit, but it is only recently in~\cite{AR_BeDe} that it was pointed out these graphs are the melonic graphs also appearing in random tensor models. This makes the Amit-Roginski model the second example of vector field theory which admits a melonic limit along with the SYK model.

\subsection{Perturbations around classical solutions}
\label{ssec:pert_sol}

Following~\cite{FaLi07}, we consider a two-dimensional perturbation over the Lie group $G$, which depends on the matter reference frame $\vec{\chi}$. The field becomes
\begin{equation}
    T_\psi(g_1,g_2,g_3;\vec{\chi})=T_f(g_1,g_2,g_3)+\xi \psi(g_1,g_3;\vec{\chi}), \label{eq:2dpert_solution}
\end{equation}
where $T_f(g_1,g_2,g_3)$ is a solution to the equations of motion of the form~\eqref{eq:Tfsolgroup} and $\psi(g_1,g_3;\vec{\chi})$ is a $2$D-perturbation with $\xi$ a real parameter $0<\xi\ll 1$. The Peter-Weyl coefficients of the perturbation are given by
\begin{align}
      \psi^{m_1m_2m_3}_{j_1j_2j_3}(\vec{\chi})&=\sum_{\{n\}}\int[d g]^3\psi(g_1,g_3;\vec{\chi})\prod_{i=1}^3\sqrt{d_{j_i}}\bar{D}^{j_i}_{m_in_i}\TJ{j_1}{j_2}{j_3}{n_1}{n_2}{n_3} \nonumber \\
  &= \delta^{j_2,0}\delta_{m_2,0}\delta^{j_1,j_3}\sqrt{d_{j_1}}\psi^{j_1}_{m_1,m_3}(\vec{\chi}).
\end{align}
 The scaling factor $\sqrt{d_{j_1}}$ is introduced anticipating on upcoming computations. The Peter-Weyl coefficients of the perturbed solution write
\begin{equation}
\label{eq:Tpsispin}
      (T_\psi)^{m_1m_2m_3}_{j_1j_2j_3}(\vec{\chi})=T^{m_1m_2m_3}_{j_1j_2j_3}+\xi \delta^{j_2,0}\delta_{m_2,0}\delta^{j_1,j_3}\psi^{j_1}_{m_1,m_3}(\vec{\chi}). 
\end{equation}

Substituting~\eqref{eq:Tpsispin} into the action~\eqref{eq:actionTchispin}, we get the action for the perturbed solution
\begin{equation}
      S_B[T_\psi(\vec{\chi})]=S_B[T] + \xi^2\cdot S_{\mathrm{eff}}[\psi]+\mathcal{O}(\xi^4),
\end{equation}
where the first order in $\xi$ vanishes since $T_\psi$ is an approximate solution to the equation of motion. The action $S_{\mathrm{eff}}[\psi]$ represents the effective action of the perturbation field $\psi^j_{mn}$ and contains corrections up to $\xi$. Therefore, $\xi^2 S_{\mathrm{eff}}[\psi]$ contains corrections up to order $\xi^3$.

\subsection{Conditions for the emergence of an Amit-Roginsky-like model} 

In this subsection, we develop each term arising from the Boulatov model in the effective action and give sufficient conditions on the coefficients of the classical solution $(T_f)^{m_1m_2m_3}_{j_1j_2j_3}$ for the effective action $S_{\mathrm{eff}}[\psi]$ to take the form of an AR-like action. The AR model involves a vector field transforming in an irreducible representation of $SU(2)$ and thus carrying only one magnetic index $m$. Hence, we already narrow the type of perturbations considered to those of the form
\begin{equation}
\label{eq:pert_spec}
\psi^{j_1}_{m_1m_3}(\vec{\chi})=\sum_{m}\sqrt{d_{j_1}}\phi^{j_1}_{m}(\vec{\chi})\TJ{j_1}{j_1}{j_1}{m_1}{m}{m_3}.
\end{equation}
and check that it can lead to an AR-like action.

\subsubsection{Quadratic terms}

Substituting perturbation~\eqref{eq:Tpsispin} into the Boulatov action~\eqref{eq:actionTchispin}, the quadratic term in $\xi$ receives three kinds of contributions.
\begin{itemize}
    \item One from the kinetic term $\sum\limits_{\{j\}}\left|(T_\psi)^{m_1,m_2,m_3}_{j_1,j_2,j_3}(\vec{\chi})\right|^2$.
    \item Two from the interaction term $T_\psi^{4_{6j}}$, depending on how the indices of the perturbation $\psi$ are contracted. Schematically, these two terms can be represented as $TT\psi\psi$ when the two perturbation fields share two magnetic indices and $T\psi T\psi$ when they only share one. There are four configurations of the form $TT\psi\psi$ and two of the form $T\psi T\psi$. The two terms give different contributions to the effective action.
\end{itemize}

\paragraph{Kinetic term\\}
\vspace{15pt}
The quadratic contribution from the kinetic term $\sum_{j_1,j_2,j_3}\left|(T_\psi)^{m_1,m_2,m_3}_{j_1,j_2,j_3}(\vec{\chi})\right|^2$ reads

\begin{align}
&\sum_{\substack{j_1,m_1,m_3\\m,m'}}(-1)^{2j_1-m_1-m_3}\phi^{j_1}_{m}\phi^{j_1}_{m'}d_{j_1}\TJ{j_1}{j_1}{j_1}{m_1}{m}{m_3}\TJ{j_1}{j_1}{j_1}{-m_1}{m'}{-m_1} \nonumber\\
  &=\sum_{j_1,m_1}(-1)^{j_1-m_1}\phi^{j_1}_{m_1}\phi^{j_1}_{-m_1}. 
\end{align}
This term is simply the quadratic term of the AR action~\eqref{eq:actionAR}. This contribution is independent of the classical solution considered and thus gives no constraints on it.

\paragraph{Terms $TT\psi\psi$\\}
\vspace{15pt}
There are four terms of type $TT\psi\psi$. Each of them contributes to the effective action as
{\small
\begin{align}
\sum_{j_1,m_1,m_6,m_4}\left[\sum_{m_3}(-1)^{-m_3-m_4}T^{-m_10-m_3}_{j_1 0 j_1}T^{m_30-m_4}_{j_1 0 j_1}\right](-1)^{2j_1-m_1-m_6}\psi^{j_1}_{m_4,-m_6}\psi^{j_1}_{m_1,m_6} 
\end{align}}%

Therefore if we require the classical solution $T^{m_1m_2m_3}_{j_1j_2j_3}$ to take the form
\begin{equation}
    \sum_{m_3}(-1)^{-m_3-m_4}T^{-m_1 0 -m_3}_{j_1,0,j_1}T^{m_3 0 -m_4}_{j_1,0,j_1}=c^{(1)}_{1} \delta_{m_1,-m_4} \label{eq:TTpsipsicon}
\end{equation}
for some coefficients $c^{(1)}_{1}$ then, for the type of perturbations~\eqref{eq:pert_spec}, the contribution of this type of term becomes
\begin{align}
\sum_{j_1,m_1}c^{(1)}_{1}(-1)^{j_1-m_1}\phi^{j_1}_{m_1}\phi^{j_1}_{-m_1}
\end{align}
which is the kinetic term of the AR model up to a scalar factor.

\paragraph{Term $T\psi T\psi$\\}
\vspace{15pt}

The remaining two quadratic contributions from the interaction term of the Boulatov model are of the form $T\psi T\psi$. Each of these terms contributes to the effective action as
\begin{align}
&\sum_{\substack{m_1,m_3,m_4,m_6\\j_1,j_3}}(-1)^{j_1+j_3-m_4-m_6}\psi^{j_3}_{m_3,-m_4}\psi^{j_1}_{m_6,m_1} \frac{1}{\sqrt{d_{j_1}d_{j_3}}} \nonumber\\
  &\times \sum_{j_2,m_2}(-1)^{\sum_{i=1}^3(2j_i-m_i)}T^{-m_1 -m_2 -m_3}_{j_1 j_2 j_3}T^{m_4 m_2 -m_6}_{j_3 j_2 j_1}. 
\end{align}
For a generic solution of the equation of motion, it leads to a non-diagonal quadratic term in $\psi$. If we impose that the classical solution satisfies the condition 
\begin{equation}
\label{eq:TpsiTpsicon}
\sum_{j_2,m_2}(-1)^{\sum\limits_{i=1}^3(2j_i-m_i)}T^{-m_1 -m_2 -m_3}_{j_1 j_2 j_3}T^{m_4 m_2 -m_6}_{j_3 j_2 j_1}=c^{(2)}_{j_1}c^{(2)}_{j_3}\delta_{m_1,-m_6}\delta_{m_3,m_4},
\end{equation}
then this contribution becomes diagonal and takes the form
\begin{align}
\left[\sum_{j_1,m_1}(-1)^{j_1-m_1}\frac{c^{(2)}_{j_1}}{\sqrt{d_{j_1}}}\psi^{j_1}_{m_1,-m_1}\right]^2.
\end{align}

Therefore, the quadratic term obtained from $T\psi T\psi$ term can also be made diagonal under a right choice of classical solution. In particular, when specializing the perturbation to those described by Equation~\eqref{eq:pert_spec} we get
\begin{align}
      \sum_{j_1,m_1}(-1)^{j_1-m_1}\frac{c^{(2)}_{j_1}}{\sqrt{d_{j_1}}}\psi^{j_1}_{m_1,-m_1}=&\sum_{j_1,m_1,m}(-1)^{j_1-m_1}c^{(2)}_{j_1}\phi^{j_1}_{m}\TJ{j_1}{j_1}{j_1}{m_1}{m}{-m_1} \nonumber\\
  =&c^{(2)}_{0}\phi^0_0.
\end{align}

\subsubsection{Cubic terms} 

There is only one type of cubic contribution which comes from the interaction term of the Boulatov model. These terms take the schematic form $T\psi\psi\psi$. There are four terms of this form and each of them contributes as
{\small
\begin{align}
 \sum_{\substack{m_1,m_3,m_4,m_6\\j_1}}(-1)^{\sum\limits_{i\neq 2,5}-m_i}T^{-m_1 0 -m_3}_{j_1 0 j_1}\frac{(-1)^{2j_1}}{d_{j_1}}\psi^{j_1}_{m_3 -m_4}\psi^{j_1}_{m_4 -m_6}\psi^{j_1}_{m_6 m_1}.
\end{align}}%
This term should take the form of the Amit-Roginsky interaction term. If we require that the classical solution satisfies
\begin{equation}
      T^{-m_1 0 -m_3}_{j_1 0 j_1} = c^{(3)}_{j_1}(-1)^{-m_3}\delta_{m_1,-m_3}, \label{eq:Tpsipsipsicon}
\end{equation}
for some coefficient $c^{(3)}_{j_1}$, this contribution becomes
\begin{align}
     &\sum_{\substack{m_1,m_3,m_4,m_6\\j_1}}(-1)^{-\sum_{i\neq 2,5}m_i}T^{-m_1 0 -m_3}_{j_1 0 j_1}\frac{(-1)^{2j_1}}{d_{j_1}}\psi^{j_1}_{m_3,-m_4}\psi^{j_1}_{m_4,-m_6}\psi^{j_1}_{m_6,m_1} \nonumber\\
 &=\sum_{\substack{m_3,m_4,m_6\\j_1}}(-1)^{2j_1-m_3-m_4-m_6}\frac{c^{(3)}_{j_1}}{d_{j_1}}\psi^{j_1}_{m_3,-m_4}\psi^{j_1}_{m_4,-m_6}\psi^{j_1}_{m_6,-m_3}.
\end{align}

And when specializing to perturbations of the form~\eqref{eq:pert_spec}, it becomes
\begin{align}
    &\sum_{\substack{m_1,m_3,m_4,m_6\\j_1}}(-1)^{-\sum_{i\neq 2,5}m_i}T^{-m_1,0,-m_3}_{j_1,0,j_1}\frac{(-1)^{2j_1}}{d_{j_1}}\psi^{j_1}_{m_3,-m_4}\psi^{j_1}_{m_4,-m_6}\psi^{j_1}_{m_6,m_1} \\
 &=\sum_{\substack{m,m',m''\\j_1}}\frac{c^{(3)}_{j_1}}{d_{j_1}} \SJ{j_1}{j_1}{j_1}{j_1}{j_1}{j_1}\phi^{j_1}_{m}\phi^{j_1}_{m'}\phi^{j_1}_{m''}\TJ{j_1}{j_1}{j_1}{m}{m'}{m''}.
\end{align}
Thus imposing a classical solution to satisfy the condition~\eqref{eq:Tpsipsipsicon}, the contribution of the cubic terms takes the form of the interaction of the Amit-Roginski model. Furthermore, as mentioned above, comparing the two conditions~\eqref{eq:TTpsipsicon} and~\eqref{eq:Tpsipsipsicon}, we can see that the two coefficients are related by the relation
\begin{equation}
      c^{(1)}_{1}={c^{(3)}_{j_1}}^2. \label{eq:c1c3rel}
\end{equation}
It follows that condition~\eqref{eq:Tpsipsipsicon} is stronger than~\eqref{eq:TTpsipsicon} as the latter is automatically satisfied when the former is.

\subsection{Emergence of the Amit-Roginsky-like model}

Now we are ready to recover the effective action for the perturbation from the Boulatov action~\eqref{eq:actionTchispin}, while requiring that the classical solution satisfies the two conditions~\eqref{eq:TpsiTpsicon} and~\eqref{eq:Tpsipsipsicon}. We first check that our classical solution~\eqref{eq:Tfsolspin} is compatible with these two conditions.

\vspace{10pt}
\paragraph{Computing coefficients $c^{(i)}$ and checking compatibility conditions.\\}

%We compute explicitly here the coefficients $c^{(1)}_{j}$,$c^{(2)}_{j}$ and $c^{(3)}_{j}$ for the homogeneous solution~\eqref{eq:Tfsolgroup} to check that these conditions are compatible with our homogeneous solution.

We start by computing $c^{(3)}_{j}$. The coefficient $T^{-m_1 0 -m_3}_{j_1 0 j_1}$ for the solution~\eqref{eq:Tfsolspin} reads
\begin{equation}
      \mu\sqrt{\frac{3!}{\lambda}}d_{j_1}f^0_{00}\TJ{j_1}{j_1}{j_1}{-m_1}{0}{-m_3}=\mu \sqrt{\frac{3!d_{j_1}}{\lambda}}f^0_{00}(-1)^{j_1+m_3}\delta_{m_1,-m_3},
\end{equation}
Therefore we can identify $c^{(3)}_j$ as
\begin{equation}
  c^{(3)}_{j}= \begin{cases} (-1)^j \mu \sqrt{\frac{3!d_{j}}{\lambda}}f^0_{00} \hspace{5pt}&\text{if} \hspace{5pt} j\in\mathbb{N} \\ 0 \hspace{5pt}&\text{otherwise}\end{cases}.
\end{equation}
And it follows that the condition~\eqref{eq:Tpsipsipsicon} is satisfied for the classical solution~\eqref{eq:Tfsolspin}. As stated before, this automatically ensures that condition~\eqref{eq:TpsiTpsicon} is satisfied due to relation~\eqref{eq:c1c3rel} between the two coefficients.

\medskip

We now turn to check condition~\eqref{eq:TpsiTpsicon}. Computing the left-hand-side for solution~\eqref{eq:Tfsolspin} yields
{\small \begin{align}
     &\frac{3!\mu^2}{\lambda} \sum_{j_2,m_2}(-1)^{\sum_{i=1}^3(4j_i-m_i)}d_{j_1}d_{j_3}\sum_{n_2,l_2}f^{j_2}_{-m_2,-n_2}f^{j_2}_{m_2,l_2}\TJ{j_1}{j_3}{j_2}{m_1}{m_3}{m_2} \TJ{j_1}{j_3}{j_2}{-m_6}{m_4}{l_2}. 
\end{align}}%
Hence to satisfy~\eqref{eq:TpsiTpsicon}, the Peter-Weyl coefficients $f^{j_2}_{m_2n_2}$ must satisfy
\begin{equation}
  \sum_{m_2}(-1)^{n_2-m_2}f^{j_2}_{-m_2,-n_2}f^{j_2}_{m_2,l_2}=d_{j_2}c_{f,j_2}^2\delta_{n_2,l_2},
\end{equation}
for some coefficients $c_{f,j_2}$. Together with the normalization condition satisfied by $f$, we get the condition that these new constants should satisfy
\begin{align}
	1&= \sum_{j_2,m_2,n_2,l_2}(-1)^{n_2-m_2}f^{j_2}_{-m_2,-n_2}f^{j_2}_{m_2,l_2} \delta_{n_2,l_2} \nonumber\\
	&=\sum_{j_2}d^2_{j_2} c_{f,j_2}^2.
\end{align}
Leaving aside the consideration due to the introduction of the heat-kernel regularization for now, we can see that this condition is satisfied by~\eqref{eq:Tfsolspin} when specialized to $f=\delta$ as its Peter-Weyl coefficients are $\delta^j_{mn} = \delta_{mn}$. 

%The heat kernel regularization doesn't spoil this feature as we shall see by computing the coefficient $c_{f,j_2}$ of the regularized solution hereafter.
%and we can get the explicit form \eqref{eq:cfepsilon} of $c_{f,j_2}$ by substituting the heat kernel regularized solution \eqref{eq:PW_eps}.

\vspace{10pt}
\paragraph{The effective action for the perturbation $\psi$\\} 

For a perturbation of the form~\eqref{eq:pert_spec}, and a classical solution satisfying both conditions~\eqref{eq:TpsiTpsicon} and~\eqref{eq:Tpsipsipsicon}, the effective action for the vector perturbation $\phi^{j}_{m}$ becomes

\begin{equation}
\label{eq:action_phi}
      S[\phi^{j}_{m}]= S_0[\phi^{0}_{0}]+\sum_{j>0}S_j[\phi^{j}_{m}],
\end{equation}
where 
\begin{equation}
      S_0[\phi^0_0]= \frac{1}{2}\int d^3\vec{\chi}\left\{(\nabla \phi^{0}_{0})^2+\left[\mu^2+\frac{\lambda}{3!}(2{c^{(3)}_{0}}^2+c_{2,0}^2)\right](\phi^{0}_{0})^2\right\}-\frac{\xi\lambda}{3!} c^{(3)}_{0}\left(\phi^0_0\right)^3,
\end{equation}
and
\begin{align}
\label{eq:actionphij}
      S_j[\phi^{j}_{m}] &=\frac{1}{2} \int d^3\vec{\chi}\left[|\nabla\phi^{j}_{n}|^2+\left(\mu^2+\frac{\lambda}{3!}{c^{(3)}_{j}}^2\right)|\phi^{j}_{n}|^2\right] \\
  &-\frac{c^{(3)}_{j_1}}{2d_{j}}\frac{\xi\lambda}{3!}\SJ{j}{j}{j}{j}{j}{j}\sum_{m_1,m_2,m_3}\phi^{j}_{m_1}\phi^{j}_{m_2}\phi^{j}_{m_3}\TJ{j}{j}{j}{m_1}{m_2}{m_3}, \nonumber
\end{align}
where $\sum_n|\phi^j_n|^2=\sum_n(-1)^{j-n}\phi^j_n\phi^j_{-n}$. 

\medskip

We can see that the modes of the perturbation with different spin labels $j$ decouple. For each spin $j$, the effective action takes the form of an Amit-Roginski action where the mass and coupling constant depend on the spin $j$.

\paragraph{For the heat kernel regularized solution\\}

The checks performed above on the classical solution~\eqref{eq:PW_eps} still hold at first order in $\epsilon$ when considering the heat kernel regularized solution~\eqref{eq:Tepsilonsolgroupreg}. Using its Peter-Weyl coefficients \eqref{eq:Deltavarepsiloncoef}, we see that the coefficient $c^{(3)}_{j}$ is simply
\begin{equation}
      c^{(3)}_{j}=(-1)^{j} \mu \sqrt{\frac{3!d_j}{\lambda}}(\Delta_\varepsilon)^0_{00}=(-1)^{j}\mu\sqrt{\frac{3!d_j}{\lambda}}N_\varepsilon.
\end{equation}
And the coefficient $c_{f,j}$ becomes
\begin{equation}
      c_{f,j}=N_\varepsilon  e^{-\varepsilon C_j}. \label{eq:cfepsilon}
\end{equation}

However, condition~\eqref{eq:TpsiTpsicon} is only satisfied approximately. Indeed at first order in $\varepsilon$ we have
\begin{equation}
    \sum_{j,m}d_je^{-2\varepsilon C_j}\TJ{j_1}{j_2}{j}{m_1}{m_2}{m}\TJ{j_1}{j_2}{j}{m'_1}{m'_2}{m} = \delta_{m_1'm_1}\delta_{m_2'm_2} + O(\varepsilon). \label{eq:3jsumdeltaappro}
\end{equation}
due to the properties of the $3j$-symbol. At first order, the coefficients $c^{(2)}_{j}$ of condition~\eqref{eq:TpsiTpsicon} are given by
\begin{equation}
 c^{(2)}_{j}=\mu d_jN_\varepsilon\sqrt{\frac{3!}{\lambda}}.
\end{equation}
Finally, the effective action we obtain for the perturbation around the heat kernel regularized solution is
\begin{align}
\label{eq:actionphi0epsilon}
       S_0[\phi^0_0]&= \int d^3\vec{\chi}\left\{\frac{1}{2}\left[(\nabla \phi^{0}_{0})^2+\mu^2\left(1+3N_\varepsilon^2\right)(\phi^{0}_{0})^2\right]-\frac{\sqrt{\lambda}\xi\mu N_\varepsilon}{\sqrt{3!}} \left(\phi^0_0\right)^3 \right\}, \\
            \label{eq:actionphijepsilon}
   S_j[\phi^{j}_{m}]  &=\int d^3\vec{\chi}\left\{\frac{1}{2}\left[|\nabla\phi^{j}_{m}|^2+\mu^2\left(1+d_jN_\varepsilon^2 \right)|\phi^{j}_{m}|^2\right]\right. \nonumber \\
  &\left.-\frac{(-1)^j}{\sqrt{3!}}\frac{\sqrt{\lambda}\xi\mu N_\varepsilon}{2\sqrt{d_{j}}}\SJ{j}{j}{j}{j}{j}{j}\sum_{m_1,m_2,m_3}\phi^{j}_{m_1}\phi^{j}_{m_2}\phi^{j}_{m_3}\TJ{j}{j}{j}{m_1}{m_2}{m_3}\right\}, 
\end{align}

Therefore the introduction of the heat kernel regularization preserves the form of the effective action, which still takes the form of an Amit-Roginski-like action. 

\medskip

This shows that the Amit-Roginski model can be obtained by considering particular perturbations around some classical solutions of the Boulatov model, provided that the classical solution satisfies some extra conditions given by Equations~\eqref{eq:TpsiTpsicon} and~\eqref{eq:Tpsipsipsicon}.

\section{Melonic dominance}\label{melonicdominance}

An important feature of the Amit-Roginski model is that its Feynman graphs are dominated by melonic graphs at large $N=2j+1$ limit as it ensures the resolvability of the theory in that limit. The dominance of the melonic graphs in the large $N$ limit holds under some conjecture on the behavior of $3nj$-symbols of $SU(2)$~\cite{Haggard_2010,CoMa15,Bonzom_2012}. The main difference between the effective action~\eqref{eq:action_phi} and the original AR action is the addition of the sum over the spin labels $j$. This additional summation could induce a different behavior of the Feynman diagrams of the theory. We study qualitatively the behavior of the Feynman amplitudes of the effective action and give different way to ensure the melonic dominance for the effective action of the perturbations.

\subsection{Amplitudes of the Feynman diagrams}

We study the Feynman diagrams associated with the action~\eqref{eq:actionphij}. As in the Amit-Roginsky model, each Feynman diagram $\gamma$ can be decomposed~\cite{AR_OG,AR_BeDe} as
\begin{equation}
    \mathcal{A}_{\gamma}=c_{\gamma}I_{\gamma}\sum_{j}\left(\frac{\lambda \{6j\}}{3!d_j}\right)^{v}a_{\gamma,j},
\end{equation}
where $c_{\gamma}$ is a combinatorial factor of the diagram, $I_{\gamma}$ is the \emph{isoscalar} which contains the usual spacetime integral, and $a_{\gamma,j}$ is the \emph{isospin} which encodes the group theoretic dependency. The isoscalar contribution $I_{\gamma}$ is the same for all spins $j$ as all the spin dependency is contained in the isospin contribution. Therefore we solely have to study the behaviour of the isospin contribution. For simplicity, we shall drop the heat kernel regularisation in this section as they do not play any role and can be restored immediately.

\medskip

Using identities of $SU(2)$ recoupling theory (see Appendix~\ref{app:su2recoupling}), the isospin part can be decomposed as a product of graphs that admit no $2$-cut. These graphs are said to be \emph{$2$-particule irreducible} ($2$-PI). On the contrary, graphs that have a $2$-cut are \emph{$2$-particule reducible} ($2$-PR). The melonic graphs are $2$-PR graphs which are \emph{fully} $2$-PR (F$2$-PR): any non-trivial melonic graph admits a $2$-cut which splits the graph into two other melonic graphs with fewer vertices. The isoscalar amplitude of a melonic graph with $v=2n$ vertices factorizes as
\begin{equation}
\label{eq:F2PR_contr}
    \mathcal{A}_{F2PR}\sim\sum_{j}d_j^{1-3n}\{6j\}^{2n}.
\end{equation}
For graphs which are not F$2$-PR, the Feynman amplitude can be factorized as a product of $k$ different $2$-particule irreducible graphs. 
\begin{equation}
    \mathcal{A}_{NF2PR}\sim \sum_{j} \{6j\}^{2n}d_j^{-n_0-2n}\prod_{i=1}^k \mathcal{A}_{\{3n_i j\}},
\end{equation}
where
\begin{equation}
    n=1+n_0-k+\sum_{i=1}^k n_i,
\end{equation}
and $\mathcal{A}_{\{3n_i j\}}$ is the amplitude of a three-particle irreducible diagrams with $2n_i$ vertices. To further evaluate this amplitude, we need to compute the $\{3nj\}$-symbols for arbitrary large $n$. Unfortunately, there is no exact expression for the evaluation of these symbols to this day and one has to resort to asymptotic results~\cite{Haggard_2010,CoMa15,Bonzom_2012} instead. The analytical and numerical analysis of the $\{3nj\}$-symbols led to the following conjecture.

\begin{conjecture}[\cite{AR_OG}]
\label{conj:3nj}
When $j$ goes to infinity, the $\{3n j\}$-symbol obeys the bound
\begin{equation}
    \label{eq:bound_3nj}
    \vert \{3n j\} \vert \leq d_j^{-n+1-\alpha},
\end{equation}
for some positive real number $\alpha>0$.
\end{conjecture}

This conjecture ensures the dominance of melonic graphs at large $j$ in the Amit-Roginsky model~\cite{AR_OG,AR_BeDe}. Here, the action~\eqref{eq:actionphijepsilon} yields an additional summation over spins $j$, hence this asymptotic result is not sufficient to ensure melonic dominance in our model. If there is a threshold value $j_t$ such that the bound~\eqref{eq:bound_3nj} holds for all $j > j_t$, then in a non-melonic graph, it would be possible for contributions of spin $j'\leq j_t$ to compensate the contribution coming from larger spins to yield a larger contribution than melonic graphs overall. Therefore, it would be possible for $\mathcal{A}_{NF2PR}$ to give a larger contribution compared to $\mathcal{A}_{F2PR}$ at large $N$ and spoil melonic dominance of the model. 

%If the bound~\eqref{eq:bound_3nj} holds all $N$, then the amplitude of graph that is not F$2$-PR is bounded order by order by the $2FPR$ contribution. However if this bound fails for $N<N_t$ for some threshold $N_t$, then the sum from $N=1$ to $N=N_t$ is still bounded order by order by the contribution of $2FPR$ graphs, but the rest of the term may not be.  Thus we see that the addition of the sum over spin labels $j$ can potentially spoil the existence of a melonic large $N$ limit.

\subsection{Ensuring melonic dominance}

From the previous discussion, we have seen that the sum over $j$ of our model can potentially change the amplitude of a Feynman graph for our model. In this section, we make different propositions to ensure that the melonic limit is preserved.

\paragraph{A stronger conjecture\\}

To ensure melonic dominance, we propose a stronger version of Conjecture~\ref{conj:3nj}, which holds for all spins $j$.
\begin{conjecture}
\label{conj:3nj_all}
For all $j \geq 0$, there is a real number $\alpha>0$ such that the $\{3n j\}$-symbol obeys the bound
\begin{equation}
    \label{eq:bound_3nj_all}
    \vert \{3n j\} \vert \leq d_j^{-n+1-\alpha},
\end{equation}
\end{conjecture}

\medskip

Under Conjecture~\eqref{conj:3nj_all}, the contribution of a non-F$2$-PR graph writes
\begin{equation}
\label{eq:bound_3nj_graph}
    \mathcal{A}_{NF2PR}\leq \sum_j d_j^{1-3n-\alpha}\{6j\}^{2n}.
\end{equation}
Therefore it is bounded term-by-term by the $2FPR$ contribution~\eqref{eq:F2PR_contr} and it follows that melonic graphs are the ones that dominate the expansion at large $N=2j+1$.

\paragraph{Recovering the original AR action\\}

Another possible way to recover a melonic large $N$ limit is to get an effective action that gives exactly the original Amit-Roginsky model, exhibiting a melonic limit under the weaker Conjecture~\ref{conj:3nj}. To do so, we can choose a perturbation of the form~\eqref{eq:Tpsispin} but enforce the selection of a single spin contribution $j$ i.e.
\begin{equation}
      (T_\psi)^{m_1m_2m_3}_{j_1j_2j_3}(\vec{\chi})=T^{m_1m_2m_3}_{j_1j_2j_3}+\delta^{j_1j}\delta^{j_2,0}\delta_{m_2,0}\psi^{j_1}_{m_1,m_3}. 
\end{equation}
The melonic limit is then ensured for $j$ large enough such that Conjecture~\ref{conj:3nj} holds.

\paragraph{Heat-kernel regularization\\}

Another way to recover melonic dominance is to use the additional parameter $\varepsilon$ of the heat-kernel regularized solution~\eqref{eq:Tepsilonsolgroupreg}. When $j_2=0$ the solution has the form
\begin{equation}
      (T_{\varepsilon})^{m_1m_2m_3}_{j_10j_3}=\mu N_{\varepsilon}\sqrt{\frac{3!}{\lambda}}e^{-2\varepsilon C_{j_1}} \sqrt{d_{j_1}}(-1)^{j_1-m_1}\delta_{j_1,j_3}\delta_{m_1,-m_3}.
\end{equation}
Now, if we fix $\displaystyle \varepsilon=\frac{1}{2N(N+1)}$, the expression above scales as $\sqrt{d_{j_1}}$ for $j_1<N$. Hence, the coefficients $(T_{\varepsilon})^{m_1m_2m_3}_{j_1j_2j_3}$ with $j_i< N_t$ for some threshold $N_t$ can be neglected, and the tensor modes are dominated by the coefficients with larger spin $j$. Then, at first order in $\epsilon$, we have
\begin{equation}
      (T_\psi)^{m_1m_2m_3}_{j_1j_2j_3}\simeq
      \left\{\begin{array}{cc}    T^{m_1m_2m_3}_{j_1j_2j_3}+\delta^{j_1j}\delta^{j_2,0}\delta_{m_2,0}\psi^{j_1}_{m_1,m_3} & N_t\leq j_i\leq N\nonumber\\
 0 &\text{otherwise}\end{array}  \right. .
\end{equation}
Such perturbation gives the following amplitude
\begin{equation}
      \mathcal{A}_{\gamma}=c_{\gamma} I_{\gamma}\sum_{j=N_t}^{j=N}\left(\frac{\lambda \{6j\}}{3!\sqrt{d_j}}\right)^{v}a_{\gamma,j},
\end{equation}
which is also dominated by melonic graphs under Conjecture~\ref{conj:3nj} for a choice of $N_t$ compatible with this conjecture.

\section{Conclusion} \label{sec:summary}

We have shown that an Amit-Roginski-like action can be obtained as a perturbation around classical solutions of the Boulatov model under mild hypothesis on the shape of the solution considered. While little is known about classical solution to the Boulatov model beyond the framework studied here, we have given sufficient conditions on the solutions to the equation of motion to recover an AR-like dynamic for the effective action of the perturbation around the classical solution. The main difference between our effective action for the perturbation and the usual AR model is the addition of the summation on the spin index $j$. While it is unclear whether this summation could spoil the dominance of melonic diagram in the most general framework, it is possible to ensure this melonic limit under the same hypothesis as in the AR-model either by specializing the type of perturbation considered or by making use of the heat kernel regularization to take a kind of double scaling limit.

\medskip

This result opens the way for at least two different generalizations. Firstly, a natural follow-up would be to find other classical solutions to the Boulatov model and to study perturbations around these solutions to see if they also admit Amit-Roginski-like perturbations or other familiar vector field theory. Secondly, as already mentioned in the Introduction, the SYK model is another type of field theory that is known to enjoy a melonic limit. Thus it would be interesting to investigate if an SYK-like action can also be obtained within a GFT setup.

%% file: Chapters/Conclusion.tex
\chapter{Conclusion}
\label{Chap:concl}

\setlength{\epigraphwidth}{0.58\textwidth}
\epigraph{\itshape Ajuste, ô fils quelconque d'Ève,\\
N'importe quel calcul à n'importe quel rêve,\\
Ajoute à l'hypothèse une lunette, et mets\\
Des chiffres l'un sur l'autre, à couvrir les sommets\\
De l'Athos, du Mont-Blanc farouche, du Vésuve,\\
Monte sur le cratère ou plonge dans la cuve,\\
Fouille, creuse, escalade, envole-toi, descends,\\
Fais faire par Gambey des verres grossissants,\\
Guette, plane avec l'aigle ou rampe avec le crabe,\\
Crois tout, doute de tout, apprends l'hébreu, l'arabe,\\
Le chinois, sois indou, grec, bouddhiste, arien,\\
Va, tu ne saisiras l'extrémité de rien.\\
Poursuivre le réel, c'est chercher l'introuvable.}{V. Hugo, \textit{A l'Homme}}

This concluding chapter provides a glimpse into potential future directions stemming from the findings presented in this thesis. The topics explored within this manuscript are part of an expansive and dynamic research domain and consequently, the propositions put forth hereafter are a mere fraction of the numerous possibilities for further investigations, reflecting my personal perspective and interests. Each section within this chapter is dedicated to a specific type of model. The first two sections focus on models studied throughout this thesis, respectively random tensor models and $b$-deformed constellations. The third section motivates the study of other multi-matrix models that seem interesting to further probe integrability features in multi-matrix models. In the years to come, I will delve into some of these questions with the hope of contributing to their resolution.

\section{Combinatorics of tensor models}
\label{sec:comb_tens_mod}

All tensor models studied in this manuscript share similar features. They admit $\frac{1}{N}$-expansion dominated by melonic graphs, they can be decomposed in schemes and their dominant schemes are decorated trees. However, the combinatorial techniques on which the implementation of these mechanisms depends on details of the model considered through its symmetry group and the choice of interactions. Hence a natural line of work is to look for results that would allow us to better understand how broadly these tools can be applied and help transfer these results to other models. 

\subsection{Extending the scheme decomposition}

\begin{figure}
    \centering
    \includegraphics[scale=0.5]{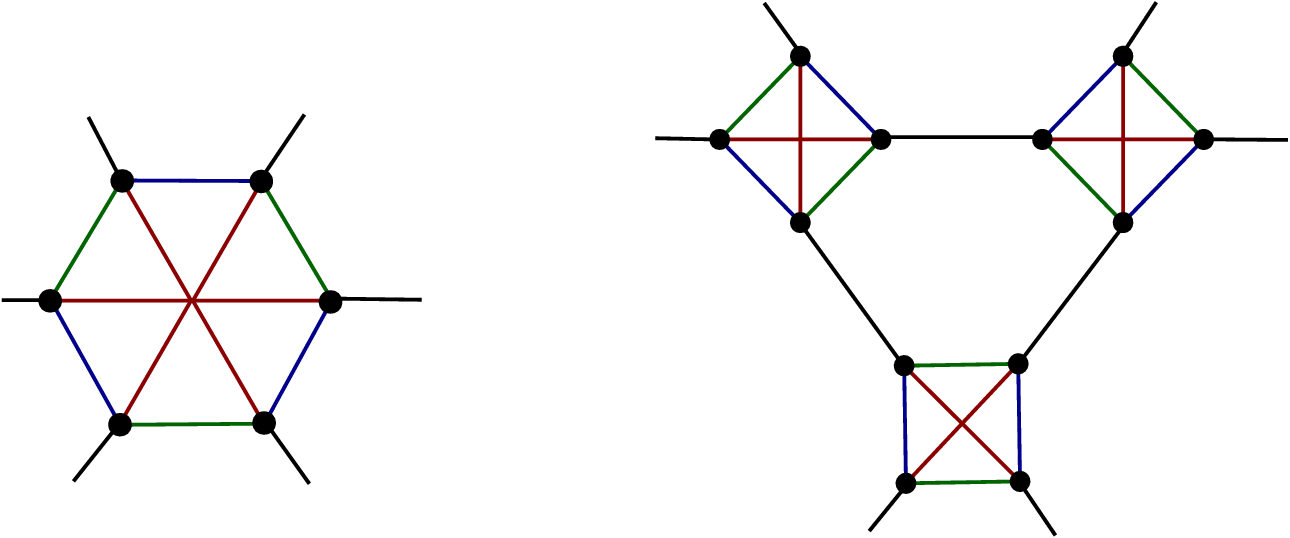}
    \caption{The wheel bubble cannot be replaced by a tetrahedral subgraph (on the right) with the correct $\frac{1}{N}$-scaling due to its internal face.}
    \label{fig:prism_to_tet}
\end{figure}

As we have seen in~Chapter~\ref{Chap:DScale} through the example of the $U(N)^2\times O(D)$ tensor model, being able to embed (through a finite mapping) schemes of a model into those of another model for which the scheme decomposition is already known allow us to bypass the derivation of the finiteness of the schemes for the first model. A similar idea can be used to include more interactions in the action within a given model. It is actually what we used for the $O(N)^3$ model with quartic interaction studied in the manuscript: the pillow interaction can be replaced by a dipole of the tetrahedral interaction. This allows to study the tensor graphs with both interactions from the simpler model with tetrahedral interaction only. This method was also applied recently in~\cite{KMT_DS} for a sextic interaction called the \emph{prismatic interaction} which can also be written using the tetrahedral bubble (and some additional decorations). It could be generalized to many other bubbles of the $O(N)^3$ model. Consider a $2n$-point subgraph $\hat{G}$ of the $O(N)^3$ model with tetrahedral interaction. By contracting all of the edges of color $0$ of $\hat{G}$, we obtain a $3$-colored graph with $2n$-edges of color $0$, that is a bubble of the $O(N)^3$ model. However, the scaling in $N^{-1}$ of the initial subgraph and the final bubble might not match, as can be illustrated via the \emph{wheel interaction} of Figure~\ref{fig:prism_to_tet}.

The only $6$-point graph that leads to this bubble has one internal face, and therefore its optimal scaling~\eqref{eq:scal_opt}) is enhanced by a factor $N$ compared to the wheel bubble. It would be interesting to characterize precisely which bubbles can be obtained as subgraphs of the tetrahedral bubble.
\begin{question}
    What are the bubbles $\hat{b}$ of the $O(N)^3$ model that can be obtained from some $2n$-point graphs $\hat{G}$ with tetrahedral bubble only such that $\hat{G}$ and $\hat{b}$ have the same scaling in $N^{-1}$?
\end{question}

Answering this question would allow us to extend the set of interactions of the $O(N)^3$ model for which the scheme decomposition applies. It would allow us to study the $O(N)^3$ model with more interactions than in the present case. In all cases studied in the literature, the double scaling limit is dominated by broken chains, which gives rise to a tree structure for the dominant schemes. Including more interactions complexifies the structure of the locus of singular for the generating function of melons. By applying the scheme decomposition to these models, we could probe whether the double scaling limit is still dominated by broken chains and if it still gives a tree structure for the dominant scheme. 
\begin{question}
Are the singular point relevant for the double scaling limit always associated with (a generalization of) broken chains? If another combinatorial structure is relevant, are the dominant schemes still in bijection with some family of decorated trees?
\end{question}
In the case where other regimes exist, a natural follow-up would be to check whether there are phase transitions between these regimes. Similar investigations could be made including interactions that cannot be written as from the tetrahedral bubble, but one would have to check first that the scheme decomposition applies. It could be interesting to see whether there are differences between the double scaling of the two types of interactions or if it is simply a shortcut to apply the scheme decomposition. 

\subsection{Tensor models with mixed index symmetry}

A different line of work concerns the tensor models with mixed index symmetry. They form a much broader class of tensor models than the ones with non-mixed index symmetry, and contain models relevant for both physics (e.g. multi-matrix models) and combinatorics (e.g. the $U(N)\times O(D)$ model whose Feynman graphs are loop-decorated maps). But these models are more difficult to study than their non-mixed index counterpart. Each interaction terms have to be described via one stranded graph, and the propagator contains different index contraction patterns hence at the combinatorial level it is given by several $2$-point stranded graphs (see e.g.~\cite{Ca18,CaHa21,TaFe}). This makes the number of stranded graphs associated with a single Feynman graph blow up very quickly, on top of potential cancellations that can occur between these configurations. This makes the stranded representation highly impractical, both at the combinatorial and algebraic levels. For these reasons, even basic results like the existence of the $\frac{1}{N}$-expansion, or the dominance of melonic graphs (when the $\frac{1}{N}$-expansion exists) are difficult to establish. They often require additional hypothesis on the tensor, typically tracelessness. Sometimes even stronger conditions are required to ensure that some configurations of the Feynman graphs that we do not how to handle cannot arise. For example, in the $U(N)\times O(D)$ model of~\cite{TaFe} studied in Section~\ref{sec:DS_U(N)xO(D)}, bipartiteness is necessary to ensure the existence of the melonic large $N$ limit as it rules out the existence of faces and O-loops of length $3$. Imposing tracelessness on the matrices instead would be a weaker hypothesis, but it proves much more difficult to establish the existence of the large $D$ limit in this case as pointed out in~\cite[Section $3$]{TaFe}. 
\begin{question}[{\cite[Section $3$]{TaFe}}]
    Can we establish the existence of the large $D$ expansion in the Hermitian traceless $U(N)\times O(D)$ multi-matrix model with tetrahedral interaction?
\end{question}
While we were writing our results on the double scaling limit of the $U(N)\times O(D)$ multi-matrix model, I tried to come up with a combinatorial proof of the existence of the large $D$ limit, working by exhaustion of all possible cases that needed to be bounded but I always ended up facing cases that I couldn't handle while they never seemed to appear in specific examples, which seem to suggest that these cases never occur.

\medskip

The complexity of the Feynman graphs of this class of model calls for the need for other tools beyond combinatorial methods to study them. Since random tensor models are defined from the action of a symmetry group, representation theory seems to be a natural framework to look for such tools. Indeed, in the cases studied in the litterature~\cite{Ca18,CaHa21,CarrozzaPozsgay2018}, the existence of the $\frac{1}{N}$ expansion is related to the irreducibility of the representation the tensor transforms in. For example, imposing the tensor to be traceless amounts to removing the vector representations when the group action on the tensor is reducible. Yet, there are as of now very few results on random tensor models that make use of representation theoretic tools. One example of results obtained through representation theory is given by~\cite[Lemma $9$]{BCGK_symm} where the authors use Schur's Lemma to show that only certain index structures can occur in the $2$-point function of a tensor model transforming in an irreducible representation of $O(N)$. 
\begin{question}
    Can we generalize Lemma $9$ of~\cite{BCGK_symm} to any $2n$-points graphs?
\end{question}
Fully exploring the relationship between the representation theory of Lie groups and the combinatorics of tensor models is an ambitious goal that goes beyond the scope of my capacities. Generalizing the result of~\cite{BCGK_symm} is a more accessible goal that would be a first step in that direction.

\section{On \texorpdfstring{$b$}{b}-deformed constellations models}
\label{sec:constr_MM}

In Chapter~\ref{Chap:const}, we derived the constraints for the two cubic $b$-deformed constellation models. The Lemma employed can in principle be applied to a wide variety of models and has the advantage of allowing for a derivation of the constraints valid for all values of $b$. However, this method relies on brute force computation of the modes of the evolution equation of the model. Not only it makes it inherently more difficult to apply this method to other models, but it doesn't bring any hindsight to the underlying algebraic or combinatorial structure at work behind this computation. 

\subsection{Understanding the \texorpdfstring{$b$}{b}-deformation}

The introduction of the $b$-deformation renders usual combinatorial techniques inapplicable since the $b$-weight depends on the order in which edges are deleted. Yet, there are hints that the results known for $b=0$ could extend to arbitrary values of $b$. First, the evolution equation for generic $b$ is similar to the case $b=0$ where it has a clear combinatorial interpretation, and second, in the quadratic and cubic models studied in this manuscript, the constraints of the generating function form an algebra for all $b$. This suggests that both the combinatorial and integrable structures should extend to arbitrary values of $b$. For the combinatorial side, this was already suggested by Goulden and Jackson in~\cite{GJ_b_conj} through the matching-Jack conjecture and the $b$-conjecture. But these different perspectives are not independent. Progress on one side would help shed light on the properties of the $b$-deformation from the other point of view. Therefore we make similar suggestions from the perspective of matrix models and integrable hierarchies. 

\medskip

The generating function of $m$-constellations can be obtained as the perturbative expansion of a multi-matrix model (with one matrix for each color of the constellation). For single matrix models, the $b$-deformation is known to be related to the $\beta$-ensembles~\cite{Desro09,BE09}. Therefore, I would expect the notion of $\beta$-ensembles to admit a generalization to these models.
\begin{question}
    Does the notion of $\beta$-ensembles generalize to constellations? 
\end{question}

On the side of integrable hierarchies, the BKP hierarchy has been shown to be related to non-oriented map enumeration~\cite{VdL01} and to the symmetric group version of the HCIZ integral~\cite{BCD22_ON_HCIZ}. Therefore, I would expect that the generating function of non-oriented constellations (i.e. evaluated at $b=1$) should also satisfy this hierarchy. This would be a first step in the study of the integrability properties of these models. More ambitiously, the fact that the constraints of the model studied here still form an algebra for generic values of $b$ suggests that integrable structure may survive the $b$-deformation and hold in all generality. However, there is to my knowledge no candidate hierarchy in this case.

\begin{question}
    Is the generating series of $b$-deformed constellation integrable for all values of $b$? If so, what is the corresponding hierarchy?
\end{question}

From the belief that the link between the combinatorics of maps and integrable hierarchy should all for all value of $b$, a tentative proposal can be made by translating the $b$-conjecture of Goulden and Jackson in this language. Their conjecture proposes combinatorial expressions for the $b$-weight as a function of the matchings. The matchings are the natural tool to count non-oriented maps at $b=1$, hence in some sense, their conjecture generalizes the case $b=1$ to more general weights. As the integrable hierarchy associated with non-oriented maps is known to be the BKP hierarchy, it would be interesting to look to incorporate a $b$-weight in the solutions of the BKP hierarchy, in a way that is coherent with the existing results on map enumeration at $b=0$.

\subsection{The algebra of constraints for orientable constellations}

Even at $b=0$, that is in the orientable case, constellations offer a rich framework to better understand the interplay between random multi-matrix models, the combinatorics of maps and integrable hierarchies. The generating series of constellations can be expressed as the generating series of a particular multi-matrix model. Different statistical physics models on a random lattice (e.g. Ising or Potts model) can be expressed as multi-matrix models. However since two random matrices do not commute, multi-matrix models are much more difficult to solve than their one-matrix counterparts. This makes constellation models particularly interesting among multi-matrix models since they are known to be integrable. Some specialization of their generating functions has been shown to satisfy the KP hierarchy and topological recursion. A different aspect of integrability lies in the algebra generated by the constraints, which forms a Virasoro algebra in maps. In this manuscript, we have explicitly checked that the constraints of cubic constellations still formed an algebra, but we did not connect it to any already known algebra. This suggests two directions to pursue this work: extracting constraints (and their algebras) in other constellations models and studying the properties of these algebras $\tilde{W}^{(m)}$ ($m=3$ in our case).

\begin{question}
    What is the algebra of constraints $\tilde{W}^{(m)}$ of $m$-constellations? 
\end{question}

To answer this question, one possible approach is to make use of the integrability properties of constellation models. In the integrability literature, the $W$-algebras are known to play a key role as they allow to define models satisfying topological recursion~\cite{BBCCD_W}. In particular in the case of maps, the Virasoro algebra $\tilde{W}^{(2)}$ is a $W$-algebra. There are natural candidate $W$-algebras denoted $W^{(m)}$ related to $r$-spin Hurwitz numbers~\cite{BCEGF} to extend this result to constellations models of order $m$ but the algebras $\tilde{W}^{(m)}$ and $W^{(m)}$ differ for $m \neq 2$.
\begin{question}
    What is (if any) the interplay between the algebra of constraints $\tilde{W}^{(m)}$ and the $W$-algebra $W^{(m)}$?
\end{question}
Albeit for a single random matrix in an external field, this question has been probed in~\cite{MMW_GKM_phases} where it has been shown that the two algebras can be obtained as different "phases" of a generalized Kontsevich model. This question would help us in both directions. It would connect the combinatorial models of constellations to the usual (in integrability) $W$-algebras, but it would also help to understand the interplay between the two, thus potentially opening the way for generalization of our results to $m \geq 4$ through a method that does not rely on such direct computation.

\section{Other multi-matrix models}
\label{sec:tensor_integ}

\subsection{The \texorpdfstring{$O(n)$}{O(n)}-loop model}

The $O(n)$-loop model~\cite{BBG12} is a multi-matrix model of $n+1$ Hermitian matrices where $n$ of the matrices appear at most quadratically and are solely coupled to the remaining matrix. It has been shown in~\cite{BGF16} to satisfy topological recursion (with a spectral curve of genus $1$) which makes it one of the few exactly solvable multi-matrix models. Its Feynman graphs resemble usual maps but are decorated with simple loops encoding the quadratic interaction of the $n$ matrices. These $n$ matrices can be explicitly integrated via Gaussian integration. After this operation is performed, the model becomes a multi-trace matrix model where $n$ is related to the weight of those multi-trace terms. This allows us to define the $O(n)$-loop model for any real value of $n$. Due to the extra parameter $n$ and the multi-trace interaction, the $O(n)$-loop model has several distinct regimes of expansion (e.g. dense vs dilute loops) and admits connection to other multi-matrix model for particular values of $n$. The connection between the $O(n)$-loop model and integrability was suggested to me by V. Bonzom~\footnote{and I thank him for allowing me to share this idea here.} through the following question.

\begin{question}[Original idea of V. Bonzom]
    Is there an integrable hierarchy satisfied by the $O(n)$-loop model for any value of $n$?
\end{question}

The first steps towards an answer to that question can be taken by studying the behavior of the partition function for particular values of $n$. For example, at $n=0$, it is the usual $1$-Hermitian random matrix model and hence it satisfies the KP hierarchy. Yet, at $n=1$ going to eigenvalues shows that the partition function can be represented as a Pfaffian, i.e. a solution of the BKP hierarchy instead of KP. This raises the following question.

\begin{question}[Original idea of V. Bonzom]
    Can we find a determinantal representation of the partition function of the $O(1)$-loop model?
\end{question}

A strategy to tackle this question is to check whether the partition function of the $O(1)$-loop model satisfies the first equations of the KP hierarchy. In particular, if it doesn't, then it hints at a possible interpolation between the KP hierarchy at $n=0$ and the BKP hierarchy at $n=1$. This would encourage looking for other values of $n$ where the model is integrable (with possibly other hierarchies than KP or BKP at play).

\medskip

Another interesting result on the connections between the $O(n)$-loop model and integrable hierarchies is given in~\cite{Ko95}. At $n=2$, the model is dual to a height model related to the Kac-Moody algebra $A^{(1)}_0$. This shows that it satisfies the KP hierarchy also for $n=2$. The two models share similar features. Their Feynmann graphs are decorated maps which gives them a nested structure and they have similar phases (dense or dilute) depending on the strength of the coupling related to this nested structure. 
\begin{question}
    Are there other duality relations between the $O(n)$-loop model and the heights model of type $A$ of~\cite{Ko95}?
\end{question}
 Establishing links with height models for other values of $n$ would be interesting as the height models of type $A^{(1)}_r$ satisfies the $r$-KdV hierarchy, that can be obtained as a reduction of the KP hierarchy~\cite{JM83}.
 
 %A strategy to investigate this question is to perform the integration over angular variables of the matrices in the $O(n)$-loop model to get its eigenvalue decomposition and compare this expression with Equation~$(28)$ of~\cite{Ko95}, rather than integrating directly over the matrices which appear quadratically.

\subsection{Tensorial HCIZ and hyperdeterminantal formulas}

Despite no results clearly relating the $b$-deformation for generic values of $b$ to hyperdeterminantal formulas, different independent results hint in that direction. This link was first suggested in~\cite{BE09} as a candidate for generalized HCIZ integral beyond the values of $\beta$ associated with the usual Lie groups. Hyperdeterminants also appear when considering the saddle-point expansion of angular integral in orientable constellations models (expressed as a multi-matrix model). A naive model of the $r$-components fermions of the KP hierarchy (i.e. the $r$-KdV hierarchy) where the fermions only interact component-wise also gives a hyperdeterminantal formula. Recently, a generalization of the HCIZ formula and Weingarten calculus for even-dimensional tensors was proposed in~\cite{CGL}. The authors also generalize monotone Hurwitz numbers to their framework and give a combinatorial interpretation of their results using constellations on nodal surfaces. The techniques developed in this article can in particular be applied to multi-matrix models by considering a $d$-matrix model as a tensor product of $d$ matrices. Therefore, the multi-matrix models can be recast as particular tensor models. The results of~\cite{CGL} are in particular applicable to the multi-matrix models relevant to the study of integrable hierarchies (e.g. constellations model or $O(n)$-loop model).
\begin{question}
    Are there hypotheses under which the tensor HCIZ of~\cite{CGL} gives an integrable model?
\end{question}
Since the techniques of~\cite{CGL} are quite general and apply to any multi-matrix model, while an arbitrary multi-matrix model has a priori no reason to be integrable, this is a natural question to ask. Conversely, using the tensor Weingarten calculus to compute the first moments for multi-matrix models related to integrable hierarchies can be interesting to conjecture on the existence of an exact integration formula for angular variables of these models. 
\begin{question}
    Can we use the tensor HCIZ of~\cite{CGL} to probe for generalization of the HCIZ formula to the multi-matrix case?
\end{question}
The (matrix) HCIZ formula is what allows to perform explicit integration of the angular variable in the Hermitian or unitary matrix model. It reduces the partition function of $1$-matrix models to the partition function of a Coulomb gas. Being able to integrate the angular degrees of freedom makes it easier to make contact with integrable properties e.g. orthogonal polynomials~\cite{EKR_Review} or to show that it is a solution of an integrable hierarchy. One difficulty of multi-matrix models is that the partition function cannot be reduced to an integral over eigenvalues of the matrices, even in multi-matrix models known to yield integrability features as with constellations. To be more explicit, in constellation models we get terms of the form  $\Tr M^{(0)}\dotsc M^{(m-1)}$ which lead to angular integral
\begin{equation}
    I^{(m)}\left[z,\Lambda^{(i)}\right] =\int_{U(N)^m} \prod\limits_{i=0}^{m-1 } dU^{(i)} \exp\left( z \Tr {U^{(0)}}^\dagger\Lambda^{(0)}U^{(0)}\dotsc {U^{(m-1)}}^\dagger\Lambda^{(m-1)}U^{(m-1)}\right)
\end{equation}
which can be seen as a generalization of the HCIZ formula to $m\geq 3$ matrices. Being able to compute this integral would help to get a deeper understanding of the interplay between the combinatorics of maps, matrix models and integrable hierarchies.

%% file: Appendices/App_cubic_constraint.tex
\chapter{Computation of the commutator for cubic constellations} % Main appendix title
\label{app:cubic_const} % Change X to a consecutive letter; for referencing this appendix elsewhere, use 

\section{Commutator for \texorpdfstring{$3$}{3}-constellations}
\label{sec:comm_3const}

This section contains the computation of the commutator involving the operator $A_i(3)$ with $A_j(s)$ with $s\leq 3$. Altogether with the computations included in Subsection~\ref{ssec:quad_const} for$b$-deformed bipartite maps, it leads to the extraction of the constraints for $b$-deformed $3$-constellations included in Subsection~\ref{sec:comm_3const}.

 \subsubsection{Computation of $[ A_i(1), A_j(3) ]$} \label{sec:CommutatorA1A3}
We start with $\left[ J^{(b)}_j, A_i(3)\right]$ and use~\eqref{eq:mod_rec_no_ui},
\begin{equation}
\left[ J^{(b)}_j, A_i(3)\right] = \sum_{n\geq 1} J_{i-n}^{(b)} [J^{(b)}_j, A_n(2)] + [J^{(b)}_j, J^{(b)}_{i-n}]A_n(2) + b(i-1)[J^{(b)}_j, A_i(2)].
\end{equation}
All commutators have been evaluated previously. This gives
\begin{multline}
\left[ J^{(b)}_j, A_i(3)\right] = j\left(\delta_{j<0} \sum_{l\geq j+1} + \sum_{l\geq \max(1, j+1)} + \delta_{i+j\geq 1}\sum_{l\geq 1}\right) J^{(b)}_{i+j-l} J^{(b)}_{l-1}\\
+bj\left( (i-1-j)\delta_{j<0} + (j+2i-2)\delta_{i+j\geq 1}\right) J^{(b)}_{i+j-1} - b^2 (i-1)^3 \delta_{i+j,1}.
\end{multline}
Therefore, by taking $j>0$ one finds
\begin{equation}
\left[p_j^*, A_i(3)\right] = j\left(\sum_{l\geq j} + \sum_{l\geq 1}\right) J^{(b)}_{i+j-l-1} p^*_{l} + bj(j+2i-2) p^*_{i+j-1},
\end{equation}
and therefore
%\begin{multline}
%\left[ A_j(1), A_i(3)\right] \\
%= (j-1)\left(\sum_{l\geq j} + \sum_{l\geq 1}\right) J^{(b)}_{i+j-1-l} A_l(1)
%+ b(j-1)(j-1+2(i-1)) A_{i+j-1}(1),
%\end{multline}
%and therefore
\begin{multline}
\left[p_i^*, A_j(3)\right] - \left[ p_j^*, A_i(3)\right] = b(i-j)(i+j-2) p^*_{i+j-1}\\
+ \left(i\sum_{l\geq i} - j\sum_{l\geq j} + (i-j)\sum_{l\geq 1}\right) J^{(b)}_{i+j-1-l} p_l^*.
\end{multline}
Introducing $M=\max(i,j)$ and $\mu=\min(i, j)$, it leads to
\begin{multline} \label{Explicitp*A3}
\left[p_i^*, A_j(3)\right] - \left[ p_j^*, A_i(3)\right] = b(i-j)(i+j-2) p^*_{i+j-1}\\
+ \left(2(i-j)\sum_{l\geq M} + (i-j)\sum_{l=1}^{M-1} - \operatorname{sgn}(i-j) \mu\sum_{l=\mu}^{M-1} \right) J^{(b)}_{i+j-1-l} p_l^*.
\end{multline}
as well as
\begin{multline} \label{ExplicitA1A3}
\left[ A_i(1), A_j(3)\right] - \left[ A_j(1), A_i(3)\right] = b(i-j)(i+j-2) A_{i+j-1}(1)\\
+ \left((i-j)\sum_{l\geq M} + (i-j)\sum_{l\geq 1} - \operatorname{sgn}(i-j) \sum_{l=\mu}^{M-1} (\mu-1)\right) J^{(b)}_{i+j-1-l} A_l(1).
\end{multline}

These quantities should be respectively given by the RHS of Equations \eqref{CommutatorpA} and \eqref{CommutatorAA}. Explicitly,
\begin{multline} \label{Thmp*A3}
\sum_{l\geq 1} D_{ij, l}(3) p_l^* = b(i-j)(i+j-2) p^*_{i+j-1}\\
+ \left(2(i-j)\sum_{l\geq M} + (i-j)\sum_{l=M}^{i+j-1} + \operatorname{sgn}(i-j)\sum_{l=\mu}^{M-1}(2l-3\mu+1)\right)J^{(b)}_{i+j-1-l} p_l^*
\end{multline}
and
\begin{multline} \label{ThmA1A3}
\sum_{l\geq 1} D_{ij, l}(1) A_l(3) + D_{ij, l}(3) A_l(1) = b(i-j)(i+j-2) A_{i+j-1}(1)\\
+ \left(2(i-j)\sum_{l\geq M} + (i-j)\sum_{l=M}^{i+j-1} + \operatorname{sgn}(i-j)\sum_{l=\mu}^{M-1}(2l-3\mu+1)\right)J^{(b)}_{i+j-1-l} A_l(1).
\end{multline}
Let us prove explicitly that the RHS of \eqref{ExplicitA1A3} and \eqref{ThmA1A3} are the same. Taking advantage of the fact that $A_l(1) = J^{(b)}_{l-1}/(1+b)$, we can make the change of summation index $l\to i+j-l$ and commute the two $J$s, so that
\begin{equation}
\sum_{l=\mu}^{M-1}(2l-3\mu+1)J^{(b)}_{i+j-1-l} J^{(b)}_{l-1} = \sum_{l=\mu+1}^M (2M-\mu-2l+1)J^{(b)}_{i+j-1-l} J^{(b)}_{l-1},
\end{equation}
and by then taking the half-sum,
\begin{equation}
\sum_{l=\mu}^{M-1}(2l-3\mu+1)J^{(b)}_{i+j-1-l} J^{(b)}_{l-1} = \left((M-\mu)\sum_{l=\mu+1}^{M-1} - (\mu-1)\sum_{l=\mu}^{M-1}\right)J^{(b)}_{i+j-1-l} J^{(b)}_{l-1}.
\end{equation}
The difference between the RHS of \eqref{ThmA1A3} and that of \eqref{ExplicitA1A3} is thus
\begin{multline}
\frac{(i-j)}{1+b}\left( \sum_{l\geq M} + \sum_{l=M}^{i+j-1} + \sum_{l=\mu+1}^{M-1} - \sum_{l\geq 1}\right)J^{(b)}_{i+j-1-l} J^{(b)}_{l-1}\\
= \frac{(i-j)}{1+b}\left(\sum_{l=\mu+1}^{i+j-1} - \sum_{l=1}^{M-1}\right)J^{(b)}_{i+j-1-l} J^{(b)}_{l-1}
\end{multline}
We consider the sum from $\mu+1$ to $i+j-1$ and change the summation index to $l\to i+j-l$ and commute the two $J$s which gives
\begin{equation}
\sum_{l=\mu+1}^{i+j-1} J^{(b)}_{i+j-1-l} J^{(b)}_{l-1} = \sum_{l=1}^{M-1} J^{(b)}_{i+j-1-l} J^{(b)}_{l-1},
\end{equation}
so that the above quantity vanishes and we have $\left[ A_i(1), A_j(3)\right] - \left[ A_j(1), A_i(3)\right] = \sum_{l\geq 1} D_{ij, l}(1) A_l(3) + D_{ij, l}(3) A_l(1)$, which is \eqref{CommutatorAA} for $s=1, s'=3$. The exact same steps can be followed to prove $\left[ p_i^*, A_j(3)\right] - \left[ p_j^*, A_i(3)\right] = \sum_{l\geq 1} D_{ij, l}(3) p_l^*$.

\subsubsection{Computation of $[A_i(2), A_j(3) ]$}
By writing again $A_j(3) = \sum_{l\geq 1} J^{(b)}_{j-l} A_l(2) + b(j-1) A_j(2)$, one finds
\begin{multline}
[ A_i(2), A_j(3) ] = \left(\sum_{n\geq j} (n-j) + \sum_{n=1}^{i+j-1} (n-j) + \sum_{n\geq i} (2i-n-1)\right) J^{(b)}_{i+j-1-n} A_n(2) \\
+ b((i-1)^2 + (j-1)(i-j)) A_{i+j-1}(2).
\end{multline}
It is then elementary to write
\begin{multline}
[ A_i(2), A_j(3) ] - [ A_j(2), A_i(3) ] = (i-j) A_{i+j-1}(3) + b(i-j)(i+j-2)A_{i+j-1}(2) \\
+ \left(2(i-j)\sum_{n\geq M} + (i-j) \sum_{n=M}^{i+j-1} + \operatorname{sgn}(i-j)\sum_{n=\mu}^{M-1}(2n-3\mu+1)\right) J^{(b)}_{i+j-1-n} A_n(2),
\end{multline}
i.e. $[ A_i(2), A_j(3) ] - [ A_j(2), A_i(3) ] = \sum_{n\geq 1} D_{ij, n}(2) A_n(3) + D_{ij, n}(3) A_n(2)$ as desired.

\subsection{Computation of $[ A_i(3), A_j(3) ]$} \label{sec:A3A3}

The last element of the proof of Proposition \ref{thm:SimplifiedCommutators} is the computation of $[ A_i(3), A_j(3) ]$, which is quite intricate. We provide a rather detailed proof in the present section. For bookkeeping purposes, we will use a notation which keeps track of both the origins and the types of all different contributions entering the calculation. The origin of a term will be denoted by a letter, $B, C, D, E, F, G, H$. The type of a term will refer to its the form of its summand. An additional number may appear if there are several terms with the same origin and of the same type. For instance,
\begin{equation*}
C^{JJA,2}_{ij} \coloneqq -(i-j)\sum_{\ell =1}^{M-\mu} \sum_{n=\mu+\ell}^{M-1} J^{(b)}_{i+j-1-n} J^{(b)}_{n-\ell} A_\ell(2)
\end{equation*}
has an origin indicated by the letter $C$, its type is $JJA$ (refering to its summand) and it is the second term of these origins and types.

We split the commutator as follows
\begin{multline}
    \left[A_i(3), A_j(3) \right] = \sum\limits_{n,\ell \geq 1} \underbrace{ \left[J^{(b)}_{i-n}A_n(2), J^{(b)}_{j-\ell} A_\ell(2)\right]}_{Q_{ij}} + \underbrace{b^2(i-1)(j-1) \left[A_i(2), A_j(2) \right]}_{H_{ij}} \\
    + \underbrace{b(i-1)\sum\limits_{\ell \geq 1}  \left[A_i(2), J^{(b)}_{j-\ell}A_\ell(2) \right]}_{G_{ij}} + \underbrace{b(j-1)\sum\limits_{n \geq 1}  \left[J^{(b)}_{i-n}A_n(2), A_j(2) \right]}_{-G_{ji}}
\end{multline}
We further expand the term $Q_{ij}$ as
\begin{multline}
    Q_{ij} = \underbrace{\sum\limits_{n,\ell \geq 1} J^{(b)}_{i-n}J^{(b)}_{j-\ell}\left[A_n(2), A_\ell(2)\right] }_{C_{ij}} + \underbrace{\sum\limits_{n,\ell \geq 1} \left[J^{(b)}_{i-n},J^{(b)}_{j-\ell}\right] A_\ell(2) A_n(2)}_{B_{ij}} \\
    + \underbrace{\sum\limits_{n,\ell \geq 1} J^{(b)}_{i-n}\left[A_n(2),J^{(b)}_{j-\ell}\right] A_\ell(2) }_{P_{ij}} + \underbrace{\sum\limits_{n,\ell \geq 1} J^{(b)}_{j-\ell}\left[J^{(b)}_{i-n}, A_\ell(2)\right]A_n(2)  }_{-P_{ji}}.
\end{multline}
We now analyze all terms so as to ``reform'' $\sum_{n\geq 1} D_{ij,n}(3) A_n(3)$.

\subsubsection{Term $B_{ij}$}
It writes
\begin{equation}
    B_{ij} = (1+b)\sum\limits_{n=1}^{i+j-1} (i-n) A_{i+j-n}(2) A_n(2).
\end{equation}
We make the change of summation index $n\to i+j-n$ and commute the two $A(2)$,
\begin{equation}
B_{ij} = (1+b) \sum\limits_{n=1}^{i+j-1} (n-j) A_{i+j-n}(2) A_n(2) + \frac{1+b}{6} (i+j)(i+j-1)(i+j-2) A_{i+j-1}(2).
\end{equation}
We then take the half-sum
\begin{equation}
     B_{ij} = \underbrace{\frac{1+b}{2} (i-j)\sum\limits_{n=1}^{i+j-1} A_{i+j-n}(2) A_n(2)}_{B^{AA}_{ij}} + \underbrace{\frac{1+b}{12} (i+j)(i+j-1)(i+j-2) A_{i+j-1}(2)}_{B^{A}_{ij}}
\end{equation}

\subsubsection{Term $C_{ij}$}
It expands as
\begin{equation}
C_{ij} = \sum_{n, l\geq 1} (n-l) J^{(b)}_{i-n} J^{(b)}_{j-l} A_{n+l-1}(2) = \sum_{k\geq 1} \sum_{n=1}^k (2n-k-1) J^{(b)}_{i-n} J^{(b)}_{j+n-k-1} A_k(2).
\end{equation}
We can either make the change $p=j+n-1$ and get
\begin{equation}
C_{ij} = \sum_{k\geq 1} \sum_{p=j}^{j+k-1} (2p-2j-k+1) J^{(b)}_{i+j-p-1} J^{(b)}_{p-k} A_k(2)
\end{equation}
or $p=i+k-n$ and get
\begin{equation}
C_{ij} = \sum_{k\geq 1} \sum_{p=i}^{i+k-1} (2i+k-2p-1)  J^{(b)}_{p-k} J^{(b)}_{i+j-p-1} A_k(2)
\end{equation}
By commuting the two $J$s the latter gives
\begin{multline}
C_{ij} = \sum_{k\geq 1} \sum_{p=i}^{i+k-1} (2i+k-2p-1) J^{(b)}_{i+j-p-1} J^{(b)}_{p-k} A_k(2) \\ \underbrace{- \frac{1+b}{6} (i+j)(i+j-1)(i+j-2) A_{i+j-1}(2)}_{2C_{ij}^{A,1}}.
\end{multline}
We then take the half-sum and split the terms as follows
{\small
\begin{multline}
C_{ij} = C_{ij}^{A,1} + \underbrace{ (i-j) \sum_{\ell \geq M-\mu+1} \sum_{n=M}^{\ell+\mu-1} J^{(b)}_{i+j-1-n} J^{(b)}_{n-\ell} A_\ell(2)}_{C^{JJA,1}_{ij}} \\
+ \frac{\operatorname{sgn}(j-i)}{2} \left[\underbrace{\sum_{\ell \geq 1} \sum_{n=\mu}^{M-1}}_{S_1} \underbrace{- \sum_{\ell =1}^{M-\mu} \sum_{n=\mu+\ell}^{M-1} }_{T_1}\right] (2\mu-2n-1+\ell) J^{(b)}_{i+j-1-n} J^{(b)}_{n-\ell} A_\ell(2) \\
- \frac{\operatorname{sgn}(j-i)}{2} \left[\underbrace{\sum_{\ell \geq 1} \sum_{n=\mu+\ell}^{M+\ell-1}}_{S_2} \underbrace{- \sum_{\ell =1}^{M-\mu} \sum_{n=\mu+\ell }^{M-1}}_{T_2}\right] (2M-2n-1+\ell) J^{(b)}_{i+j-1-n} J^{(b)}_{n-\ell} A_\ell(2)
\end{multline}}%
where we can already notice $C_{ij}^{A,1} = -B_{ij}^{A}$.

The sums $T_1+T_2$ get together as
\begin{equation}
C^{JJA,2}_{ij} \coloneqq -(i-j)\sum_{\ell =1}^{M-\mu} \sum_{n=\mu+\ell}^{M-1} J^{(b)}_{i+j-1-n} J^{(b)}_{n-\ell} A_\ell(2).
\end{equation}
In order to pack together $S_1$ and $S_2$, we relabel $n \rightarrow i+j-1-n-\ell$ in $S_2$ and commute the two $J$s. It comes
\begin{multline}
C_{ij} = C^{JJA,1}_{ij} + C^{JJA,2}_{ij} + C_{ij}^{A,1} + \underbrace{\frac{\operatorname{sgn}(j-i)}{2}\sum_{n=\mu}^{M-1} (2n-3\mu-M+2)(i+j-1-n) A_{i+j-1}(2)}_{C^{A,2}_{ij}} \\
 \underbrace{-\operatorname{sgn}(j-i) \sum_{\ell \geq 1} \sum_{n=\mu}^{M-1} (2n-2\mu+1-\ell) J^{(b)}_{i+j-1-n} J^{(b)}_{n-\ell} A_\ell(2)}_{\tilde{C}^{JJA}_{ij}}.
\end{multline}

\subsubsection{Term $P_{ij}$}
This contribution splits into three terms,
\begin{multline}
P_{ij} = \underbrace{\sum\limits_{n \geq j} b(n-j)^2J^{(b)}_{i+j-1-n} A_{n}(2)}_{F_{ij}} 
+ \underbrace{\sum\limits_{n, \ell \geq j}  (\ell-j) J^{(b)}_{i+j-1-n} J^{(b)}_{n-\ell} A_\ell(2)}_{D_{ij}} \\
+ \underbrace{\sum\limits_{n \geq j} \sum_{\ell=1}^{n} (\ell-j) J^{(b)}_{i+j-1-n} J^{(b)}_{n-\ell} A_\ell(2)}_{E_{ij}}.
\end{multline}
Recall that we do not need $P_{ij}$ itself but $P_{ij}-P_{ji}$, so we can compute directly the antisymmetrized contribution for each of the terms above. We have
\begin{multline}
    D_{ij}-D_{ji} = (i-j)\sum\limits_{n, \ell \geq M} J^{(b)}_{i+j-1-n} J^{(b)}_{n-\ell}A_\ell(2) \\
    +  \underbrace{\sgn(j-i) \left[ \sum\limits_{n \geq M}  \sum\limits_{\ell = \mu}^{M-1} +  \sum\limits_{n = \mu}^{M-1} \sum\limits_{\ell = \mu}^{M-1} +  \sum\limits_{n = \mu}^{M-1}  \sum\limits_{\ell \geq M}\right] (\mu-\ell) J^{(b)}_{i+j-1-n} J^{(b)}_{n-\ell} A_\ell(2)}_{\Omega^{(1)}_{ij}}.
\end{multline}
We use the first term above to form some $A_{n}(3)$,
\begin{multline}
(i-j)\sum\limits_{n, \ell \geq M} J^{(b)}_{i+j-1-n} J^{(b)}_{n-\ell}A_\ell(2) = (i-j)\sum\limits_{n \geq M}  J^{(b)}_{i+j-1-n} A_n(3) \\
\underbrace{-(i-j) \sum\limits_{n \geq M} \sum\limits_{\ell =1}^{M-1} J^{(b)}_{i+j-1-n} J^{(b)}_{n-\ell} A_\ell(2)}_{\alpha_{ij}}
\underbrace{-b(i-j)\sum\limits_{n \geq M} (n-1) J^{(b)}_{i+j-n-1} A_n(2)}_{D^{bJA,1}_{ij}}
\end{multline}

Similarly the sum over $n \geq M$ in $\alpha_{ij}$ can be completed to give a contribution of the form $A(2)A(2)$,
\begin{multline}
    \alpha_{ij} = - (1+b)(i-j) \sum_{\ell=1}^{M-1} A_{i+j-\ell}(2) A_\ell(2) + \underbrace{(i-j) \sum_{\ell=1}^{M-1} \sum_{n = \ell}^{M-1} J^{(b)}_{i+j-n-1} J^{(b)}_{n-\ell} A_\ell(2)}_{D^{JJA}_{ij}} \\
    + \underbrace{b(i-j) \sum_{\ell=1}^{M-1} (i+j-1-\ell) J^{(b)}_{i+j-1-\ell} A_\ell(2)}_{D^{bJA,2}_{ij}}
\end{multline}
In the first term, we commute the two $A(2)$, relabel $\ell \rightarrow i+j-\ell$ and take the half-sum, leading to 
\begin{multline}
    \alpha_{ij} = D^{JJA,1}_{ij} + D^{bJA,2}_{ij} \underbrace{-\frac{1+b}{2} (i-j) \sum_{\ell=1}^{i+j-1}A_{i+j-\ell}(2) A_\ell(2)}_{D^{AA,1}_{ij}} \\
     \underbrace{-\frac{1+b}{2} (i-j) \sum_{\ell=\mu+1}^{M-1} A_{i+j-\ell}(2) A_\ell(2)}_{D^{AA,2}_{ij}} + \underbrace{ \frac{1+b}{2}(i-j) \sum_{\ell=1}^{M-1} (2\ell-i-j) A_{i+j-1}(2)}_{D^A_{ij}}
\end{multline}

Finally, we split $\Omega_{ij}^{(1)}$ as 
\begin{multline}
\Omega_{ij}^{(1)} = \underbrace{\sgn(j-i) \sum\limits_{n \geq M}  \sum\limits_{\ell = \mu}^{M-1} (\mu-\ell) J^{(b)}_{i+j-1-n} J^{(b)}_{n-\ell} A_\ell(2)}_{\tilde{D}^{JJA,1}_{ij}} \\
+ \underbrace{\sgn(j-i) \sum\limits_{n = \mu}^{M-1} \sum\limits_{\ell \geq \mu} (\mu-\ell) J^{(b)}_{i+j-1-n} J^{(b)}_{n-\ell} A_\ell(2)}_{\tilde{D}^{JJA,2}_{ij}}
\end{multline}
All in all, we find the decomposition
\begin{multline}
D_{ij}-D_{ji}
= (i-j)\sum\limits_{n \geq M}  J^{(b)}_{i+j-1-n} A_n(3) \\
+ D^{AA,1}_{ij} + D^{AA,2}_{ij} + D^A_{ij} + D^{JJA}_{ij} + \tilde{D}^{JJA,1}_{ij} + \tilde{D}^{JJA,2}_{ij} + D^{bJA,1}_{ij} + D^{bJA,2}_{ij}
\end{multline}

For the term $E_{ij}$ we get
\begin{multline}
    E_{ij} - E_{ji} = \underbrace{(i-j)\sum\limits_{n \geq M} \sum\limits_{\ell = 1}^{n-1} J^{(b)}_{i+j-1-n}J^{(b)}_{n-\ell} A_\ell(2) }_{X_{ij}}\\
    + \underbrace{\sgn(j-i) \sum\limits_{n =\mu}^{M-1}\sum\limits_{\ell=1}^{n-1} (\mu-\ell) J^{(b)}_{i+j-1-n}J^{(b)}_{n-\ell} A_\ell(2)}_{\Omega^{(2)}_{ij}}
\end{multline}
The latter term is split as follows
\begin{multline}
\Omega^{(2)}_{ij} = \underbrace{\sgn(j-i) \sum\limits_{n =\mu}^{M-1}\sum\limits_{\ell=1}^{\mu-1} (\mu-\ell) J^{(b)}_{i+j-1-n}J^{(b)}_{n-\ell} A_\ell(2)}_{E^{JJA,3}_{ij}} \\
+ \underbrace{\sgn(j-i) \sum\limits_{n =\mu}^{M-1}\sum\limits_{\ell=\mu}^{n} (\mu-\ell) J^{(b)}_{i+j-1-n}J^{(b)}_{n-\ell} A_\ell(2)}_{E^{JJA,2}_{ij}}
\end{multline}

We complete the sum on $\ell$ to get a $A_n(3)$,
\begin{multline}
X_{ij} = (i-j) \sum\limits_{n \geq M} J^{(b)}_{i+j-1-n} A_{n}(3) \underbrace{- (i-j) \sum\limits_{n \geq M} \sum\limits_{\ell \geq n} J^{(b)}_{i+j-1-n} J^{(b)}_{n-\ell} A_{\ell}(2)}_{E^{JJA,1}_{ij}} \\
\underbrace{-b(i-j) \sum\limits_{n \geq M} (n-1) J^{(b)}_{i+j-1-n} A_{n}(2)}_{E^{bJA}_{ij}}
\end{multline}
hence
\begin{equation}
E_{ij} - E_{ji} = (i-j) \sum\limits_{n \geq M} J^{(b)}_{i+j-1-n} A_{n}(3) + E^{JJA,1}_{ij} + E^{JJA,2}_{ij} + E^{JJA,3}_{ij} + E^{bJA}_{ij}.
\end{equation}

And $F_{ij}$ gives
\begin{multline}
    F_{ij} - F_{ji} = \underbrace{b\sum\limits_{n \geq M} (i-j)(2n-i-j) J^{(b)}_{i+j-1-n} A_{n}(2)}_{F^{bJA}_{ij}} \\
    + \underbrace{b\sgn(i-j) \sum\limits_{n=\mu}^{M-1} (n-\mu)^2 J^{(b)}_{i+j-1-n} A_{n}(2)}_{\tilde{F}^{bJA}_{ij}}.
\end{multline}

\subsubsection{Terms $G_{ij}$ and $H_{ij}$}
The term $H_{ij} = b^2(i-1)(j-1) \left[A_i(2), A_j(2) \right]$ is simply
\begin{equation}
    H_{ij} = b^2(i-1)(j-1)(i-j) A_{i+j-1}(2).
\end{equation}

As for $G_{ij} = b(i-1)\sum_{l\geq 1} [A_i(2), J^{(b)}_{j-l} A_l(2)]$,
\begin{multline}
    G_{ij} = \frac{b^2}{1+b}(i-1)^3 A_{i+j-1}(2) \\
    + b(i-1) \left(\sum_{l\geq i} (2i-l-1) + \sum_{l\geq j} (l-j) + \sum_{l=1}^{i+j-1} (l-j)\right) J^{(b)}_{i+j-l-} A_l(2).
\end{multline}
Antisymmetrization leads to
\begin{multline}
    G_{ij}-G_{ji} = \underbrace{b^2\left((i-1)^3-(j-1)^3\right) A_{i+j-1}(2)}_{G^{A}_{ij}} \\
    + \underbrace{b(i-j)\sum\limits_{n=1}^{i+j-1} (n-1) J^{(b)}_{i+j-1-n} A_{n}(2)}_{G^{bJA,1}_{ij}}
    + \underbrace{2b(i-j)(i+j-2)\sum\limits_{n \geq M}  J^{(b)}_{i+j-1-n} A_{n}(2)}_{G^{bJA,2}_{ij}} \\
    + \underbrace{b\sgn(i-j) \sum\limits_{n=\mu}^{M-1} \left[(M-1)(n-\mu) - (\mu-1)(2\mu-n-1)\right] J^{(b)}_{i+j-1-n} A_{n}(2)}_{\tilde{G}^{bJA}_{ij}}
\end{multline}

\subsubsection{Repackaging the contributions}
We start with adding together the three terms $\tilde{C}^{JJA}_{ij}, \tilde{D}^{JJA,2}_{ij}, \tilde{E}^{JJA,1}_{ij}$.
\begin{equation}
\begin{aligned}
&\tilde{C}^{JJA}_{ij}+ \tilde{D}^{JJA}_{ij}+ \tilde{E}^{JJA,1}_{ij}
= \sgn(j-i)\\
& \sum_{n=\mu}^{M-1} \left(\sum_{l\geq \mu} (\mu-l) + \sum_{l=1}^{\mu-1} (\mu-l) - \sum_{l\geq 1}(2n-2\mu-l+1)\right)J^{(b)}_{i+j-1-n} J^{(b)}_{n-l} A_l(2)\\
&= \sgn(j-i) \sum_{n=\mu}^{M-1} \sum_{l\geq 1} (3\mu-2n-1) J^{(b)}_{i+j-1-n} J^{(b)}_{n-l} A_l(2)\\
&= \sgn(j-i) \sum_{n=\mu}^{M-1} (3\mu-2n-1) J^{(b)}_{i+j-1-n} (A_n(3) - b (n-1) A_n(2))\\
&= \sgn(j-i) \sum_{n=\mu}^{M-1} (3\mu-2n-1) J^{(b)}_{i+j-1-n} A_n(3) + \widetilde{{CE}}^{bJA}_{ij}
\end{aligned}
\end{equation}
with $\widetilde{{CE}}^{bJA}_{ij} = - b\sgn(j-i) \sum_{n=\mu}^{M-1} (3\mu-2n-1)(n-1) J^{(b)}_{i+j-1-n} A_n(2)$. Then
\begin{equation}
\begin{aligned}
D^{JJA,2}_{ij} + E^{JJA,2}_{ij} &= \sgn(j-i) \sum_{l=\mu}^{M-1}(\mu-l) \sum_{n\geq l} J^{(b)}_{i+j-1-n} J^{(b)}_{n-l} A_l(2)\\
&= \sgn(j-i) \sum_{l=\mu}^{M-1}(\mu-l) \sum_{p\geq 1} J^{(b)}_{i+j-1-p-l} J^{(b)}_{p-1} A_l(2)\\
&= \sgn(j-i) (1+b) \sum_{l=\mu}^{M-1}(\mu-l) (A_{i+j-l}(2) - b (i+j-l-1) A_{i+j-l}(1)) A_l(2).
\end{aligned}
\end{equation}
The term with $A_{i+j-l}(2) A_l(2)$ can be rewritten by making the change of index variable $l\to i+j-l$ and commuting the two operators,
\begin{multline}
\sgn(j-i) (1+b) \sum_{l=\mu+1}^{M-1}(\mu-l) A_{i+j-l}(2) A_l(2)
= \sgn(j-i) (1+b) \sum_{l=\mu+1}^{M-1}(l-M) A_{i+j-l}(2) A_l(2) \\+ \sgn(j-i) (1+b) \sum_{l=\mu+1}^{M-1}(l-M)(2l-i-j) A_{i+j-1}(2).
\end{multline}
then by taking the half-sum, we get
\begin{multline}
D^{JJA,2}_{ij} + E^{JJA,2}_{ij} = \underbrace{\frac{1+b}{2} (i-j) \sum_{l=\mu+1}^{M-1} A_{i+j-l}(2) A_l(2)}_{(DE)^{AA}_{ij}}\\
+ \underbrace{\sgn(i-j) b\sum_{l=\mu}^{M-1} (\mu-l)(i+j-l-1) J^{(b)}_{i+j-l-1} A_l(2)}_{\widetilde{(DE)}^{bJA}_{ij}}\\
+ \underbrace{\sgn(j-i) \frac{1+b}{2} \sum_{l=\mu+1}^{M-1}(l-M)(2l-i-j) A_{i+j-1}(2)}_{(DE)^A_{ij}}
\end{multline}

Next up, we look at $C^{JJA,1}_{ij} + E^{JJA,1}_{ij}$,
\begin{equation}
\begin{aligned}
C^{JJA,1}_{ij} + E^{JJA,1}_{ij} &= (i-j) \sum_{n\geq M} \sum_{l=n-\mu+1}^{n-1} J^{(b)}_{i+j-1-n} J^{(b)}_{n-l} A_{l}(2)\\
&= (i-j) \Biggl( \underbrace{\sum_{l\geq M-1} \sum_{n=l+1}^{l+\mu-1}}_{z_{ij}} + \underbrace{\sum_{l=M-\mu+1}^{M-2} \sum_{n=M}^{l+\mu-1}}_{(CE)^{JJA,1}_{ij}}\Biggr) J^{(b)}_{i+j-1-n} J^{(b)}_{n-l} A_{l}(2)
\end{aligned}
\end{equation}
then
\begin{equation}
\begin{aligned}
z_{ij} &= (i-j) \sum_{l\geq M-1} \sum_{p=1}^{\mu-1} J^{(b)}_{i+j-1-l-p} J^{(b)}_{p} A_{l}(2)\\
&= \underbrace{(i-j) \sum_{l\geq M-1} \sum_{p=1}^{\mu-1}  J^{(b)}_{p} J^{(b)}_{i+j-1-l-p} A_{l}(2)}_{\zeta_{ij}} \underbrace{- (i-j) (1+b)\frac{\mu(\mu-1)}{2} A_{i+j-1}(2)}_{(CE)^A_{ij}}
\end{aligned}
\end{equation}
In the first line we have set $p=n-l$, in the second line we have commuted the two $J$s. In $\zeta^{(1)}_{ij}$ we form $A_{i+j-1-p}(3)$
\begin{multline}
\zeta_{ij} = (i-j) \sum_{p=1}^{\mu-1} J^{(b)}_{p} (A_{i+j-p-1}(3) - b(i+j-p-2) A_{i+j-1-p}(2))\\
- (i-j) \sum_{p=1}^{\mu-1} \sum_{l=1}^{M-2} J^{(b)}_{p} J^{(b)}_{i+j-1-l-p} A_{l}(2)
\end{multline}
and making the change $n=i+j-1-p$,
\begin{multline}
\zeta_{ij} = (i-j) \sum_{n=M}^{i+j-1} J^{(b)}_{i+j-1-n} A_{n}(3) \underbrace{- b(i-j) \sum_{n=M}^{i+j-1} (n-1) J^{(b)}_{i+j-1-n} A_{n}(2)}_{(CE)^{bJA}_{ij}}\\
\underbrace{- (i-j) \sum_{n=M}^{i+j-2} \sum_{l=1}^{M-2} J^{(b)}_{i+j-1-n} J^{(b)}_{n-l} A_{l}(2)}_{(CE)^{JJA,2}_{ij}}.
\end{multline}
Overall
\begin{multline}
C^{JJA,1}_{ij} + E^{JJA,1}_{ij} = (i-j) \sum_{n=M}^{i+j-1} J^{(b)}_{i+j-1-n} A_{n}(3) \\
+ (CE)^{JJA,1}_{ij} + (CE)^{JJA,2}_{ij} + (CE)^A_{ij} + (CE)^{bJA}_{ij}.
\end{multline}

%%%%%%%%%%%%%%%%%%%%%%%%%%
%\begin{align}
%\Omega^{(1)}+\Omega^{(2)}+\Omega^{(3)} &= \frac{\sgn(j-i)}{1+b}\sum\limits_{n=\mu}^{M-1}(3\mu-2n-1)J^{(b)}_{i+j-1-n}M^{(k-1,1,3)}_{n} \nonumber \\ 
%&\underbrace{-\frac{b\sgn(j-i)}{1+b}\sum\limits_{n=\mu}^{M-1}(3\mu-2n-1)(n-1)J^{(b)}_{i+j-1-n} M^{(k-1,1,2)}_{n}}_{r^{(4)}} \nonumber\\
%&+ \underbrace{\frac{\sgn(j-i)}{1+b} \frac{1}{2} \sum\limits_{n=\mu+1}^{M-1}(\mu-M) M^{(k-1,1,2)}_{i+j-n}M^{(k-1,1,2)}_n}_{-\alpha^{(4)}} \nonumber\\
%&+ \underbrace{\frac{\sgn(j-i)}{1+b} \frac{1}{2} \sum\limits_{\ell=\mu+1}^{M}(\mu-\ell)(M+\mu-2\ell)M^{(k-1,1,2)}_{i+j-1}}_{\omega} \nonumber\\
%&\underbrace{-\frac{\sgn(j-i)}{1+b} \sum\limits_{n=\mu}^{M-1}b(n-\mu) (i+j-1-n)J^{(b)}_{i+j-1-n}M^{(k-1,1,2)}_n}_{r^{(5)}} \nonumber
%\end{align}
%%
%Similarly we bring together the terms $Y_{ij}$ and $Z_{ij}$. We form a $M^{(k-1,1,3)}$ contribution and rewrite the other term by
%\begin{align}
%Y_{ij} + Z_{ij} &= \frac{(i-j)}{1+b}\sum\limits_{n=1}^{\mu-1} J^{(b)}_n M^{(3,1)}_{i+j-1-n} \underbrace{-\frac{(i-j)}{1+b}\sum\limits_{\ell=1}^M \sum\limits_{\substack{n=M \\ n \neq \ell }}^{\ell} J^{(b)}_{i+j-1-n} J^{(b)}_{n-\ell} M^{(k-1,1,2)}_\ell }_{\alpha^{(5)}} \\
%&\underbrace{-\frac{(i-j)}{1+b}\sum\limits_{n=1}^{\mu-1} b(i+j-2-n)J^{(b)}_nM^{(k-1,1,2)}_{i+j-1-n}}_{r^{(6)}} .\nonumber
%\end{align}

We package together some of the contributions which receive some $b$ during the calculation. First,
\begin{equation}
    H_{ij} + G^{A}_{ij} = b^2(i-j)(i+j-2)^2 A_{i+j-1}(2).
\end{equation}
Moreover,
\begin{align}
D^{bJA,1}_{ij} + E^{bJA}_{ij} + F^{bJA}_{ij} + G^{bJA,2}_{ij} &= b(i-j)(i+j-2) \sum_{n\geq M} J^{(b)}_{i+j-1-n} A_n(2)\\
(CE)^{bJA}_{ij} + G^{bJA,1}_{ij} + D^{bJA,2}_{ij} &= b(i-j)(i+j-2) \sum_{n=1}^{M-1} J^{(b)}_{i+j-1-n} A_n(2)
\end{align}
hence the sum of the above three contributions reduces to
\begin{multline}
(D^{bJA,1}_{ij} + E^{bJA}_{ij} + F^{bJA}_{ij} + G^{bJA,2}_{ij}) + ((CE)^{bJA}_{ij} + G^{bJA,1}_{ij} + D^{bJA,2}_{ij}) + (H_{ij} + G^{A}_{ij}) \\
= b(i-j)(i+j-2) A_{i+j-1}(3).
\end{multline}

\subsubsection{Cancellations}
We notice the following cancellations,
\begin{equation}
\begin{aligned}
C^{A,1}_{ij} + B^{A}_{ij} &= 0, \\
D^{AA,1}_{ij} + B^{AA}_{ij} &= 0, \\
D^{AA,2}_{ij} + (DE)^{AA}_{ij} &=0, \\
C^{A,2}_{ij} + D^A_{ij} + (DE)^A_{ij} + (CE)^A_{ij}  &= 0, \\
C^{JJA,2}_{ij} + D^{JJA,1}_{ij} + (CE)^{JJA,1}_{ij} + (CE)^{JJA,2}_{ij} &= 0, \\
\widetilde{(CE)}^{bJA}_{ij} + \widetilde{(DE)}^{bJA}_{ij} + \tilde{G}^{bJA}_{ij} + \tilde{F}^{bJA}_{ij}  &= 0.
\end{aligned}
\end{equation}

\subsubsection{Final expression}
Bringing together the remaining contributions, we get
\begin{multline}
\left[ A_i(3) , A_j(3) \right] = 2(i-j) \sum\limits_{n \geq M} J^{(b)}_{i+j-1-n} A_n(3)
+ (i-j) \sum\limits_{n=M}^{i+j-1} J^{(b)}_{i+j-1-n} A_n(3) \\
+ \sgn(i-j) \sum\limits_{n=\mu}^{M-1}(2n-3\mu+1) J^{(b)}_{i+j-1-n} A_n(3) + b(i-j)(i+j-2)A_{i+j-1}(3)
\end{multline}
which concludes the proof of Proposition~\ref{thm:SimplifiedCommutators}.

\section{Commutator for at most trivalent bipartite maps}
\label{sec:comm_triv_bip}

This section is similar to the previous one. It contains the computation of all commutator of operators $M^{(1,m)}$ up to $m=3$. We recall that for this family of operators, we set $J_0 = u$.  The general outline of the computation is similar to the one of $3$-constellations but the shift due to the extra operator $Y_+$ induce minor changes in some of the terms.

\subsubsection{Computation of $\left[M_i^{(1,1)}, M_j^{(1,1)}\right]$}
We have $M^{(1,1)}_i = \frac{J^{(b)}_{i-1}}{1+b}$, then $\left[J^{(b)}_j, M_i^{(1,1)} \right] = j\delta_{i+j,1}$. Therefore $[M_i^{(1,1)}, M_j^{(1,1)}] = [M_i^{(1,1)}, p_j^*] = 0 $, which are indeed compatible with $\tilde{D}_{ij,l}(1)=0$.

\subsubsection{Computation of $\left[ M_i^{(1,m)}, M_j^{(1,m')}\right]$ for $m, m'\leq 2$}
We use relation~\eqref{eq:mod_rec_no_ui'} to compute the various commutators. We have
\begin{equation}
    \left[ J^{(b)}_j, M_i^{(1,2)}  \right] = (1+b) j M^{(1,1)}_{i+j-1}\left( \delta_{i+j \geq 2} +  \delta_{j \leq 0} \right) - b (i-1)(i-2) \delta_{i+j,2}.
\end{equation}

This $[p_j^*, M^{(1,2)}_i] = jJ^{(b)}_{i+j-2}$ for $i, j\geq 1$ and $[M^{(1,1)}_j, M^{(1,2)}_i] = (j-1) M^{(1,1)}_{i+j-2}$, hence
\begin{equation}
\begin{aligned}
    [M^{(1,2)}_i, p_j^*] - [M^{(1,2)}_j, p_i^*] &= (i-j) p_{i+j-2}^*\\ 
    [M^{(1,2)}_i, M^{(1,1)}_j] - [M^{(1,2)}_j,M^{(1,1)}_i] &= (i-j) M^{(1,1)}_{i+j-2}
\end{aligned}
\end{equation}
which correspond to \eqref{CommutatorMM} and \eqref{CommutatorpM} with $m'=1, m=2$. Furthermore, the calculation of $[M^{(1,2)}_i, M^{(1,2)}_j]$ follows exactly the same steps as $[A_i(2), A_j(2)]$ in Section \ref{sec:A2A2} and again produces a half-Virasoro algebra,
\begin{equation}
[M^{(1,2)}_i, M^{(1,2)}_j] = (i-j) M^{(1,2)}_{i+j-2}.
\end{equation}
which is also compatible with $\tilde{D}_{ij,l}(2) = (i-j) \delta_{l,i+j-2}$.

At this stage, we have Lemma~\ref{lemma:SimplifiedCommutators2} for bipartite maps with vertices of valency at most $2$, which in particular includes $b$-deformed maps when $q_1=0$. This leads to Theorem~\ref{thm:constr_bmaps}.

\subsubsection{Computation of $\left[ M^{(1,1)}_i, M^{(1,3)}_j\right]$}
Again, we use~\eqref{eq:mod_rec_no_ui'} starting with
\begin{multline}
 \left[ J^{(b)}_j, M^{(1,3)}_i  \right] = -b^2 (i-1)(i-2)(i-3) \delta_{i+j,3} - b j(j-1) J^{(b)}_{i+j-3}\delta_{j \leq 1} \\ 
+ b(i-1)j J^{(b)}_{i+j-3} \left( \delta_{i+j \geq 3 } + \delta_{j \leq 0 } \right) + (1+b)j M^{(1,2)}_{i+j-1}\delta_{i+j\geq 2} \\
+ \sum\limits_{n\geq j} (1+b) j J^{(b)}_{i+j-n-2} M^{(1,1)}_{n}\left( \delta_{n \geq1} +  \delta_{j\leq 0} \right).
\end{multline}
For $j\geq 1$, this leads to
\begin{equation}
[p_j^*,  M^{(1,3)}_i] = j\left(\sum_{l\geq 0} + \sum_{l\geq j-1}\right) J^{(b)}_{i+j-3-l} J^{(b)}_l + \frac{b}{1+b}j(2i+j-3) J^{(b)}_{i+j-3} \delta_{i+j\geq 3}
\end{equation}
and therefore (recall that $p^*_n = 0$ if $n\leq 0$)
\begin{multline} \label{Commutatorp*M13}
[p_i^*,  M^{(1,3)}_j] - [p_j^*,  M^{(1,3)}_i] = b(i-j)(i+j-3) p^*_{i+j-3}\\
+ \frac{1}{1+b}\left(2(i-j) \sum_{l\geq M-1} -\sgn(i-j) \mu \sum_{l=\mu-1}^{M-2} + (i-j)\sum_{l=0}^{M-2}\right)J^{(b)}_{i+j-3-l} J^{(b)}_l
\end{multline}
and (recall that $M^{(1,1)}_n = 0$ if $n\leq 0$)
\begin{multline} \label{CommutatorM11M13}
[M^{(1,1)}_i,  M^{(1,3)}_j] - [M^{(1,1)}_j,  M^{(1,3)}_i] = b(i-j)(i+j-3) M^{(1,1)}_{i+j-3}\\
\left(2(i-j) \sum_{l\geq M-1} -\sgn(i-j) (\mu-1) \sum_{l=\mu-1}^{M-2} + (i-j)\sum_{l=1}^{M-2}\right)J^{(b)}_{i+j-3-l} M^{(1,1)}_l.
\end{multline}
The RHS of those equations should be respectively equal to 
\begin{multline} \label{DTilde3p*}
\sum_{l\geq 1} \tilde{D}_{ij, l}(3) p_l^* = b(i-j)(i+j-3) p_{i+j-3}^*\\
+ \left(2(i-j)\sum_{l\geq M-1} + (i-j) \sum_{l=M-1}^{i+j-3} + \sgn(i-j) \sum_{l=\mu-1}^{M-2} (2l-3\mu+3)\right) J^{(b)}_{i+j-3-l} p_l^*
\end{multline}
and
\begin{multline}
\sum_{l\geq 1} \tilde{D}_{ij, l}(3) M^{(1,1)}_l + \tilde{D}_{ij, l}(1) M^{(1,3)}_l = b(i-j)(i+j-3) M^{(1,1)}_{i+j-3}\\
+ \left(2(i-j)\sum_{l\geq M-1} + (i-j) \sum_{l=M-1}^{i+j-3} + \sgn(i-j) \sum_{l=\mu-1}^{M-2} (2l-3\mu+3)\right) J^{(b)}_{i+j-3-l} M^{(1,1)}_l
\end{multline}
This can be proven exactly as we did for 3-constellations in Section \ref{sec:CommutatorA1A3}, and therefore we do not repeat it. There is however a subtlety in equating \eqref{Commutatorp*M13} with \eqref{DTilde3p*}, because the former is written with $J^{(b)}_l$, where $l=0$ is not vanishing, instead of $p_l^*$ for \eqref{DTilde3p*}. To deal with this, one can separate three cases: $i)$ $\mu>1$, $ii)$ $\mu=1$ and $M>2$, $iii)$ $\mu=1, M=2$ and prove the equality in all of them.

\subsubsection{Computation of $[M^{(1,2)}_i, M^{(1,3)}_j]$}
By writing again $M^{(1,3)}_j = \sum_{l\geq 1} J^{(b)}_{j-l-1} M^{(1,2)}_l + b(j-1) M^{(1,2)}_{j-1}$, one finds
\begin{multline}
[M^{(1,2)}_i, M^{(1,3)}_j] =  b((i-1)(i-2) + (j-1)(i-j+1)) M^{(1,2)}_{i+j-3}\\
+\left(\sum_{n\geq j-1} (n-j+1) + \sum_{n=1}^{i+j-3} (n-j+1) + \sum_{n\geq i-1} (2i-n-2)\right) J^{(b)}_{i+j-3-n} M^{(1,2)}_{n}.
\end{multline}
It is then elementary to write
\begin{multline}
[M^{(1,2)}_i, M^{(1,3)}_j] - [M^{(1,2)}_j, M^{(1,3)}_i] = (i-j) M^{(1,3)}_{i+j-2} + b(i-j)(i+j-3) M^{(1,2)}_{i+j-3} \\
+ \left(2(i-j)\sum_{n\geq M-1} + (i-j) \sum_{n=M-1}^{i+j-3} + \operatorname{sgn}(i-j)\sum_{n=\mu-1}^{M-2}(2n-3\mu+3)\right) J^{(b)}_{i+j-3-n} M^{(1,1)}_n,
\end{multline}
i.e. $[M^{(1,2)}_i, M^{(1,3)}_j] - [M^{(1,2)}_j, M^{(1,3)}_i] = \sum_{n\geq 1} \tilde{D}_{ij, n}(2) M^{(1,3)}_n + \tilde{D}_{ij, n}(3) M^{(1,2)}_n$ as desired.

\subsection{Computation of $\left[ M^{(1,3)}_i, M^{(1,3)}_j\right]$}
This computation essentially follows the same steps as $[A_i(3), A_j(3)]$ for 3-constellations in Section~\ref{sec:A3A3}. In order to emphasize the similarty with $[A_i(3), A_j(3)]$, we decompose the calculation of $[M^{(1,3)}_i, M^{(1,3)}_j]$ in the same way and to avoid cumbersome notation, we reallocate the notation of Section \ref{sec:A3A3} to the analogous quantities they represent in our decomposition of $[M^{(1,3)}_i, M^{(1,3)}_j]$. The difference in the notation reallocation is that the summands are now of the form $JJM, JM, bJM, \dotsc$

The commutator $[ M^{(1,3)}_i, M^{(1,3)}_j]$ expands as

\begin{multline}
[ M^{(1,3)}_i, M^{(1,3)}_j] = \underbrace{\sum\limits_{\substack{\ell \geq 1 }} \sum\limits_{\substack{n \geq 1 }} \left[ J^{(b)}_{i-1-\ell}M^{(1,2)}_\ell,J^{(b)}_{j-1-l}M^{(1,2)}_n\right]}_{Q_{ij}} + \underbrace{b^2(j-1)(i-1) \left[M^{(1,2)}_{i-1},M^{(1,2)}_{j-1} \right]}_{H_{ij}} \nonumber\\
+ \underbrace{b(j-1)\sum\limits_{\substack{\ell \geq 1 }} \left[ J^{(b)}_{i-1-\ell}M^{(1,2)}_\ell,M^{(1,2)}_{j-1}\right] }_{G_{ij}} +\underbrace{b(i-1)\sum\sum\limits_{\substack{n \geq 1 }} \left[M^{(1,2)}_{i-1},J^{(b)}_{j-1-n}M^{(1,2)}_n\right] }_{-G_{ji}} .
\end{multline}%
We decompose $Q_{ij}$ by expanding the commutator
\begin{multline}
Q_{ij} = \underbrace{\sum\limits_{\substack{\ell \geq 1 }} \sum\limits_{\substack{n \geq 1 }} J^{(b)}_{i-1-\ell} [M^{(1,2)}_\ell, J^{(b)}_{j-1-n}] M^{(1,2)}_n}_{P_{ij}} + \underbrace{\sum\limits_{\substack{\ell \geq 1 }} \sum\limits_{\substack{n \geq 1 }} J^{(b)}_{j-1-n} [J^{(b)}_{i-1-\ell} , M^{(1,2)}_n] M^{(1,2)}_\ell}_{-P_{ji}} \nonumber\\
+ \underbrace{\sum\limits_{\substack{\ell \geq 1 }} \sum\limits_{\substack{n \geq 1 }} [J^{(b)}_{i-1-\ell} , J^{(b)}_{j-1-n}] M^{(1,2)}_n M^{(1,2)}_\ell}_{B_{ij}} + \underbrace{\sum\limits_{\substack{\ell \geq 1 }} \sum\limits_{\substack{n \geq 1 }} J^{(b)}_{i-1-\ell} J^{(b)}_{j-1-n} [ M^{(1,2)}_\ell, M^{(1,2)}_n]}_{C_{ij}}. 
\end{multline}
We now rewrite each of these contributions individually before assembling them to show that the commutator closes.

\subsubsection{Term $B_{ij}$} This term writes
\begin{equation}
B_{ij} = (1+b) \sum\limits_{\ell = 1}^{i+j-3} (i-1-\ell) M^{(1,2)}_{i+j-2-\ell} M^{(1,2)}_\ell .
\end{equation}
Now commuting the two factors $M^{(1,2)}$, relabelling the sum using $\ell\to i+j-2-\ell$ and taking half-sums of the two expressions gives
\begin{multline}
B_{ij} = \underbrace{\frac{1+b}{12} (i+j-2)(i+j-3)(i+j-4) M^{(1,2)}_{i+j-4}}_{B^{M}_{ij}} \\+ \underbrace{\frac{1+b}{2}(i-j) \sum\limits_{\ell=1}^{i+j-3} M^{(1,2)}_{i+j-2-\ell} M^{(1,2)}_\ell}_{B^{MM}_{ij}}.
\end{multline}

\subsubsection{Term $C_{ij}$} It is given by
\begin{equation}
    C_{ij} = \sum\limits_{\substack {\ell \geq 1 }} \sum\limits_{\substack {n \geq 1 }} (\ell-n) J^{(b)}_{i-1-\ell}  J^{(b)}_{j-1-n} M^{(1,2)}_{\ell+n-2} = \sum\limits_{\substack{ n \geq 0 }} \sum\limits_{\substack{\ell = j-1 }}^{n+j-1}  (2\ell-2j-n+2) J^{(b)}_{i+j-3-\ell} J^{(b)}_{\ell-1-n} M^{(1,2)}_n.
\end{equation}

Following the exact same steps as in Section~\ref{sec:comm_3const} we arrive at
\begin{equation}
C_{ij} = C^{JJM,1}_{ij} + C^{JJM,2}_{ij} + C_{ij}^{M,1} + C^{M,2}_{ij} + \tilde{C}^{JJM}_{ij}
\end{equation}
with
\begin{equation}
\begin{aligned}
C^{JJM,1}_{ij} &= (i-j) \sum_{ \substack{\ell \geq M-\mu}} \sum_{\substack{n=M-1 }}^{\ell+\mu-1} J^{(b)}_{i+j-3-n} J^{(b)}_{n-1-\ell} M^{(1,2)}_\ell \\
C^{JJM,2}_{ij} &= (i-j)\sum_{\ell =0}^{M-\mu-1} \sum_{\substack{ n=\mu+\ell}}^{M-2} J^{(b)}_{i+j-3-n} J^{(b)}_{n-1-\ell}M^{(1,2)}_\ell\\
C_{ij}^{M,1} &= - B^{M}_{ij}\\
C^{M,2}_{ij} &= (1+b)\frac{\operatorname{sgn}(j-i)}{2}\sum_{n=\mu-1}^{M-2} (2n+6-3\mu-M)(i+j-3-n) M^{(1,2)}_{i+j-4}\\
\tilde{C}^{JJM}_{ij} &= \operatorname{sgn}(j-i) \sum_{\substack{\ell \geq 1 }} \sum_{n=\mu-1}^{M-2} (2\mu-2n+\ell-2) J^{(b)}_{i+j-3-n} J^{(b)}_{n-1-\ell} M^{(1,2)}_\ell
\end{aligned}
\end{equation}

\subsubsection{Term $P_{ij}$}
This term expands as 
\begin{multline}
    P_{ij} = \underbrace{ b \sum\limits_{n \geq j-1} (n-j)(n-j+1) J^{(b)}_{i+j-3-n} M^{(1,2)}_{n-1}}_{F_{ij}}
    + \underbrace{\sum\limits_{n,\ell \geq j-1} (\ell-j+1) J^{(b)}_{i+j-3-n} J^{(b)}_{n-1-\ell}  M^{(1,2)}_\ell}_{D_{ij}} \nonumber\\
    + \underbrace{\sum\limits_{n \geq j-1} \sum\limits_{ \substack{ \ell =1 }}^{n-1} (\ell-j+1) J^{(b)}_{i+j-3-n} J^{(b)}_{n-1-\ell}  M^{(1,2)}_\ell}_{E_{ij}}. \nonumber
\end{multline}
Following the same steps as for the previous model in Section~\ref{sec:comm_3const}, we arrive at
\begin{multline}
D_{ij}-D_{ji}
= (i-j) \sum\limits_{\substack{n \geq M-1 }} J^{(b)}_{i+j-3-n} M^{(1,3)}_n \\
+ D^{MM,1}_{ij} + D^{MM,2}_{ij} + D^M_{ij} + D^{JJM}_{ij} + \tilde{D}^{JJM,1}_{ij} + \tilde{D}^{JJM,2}_{ij} + D^{bJM,1}_{ij} + D^{bJM,2}_{ij}
\end{multline}
with
\begin{equation}
\begin{aligned}
D^{MM,1}_{ij} &= -(1+b)\frac{(i-j)}{2}\sum\limits_{\ell = 1}^{M+\mu-3}  M^{(1,2)}_{i+j-2-\ell} M^{(1,2)}_\ell\\
D^{MM,2}_{ij} &= -(1+b)\frac{(i-j)}{2}\sum\limits_{\ell = \mu}^{M-2}  M^{(1,2)}_{i+j-2-\ell} M^{(1,2)}_\ell\\
D^M_{ij} &= (1+b)\frac{(i-j)}{2}\sum\limits_{\substack{\ell =1}}^{M-2} (i+j-2-2\ell) M^{(1,2)}_{i+j-4}\\
D^{JJM}_{ij} &= (i-j)\sum\limits_{\substack{\ell =1}}^{M-2} \sum\limits_{\substack{ n=\ell+1 }}^{M-2}  J^{(b)}_{i+j-3-n} J^{(b)}_{n-1-\ell} M^{(1,2)}_\ell\\
\tilde{D}^{JJM,1}_{ij} &= \sgn(j-i)\sum\limits_{\substack{n = \mu-1}}^{M-2} \sum\limits_{\substack{\ell \geq \mu-1 }} (\mu-1-\ell) J^{(b)}_{i+j-3-n} J^{(b)}_{n-1-\ell} M^{(1,2)}_\ell\\
\tilde{D}^{JJM,2}_{ij} &= \sgn(j-i) \sum\limits_{\substack{n \geq M-1 }} \sum\limits_{\substack{\ell = \mu-1 }}^{M-2} (\mu-1-\ell) J^{(b)}_{i+j-3-n} J^{(b)}_{n-1-\ell} M^{(1,2)}_\ell\\
D^{bJM,1}_{ij} &= - b(i-j) \sum\limits_{\substack{n \geq M-1 }} (n-1) J^{(b)}_{i+j-3-n} M^{(1,2)}_{n-1}\\
D^{bJM,2}_{ij} &= b(i-j) \sum\limits_{\substack{\ell =1}}^{M-2} (i+j-3-\ell) J^{(b)}_{i+j-4-\ell-n} M^{(1,2)}_\ell
\end{aligned}
\end{equation}
for $D_{ij}-D_{ji}$. As for $E_{ij}-E_{ji}$,
\begin{equation}
E_{ij} - E_{ji} = (i-j)\sum\limits_{\substack{ n \geq M-1 }} J^{(b)}_{i+j-3-n} M^{(1,3)}_{n} + E^{JJM}_{ij} + \tilde{E}^{JJM,1}_{ij} + \tilde{E}^{JJM,2}_{ij} + E^{bJM}_{ij},
\end{equation}
with
\begin{equation}
\begin{aligned}
E^{JJM}_{ij} &= -(i-j) \sum\limits_{\substack{ n \geq M-1 }} \sum\limits_{ \ell \geq n-1} J^{(b)}_{i+j-3-n} J^{(b)}_{n-1-\ell} M^{(1,2)}_\ell\\
\tilde{E}^{JJM,1}_{ij} &= \sgn(j-i) \sum\limits_{\ell=1}^{\mu-2} \sum\limits_{\substack{n = \mu-1 }}^{M-2} (\mu-1-\ell) J^{(b)}_{i+j-3-n} J^{(b)}_{n-1-\ell} M^{(1,2)}_\ell\\
\tilde{E}^{JJM,2}_{ij} &=  \sgn(j-i) \sum\limits_{\ell=\mu-1}^{M-3} \sum\limits_{\substack{n=\ell+1}}^{M-2} (\mu-1-\ell)   J^{(b)}_{i+j-3-n} J^{(b)}_{n-1-\ell} M^{(1,2)}_\ell\\
E^{bJM}_{ij} &= -b(i-j)\sum\limits_{\substack{ n \geq M-1 }} (n-1) J^{(b)}_{i+j-3-n} M^{(1,2)}_{n-1}
\end{aligned}
\end{equation}
Finally,
\begin{multline}
F_{ij}-F_{ji} = \underbrace{b(i-j)\sum\limits_{n \geq M-1} (2n-i-j+1)J^{(b)}_{i+j-3-n} M^{(1,2)}_{n-1}}_{F^{bJM}_{ij}} \\
+ \underbrace{b\sgn(i-j)\sum\limits_{n = \mu-1}^{M-2} (\mu-n)(\mu-1-n)J^{(b)}_{i+j-3-n} M^{(1,2)}_{n-1}}_{\tilde{F}^{bJM}_{ij}}.
\end{multline}

\subsubsection{Terms $G_{ij}$ and $H_{ij}$} The term $H_{ij} = b^2(j-1)(i-1) [M^{(1,2)}_{i-1},M^{(1,2)}_{j-1}]$ is direct to compute and gives
\begin{equation}
H_{ij} = b^2(i-1)(j-1)(i-j) M^{(1,2)}_{i+j-4}.
\end{equation}

The term $G_{ij} = b(j-1)\sum\limits_{\substack{\ell \geq 1 }} [ J^{(b)}_{i-1-\ell}M^{(1,2)}_\ell,M^{(1,2)}_{j-1}]$ gives
\begin{multline}
    G_{ij}-G_{ji} = \underbrace{b^2\left[(i-1)(i-2)(i-3)-(j-1)(j-2)(j-3)\right] M^{(1,2)}_{i+j-4}}_{G^{M}_{ij}} \\
    + \underbrace{b(i-j)\sum\limits_{n=1}^{i+j-4} n J^{(b)}_{i+j-4-n} M^{(1,2)}_{n}}_{G^{bJM,1}_{ij}}
    + \underbrace{2b(i-j)\sum\limits_{n\geq M-2} (i+j-3)J^{(b)}_{i+j-4-n}M^{(1,2)}_n}_{G^{bJM,2}_{ij}} 
    + \tilde{G}^{bJM}_{ij}
\end{multline}
with 
\begin{multline}
\tilde{G}^{bJM}_{ij} = b\sgn(i-j)\sum\limits_{n=\mu-2}^{M-3} \left((\mu-1)(n-2\mu+4)-(M-1)(\mu-1-n)\right)J^{(b)}_{M+\mu-4-n}M^{(1,2)}_n \\
+ b\sgn(i-j)\left((M-1) J^{(b)}_{M-2} M^{(1,2)}_{\mu-2} - (\mu-1) J^{(b)}_{\mu-2} M^{(1,2)}_{M-2}\right)
\end{multline}

\subsubsection{Repackaging contributions}
We start by adding together $\tilde{C}_{ij}^{JJM}$, $\tilde{D}_{ij}^{JJM,1}$ and $\tilde{E}_{ij}^{JJM,1}$ 
\begin{equation}
\begin{aligned}
\tilde{C}_{ij}^{JJM} + \tilde{D}_{ij}^{JJM,1} + \tilde{E}_{ij}^{JJM,1} &= \sgn(j-i) \sum\limits_{n=\mu-1}^{M-2} \sum\limits_{\ell \geq 1} (3\mu-3-2n) J^{(b)}_{i+j-3-n} J^{(b)}_{n-1-\ell} M^{(1,2)}_\ell \\
&\begin{multlined} = \sgn(j-i) \sum\limits_{n=\mu-1}^{M-2}(3\mu-3-2n) J^{(b)}_{i+j-3-n} M^{(1,3)}_n \\
 \underbrace{-b\sgn(j-i) \sum\limits_{n=\mu-1}^{M-2}(3\mu-3-2n) (n-1) J^{(b)}_{i+j-3-n} M^{(1,2)}_{n-1}}_{{\widetilde{(CE)}}^{bJM}_{ij}}.\end{multlined}
\end{aligned}
\end{equation}

Now we add together the contributions $\tilde{D}_{ij}^{JJM,2}$ and $\tilde{E}_{ij}^{JJM,2}$. They have similar range except for the $\ell = M-2$ contribution of  $\tilde{D}_{ij}^{JJM,2}$ which is
\begin{equation}
    \lambda_{ij} \coloneqq \sgn(j-i)\sum\limits_{n\geq M} (\mu-M+1)J^{(b)}_{M+\mu-3-n} J^{(b)}_{n+1-M} M^{(1,2)}_{M-2}.
\end{equation}
Now the sum over $n$ in $\tilde{D}_{ij}^{JJM,2}+\tilde{E}_{ij}^{JJM,2}$ can be encoded into some $M^{(1,2)}$ as follows
\begin{align}
\tilde{D}_{ij}^{JJM,2}+\tilde{E}_{ij}^{JJM,2} &=  \lambda_{ij} + \sgn(j-i)(1+b)\sum\limits_{\ell=\mu}^{M-3} (\mu-1-\ell) M^{(1,2)}_{M+\mu-2-\ell} M^{(1,2)}_{\ell} \\
&-b(1+b)\sgn(j-i)\sum\limits_{\ell=\mu-1}^{M-3} (\mu-1-\ell)(M+\mu-3-\ell) M^{(1,1)}_{M+\mu-3-\ell}M^{(1,2)}_{\ell}. \nonumber
\end{align}
Commuting the two $M^{(1,2)}$s and relabelling the sum as $\ell \rightarrow M+\mu-2-\ell$ then taking the half-sum gives
\begin{multline}
\sgn(j-i) \sum\limits_{\ell=\mu}^{M-3} (\mu-1-\ell) M^{(1,2)}_{M+\mu-2-\ell} M^{(1,2)}_{\ell} = \frac{i-j}{2} \sum_{\ell=\mu+1}^{M-3} M^{(1,2)}_{i+j-2-\ell} M^{(1,2)}_\ell\\
+ \frac{\sgn(j-i)}{2} \left(\sum_{\ell=\mu+1}^{M-2} (\ell+1-M)(2\ell +2-i-j) M^{(1,2)}_{i+j-4} - M^{(1,2)}_\mu M^{(1,2)}_{M-2} -  M^{(1,2)}_{M-2} M^{(1,2)}_\mu\right)
%&= \sum\limits_{\ell=\mu+1}^{M-2} (M-1-\ell)(2\ell-M-\mu+2)A^2_{M+\mu-4} \\
%&+\sum\limits_{\ell=\mu+1}^{M-2} (M-1-\ell) A^2_{M+\mu-2-\ell}A^2_{\ell} \nonumber.
\end{multline}
Similarly, the sum over $n$ in $\lambda_{ij}$ can also be written in terms of some $M^{(1,2)}$,
\begin{multline}
\lambda_{ij} = \frac{1}{2} (i-j+\sgn(j-i)) (M^{(1,2)}_\mu M^{(1,2)}_{M-2} + M^{(1,2)}_{M-2} M^{(1,2)}_\mu + (\mu-M+2) M^{(1,2)}_{i+j-4}) \\
-b(1+b)\sgn(j-i)(M-\mu-1)(\mu-1) M^{(1,1)}_{\mu-1} M^{(1,2)}_{M-2}
%(1+b)\sgn(j-i)(M-\mu-1)A^2_{\mu}A^2_{M-2} \\
%\underbrace{-b(1+b)\sgn(j-i)(M-\mu-1)(\mu-1)A^1_{\mu-1}A^2_{M-2}}_{\theta_{ij}^{(1)}}. \nonumber
\end{multline}
%Commuting the $A^2A^2$ term in $\lambda$ gives
%\begin{align}
%\lambda_{ij} &= (1+b)\sgn(j-i)(M-\mu-1)A^2_{M-2}A^2_{\mu} \\
%&+ (1+b)\sgn(j-i)(M-\mu-1)(\mu-M+2)A^2_{M+\mu-4} + \theta_{ij}^{(1)}. \nonumber
%\end{align}
%
%Taking one-half of both expressions, and adding the terms together, we get
Altogether, one gets
\begin{multline}
\tilde{D}_{ij}^{JJM,2}+\tilde{E}_{ij}^{JJM,2} = \underbrace{\frac{1+b}{2}(i-j) \sum\limits_{\ell=\mu}^{M-2} M^{(1,2)}_{i+j-2-\ell} M^{(1,2)}_{\ell}}_{(DE)_{ij}^{MM}}\\
+ \underbrace{(1+b)\frac{\sgn(j-i)}{2}\sum\limits_{\ell=\mu}^{M-2} (\ell-M+1)(2\ell-i-j+2) M^{(1,2)}_{i+j-4}}_{(DE)^M_{ij}}\\
+ \underbrace{ b(1+b)\sgn(i-j)\sum\limits_{\ell=\mu-1}^{M-2} (\mu-1-\ell)(i+j-3-\ell) M^{(1,1)}_{i+j-3-\ell}M^{(1,2)}_{\ell}}_{\widetilde{(DE)}_{ij}^{bJM}}
\end{multline}

Next, we add together the terms $C^{JJM,1}_{ij}$ and $E^{JJM}_{ij}$. The same treatment as for $C^{JJA,1}_{ij} + E^{JJA}_{ij}$ leads to
% TEMPLATE FOR ONCE THE CALCULATION IS DONE
\begin{multline}
C^{JJM,1}_{ij} + E^{JJM}_{ij} = (i-j)\sum\limits_{n=M-1}^{i+j-3}  J^{(b)}_{i+j-3-n} M^{(1,3)}_{n} + (CE)^{JJM,1}_{ij} \\ + (CE)^{JJM,2}_{ij} + (CE)^M_{ij} + (CE)^{bJM}_{ij}
\end{multline}
with
\begin{equation}
\begin{aligned}
(CE)^{JJM,1}_{ij} &= (i-j)\sum\limits_{n=M-1}^{i+j-3} \sum\limits_{\ell = n-\mu+1}^{M-2}  J^{(b)}_{i+j-3-n} J^{(b)}_{n-1-\ell} M^{(1,2)}_{\ell}\\
(CE)^{JJM,2}_{ij} &= -(i-j) \sum\limits_{n=M-1}^{i+j-3} \sum\limits_{\ell = 1}^{M-2} J^{(b)}_{i+j-3-n} J^{(b)}_{n-1-\ell} M^{(1,2)}_{\ell}\\
(CE)^M_{ij} &= -(1+b)(i-j)\frac{(\mu-1)(\mu-2)}{2} M^{(1,2)}_{i+j-4}\\
(CE)^{bJM}_{ij} &= -b(i-j)\sum\limits_{n=M-1}^{i+j-3} (n-1)J^{(b)}_{i+j-3-n} M^{(1,2)}_{n-1}
\end{aligned}
\end{equation}
% VICTOR'S CALCULATION
%\begin{align}
%C^{JJA,1}_{ij} + E^{JJA}_{ij} &= (i-j)\sum\limits_{n=M-1}^{M+\mu-4}  J^{(b)}_{M+\mu-3-n} M^{(1,3)}_{n} \underbrace{-(i-j)\sum\limits_{\substack{n=M-1 \\ n \neq i+j-3}}^{M+\mu-4} \sum\limits_{\substack{\ell = n-\mu+1 \\ \ell \neq n-1 }}^{M-3}  J^{(b)}_{M+\mu-3-n} J^{(b)}_{n-1-\ell} M^{(1,3)}_{\ell} }_{(CE)_{ij}^{JJA,1}}\\
%& \underbrace{-(i-j)\sum\limits_{n=M-1}^{M+\mu-4} b(n-1)J^{(b)}_{i+j-3-n} M^{(1,3)}_{n-1}}_{(CE)_{ij}^{bJA}} \underbrace{-(1+b)(i-j)\frac{\mu(\mu-1)}{2} M^{(1,2)}_{i+j-4}}_{(CE)_{ij}^{bA}} \nonumber\\
%&+ \underbrace{(i-j) \sum\limits_{n=M-1}^{M+\mu-4} \sum\limits_{\substack{\ell = 1 \\ \ell \neq n-1}}^{M-1} J^{(b)}_{M+\mu-3-n} J^{(b)}_{n-1-\ell} M^{(1,3)}_{\ell}}_{(CE)_{ij}^{JJA,2}}. \nonumber
%\end{align}

Finally, we repackage together
\begin{align}
G^{M}_{ij} + H_{ij} = b^2(i-j)(i+j-3)(i+j-4) M^{(1,2)}_{i+j-4}.
\end{align}
Moreover,
\begin{multline}
D^{bJM,1}_{ij} + E^{bJM}_{ij} + F^{bJM}_{ij} + G^{bJM,2}_{ij} = b(i-j)(i+j-3) \sum_{n\geq M-1} J^{(b)}_{i+j-3-n} M^{(1,2)}_{n-1} 
%\\+ \underbrace{b\sgn(j-i)(M-1)J^{(b)}_{M-2}A^2_{\mu-2}}_{\theta_{ij}^{(3)}}. 
\end{multline}
and
\begin{multline}
(CE)^{bJM}_{ij} + G^{bJM,1}_{ij} + D^{bJM,2}_{ij} = b(i-j)(i+j-3) \sum_{n=1}^{M-2} J^{(b)}_{i+j-3-n} M^{(1,2)}_{n-1} \\ + \underbrace{b(i-j)(\mu-1) J^{(b)}_{\mu-2} M^{(1,2)}_{M-2}}_{\theta_{ij}}. 
\end{multline}
Hence, adding together the last three equations gives a contribution
\begin{multline}
(D^{bJM,1}_{ij} + E^{bJM}_{ij} + F^{bJM}_{ij} + G^{bJM,2}_{ij}) + ((CE)^{bJM}_{ij} + G^{bJM,1}_{ij} + D^{bJM,2}_{ij}) + (H_{ij} + G^{M}_{ij}) \\
= b(i-j)(i+j-3) M^{(1,3)}_{i+j-3} + \theta_{ij}.
\end{multline}

\subsubsection{Cancellations}
We notice the following cancellations,
\begin{equation}
\begin{aligned}
C^{M,1}_{ij} + B^{M}_{ij} &= 0, \\
D^{MM,1}_{ij} + B^{MM}_{ij} &= 0, \\
D^{MM,2}_{ij} + (DE)^{MM}_{ij} &=0, \\
C^{M,2}_{ij} + D^M_{ij} + (DE)^M_{ij} + (CE)^M_{ij}  &= 0, \\
C^{JJM,2}_{ij} + D^{JJM,1}_{ij} + (CE)^{JJM,1}_{ij} + (CE)^{JJM,2}_{ij} &= 0, \\
\widetilde{(CE)}^{bJM}_{ij} + \widetilde{(DE)}^{bJM}_{ij} + \tilde{G}^{bJM}_{ij} + \tilde{F}^{bJM}_{ij} + \theta_{ij} &= 0.%,\\
%\theta_{ij}^{(1)} + \theta_{ij}^{(2)} + \theta_{ij}^{(3)} + G^{bJA,3}_{ij} &= 0.
\end{aligned}
\end{equation}

\subsubsection{Final expression} Bringing together the remaining terms, we get the final expression for the commutator
\begin{align}
\label{eq:comm_triv_bip}
[M^{(1,3)}_i, M^{(1,3)}_j] &= 2(i-j) \sum\limits_{n \geq M-1} J^{(b)}_{i+j-3-n} M^{(1,3)}_n + (i-j) \sum\limits_{n=M-1}^{i+j-3} J^{(b)}_{i+j-3-n} M^{(1,3)}_n \\
&+ \sgn(j-i)\sum\limits_{n=\mu-1}^{M-2} (3\mu-2n-3) J^{(b)}_{i+j-3-n} M^{(1,3)}_n 
+b(i-j)(i+j-3) M^{(1,3)}_{i+j-3} \\
&= \tilde{D}_{ij,\ell}(3) M^{(1,3)}_\ell \nonumber
\end{align}
which concludes the computation and the proof of Proposition~\ref{thm:SimplifiedCommutators2}.

%% file: Appendices/App_UN2_OD.tex
\chapter{Double scaling limit for the quartic \texorpdfstring{$U(N)^2\times O(D)$}{U(N)2XO(D)} model}
\label{app:DS_UN2OD} 

%%%%%%%%%%%%%%%%%%%
\section{Definition of the model and its large \texorpdfstring{$N$}{N}, large \texorpdfstring{$D$}{D} expansion}

%%%%%%%%%%%%%%%%%%%
\subsection{Feynman graphs, genus and grade}

The $U(N)^2 \times O(D)$ multi-matrix model is a model involving a vector of $D$ complex matrices of size $N \times N$, denoted $(X_\mu)_{\mu=1, \dotsc, D} = (X_1, \dotsc, X_D)$. The model is required to be invariant under unitary actions on the left and on the right of each matrix $X_\mu$,
\begin{equation}
X_\mu \rightarrow X'_\mu = U_1 X_{\mu} U_2^\dagger,
\end{equation}
with $U_1,U_2 \in U(N)$. The ring of polynomials which are invariant under this transformation is generated by products of traces of the form
\begin{equation} \label{U(N)2Invariance}
\Tr X_{\mu_1} X^\dagger_{\nu_1} \dotsm X_{\mu_n} X^\dagger_{\nu_n}.
\end{equation}
The model is further required to be invariant orthogonal transformations on the vector $(X_\mu)_{\mu=1, \dotsc, D}$,
\begin{equation}
X_\mu \rightarrow X'_\mu = \sum_{\mu'=1}^D O_{\mu\mu'}\ X_{\mu'},
\end{equation}
for any $O \in O(D)$. To enforce this on polynomials which are products of the traces of the type \eqref{U(N)2Invariance}, each vector index $\mu_i$ and $\nu_i$ must be identified with another vector index as follows,
\begin{equation} \label{O(D)Invariance}
\sum_{\mu=1}^D X_\mu \dotsb X_{\mu}\dotsb \quad \text{or} \quad \sum_{\mu=1}^D X_\mu \dotsb X^\dagger_{\mu}\dotsb \quad \text{or} \quad \sum_{\mu=1}^D X^\dagger_\mu \dotsb X^\dagger_{\mu}\dotsb
\end{equation}
%which transforms under $U(N)^2\times O(D)$~\cite{BeCa} as

Graphically, the bubbles of the model are represented as follows 
\begin{itemize}
\item We represent $(X_\mu)_{ab}$ as a white vertex, $(X^\dagger_\nu)_{b'a'}$ as a black vertex.
\item Their matrix indices $a, a'$ are half-edges of color 1, and their matrix indices $b, b'$ are half-edges of color 2, and the vector indices $\mu, \nu$ are half-edges of color 3. 
\item One then forms graphs by connecting half-edges of the same color, with the only constraint being that the subgraph obtained by removing all edges of color 3 is bipartite. 
\end{itemize}

%The product of traces like in~{\eqref{U(N)2Invariance} implies graphically that the subgraph with colors 1 and 2 only is a disjoint union of bipartite cycles whose edges alternate the colors 1 and 2. The identification of vector indices from \eqref{O(D)Invariance} then means that the edges of color 3 form a perfect matching on the vertices (black and white undifferently).
The quadratic bubble is
\begin{equation}
I_k(X, X^\dagger) = \sum_{\mu=1}^D \Tr X_\mu X^\dagger_\mu.
\end{equation}
We will only consider connected, quadratic and quartic interactions. The four connected, quartic interactions are
\begin{align}
I_t(X, X^\dagger) = \sum_{\mu, \nu} \Tr \bigl(X_\mu X^\dagger_\nu X_\mu X^\dagger_\nu\bigr) &= \begin{array}{c}\includegraphics[scale=.28]{TetrahedralComplex.pdf}\end{array}\\
I_{1}(X, X^\dagger) = \sum_{\mu, \nu} \Tr \bigl(X_\mu X^\dagger_\mu X_\nu X^\dagger_\nu\bigr) &= \begin{array}{c}\includegraphics[scale=.28]{PillowColor1.pdf}\end{array}\\
I_{2}(X, X^\dagger) = \sum_{\mu, \nu} \Tr \bigl(X_\mu X^\dagger_\nu X_\nu X^\dagger_\mu\bigr) &= \begin{array}{c}\includegraphics[scale=.28]{PillowColor2.pdf}\end{array}\\
I_{3b}(X, X^\dagger) = \sum_{\mu, \nu} \Tr \bigl(X_\mu X^\dagger_\nu\bigr)\ \Tr \bigl(X_\nu X^\dagger_\mu\bigr) &= \begin{array}{c}\includegraphics[scale=.28]{PillowColor3.pdf}\end{array}\\
I_{3nb}(X, X^\dagger) = \sum_{\mu, \nu} \Tr \bigl(X_\mu X^\dagger_\nu\bigr)\ \Tr \bigl(X_\mu X^\dagger_\nu\bigr) &= \begin{array}{c}\includegraphics[scale=.28]{PillowColor3NonBip.pdf}\end{array}
\end{align}

The action of the model writes 
\begin{multline}
S_{U(N)^2\times O(D)}(X_{\mu},X_{\mu}^{\dagger }) = -ND \sum_{\mu=1}^D \Tr(X_{\mu} X^\dagger_{\mu}) + N D^{3/2}\frac{\lambda_1}{2} I_t(X_{\mu},X_{\mu}^{\dagger }) \\
+ ND\frac{\lambda_2}{2} \bigl(I_{p;1}(X_{\mu},X_{\mu}^{\dagger })+I_{p;2}(X_{\mu},X_{\mu}^{\dagger })\bigr) +  D^2\frac{\lambda_2}{2} \bigl(I_{p;3b}(X_{\mu},X_{\mu}^{\dagger }) + I_{p;3nb}(X_{\mu},X_{\mu}^{\dagger })\bigr),
\end{multline}
and thus the partition function is
\begin{equation}
Z_{U(N)^2\times O(D)}(\lambda_1, \lambda_2) = e^{F_{U(N)^2\times O(D)}(\lambda_1, \lambda_2)} = \int \prod_{\mu=1}^D \prod_{a,b=1}^N d(X_\mu)_{ab} d(\overline{X_\mu})_{ab}\ e^{S_{N,D}(X_\mu, X^\dagger_\mu)}.
\end{equation}

Feynman graphs are obtained by taking any collection of interactions and performing Wick contractions. The latter can be represented graphically as edges pairing every white vertex to a black vertex. As usual in the literature, we give those edges the fictitious color 0. As a result the set of connected Feynman graphs, denoted $\bar{\mathbb{G}}$ is the set of connected, 4-regular, properly-edge-colored graphs such that the subgraph obtained by removing all edges of color 0 is a disjoint union of quartic bubbles $I_{p;1}, I_{p;2}, I_{p;3b}, I_{p;3nb}, I_t$.

\medskip

Let $\bar{\cG}\in\bar{\mathbb{G}}$, then denote
\begin{itemize}
\item $n_t, n_{1}, n_{2}, n_{3b}, n_{3nb}$ the number of interactions $I_t, I_{1}, I_{2}, I_{3b}, I_{3nb}$ respectively,
\item $E$ the number of edges of color 0 (with $E = 2(n_t+n_{1}+n_{2}+n_{3b}+n_{3nb})$),
\item $F_{0a}$ for $a=1,2,3$, the number of bicolored cycles which alternate the colors $0, a$. For $a=1,2$, we will also call them faces of colors 1 and 2 (the reason for this will be clear below).
\end{itemize}

The free energy has $1/N, 1/D$ expansions \cite{FeVa}. The $1/N$ expansion is governed by $h(\bar{\cG})$ and the $1/D$ expansion is governed by a non-negative integer, the \emph{grade} denoted $l(\bar{\cG})$.

\begin{theorem}\label{thm:free_energy}
The free energy expands as
\begin{equation}
F_{U(N)^2\times O(D)}(\lambda_1, \lambda_2) = \sum_{\bar{\cG}\in\bar{\mathbb{G}}} \biggl(\frac{N}{\sqrt{D}}\biggr)^{2-2h(\bar{\cG})} D^{2-l(\bar{\cG})/2}\ \mathcal{A}_{\bar{\cG}}(\lambda_1, \lambda_2)
\end{equation}
where $h(\bar{\cG})$ is a non-negative integer which we call the genus and $l(\bar{\cG})$ is a non-negative integer called the \emph{grade}.
\end{theorem}
Notice that in fact, the large $N$, large $D$ expansion is an expansion in $D$ and $L := \frac{N}{\sqrt{D}}$.

\begin{proof}
Here we give an alternative proof to \cite{FeVa}. From the Feynman expansion, one has
\begin{align} \label{GenusExpansion}
2-2h(\bar{\cG}) &= F_{01} + F_{02} - E + n_t + n_{1} + n_{2}\\
1+h(\bar{\cG}) - \frac{l(\bar{\cG})}{2} &= F_{03} - E +\frac{3}{2}n_t + n_{1} + n_{2} + 2(n_{3b}+n_{3nb}), \label{Grade}
\end{align}
and we want to show that both are non-negative.

\medskip

Since we are primarily dealing with matrices, one expects the Feynman graphs to involve ribbon graphs, and the $1/N$ expansion to involve their genera, and thus be related to $h(\bar{\cG})$ and $l(\bar{\cG})$. Let us explain how ribbon graphs are encoded in our colored graphs. From $\bar{\cG}\in \bar{\mathbb{G}}$, remove all edges of color 3. Since they correspond to identifications of vector indices, it means that we are only keeping the information about the matrix part of the model. The resulting graph is a not-necessarily connected, 3-regular, properly-edge-colored graph, with colors 0, 1, 2.

\medskip

Such graphs are known to be equivalent to ribbon graphs (more precisely of models for a complex matrix). One way to see this is to choose a cyclic order for the three colors meeting at every white vertex and the reverse cyclic order at the black vertices. For instance, one draws the edges of colors $(0,1,2)$ in the counter-clockwise order around the white vertices and in the clockwise order around black vertices. This provides each cycle of colors $\{1,2\}$ with a cyclic order of its incident edges of color 0. It can thus be replaced with a ribbon vertex,
\begin{equation}
\includegraphics[valign=c,scale=.5]{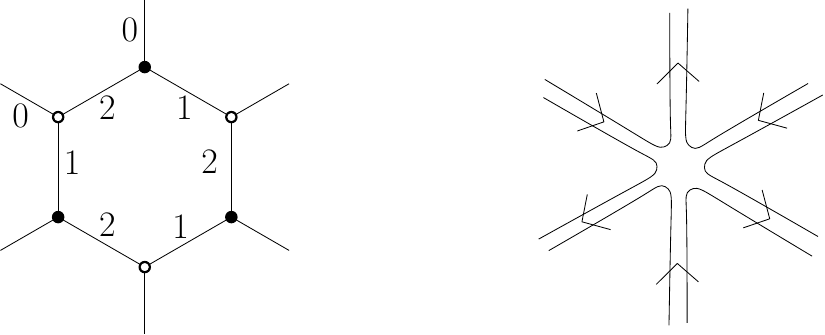}
\end{equation}
(of course in our case all cycles of colors $\{1,2\}$ are of length 4, but the correspondence is more general). Only edges of color 0 remain and they can be thickened (without twist) to obtain a ribbon graph denoted $\cG_{012}$. Arrows can be used to recover the vertex coloring of $\bar{\cG}$: for instance an edge of color 0 in $\bar{\cG}$ is oriented from black to white and those orientations are inherited in $\cG_{012}$. 

\medskip

An alternative representation is as a combinatorial map, where one collapses the ribbon but keep the cyclic ordering of the edges around each vertex. This is enough to reproduce the ribbon graph. An example of a graph $\bar{\cG} \in \mathbb{G}_{U(N)^2\times O(D)}$ and its associated combinatorial map is given on Figure~\ref{fig:un2xod_to_map}. All vertices have degree four (because all bubbles are quartic).
\begin{figure}
\centering
\includegraphics[scale=0.63]{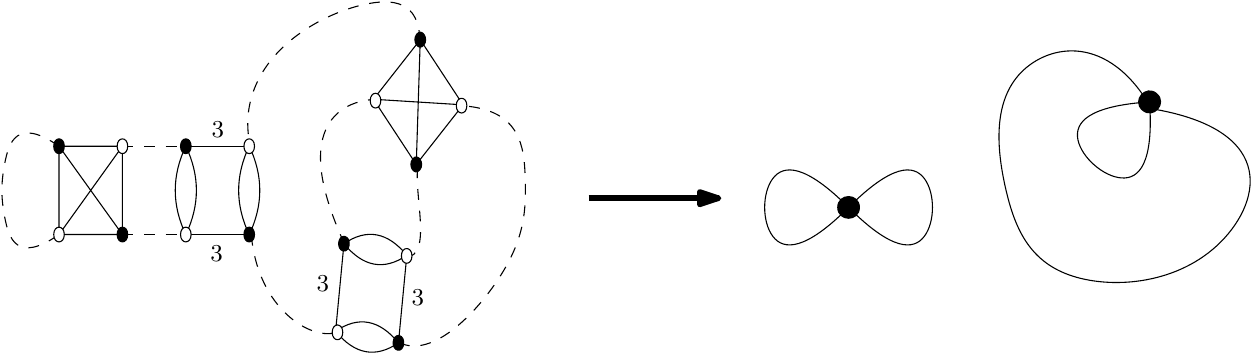}
\caption{A graph of the $U(N)^2\times O(D)$ model and its ribbon graph.}
\label{fig:un2xod_to_map}
\end{figure}

\medskip

In addition to vertices and edges, a ribbon graph has faces, obtained by following the border of the ribbon. It is well-known that Euler's relation applies in this context: if $F, E, V$ are respectively the number of faces, edges and vertices of a connected ribbon graph, then
\begin{equation}
F-E+V=2-2g
\end{equation}
where $g$ is a non-negative integer called the genus of the ribbon graph.

\medskip

In order to apply this to $\cG_{012}$, we first need to identify the faces and vertices, and to discuss connectedness. The faces of $\cG_{012}$ are the faces of colors 1 and 2 of $\bar{\cG}$, as seen from the bijection described above. Moreover,
\begin{equation} 
V(\cG_{012}) = n_t+n_{1}+n_{2}+2(n_{3b}+n_{3nb}).
\end{equation}
Note that $\cG_{012}$ may not be connected if $G$ contains some interactions $I_{p;3b}$ and $I_{p;3nb}$ because the latter are disconnected by the removal of their edges of color 3 (corresponding to the fact that they are double-trace interactions). The genus and number of connected components of this ribbon graph are denoted $g_{012}$ and $c_{012}$. Here the genus is the sum of the genera of its connected components. Euler's formula thus gives
\begin{equation}
2c_{012} - 2g_{012} = F_{01} +F_{02} - E_0 + n_t + n_{1} + n_{2} + 2(n_{3b}+n_{3nb}).
\end{equation}

As a result, Equation~\eqref{GenusExpansion} for $h(\bar{\cG})$ can be expressed in terms of the genus of the ribbon graph $\cG_{012}$ as
\begin{equation} \label{GenusExpansion2}
h(\bar{\cG}) = g_{012} + n_{3b}+n_{3nb}+1-c_{012}.
\end{equation}
Crucially, $g_{012}\geq 0$ and $n_{3b}+n_{3nb}+1-c_{012}\geq 0$ which is easily proved by induction for example, so that $h(\cG)\geq 0$. 

\medskip

In the model with tetrahedral bubble only, $n_{3b}=n_{3nb}=0$, then $c_{012}=1$ and $h(\bar{\mathcal{G}})$ reduces to the genus of the ribbon graph $g_{012}$ (and more generally, whenever only single trace interaction are considered, as often in the literature \cite{BeCa}). In the rest of this section, we will use $h$ when working with Feynman graph while $g$ will be used only when working with maps.
Another ribbon graph can be obtained by following the same procedure as the one leading to $\cG_{012}$, but starting with removing the edges of color 1 instead of color 3. The corresponding ribbon graph is denoted $\cG_{023}$. Similarly one obtains $\cG_{013}$ when one starts by removing the edges of color 2. The only difference with the procedure leading to $\cG_{012}$ is that now edge-twists have to be allowed and the genera may be half-integers. Euler's formulas are
\begin{equation}
\begin{aligned}
2c_{013} - 2g_{013} &= F_{01} + F_{03} - E + n_t + n_1 + 2n_2 + n_{3b} + n_{3nb}\\
2c_{023} - 2g_{023} &= F_{02} + F_{03} - E + n_t + 2n_1 + n_2 + n_{3b} + n_{3nb}.
\end{aligned}
\end{equation}
Equation \eqref{Grade} for $l(\bar{\cG})$ can thus be written as
\begin{equation}
\frac{l(\bar{\cG})}{2} = g_{023} + g_{013} + (n_{1}+1-c_{023}) + (n_{2}+1-c_{013}).
\end{equation}
All the quantities into brackets are non-negative integers, and the genera may be half-integers.
\end{proof}

For reference, Equations~\eqref{GenusExpansion} and~\eqref{Grade} give the following combinatorial expression of the grade
\begin{equation}
\label{eq:grade_UN2_comb}
2-\frac{l(\bar{\cG})}{2} = \frac{F_1+F_2}{2} + F_3 - \frac{3}{2}E + 2n_t + 2(n_{3b}+n_{3nb}) + \frac{3}{2}(n_1+n_2).
\end{equation}

\paragraph{2-point graphs.\\}

Due to the $U(N)^2\times O(D)$ invariance, the 2-point function is
\begin{equation}
\langle (X_\mu)_{ab} (\overline{X}_\nu)_{cd}\rangle = \frac{1}{N^2 D} G_{N,D}(\lambda_1, \lambda_2) \delta_{\mu, \nu} \delta_{ac} \delta_{bd},
\end{equation}
with $G_{N, D}(\lambda_1, \lambda_2) = \langle \sum_{\mu=1}^D \Tr X_\mu X^\dagger_\mu\rangle$. It has an expansion on 2-point graphs, whose set is denoted $\mathbb{G}$. A 2-point graph is like a vacuum graph with one white and one black vertex having exactly one half-edge of color 0 incident on them. 

\medskip

From a graph $\cG\in\mathbb{G}$ one can obtain a unique vacuum graph $\bar{\cG}\in \bar{\mathbb{G}}$ called its \emph{closure}, by connecting the two half-edges of color 0. This vacuum graph is furthermore equipped with a marked edge (that obtained by connecting the two half-edges) called a \emph{root}. The set of rooted graphs is the set of Feynman graphs for the expansion of $G_{N,D}(\lambda_1, \lambda_2)$.

\medskip

The amplitudes of $\cG$ and its closure $\bar{\cG}$ are the same, up to a factor $N^2 D$, since the difference between those graphs is one face of color 1, one face of color 2 (both contributing to a power of $N$), and one bicolored cycle of colors $\{0,3\}$ (contributing to $D$). We therefore simply define the genus and the grade of a 2-point graph to be those of its closure, i.e. $h(\cG):= h(\bar{\cG})$ and $l(\cG) := l(\bar{\cG})$. 

%%%%%%%%%%%%%%%%%%
\subsection{Melons, dipoles and chains}

As with the quartic $O(N)^3$ tensor model studied in the previous section, the scheme decomposition of the $U(N)^2 \times O(D) $ relies on the same families of subgraph which we define here.

\paragraph{Melons\\}

The elementary melons of our model are the following 2-point graphs:
\begin{equation}
\includegraphics[scale=.5, valign=c]{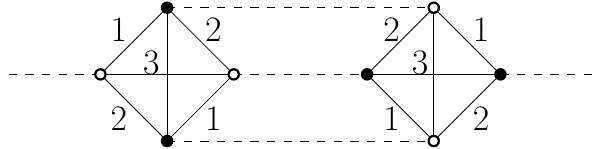}, \quad \includegraphics[scale=.5, valign=c]{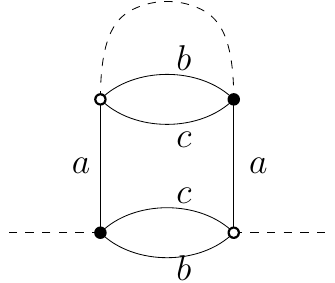}, \quad \includegraphics[scale=.5, valign=c]{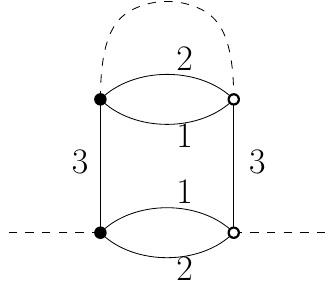}
\end{equation}
where $(a,b,c)$ is a permutation of $(1,2,3)$. Melonic graphs are obtained iteratively. First consider the closure of the elementary melons. Then on any edge of color 0, cut and insert any one of the elementary melons (still so that edges of color 0 connect white to black vertices), and repeat.

\begin{proposition}
Inserting a melon on an edge of a graph $\bar{\cG}$ leaves $h(\bar{\cG})$ and $l(\bar{\cG})$ invariant. Melonic graphs are the only graphs of vanishing genus and grade. This is also true for 2-point graphs.
\end{proposition}

From their recursive structure, it can be seen that the generating function of melonic 2-point graphs thus satisfies the equation
\begin{equation}
M(\lambda_1,\lambda_2) = 1 + \lambda_1^2M(\lambda_1,\lambda_2)^4 + 4\lambda_2M(\lambda_1,\lambda_2)^2
\label{eq:mel_un2}
\end{equation}
Performing the change of variable $(t,\mu) = \left(\lambda_1^2,\frac{4\lambda_2}{\lambda_1^2}\right)$ gives the same function $M(t,\mu)$ as in equation~\eqref{eq:mel} of the $O(N)^3$-invariant tensor model.

\vspace{5pt}
\paragraph{Dipoles\\}

Dipoles are defined similarly as for the $O(N)^3$ tensor model, see Def.~\ref{def:dip}. The set of dipoles of color $1$ and $2$ is similar to the $O(N)^3$ tensor model while adding bipartite coloring of their verties
\begin{equation}
D_{1,2} = \left\{\includegraphics[scale=0.22,valign=c]{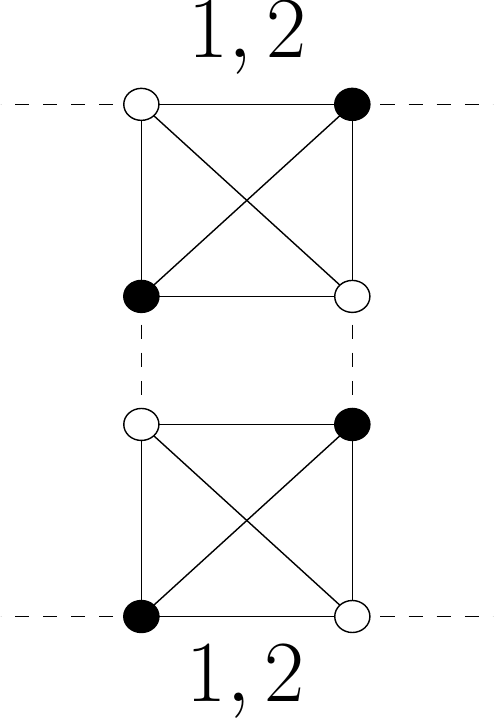},\hspace{.5cm} \includegraphics[scale=0.35,valign=c]{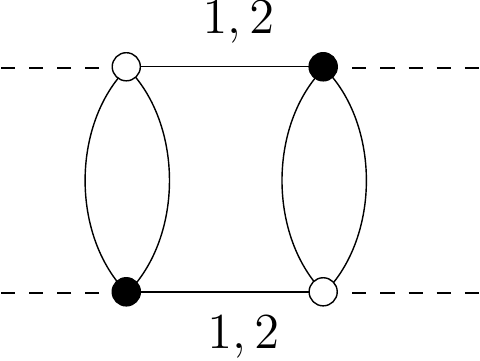}\right\}
\label{eq:d12}
\end{equation}

But we now have two types of dipoles of color $3$ according to the coloring of their vertices.
\begin{itemize}
\item One involving the tetrahedral interaction, and the non-bipartite pillow of color $3$, %denoted $D_3$
\begin{equation}
D_3 = \left\{\includegraphics[scale=0.25,valign=c]{d3_m.pdf}, \hspace{.5cm} \includegraphics[scale=0.26,valign=c]{d3_pnb.pdf}\right\}
\end{equation}
\item One for the bipartite pillow of color $3$, %denoted $D'_3$
\begin{equation}
D_3' = \left\{\includegraphics[scale=0.3,valign=c]{d3_pb.pdf}\right\}
\end{equation}
\end{itemize}

Again, dipoles may be non-isolated. When they are isolated, it is useful to represent dipoles from the same group as a \emph{dipole-vertex}. In terms of generating series, the series associated to a dipole-vertex is the sum of the series of the dipoles in that group. We represent a dipole-vertex as a box with each pair on either sides. To distinguish the sides, we separate them with some thickened edges along the box, so that
\begin{equation}
    \begin{aligned}
    \includegraphics[scale=0.37,valign=c]{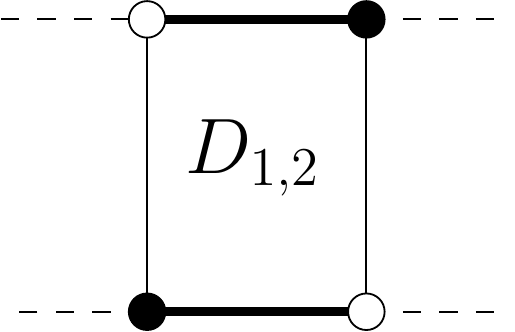} &= \includegraphics[scale=0.22,valign=c]{d12_m.pdf} + \includegraphics[scale=0.35,valign=c]{d12_pb.pdf}\\
    \includegraphics[scale=0.3,valign=c]{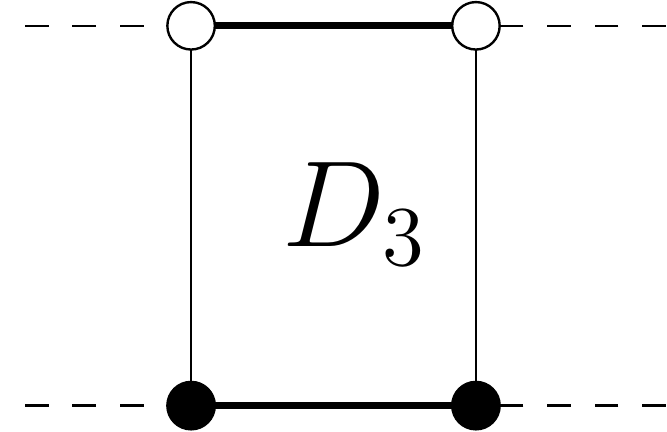} &= 
    \includegraphics[scale=0.25,valign=c]{d3_m.pdf} +  \includegraphics[scale=0.26,valign=c]{d3_pnb.pdf}\\
    \includegraphics[scale=0.45,valign=c]{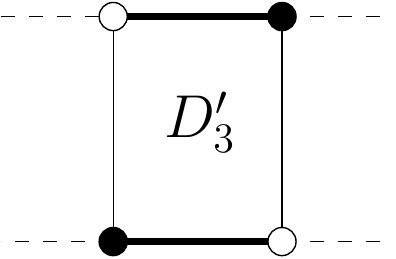} &= \includegraphics[scale=0.3,valign=c]{d3_pb.pdf}
    \end{aligned}
\end{equation}
Let us denote $U$ the generating function for the dipoles $D_{1,2,3}$ decorated with melons on one side and on the internal edges, and $V$ the generating function for $D'_3$ decorated with melons on one side. Then we have
\begin{align}
U(t,\mu) &= tM(t,\mu)^4 + \frac{1}{4}t\mu M(t,\mu)^2 \underset{\eqref{eq:mel_un2}}{=} M(t,\mu)-1-\frac{3}{4}t\mu M(t,\mu)^2 \\
V(t,\mu) &= \frac{1}{4}t\mu M(t,\mu)^2 
\end{align}

\paragraph{Chains\\}

At a purely combinatorial level, chains follow definition~\ref{def:chains} and share all the subsequent properties introduced for the $O(N)^3$ tensor model. The difference with the $O(N)^3$ tensor model lie in their generating functions. Since the model considered here allows for more dipoles of color $3$, the generating function of chains is also impacted by this change. Chains of color 1 and 2 are geometric series of dipoles of color 1 and 2. Chains of color 3 can have at each step either $D_3$ or $D'_3$. The corresponding generating functions are
\begin{equation}
C_1 = C_2 = U \sum_{n\geq0} U^n = \frac{U}{1-U} \qquad
C_3 = (U+V) \sum\limits_{n\geq0} (U+V)^n = \frac{(U+V)}{1-(U+V)}
\label{eq:chain_un2}
\end{equation}
We can further distinguish two types of chains of color $3$, depending on the coloring on the vertex at their ends. We denote the respective generating series $C_{3,\hspace{1pt}\includegraphics[scale=0.2,valign=c]{wbbw.pdf}}$ and $C_{3,\hspace{1pt}\includegraphics[scale=0.2,valign=c]{wwbb.pdf}}$. To evaluate them, notice that $C_{3,\hspace{1pt}\includegraphics[scale=0.2,valign=c]{wbbw.pdf}}$ (resp. $C_{3,\hspace{1pt}\includegraphics[scale=0.2,valign=c]{wwbb.pdf}}$) is a series in $U$ and $V$ where each term has an even (resp. odd) number of $U$. Thus, we introduce a dummy variable $x$ which counts the number of $U$ in each term and then extract the even powers (resp. odd powers) of $x$ to get $C_{3,\hspace{1pt}\includegraphics[scale=0.2,valign=c]{wbbw.pdf}}$ (resp. $C_{3,\hspace{1pt}\includegraphics[scale=0.2,valign=c]{wwbb.pdf}}$~). We have
\begin{align}
C_{3,\hspace{1pt}\includegraphics[scale=0.2,valign=c]{wbbw.pdf}} &= \sum_{k\geq 0} [x^{2k}] \frac{xU+V}{1-(xU+V)} = \sum_{k\geq 0} [x^{2k}] \frac{1}{1-V}\frac{xU+V}{1-\frac{xU}{1-V}} \nonumber\\
&= \frac{1}{1-V}\sum_{k\geq 0} [x^{2k}] (xU+V)\sum\limits_{n\geq0}\left(\frac{xU}{1-V}\right)^{n} = \frac{U^2+V-V^2}{(1-V)^2-U^2} \nonumber\\
&= \frac{U^2+V-V^2}{(1-V-U)(1-V+U)}
\end{align}
And similarly
\begin{align}
C_{3,\hspace{1pt}\includegraphics[scale=0.2,valign=c]{wwbb.pdf}} &= \sum_{k\geq 0} [x^{2k+1}] \frac{xU+V}{1-(xU+V)} = \frac{U}{(1-V-U)(1-V+U)}
\end{align}

Finally the generating function of broken chains is that of all remaining chains.
\begin{equation} \label{BrokenChainsUN2OD}
\begin{aligned}
B &= \left( 3U+V \right) \sum_{n\geq 0} \left( 3U+V \right)^n - \sum\limits_{i=1}^{3} C_i = \frac{\left(3U+V\right)}{1-3U-V} - 2\frac{U}{1-U} - \frac{U+V}{1-U-V} \\
&= \frac{-6U^3 -8U^2V+6U^2-2UV^2 +4 UV}{\left(1-3U-V\right)\left(1-U\right)\left(1-U-V\right)}
\end{aligned}
\end{equation}
Similarly as for chains of color $3$, we could distinguish two types of broken chains depending on the coloring of the vertices at its boundary. However, that will not be necessary for our analysis.

%%%%%%%%%%%%%%%%%%%%
\subsection{Schemes of the quartic \texorpdfstring{$U(N)^2 \times O(D)$}{U(N)2xO(D)} multi-matrix model}

Let $\cG\in\bar{\mathbb{G}}$. The scheme of $\cG$ is obtained by 
\begin{enumerate}
\item Replacing every melonic 2-point subgraph with an edge of color 0,
\item Replacing every maximal chain with a chain-vertex of the same type,
\end{enumerate}
All graphs which reduce to the same scheme have the same genus $h$ and grade $l$, so it makes sense to define the genus and grade of a scheme as those values. We denote $\mathbb{S}_{h,l}$ the set of schemes at fixed values of $h, l$. We will prove the following result

\begin{theorem}
$\mathbb{S}_{h,l}$ is a finite set.
\end{theorem}

In other words, there is a finite number of schemes of fixed genus and grade. All graphs of fixed genus and grade $(h,l)$ can be obtained starting with a scheme in $\mathbb{S}_{h,l}$ and replacing chain-vertices with chains and edges of color 0 with melons.

%%%%%%%%%%%%%%%
\section{Finiteness of schemes of genus \texorpdfstring{$g$}{g} and grade \texorpdfstring{$l$}{l}}

Graphically, the only difference between the graphs of $\mathbb{G}$ and those of $\mathbb{G}_{O(N)^3}$ is the vertex coloring. There are two vertex colors in $\mathbb{G}$ because there are two fields $X, X^\dagger$, but no coloring in $\mathbb{G}_{O(N)^3}$ because $\phi$ is real. Hence there is a map 
\begin{equation}
\theta: \mathbb{G} \to \mathbb{G}_{O(N)^3}
\end{equation}
which simply consists in forgetting the vertex coloring.

\medskip

If $\mathcal{G}\in\mathbb{G}_{U(N)^2\times O(D)}$, we define the degree of $\cG$ as $\omega(\mathcal{G}) := \omega(\theta(\mathcal{G}))$ the degree of the graph it is mapped to in the $O(N)^3$ model. Then
\begin{equation}
\omega(\mathcal{G}) = h(\mathcal{G}) + \frac{l(\mathcal{G})}{2}
\end{equation}
which can be established by setting $D=N$ in \eqref{UN2ODAction}.

\medskip

The map $\theta$ maps melons to melons, dipoles to dipoles and chains to chains. We denote $\mathbb{S}_{O(N)^3}(\omega)$ be the set of schemes of degree $\omega$ in the $O(N)^3$ model. We can therefore descend the map $\theta$ to schemes by
\begin{equation}
\tilde{\theta}_{h,l}: \mathbb{S}_{h,l} \to \mathbb{S}_{O(N)^3}(h+l/2)
\end{equation}
from the schemes of genus $h$ and grade $l$ to those of the $O(N^3)$ model with degree $h+l/2$. It simply consists in forgetting the coloring of the vertices on the vertices, dipoles and chains.

\medskip

Theorem~\ref{th:sch} has been established for the quartic $O(N)^3$ tensor model (see Section~\ref{ssec:scheme_O(N)^3}). It is sufficient to show that each scheme $\mathcal{S}\in \mathbb{S}_{O(N)^3}(h+l/2)$ has a finite fiber via $\tilde{\theta}_{h,l}$. This is clear since the fiber of $\mathcal{S}\in \mathbb{S}_{O(N)^3}(\omega)$ is found by considering all colorings of the vertices and chain-vertices. This proves Theorem \ref{thm:FiniteSchemesI}.

\begin{figure}[!ht]
\begin{center}
\includegraphics[scale=0.55]{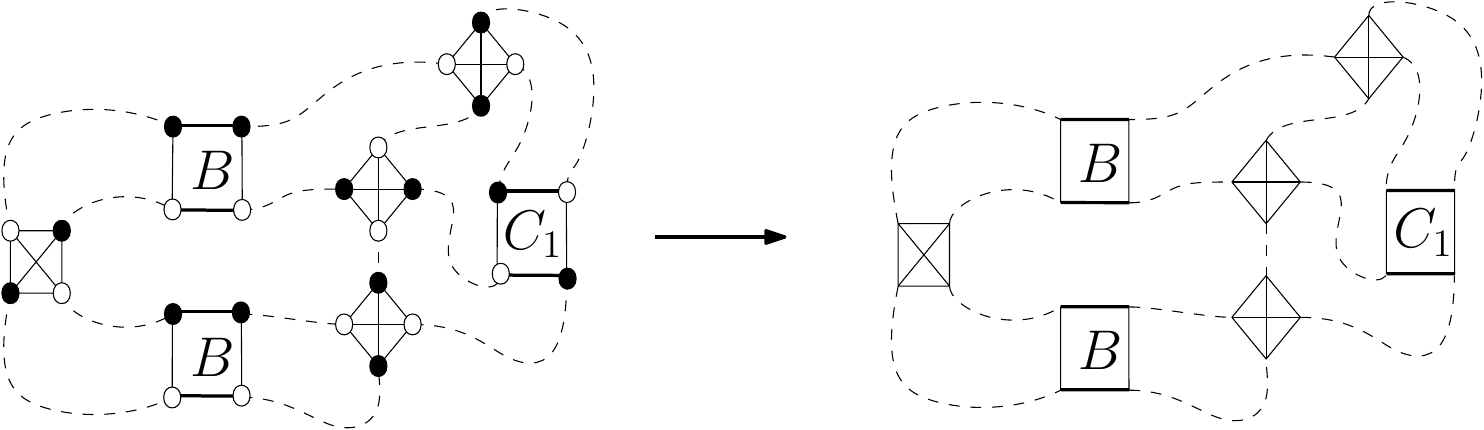}
\caption{A scheme of the $U(N)^2\times O(D)$-invariant model and its corresponding scheme on the side of the $O(N)^3$-invariant model.}
\end{center}
\end{figure}

For a fixed degree $\omega$, schemes of $\mathbb{S}_{O(N)^3}(\omega)$ have a bounded number of chains. Therefore, for a fixed value of $(h,l)$, the number of chains of schemes in $\mathbb{S}_{h,l}$ are also bounded. However, since the nature of the expansion is different in both models, the map $\tilde{\theta}$ cannot be used to identify directly dominant scheme for the $U(N)^2\times O(D)$ model from the dominant schemes of the $O(N)^3$ model. Indeed, in this case the double scaling limit selects graphs with vanishing grade while some dominant schemes of the $O(N)^3$ model have non-zero grade. As we will see later in section~\ref{sssec:id_dom_scheme}, the structure of the dominant schemes for the $U(N)^2\times O(D)$ model will be different from the one of Theorem~\ref{thm:DominantSchemesO(N)3} of the quartic $O(N)^3$ model.

%%%%%%%%%%%%%%%%%%
\section{Identification of the dominant schemes and double scaling limit}

%%%%%%%%%%%%%%%
\subsection{Singularity analysis}  \label{sec:UN2ODSingularity}

Looking at the form of the generating series above, there are four different types of possible singular points:
\begin{itemize}
\item Singular points of $M(t,\mu)$,
\item Points where $U(t,\mu)=1$, these points are singular for chains of colors 1 and 2, and for broken chains,
\item Points where $1-U(t,\mu)-V(t,\mu)=0$, which are singular for chains of color 3 and for broken chains,
\item Points such that $1-V(t,\mu)+U(t,\mu)=0$, which are singular for even and odd chains of color 3,
\item Points such that $1-3U(t,\mu)-V(t,\mu)=0$, which are singular for broken chains.
\end{itemize}
Since $M(t, \mu)$ is an increasing function of $t$, so are $U(t,\mu)$ and $V(t,\mu)$, which are also non-negative since $M(0,\mu) = 1$. Therefore, at fixed $\mu$, the singular points $1-3U(t,\mu)-V(t,\mu)=0$ are always encountered before any other point, except for potential singular points of $M(t,\mu)$. Noticing that $3U(t,\mu)+V(t,\mu) = tM(t,\mu)^4 + t\mu M(t,\mu)^2$, the analysis performed for the $O(N)^3$-invariant tensor model can be applied here. Thus, points where $3U(t,\mu)+V(t,\mu) = 1$ correspond to the locus of singular points of $M(t,\mu)$. These points are plotted on Figure~\ref{fig:plot_tc}.

\medskip 
Using Equation~\eqref{eq:crit_behav}, near a singular point $(t_c(\mu),\mu)$ of $M(t,\mu)$ we have:
{\small
\begin{align}
V(t,\mu) = \frac{1}{4}t_c(\mu)\mu \left( M_c(\mu) + 2M_c\mu)K(\mu)  \sqrt{1 - \frac{t}{t_c(\mu)}} \right) + \mathcal{O}( t_c(\mu) - t), \\
U(t,\mu) = M_c(\mu)-1 - \frac{3}{4} t_c(\mu) \mu M_c(\mu)^2 -\frac{3}{2} t_c(\mu) \mu M_c(\mu) K(\mu) \sqrt{1- \frac{t}{t_c(\mu)}} + \mathcal{O}( t_c(\mu) - t),
\end{align}}%
which leads to
{\small 
\begin{align}
B(t,\mu) &\underset{t\rightarrow t_c(\mu)}{\sim}  \frac{1}{\left(1-\frac{4}{3} t_c(\mu) \mu M_c(\mu) \right) \sqrt{1- \frac{t}{t_c(\mu)}}} \\ &\hspace{-20pt}\times \left[\frac{-6U(t,\mu)^3 -8U(t,\mu)^2V(t,\mu)+6U(t,\mu)^2-2U(t,\mu)V(t,\mu)^2 +4 U(t,\mu)V(t,\mu)}{\left(1-U(t,\mu)\right)\left(1-U(t,\mu)-V(t,\mu)\right)}\right]_{\big\rvert_{t_c(\mu)}} \nonumber
\end{align}}%

%%%%%%%%%%%%%%%
\subsection{Dominant schemes}
\label{sssec:dom_sch_U(N)2}

There are finitely many schemes for a given genus and grade $(h,l)$, thus all singularities of their generating function come from the generating series of melons and chains. As seen before, the leading singularity is that of broken chains, which we will show are in bounded number in a scheme of fixed genus and grade. Moreover, we restrict attention to schemes of vanishing grade as they are the ones dominating the large $N$, large $D$ expansion.

\medskip

Let us stress again that those dominant schemes cannot be found from the equivalent result for the $O(N)^3$ model. Indeed, since we have embedded our model into the $O(N)^3$ model and used that to prove the finiteness of the number of schemes at fixed genus and grade, it is a natural question to ask whether the latter model can also be used to identify the dominant schemes. However, dominant schemes of the $O(N)^3$ have sub-schemes of degree $\omega=1/2$ (which correspond to leaves from the tree in the representation explained below). A graph of the $U(N)^2\times O(D)$ model which is mapped to a graph of degree 1/2 of the $O(N)^3$ model has non-vanishing grade. It is thus not possible to immediately deduce the dominant schemes of the $U(N)^2\times O(D)$ model from those of the $O(N)^3$ model.

\medskip

The combinatorial analysis of dominant schemes for this model with tetrahedral interaction only has been done in~\cite{BeCa} and is similar to the analysis performed for the $O(N)^3$ tensor model before. It can be adapted to our situation with minor adjustments. We sketch the main points of the derivation.

\subsection{Removals of chains and dipoles}
\label{sssec:DipoleRemovals}

The chains and dipoles are at the heart of scheme decomposition. As such, chain or dipole removal helps understand the structure of schemes as it allows to decompose a scheme into a \textit{skeleton graph} which has no chains nor dipoles and turns out to be one of the key ingredients in the identification of dominant schemes. The removal of a dipole or chain-vertex is the following move
\begin{equation}
    \begin{array}{c}\includegraphics[scale=.5]{DipoleOrChainVertex.pdf}\end{array}
\end{equation}
A dipole-vertex is said to be separating if its removal disconnects $\cG$ into two graphs $\cG_1, \cG_2$ and \emph{non-separating} otherwise, and similarly for chain-vertices. For a quantity $\mathcal{O}$, let us denote $\Delta\mathcal{O} = \mathcal{O}(\cG_1) + \mathcal{O}(\cG_2) - \mathcal{O}(\cG)$ in the case of a separating removal and $\Delta \mathcal{O} = \mathcal{O}(\cG') - \mathcal{O}(\cG)$ in the case of a non-separating removal, with $\cG'$ being the (connected) graph obtained after the removal.

\paragraph{Separating dipole removals.} It can be checked that for all types of dipoles
\begin{equation}
    h(\cG) = h(\cG_1) + h(\cG_2),\qquad l(\cG) = l(\cG_1) + l(\cG_2)
\end{equation}
This is proved in \cite{BeCa} for tetrahedral interactions. Let us give an example in the case of a dipole pillow of color 1,
\begin{equation}
    \begin{array}{c} \includegraphics[scale=.3]{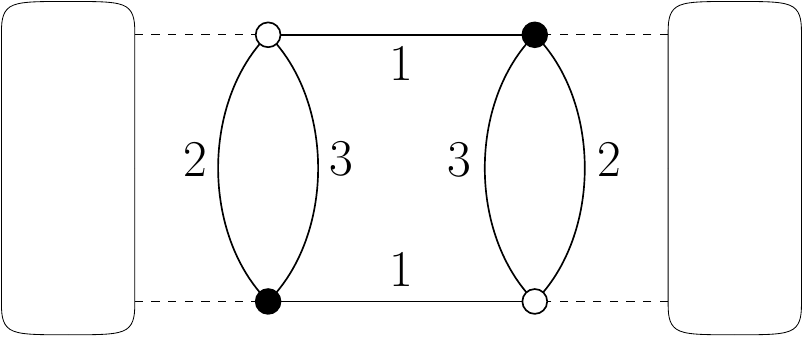} \end{array} \to \begin{array}{c} \includegraphics[scale=.3]{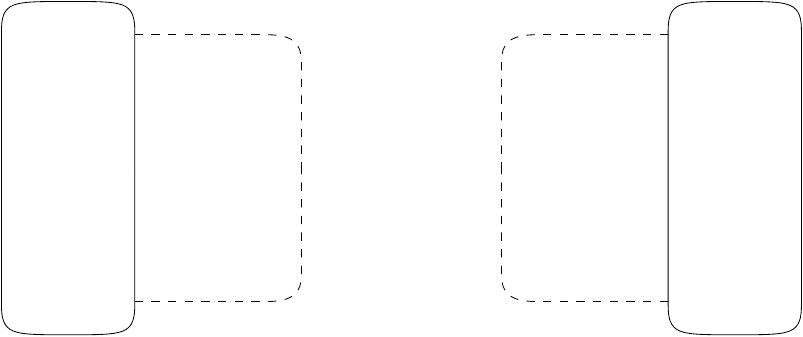} \end{array}
\end{equation}
Then, it is clear that 
\begin{equation} \label{VariationsPillowColor1Removal}
    \Delta F_{02} = \Delta F_{03} = 0, \Delta n_1=-1, \Delta E_0 = -2, \Delta F_{01} = 1.
\end{equation}
Then one finds $\Delta h = \Delta l =0$ from \eqref{GenusExpansion} and \eqref{Grade} (which hold for each connected component independently).

\paragraph{Non-separating dipole removals.}
\subparagraph{Non-separating dipoles of color $a\in\{1,2\}$.} Then
\begin{equation} \label{DipoleRemoval12}
    \Delta h =-1\text{ or }0, \qquad \Delta l = -2 \text{ or } -4.
\end{equation}
The case of dipoles built on tetrahedral interactions is proved in \cite{BeCa}. Let us consider the case of a dipole made of a pillow of color 1. Then, the variations are the same as in \eqref{VariationsPillowColor1Removal} except for $\Delta F_{01} \geq -1$ since $\cG'$ can have either one more or one less face of color 1 than $\cG$. One concludes by taking the variations in \eqref{GenusExpansion} and \eqref{Grade} again.

\subparagraph{Non-separating dipoles of color 3.} Then 
\begin{equation} \label{DipoleRemoval3}
    \Delta h =-1, \qquad \Delta l = 0 \text{ or } -4.
\end{equation}
The case of dipoles of color 3 built on tetrahedral interactions is proved in \cite{BeCa}, therefore we only have to consider dipoles made of a pillow of color 3. Then it is found that
\begin{equation}
    \Delta F_{01} = \Delta F_{02} = 0, \Delta F_{03} = \pm1, \Delta E_0 = -2, \Delta (n_{3b} + n_{3nb}) = -1,
\end{equation}
since $\cG'$ can have either one more or one less cycle of colors $\{0,3\}$. One concludes with \eqref{GenusExpansion} and \eqref{Grade} once again.

\paragraph{Chain removals.}

A (maximal) chain removal can be performed by removing a dipole from the chain, then removing the melonic 2-point functions this creates. The removal of a separating chain is obviously the same as that of a separating dipole.

\subparagraph{Non-separating chains of color $i$.}
Chains of color $i \in \{1,2,3\}$ are sequences of isolated dipoles of color $i$. Thus their removal is exactly the same as the removal of a non-separating dipole of color $i$.

\subparagraph{Non-separating broken chains.} By definition a broken chain has at least two dipoles of different colors. If it has no dipoles of color 3, then it must be a special case of \eqref{DipoleRemoval12} and if it has one it must also be a special case of \eqref{DipoleRemoval3}. By checking all cases, it is found that
\begin{equation} \label{BrokenChainRemoval}
    \Delta h=-1,\qquad \Delta l=-4.
\end{equation}
This is proved for tetrahedral interactions in \cite{BeCa}. Let us give an example with a pillow bubble of color 1 and a pillow bubble of color 3. Then $\Delta F_1=\Delta F_3=0$ and $\Delta F_2 =-1$. Moreover $\Delta E_0=-4$ and $\Delta n_1=\Delta (n_{3b}+n_{3nb})=-1$. One concludes with \eqref{GenusExpansion} and \eqref{Grade}.

\subsection{Skeleton graphs}

The definition of skeleton graphs is identical to the one of the $O(N)^3$ tensor model, see Def.~\ref{def:skel_graph}. Its yields similar properties as Lemma~\ref{thm:SkeletonGraph}, which is simply adapted to take into account the difference in the nature of the expansion in both model, i.e. replacing the degree by the pair genus and grade.

\begin{lemma}\label{thm:SkeletonGraph_U(N)2}
The skeleton graph $\mathcal{I}(\cS)$ of a scheme $S$ satisfies the following properties:
\begin{enumerate}
\item Any vertex of $\mathcal{I}(\cS)$ corresponding to a component of vanishing genus and grade, and not carrying the external legs of $\cS$, has degree at least $3$.
\item If $\mathcal{I}(\cS)$ is a tree, then the genus and grade of $\cS$ are split among the components.
\end{enumerate}
\end{lemma}

\subsubsection{Identifying the dominant schemes}
\label{sssec:id_dom_scheme}

We can now adapt the proof of Theorem~\ref{thm:DominantSchemesO(N)3} to this model, allowing us to obtain the combinatorial structure of the dominant scheme. We first prove that if $\mathcal{S}$ is a dominant scheme, its skeleton graph $\mathcal{I}(\mathcal{S})$ is a tree. Assume that $\mathcal{I}(\mathcal{S})$ is a not a tree, i.e. $\mathcal{S}$ has a non-separating chain.

\begin{itemize}
    \item If it has a non-separating broken chain-vertex, then removing it decreases the grade, see Equation \eqref{BrokenChainRemoval}, so $\mathcal{S}$ cannot be dominant. The same holds if it has a non-separating chain-vertex of color $a\in\{1,2\}$, see Equation \eqref{DipoleRemoval12}.
    \item If it has a non-separating chain or dipole of color 3, there two possibilities. Either the grade decreases, then $\mathcal{S}$ cannot be dominant, or it does not but the genus decreases anyway, see Equation \eqref{DipoleRemoval3}. We then perform the following move
    \begin{equation}
        \begin{array}{c}\includegraphics[scale=.4]{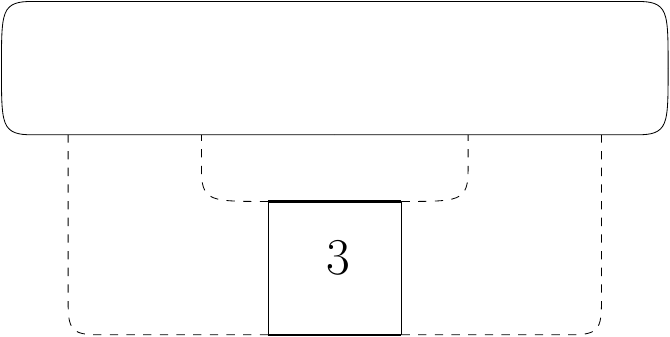}\end{array} \quad \to \quad \begin{array}{c}\includegraphics[scale=.4]{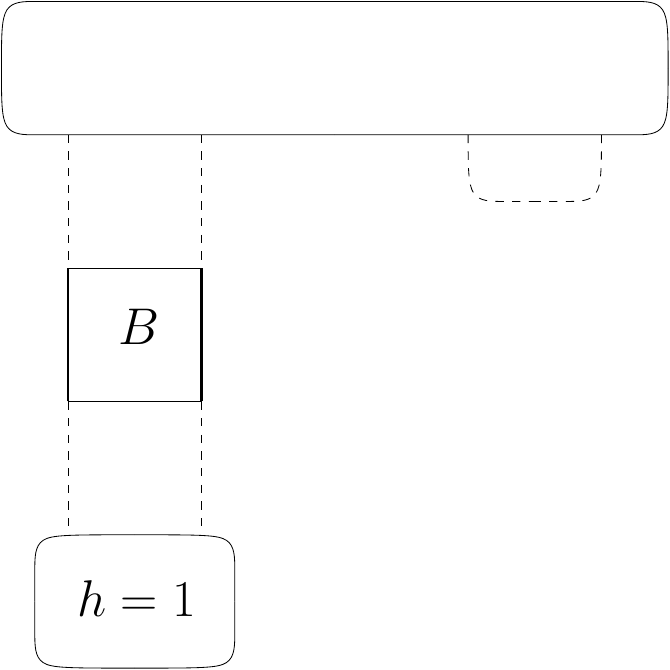}\end{array}
    \end{equation}
    where the component which is added has genus 1. On the RHS we get a scheme with the same genus as on the LHS but one more broken chain, hence it is more singular and as a result $\mathcal{S}$ cannot be dominant.
\end{itemize}

Thus, the skeleton graph must be a tree $\mathcal{T}$. Using Lemma~\ref{lemma:fst}, we know that the leaves of $\mathcal{T}$ cannot have vanishing genus and grade, therefore all leaves of $\mathcal{T}$ must have genus at least $1$ (and grade 0), hence there are at most $g$ leaves in $\mathcal{T}$. Finally, in order to be dominant, $\mathcal{T}$ must have as many broken chains as possible. We thus have an optimization problem: maximizing the number of chains with an upper bound on the number of leaves, and involving as variables the degrees and the genus of the internal vertices of $\mathcal{T}$. The solution is to have as many leaves as possible, here $g$, and inner vertices of degree exactly 3 so that $\mathcal{T}$ is a plane binary tree.

\medskip

In terms of schemes, the leaves of $\mathcal{T}$ correspond to 2-point schemes of genus 1 (and vanishing grade) with no separating chains, while the internal vertices of $\mathcal{T}$ correspond to 6-point functions of vanishing genus (and grade). By repeating the analysis of~\cite{BeCa}, one finds the structure of those objects. This gives the following proposition.

\begin{proposition} 
A dominant scheme of genus $h>0$ has $2h-1$ broken chain-vertices, all separating. Such a scheme has the structure of a rooted binary plane tree where
\begin{itemize}
\item Edges correspond to broken chain-vertices.
\item The root of the tree corresponds to the two external legs of the 2-point function.
\item The $h$ leaves are one of the two following graphs
\begin{equation}
\includegraphics[scale=0.5]{leaf_tree_un2.pdf}
\end{equation}
Observe that the chains of color $3$ have different boundary vertices in the two graphs.
\item Each internal vertex corresponds to one of four $6$-point subgraphs:
\begin{equation}
\includegraphics[scale=0.5]{node_tree_un2.pdf}
\end{equation}
Here we have used the embedding convention that the edge coming out of the bottom left is connected to a white vertex.
\end{itemize}
\end{proposition}

Notice that the two types of chains of color 3, which differ from their vertex colorings, play different roles (at the leaves of the tree). Remarkably, it is not necessary to also distinguish two types of broken chains according to their vertex colorings. Indeed, the only constraint imposed by the above Proposition is that every broken chain (corresponding to the edges of the tree) has one black and one white vertex on each side.

\medskip

We now give the generating series of graphs associated to a dominant scheme. Recall that if $\mathcal{T}$ is a rooted  binary tree with $n$ edges, then it has $\frac{n-1}{2}$ internal vertices, $\frac{n+1}{2}$ leaves (not counting the root).

\medskip

For a dominant scheme corresponding to a rooted plane binary tree $\mathcal{T}$, every edge of $\mathcal{T}$ contributes with the generating series of broken chains, every inner node with $(1+6t)$ (1 for the leftmost case of internal nodes of Proposition \ref{prop:dom_scheme_un2}, all 6 cases come with a factor $t$). Since at genus $h>0$, a dominant scheme has are $2h-1$ edges, its generating series is
\begin{align}
G_{\mathcal{T}}^{h}(t,\mu) = L^{2-2h}\bigl(C_{3,\hspace{1pt}\includegraphics[scale=0.2,valign=c]{wbbw.pdf}}(t,\mu)+tC_{3,\hspace{1pt}\includegraphics[scale=0.2,valign=c]{wwbb.pdf}}(t,\mu)\bigr)^{h}\left(1+6t\right)^{h-1}B(t,\mu)^{2h-1}
\end{align}
where $B(t,\mu)$ is given by Equation~\eqref{BrokenChainsUN2OD}. We recall that $L=N/\sqrt{D}$ is the parameter of the genus expansion, see Theorem \ref{thm:free_energy}.

\medskip

Since it only depends on $h$, the sum over all rooted binary trees with $h$ leaves (not counting the root) is 
{\small 
\begin{equation}
G_{\text{dom}}^{h}(t,\mu) = \sum_{\substack{\mathcal{T}\\ \text{$h$ leaves}}} G_{\mathcal{T}}^{h}(t,\mu) = L^{2-2h}\frac{\Cat_{h-1}}{(1+6t) B(t,\mu)} \Bigl(\bigl(C_{3,\hspace{1pt}\includegraphics[scale=0.2,valign=c]{wbbw.pdf}}+tC_{3,\hspace{1pt}\includegraphics[scale=0.2,valign=c]{wwbb.pdf}}\bigr)(1+6t)B(t,\mu)^{2}\Bigr)^h.
\end{equation}}%
Here $\Cat_{h-1} = \frac{1}{h}\binom{2h-2}{h-1}$ the number of rooted trees with $h$ leaves. $G_{\text{dom}}^{h}(t,\mu)$ is the total contribution of dominant schemes of genus $h$.

\subsection{Double scaling limit for the 2-point function}

As higher order graphs in the genus expansion can have more broken chains, we can tune the way we approach the critical point such that the singular behaviour of broken chains makes up for the loss of scaling in $L$ of the graphs, enhancing their contribution. Thus the double scaling parameter $\kappa(\mu)$ is defined graphs of all genus contribute. The dominant schemes of genus $h$ scale with $L$ as $L^{2-2h}$, therefore we defined $\kappa(\mu)$ as
\begin{equation}
\kappa(\mu)^{-1} = L^2 \frac{1}{(1+6t)(C_{3,\hspace{1pt}\includegraphics[scale=0.2,valign=c]{wbbw.pdf}}(t,\mu)+tC_{3,\hspace{1pt}\includegraphics[scale=0.2,valign=c]{wwbb.pdf}}(t,\mu))}\left(\frac{1}{B(t,\mu)}\right)^2.
\label{eq:kappa_un2A}
\end{equation}
Using Equation~\eqref{eq:kappa_un2A} we have
\begin{equation}
 \frac{1}{(1 + 6t)B(t,\mu)}  = \frac{\kappa(\mu)^{-\frac{1}{2}}}{L} \left(\frac{C_{3,\hspace{1pt}\includegraphics[scale=0.2,valign=c]{wbbw.pdf}}+tC_{3,\hspace{1pt}\includegraphics[scale=0.2,valign=c]{wwbb.pdf}}}{1+6t}\right)^{\frac{1}{2}}
\end{equation}

Therefore, the contribution of the dominant scheme of genus $h>0$ in the double scaling limit reads
\begin{equation}
G_{dom}^{h}(t,\mu) = L\Cat_{h}\left(\frac{(C_{3,\hspace{1pt}\includegraphics[scale=0.2,valign=c]{wbbw.pdf}}+tC_{3,\hspace{1pt}\includegraphics[scale=0.2,valign=c]{wwbb.pdf}})}{1+3t}\bigg\rvert_{(t_c(\mu),\mu)} \right)^{\frac{1}{2}}\kappa(\mu)^{h-\frac{1}{2}}.
\end{equation}

For $h=0$, there is only one scheme with one edge and no vertices, corresponding to the generating function of melonic graphs. It contributes as $M_c(t_c(\mu),\mu)$.

\medskip

The leading order of the $2$-point function in the double scaling limit is obtained by resumming the contribution of all genus $h$. Therefore we get
\begin{align}
G_{2}^{DS}(\mu) &= \sum\limits_{h\in\mathbb{N}} G_{dom}^h(\mu) \\
&= M_c(t_c(\mu),\mu)\left(1+ L \left(\frac{(C_{3,\hspace{1pt}\includegraphics[scale=0.2,valign=c]{wbbw.pdf}}+tC_{3,\hspace{1pt}\includegraphics[scale=0.2,valign=c]{wwbb.pdf}})}{1+6t}\bigg\rvert_{(t_c(\mu),\mu)} \right)^{\frac{1}{2}} \frac{1-\sqrt{1-4\kappa(\mu)}}{2\kappa(\mu)^\frac{1}{2}}\right) \nonumber
\end{align}
The sum converges for $\kappa(\mu) \leq \frac{1}{4}$, similarly to other tensor models. The parameter $\kappa(\mu)$ encodes a balance between the large $N$ limit and the criticality (in the large $D$ limit). In particular, for $\kappa(\mu) = 0$, we obtain the usual large $N,D$ limit where only melons contribute, and at the other end when $\kappa(\mu) = \frac{1}{4}$ analyticity is lost. 

\medskip

Note however that the scaling with $L$ differs from that of~\ref{eq:GDS_ON3}. This is due to differences in the dominant schemes between the two models (although both are indeed mapped to rooted binary plane trees). In the $U(N)^2 \times O(D)$-invariant multi-matrix model, the leaves of the trees associated to a dominant scheme all have genus $1$ and vanishing grade. However, they have degree $\omega = \frac{1}{2}$ in the $O(N)^3$-invariant tensor model. Thus, the dominant schemes of the $U(N)^2 \times O(D)$-invariant multi-matrix models are not mapped to the dominant schemes of the $O(N)^3$ tensor model by the map $\theta$ defined in Equation~\eqref{eq:theta}.

%% file: Appendices/App_prop_3j.tex
\chapter{Some identities from \texorpdfstring{$SU(2)$}{SU(2)} recoupling theory}
\label{app:su2recoupling}

This Appendix contains some definitions and properties related to $SU(2)$ recoupling theory used in Chapter~\ref{Chap:melons_th}. These properties are classical results on the domain and we will simply state them here. Their proof can be found in the litterature on the topic e.g. Ilkka M\"akinen's introduction~\cite{Mak19}. 

\paragraph{Haar measure and Wigner matrices\\}

From the Peter-Weyl theorem, the Wigner matrices $D^j_{mn}(g)$ form an orthogonal basis of the functions $f:SU(2) \rightarrow \mathbb{C}$. This orthogonality relation is encoded at the level of the Haar measure via the relation
\begin{equation}
      \int d g D^j_{mn}(g)\bar{D}^{j'}_{m'n'}(g) =\frac{1}{(2j+1)}\delta^{jj'}\delta_{mm'}\delta_{nn'},
\end{equation}
where the Wigner matrices satisfy
\begin{equation}
      D^j_{mn}(g)=(-1)^{m-n}\bar{D}^j_{-m,-n}(g).
\end{equation}

\paragraph{The \texorpdfstring{$3j$}{3j}-symbol and its properties\\}

The $3j$ symbol is invariant under the action of $SU(2)$ group, 
\begin{equation}
D^{j_1}_{m_1n_1}(g_1)D^{j_2}_{m_2n_2}(g_2)D^{j_3}_{m_3n_3}(g_3)\TJ{j_1}{j_2}{j_3}{n_1}{n_2}{n_3}=\TJ{j_1}{j_2}{j_3}{m_1}{m_2}{m_3}.
\end{equation}

It is also invariant under the even permutations of indices, and it acquires an additional phase under odd permutations
\begin{equation}
       \TJ{j_1}{j_2}{j_3}{m_1}{m_2}{m_3}=(-1)^{j_1+j_2+j_3}\TJ{j_2}{j_1}{j_3}{m_2}{m_1}{m_3}. \label{eq:3jpermu}
\end{equation}

Similar relations holds when changing the sign of all magnetic indices simultanesously.
\begin{equation}
       \TJ{j_1}{j_2}{j_3}{-m_1}{-m_2}{-m_3}=(-1)^{j_1+j_2+j_3}\TJ{j_1}{j_2}{j_3}{m_1}{m_2}{m_3}. 
\end{equation}

The $3j$ symbols satisfy two orthonormal relations
\begin{align}
     &(2j_3+1)\sum_{m_1,m_2}\TJ{j_1}{j_2}{j_3}{m_1}{m_2}{m_3}\TJ{j_1}{j_2}{j'_3}{m_1}{m_2}{m'_3}=\delta_{j_3,j_3'}\delta_{m_3,m_3'}, \\
  &\sum_{j_3,m_3}(2j_3+1)\TJ{j_1}{j_2}{j_3}{m_1}{m_2}{m_3}\TJ{j_1}{j_2}{j_3}{m'_1}{m'_2}{m_3}=\delta_{m_1,m_2'}\delta_{m_2,m_2'},  \label{eq:3jsumdelta}
\end{align}

Finally, when one of the magnetic moment (say $m_3$) vanishes, then the $3j$ symbol vanishes unless $m_1 = -m_2$ and we have
\begin{equation}
       \sum_m (-1)^{j-m}\TJ{j}{j}{k}{m}{-m}{0}=\sqrt{2j+1}\delta_{k,0}. \label{eq:3jm3eq0}
\end{equation}

In particular for $k=0$ it gives
\begin{equation}
      \TJ{j_1}{0}{j_3}{m_1}{0}{m_3}=\delta^{j_1,j_3}\frac{1}{\sqrt{2j_1+1}}(-1)^{j_1+m_1}\delta_{m_1,-m_3} \label{eq:3jj2eq0} \\
\end{equation}

\paragraph{The \texorpdfstring{$6j$}{6j}-symbol and its properties\\}

The $6j$ symbol is defined as
\begin{align}
      \SJ{j_1}{j_2}{j_3}{j_4}{j_5}{j_6} =&\sum_{j_i,m_i} (-1)^{\sum_{a=1}^6(j_a-m_a)} \TJ{j_1}{j_2}{j_3}{-m_1}{-m_2}{-m_3}\TJ{j_1}{j_5}{j_6}{m_1}{-m_5}{m_6} \nonumber \\
  &\cdot \TJ{j_4}{j_2}{j_6}{m_4}{m_2}{-m_6}\TJ{j_4}{j_5}{j_3}{-m_4}{m_5}{m_3}.   \label{eq:6jdef}
\end{align}

It satisfies the following relation
\begin{align}
       &\sum_{n_1,n_2,n_3}(-1)^{\sum_{a=1}^{3}(k_a-n_a)} \TJ{j_1}{k_2}{k_3}{m_1}{-n_2}{n_3}\TJ{k_1}{j_2}{k_3}{n_1}{m_2}{-n_3}\TJ{k_1}{k_2}{j_3}{-n_1}{n_2}{m_3} \nonumber \\
   &=\SJ{j_1}{j_2}{j_3}{k_1}{k_2}{k_3} \TJ{j_1}{j_2}{j_3}{m_1}{m_2}{m_3}. \label{eq:3jsum6j} 
\end{align}

Finally when one of the spin index, say $j_6$, vanishes we have
\begin{equation}
     \SJ{j_1}{j_2}{j_3}{j_4}{j_5}{0}=\frac{\delta_{j_1,j_5}\delta_{j_2,j_4}}{\sqrt{d_{j_1}d_{j_2}}}(-1)^{j_1+j_2+j_3}\{j_1~j_2~j_3\}. \label{eq:6jj6eq0}
\end{equation}
where $\{j_1~j_2~j_3\}$ is one if $(j_1,j_2,j_3)$ satisfy the triangular inequality and zero otherwise.